 \documentclass[reqno,table]{amsart}
\usepackage[table]{xcolor}
\usepackage{array}
\usepackage{booktabs}
\usepackage{amssymb,amsmath,amsfonts,amsthm,epsfig,amscd,stmaryrd}
 \usepackage[textsize=tiny]{todonotes}
 \usepackage{pifont}
 \usepackage{bbding}
\usepackage{url}

\newcommand{\Zptimes}{\Z_{p}^{\times}}
\newcommand{\Qptimes}{\Q_{p}^{\times}}
\newcommand{\Sen}{\operatorname{Sen}}
\newcommand{\cycl}{\operatorname{cycl}}
\newcommand{\arit}{\operatorname{arit}}
\newcommand{\Mmw}{M(\mf{m}_w)} 
\newcommand{\Lmw}{L(\mf{m}_w)}
\newcommand{\CwQ}{C_{w,Q}}
\newcommand{\CwQdag}{C_{w,Q}^\dag}
\newcommand{\StabQ}{\Stab_{Q}}
\newcommand{\StabQw}{\Stab_{Q}(w)}
 \renewcommand{\setminus}{\smallsetminus}
\chardef\bslash=`\\

\newcommand{\Adm}{\operatorname{Adm}}
\newcommand{\Loc}{\operatorname{Loc}}
\newcommand{\proket}{\operatorname{prok\'et}}

\newcommand{\cha}{\operatorname{char}}
\newcommand{\tv}{{\widetilde{v}}}
\newcommand{\ilim}{\mathop{\varinjlim}\limits}
\newcommand{\plim}{\mathop{\varprojlim}\limits}

\newcommand{\disjointunion}{\sqcup}

\newcommand{\Qtwo}{\Q_2}
\newcommand{\cCbar}{\overline{\cC}}
\newcommand{\Hbar}{\overline{H}}

\newcommand{\Gbar}{\overline{G}}
\newcommand{\corank}{\operatorname{corank}}
\newcommand{\patch}{\operatorname{patch}}

\newcommand{\cotimes}{\widehat{\otimes}}
\newcommand{\Pone}{\mathbb{P}^1}

\newcommand{\vzero}{r}

\def\top{\mathrm{top}}
\def\alg{\mathrm{alg}}
\def\dalg{\text{-}\mathrm{alg}}
\def\ev{\mathrm{ev}}
\def\coad{\mathrm{coad}}
\def\fing{\mathrm{fg}}
\def\la{\mathrm{la}}
\def\Stab{\mathrm{Stab}}
\def\sm{\mathrm{sm}}
\def\tor{\mathrm{tor}}
\def\Rla{\mathrm{Rla}}
\def\nfs{f\kern-0.06em{s}}
\def\Sh{\mathrm{Sh}}
\def\hor{\mathrm{hor}}
\def\VB{VB}
\def\VBzero{VB^0}
\def\VBred{VB^{\red}}
\def\LB{LB}
\def\LS{LS}
\def\QQQ{\mathcal{Q}}

\def\Wbar{\overline{W}}
\def\Alt{\mathrm{Alt}}
\def\WW{\mathcal{W}}
  
\newcommand{\et}{\mathrm{\acute{e}t}}

\def\aaaa{\psi}
\def\bbbb{b}
\def\taaaa{\widetilde{\aaaa}}
\def\aaaabar{\overline{\aaaa}}
\newcommand{\xmark}{\ding{55}}
\def\Vbar{\overline{V}}
\def\nubar{\overline{\nu}}
\def\FFFF{\xmark}
\def\FFFFF{\FFFF}
\def\TTTT{\checkmark}

\def\MMM{\mathcal{M}}

\def\franksp{\mathfrak{sp}}
\def\frankgsp{\mathfrak{gsp}}

\def\gl{\mathfrak{gl}}
\def\sl{\mathfrak{sl}}

\def\sss{\mathrm{ss}}

\newcommand{\TW}{\mathrm{TW}}

\newcommand{\chit}{\widetilde{\chi}}

\newcommand{\Lift}{\mathrm{Lift}}
\newcommand{\Def}{\mathrm{Def}}

\newcommand\llb{\llbracket}
\newcommand\rrb{\rrbracket}

 \def\Ghat{\widehat{G}}

\usepackage[backref=page,bookmarksopen,bookmarksdepth=2]{hyperref}
\renewcommand\j{\Don't use this because it conflicts with another style file}
\usepackage{etex} 
\usepackage{phoenician}
\usepackage{linearb}
\usepackage{linearA}
\usepackage{scalerel}
\usepackage{graphicx}

\usepackage[utf8]{inputenc}
\usepackage[T1]{fontenc}
\usepackage{fixltx2e}
\usepackage{mathrsfs}

\def\AFF{AF\kern-0.12em{F}}
\def\KZ{Z\kern-0.1em{K}}

\usepackage[OT2,T1]{fontenc}
\DeclareSymbolFont{cyrletters}{OT2}{wncyr}{m}{n}
\DeclareMathSymbol{\Sha}{\mathalpha}{cyrletters}{"58}

\def\HH{H}

\renewcommand{\star}{*}
\def\RGamma{\mathrm{R}\Gamma}
\def\varrhobar{\overline{\varrho}}

\def\PSp{\mathrm{PSp}}
\def\PGSp{\mathrm{PGSp}}

\newcommand{\oscr}{\mathcal{O}}

\newcommand{\an}{\textrm{an}}

\newcommand{\cGhat}{\widehat{\cG}}

\def\doubleslash{/\kern-0.2em{/}}

\def\tw{\widetilde{w}}

\def\sl{\mathfrak{sl}}

\def\af{a \kern-0.05em{f}}

\def\cusp{\mathrm{cusp}}

\def\OL{\mathcal{O}}

\def\act{\mathrm{act}}

\newcommand{\Jac}{\operatorname{Jac}}
 
\newcommand{\Iw}{\mathrm{Iw}}

\newcommand{\Par}{\mathrm{Par}}
\newcommand{\Lambdat}{\widetilde{\Lambda}}

\newcommand{\Spa}{\mathrm{Spa}}

\renewcommand{\mathbb}{\mathbf}

\newcommand{\alphabeta}{\scaleobj{0.7}{\textlinb{\LinearAC}}}
\newcommand{\alphabar}{\overline{\alpha}}
\newcommand{\betabar}{\overline{\beta}}
\newcommand{\alphabetabar}{\overline{\alphabeta}}
\newcommand{\recGT}{\rec_{\operatorname{GT}}}
\newcommand{\recGTp}{\rec_{\operatorname{GT},p}}

\newcommand{\chibar}{\overline{\chi}}
\newcommand{\rbar}{\overline{r}}

\newcommand{\gr}{\operatorname{gr}}

\usepackage{tikz}
\usepackage{tikz-cd}
\usetikzlibrary{arrows, matrix}
\usetikzlibrary{arrows.meta,
                positioning}

\usepackage[shortlabels]{enumitem}\usetikzlibrary{decorations.pathmorphing}

\newcommand{\Lie}{{\operatorname{Lie}\,}}
\newcommand{\semis}{{\operatorname{ss}}}

\newcommand{\fg}{{\mathfrak{g}}}

\newcommand{\mf}{\mathfrak}

\newcommand{\rec}{{\operatorname{rec}}}

\newcommand{\T}{\mathbf{T}}

\def\PGL{\mathrm{PGL}}

\def\tS{\widetilde{S}}

\usepackage[all,cmtip,poly]{xy}
\usepackage{color}

\newcommand{\loc}{\operatorname{loc}}
\newcommand{\ad}{\operatorname{ad}}
\newcommand{\diag}{\operatorname{diag}}
\newcommand{\tr}{\operatorname{tr}}
\newcommand{\Tr}{\operatorname{Tr}}

\newcommand{\CNL}{\operatorname{CNL}}
\newcommand{\Sets}{\operatorname{Sets}}

\newcommand{\red}{\operatorname{red}}
\newcommand{\colim}{\operatorname{colim}}

\newcommand{\To}{\longrightarrow}
\newcommand{\isoto}{\stackrel{\sim}{\To}}

\newcommand{\pr}{\operatorname{pr}}

\newcommand{\Rep}{\operatorname{Rep}}

\newcommand{\BT}{\operatorname{BT}}

\def\iso{\buildrel \sim \over \longrightarrow}

\newcommand{\RHom}{{\mathrm{RHom}}}

\makeatletter
\let\c@figure\c@equation
\let\c@table\c@equation
\makeatother

\def\numfigure{\addtocounter{subsubsection}{1}\begin{figure}}
\def\numtable{\addtocounter{subsubsection}{1}\begin{table}}

\newtheorem{theorem}[subsubsection]{Theorem}

\newtheorem{claim}[subsubsection]{Claim}
\newtheorem{thm}[subsubsection]{Theorem}
\newtheorem{lemma}[subsubsection]{Lemma}
\newtheorem{lem}[subsubsection]{Lemma}

\newtheorem{ithm}{Theorem}

\newtheorem{cor}[subsubsection]{Corollary}

\newtheorem{prop}[subsubsection]{Proposition}

\theoremstyle{definition}
\newtheorem{df}[subsubsection]{Definition}
\newtheorem{defn}[subsubsection]{Definition}

\newtheorem{remark}[subsubsection]{Remark}
\newtheorem{rem}[subsubsection]{Remark}

\newtheorem{example}[subsubsection]{Example}

\newtheorem{hypothesis}[subsubsection]{Hypothesis}

\def\numequation{\addtocounter{subsubsection}{1}\begin{equation}}
\def\numalign{\addtocounter{subsubsection}{1}\begin{align}}
\def\nummultline{\addtocounter{subsubsection}{1}\begin{multline}}

\newif\iffinalrunBANG
\iffinalrunBANG
\else
\fi

\iffinalrunBANG
  \newcommand{\need}[1]{}
  \newcommand{\mar}[1]{}
\else
  \newcounter{margcount}
  \newcommand{\need}[1]{{\tiny *** #1}}
  \newcommand{\mar}[1]{\stepcounter{margcount}\marginpar{\raggedright\tiny \themargcount.FIXME #1}}
\fi

\newcommand{\A}{\mathbf{A}}
\newcommand{\C}{\CC}
\newcommand{\F}{\FF}
\newcommand{\N}{\mathbf{N}}
\newcommand{\Q}{\QQ}
\newcommand{\qq}{\QQ}
\newcommand{\R}{\RR}
\newcommand{\Z}{\ZZ}

\renewcommand{\O}{\cO}

\newcommand{\m}{\frakm}

\newcommand{\q}{\frakq}

\newcommand{\CC}{{\mathbb C}}
\newcommand{\FF}{{\mathbb F}}
\newcommand{\GG}{{\mathbb G}}
\newcommand{\RR}{{\mathbb R}}
\newcommand{\QQ}{{\mathbb Q}}
\newcommand{\TT}{{\mathbb T}}

\newcommand{\ZZ}{{\mathbb Z}}

\newcommand{\cA}{{\mathcal A}}

\newcommand{\cC}{{\mathcal C}}
\newcommand{\cD}{{\mathcal D}}

\newcommand{\cF}{{\mathcal F}}
\newcommand{\cG}{{\mathcal G}}
\newcommand{\cH}{{\mathcal H}}

\newcommand{\cL}{{\mathcal L}}

\newcommand{\cO}{{\mathcal O}}
\newcommand{\ocal}{{\mathcal O}}

\newcommand{\cQ}{{\mathcal Q}}

\newcommand{\cS}{{\mathcal S}}
\newcommand{\cT}{{\mathcal T}}
\newcommand{\cU}{{\mathcal U}}

\newcommand{\cX}{{\mathcal X}}

\newcommand{\cZ}{{\mathcal Z}}
\newcommand{\frakm}{\mathfrak{m}}

\newcommand{\frakq}{\mathfrak{q}}

\newcommand{\Fbar}{\overline{\F}}
\newcommand{\Qbar}{\overline{\Q}}
\newcommand{\Zbar}{\overline{\Z}}

\newcommand{\Fp}{\F_p}

\newcommand{\Fpbar}{\Fbar_p}
\newcommand{\Ftwobar}{\Fbar_2}
\newcommand{\Qtwobar}{\overline{\Q}_2}

\newcommand{\Fpbartimes}{\Fpbar^{\times}}

\newcommand{\Zp}{\Z_p}

\newcommand{\Zpbar}{\Zbar_p}

\newcommand{\WM}{{}^{M}W}
\newcommand{\Ql}{\Q_{l}}
\newcommand{\Qp}{\Q_p}

\newcommand{\Qpbar}{\Qbar_p}

\newcommand{\Qpbartimes}{\Qpbar^{\times}}

\newcommand{\Cp}{\C_p}
\newcommand{\kbar}{\overline{k}}

\DeclareMathOperator{\Aut}{Aut}

\DeclareMathOperator{\depth}{depth}
\DeclareMathOperator{\disc}{disc}
\DeclareMathOperator{\End}{End}
\DeclareMathOperator{\Ext}{Ext}
\DeclareMathOperator{\Fil}{Fil}
\DeclareMathOperator{\Gal}{Gal}
\DeclareMathOperator{\GL}{GL}

\DeclareMathOperator{\GSp}{GSp}
\DeclareMathOperator{\GSpin}{GSpin}

\DeclareMathOperator{\Hom}{Hom}
\DeclareMathOperator{\im}{im}
\DeclareMathOperator{\Ind}{Ind}

\DeclareMathOperator{\Mod}{Mod}

\DeclareMathOperator{\ord}{ord}

\DeclareMathOperator{\SL}{SL}

\DeclareMathOperator{\Sp}{Sp}
\DeclareMathOperator{\Spec}{Spec}
\DeclareMathOperator{\AnSpec}{AnSpec}

\DeclareMathOperator{\Spf}{Spf}

\DeclareMathOperator{\Sym}{Sym}

\DeclareMathOperator{\WD}{WD}

\newcommand{\ab}{\mathrm{ab}}

\newcommand{\Frob}{\mathrm{Frob}}

\newcommand{\HT}{\mathrm{HT}}
\newcommand{\Id}{\mathrm{Id}}

\newcommand{\rhobar}{\overline{\rho}}

                               \newcommand{\into}{\hookrightarrow}
\newcommand{\onto}{\twoheadrightarrow}

\newcommand{\Gm}{\GG_m}

\newcommand{\etabar}{\overline\eta}

\newcommand{\Art}{{\operatorname{Art}}}

\newcommand{\varepsilonbar}{\overline{\varepsilon}}

\newcommand{\Res}{\operatorname{Res}}

\title[Modularity theorems for abelian surfaces]{Modularity theorems for abelian surfaces}%

\author[G.~Boxer]{George Boxer}  \email{g.boxer@imperial.ac.uk} \address{Department of
  Mathematics, Imperial College London,
  London SW7 2AZ,~UK}
  
\author[F.~Calegari]{Frank Calegari}  \email{fcale@math.uchicago.edu} \address{The University of Chicago,
5734 S University Ave,
Chicago, IL 60637, USA}

\author[T.~Gee]{Toby Gee} \email{toby.gee@imperial.ac.uk} \address{Department of
  Mathematics, Imperial College London,
  London SW7 2AZ,~UK}
  
\author[V.~Pilloni]{Vincent Pilloni}\email{vincent.pilloni@universite-paris-saclay.fr}
\address{D\'{e}partement de Math\'{e}matique, Universit\'{e} Paris-Sud, B\^atiment 307,
81405 Orsay cedex, France}

\thanks{G.B. was supported by a Royal Society University Research Fellowship.
F.C. was supported in part by NSF Grant DMS-2001097. T.G.\ was
  supported in part by an ERC Advanced grant and by the Simons Collaboration on Perfection in Algebra, Geometry, and Topology.  %
  This project has received funding from the European Research Council (ERC) under the European Union’s Horizon 2020 research and innovation programme (grant agreement No. 884596). V.P was supported in part by  the ERC-2018-COG-818856-HiCoShiVa. }
  
\begin{document}

\begin{abstract}We prove the modularity of a positive proportion of
  abelian surfaces over~$\Q$. More precisely, we prove the modularity
  of abelian surfaces which are ordinary at~$3$ and are
  $3$-distinguished, subject to some assumptions on the $3$-torsion
  representation (a ``big image'' hypothesis, and a technical
  hypothesis on the action of a decomposition group at~$2$). %
  We employ a 2--3 switch and a new
  classicality theorem (in the style of Lue Pan) for ordinary $p$-adic Siegel modular forms. 
 \end{abstract}

 \maketitle

\setcounter{tocdepth}{2}
{\footnotesize
\tableofcontents
}

\section{Introduction}

\subsection{The main theorems}

Our main theorem is as follows (see~\S\ref{sec:proofsofAandB}).

\begin{ithm} \label{first}
Let~$A/\Q$ be an abelian surface with a polarization of degree prime to~$3$. Suppose 
the following  holds:
\begin{enumerate}
\item \label{surjective} The mod~$3$ representation:
$$\rhobar_{A,3}: \Gal(\Qbar/\Q) \rightarrow \GSp_4(\F_3)$$
is surjective.
\item  \label{conditionattwo}  $\rhobar_{A,3}|_{G_{{\Q_2}}}$ is unramified, and the  characteristic polynomial 
of $\rhobar_{A,3}(\Frob_2)$ is not $(x^2\pm x+2)^2$. %
\item \label{conditionatthree} $A$ has good ordinary reduction at~$3$ and
 the characteristic polynomial %
of %
Frobenius at~$3$ does not have repeated roots.
\end{enumerate}
Then~$A$ is modular. More precisely, there exists a cuspidal automorphic representation~$\pi$
for~$\GL_4/\Q$ \emph{(}the transfer of a cuspidal
automorphic representation of~$\GSp_4/\Q$ of weight~$2$\emph{)} %
such that~$L(s,H^1(A)) = L(s,\pi)$, and hence ~$L(s,H^1(A))$ has a holomorphic continuation to~$\C$ satisfying the expected functional equation.
\end{ithm}

In our previous paper~\cite{BCGP}, we proved the \emph{potential} modularity of all abelian surfaces over  
totally real fields.
(We refer the reader to the introduction to~\cite{BCGP} for a history of the modularity conjecture for abelian surfaces.)
As a consequence, the main result of~\cite{BCGP} implies that~$L(s,H^1(A))$ has a \emph{meromorphic} continuation
to all of~$\C$, but it does not suffice to prove the conjectured holomorphicity, for essentially the same reason that Brauer~\cite{Brauer} was able to prove the meromorphic continuation of Artin $L$-functions but their holomorphicity remains conjectural.
The results of~\cite{BCGP} also allowed one to establish the \emph{modularity} of an abelian surface under extremely restrictive conditions, and in particular to produce~\cite[Thm~10.2.6]{BCGP} infinite (thin) sets of modular abelian surfaces~$A/\Q$ (up to twist) with 
$\End(A_{\Qbar}) = \Z$. These sets, however, 
account
for 0\% of all abelian surfaces over~$\Q$ counted in any reasonable way. Indeed, even producing \emph{any} explicit examples
where our modularity theorems applied was somewhat of a challenge~\cite{CGS}. 
In contrast, we expect that Theorem~\ref{first} applies to a positive proportion
of all abelian surfaces over~$\Q$ counted in any reasonable way.\footnote{It is very hard to say anything rigorous
(even for elliptic curves) if one orders by conductor.}
For example, conditions (\ref{surjective})--(\ref{conditionatthree}) can be guaranteed by imposing congruence conditions at finitely many primes (including~$2$ and~$3$).
See Section~\ref{sub:examples} for some more precise heuristics and examples; in particular we show that Theorem~\ref{first} applies to the Jacobians of precisely~$11743$ of the~$66158$ genus two curves in~\cite{LMFDB,database}.

The hypothesis~(\ref{surjective}) (which comes
from the Taylor--Wiles method) on the residual image can be weakened;   the allowable subgroups are precisely those
listed in Lemma~\ref{listofsubgroups} (they are all absolutely irreducible). 
Although there is some scope for 
marginal improvement on the local conditions~\eqref{conditionattwo} and~\eqref{conditionatthree} (as a direct
consequence of the modularity theorems proved in this paper), 
our expectation is that  the best way to relax the local assumptions is to make use of base change by generalizing
our main results to totally real fields, which  we hope to return  in the future. 
While some of our arguments will generalize straightforwardly, the proof of the main classicality theorem will require new ideas.
In ~\cite{BCGP}, we were able to work over totally real fields~$F$ in which~$p$
splits completely, and additionally there was considerable freedom to choose the prime~$p$. In contrast, in the current paper,  we are often forced to take~$p=3$ or~$p=2$, and in order to relax~\eqref{conditionattwo} and~\eqref{conditionatthree} it will be necessary to allow these primes to behave arbitrarily in the totally real field~$F$.

We also note the following easy to formulate corollary of Theorem~\ref{first} (again, see~\S\ref{sec:proofsofAandB}).

\begin{ithm} \label{second}
Let~$X: y^2 = f(x)$ with~$\deg(f(x)) = 5$
be a genus two curve over~$\Q$.
Suppose that:
\begin{enumerate}
\item ~$\rhobar_{X,3}: \Gal(\Qbar/\Q) \rightarrow \GSp_4(\F_3)$ is surjective.
\item $X$ has good ordinary reduction at~$2$.
\item $X$ has good ordinary reduction at~$3$.
\end{enumerate}
Then~$X$ is modular.
\end{ithm}

\subsection{The~$2$-$3$ switch} \label{sub:intro} 

The starting point of this paper is an analogue of 
the~$3$-$5$ switch used by Wiles~\cite{MR1333035}
 to prove residual modularity, %
which exploited  the rationality of certain twists of the modular curve~$X(5)/\Q$.  In our case, we use a rational moduli
space of abelian surfaces to carry out a~$2$-$3$ switch. %
This space was defined in~\cite{BCGP} as follows:
given an abelian surface~$A$ with a prime to~$3$ polarization, %
one may consider
the moduli space~$P(A[3])$ of genus two curves~$X$ equipped with a symplectic isomorphism~$\Jac(X)[3] \simeq A[3]$  and a fixed rational Weierstrass point.
By forgetting the Weierstrass point, the variety $P(A[3])$ admits a degree~$6$ map to a twist %
of the Siegel~$3$-fold of full level three, which is rational over~$\C$ but almost never over~$\Q$
~\cite{CC}.
However,~$P(A[3])$ is always rational by~\cite[Thm~10.2.1]{BCGP}.
  In particular, %
  we may find another
 abelian surface~$B/\Q$ with~$B[3] \simeq A[3]$ and such that  the map
$$\rhobar_{B,2}: G_{\Q} \rightarrow \GSp_4(\F_2) \simeq S_6$$
has image %
isomorphic to~$S_5$ (because of the rational Weierstrass point). 
Condition~\eqref{conditionattwo} of Theorem~\ref{first} ensures that we can find
 such a~$B$ with good ordinary reduction at~$2$.
After 
restricting to a quadratic extension~$F^+/\Q$ (which we can arrange to be totally real), we may assume
that~$\rhobar_{B,2} |_{G_{F^{+}}}$ is absolutely irreducible with image~$A_5$ in~$\GSp_4(\F_2)$. %
Known cases of the Artin conjecture in dimension two~\cite{MR3581178, MR3904451} allow
us to  identify this representation
with the mod~$2$ reduction of the symmetric cube of the $2$-adic Galois
representation associated to a Hilbert modular form of parallel weight~$2$, and thus (via known functorialities) 
to the mod~$2$ representation associated to a Hilbert--Siegel
eigenform %
(see also~\cite[Thm.\ 4.7]{tsuzuki2022automorphy}).
The goal is now to use modularity lifting theorems to go from the modularity of~$B[2]$
to the modularity of~$B$ and thus to the modularity of~$A[3]\simeq B[3]$,
and finally to the modularity of~$A$.

One difficulty that we encounter is that we need to prove modularity lifting theorems which apply when~$p=2$ and the residual image is rather small.
However by far the most serious difficulty compared to our previous work is that the argument above gives modularity of the residual representation~$\rhobar_{B,2}$ in regular weight, but the representation~$\rho_{B,2}$ has irregular weight.
The main innovation in our earlier work~\cite{BCGP} was to prove a modularity lifting theorem in irregular weight; however, this theorem crucially depended on having residual modularity in irregular weight as an input.

Deducing residual modularity in irregular weight from regular weight would be a higher dimensional analogue 
of showing (for modular forms) that a residual modular representation~$\rhobar: G_{\Q} \rightarrow \GL_2(\Fbar_p)$ which is unramified at~$p$ arises from a Katz modular
form of weight one~\cite{MR1074305,MR1185584}, and we do not know how to do this. %
Instead we use modularity lifting theorems to prove the existence of a $p$-adic Siegel modular form associated to the $p$-adic Tate module of our abelian surface, %
and we then  prove a classicality theorem for ordinary $p$-adic Siegel forms in irregular weight. Such $p$-adic modular forms are not necessarily classical; indeed their associated Galois representations need not be de Rham. However, we prove (under mild technical hypotheses, see Theorem~\ref{thm:multiplicity-one-implies-classical})  that if the Galois representation is (ordinary and) de Rham, then the form is indeed
classical. (Condition~\eqref{conditionatthree} in Theorem~\ref{first} guarantees that we can apply Theorem~\ref{thm:multiplicity-one-implies-classical} to~$\rho_{A,3}$.)

\subsection{Classicality for ordinary $p$-adic Siegel modular
  forms}\label{subsec:intro-classicality}Our classicality theorem
follows the strategy introduced by Lue Pan in his
paper~\cite{MR4390302}, as reinterpreted in~\cite{PilloniVB}. Both of
these papers (and Pan's sequel \cite{pan2022locally}) are in the
setting of the modular curve. While Pan gives a complete treatment of
arbitrary de Rham representations, we restrict to the ordinary
setting, which is considerably simpler; but there are still many
additional complications for higher-dimensional Shimura varieties. We prove our classicality theorem in the case of the Siegel threefold  (the Shimura variety
for $\GSp_4/\Q$) and it  asserts that if an ordinary $p$-adic modular eigenform of weight $2$ has an associated  Galois representation which is irreducible and de Rham (and satisfies a few more technical hypothesis), then it is classical (see Theorem \ref{thm:multiplicity-one-implies-classical}). 
The strategy for proving this theorem is to realize the Galois representation in
the completed cohomology of the Siegel threefold (it does not contribute to the classical \'etale cohomology because its Hodge--Tate
 weights are singular) and to relate the Sen operator of this Galois representation to a Cousin map which measures the obstruction for a $p$-adic modular form to be a classical modular form. In our ordinary case, the de Rhamness is equivalent to the semi-simplicity of the Sen operator which translates into the vanishing of the Cousin map and therefore implies the classicality of the $p$-adic  modular form.

We now give a more precise account of our strategy.  Let $\Sh_{K_pK^p}^{\tor}$ be a toroidal compactification of the  Siegel threefold  of  level $K^pK_p$ over $\Spa(\C_p, \ocal_{\C_p})$, and let  $\Sh_{K^p}^{\tor} = \lim_{K_p} \Sh_{K_pK^p}^{\tor} $ be   the perfectoid Siegel threefold over $\Spa(\C_p, \ocal_{\C_p})$, of prime-to-$p$ level $K^p$. Let $\omega^2_{K_p}$ be the sheaf of weight $2$ Siegel modular forms over $\Sh_{K_pK^p}^{\tor}$ and let $\omega^{2,\sm} = \colim_{K_p} \omega^{2}_{K_p}$, viewed as a sheaf over  $\Sh_{K^p}^{\tor}$, whose   cohomology is  $\colim_{K_p} \mathrm{R}\Gamma(\Sh_{K_pK^p}^{\tor}, \omega^{2}_{K_p})$. Thus, an element of the degree $0$ cohomology of $\omega^{2, \sm}$ is a weight $2$ Siegel modular form of   level $K^pK_p$ for some~$K_p$. 

The ordinary part $\mathrm{R}\Gamma(\Sh_{K^p}^{\tor}, \omega^{2,\sm})^{\ord}$ is computed by the following complex (more precisely, this is only true for cuspidal cohomology but we ignore this subtlety in the introduction) in degrees $0$ and $1$: 
\numequation
\label{cousinintro}
[\HH^0_{Id}(\Sh_{K^p}^{\tor}, \omega^{2,\sm})^{\ord} \stackrel{Cous}\rightarrow \HH^1_{^1w}(\Sh_{K^p}^{\tor}, \omega^{2,\sm})^{\ord}]
\end{equation}
where the module in degree $0$ in the complex is the space of ordinary $p$-adic modular forms of weight $2$, and the module in degree $1$ is a space of ordinary  higher $p$-adic modular forms (studied in higher Coleman theory). The differential is the Cousin map.

 Let $\mathrm{R}\Gamma(\Sh_{K^p}^{\tor}, \qq_p)$ denote the complex of completed cohomology.  We prove (under technical assumptions, Theorem \ref{thm-ECpadic}) that the ordinary part of the  $\mathfrak{b}$-cohomology  
\numequation
 \label{acomplex}
 \mathrm{R}\Hom_{\mathfrak{b}}(\lambda, \mathrm{R}\Gamma(\Sh_{K^p}^{\tor}, \qq_p)^{\la})^{\ord}
 \end{equation}
  of locally analytic vectors in completed cohomology is concentrated in degree
  $3$. Here $\lambda $ is  a non-dominant %
  character of the torus of $\mathrm{GSp}_4$ where we expect to see (by interpolation from what happens for dominant characters) the Galois representation of weight $2$ modular forms. %

After we tensor \eqref{acomplex} with $\C_p$, by  the $p$-adic Eichler--Shimura theory developed in this paper,  the cohomology  admits a $4$ step filtration and the graded pieces are  given  by the various relevant higher Coleman theories. If we denote by $V$ this degree $3$ cohomology group, we prove that $V_{\C_p}$ has a decreasing  filtration with $\mathrm{Gr}^0 V_{\C_p} = \HH^0_{Id}(\Sh_{K^p}^{\tor}, \omega^{2,\sm})^{\ord}$ and $\mathrm{Gr}^1 V_{\C_p} = \HH^1_{^1w}(\Sh_{K^p}^{\tor}, \omega^{2,\sm})^{\ord}$. 
We use this to verify that the Galois representation of our $p$-adic modular eigenform of weight $2$ is realized in  completed cohomology. 

The Sen operator respects the filtration and acts by the scalars $1, 1,0,0$ %
  on the respective graded pieces,  but possibly acts non semi-simply on $V_{\C_p}$. Looking at the generalized $0$ eigenspace for the Sen operator,  we get an induced map  
\numequation
\label{senintro}
\Sen: \mathrm{Gr}^0 V_{\C_p} = \HH^0_{Id}(\Sh_{K^p}^{\tor}, \omega^{2,\sm})^{\ord} \rightarrow
 \HH^1_{^1w}(\Sh_{K^p}^{\tor}, \omega^{2,\sm})^{\ord}=  \mathrm{Gr}^1 V_{\C_p} 
 \end{equation}
  measuring the failure of semi-simplicity of the Sen operator. 

The main result from which we deduce our classicality theorem  is the property that the Cousin map~\eqref{cousinintro} and the Sen map~\eqref{senintro} agree up to a non-zero scalar (Theorem \ref{thm-Sen-Cousin}).
 The key idea behind the proof is that the Sen operator which acts on the cohomology arises from an operator  defined on the complex of sheaves  $\mathrm{RHom}_{\mathfrak{b}}(\lambda, \oscr_{\Sh^{\tor}_{K^p}}^{\la})$ on the perfectoid Shimura varieties  whose cohomology is the $\mathfrak{b}$-cohomology of locally analytic vectors in completed cohomology. The complex  $\mathrm{RHom}_{\mathfrak{b}}(\lambda, \oscr_{\Sh^{\tor}_{K^p}}^{\la})$ is closely related to  twisted $D$-modules on the flag variety (a form of the Beilinson--Bernstein localization  of the Verma module of weight $\lambda$), and the Sen operator to a certain horizontal Cartan action. It turns out that the relation  between the Sen and Cousin maps can already be studied  and established at this more ``explicit'' geometric representation level (Theorem \ref{thm-extequal}). 

\subsection{$p$-adic Eichler--Shimura theory} \label{sec:babyintro} 
  A substantial part of this paper is dedicated to  $p$-adic Eichler--Shimura theory. In this part of the work, we are able to work in the generality of  Hodge type Shimura varieties. %
Let us recall first the classical Hodge--Tate decomposition of  modular
curves. Let $G = \mathrm{GL}_2$, and
let $K \subseteq G(\mathbb{A}_f)$ be a compact open subgroup. Let $\Sh_{K}$ be the modular curve over $\C_p$ of level $K$, with compactification $\Sh_K^{\tor}$. Let $E \rightarrow \Sh_K$ be the universal elliptic curve. In \cite{MR0909221},  Faltings proved the following Hodge--Tate decomposition for the \'etale cohomology of modular curves:  
\numequation
\label{faltingsintro}
\HH^1(\Sh_{K}, \mathrm{Sym}^k T_p E) \otimes_{\ZZ_p} \C_p = \HH^0(\Sh_K^\tor, \omega^{k+2}) (-k-1) \oplus \HH^1(\Sh_K^\tor, \omega^{-k})
\end{equation}
This isomorphism is equivariant for the Hecke action and the  local Galois action, and $(-k-1)$ indicates a Tate twist.
This kind of  statement  has been generalized to all Shimura varieties (see for example \cite{MR4552642}). 

Both sides of  \eqref{faltingsintro} 
are classical instances of bigger $p$-adic objects. On the left hand side we may consider completed cohomology, and on the right hand side we can consider higher Coleman theory \cite{boxer2021higher}. The main
goal of $p$-adic Eichler--Shimura theory (as taken up
in~\S\ref{sec:application-to-Shimura}) is to express some relation between both big $p$-adic spaces, generalizing Faltings's theory to non-classical cohomologies. 
We note that  $p$-adic Eichler--Shimura theory was initiated in
\cite{MR3315057} and completely transformed after the work of Pan~
\cite{MR4390302} followed by that of Rodríguez Camargo \cite{camargo2022locally,camargo2023geometric}. 
In order to state our main results, we need to introduce a certain amount
of notation as well as recall a number of facts from higher Coleman theory.
For this reason,
we defer any further discussion
to the more technical introduction given in~\S\ref{sec:eichlerintro}.

\subsection{An outline of the paper}
Here is a brief synopsis of the sections in our paper; see also the
introductions to the individual sections for more details. %

\S\ref{sec:Lie-algebra-homology} is concerned with Lie algebras. We consider the enveloping algebra $U(\mathfrak{g})$ of a finite dimensional  Lie algebra $\mathfrak{g}$ and  its Fr\'echet completion $\hat{U}(\mathfrak{g})$, as well as modules over them. Our main result (of independent interest) is Theorem \ref{thm-STrict} which compares the Lie algebra  cohomology of  a unipotent radical of a parabolic of $\mathfrak{g}$ of certain algebraic $U(\mathfrak{g})$-modules and of their completions. 

\S\ref{sec: classicality} is about equivariant twisted $D$-modules on flag varieties. We develop a somewhat ad hoc language to describe them. One of the main difficulties is to keep track of the various topologies and finiteness conditions we want to impose. We use the language of condensed mathematics to deal with functional analysis. We introduce a version of Beilinson--Bernstein localization and describe it using the results of Section \ref{sec:Lie-algebra-homology}.

\S\ref{sec:application-to-Shimura} contains our main classicality result. We also give some complements on higher Coleman theory and establish the $p$-adic Eichler--Shimura theory.

\S\ref{sec:
  R=T unitary} proves an~$R = \T$ theorem  in regular weight
  when~$p=2$ under a suitable
  oddness hypothesis, following~\cite{MR3598803}.
   For technical reasons (due to the small residual image of our representations), we need to work with unitary groups rather than symplectic groups.

  \S\ref{sec: R=T GSp4 preliminaries all p} proves an~$R=\T$ theorem
  in regular weight for~$p > 2$ 
 for symplectic representations.
Curiously enough,  when~$p=3$ and
the image of~$\rhobar$ is~$\GSp_4(\F_3)$ (the main case of interest),   technical reasons
now mandate that we work
with symplectic groups rather than unitary groups; see Remark~\ref{rem:curious}.

\S\ref{sec: SW type argument} proves a multiplicity one
result for certain Hida families, which  (once again) for technical reasons is necessary for our
classicality argument. This is where the main modularity
theorems Proposition~\ref{prop: the modularity result for p equals 2 or 3} and
Theorem~\ref{thm:rho-is-modular-from-mult-one-and-classicity} are proved,
using the  classicality result~Theorem~\ref{thm:multiplicity-one-implies-classical}.

\S\ref{sec:  modularity abelian surfaces} beings by recalling the basic theory of~$2$-torsion points on an abelian surface, and then establishes some basic but necessary facts
concerning the modular representation theory of~$A_5$.
This section also addresses the residual modularity
of mod-$2$ representations with image~$A_5$ using known cases of the Artin
conjecture for~$n = 2$.

\S\ref{sec:twothreeswitchsection} studies the representations~$\rhobar: G_{\Q_p} \rightarrow \GSp_4(\F_3)$
such that~$\rhobar^{\vee} \simeq \Jac(X)[3]$ for a genus two curve~$X/\Q_p$ with
 good ordinary reduction and a rational Weierstrass point when~$p=2$ or~$3$. We also study the related question of when~$\rhobar^{\vee} \simeq A[3]$
 where~$A/\Q_p$ is an abelian surface with good ordinary reduction and  a rational odd theta characteristic, as well as variants
 in which ordinary semistable reduction is allowed --- note that even when~$A = \Jac(X)$, it is possible that~$A$ has good reduction even
 when~$X$ does not.  This analysis is then used in~\S\ref{subsec: 2 3 switch} to carry out the~$2$-$3$ switch and
  then in~\S\ref{sec:proofsofAandB} to complete the  proofs of our main modularity theorems (Theorem~\ref{first} and~\ref{second}).

\S\ref{sec:complements} gives some examples and complements to our main theorem,
proving a residual modularity theorem for mod~$2$ representations with image~$A_6$
or~$S_6$, and proving the automorphy of any abelian surface~$A/\Q$
which  neither 
has~$\End(A_{\Qbar}) = \Z$ nor satisfies~$\End(A_K) \otimes \R = \End(A_{\Qbar})  \otimes \R =  \R \oplus \R$
for some quadratic field~$K$ (this excluded case includes the restriction of scalars of a general elliptic curve over~$K$).
We also explain why the full modularity theorem for all abelian surfaces over $\Q$ would follow
from a version of Serre's conjecture for~$\GSp_4(\F_p)$ in \emph{regular} weight.

\subsection{The work of Arthur} It should be noted that this paper,
as with the paper~\cite{BCGP} (see~\cite[1.4.1]{BCGP}), relies on results stated
by Arthur in~\cite{MR2058604} which ultimately rely on 
references [A24], [A25], [A26], and [A27] which have not (still) yet appeared,
as well as cases of the twisted weighed fundamental lemma announced
in~\cite{MR2735371}. However,  the situation has improved
remarkably in recent times. As a result of the recent preprint~\cite{A25A26A27}
of Atobe, Gan, Ichino, Kaletha, M\'{\i}nguez, and Shin,
a complete proof of  all the missing ingredients from Arthur's
papers is now  available, and thus the only result
we use for which a proof is not yet
available is 
the twisted weighted fundamental lemma.

\subsection{Acknowledgments}We would like to thank Jack Thorne for
several helpful conversations about his paper~\cite{MR3598803}, and
about $2$-adic automorphy lifting theorems in general.
We would also like to thank
Shiva Chidambaram, Lue Pan, Dave Roberts, Juan Esteban Rodríguez Camargo,
Travis Scrimshaw, and Andrew Sutherland for helpful conversations.
 Several of the ideas of this paper were discovered when all four authors
  were visiting the DFG-funded Hausdorff Research Institute for Mathematics as part of the
   trimester program ``The Arithmetic of the Langlands Program'' in~$2023$.

\subsection{Notation and conventions}\label{notn}

\subsubsection{Assorted notation}
We write~$\Z^n_+\subset\Z^n$ for the subset of
tuples~$(\lambda_1,\dots,\lambda_n)$ with $\lambda_1\ge\lambda_2
\ge\dots\ge\lambda_n$. If~$L/\Qp$ is a finite extension, we write $L_{\cycl}:=L(\zeta_{p^{\infty}})$
for the cyclotomic extension.%

\subsubsection{Coefficients}
We let~ $E$ be a finite extension of $\Qp$ with ring of integers
$\cO$, uniformizer~$\varpi$ and residue field $k$. We will always
assume that $E$ is chosen to be large enough such that all irreducible
components of all deformation rings that we consider, and all
irreducible components of their special fibres, are geometrically
irreducible. (We are always free to enlarge~$E$ in all of the
arguments that we make, so this is not a serious assumption.) Given a
complete Noetherian local $\cO$-algebra $\Lambda$ with residue field
$k$, we let $\CNL_{\Lambda}$ denote the category of complete
Noetherian local $\Lambda$-algebras with residue field $k$.  We refer
to an object in $\CNL_{\Lambda}$ as a
$\CNL_{\Lambda}$-algebra. If~$G$ is a group functor on
$\CNL_{\Lambda}$ then we write $\Ghat$ for the group functor
on~$\CNL_{\Lambda}$ given by $\Ghat(R):=\ker(\Ghat(R)\to\Ghat(k))$. %

\subsubsection{Galois representations and $p$-adic Hodge theory}We assume
without further comment that all Galois representations are continuous with
respect to the natural topologies. We normalize Hodge--Tate weights so that
the cyclotomic character has Hodge--Tate weight $-1$, and the Sen operator acts
via $1$ on the Sen module of $\qq_p(1)$, so that the (generalized) Hodge--Tate
weights are the negatives of the eigenvalues of the Sen operator. We
write~$\varepsilon$ for the $p$-adic cyclotomic character.

Let~$K/\Q_l$ be a finite extension for some~$l$ (possibly equal to~$p$). As
in~\cite[\S 2.8]{BCGP} we say that a representation $G_K\to\GL_n(\Qpbar)$ is
pure if the corresponding Weil--Deligne representation is pure; in the
case~$l=p$, this presupposes that the representation is de Rham. We say that a
representation $G_K\to\GSp_4(\Qpbar)$ is pure if the corresponding
representation  $G_K\to\GL_4(\Qpbar)$ is pure. If~$F$ is a number field then we
say that a representation $G_F\to\GL_n(\Qpbar)$ (or $G_F\to\GSp_4 (\Qpbar)$) if
it is pure at all finite places of~$F$.

\subsubsection{Notation for reductive groups}  \label{subsubsec-reductive-groups-notation}%
We consider a split reductive group $G$ over $\Spec~E$, with Borel $B$ and torus $T$. We denote by $\mathfrak{g}$, $\mathfrak{b}$,
$\mathfrak{h}$ their Lie algebras.  We write $\bar{B}$ for the opposite Borel of $B$ and write $\overline{\mf{b}}$ for
its Lie algebra. We let~$\Phi$ be the set of roots
of~$G$, with positive roots~$\Phi^+$ and negative
roots~$\Phi^-=-\Phi^+$. For~$\alpha\in\Phi$ choose standard basis elements~
$X_{\alpha}$ so that $(X_{\alpha},X_{-\alpha},H_{\alpha})$ is an
$\mf{sl}_2$-triple, where~$H_{\alpha}:=[X_{\alpha},X_{-\alpha}]$. We write~$\Delta$ for the set of simple roots.  We let $W$ be the Weyl group of $G$, with length function
$\ell:W\to\Z_{\ge 0}$, and write~$w_0$ for the longest element of~$W$. The Weyl group acts on the left on the
character group~$X^{*}(T)$ via $(w\lambda)(t):=\lambda(w^{-1}tw)$. It
also acts on the left on~$X_{*}(T)$, and the natural pairing
$\langle,\rangle$ between $X^{*}(T)$ and~$X_{*}(T)$ is $W$-equivariant.

Let $ P\supseteq B$ be a standard
parabolic with Levi quotient~$M$. We let $\mathfrak{p}$ be its Lie algebra, with
 unipotent radical $\mathfrak{u}_{\mathfrak{p}}$ and Levi
 $\mathfrak{m}$, and we write~$\mf{z}_{\mf{m}}$ for the centre
 of~$\mf{m}$. We have a Borel $\mathfrak{b}_\mathfrak{m}=\mathfrak{m}\cap\mathfrak{b}\subseteq\mathfrak{m}$.  We write $\mathfrak{u}$ for the unipotent radical of $\mathfrak{b}$ and $\mathfrak{u}_{\mathfrak{p}}$ for the unipotent radical of $\mathfrak{p}$.
 Let~$\Phi^+_M$ be the subset of~$\Phi^+$ which lie in
the Lie algebra of~$M$, and set $\Phi^{+,M}:=\Phi^+\setminus
\Phi^+_M$; and write  $\Phi^-_M:=-\Phi^+_M$, $\Phi^{-,M}:=-\Phi^{+,M}$. We let  $W_M$  be the Weyl
group of $M$, with longest element~$w_{0,M}$, and we let $\WM \subseteq W$ be the set of Kostant %
representatives of $W_M \backslash W$ (i.e.\
those~$w\in W$ with $\Phi^+_M\subseteq w\Phi^+$; this is a set of
coset representatives of minimal length). There is an involution of
$\WM$ given by $w\mapsto w_{0,M}ww_0$, and we have
$l(w_{0,M}ww_0)+l(w)=|\Phi^{+,M}|$. In particular the Kostant
representative of maximal length is $w_0 ^M:=w_{0,M}w_0$.

We let $\rho$ be half the sum of the
positive roots, and write $\rho = \rho^M+\rho_M$ where $\rho^M$ is half the
sum of the roots in~$\Phi^{+,M}$ and
$\rho_M$ is half the sum of the  roots in
$\Phi^+_M$. 

We define the ``dot action'' of~$W$ on~$X^{*}(T)$ by
$w\cdot\lambda:=w(\lambda+\rho)-\rho$. We say that $\lambda\in
X^{*}(T)$ is \emph{regular} if the stabilizer of~$\lambda$ for the dot
action is trivial, and otherwise we say that~$\lambda$ is
\emph{singular} or \emph{irregular}.

Let $w \in W$. Then we let~$P_w:=w^{-1}Pw$, with Lie algebra $\mathfrak{p}_w:=w^{-1}
\mathfrak{p} w$, and similarly we define
$\mathfrak{u}_{\mathfrak{p}_w}$, $\mathfrak{m}_w$, and so on.

\subsubsection{Notations for a $p$-adic torus}\label{notation-torus}
Let $T \rightarrow \Spec~\qq_p$ be a torus. We let $T^d$ be its maximal split subtorus. 
The group $T(\qq_p)$ has a unique maximal compact subgroup that we denote (abusing notation) by $T(\ZZ_p)$. 
We have an exact sequence \[0 \rightarrow \bar{\ZZ}_p^\times \rightarrow \Qpbar^\times \stackrel{v}\rightarrow \qq \rightarrow 0\] given by the $p$-adic valuation, normalized by $v(p)=1$. Tensoring this sequence by $X_\star(T)$ and taking invariants under the Galois group $\mathrm{Gal}(\Qpbar/\qq_p)$ yields an exact sequence \numequation\label{eqn:slope-left-exact-tori}0 \rightarrow T(\ZZ_p) \rightarrow T(\qq_p) \stackrel{v}\rightarrow X_\star(T^d)\otimes \qq\end{equation} where the image of $v$ is a lattice. 

Let $\chi : T(\qq_p) \rightarrow  \Qpbar^\times$ be a character. We can compose it with the $p$-adic valuation and get a map $v (\chi) : T(\qq_p) \rightarrow {\qq}$, which factors as $T(\qq_p) \stackrel{v}\rightarrow X_\star(T^d)\otimes \qq \rightarrow \qq$. We can therefore think of $v(\chi)$ as an element of $X^\star(T^d)_{\qq}$. 

If $T$ is a maximal torus contained in a Borel of a quasi-split reductive group $G$ defined over $\qq_p$ and if $\Phi = \Phi^+ \cup \Phi^{-}$ is the set of absolute roots, we let 
$T^+(\qq_p) = \{ t \in T(\qq_p), v( \alpha(t)) \geq 0~\forall \alpha \in \Phi^+\}$ and $T^{++}(\qq_p) = \{ t \in T(\qq_p), v( \alpha(t)) > 0~\forall \alpha \in \Phi^+\}$.

 \subsubsection{Notations in the symplectic case}\label{notn:GSp2g}%

We will often consider the case where $G = \mathrm{GSp}_{2g}$ and $P$ is the Siegel parabolic. In this case we make some more explicit choices. The group $G$ has a natural model over $\Spec~\ZZ$,  namely we realize $G$ as the subgroup of $\mathrm{GL}_{2g}$ acting on the free $\ZZ$-modules of rank $2g$, with basis $e_1, \cdots, e_{2g}$ and  preserving up to a similitude factor the symplectic form with matrix    $$J=\begin{pmatrix} %
      0 & S \\
      -S & 0 \\
   \end{pmatrix}$$ where $S$ is the $g \times g$ anti-diagonal matrix with only $1$'s on the anti-diagonal. We denote by $\nu: \mathrm{GSp}_{2g} \rightarrow \mathbb{G}_m$ the similitude factor. 
   
   We let $P$ be the stabilizer   of $\langle e_{g+1}, \cdots, e_{2g} \rangle$.
   We choose $B \subseteq P$ to be upper triangular on each diagonal block. We
   let $T$ be the diagonal torus. An element of $T$ is labelled %
 $t=\mathrm{diag} (z t_1, \cdots, zt_g, zt_{g}^{-1}, \cdots, zt_1^{-1})$. Characters $X^\star(T)$ of $T$ are tuples: 
   $\kappa = ( k_1, \cdots, k_g; w) \in \ZZ^g \times \ZZ$ with $w = \sum
   k_i\pmod 2$, and $\kappa (t) =z^w \prod_{i=1}^g t_i^{k_i}$. 
   A character is $M$-dominant if $k_1 \geq \cdots \geq k_g$. The set of $M$-dominant characters is denoted by $X^\star(T)^{M,+}$.  A character is $G$-dominant if $0 \geq k_1 \geq \cdots \geq k_g$. The set of $G$-dominant characters is denoted by $X^\star(T)^{G,+}$.

\subsubsection{$\GSp_4$}\label{notn:GSp4}
We now specialize further to the case~$G=\GSp_4$. We continue to take $P$ to be
the (``block lower-triangular'') Siegel parabolic stabilizing $e_3
,e_4 $, and  $B$ the %
Borel inside it which is upper-triangular in each of the diagonal
$2\times 2$ blocks.
We let~$Q$ be the Klingen parabolic containing $B$ (this is the other maximal parabolic in $\GSp_4 $) with Levi $M_Q$. %

Let $\kappa = (k_1, k_2; w) \in X^\star(T)^+$ be a  dominant weight
for $\mathrm{GSp}_4$, so that $0 \geq k_1 \geq
k_2$. Given our choice of Borel, the positive roots~$\Phi^+$ are $\alpha = (1,-1;0),
\beta = (-2,0;0), \gamma = \alpha + \beta  = (-1,-1;0), \delta =
2\alpha +\beta = (0,-2;0)$. We have~$\rho=(-1,-2;0)$.

The Weyl group~$W$ is generated by $s_\alpha$ and $s_\beta$ where
$s_\alpha(k_1,k_2;w) = (k_2,k_1;w)$ and $s_\beta(k_1,k_2;w) = (-k_1,k_2;
w)$, so that $w_0=s_\alpha s_\beta s_\alpha s_\beta $, and $w_0 (k_1 ,k_2 ;w)=(-k_1 ,-k_2 ;w)$. We have~$W_M=\{\Id,s_\alpha\}$ and~$W_{M_{Q}}=\{\Id,s_\beta\}$. The elements of $\WM$ are $\Id, s_\beta, s_\beta s_\alpha, s_\beta s_\alpha s_\beta$. We label
them ${^0w},{^1w},{^2w}, {^3w}$; they respectively have length
$0,1,2,3$. In particular, $^3w=w_0^M$  is the length three element. 
We  use the
pairing between characters and cocharacters coming from the standard pairing on
~$\qq^3$. Thus, we label cocharacters $X_\star(T)$ as triples $(a,b;c) \in \frac{1}{2}\mathbb{Z}^3$, with $a+c, b+c  \in \ZZ$. To $(a,b;c)$ we attach the cocharacter   $t\mapsto \mathrm{diag}(t^{a+c}, t^{b+c}, t^{-b+c}, t^{-a+c})$. 
We let $\mu = (-1/2,  -1/2; 1/2) \in X_\star(T)$. %
We sometimes view $\mu$
as an element of $\mathfrak{z}_{\mathfrak{m}}$ (the centre of the Lie algebra~$\mf{m}$). 

We let $\widehat{\mathrm{GSp}_4}$ be the dual group of $\mathrm{GSp}_4$. Our choice of Borel $B$ and torus $T$ in $\mathrm{GSp}_4$ gives a Borel $\widehat{B}$ and torus $\widehat{T}$ in $\widehat{\mathrm{GSp}_4}$. We use the spin representation to identify $\widehat{\mathrm{GSp}_4}$, the Borel $\widehat{B}$ and torus $\widehat{T}$ with the group $\mathrm{GSp}_4$, its usual upper triangular Borel and diagonal torus. 
In particular, this fixes  an isomorphism $X^\star(T) =  X_\star(\widehat{T}) \simeq X_\star(T)$, given by $$(\lambda_1, \lambda_2;w) \mapsto [t \mapsto \mathrm{diag}(t^{\frac{-\lambda_1-\lambda_2+w}{2}}, t^{\frac{\lambda_1-\lambda_2+w}{2}}
, t^{\frac{-\lambda_1+\lambda_2+w}{2}},
t^{\frac{\lambda_1+\lambda_2+w}{2}})]$$
Dually, there is  an isomorphism $X_\star(T) = X^\star(\widehat{T}) \simeq X^\star(T)$ for which $\mu$ corresponds to the dominant character $(1,0;1)$ of $X^\star(\widehat{T})$. 
When we work on the dual side (typically when we consider Galois representations), we will also denote by $B$ the upper triangular  Borel in $\mathrm{GSp}_4 \simeq \widehat{\mathrm{GSp}_4}$. This should not cause any confusion.

\subsubsection{Ordinary Galois
  representations}\label{subsec-ordinary-Galois-definitions}  %
\begin{defn}\label{defn:ordinary-for-Galois-reps}Let~$K/\Qp$ be a finite extension, and let
$\rho:G_K\to\GSp_4 (\Qpbar)$ be a representation with similitude
factor~$\varepsilon^{-1}$.
  We say that ~$\rho$ is \emph{ordinary} if there are characters
  $\chi_1 ,\chi_2:G_K\to\Qpbartimes$ with\[\rho\cong \begin{pmatrix}
    \chi_{1}&*&*&*\\
    0 &\chi_{2}&*&*\\
    0&0&\varepsilon^{-1}\chi_{2}^{-1}&*\\
    0&0&0&\varepsilon^{-1}\chi_{1}^{-1}
  \end{pmatrix}.\]We say that the ordered pair
    $(\chi_1 ,\chi_2 )$ is a \emph{$p$-stabilization} of~$\rho$. We say
    that~$\rho$ is \emph{$p$-distinguished} if the 4
characters
$\chi_1,\chi_2,\varepsilon^{-1}\chi_2^{-1},\varepsilon^{-1}\chi_1^{-1}$
are pairwise distinct. We say that~$\rho$ is \emph{semistable of weight 2} if
the subrepresentation \[
  \begin{pmatrix}
    \chi_{1}&*\\
    0 &\chi_{2}
  \end{pmatrix}
\]is unramified. (Such a representation is automatically semistable in the usual
sense.) In this case we will sometimes denote the $p$-stabilization $(\chi_1
,\chi_2 )$ by $(\alpha,\beta)$ with $\alpha=\chi_1 (\Frob_K)$, $\beta=\chi_2 (\Frob_K)$.

Similarly, we say that a representation
$\rhobar:G_K\to\GSp_4 (\Fpbar)$ with similitude factor~$\varepsilonbar$  is \emph{ordinary} if there are characters
  $\chibar_1 ,\chibar_2:G_K\to\Fpbartimes$ with\[\rhobar\cong \begin{pmatrix}
    \chibar_{1}&*&*&*\\
    0 &\chibar_{2}&*&*\\
    0&0&\varepsilon^{-1}\chibar_{2}^{-1}&*\\
    0&0&0&\varepsilon^{-1}\chibar_{1}^{-1}
  \end{pmatrix}.\]We say that the ordered pair
    $(\chibar_1 ,\chibar_2 )$ is a \emph{$p$-stabilization} of~$\rhobar$.
    We say
    that~$\rho$ is \emph{residually $p$-distinguished} if the 4
characters
$\chibar_1,\chibar_2,\varepsilon^{-1}\chibar_2^{-1},\varepsilon^{-1}\chibar_1^{-1}$
are pairwise distinct.  We say that~$\rhobar$ is \emph{of weight 2} if the
subrepresentation \[
  \begin{pmatrix}
    \chibar_{1}&*\\
    0 &\chibar_{2}
  \end{pmatrix}
\]is unramified; in particular
the characters~$\chibar_1 ,\chibar_2 $ are  unramified. (Conversely,
if~$\chibar_1,\chibar_2 $ are distinct and unramified, then~$\rhobar$ is of
weight 2.) If $\rhobar$ is of weight~$2$, then we will usually denote the
$p$-stabilization $(\chibar_1 ,\chibar_2 )$ by~$(\alphabar,\betabar)$, where
$\alphabar=\chibar_1 (\Frob_K)$, $\betabar=\chibar_2 (\Frob_K)$.

We make the same definitions for integral representations
$\rho:G_K\to\GSp_4(\Zpbar)$, in which case a $p$-stabilization $(\chi_1 ,\chi_2
)$ induces a $p$-stabilization $(\chibar_1 ,\chibar_2 )$ of the mod~$p$
representation $\rhobar:G_K\to\GSp_4 (\Fpbar)$. Note that if~$\rho$ is
semistable of weight~$2$, then~$\rhobar$ is of weight~$2$. If we regard~$\rho$ as a lift
of~$\rhobar$, then we say that $(\chi_1 ,\chi_2 )$ is \emph{compatible} with
$(\chibar_1 ,\chibar_2 )$. Given two ordinary lifts~$\rho_1 ,\rho_2 $
of~$\rhobar$, we say that $p$-stabilizations of~$\rho_1 $ and~$\rho_2 $ respectively are
compatible if they induce the same $p$-stabilization of~$\rhobar$.
\end{defn}
\begin{rem}
  \label{rem: apologize for weight 2}We again (see~\cite[Rem.\
  7.3.2]{BCGP}) apologize for the terminology ``of weight 2''; these definitions
  are convenient later in the paper when we wish to appeal to results from~\cite{BCGP}. In particular we
  caution the reader that if~$\rho$ is of weight~$2$ and pure, then
  it is pure of weight~$1$ in the usual sense. Since we will never use
  the terminology ``pure of weight~$1$'' (or ``pure of weight~$2$'',
  for that matter), we hope that this will not lead to any confusion.
\end{rem}

\subsubsection{Galois representations associated to automorphic representations}\label{subsec-GalAAR}

We for the most part follow the conventions of our earlier paper~\cite{BCGP}, to which we refer
for further details. We begin with some brief recollections from~\cite[\S
2.3]{BCGP}.  If~$K/\Ql$ is a finite extension for some~$l$, then we let $\rec_K$
be the local Langlands correspondence of~\cite{ht}, which assigns to an
irreducible complex admissible representation $\pi$ of $\GL_n(K)$ a Frobenius
semi-simple Weil--Deligne complex representation $\rec_K(\pi)$ of the Weil group
$W_K$. We will write $\rec$ for $\rec_K$ when the choice of $K$ is
clear. In the case~$n=1$, $\rec_K$ is obtained from the Artin map
$\Art_K:K^\times \iso W_K^\ab$, which we  normalize to
send uniformizers to geometric Frobenius elements. Similarly, we denote the local Langlands correspondence
of~\cite{gantakeda} by $\recGT$; this assigns a $\GSp_4$-conjugacy classes of
$\GSp_4(\C)$-valued Weil--Deligne representation of~$W_K$ to each irreducible
smooth complex representation of~$\GSp_4(K)$. If $(r,N)$ is a Weil--Deligne
representation of $W_K$ we will write $(r,N)^{F-\semis}$ for its Frobenius
semi-simplification.

We fix once and for all for each prime~$p$ an isomorphism
$\imath=\imath_p:\C\cong\Qpbar$. We will sometimes omit these
isomorphisms from our notation, in order to avoid clutter. In
particular, we will frequently use that~$\imath$ determines a square
root of~$p$ in~$\Qpbar$ (corresponding to the positive square root
of~$p$ in~$\C$). We will often regard  automorphic representations as
  being defined over~$\Qpbar$, rather than~$\C$, by means of the fixed
  isomorphism $\imath:\C\cong\Qpbar$.  We 
write~$\rec_p$ and~$\recGTp$ for the local Langlands correspondences
for~$\Qpbar$-representations given by conjugating by~$\imath$. 

Suppose that~$F^+$ is a
  totally real field and that~$\pi$ is a cuspidal automorphic
  representation of~$\GSp_4/F^+$. We will always assume that such a~$\pi$ has
  central character~$|\cdot|^2$. (We apologize for this assumption, which seemed
  helpful at some points when writing~\cite{BCGP}, and suffices for applications to abelian surfaces.) We say that~$\pi$ is \emph{algebraic} if it is
  $C$-algebraic, and we say that it is \emph{regular algebraic}
  if~$\pi_{\infty}$ is an (essentially) discrete series representation.
  Suppose that~$\pi$ is algebraic. We say that it has weight $(\lambda_v)_{v|\infty}$
  where~$\lambda_v \in (X^{*}(T)_{\qq}^{+}-\rho )\cap X^\star(T)$, if $\pi_v$ has infinitesimal character $-\lambda_v- \rho$. If $\pi$ is regular algebraic then $\lambda_v \in X^\star(T)^{+}$, and  we know  that $\pi \otimes \bigotimes_{v \mid \infty} V_{\lambda_v}$ has non-trivial $(\mathfrak{g} ,K_\infty)$-cohomology where $V_{\lambda_v}$ is the highest weight $\lambda_v$-representation.

  We now come to the definition of ordinarity. Assume furthermore that~$p$
  splits completely in~$F^+$ (this is sufficient to us). Our fixed isomorphism
  $\C\cong\Qpbar$ identifies $\{ w\mid p\}$ and $\{v \mid \infty\}$.  Suppose
  that $w \mid p$. We say that $\pi_w$ is \emph{finite slope} if it has
  non-trivial Jacquet module. The Jacquet module of $\pi_w$ is then a direct sum
  of characters $\chi_w : T(\qq_p) \rightarrow \bar{\qq}_p^\times$.  We say that
  $\pi_w$ is \emph{ordinary} if there is a character $\chi_w$ occurring in the
  Jacquet module such that $v(\chi_w) =-\lambda_w$ (see
  Section~\ref{notation-torus} for the definition of $v(\chi_w)$). We refer to a
  choice of such a character as an (ordinary) \emph{$p$-stabilization}
  of~$\pi_w$.  %
  We say that~$\pi$
  is ordinary if ~$\pi_w$ is ordinary for all~$w\mid p$. If~$\pi$ is ordinary
  and regular algebraic, then each ~$\pi_w$ has a unique  (ordinary)
  $p$-stabilization.

  \begin{thm}\label{thm:regular-weight-Galois-rep-recalled}Suppose that~$F^+$ is
    totally real and that~$p$ splits completely in~$F^+$.
    If~$\pi$ is regular algebraic of weight
    $\lambda = ((k_v, l_v;2))_{v\mid \infty}$, then for each prime~$p$ there is
    (see e.g.\ \cite[Thms.\ 2.7.1, 2.7.2]{BCGP})
    a semi-simple representation $\rho_{\pi,p}:G_{F^+}\to\GSp_4(\Qpbar)$
    satisfying the following properties. %
    \begin{itemize}
    \item $\nu\circ\rho_{\pi,p}=\varepsilon^{-1}$. \item For each finite
      place~$v\nmid p$, we
      have
      \[\WD(\rho_{\pi,p}|_{G_{F^+_v}})^{\semis}\cong\recGTp(\pi_v\otimes|\nu|^{-3/2})^{\semis}.\]

    \item If $\rho_{\pi,p}$ is irreducible, then for each finite place~$v$
      of~$F$, $\rho_{\pi,p}|_{G_{F^+_v}}$ is pure
      and
      \[\WD(\rho_{\pi,p}|_{G_{F^+_v}})^{F-\semis}\cong\recGTp(\pi_v\otimes|\nu|^{-3/2}).\]%

    \item For each~ $v|p$, $\rho_{\pi,p}|_{G_{F^+_v}}$ is de Rham with
      Hodge--Tate weights 
      $((k_v+l_v)/2-1,-(k_v-l_v)/2, (k_v-l_v)/2+1, 2-(k_v+l_v)/2)$. %

    \item If~ $p$ splits completely in $F$, $v \mid p$, and $\pi_v$ is ordinary,
      then there are potentially unramified characters $\alpha, \beta$ such
      that:
      \numequation\label{eqn:rho-pi-p-Borel-explicit}\rho_{\pi,p}|_{G_{F^+_v}}\cong \begin{pmatrix}
        \alpha\varepsilon^{1-(k_{v}+l_{v})/2}&*&*&*\\
        0 &\beta\varepsilon^{(k_{v}-l_{v})/2}&*&*\\
        0&0&\beta^{-1}\varepsilon^{-1-(k_{v}-l_{v})/2}&*\\
        0&0&0&\alpha^{-1}\varepsilon^{(k_{v}+l_{v})/2-2}
      \end{pmatrix}.\end{equation}
  \end{itemize}
\end{thm}

\begin{rem}\label{rem-recipe-diag} We can spell out more precisely the
  characters on the diagonal in~\eqref{eqn:rho-pi-p-Borel-explicit}. Let $\chi_v
  : T(\qq_p) \rightarrow \bar{\qq}_p^\times$ be an ordinary
  $p$-stabilization of~$\pi_v$. This induces a
  $p$-stabilization of~$\rho_{\pi,p}|_{G_{F^+_v}}$ in the sense of
  Definition~\ref{defn:ordinary-for-Galois-reps} as follows. Let \[\widetilde{\chi}_v = \chi_v \lambda_v \rho
  \nu^{-\frac{3}{2}} \vert \rho \nu^{-\frac{3}{2}} \vert^{-1} : T(\qq_p)
  \rightarrow \bar{\Z}_p^\times.\] This character is valued in $\bar{\Z}_p^\times$ by the ordinarity assumption. We can identify    $\widetilde{\chi}_v$ with an  homomorphism  $\qq_p^\times \rightarrow \widehat{{T}}(\bar{\Z}_p)$, where $\widehat{T}$ is the dual torus, which we identify with $T$ by using the isomorphism $\mathrm{GSp}_4 \simeq \widehat{\mathrm{GSp}}_4$ of Section \ref{notn:GSp4}. Then by class field theory, we interpret $\widetilde{\chi}_v : G_{\qq_p} \rightarrow {{T}}(\bar{\Z}_p)$. This is  the  character on the diagonal of $\rho_{\pi,p}|_{G_{F^+_v}}$. 
\end{rem}

We will also need to use the Galois representations associated to certain
irregular weight algebraic cuspidal automorphic representations for~$\GSp_4 /F^+$. 

 \begin{defn} We say that $\pi$   has weight~$2$ if it is algebraic of weight $\lambda = (1,1;2)_{v\mid \infty}$ (remember that by our convention, this $\lambda$ is not $G$-dominant) and $\pi_\infty$ is a non-degenerate limit of discrete series. 
 \end{defn}
The following theorem is well known. We provide a sketch of proof since we
couldn't find a precise reference in the literature. %

\begin{thm}\label{thm:irregular-ordinary-pi-Galois-rep}Suppose that~$F^+$ is
    totally real and that~$p$ splits completely in~$F^+$. Let  $\pi$ be an
    ordinary weight $2$ automorphic representation for~$\GSp_4 /F^+$.  There is 
a semi-simple representation
$\rho_{\pi,p}:G_{F^+}\to\GSp_4(\Qpbar)$ satisfying the following properties. %
\begin{itemize}
\item 
  $\nu\circ\rho_{\pi,p}=\varepsilon^{-1}$. \item For each finite place~$v\nmid p$, we have \numequation\label{eqn:semi-simplified-LGC}\WD(\rho_{\pi,p}|_{G_{F^+_v}})^{\semis}\cong\recGTp(\pi_v\otimes|\nu|^{-3/2})^{\semis}.\end{equation}

  \item There are potentially unramified characters $\alpha, \beta$ such that:
 \numequation\label{eqn:rho-pi-p-Borel-explicit}\rho_{\pi,p}|_{G_{F^+_v}}\cong  \begin{pmatrix}
    \alpha&*&*&*\\
0 &\beta&*&*\\
0&0&\beta^{-1}\varepsilon^{-1}&*\\
0&0&0&\alpha^{-1}\varepsilon^{-1}
  \end{pmatrix}.\end{equation}
  In fact, the character on the diagonal is described by the recipe explained in
  Remark~ \ref{rem-recipe-diag}.%
\end{itemize}
\end{thm}

\begin{proof} The representation $\pi_f$ will realize in the interior coherent
  cohomology of the Hilbert--Siegel Shimura variety by \cite[Thm.\ 2.7, Thm.\
  3.6.2] {harris-ann-arb}. By \cite[Thm.\ 1.4.3  (1) (4)]{boxer2021higher},
  $\pi_f$ defines a point $x$ on an equidimensional eigenvariety which dominates
  weight space. Let $\Spa(A,A^+)$ be an affinoid open subset of the eigenvariety
  containing $x$.   By \cite[Thm.\ 1.4.3 (2)]{boxer2021higher}, there is a
  Zariski dense set of classical points in $\Spa(A,A^+)$,  with regular
  algebraic weight. %
  Let~$X$ be the space of $\GSp_4$-valued pseudorepresentations 
  of~$G_{F^+}$ (in the sense of Lafforgue, see~\cite{quast2023deformationsgvaluedpseudocharacters}). Then by interpolation of the representations in Theorem~\ref{thm:regular-weight-Galois-rep-recalled} there is a map $\Spec A\to X$.  
  Specializing at $x$
  produces the semi-simple representation $\rho_{\pi,p}$. By interpolation the
  representation $\wedge^2\rho_{\pi,p}$ contains the
  character~$\varepsilon^{-1}$, so $\rho_{\pi,p}$ admits a symplectic pairing
  with multiplier~$\varepsilon^{-1}$. The statement
  regarding local-global compatibility away from $p$   follows by a standard
  argument from $p$-adic interpolation (note that the Weil--Deligne
  representations are only considered up to semi-simplification). If we assume that $\pi$ is ordinary, then we can assume that $\Spa(A,A^+)$ is an ordinary component of the eigenvariety, and by interpolation our local-global compatibility statement  at $p$ follows (see also \cite[Thm. 1.4.8]{boxer2021higher}, for a more general statement in the finite slope case). 
\end{proof}
\begin{rem}
  \label{rem:purity-implies-full-local-global}In the situation of Theorem~\ref{first}, we can
  upgrade the semi-simplified local-global compatibility~\eqref{eqn:semi-simplified-LGC} in
  Theorem~\ref{thm:irregular-ordinary-pi-Galois-rep} to full local-global
  compatibility. More precisely, if~$\pi$ is of general type and~$\rho_{\pi,p}$ in
   is pure,
  then \[\WD(\rho_{\pi,p}|_{G_{F^+_v}})^{F-\semis}\cong\recGTp(\pi_v\otimes|\nu|^{-3/2})\]for
  all~$v$; that is, in addition to~\eqref{eqn:semi-simplified-LGC}, the
  monodromy operators~$N$ on each side agree. To see this, note firstly that  since~$\pi$ is of general type, and cuspidal automorphic
  representations of~$\GL_n$ are generic, the $L$-packet containing
  $\pi_v\otimes|\nu|^{-3/2}$ is generic; so by part vii of the main theorem of~\cite{gantakeda}, the adjoint $L$-factor
$L(s,\ad(\recGTp(\pi_v\otimes|\nu|^{-3/2})))$ is holomorphic at
$s=1$. Equivalently, \numequation\label{eqn:adjoint-invariants-vanish-generic}(\ad(\recGTp(\pi_v\otimes|\nu|^{-3/2}))(1))^{\varphi=1,N=0}=0.\end{equation}
On the other hand, since $\rho_{\pi,p}|_{G_{F^+_v}}$ is pure, so is
$\WD(\rho_{\pi,p}|_{G_{F^+_v}})^{F-\semis}$ (by~\cite[Lem.\
1.4(1)]{ty}). By~\cite[Lem.\ 1.4(4)]{ty} and its proof, this means
that~$\WD(\rho_{\pi,p}|_{G_{F^+_v}})^{F-\semis}$ is equipped with the unique choice of~$N$ 
satisfying~\eqref{eqn:adjoint-invariants-vanish-generic}, as required.
\end{rem}

Finally, we will  need to use the Galois representations associated to certain
automorphic representations of~$\GL_n$, which we now very briefly recall. Let~$F$ be an imaginary CM field. Recall that an automorphic representation~$\pi$ of
$\GL_n/F$ is RACSDC if it is regular algebraic, conjugate
self-dual, (i.e.\ $\pi^c\cong\pi^\vee$), and cuspidal.
(See e.g.\ \cite[\S 2]{BLGGT}) for more details.)  Associated to a RACSDC automorphic representation~$\pi$ is a continuous semi-simple
representation $r_{\pi,p}:G_F\to\GL_n(\Qpbar)$ such that $r_{\pi,p}|_{G_{F_v}}$
is de Rham for all~$v|p$, and for each finite place~$v$ of~$F$ we have  \[ \imath \WD(r_{\pi,p}|_{G_{F_v}})^{F-\semis} \cong \rec(\pi_v \otimes
     |\det|_v^{(1-n)/2})\]  %
   (see e.g.\ \cite[Thm.\ 2.1.1]{BLGGT} and~\cite[Thm.\ 1.1]{MR3272276}).
In
particular, we have
$r_{\pi,p}^c\cong r_{\pi,p}^\vee\varepsilon^{1-n}$.

\subsubsection{Transfer between~$\GL_4 $ and~$\GSp_4 $}\label{sec:Arthur}%
We firstly very briefly recall some results on Arthur's
   classification of discrete automorphic representations of~$\GSp_4$;
   see~\cite[\S 2.9]{BCGP} for a slightly longer treatment with precise
   references to the literature.  Suppose that~$F$ is a number field,
 that~$\Pi$ is a cuspidal
   automorphic representation of~$\GL_4/F$, and that
   ~$\chi:\A_F^{\times}/F^{\times}\to \C^{\times}$ is unitary. Then we
   say that~$\Pi$ is $\chi$-self dual if 
   $\Pi\cong\Pi^\vee\otimes\chi\circ\det$, in which case the pair
   ~$(\Pi,\chi)$ is of \emph{symplectic type} if the partial
   $L$-function $L^S(s,\Pi,\bigwedge^2\otimes\chi^{-1})$ has a pole at
   $s=1$ (where $S$ is any finite set of places of $F$) or of
   \emph{orthogonal type} if $L^S(s,\Pi,\Sym^2\otimes\chi^{-1})$ has a
   pole at $s=1$. Exactly one of these alternatives holds, and
   if~$(\Pi,\chi)$ is of symplectic (resp.\ orthogonal) type then it
   descends to a discrete automorphic representation~$\pi$
   of~$\GSp_4/F$ (resp.\ ~$\GSpin_4^{\alpha}/F$ for some inner
   form~$\GSpin_4^{\alpha}$ of~$\GSpin_4$) with central
   character~$\omega_{\pi}=\chi$. (See for example~\cite[Prop.\
   6.1.7]{GeeTaibi}.)
   We say that a discrete automorphic representation~$\pi$
   of~$\GSp_4/F$ is \emph{of general type} if it arises in this
   way for some~$(\Pi,\chi)$, in which case we say that~$\Pi$ is the
   \emph{transfer} of~$\pi$, and that~$\pi$ is a \emph{descent}
   of~$\Pi$. For each place~$v$ of~$F$, the $L$-parameter obtained
   from $\recGT(\pi_v)$ by composing with the usual embedding
   $\GSp_4\into\GL_4$ is $\rec(\Pi_v)$.  In this case~$\pi$ is
   necessarily cuspidal, and it is stable. In fact
   if~$\pi':=\otimes'\pi'_v$ with $\pi_v,\pi'_v$ in the same
   $L$-packet for all~$v$, then~$\pi'$ is automorphic, and occurs with
   multiplicity one in the discrete spectrum. If~$\pi$ is
   (regular) algebraic then ~$\Pi$ is also (regular)
   algebraic. %

      If~$F$ is totally real, and ~$\pi$ is regular algebraic and \emph{not} of general type, then the
   Galois representations~$\rho_{\pi,p}$ associated to~$\pi$ are reducible
   by~\cite[Lem.\ 2.9.1]{BCGP}. Since we will always be in a
   situation where our Galois representations are irreducible (even
   irreducible modulo~$p$), we will only need to consider~$\pi$ of
   general type in this paper.

  \subsubsection{Galois representations associated to abelian surfaces} 
\label{sec:galoisintro}
Let~$F$ be a number field, and let~\(A/F\) be an abelian surface. 
For each prime~\(p\), we may write~\(\rho_{A,p}\)
for the Galois representation associated to~\(H^1(A_{\overline{F}},\Z_p)\).
We often think of \(\rho_{A,p}\) as a representation
\[\rho_{A,p}: G_{F} \to \GSp_4(\Q_p)\]
with multiplier given by the inverse cyclotomic character~$\varepsilon^{-1}$
(compare~\cite[Defn.\
2.8.2]{BCGP}). We also let \(\rhobar_{A,p}\) denote the Galois
 representation associated to \(H^1(A_{\overline{F}},\F_p)\). If \(A\) admits a principal polarization
 of degree prime to \(p\), then we can and do think of \(\rhobar_{A,p}\) as a representation
 \[\rhobar_{A,p}: G_{F} \to \GSp_4(\Fbar_p).\]
 We take the coefficient field of \(\rho_{A,p}\) (respectively, of \(\rhobar_{A,p}\))
  to be
 \(\Q_p\) or \(\Qbar_p\)
 (resp. \(\F_p\) or \(\Fbar_p\))
 depending on what is most convenient.
  If \(T_p(A)\)
denotes the \(p\)-adic Tate module of \(A\), then (in our conventions)
the Galois
representations associated to \(T_p(A)\)  and \(A[p]\) are   the dual representations
\(\rho^{\vee}_{A,p} \simeq \rho_{A,p} \otimes \varepsilon\)
and
\(\rhobar^{\vee}_{A,p} \simeq \rhobar_{A,p} \otimes \varepsilon\)
respectively.
   The representation~$\rho_{A,p}$ is unramified at all but finitely many
  places \(v\) of \(F\), and if~ 
  $v|p$ then $\rho_{A,p}|_{G_{F_v}}$ is de Rham with Hodge--Tate weights
  0,0,1,1 for every choice of embedding \(F \rightarrow \Qbar_p\). Furthermore $\rho_{A,p}|_{G_{F_v}}$  is
  pure at all finite places~\(v\)
  (see e.g.\ \cite[Prop.\ 2.8.1]{BCGP}).
  If \(A/F_v\) has good ordinary reduction for some \(v|p\),
  then  $\rho_{A,p}|_{G_{F_v}}$ is crystalline and ordinary
  of weight 2. %

\subsubsection{Notions of modularity}%
\label{sec:definitionofmodularity}
Let~$F$ be a number field.

\begin{df} \label{definitionautomorphic} An abelian surface~$A/F$ 
is \emph{modular}, or equivalently, \emph{automorphic}  if
there exist %
$C$-algebraic
cuspidal automorphic representations~$\pi_i$ for~$\GL_{n_i}/F$
with~$4 = \sum n_i$ 
such that
\[L(s,H^1(A)) = \prod L(s,\pi_i \otimes | \det |^{(1-n_i)/2}).\]
A genus two curve~$X/F$ is modular if~$A = \Jac(X)/F$ is modular.
\end{df}

If~$A$ is modular, 
then~$L(s,H^i(A)) = L(s,\wedge^i H^1(A)) 
 = L(s,\Pi_i)$ for some automorphic representation~$\Pi_i$.
 This follows from known functorialities in small degrees, most 
 notably~\cite{MR1937203,MR2567395} (cf.\ the proof of~\cite[Thm~9.3.1]{BCGP}).

\begin{rem}[Warning] In~\cite{BCGP}, particularly~\cite[Defn.\ 9.1.8]{BCGP}, we
  reserved the term \emph{modular} to specifically refer to the stronger
  statement that~\(F\) was totally real and that \(A\) was associated to a
  cuspidal automorphic representation of \(\GSp_4/F\) with certain
  properties. With such a restriction, there are abelian surfaces and genus two
  curves which fail to be modular, for example, when~$A/\Q = \Jac(X)/\Q$ is
  isogenous to a product of two elliptic curves or an abelian surface
  of~$\GL_2$-type. In retrospect, we feel that this distinction is unhelpful. In
  the main theorems of this paper, we (under certain hypotheses) establish the
  modularity of~$A/\Q$ by proving the modularity of~$\rho_{A,p}$ for
  some~$p$. More precisely, we assume that~$\rho_{A,p}$ is absolutely
  irreducible, and show that~$\rho_{A,p}\cong\rho_{\pi,p}$ for some weight~$2$
  cuspidal automorphic representation~$\pi$ for~$\GSp_4 /\Q$. This~$\pi$ will be
  of general type, and thus transfers to a $C$-algebraic cuspidal automorphic
  representation of~$\GL_4 $.
\end{rem}

\subsubsection{The Eichler--Shimura relation}\label{sec-ESR}
We let $\ell$ be a prime. We let 
\[\mathbf{T}_\ell=\Z[\GSp_4(\Q_{\ell})/\kern-0.2em{/}{\GSp_4(\Z_{\ell})]}\]
 be the spherical Hecke
algebra for~$\GSp_4(\Q_{\ell})$ with $\ZZ$-coefficients. As a $\ZZ$-module, it
has a basis consisting of  the characteristic functions of the double cosets $T_\lambda =
[\GSp_4(\Z_{\ell}) \lambda(\ell) \GSp_4(\Z_{\ell})]$ where $\lambda \in
X_\star(T)^+$. 
In particular, we define $T_{\ell,i}=[\GSp_4(\Z_l)\beta_{\ell,i}\GSp_4(\Z_l)]$, where %
\[\beta_{\ell,0}=\diag(\ell,\ell,\ell,\ell),\]
\[\beta_{\ell,1}=\diag(\ell,\ell,1,1),\] \[\beta_{\ell,2}=\diag(\ell^2,\ell,\ell,1).\]  We
write $Q_\ell(X)\in \TT[X]$ for the
polynomial \[X^4-T_{\ell,1}X^3+(\ell T_{\ell,2}+(\ell^3+l)T_{\ell,0})X^2-\ell^3T_{\ell,0}T_{\ell,1}X+\ell^6T_{\ell,0}^2.\]

We have the Satake isomorphism $$S : \qq(\sqrt{\ell}) [X^\star(\widehat{T})]^{W}
\isoto  \qq(\sqrt{\ell}) \otimes_{\ZZ} \mathbf{T}_\ell.$$ %
For each representation $V$ of $\widehat{\GSp_4}$
we let
$[V]$ be the character of $\widehat{T}$ on~$V$. This defines an element of $ \qq(\sqrt{\ell}) [X^\star(\widehat{T})]^{W}$.
To each $\lambda \in X_\star(T)$ we can associate a representation $V_\lambda$ of $\widehat{\mathrm{GSp}_4}$ with highest weight $\lambda$. The $[V_\lambda]$ form a basis of $\qq(\sqrt{\ell}) [X^\star(\widehat{T})]^{W}$. 
We consider in particular  $\lambda = (-\frac{1}{2},-\frac{1}{2};\frac{1}{2})
\in X_\star({T})$, %
corresponding to the Spin representation
$\widehat{\mathrm{GSp}}_4 \rightarrow \mathrm{GL}_4$ (which, as explained above, we use to identify
$\widehat{\mathrm{GSp}}_4$ and $\mathrm{GSp}_4$, so that via the isomorphism
$X_\star(T) = X^\star(\hat{T}) = X^\star(T)$,
$(-\frac{1}{2},-\frac{1}{2};\frac{1}{2})$ goes to $(1,0;1)$). We also consider
the dual of the Spin representation, corresponding to $\lambda =
(-\frac{1}{2},-\frac{1}{2};-\frac{1}{2})$. We write $P_{\ell,\lambda}(X) = X^4 -
 [V_\lambda]X^3 + [\Lambda^2V_\lambda]X^2 - [\Lambda^3 V_\lambda]X + [\Lambda^4
 V_\lambda]$  for the  characteristic polynomial of the
 representation $V_\lambda$ in either of these cases. We write $Q_\ell(X)=\ell^6P_{\ell,
   (-\frac{1}{2},-\frac{1}{2};\frac{1}{2})}(\ell^{-\frac{3}{2}}X)$. Then~$Q_{\ell}$
 is the usual Hecke polynomial in $\mathbf{T}_\ell[X]$ whose definition was recalled above.  The coefficient of $X^3$ is
 $-T_{(-\frac{1}{2},-\frac{1}{2};\frac{1}{2})}$.  We also let
 \begin{equation}
   \label{eq:P-ell}   P_{\ell}(X):=\ell^6P_{\ell,(-\frac{1}{2},-\frac{1}{2};-\frac{1}{2})}(\ell^{-\frac{3}{2}}X).
 \end{equation}
 The coefficient of $X^3$ in~$P_{\ell}$ is $- T_{(-\frac{1}{2},-\frac{1}{2};-\frac{1}{2})}$.

 Let $\pi$ be a $C$-algebraic automorphic representation of $\mathrm{GSp}_4/\qq$
 whose component at infinity is a non-degenerate limit of discrete series. Let
 $S$ be the set of finite places at which $\pi$ is not spherical. Let
 $\mathbb{T}^{S}:=\otimes_{\ell\notin S}\mathbb{T}_{\ell}$ be the spherical Hecke algebra away from $S$. %
 We let $\Theta_\pi : \mathbb{T}^{S} \rightarrow \C$ be the character describing
 the action of~$\mathbb{T}^{S}$ on the
 one dimensional $\C$-vector space of spherical vectors of $\otimes_{v\notin
   S}\pi_v$. %
 Then by the definition of~$\rec$, the Galois representation
 $\rho_{\pi,p} : G_{\qq} \rightarrow \mathrm{GSp}_4(\Qpbar)$ has the property
 for all primes $\ell \notin S \cup \{p\}$,
 $\imath(\Theta_{\pi}(Q_\ell))(X)$ is the characteristic
 polynomial of  $\rho_{\pi,p} (\mathrm{Frob}_\ell)$ (here~$\Frob_{\ell}$
 denotes a geometric
 Frobenius element). %
 
 \begin{lem}  For all primes $\ell \notin S \cup \{p\}$,
   $\imath(\Theta_{\pi}(P_\ell))(X)$ is the characteristic
   polynomial of  $(\rho^\vee_{\pi,p}\otimes\varepsilon^{-3})
   (\mathrm{Frob}_\ell)$. %
 \end{lem}
 \begin{proof} Unraveling the definitions, we find that $\ell^{-6}\imath(\Theta_{\pi}(P_{\ell,(-\frac{1}{2},-\frac{1}{2};-\frac{1}{2})}))(\ell^{\frac{3}{2}}X)$ is the characteristic polynomial of $\rho^\vee_{\pi,p} (\mathrm{Frob}_\ell)$. The characteristic polynomial of $\ell^{3}\rho^\vee_{\pi,p} (\mathrm{Frob}_\ell)$ is therefore $\ell^6\imath(\Theta_{\pi}(P_{\ell,(-\frac{1}{2},-\frac{1}{2};-\frac{1}{2})}))(\ell^{-\frac{3}{2}}X)$. 
 \end{proof}

 For any neat compact open subgroup $K = \prod_\ell K_\ell \subseteq
 \mathrm{GSp}_4(\mathbb{A}_f)$, let $\Sh^{alg}_{K} \rightarrow \Spec~\qq$ denote
 the Siegel threefold of level $K$. %
 We will make use of the following Eichler--Shimura relation.
 
 \begin{thm}\label{thm-ESRFC} On $\mathrm{R}\Gamma(\Sh^{alg}_{K, \bar{\qq}}, \ZZ/p^n\ZZ)$ and $\mathrm{R}\Gamma_c(\Sh^{alg}_{K, \bar{\qq}}, \ZZ/p^n\ZZ)$, for each place $\ell\neq p$ at which $K_\ell$ is hyperspecial, the local Galois representation of $G_{\qq_\ell}$ is unramified at $\ell$ and $P_\ell(\mathrm{Frob}_\ell)=0$. 
 \end{thm}
 \begin{proof}  Let $\ell \neq p$ be a place at which $K_\ell$ is hyperspecial. We have a natural smooth integral model $\Sh^{alg}_{K, \ZZ_\ell} \rightarrow \Spec~\ZZ_\ell$.  We first claim that $\mathrm{R}\Gamma(\Sh^{alg}_{K, \bar{\qq}_\ell}, \ZZ/p^n\ZZ) = \mathrm{R}\Gamma(\Sh^{alg}_{K, \bar{\FF}_\ell}, \ZZ/p^n\ZZ)$.   By  \cite[Coro. 5.20]{Lan_Stroh_2018}\footnote{whose proof considerably simplifies in our case, due to the existence of smooth toroidal compactifications, with normal crossing boundary divisor.}, 
 $$\mathrm{R}\Gamma(\Sh^{alg}_{K, \bar{\qq}_\ell}, \ZZ/p^n\ZZ) = \mathrm{R}\Gamma(\Sh^{alg}_{K, \bar{\FF}_\ell}, \mathrm{R}\Psi\ZZ/p^n\ZZ).$$ Since  $\Sh^{alg}_{K, \ZZ_\ell} \rightarrow \Spec~\ZZ_\ell$ is smooth,  the map $\ZZ/p^n\ZZ \rightarrow  \mathrm{R}\Psi\ZZ/p^n\ZZ$ is an isomorphism.  By Poincar\'e duality, we deduce  that 
 \[\mathrm{R}\Gamma_c(\Sh^{alg}_{K, \bar{\qq}_\ell}, \Z/p^n\Z) = \mathrm{R}\Gamma_c(\Sh^{alg}_{K, \bar{\FF}_\ell}, \ZZ/p^n\ZZ).\]
 We now use the Eichler--Shimura relation of \cite[VII, Thm 4.2]{MR1083353}, to
 deduce that $P_\ell(\mathrm{Frob}_\ell)=0$. It only remains to explain why it is the polynomial $P_\ell$ and not $Q_\ell$ that we need to use. This all boils down to understanding how we attach to a characteristic function of a  double coset in the Hecke algebra, a Hecke correspondence. 
 Using our conventions (which we think are  standard, but are the transpose of that of \cite{MR1083353}), to the double coset $T_{\lambda}$ is associated  the Hecke correspondence:  
 \begin{eqnarray*}
 \xymatrix{  \Sh^{alg}_{\lambda(\ell)^{-1} K \lambda(\ell) \cap K} \ar[d]^{p_2} &   \Sh^{alg}_{\lambda(\ell) K \lambda(\ell)^{-1} \cap K} \ar[l]^{\lambda(\ell)} \ar[d]^{p_1} \\
  \Sh^{alg}_{K} &  \Sh^{alg}_{ K}  }
 \end{eqnarray*}
 For example, for $\lambda = (-\frac{1}{2},-\frac{1}{2};-\frac{1}{2})$, this is the moduli space parametrizing abelian surfaces 
 $p_1^\star A$ and $p_2^\star A$, with certain  prime-to-$\ell$ level structure  and prime-to-$\ell$ polarization, together with an isogeny (compatible with level structure and polarization) $p_1^\star A \rightarrow p_2^\star A$ whose kernel  is a  maximal isotropic subgroup of $p_1^\star A[\ell]$. The reduction of the natural integral model of this correspondence modulo $\ell$  contains the Frobenius correspondence. 
 \end{proof}

\section{Lie algebra homology}%
\label{sec:Lie-algebra-homology}
\subsection{Introduction}%
Let~$\mf{g}$ be a reductive Lie algebra over~$E$, let $\mf{p}$ be a parabolic
subalgebra of~$\mf{g}$ with Levi~$\mf{m}$ and unipotent
radical~$\mf{u}_{\mf{p}}$, and let~$\mf{b}$ be a Borel of~$\mf{g}$ containing a
Cartan~$\mf{h}$, which we assume is also contained in~$\mf{m}$. In this section
we study the~$\mf{u}_{\mf{p}}$-cohomology of objects of category~$\cO$ and of
category~$\hat{\cO}$, a $p$-adic analytic version of the BGG category~$\cO$. The categories~$\cO$ and~$\hat{\cO}$ are equivalent,
via base change from the universal enveloping algebra~$U(\mf{g})$ to
its completion, the Fr\'{e}chet--Stein algebra~$\hat{U}(\mf{g})$. We establish in particular the key Theorem~\ref{thm-STrict}, which shows that in
a fixed $p$-adically non-Liouville weight, the operation of taking
$\mf{u}_{\mf{p}}$-cohomology is compatible with completion, i.e.\ with passage
from category~$\cO$ to category~$\hat{\cO}$.

We use the language of condensed mathematics throughout, and we begin
in Section~\ref{subsec:solid-K-v-s} with an overview of the results
that we need (mostly from~\cite{MR4475468}) on solid $E$-vector
spaces, together with a summary of some results from~\cite{MR3072116}
on category~$\hat{\cO}$.

\subsection{Solid functional analysis and representations}%

\subsubsection{Solid $E$-vector spaces}\label{subsec:solid-K-v-s}%
Rather than use the classical theory of
topological vector spaces, we work throughout with the condensed
mathematics of Clausen--Scholze~\cite{Clausen-Scholze}; for the
convenience of the reader, here and below we recall some of the comparisons to
the classical definitions.
Let $E$ be a finite extension of $\qq_p$.
By \cite[lecture 7]{Clausen-Scholze}, the non-archimedean field $E$
can be viewed as a solid abelian group. It follows that $E$ can be equipped with a structure of an analytic ring,  where for any profinite set $S$, $E_{\blacksquare}[S] = E \otimes_{\ZZ} \ZZ_{\blacksquare}[S]$. 
We let $\operatorname{Mod}(E)$ %
be the abelian category of solid $E$-vector spaces; this has a tensor product,
which we denote by~$\otimes$, and an internal~$\Hom$, which we denote
by~$\underline{\Hom}(-,-)$.  We refer to \cite[\S3]{MR4475468},  for a  complete
treatment of non-archimedean functional analysis from the condensed
perspective. We simply recall what is strictly necessary for us.\footnote{In
  order to fix set-theoretical issues,  we choose a strongly inaccessible cardinal $\kappa$ and  we only consider $\kappa$-small profinite sets. See \cite[Lecture 1, rem. 1.3]{Clausen-Scholze}}
  
  We have a functor $V \mapsto \underline{V}$ from topological spaces to condensed sets,  where $\underline{V}$ is the condensed set defined by $\underline{V}(S) = \mathcal{C}^0(S, V)$ for any profinite set $S$. 
  This functor has a left adjoint $X \mapsto X(\star)_{\top}$ from condensed sets to topological spaces, given by evaluating a condensed set $X$ on the point $\star$ and endowing $X(\star)$ with the quotient topology of the map $\coprod_{S, x \in X(S)} S \rightarrow X(\star)$, where $S$ runs through all profinite sets. 
  The restriction of the  functor $V \mapsto \underline{V}$ to the
  category of  compactly generated topological spaces is fully
  faithful, and if~ $V$ is compactly generated then $V = \underline{V}(\star)_{\top}$ (more
  precisely, the counit $ \underline{V}(\star)_{\top}\to V$ of the adjunction restricts to the
  identity functor on compactly generated topological spaces, see
  \cite[Prop.\ 1.7]{Clausen-Scholze}).

  By  \cite[Proposition 3.7]{MR4475468}, the  functor $V \mapsto \underline{V}$ restricts to a functor from the category of complete locally convex $E$-vector spaces to the category of solid $E$-vector spaces.
 All the complete locally convex $E$-vector spaces that we will
 encounter will be considered as solid $E$-vector spaces unless
 explicitly specified otherwise.

  We introduce  certain  full
subcategories of $\Mod(E)$. 

\begin{defn} \leavevmode \begin{enumerate}
\item  A Banach space is a solid $E$-module of the form 
$ (\lim_n (\oplus_I \ocal_E/p^n\ocal_E))[1/p]$
 for some set $I$.
\item A Smith space is a solid $E$-module which has the form $(\prod_I \ocal_E)[1/p]$ for some set $I$.
\end{enumerate}
We let $B(E)$ be the category of Banach spaces  and $S(E)$ be the category of Smith spaces.
\end{defn}

\begin{rem}\label{rem-exact-Banach-classical-solid} The categories of
  solid and classical Banach spaces (resp.\ Smith spaces) are
  equivalent via the functors $V \mapsto V(\star)_{\top}$ and $V
  \mapsto \underline{V}$. The essential surjectivity follows from the
  explicit description of the objects. The full faithfulness is a consequence of the fact that classical Banach spaces and Smith spaces are compactly generated. (See for example~\cite[Prop.\ 3.5]{MR4475468}.)
\end{rem}

\begin{prop}\label{prop: antieq Smith Banach} \cite[Lem.\ 3.10]{MR4475468} There is an anti-equivalence of categories    between Smith and Banach spaces given by $V \mapsto V^\vee:=\underline{\mathrm{Hom}}(V, E)$. Moreover, $(V^\vee)^\vee = V$. 
\end{prop}

\begin{rem}\label{rem: exactness of antieq Smith Banach} The functor $V \mapsto V^\vee$ is exact in the sense that it sends short exact sequences of Banach spaces (resp.\ Smith spaces) to short exact sequence of Smith spaces (resp.\ Banach spaces). In fact, any short exact sequence is split. 
\end{rem}

\begin{rem} If $V$ is in $B(E)$,  then $V^\vee(\star)_{\top}$ is the classical Smith space equal to the continuous dual $\mathrm{Hom}(V(\star)_{\top}, E)$ equipped with the compact open topology. 
\end{rem}

\begin{defn}  \leavevmode
\begin{enumerate}
\item A Fr\'echet space is a solid $E$-module which can be written as a sequential limit of Banach spaces. 
\item An $\LS$-space is a solid $E$-module which can be written as a
  sequential colimit of Smith spaces with injective transition maps. 
\item An $\LB$-space is a solid $E$-module which can be written as a
  sequential colimit of Banach spaces with injective transition maps.
\end{enumerate}
We let $F(E)$ be the category of Fr\'echet spaces, we let $\LS(E)$ be
the category of $\LS$-spaces, and we let  $\LB(E)$ be the category of $\LB$-spaces. 
\end{defn}

\begin{rem}\label{rem:solid-classical-Frechet} The categories of solid and classical Fr\'echet spaces are
  equivalent under the functors $V \mapsto V(\star)_{\top}$ and $V
  \mapsto \underline{V}$, \cite[Lem.\ 3.24(1)]{MR4475468}. To see
  this, we claim that it  suffices to show that  $V \mapsto
\underline{V}$ is an essentially surjective functor from classical to
solid  Fr\'echet spaces. Indeed, since the counit
$\underline{V}(\star)_{\top}\to V$ is an isomorphism (because  Fr\'echet
  spaces are in particular compactly generated), we will then know
  that $V \mapsto
\underline{V}$ is fully faithful and essentially surjective, and thus
an equivalence; it follows formally from this that 
the unit $V\mapsto
  \underline{V(\star)_{\top}}$ of the adjunction is also an
  isomorphism of Fr\'echet spaces, as required.

 Now, if $V =\lim_r V_r$ is a classical Fr\'echet space (where the~$V_r$
  are classical Banach spaces), then by Remark~\ref{rem-exact-Banach-classical-solid},
$\underline{V} = \lim_r \underline{V_r}$ is a solid Fr\'echet space
(note that $V \mapsto \underline{V}$ commutes with limits, being
a right adjoint). Conversely, since by
definition a solid Fr\'echet space can be written as $V = \lim_r V_r$
where the~$V_r$ are Banach spaces, we have \[V=\lim_r V_r\isoto \lim_r
  \underline{V_r(\star)_{\top}},\]
which gives the essential surjectivity.

Note in particular that as a consequence of this equivalence, any (solid) Fr\'echet space admits a presentation where $V = \lim_r V_r$ with $V_{r+1}(\star)_{\top} \rightarrow V_r(\star)_{\top}$ has dense image. 
\end{rem}

\begin{prop}\cite[Thm.\ 3.40]{MR4475468}  We have an anti-equivalence of categories  $V \mapsto V^\vee: = \underline{ \mathrm{Hom}}(V, E)$ between $F(E)$ and $\LS(E)$ extending the biduality between $B(E)$ and $S(E)$. Moreover, $(V^\vee)^\vee = V$. The functor $V \mapsto V^\vee $ is exact. 
\end{prop}

 \begin{defn} \leavevmode \begin{enumerate}
 \item  A map $f: V \rightarrow W$ of Smith spaces is \emph{trace class} if there exists a map $g: E \rightarrow V^\vee \otimes W$ such that $f$ is the composite $V \stackrel{\Id_V \otimes g}\rightarrow V \otimes V^\vee \otimes W \stackrel{\ev \otimes \Id_W}\rightarrow W$. 
 \item A map $f: V \rightarrow W$ of Banach spaces is \emph{compact} if its dual is trace class. 
 \end{enumerate}
 \end{defn}
 
 \begin{example} Let $I$ be a set and let $ (a_i)_{i \in I} \in
   E^I:=\prod_{i\in I}E$ be a family converging to zero with respect to the net of the complements of finite subsets of $I$. 
 Let $ f: \ocal_E^I[1/p] \rightarrow \ocal_E^I[1/p]$ be the map sending $(x_i)_{i\in I}$ to $(a_i x_i)_{i\in I}$. Then one sees that $f$ is trace class, represented by the tensor $\sum_i a_i e^\vee_i \otimes e_i$ in $(\lim_n (\oplus_I \ocal_E/p^n\ocal_E))[1/p] \otimes \ocal_E^I[1/p]$ (where $e_i$ is $i$-th basis vector of  $ \ocal_E^I[1/p]$).  \end{example}

 \begin{defn} \leavevmode \begin{enumerate}
 \item An object $V$ of $\LS(E)$ is of compact type if it has a presentation $V = \colim_n V_n$ where the maps $V_n \rightarrow V_{n+1}$ are trace class.  \item An object $V $ of $F(E)$ is of compact type if it has a presentation $V= \lim_n V_n$ where the maps $V_n \rightarrow V_{n-1}$ are compact.
\item  An object $V$ of $\LB(E)$ is of compact type if it has a
  presentation $V= \colim_n V_n$ where the maps $V_n \rightarrow V_{n+1}$ are compact.  \end{enumerate}
 \end{defn}
 
 \begin{prop}\cite[Cor.\  3.38]{MR4475468} A solid $E$-module is an $\LB$-space of compact type if and only if it is a $\LS$-space of compact type. 
 \end{prop}

 We will use the following lemma in Remark~\ref{rem:flatness-of-U-p-hat}.
 \begin{lem}\label{lem-flat-frechet}  A Smith, Banach, $\LB$ or $\LS$-space is flat. A Fr\'echet space of compact type is flat over $E$.
\end{lem}
\begin{proof} The flatness of Smith, Banach, $\LB$ or $\LS$-spaces is \cite[Lem.\ 3.21]{MR4475468}.   Let~$V$ be a Fr\'echet space of compact type. By~\cite[Cor.\ 3.38(1)]{MR4475468}, we can write $V = \lim_n V_n$ as an inverse limit of
  Smith spaces. Following the proof of \cite[Lem.\ 3.21]{MR4475468}, it suffices
  to show that if $W' \rightarrow W$ is an injection of Smith spaces, then $W'
  \otimes V \rightarrow W \otimes V$ is injective.

  Since we have an injection $ \lim_n V_n \into \prod_n V_n $, and since Smith
  spaces are flat over~$E$ (by~\cite[Prop.\ 3.20, Lem.\ 3.21]{MR4475468}), it
  suffices to show that  $(\prod_n V_n) \otimes_E W' \rightarrow (\prod_n V_n)
  \otimes_E W$ is injective. For any Smith space~$X$, we have (see~\cite[Prop.\
  3.12]{MR4475468}) \[(\prod_n
  V_n) \otimes_E X = \prod_n ( V_n \otimes_E X),\]  so it suffices in turn to show that $\prod_n (V_n \otimes_E
  W') \rightarrow \prod_n (V_n  \otimes_E W)$  is injective. Since the Smith
  spaces~$V_n$ are flat, each morphism $V_n \otimes_E
  W' \rightarrow V_n  \otimes_E W$ is injective, and we are done.\end{proof}

 \subsubsection{Representations of  algebraic groups} In this section we recall the classical notion of representation of an algebraic group, before moving to representations of analytic groups.  We let   $\Mod^{\delta}(E)$ be the usual category of $E$-vector spaces (the superscript $\delta$ stands for discrete).
 Let $G = \Spec \oscr_G$ be an affine group scheme  over $\Spec E$. The algebra $\oscr_G$ is a Hopf algebra with comultiplication $\mu: \oscr_G \rightarrow \oscr_G \otimes \oscr_G$ and augmentation $e: \oscr_G \rightarrow E$.  
 We let $\mathrm{Mod}^{\delta}_G(E)$ be the category of algebraic representations of $G$. Its objects are vector spaces $V$ over $E$, equipped with a co-action map $ c: V \rightarrow V \otimes \oscr_G$ such that:
 \begin{enumerate}
 \item (associativity)  The maps $(c \otimes \Id) \circ c $ and $(\Id \otimes \mu) \circ c$:   $V \rightarrow V \otimes \oscr_G \rightarrow  V \otimes \oscr_G \otimes \oscr_G$ agree. 
 \item (neutral element) The map $(\Id\otimes  e) \circ c$: $V \rightarrow V \otimes \oscr_G
   \rightarrow V$ is the identity. 
 \end{enumerate}
 
 \subsubsection{Representations of analytic groups}
We recall the following standard definition.
 \begin{defn}\label{defn:quasi-Stein-space}
  An adic space $\mathcal{X}$ is called quasi-Stein if it has an open cover
  given by an increasing countable union of affinoid spaces of finite type
  $\mathcal{X} = \cup_n \mathcal{X}_n$ where
  $\HH^0(\mathcal{X}_{n+1}, \oscr_{\mathcal{X}_{n+1}}) \rightarrow
  \HH^0(\mathcal{X}_n, \oscr_{\mathcal{X}_n})$ has dense image.  A quasi-Stein
  space is Stein if it admits a covering as before having the property that
  $\mathcal{X}_n$ is relatively compact in $\mathcal{X}_{n+1}$
  (\cite[2.4]{MR1032938}); equivalently, if the closure
  $\overline{\mathcal{X}_n}$ of $\mathcal{X}_n$ in $\mathcal{X}_{n+1}$ is proper
  over $\Spa (E, \ocal_E)$.
 \end{defn}
We now let $G$ be a Stein  analytic group over $\Spa(E, \ocal_E)$.%
\begin{rem}We have two cases in mind: either $G$ is the
  analytification of an affine group scheme over $\Spec E$, or $G$ is
  a quasi-compact affinoid open subgroup of such an analytification. 
\end{rem}
We let $\oscr_G$ be the algebra of functions on $G$, which is an
object of $\Mod(E)$ (it is a Fr\'echet space). It has a structure of a
Hopf algebra.  We define the category $\Mod_G(E)$ %
of representations of $G$ on solid $E$-vector spaces. Its objects are solid vector spaces equipped with  a co-action map 
$c: V \rightarrow V \otimes \oscr_G$ satisfying the same conditions as
before. Similarly, we let $B_G(E)$ be the category of representations of $G$ on Banach modules. We let $\LB_G(E)$ be the category of representations of $G$ on $\LB$-spaces. 

We let $D(G) = \underline{\mathrm{Hom}}(\oscr_G, E) = \oscr_G^\vee$ be the distribution algebra of $G$. The dual  of the comultiplication $\mu: \oscr_G \rightarrow \oscr_G \otimes \oscr_G$  induces the algebra structure on $D(G)$. 
If $V$ is an object of $\Mod_G(E)$, then it is naturally a $D(G)$-module (via $V \otimes D(G) \stackrel{c}{\rightarrow} V \otimes \oscr_G  \otimes D(G) \rightarrow V$). 
We therefore have a natural functor $\Mod_G(E) \rightarrow \Mod(D(G))$
from the category of solid $G$-representations to the category of solid $D(G)$-modules.

\begin{rem} In some cases, one can go backwards. For example if $G$ is quasi-compact and $V$ is a Banach space, we have that $\underline{\mathrm{Hom}}( D(G), V) = \oscr_G \otimes V$ by \cite[Cor.\  3.17]{MR4475468} so that any $D(G)$-module structure on $V$ can be turned into an action of $G$ on $V$. 
If we denote by $B(D(G))$ the category of $D(G)$-modules which are Banach spaces, then the categories $B_G(E)$ and $B(D(G))$ are equivalent. 
\end{rem}
\subsubsection{Representations of locally profinite groups} %
We now let $M$ be a locally profinite group. We view $M$ as a condensed group.  We let $E[M]$ be the associated condensed ring and we let $E_{\blacksquare}[M]$ be its solidification.
If $M$ is compact, then $E_{\blacksquare}[M] = (\lim_N \ocal_E [M/N]
[1/p]$ where $N$ runs through the compact open subgroups of $M$. In
general, if $M_0 \subseteq M$ is a compact open subgroup, then we have
the formula $E_{\blacksquare}[M] = E[M] \otimes_{E[M_0]}
E_{\blacksquare}[M_0]$. %

\begin{defn}\label{def-locally-pro-rep} A  representation of $M$ over a  solid $E$-vector space is a solid $E_{\blacksquare}[M]$-module. The category of $M$-representations is denoted by $\Mod_M(E)$.  
\end{defn}
  
 \begin{rem} Equivalently, a representation of $M$ is a solid $E$-vector space $V$ and an action map $M \times V \rightarrow V$ of condensed sets  satisfying the usual group action axioms.
 \end{rem}
 \subsubsection{Smooth representations}
 Let $V \in \Mod_M(E)$. We let $V^{\sm}  = \colim_{N \subseteq M} V^N$ where $N$ runs through all compact open subgroups of $M$. We say that $V$ is smooth if the natural map $V^{\sm} \rightarrow V$ is an isomorphism. 
 We let $\Mod^{\sm}_M(E)$ be the category of smooth representations. 
 
 We let $M_{\disc}$ be the group $M$ equipped with the discrete topology. There is a natural map $M_{\disc} \rightarrow M$ of condensed sets. One can consider the category $\Mod_{M_{\disc}}(E)$ of $E[M_{\disc}]$-modules. 
 We can define the subcategory $\Mod_{M_{\disc}}^{\sm}(E)$ of smooth representations of $M_{\disc}$. Its objects are representations $V$ of $M_{\disc}$ such that $V = \colim V^{N_{\disc}}$ where $N$ goes through all compact open subgroups of $M$. 
 
 \begin{lem} The categories $\Mod^{\sm}_M(E)$ and $\Mod^{\sm}_{M_{\disc}}(E)$ are equivalent. 
 \end{lem}
 \begin{proof} We have a natural functor $\Mod_M(E) \rightarrow
   \Mod_{M_{\disc}}(E)$, induced by the map $M_{\disc} \rightarrow M$. This
   induces a functor $\Mod^{\sm}_M(E) \rightarrow
   \Mod^{\sm}_{M_{\disc}}(E)$.  We now construct a functor 
 $\Mod^{\sm}_{M_{\disc}}(E) \rightarrow \Mod^{\sm}_{M}(E)$. Let $V \in
 \Mod^{\sm}_{M_{\disc}}(E)$. Let $M_0$ be a compact open subgroup of
 $M$ and let $M_n$ be a system of normal compact open subgroups of $M_0$.  %
 We see that $V^{M_n}$ is an $E[(M_0)_{\disc}/(M_n)_{\disc}]$-module. Since $(M_0)_{\disc}/(M_n)_{\disc} = M_0/M_n$, we deduce that the $(M_0)_{\disc}$-module structure on $V^{M_n}$ extends uniquely to an $M_0$-module structure. Passing to the colimit over $n$, we deduce that $V$ is an $M_0$-module. Since $E[M] = E[M_{\disc}] \otimes_{E[(M_0)_{\disc}]} E[M_0]$, we are done. 
 \end{proof}

 Recall that an abelian category is an \emph{Grothendieck abelian category} if
 it has arbitrary colimits, it has a generator, and filtered colimits are exact
 (AB5). By \cite[\href{https://stacks.math.columbia.edu/tag/079H}{Tag
   079H}]{stacks-project}, any Grothendieck abelian category has enough
 injectives. We will now show that  $\Mod^{\sm}_M(E)$ is a Grothendieck abelian
 category; note that this relies on our set-theoretic assumption that we only consider
 $\kappa$-small profinite sets for some fixed~$\kappa$.
 \begin{lem} The category $\Mod^{\sm}_M(E)$ is a Grothendieck abelian category,
   and in particular it has enough injectives. %
 \end{lem}
\begin{proof} This is obvious, except for the existence of a generator. For a totally disconnected  $S$ and compact open subgroup $N \subseteq M$, we consider
$V_{S,N} = E_{\blacksquare}[M]\otimes_{E_{\blacksquare}[N]}
E_{\blacksquare}[S]$ with $N$ acting trivially on
$E_{\blacksquare}[S]$. 
We claim that $\oplus_{S, N} V_{S,N}$ is a generator. This follows from the property that for any $V \in \Mod^{\sm}_M(E)$, $\mathrm{Hom}_{ \Mod^{\sm}_M(E)}(V_{S,N}, V) =
V^N(S)$. %
\end{proof}

Here is a slight generalization of the concept of smooth.  Let
$\lambda: M \rightarrow E^\times$ be a character, and write
$E(-\lambda)$ for the corresponding representation of~$M$ (with
underlying vector space~$E$).
\begin{defn}
  We say that $V$ is $\lambda$-smooth if $V \otimes E(-\lambda)$
  is
  smooth. We let $\Mod_M^{\lambda-\sm}(E)$ be the category of
  $\lambda$-smooth $M$-modules.
\end{defn}
Note that if $\lambda, \lambda'$ are two characters such that $\lambda \otimes (\lambda')^{-1}$ is smooth, then the categories $\Mod_M^{\lambda-\sm}(E)$ and $\Mod_M^{\lambda'-\sm}(E)$ are canonically equivalent.

\subsubsection{Locally analytic representations}
We now assume that $M$ arises as the set of $\qq_p$-points of an analytic group $G$ over $\Spa(\qq_p, \ZZ_p)$ and we also assume that we have a fundamental system of quasi-compact open subgroups   $\{G_r\}_{r \geq 0}$ of $G$, where $G_r$ is a polydisc.  We let  $G_r(\qq_p) = M_r$. The $\{M_r\}$ form  a fundamental system of compact open subgroups in $M$. 
We now define the locally analytic vectors of $V \in \Mod_M(E)$. %
Note that $\oscr_{G_{r}} \otimes_E V$ has
three commuting left actions of $M_r$: %
$\star_l$, $\star_r$ and $\star_V$
(induced respectively by left translation  on the group, right translation on the group, and the original action on $V$).    The action $\star_V$  comes from an action of $M$.  Moreover, the group $M$ acts by conjugations $\star_{l,r}$ on its system of neighborhoods of identity $\{M_r\}$.  %
We set $V^{M_r-\an} = \HH^0(M_r, (\oscr_{G_{r}} \otimes V))$, where the
invariants are taken for  the action $\star_l \otimes
\star_V$. %
The space $V^{M_r-\an}$  still carries a $\star_{l,r} \otimes \star_V$-action of $M_r$. The evaluation map at $e$,  $ \oscr_{G_{r'}} \otimes V \rightarrow V$  induces an injective map  $V^{M_r-\mathrm{an}} \rightarrow V$. We let $V^{\la} = \colim_{r} V^{M_r-\mathrm{an}}$. This is an $M$-representation.  We thus have inclusions $V^{\sm} \rightarrow V^{\la} \rightarrow V$. 
\begin{rem} The functor $V \mapsto V^{\la}$ can naturally be derived into a functor $V \mapsto V^{\mathrm{Rla}}$. See \cite[sect. 4.4]{MR4475468}.
\end{rem}

\subsubsection{The algebra $\hat{U}(\mathfrak{g})$}\label{sect-dist-alg}
Let us assume now that we have an algebraic group $G^{\alg} \rightarrow \Spec \ocal_E$. Its analytification defines a quasi-compact affinoid analytic group $G=G_0 \rightarrow \Spa(E, \ocal_E)$. For any $r \in
\qq_{\geq 0}$, we let $G_r$ be the quasi-compact analytic subgroup
of $G_0$  of elements reducing to the identity $e$  modulo
$p^r$.  We have that $\oscr_{G,e} = \colim \oscr_{G_r}$ is an $\LB$-space of compact type. We let $\mathfrak{g}$ be the Lie-algebra of~ $G$. We define 
$\hat{U}(\mathfrak{g}) = \oscr_{G,e}^\vee$. This is a Fr\'echet space of compact type.

Since the categories of solid and classical Fr\'echet spaces are
equivalent, we will freely write $\hat{U}(\mathfrak{g})$ for the
underlying classical  $E$-algebra $\hat{U}(\mathfrak{g})(\star)_{\top}$ of the solid $E$-algebra $\hat{U}(\mathfrak{g})$. %
We have a natural  
map $\mathfrak{g} \rightarrow \hat{U}(\mathfrak{g})$ given by $X
\mapsto [f(g) \mapsto f'(g \operatorname{exp}(-tX))\vert_{t=0}]$, 
 which extends to a map from the enveloping algebra  $U(\mathfrak{g})
 \rightarrow \hat{U}(\mathfrak{g})$ with dense image.

 One can describe $\hat{U}(\mathfrak{g})$ as a
completion of $U(\mathfrak{g})$ %
as follows (following \cite[lemma 2.4]{MR1887640}   and  \cite[sect. 3.2]{MR3072116}). If we fix a basis $x_1, \cdots, x_n$ of $\mathfrak{g}$, then $U(\mathfrak{g}) = \oplus_{(n_i)} E \prod x_i^{n_i}$ by the PBW theorem. For each $r \in \R_{>0}$, we define a norm $\vert ~\vert_r$ on $U(\mathfrak{g})$ by putting $\vert  \sum a_{\underline{n}} x^{\underline{n}}\vert_r = \mathrm{sup}_{\underline{n}} \vert a_{\underline{n}} \vert r^{\sum n_i}$. We let $U(\mathfrak{g})_r$ be the completion of $U(\mathfrak{g})$ for $\vert-\vert_r$ and we have $\hat{U}(\mathfrak{g}) = \lim_{r \geq 1} U(\mathfrak{g})_r$.  %

\begin{rem} For any $r$, there exists $r'$ such that
$U(\mathfrak{g})_{r'} \rightarrow D(G_r)$ and  $D(G_{r'}) \rightarrow
U(\mathfrak{g})_{r}$. %
We thus have two presentations of $\hat{U}(\mathfrak{g}) = \lim_r
U(\mathfrak{g})_r = \lim_r D(G_r)$, as an inverse limit of Banach
spaces with compact transition maps and as an inverse limit of Smith spaces with trace class transition maps. 
\end{rem}
 
 Since $\hat{U}(\mathfrak{g})$ is a Fr\'echet--Stein algebra (see
 \cite[sect. 3]{MR1990669}), there is an associated  abelian category
 of coadmissible modules
 $\mathrm{Mod}^{\coad}(\hat{U}(\mathfrak{g}))$, which is defined as follows.
 \begin{defn} %
   A (left) $\hat{U}(\mathfrak{g})$-module $M$ is coadmissible if  it has a presentation $M = \lim M_r$ where $M_r$ is a finitely generated $\hat{U}(\mathfrak{g})_r$-module and  $M_{r+1} \otimes_{\hat{U}(\mathfrak{g})_{r+1}} \hat{U}(\mathfrak{g})_r = M_r$. 
 \end{defn} 
 We let $\Mod(\hat{U}(\mathfrak{g}))$ be the category of solid $\hat{U}(\mathfrak{g})$-modules. 
 \begin{thm}\label{thm-coad-to-sol} We have a fully faithful exact functor $\mathrm{Mod}^{\coad}(\hat{U}(\mathfrak{g})) \rightarrow \Mod(\hat{U}(\mathfrak{g}))$. 
 \end{thm}
 \begin{proof}By for example \cite[Prop.\ 3.1.1]{MR3072116}, any coadmissible module is canonically an  object of $F(E)$ of compact type. 
 \end{proof}
 \begin{defn}\label{defn:admissible-module}
   An \emph{admissible module} is the dual of a coadmissible module. Admissible
   modules are objects of $\LB(E)$ of compact type.
 \end{defn}

  \subsubsection{Categories of $\hat{U}(\mathfrak{g})$ and $U(\mathfrak{g})$-modules}
  
  We recall the the maps (of classical
  rings)  $U(\mathfrak{g}) \rightarrow U(\mathfrak{g})_r$ and $U(\mathfrak{g}) \rightarrow
  \hat{U}(\mathfrak{g})$ %
  are flat (see for
example \cite[Thm.\ 4.3.3]{MR3072116}).  Let
$\mathrm{Mod}^{\fing}(U(\mathfrak{g}))$ be the category of finitely generated left
$U(\mathfrak{g})$-modules. %

  \begin{prop}\label{prop-U(g)-to-hat} We have an exact functor: 
  
  \begin{eqnarray*}
\hat{U}(\mathfrak{g}) \otimes_{U(\mathfrak{g})}-: \mathrm{Mod}^{\fing}(U(\mathfrak{g})) & \rightarrow & \mathrm{Mod}^{\coad}(\hat{U}(\mathfrak{g})) \\
M & \mapsto & \hat{U}(\mathfrak{g}) \otimes_{U(\mathfrak{g})} M
\end{eqnarray*}
  \end{prop}
  \begin{proof} This follows from the flatness of $U(\mathfrak{g})  \rightarrow \hat{U}(\mathfrak{g})$.
  \end{proof}

 \begin{cor}\label{cor-exact-U(g)-to-hat}  We have  an exact  functor: 
  
  \begin{eqnarray*}\hat{U}(\mathfrak{g}) \otimes_{U(\mathfrak{g})}-:
 \mathrm{Mod}^{\fing}(U(\mathfrak{g})) & \rightarrow & \Mod( \hat{U}(\mathfrak{g}) ) \\
M & \mapsto & \widehat{M}:= \hat{U}(\mathfrak{g}) \otimes_{U(\mathfrak{g})} M
\end{eqnarray*}

  \end{cor}%

  \begin{proof} By combining Theorem \ref{thm-coad-to-sol} and Proposition \ref{prop-U(g)-to-hat}, we obtain an exact functor $\mathrm{Mod}^{\fing}(U(\mathfrak{g}))  \rightarrow \Mod(\hat{U}(\mathfrak{g}))$. 
  \end{proof}

 \subsubsection{Category $\ocal$ and category $\hat{\ocal}$}
  We assume that $ \mathfrak{g}$ is a reductive Lie algebra with Borel
  $\mathfrak{b}$ and Cartan $\mathfrak{h}$, and as usual we write~$\Phi^+$ (resp.\
  $\Phi^{-}$) for the positive (resp.\ negative) roots determined by our
fixed Borel subgroup~$\mf{b}$.    We can consider the abelian
  category $\ocal( \mathfrak{g}, \mathfrak{b})$ (simply denoted
  $\ocal$ if $\mathfrak{g}$ and $\mathfrak{b}$ are clear from the
  context) whose objects are finitely generated left $U(\mathfrak{g})$-modules for which the $\mathfrak{b}$-action is locally finite and the $\mathfrak{h}$-action is semi-simple (see \cite{MR2428237}, chapter 1). 

\begin{defn}
  Following \cite[Defn.\ 3.6.2]{MR3072116}, we let $\hat{\ocal}$ be
  the category whose objects are coadmissible
  $\hat{U}(\mathfrak{g})$-modules~$M$ for which the action of
  $\hat{U}(\mathfrak{h})$ is diagonalizable and the following
  properties hold:

  \begin{enumerate}
  \item All weights of~$M$ are contained in finitely many subsets of
    the form $\lambda + \mathbb{N}[\Phi^{-}]$, and
  \item all weight spaces of~$M$ are finite dimensional.
  \end{enumerate}
\end{defn}

\begin{thm}[\cite{MR3072116}, Thm.\ 4.3.1]\label{thm-schmidt}  We have an equivalence of categories: 
\begin{eqnarray*}
\ocal & \rightarrow &\hat{\ocal} \\
M & \mapsto &  \hat{U}(\mathfrak{g})\otimes_{U(\mathfrak{g})}M.
\end{eqnarray*}
A quasi-inverse to this functor is given by the functor: $M \mapsto M^{\sss}$, which takes $M$ to the direct sum of its weight spaces. 
\end{thm}

\subsection{Lie algebra cohomology and homology} \label{subsec:Lie-algebra-cohomology}
\subsubsection{Definitions}
Recall that    $\Mod^{\delta}(E)$ is the usual category of $E$-vector spaces,
and $\Mod(E)$ is the category of solid $E$-vector spaces. We let
$D(\mathrm{Mod}^\delta(E))$ be the derived category of $\Mod^{\delta}(E)$, and
we let $D(\Mod(E))$ be the derived category of $\Mod(E)$. Let $\mathfrak{g}$ be
a    Lie algebra (not necessarily reductive)  with enveloping Lie algebra
$U(\mathfrak{g})$. We  let  $\Mod(U(\mathfrak{g}))$ be the category of
(discrete) $U(\mathfrak{g})$-modules, and let $D(\Mod(U(\mathfrak{g}))$ be its derived category. %

  We have a   functor  ``homology of $\mathfrak{g}$'':
\begin{eqnarray*}
E\otimes^L_{U(\mathfrak{g})} -: D( \Mod(U(\mathfrak{g}))) & \rightarrow& D(\mathrm{Mod}^{\delta}(E))\\
M & \mapsto & E\otimes^L_{U(\mathfrak{g})} M
\end{eqnarray*}
We let $\HH_i(\mathfrak{g}, M):= \HH^{-i}( E\otimes^L_{U(\mathfrak{g})} M)$.

We also have a functor  ``cohomology of $\mathfrak{g}$'':

\begin{eqnarray*}
\mathrm{R}\Gamma(\mathfrak{g}, -): D(  \Mod(U(\mathfrak{g})) ) & \rightarrow& D(\mathrm{Mod}^{\delta}(E))\\
M & \mapsto & \RHom_{\mathfrak{g}}(E, M)
\end{eqnarray*}%

These functors can be  computed by taking the Chevalley--Eilenberg resolution $CE(E)$ of $E$, in cohomological degrees $[-d, 0]$ with $d = \mathrm{dim}(\mathfrak{g})$ (see \cite[sect. 7]{weibel}):
$$ 0 \rightarrow U(\mathfrak{g}) \otimes \Lambda^d\mathfrak{g} \rightarrow \cdots \rightarrow U(\mathfrak{g}) \rightarrow  0.$$ 
We can also define functors:

\begin{eqnarray*}
E\otimes^L_{U(\mathfrak{g})} -: D(\Mod(\hat{U}(\mathfrak{g})) ) & \rightarrow& D(\Mod(E))\\
M & \mapsto & E\otimes^L_{U(\mathfrak{g})} M
\end{eqnarray*}

\begin{eqnarray*}
\mathrm{R}\Gamma(\mathfrak{g}, -): D( \Mod(\hat{U}(\mathfrak{g})) ) & \rightarrow& D(\Mod(E))\\
M & \mapsto & \RHom_{\mathfrak{g}}(E, M)
\end{eqnarray*}

These functors can   also be computed by taking the solid Chevalley--Eilenberg resolution $\hat{U}(\mathfrak{g}) \otimes_{U(\mathfrak{g})}CE(E)$ of $E$ (which remains a resolution of $E$ by Corollary \ref{cor-exact-U(g)-to-hat}):

$$ 0 \rightarrow \hat{U}(\mathfrak{g}) \otimes \Lambda^d\mathfrak{g} \rightarrow \cdots \rightarrow \hat{U}(\mathfrak{g}) \rightarrow  0.$$ 

We sometimes write respectively $\RHom_{U(\mathfrak{g})}$ or  $\RHom_{\hat{U}(\mathfrak{g})}$
in place of  $\RHom_{\mathfrak{g}}$.
\begin{rem}\label{rem:trivial-reln-homology-cohomology} We have the following trivial relation between homology and cohomology:
$\underline{\mathrm{RHom}}(E\otimes^L_{U(\mathfrak{g})} M,E) = \RHom_{U(\mathfrak{g})}(E, \underline{\mathrm{RHom}}(M,E))$. 
We also have the following relation between homology and cohomology (\cite{MR0276285}):
\numequation\label{eqn:Poincare-duality-Lie-algebra}E\otimes^L_{U(\mathfrak{g})} (M[-d] \otimes_E
\Lambda^d\mathfrak{g}^\vee) = \RHom_{U(\mathfrak{g})}(E,
M).\end{equation} %
\end{rem}
We have the following well-known lemma.
\begin{lem}\label{lem:OGe-is-acyclic}$\cO_{G,e}$ is an acyclic $\hat{U}(\mf{g})$-module.  
\end{lem}
\begin{proof}
This is a simple consequence of the standard relationship between the
Chevalley--Eilenberg resolution and the de Rham complex, \emph{cf.}\ \cite[Prop.\ 5.12]{MR4475468}.
\end{proof}

\subsubsection{Homology and cohomology of $\mathfrak{g}$ and $\mathfrak{m}$}\label{sect-homo-cohom}
For the rest of this section we put ourselves in the situation of
Section~\ref{subsubsec-reductive-groups-notation}, so that in particular
$\mathfrak{g}$ is a reductive Lie algebra with Cartan and Borel
$\mathfrak{h}\subseteq\mathfrak{b}\subseteq\mathfrak{g}$, and
 $\mathfrak{p}\supseteq \mathfrak{b}$ is a standard parabolic with Levi $\mathfrak{m}\supseteq\mathfrak{h}$.  %
 We let $d=\dim\mathfrak{u}_\mathfrak{p}$.  %
 From now on until the end of Section~\ref{subsec:Lie-algebra-cohomology}, we fix a $w\in \WM$ and consider the conjugates $\mathfrak{p}_w,\mathfrak{m}_w,\mathfrak{u}_{\mathfrak{p}_w},\mathfrak{b}_{\mathfrak{m}_w}$.  We note that because $w\in\WM$, $\mathfrak{b}_{\mathfrak{m}_w}=\mathfrak{b}\cap\mathfrak{m}_w$.

We can define the homology functor of $\mathfrak{u}_{\mathfrak{p}_w}$: 
\begin{eqnarray*}
E\otimes^L_{U(\mathfrak{u}_{\mathfrak{p}_w})} -: D^-( \Mod(U(\mathfrak{g})) ) & \rightarrow& D^-(\Mod(U(\mathfrak{m}_w)))\\
M & \mapsto & E\otimes^L_{U(\mathfrak{u}_{\mathfrak{p}_w})} M
\end{eqnarray*}

\begin{rem} \label{rem: proj res as Up module}  This functor is defined by taking a projective resolution of $M$ as a $U(\mathfrak{p}_w)$-module. By the PBW theorem, $U(\mathfrak{p}_w)$ is free over $U(\mathfrak{u}_{\mathfrak{p}_w})$, and so this is also a projective resolution of $M$ as a $U(\mathfrak{u}_{\mathfrak{p}_w})$-module. 
We also have a natural functor $D(\Mod(U(\mathfrak{m}_w))) \rightarrow D(\Mod(U(\mathfrak{p}_w)))$. If we resolve $E$ via the Chevalley--Eilenberg resolution 
\[U(\mathfrak{u}_{\mathfrak{p}_w}) \otimes \Lambda^i \mathfrak{u}_{\mathfrak{p}_w}  \rightarrow \cdots \rightarrow U(\mathfrak{u}_{\mathfrak{p}_w}) \rightarrow E\]
then we get a complex which computes $ E\otimes^L_{U(\mathfrak{u}_{\mathfrak{p}_w})} M$ in $D(\Mod(U(\mathfrak{p}_w)))$ (but on the cohomology groups, the action of $\mathfrak{p}_w$ factors through an action of $\mathfrak{m}_w$). 
\end{rem}

We similarly have a cohomology functor of $\mathfrak{u}_{\mathfrak{p}_w}$:

\begin{eqnarray*}
\mathrm{R}\Gamma(\mathfrak{u}_{\mathfrak{p}_w}, -): D^+(\Mod( U(\mathfrak{g})) ) & \rightarrow& D^+(\Mod(U(\mathfrak{m}_w)))\\
M & \mapsto & \RHom_{\mathfrak{u}_{\mathfrak{p}_w}}(E, M)
\end{eqnarray*}

\begin{rem}\label{rem: proj res as Up module2} Similarly to Remark~\ref{rem: proj res as Up module}, this functor is obtained by taking an injective resolution of $M$ as a $U(\mathfrak{p}_w)$-module. If one uses the Chevalley--Eilenberg resolution of $E$ instead,  then we obtain a complex which computes the composition of this functor with the natural functor $D(\Mod(U(\mathfrak{m}_w))) \rightarrow D(\Mod(U(\mathfrak{p}_w)))$. 
\end{rem}

We can also define functors:

\begin{eqnarray*}
E\otimes^L_{U(\mathfrak{u}_{\mathfrak{p}_w})} -: D^-( \Mod(\hat{U}(\mathfrak{g}) )) & \rightarrow& D^-(\Mod(\hat{U}(\mathfrak{m}_w)))\\
M & \mapsto & E\otimes^L_{U(\mathfrak{u}_{\mathfrak{p}_w})} M
\end{eqnarray*}

\begin{eqnarray*}
\mathrm{R}\Gamma(\mathfrak{u}_{\mathfrak{p}_w}, -): D^+( \Mod(\hat{U}(\mathfrak{g})) ) & \rightarrow& D^+(\Mod(\hat{U}(\mathfrak{m}_w)))\\
M & \mapsto & \RHom_{\mathfrak{u}_{\mathfrak{p}_w}}(E, M)
\end{eqnarray*}

\begin{rem}\label{rem:flatness-of-U-p-hat} Following Remarks \ref{rem: proj res as Up module} and \ref{rem: proj res as Up module2}, these functors are well defined because $\hat{U}(\mathfrak{p}_w)$ is flat over $\hat{U}(\mathfrak{u}_{\mathfrak{p}_w})$ so that $\hat{U}(\mathfrak{m}_w) = \hat{U}({\mathfrak{p}_w}) \otimes^L_{\hat{U}(\mathfrak{u}_{\mathfrak{p}_w})} E$.
 In order to see the flatness,  by the PBW theorem and the description of
 $\hat{U}(\mathfrak{p}_w)$ given in section  \ref{sect-dist-alg}, we find that
 $\hat{U}(\mathfrak{p}_w) = \hat{U}(\mathfrak{u}_{\mathfrak{p}_w}) \otimes_E
 \hat{U}(\mathfrak{m}_{w})$ and so it remains to note that the Fr\'echet space
 of compact type $\hat{U}(\mathfrak{m}_{w})$ is flat over $E$ %
 by Lemma \ref{lem-flat-frechet}.
\end{rem}

\subsubsection{Finiteness of the algebraic cohomology}\label{sect-finitenessliealgebracohomology}
Let $Z(\mathfrak{g})$ and $Z(\mathfrak{m}_w)$ denote the centres of
$U(\mathfrak{g})$ and $U(\mathfrak{m}_w)$ respectively, and let $W$ and
$W_{\mathfrak{m}_w}$ be the Weyl groups of~ $\mathfrak{g}$ and~
$\mathfrak{m}_{w}$. If~$M$ is a $U(\mathfrak{g})$-module on which~$Z(\mathfrak{g})$
acts via a character~$\chi$, then we say that~$\chi$ is the
\emph{infinitesimal character} of~$M$.
Recall the Harish-Chandra isomorphism $HC_{\mathfrak{g}}:
Z(\mathfrak{g}) \rightarrow U(\mathfrak{h})^{W,\cdot}$ (where the
target is the invariants for the dotted action of~$W$), determined by the property that $z \otimes 1 = 1 \otimes HC_{\mathfrak{g}}(z)$ in $U(\mathfrak{g}) \otimes_{U({\mathfrak{b}})} U(\mathfrak{h})$.
Using this isomorphism, any character $\lambda: Z(\mathfrak{g})
\rightarrow E$ is identified with an element of $\mathfrak{h}^\vee$,
well defined up to the  dotted action of $W$. 

\begin{rem}\label{rem:iota-on-HC} Let $\iota: Z(\mathfrak{g}) \rightarrow Z(\mathfrak{g})$ be the map induced by the inverse map on~ $\mathfrak{g}$, $X \mapsto -X$. We have $HC \circ \iota = -w_0 \circ HC$. If $\lambda \in \mathfrak{h}^\vee$ represents  a character of~ $Z(\mathfrak{g})$, $-w_0\lambda$ represents the character $\lambda \circ \iota$. 
\end{rem}

We similarly have a natural Harish-Chandra isomorphism $HC_{\mathfrak{m}_w}: Z(\mathfrak{m}_w)
\rightarrow  U(\mathfrak{h})^{W_{\mathfrak{m}_w},\cdot}$, %
and we deduce that  there is a natural Harish-Chandra  map \numequation\label{eqn:HCmap}HC:
  Z(\mathfrak{g}) \rightarrow Z(\mathfrak{m}_w).\end{equation} This map is
characterized by the property that in $U(\mathfrak{g})
\otimes_{U({\mathfrak{p}_w})} U(\mathfrak{m}_w)$, we have $z\otimes 1 = 1 \otimes HC(z)$. 

We also have a map from $Z(\mathfrak{g})$ (resp.\ $Z(\mathfrak{m}_w)$) to
the centre of the derived categories $D(\Mod(U(\mathfrak{g})))$ (resp.\
$D(U(\mathfrak{m}_w))$) (i.e.\ the $t$-centre in the sense of  \cite{MR3308788}). %
The following is known as the Casselman--Osborne theorem. 
\begin{thm}[\cite{MR0396704},  \cite{MR3308788}] \label{thm: CO thm}  The functor $\mathrm{R}\Gamma(\mathfrak{u}_{\mathfrak{p}_w},-): D^+(\Mod(U(\mathfrak{g}))) \rightarrow D^+(\Mod(U(\mathfrak{m}_w)))$ is $Z(\mathfrak{g})$-homogeneous, in the sense that for $z \in Z(\mathfrak{g})$, we have $$\mathrm{R}\Gamma(\mathfrak{u}_{\mathfrak{p}_w}, z) = HC(z).$$
In particular, if  $M$ is a $\mathfrak{g}$-module with infinitesimal character $\lambda \in \mathfrak{h}^\vee$, then $\HH^i(\mathfrak{u}_{\mathfrak{p}_w}, M)$ is a  $Z(\mathfrak{m}_w) \otimes_{Z(\mathfrak{g})} \lambda$-module. 
\end{thm}

\begin{thm}\label{thm:u-cohom-in-cat-O} Let $M \in \ocal(\mathfrak{g}, \mathfrak{b})$. 
Then $\HH^i(\mathfrak{u}_{\mathfrak{p}_w}, M) \in \ocal(\mathfrak{m}_w, \mathfrak{b}_{\mathfrak{m}_w})$. 
\end{thm}
\begin{proof}Using the Chevalley--Eilenberg resolution, the cohomology is computed by the complex $0 \rightarrow M \rightarrow M \otimes \mathfrak{u}_{\mathfrak{p}_w}^\vee \rightarrow \cdots $.
We see that all modules occurring in this complex have  locally
nilpotent action of $\mathfrak{u} \cap \mathfrak{m}_w$ (the unipotent radical of $\mathfrak{b}_{\mathfrak{m}_w}$) and
semi-simple $\mathfrak{h}$-action, and %
furthermore each $\mathfrak{h}$-eigenspace has finite
dimension. This also holds for the cohomology groups. It follows
that the cohomology groups admit a (possibly infinite) increasing
filtration where each graded  is a simple object in category
$\ocal(\mathfrak{m}_w, \mathfrak{b}_{\mathfrak{m}_w})$.  Indeed, if a
cohomology group is non-zero, we can find a highest weight vector
since $\mathfrak{u}\cap\mathfrak{m}_w$ acts locally nilpotently,
so we get a map from a Verma module. We can repeat the process with
the quotient.  Since all simple objects of
category~$\cO$ have a generalized infinitesimal character, it follows from the Casselman--Osborne
Theorem~\ref{thm: CO thm} that there are only a finite number of
possible infinitesimal characters of the simple  subquotients, and
therefore only finitely many
possible highest weight vectors of all irreducible subquotients. Since~$\mf{h}$
acts semi-simply with finite-dimensional eigenspaces, we
deduce that the filtration is finite and that the  cohomology groups
belong to $\ocal(\mathfrak{m}_w, \mathfrak{b}_{\mathfrak{m}_w})$, as required. 
\end{proof}
We write $D^{\mathrm{perf}}(U(\mathfrak{m}_w))$ for the category of perfect
complexes of~$U(\mf{m}_{w})$-modules. Since~$U(\mf{m}_{w})$ is Noetherian and
has finite global dimension (see e.g.\ \cite[Ex.\ 7.7.2]{weibel}), %
these are equivalently the complexes whose
cohomologies are finitely generated $U(\mf{m}_{w})$-modules, and are nonzero in
only finitely many degrees.
\begin{cor}\label{cor:coh-is-perfect-over-U} If $M \in \ocal(\mathfrak{g}, \mathfrak{b})$, then
  $\mathrm{R}\Gamma(\mathfrak{u}_{\mathfrak{p}_w}, M)$ and $E
  \otimes^L_{U(\mathfrak{u}_{\mathfrak{p}_w})} M$ %
  belong to $D^{\mathrm{perf}}(U(\mathfrak{m}_w))$, and their cohomologies
  belong to $\ocal(\mathfrak{m}_w, \mathfrak{b}_{\mathfrak{m}_w})$. %
\end{cor}%
\begin{proof}
  This is immediate from Theorem~\ref{thm:u-cohom-in-cat-O} and
  Remark~\ref{rem:trivial-reln-homology-cohomology}.
\end{proof}
  \subsubsection{Cohomology of Verma modules}

  \begin{defn}
    \label{defn:Verma-module-notation}For any~ $\lambda \in
    \mathfrak{h}^\vee$, we write $M_{\lambda}:=U(\mathfrak{g})
  \otimes_{U({\mathfrak{b}})} \lambda$ for the corresponding Verma
  module for~$\mf{g}$, and $\Mmw_{\lambda}:=U(\mf{m}_{w})
  \otimes_{U({\mathfrak{b}_{\mathfrak{m}_w}})} \lambda$ for the Verma module
  for~$\mf{m}_{w}$.  We write $L(\m_w)_\lambda$ for the simple object in $\O(\mathfrak{m}_w,\mathfrak{b}_{\mathfrak{m}_w})$ with highest weight $\lambda$.%
  \end{defn}
  
  \begin{defn} Let $\lambda \in\mathfrak{h}^\vee$. We let $W_{<\lambda}$ be the subset of $W$ consisting of elements $w'$ which satisfy:
 $w'\cdot\lambda  = \lambda - \sum_{\alpha \in \Phi^+} n_\alpha \alpha$,
 where  $n_\alpha\in\Z_{ \geq 0}$, and  $ n_\alpha >0$ for some~$\alpha$. 
\end{defn}

\begin{rem} If $\lambda \in \mathfrak{h}^\vee$ is such that $w'\cdot\lambda-\lambda \notin \Z\Phi$ for $w' \neq 1$ (a generic condition on $\lambda$), then $W_{<\lambda} = \emptyset$. 
\end{rem}

\begin{thm}\label{thm-algebraic-computationStrict} Let $\lambda \in
  \mathfrak{h}^\vee$, and let $w\in{}^MW$. Assume that~$\mf{u}_\mathfrak{p}$ is
  abelian. %
\begin{enumerate}%
\item The groups $\HH_i(\mathfrak{u}_{\mathfrak{p}_w}, M_\lambda)$ belong to the category $\ocal(\mathfrak{m}_w,{\mathfrak{b}}_{\mathfrak{m}_w})$. 
\item These homology groups vanish if $i>d-\ell(w)$.  %
\item There is an injective  ``highest weight'' map  $$ M(\mf{m}_w)_{w^{-1}(w\cdot\lambda + 2 \rho^M)} \hookrightarrow \HH_{d-\ell(w)}(\mathfrak{u}_{\mathfrak{p}_w}, M_\lambda).$$ 
\item The cokernel of the  map $ M(\mf{m}_w)_{w^{-1}(w\cdot\lambda + 2 \rho^M)} \hookrightarrow \HH_{d-\ell(w)}(\mathfrak{u}_{\mathfrak{p}_w}, M_\lambda)$ and the homology groups $\HH_{i}(\mathfrak{u}_{\mathfrak{p}_w}, M_{\lambda})$ for $i < d-\ell(w)$  have
Jordan-- H\"older factors among the
$L(\mathfrak{m}_w)_{w^{-1}(w'\cdot\lambda + 2 \rho^M)} $
 with $w' \in w W_{<\lambda}$.
\end{enumerate}
\end{thm}

\begin{rem} In particular, if $W_{< \lambda} = \emptyset$, then $\HH_i(\mathfrak{u}_{\mathfrak{p}_w}, M_\lambda)$ is concentrated in degree $d-\ell(w)$ and $\HH_{d-\ell(w)}(\mathfrak{u}_{\mathfrak{p}_w}, M_{\lambda})  = M(\mf{m}_w)_{w^{-1}(w\cdot\lambda + 2 \rho^M)}$. 
\end{rem}

\begin{proof}[Proof of Theorem~\ref{thm-algebraic-computationStrict}]By Corollary~\ref{cor:coh-is-perfect-over-U} ,
    the homologies belong to  the category
  $\ocal(\mathfrak{m}_w,{\mathfrak{b}}_{\mathfrak{m}_w})$. %
Since $M_\lambda$ is free as a $U(\bar{\mathfrak{u}})$-module and so also as a $U(\mathfrak{u}_{\mathfrak{p}_w}\cap\bar{\mathfrak{u}})$-module,  there is an isomorphism \numequation\label{eqn:homology-of-Verma-b-and-b-op}\HH_\star(\mathfrak{u}_{\mathfrak{p}_w},
M_\lambda) =  \HH_\star
(\mathfrak{u}_{\mathfrak{p}_w} \cap \mathfrak{u},  \HH_0(\mathfrak{u}_{\mathfrak{p}_w} \cap
\bar{\mathfrak{u}},  M_\lambda)).\end{equation}(Here we used
that~$\mf{u}_{\mf{p}_w}$ is abelian, since $\mf{u}$ is abelian by assumption.)
Write $N:= \HH_0(\mathfrak{u}_{\mathfrak{p}_w} \cap
\bar{\mathfrak{u}},  M_\lambda)$. %
 The homology~\eqref{eqn:homology-of-Verma-b-and-b-op} is computed by  the Chevalley--Eilenberg complex: 
$$ 0 \rightarrow \wedge^{d-\ell(w)} (\mathfrak{u}_{\mathfrak{p}_w} \cap \mf{u}) \otimes N \rightarrow \cdots \rightarrow N \rightarrow 0$$
and so in particular vanishes in degree bigger than $d-\ell(w)=\dim(\mathfrak{u}_{\mathfrak{p}_w} \cap \mf{u})$.  Moreover the highest weight occurring in this complex is
\numequation\label{eqn:two-ways-to-write-the-weights}\lambda + w^{-1} w_{0,M} \rho + \rho=w^{-1}(w\cdot\lambda + 2
                                               \rho^M)\end{equation} (the equality
                                               holding because we have
  assumed that~$w\in{}^MW$), which occurs exactly once in  $\wedge^{d-\ell(w)} (\mathfrak{u}_{P_w} \cap \mf{b}) \otimes N$. It follows that  there is a natural map:
$$M(\mathfrak{m}_w)_{w^{-1}(w\cdot\lambda + 2 \rho^M)} \rightarrow \wedge^{d-\ell(w)}
(\mathfrak{u}_{\mathfrak{p}_w} \cap \mf{u}) \otimes N$$ which induces the
                                               (necessarily injective) map of
                                               part (3) \[
M(\mf{m}_w)_{w^{-1}(w\cdot\lambda + 2 \rho^M)} \hookrightarrow
\HH_{d-\ell(w)}(\mathfrak{u}_{\mathfrak{p}_w}, M_\lambda).\]

It remains to prove~(4). By the
  Casselman--Osborne Theorem~\ref{thm: CO thm}, together with~\eqref{eqn:Poincare-duality-Lie-algebra}, we see that the
  possible infinitesimal characters of simple subquotients in the
  homology belong to the set $\{w^{-1}(w'\cdot\lambda + 2\rho^M), w '
  \in \WM\}$.  %
Thus the simples which can occur are the $L(\mathfrak{m}_w)_{w^{-1}(w'\cdot\lambda+2\rho^M)}$.
 Moreover in order for $L(\mathfrak{m}_w)_{w^{-1}(w'\cdot\lambda + 2 \rho^M)} $ to
 occur as a subquotient of $\HH_i(\mathfrak{u}_{\mathfrak{p}_w}, M_\lambda)$ for $ i < d-\ell(w)$ or of
 $\HH_{d-\ell(w)}(\mathfrak{u}_{\mathfrak{p}_w}, M_\lambda)/ M(\mf{m}_w)_{w^{-1}(w\cdot\lambda + 2
   \rho^M)} $ we must have that $w^{-1}(w'\cdot\lambda + 2 \rho^M)$  is of the form $w^{-1}(w\cdot\lambda + 2 \rho^M)-\sum n_\alpha\alpha$ with $n_\alpha\geq0$ and some $n_\alpha\not=0$.   Therefore $w' \in
 wW_{<\lambda}$. %
\end{proof}

\begin{rem} One can show that the highest weight map computes 
  the Euler characteristic of $\HH_{\star}(\mathfrak{u}_{\mathfrak{p}_w},
  M_\lambda)$ in the Grothendieck group, see Proposition \ref{prop-Groth-group}.%
\end{rem}

\subsubsection{Strictness}

We now state our main theorem on the comparison of algebraic and solid
cohomology of Lie algebras (Theorem~\ref{thm-STrict}). 
\begin{defn}\label{defn-number-non-Liouville} We say $x\in E$ is \emph{$p$-adically non-Liouville}, or
  simply \emph{non-Liouville}, if $\liminf_{r\in\Z_{>0}}|x-r|^{1/r}\not=0$.\end{defn}

\begin{remark}[Inconsistencies in the literature concerning the definition of non-Liouville]
There are  a  number of conflicting definitions in the literature of 
what it means for~$\alpha \in \Qbar_p$ to be $p$-adically non-Liouville.
The original definition in~\cite[Def~1]{Clark} is equivalent to the existence of
a real number~$d$ such that
\numequation \label{clarkversion}
|\alpha - r| \gg |r|^{-d},  \quad r \in \Z, \ r \rightarrow \infty.
\end{equation}
This is the most direct analogue
of the definition over~$\R$ --- Liouville's original argument shows that 
any~$\alpha \in \Qbar \subset \Qbar_p$ satisfies~(\ref{clarkversion}) with~$d = [\Q(\alpha):\Q]$.
There is a weaker definition of non-Liouville
given in~\cite{Adolphson}, which is
equivalent to
\numequation
\forall B \in (0,1), \quad  |\alpha - r| \gg B^{|r|},  \quad r \in \Z, \ r \rightarrow \infty. \end{equation}
Definition~\ref{defn-number-non-Liouville}, following Pan~\cite[Rem~5.2.11]{MR4390302}
 (see also~\cite[Def.~1]{MR3053445})
 is weaker still, and is equivalent to
\numequation \label{formulaL}
\exists
B \in (0,1), \quad  |\alpha - r| \gg B^{|r|},  \quad r \in \Z, \ r \rightarrow \infty.\end{equation} 
Despite these differences, 
both~\cite{Adolphson} and~\cite{MR4390302} attribute their respective definitions
to~\cite{Clark}.
 The definition given in Kedlaya's book~\cite[\S13]{Kedlaya}
is equivalent to the one in~\cite{Adolphson}.
The reason we use the definition in~\cite{MR4390302} is because (for a number of arguments) our
results are true \emph{if and only if} $\lambda$ is non-Liouville in the sense of
 Definition~\ref{defn-number-non-Liouville}, although we do not stress this point.
 As a practical matter, however, the reader should feel free to take any definition they like, since:
\begin{enumerate}
\item In applications, the strongest assumption from~\cite{Clark} is satisfied.
\item In  proofs, only the weakest assumption from~\cite{MR4390302} will be assumed.
\end{enumerate}
 \end{remark}

\begin{defn}\label{defn-weight-non-Liouville}
We say that a weight $\lambda\in \mathfrak{h}^\vee$ is $p$-adically
non-Liouville if $\langle\lambda+\rho,\alpha^\vee\rangle$ is
$p$-adically non-Liouville for all $\alpha\in\Phi^+$.  Again, we
will often abbreviate ``$p$-adically non-Liouville'' to
``non-Liouville''.
\end{defn}

\begin{rem} Any integer is $p$-adically non-Liouville, and if
  $(\mathfrak{g}, \mathfrak{h}, \mathfrak{b})$ arises as the Lie
  algebras of $(G,T,B)$ with $G$  a reductive group, then any
  algebraic weight $\lambda \in X^\star(T) \subseteq \mathfrak{h}^\vee
  = X^\star(T)_E$ is $p$-adically non-Liouville. %
\end{rem}

We let $\ocal_{nL}(\mathfrak{g}, \mathfrak{b}) \subseteq
\ocal(\mathfrak{g}, \mathfrak{b})$ be the direct factor abelian
subcategory consisting of objects whose weights are $p$-adically non-Liouville.

\begin{thm}\label{thm-STrict}  Let $M \in \ocal_{nL}(\mathfrak{g}, \mathfrak{b})$. The canonical map 
$$\hat{U}(\mathfrak{m}_w)\otimes_{U(\mathfrak{m}_w)}(E \otimes^L_{U(\mathfrak{u}_{\mathfrak{p}_w})} M)  = \hat{U}(\mathfrak{p}_w) \otimes_{U(\mathfrak{p}_w)}(E \otimes^L_{U(\mathfrak{u}_{\mathfrak{p}_w})} M) \rightarrow E \otimes^L_{\hat{U}(\mathfrak{u}_{\mathfrak{p}_w})} ( \hat{U}(\mathfrak{g})\otimes_{U(\mathfrak{g})} M)$$is a quasi-isomorphism.
\end{thm}%
\begin{proof}
  This is proved below as Theorem~\ref{thm-STrict2}.
\end{proof}

\begin{cor}\label{Main-cor-strictness} For all $i$,
  $\HH^i(\mathfrak{u}_{\mathfrak{p}_w}, \widehat{M}^\vee) =
  \widehat{\HH_i(\mathfrak{u}_{\mathfrak{p}_w}, M)}^\vee$ is an
  admissible $\hat{U}(\mathfrak{m}_w)$-module (where as usual we write $\widehat{M}:= M \otimes_{U(\mathfrak{g})} \hat{U}(\mathfrak{g}$). %
\end{cor}%
\begin{proof} The Chevalley--Eilenberg complex $CE(E) \otimes_E \widehat{M}$ which computes $E \otimes^L_{\hat{U}(\mathfrak{u}_{\mathfrak{p}_w})}\widehat{M}$ has the shape:
$0 \rightarrow \Lambda^d\otimes\widehat{ M}  \mathfrak{u}_{\mathfrak{p}_w} \rightarrow \cdots \rightarrow  \widehat{M} \rightarrow 0$. 
This is a complex of Fr\'echet  spaces, and its cohomology groups are also Fr\'echet spaces by  Theorem \ref{thm-STrict}. 
Moreover, we have that $\HH^{-i}(CE(E) \otimes_E \widehat{M}) = \widehat{\HH_i(\mathfrak{u}_{\mathfrak{p}_w}, M)}$.
We recall from Remark~\ref{rem: exactness of antieq Smith Banach} that  $\underline{\mathrm{Hom}}(-,E)$  is an exact functor between Fr\'echet and $\LB$-spaces.  We see in the first place that  
  $\underline{\mathrm{Hom}}(CE(E) \otimes_E \widehat{M},E) =
  \underline{\mathrm{Hom}}(CE(E), \widehat{M}^\vee)$ computes
  $\mathrm{R}\Gamma(\mathfrak{u}_{\mathfrak{p}_w}, \widehat{M}^{\vee})$ %
  and that $\HH^i (\underline{\mathrm{Hom}}(CE(E) \otimes_E \widehat{M},E)) =
  \underline{\mathrm{Hom}}( \HH^{-i}(CE(E) \otimes_E \widehat{M}),E)$. Since
  $\underline{\mathrm{Hom}}( \HH^{-i}(CE(E) \otimes_E \widehat{M}),E)  =
  \widehat{\HH_i(\mathfrak{u}_{\mathfrak{p}_w}, M)}^\vee$ by Theorem
  \ref{thm-STrict}, we are done by Theorem~\ref{thm:u-cohom-in-cat-O}. %
  \end{proof}

\subsection{Algebraic local cohomology and twisted Verma modules}  Now we fix a
split reductive group $G/E$ with Lie algebra $\mathfrak{g}$. We work over the field $E$, viewed as a discrete field (we ignore its natural $p$-adic  topology for the moment). We want to introduce twisted Verma modules as local cohomology on the flag variety.  

\begin{rem} We will make a small variation on the classical presentation since we will use the six-functor formalism in coherent cohomology of Clausen and Scholze \cite{Clausen-Scholze}, and endow the Bruhat cells with the structure of analytic stacks. We feel that this perspective clarifies the discussion. We note that most of  our statements are classical (see for example \cite{MR1985191}),  and  our   proofs can easily be translated into more classical language. 
\end{rem}

To any affine scheme $\Spec~A$, we can attach an analytic stack  $\AnSpec(A, \Mod(A))$ where $\Mod(A)$   is the category of condensed $A$-modules which are solid $\Z$-modules.  This procedure glues to  define a functor from the category of schemes to the category analytic stacks, which we denote  by $X \mapsto \tilde{X}$. 

For any  Zariski open  subset $\Spec~A[1/f]$, the corresponding map $$i:
\AnSpec(A[1/f], \Mod(A[1/f])) \rightarrow \AnSpec(A, \Mod(A))$$ is proper (!)
and can be regarded as a closed immersion.  The inclusion $i$ has an open
complement $j: U \rightarrow \AnSpec(A, \Mod(A))$. We can describe $U$ as the
ind-scheme equal to the formal completion of $A$ along the ideal $(f)$. Via the
morphism $j_\star$, the (derived) category of quasi-coherent sheaves on $U$
identifies with the subcategory of $D(\Mod(A))$ of modules which are derived
$(f)$-complete (i.e.\ modules $M$ which satisfy $\lim_{\times f} M = 0$).

We now  let $FL:=B \backslash G$. %
We  consider the classical Bruhat stratification $FL =
\coprod_{w \in W} C_{w}$ where $C_{w} = B \backslash Bw B$.  We let $X_w =
\overline{C_w}$ be the Schubert variety. %
We now equip each $X_{w}$ with the
 structure of an  analytic stack $\tilde{X}'_{w}$ which admits a map to
 $\tilde{X}_{w}$.  Let $Y_{w}$ be the open complement of $X_{w}$ in $FL$. Then $Y_w$ defines an analytic stack $\tilde{Y}_{w}$,   the map $\tilde{Y}_{w} \rightarrow \tilde{FL}$ is proper, and we let  $\tilde{X}'_{w}$ be its open complement. 
Its structure sheaf is the completion
$\widehat{\oscr_{FL}}^{I_{X_w}}$ where $I_{X_{w}}$ is the  ideal of $X_{w}$ in
$\oscr_{FL}$.  The corresponding category of  modules are the
$\oscr_{FL}$-modules which are  $\ZZ$-solid %
and derived complete modules for the $I_{X_{w}}$-adic topology. In other words, we are considering the formal scheme equal to the formal completion of $FL$ along $X_{w}$, which we view naturally as an object of the category of analytic stacks.  From that perspective, the map $\tilde{X}'_w \rightarrow  \tilde{FL}$ is an open immersion.  %

This induces a  structure of analytic stack $\tilde{C}'_{w}$ on  each Schubert cell $C_{w}$. Indeed, we have classically $C_{w} = X_{w} \setminus \cup_{w' \leq w} X_{w'}$ and we let $\tilde{C}'_{w} = \tilde{X}'_{w} \setminus \cup_{w' \leq w} \tilde{X}'_{w'}$. Note that  $\tilde{C}'_{w}$ is naturally closed in $\tilde{X}'_{w}$. We now simplify our notations, and denote by $FL$, $C_w$, $X_w$  the analytic stacks we just defined. 

\begin{example} We can illustrate how this works for $\mathrm{SL}_2$. In this case, we have $FL = \mathbb{P}^1$. 
We have that $C_{Id}$ has structure sheaf $E \llb T^{-1}\rrb$ and category of modules the solid $E\llb T^{-1}\rrb$-modules which are derived complete for the  $T^{-1}$-adic topology. 
We have that $C_{w_0}$ has the structure sheaf $E[T]$ and modules are the $E[T]$-modules which are solid $\Z$-modules.
\end{example}
\begin{example}\label{ex:an-stack-str-on-Cw} We can also describe $C_w$ in general. Write
  $E[T_\alpha]$ for the  underlying ring of the root group~$U_{\alpha}$.  Then $C_w$ has structure sheaf $E[T_\alpha, \alpha \in w^{-1}\Phi^- \cap \Phi^+]\llb T_\alpha, \alpha \in w^{-1}\Phi^- \cap \Phi^-\rrb$ and category of modules the solid $\Z$-modules which are $E[T_\alpha, \alpha \in w^{-1}\Phi^- \cap \Phi^+]\llb T_\alpha, \alpha \in w^{-1}\Phi^- \cap \Phi^-\rrb$-modules and are $(T_\alpha, \alpha \in w^{-1}\Phi^- \cap \Phi^-)$-derived complete. 
\end{example}%

We can in fact consider the action of  an analytic stack in groups $\hat{G} \rtimes B$  on $FL$ (via the the product $\hat{G} \rtimes B \rightarrow G$ and the obvious $G$-action on $FL$),
which is such that the $C_w$ (with their analytic structures) are the $\hat{G}
\rtimes B$-orbits. %
We begin with some definitions.

Whenever we have a classical affine algebraic group $H$ we view it as an analytic stack using the functor $H \mapsto \tilde{H}$. In other words, it is equipped with the structure sheaf $\oscr_H$ and category of modules the condensed $\oscr_H$-modules which are   solid $\Z$-modules.  Similarly, we define $\hat{H}$ (the completion at identity) with
structure sheaf $\widehat{\oscr_{H,e}}$ (the completed structure sheaf
at the identity) and modules the solid $\ZZ$-modules which are
$\widehat{\oscr_{H,e}}$-modules  and are  derived complete modules for the
$\mathfrak{m}_{\oscr_{H,e}}$-adic topology. Note that $\hat{H}$ is the complement of $H \setminus \{e\}$  (where again $H
\setminus \{e\}$ is equipped with its structure sheaf $\oscr_{H \setminus
  \{e\}}$ and category of modules all the solid $\ZZ$-modules which are
$\oscr_{H \setminus \{e\}}$-modules).  Note also that $\hat{H}$ is open in $H$ and $H \setminus \{e\}$ is closed. 

The relevance of the group $\hat{G}$ is clarified by the following lemma. 

\begin{lem}\label{lem-repofhatG} The category of  representations of
  $\hat{G}$ is naturally equivalent to the category of
  $U(\mathfrak{g})$-modules on solid $E$-vector spaces. %
\end{lem}
\begin{proof}Write $\pi: \hat{G} \rightarrow \Spec~E$, so that a representation
  of $\hat{G}$ is a  solid $E$-module~$M$, together with a comodule map $M
  \rightarrow \pi^\star M$ satisfying the usual cocycle condition. We have
  $\pi^\star M = \mathrm{R}\lim_n M \otimes_E \oscr_G /
  \mathfrak{m}_{\oscr_{G,e}}^n$, so that the map $M \rightarrow \pi^\star M$ is equivalent to the data of compatible maps $M \rightarrow M \otimes_E \oscr_G / \mathfrak{m}_{\oscr_{G,e}}^n$, which dually corresponds to a map $U(\mathfrak{g}) \otimes_{E} M \rightarrow M$. 
\end{proof}

The following computation will be used repeatedly.

\begin{lem}\label{lemkey-comput} Let $U_\alpha$ be a root group, with underlying ring $E[T_\alpha]$. 
Then we have $\mathrm{R}\Gamma_c(U_\alpha, \oscr_{U_\alpha}) = E[T_\alpha]$ and $\mathrm{R}\Gamma_c(\hat{U}_\alpha, \oscr_{U_\alpha}) = E[T_\alpha, T_\alpha^{-1}]/E[T_\alpha][-1]$.
Moreover, $E \otimes^L_{U(\mathfrak{u}_{\alpha})} \mathrm{R}\Gamma_c(U_\alpha, \oscr_{U_\alpha}) = E(\alpha)[1]$ and $E \otimes^L_{U(\mathfrak{u}_{\alpha})} \mathrm{R}\Gamma_c(\hat{U}_\alpha, \oscr_{U_\alpha}) = E(-\alpha)[-1]$.

\end{lem}
\begin{proof} When regarded as above as an analytic stack, %
$U_\alpha$ is
proper over $\Spec~E$ so that  $\mathrm{R}\Gamma_c(U_\alpha, \oscr_{U_\alpha}) =
E[T_\alpha]$.  We can then compute 
$\mathrm{R}\Gamma_c(\hat{U}_\alpha, \oscr_{U_\alpha}) = E[T_\alpha, T_\alpha^{-1}]/E[T_\alpha][-1]$ by using the triangle:
$$\mathrm{R}\Gamma_c(\hat{U}_\alpha, \oscr_{U_\alpha})  \rightarrow \mathrm{R}\Gamma({U}_\alpha, \oscr_{U_\alpha}) 
\rightarrow \mathrm{R}\Gamma (\hat{U}_\alpha \setminus \{e\}, \oscr_{U_\alpha}) \stackrel{+1}\rightarrow .$$
Next, $E \otimes^L_{U(\mathfrak{u}_{\alpha})} \mathrm{R}\Gamma_c(U_\alpha,
\oscr_{U_\alpha})$ is computed by the (Chevalley--Eilenberg) complex in degrees $-1$ and $0$: $E[T_\alpha] \otimes \mathfrak{u}_\alpha \rightarrow E[T_\alpha]$ with basis vector $u_\alpha$ of $\mathfrak{u}_\alpha$ acting by the derivation $\partial_{T_\alpha}$. It is thus $E(\alpha)[1]$. 
Similarly, $E \otimes^L_{U(\mathfrak{u}_{\alpha})} \mathrm{R}\Gamma_c(\hat{U}_\alpha, \oscr_{U_\alpha})$ is computed by the complex in degrees $0$ and $1$: $E[T_\alpha,T_\alpha^{-1}]/E[T_\alpha] \otimes \mathfrak{u}_\alpha \rightarrow E[T_\alpha, T_\alpha^{-1}]/E[T_\alpha]$. It is thus $E(-\alpha)[-1]$. 
\end{proof}

We check that the semi-direct product $\hat{G} \rtimes B$ is well defined (i.e.\ there is an action of $B$ on $\hat{G}$). 
First, there is an action of $B$ on $G$ by conjugation (with $G$ equipped with its structure sheaf $\oscr_G$ and category of modules the solid $\ZZ$-modules which are also $\oscr_G$-modules). We observe next  that $\hat{G}$ is the complement of $G \setminus \{e\}$.  It is clear that $B$ preserves $G \setminus \{e\}$ and thus it also acts on its open complement $\hat{G}$. 

We see that each $C_w$ is a $\hat{G}\rtimes B$-orbit in $FL$. %
Therefore, we have an equivalence of categories between $\hat{G}\rtimes B$-equivariant sheaves on $C_w$  and  representations of the stabilizer $Stab(w)$ of $w$, given by the fiber functor $\mathscr{F} \mapsto \mathscr{F}\vert_w$. An inverse of this functor is given by $V \mapsto \pi_\star( \oscr_{\hat{G} \rtimes B} \otimes V)^{Stab(w)}$ where $\pi: \hat{G} \rtimes B \rightarrow C_w$ is the uniformization map. One can describe the stabilizer $Stab(w)$ of the point $w$ under this map. 

\begin{lem}\label{lem:stab-w-analytic-stack} We have $Stab(w) = \widehat{B \cap B_w} \backslash [(\hat{B}_w \times \hat{B}) \rtimes B \cap B_w]$ where the map $(\hat{B}_w \times \hat{B}) \rtimes B \cap B_w \rightarrow \hat{G} \rtimes B$ is given by $(b,b', b'') \mapsto (b(b')^{-1}, b'b'')$ and the map $\widehat{B \cap B_w} \rightarrow (\hat{B}_w \times \hat{B}) \rtimes B \cap B_w$ is given by $b \mapsto (b^{-1}, b^{-1}, b)$. 
\end{lem}%
\begin{proof} This is straightforward; see 
  Lemma~\ref{lem-decompose-stab(w)-Q} for the proof of a very similar statement.%
\end{proof} 
\begin{lem}\label{lem-triv-of-sheaf} Let $\mathscr{F}$ be a
  $\hat{G}\rtimes B$ equivariant sheaf on $C_w$. Then  $\mathscr{F}$
  is isomorphic as a $\hat{B}_{w_0w}$-equivariant sheaf to %
  $\mathscr{F}\vert_w \otimes_E \oscr_{C_w}$ (where
  the $\hat{B}_{w_0w}$-equivariant sheaf structure on~
  $\mathscr{F}$ arises from regarding  $\hat{B}_{w_0w}$ as a subgroup
  of~$\widehat{G}$; and $\hat{B}_{w_0w}$ acts on $\mathscr{F}\vert_w$ through $\hat{T}$ and
  via its natural action on $\oscr_{C_w}$). 
\end{lem}
\begin{proof}By Example \ref{ex:an-stack-str-on-Cw} $ \widehat{{B}_{w_0}\cap
    {U}_{w_0w}} \times B \cap {U}_{w_0w}$ (viewed as a substack -- but not a
  subgroup -- of $\hat{G}
  \rtimes B$) maps isomorphically to $C_w$ via the uniformization map $x\mapsto wx$. It
  follows that the product map gives an isomorphism of analytic stacks (not of groups): %
  $$Stab(w) \times ( \widehat{{B}_{w_0}\cap {U}_{w_0w}} \times B \cap {U}_{w_0w}) \rightarrow \hat{G} \rtimes B.$$ This isomorphism is equivariant for the $\hat{U}_{w_0w}$-action by translation on the right on  $\widehat{{B}_{w_0}\cap {U}_{w_0w}} \times B \cap {U}_{w_0w} $  and  $\hat{G}\rtimes B$. It is also equivariant for the $\hat{T}$-action (the one by right translation on $\hat{G}\rtimes B$, and by right translation on $Stab(w)$ and conjugation on $ \widehat{{B}_{w_0}\cap {U}_{w_0w}} \times B \cap {U}_{w_0w}$). 
We construct a map $\mathscr{F}\vert_w \rightarrow \mathscr{F}$, by sending $v \in \mathscr{F}\vert_w$ to $(ss' \mapsto sv)$ viewed as an element of  $\pi_\star( \oscr_{\hat{G} \rtimes B} \otimes \mathscr{F}\vert_w)^{Stab(w)} = \mathscr{F}$
for $s \in Stab(w)$ and $s' \in   \widehat{{B}_{w_0}\cap {U}_{w_0w}} \times B \cap {U}_{w_0w}$. This induces an isomorphism
$\mathscr{F}\vert_w \otimes_E \oscr_{C_w} \rightarrow \mathscr{F}$ which satisfies the expected properties. 
\end{proof}

Let $(\kappa,\nu) \in X^\star(T)^2_E$ with the property that $\kappa +
\nu \in X^\star(T)$. We define a character of %
$Stab(w)$ as follows: we let $\hat{B}_w$ act via $\kappa$, we let $\hat{B}$ act via $\nu$, and we let  $B \cap B_w$ act via $\kappa + \nu$. This defines a $\hat{G} \rtimes B$ equivariant sheaf $\mathcal{L}_\kappa(\nu)$ over $C_w$. We sometimes drop $\nu$ from the notation since we are mostly interested in the $\hat{G}$-equivariant action (the $B$-action rigidifies the construction and will be used in the construction of intertwining maps). 
We let $j_w: C_w \rightarrow {FL}$ be the inclusion. Let $d_{FL}$ be the dimension of $FL$.  We can now define the twisted Verma modules:

 \begin{defn}\label{defn:twisted-Verma} We define the twisted Verma module $M(\mf{g})_\lambda^w =
   \HH^{d_{FL}-\ell(w)}(FL, (j_w)_! \mathcal{L}_{\lambda + w^{-1} \rho +
     \rho}(\nu))$. 
 \end{defn}
 
 This is a representation of $\hat{G} \rtimes B$. By Lemma \ref{lem-repofhatG},
 the $\hat{G}$-action amounts to a $\mathfrak{g}$-module structure. We will
 usually write  $M_\lambda^w$ for $M(\mf{g})_\lambda^w$.
 
 \begin{prop}\label{prop-computation-of-local-coho} The $\mathfrak{g}$-module $M_\lambda^w$ belongs to $\ocal(\mathfrak{g},\mathfrak{b})$. It has the following properties: 
    \begin{enumerate}
 \item Its highest weight is $\lambda$. 
 \item It is isomorphic to the direct sum
 $$ \oplus_{ k_\alpha \geq 0, \alpha \in w^{-1} \Phi^- \cap \Phi^+, k_\alpha < 0, \alpha \in w^{-1} \Phi^{-} \cap \Phi^-} E(\lambda+w^{-1}\rho + \rho)\prod T_{\alpha}^{k_\alpha}$$
 and in particular the  action of $\mathfrak{b}_{w_0 w}$ is completely explicit. 
 \item $M_\lambda^{Id}$ is the Verma module of highest weight $\lambda$. 
 \item $M_\lambda^{w_0}$ is the dual Verma module of highest weight $\lambda$.
 \item\label{item:character-of-twisted-Verma} The elements
   $[M_\lambda^w]$ of the Grothendieck group  are independent of~$w$. 
 \end{enumerate}
 \end{prop}
 \begin{proof} This is proven in the course of the proof of
   \cite[Lem.\ 3.2.2]{boxer2021higher}. %
   Let us give some details.  Given Lemma \ref{lem-triv-of-sheaf}, and the projection formula,  the key computation is that 
 $\mathrm{R}\Gamma_c(C_w, \oscr_{C_w}) =  \oplus_{ k_\alpha \geq 0, \alpha \in w^{-1} \Phi^- \cap \Phi^+, k_\alpha < 0, \alpha \in w^{-1} \Phi^{-} \cap \Phi^-} E\prod T_{\alpha}^{k_\alpha}.$
 This follows from  Lemma \ref{lemkey-comput}. We deduce that $M_\lambda^w$ is a
 finitely generated $U(\mathfrak{g})$-module and that the action of
 $\mathfrak{h}$ is semi-simple with the same character as that of the  Verma module $M^{Id}_\lambda$. 
 This implies that $M_\lambda^w$ belongs to category $\ocal$ and that $[M_\lambda^w]$ is independent of $w$ (see \cite[1.15]{MR2428237}).
 \end{proof}

 One can compute easily the homology of
 $\mathfrak{u}_{\mathfrak{p}_w}$ on $M_\lambda^{w_0w}$ as follows. We continue to
 write $d = \mathrm{dim}_E
(\mathfrak{u}_{\mathfrak{p}})$.%

 \begin{prop}\label{prop-computing-coho-twistedverma} We have \[E
   \otimes^L_{U(\mathfrak{u}_{\mathfrak{p}_w})} M_\lambda^{w_0w} =
   M(\mathfrak{m}_w)^{(w^M)^{-1}w_{0,M}w}_{(w^M)^{-1}(w^M\cdot\lambda
     + 2\rho^M)}[d-\ell(w^M)]\] where $w=w_Mw^M$ for $w_M\in W_M$ and $w^M\in\WM$. %
 \end{prop}
 \begin{proof} %
   We have a
   map \[\pi:\prod_{ \alpha \in w^{-1} \Phi^+ \cap \Phi^+} U_\alpha
     \prod_{\alpha \in w^{-1} \Phi^+ \cap \Phi^-} \hat{U}_\alpha
     \rightarrow \prod_{\alpha \in w^{-1} \Phi^+_M \cap \Phi^+}
     U_\alpha \prod_{\alpha \in w^{-1} \Phi^+_M \cap \Phi^-}
     \hat{U}_\alpha\] of analytic stacks with fiber
 $$F := \prod_{\alpha \in w^{-1}\Phi^{+,M} \cap \Phi^+} U_\alpha \prod_{\alpha \in w^{-1} \Phi^{+,M} \cap \Phi^{-}} \hat{U}_{\alpha}.$$
 This map is $\hat{P}_{w} \rtimes B\cap P_w$-equivariant (the action of $\hat{P}_w \rtimes B\cap P_w$ factors through an action $\hat{M}_w \rtimes B_{M_w}$ on the target). 
 The space $\prod_{ \alpha \in w^{-1} \Phi^+ \cap \Phi^+} U_\alpha \prod_{\alpha \in w^{-1} \Phi^+ \cap \Phi^-} \hat{U}_\alpha$ is the  Bruhat cell  $C_{w_0w}$ in  $B \backslash G$. The space  $\prod_{\alpha \in w^{-1} \Phi^+_M \cap \Phi^+} U_\alpha \prod_{\alpha \in w^{-1} \Phi^+_M \cap \Phi^-} \hat{U}_\alpha$ is the Bruhat cell $C'_{(w^M)^{-1}w_{0,M}w_Mw^M}$  in $B_{M_w} \backslash M_w$. We deduce that $$E \otimes^L_{U(\mathfrak{u}_{\mathfrak{p}_w})} M^{w_0 w}_\lambda = \mathrm{R}\Gamma_c ( C'_{(w^M)^{-1}w_{0,M}w_Mw^M}, \pi_! \mathcal{L}_{\lambda + w^{-1}w_0^{-1} \rho + \rho} \otimes^L_{U(\mathfrak{u}_{\mathfrak{p}_w})} E).$$
It therefore suffices to compute the sheaf $\pi_! \mathcal{L}_{\lambda
  + w^{-1}w_0^{-1} \rho + \rho}
\otimes^L_{U(\mathfrak{u}_{\mathfrak{p}_w})} E$. This is a $\hat{M}_w
\rtimes B_{M_w}$-equivariant sheaf, so it is determined by its fiber at~$w$.

It follows from the basic computations of Lemma
\ref{lemkey-comput}  that $\mathrm{R}\Gamma_c(F, \oscr_{F})$ is concentrated in degree $\ell(w^M)$ and equals 
$$ \oplus_{ k_\alpha \geq 0, \alpha \in  w^{-1}\Phi^{+,M} \cap \Phi^+ , k_\alpha < 0, \alpha \in  w^{-1}\Phi^{+,M} \cap \Phi^-}  E \prod_\alpha T_{\alpha}^{k_\alpha}.$$
 We then compute that

\begin{eqnarray*}
E \otimes^L_{\mathfrak{u}_{P_w}} \mathrm{R}\Gamma_c(F, \oscr_{F}) &  = & \HH_0( \mathfrak{u}_{P_w} \cap \bar{\mathfrak{b}}, \HH_{d-\ell(w^M)}(\mathfrak{u}_{P_w} \cap \mathfrak{b}, \HH_c^{\ell(w)}(F, \oscr_F))[d-2\ell(w^M)] \\
&=& E( \sum_{w^{-1}\Phi^{+,M} \cap \Phi^+} \alpha - \sum_{w^{-1}\Phi^{+,M} \cap \Phi^-} \alpha)[d-2\ell(w^M)].
\end{eqnarray*} 
It follows that $\pi_! \mathcal{L}_{\lambda + w^{-1}w_0^{-1} \rho + \rho} \otimes^L_{U(\mathfrak{u}_{\mathfrak{p}_w})} E$
 is an invertible sheaf in degree $2\ell(w^M)-d$ of  weight $$\lambda + w^{-1}w_0^{-1} \rho + \rho + \sum_{w^{-1}\Phi^{+,M} \cap \Phi^+} \alpha - \sum_{w^{-1}\Phi^{+,M} \cap \Phi^-} \alpha.$$ 
It follows that  $M^{w_0 w}_\lambda
\otimes^L_{U(\mathfrak{u}_{\mathfrak{p}_w})} E$ is concentrated in
degree $2\ell(w^M)-d + \ell(w^M)- \ell(w) = \ell(w^M)-d$, and is the twisted Verma of weight 
\[\lambda + w^{-1}w_0^{-1} \rho + \rho +  \sum_{w^{-1}\Phi^{+,M} \cap \Phi^+} \alpha - \sum_{w^{-1}\Phi^{+,M} \cap \Phi^-} \alpha - \sum_{w^{-1}\Phi^{+}_M \cap \Phi^-} \alpha$$  
$$ =\lambda + \sum_{w^{-1}\Phi^{+,M} \cap \Phi^+} \alpha 
= (w^M)^{-1}(w^M\cdot (\lambda + \rho) + 2 \rho^M).\qedhere \]
 \end{proof}%
 
 \begin{prop}\label{prop-Groth-group}  Let $w \in \WM$. In the Grothendieck group of category $\ocal$ for $\mathfrak{m}_w$, we have 
$$ \big[E \otimes^L_{U(\mathfrak{u}_{\mathfrak{p}_w})} M_\lambda^{Id}\big]=(-1)^{d-\ell(w)}\big[M(\mathfrak{m}_w)_{w^{-1}(w\cdot\lambda + 2 \rho^M)}^{Id}\big].$$
\end{prop}
\begin{proof} In the Grothendieck group of $\ocal(\mathfrak{g},\mathfrak{b})$ we have $[M_\lambda^{Id} ] =[ M^{w_0w_{0,M}w}_\lambda]$ and by Proposition \ref{prop-computing-coho-twistedverma} we have $E\otimes_{U(\mathfrak{u}_{\mathfrak{p}_w})}^LM_\lambda^{w_0w_{0,M}w}=M(\mathfrak{m}_w)_{w^{-1}(w\cdot\lambda + 2 \rho^M)}^{Id}[d-\ell(w)]$.
\end{proof}

\subsection{Some $\mathrm{SL}_2$-computations}
We now make some explicit calculations in the case $G = \mathrm{SL}_2$.  %
We let $H, X, \bar{X}$ be the standard basis of $\mathfrak{g}$  with $E\cdot H
\oplus E\cdot X = \mathfrak{b}$ and $[X,\bar{X}] = H$.   Let $\lambda :
\mathfrak{h} \rightarrow E$ be a character (identified with its value
$\lambda(H) \in E)$, and write $E(\lambda)$ for the underlying representation. %
The Verma module~$M_{\lambda}=U(\mathfrak{g}) \otimes_{U(\mathfrak{b})} E(\lambda)$  has basis the $\{\bar{X}^n\}_{n \geq 0}$ (or more precisely, $\bar{X}^n \otimes 1$  where $1 \in E(\lambda)$ is a basis  vector). 

We let $M_\lambda^\vee$ be the dual Verma module (the dual  in category $\ocal$).  Concretely, $M_\lambda^\vee$ is the subspace of the algebraic dual $\mathrm{Hom}_E(M_\lambda, E)$ which has basis the vectors $\{(\bar{X}^n)^\star\}_{n \geq 0}$ where $(\bar{X}^n)^\star( \bar{X}^{m}) = 1 $ if $m=n$ and $(\bar{X}^n)^\star( \bar{X}^{m})=0$ if $m \neq n$. 
For $g \in \mathfrak{g}$ and $f \in M_\lambda^\vee$, we have that $g f = f(^tg-)$.

\begin{lem}\label{lem-actingVermaSL2} We have that $\bar{X}\cdot (\bar{X}^{n})^\star = (n+1)(\lambda-n) (\bar{X}^{n+1})^\star$ and ${X}\cdot (\bar{X}^{n})^\star =  (\bar{X}^{n-1})^\star$.
\end{lem}
\begin{proof} These follow from the corresponding formulas in $M_\lambda$: $X\cdot \bar{X}^{n} = n(\lambda-n+1) \bar{X}^{n-1}$ and $\bar{X}\cdot \bar{X}^n = \bar{X}^{n+1}$. 
\end{proof}

\begin{cor}\label{coro-intertwining}\leavevmode \begin{enumerate}
\item There is a  unique  map of $U(\mathfrak{g})$-modules,  $ I : M_\lambda \rightarrow M_\lambda^\vee$ which sends $\bar{X}^n$ to $ n! \lambda(\lambda-1)\cdots (\lambda-(n-1)) (\bar{X}^n)^\star$.  Any other map of $U(\mathfrak{g})$-modules is a $E$-multiple of this map. 
\item If $\lambda \notin \Z_{\geq 0}$, the map  $M_\lambda \rightarrow M_\lambda^\vee$ is an isomorphism. 
If $\lambda \in \Z_{\geq 0}$, we have a long exact sequence:
$$ 0 \rightarrow M_{-2-\lambda} \rightarrow M_{\lambda} \rightarrow M^\vee_\lambda \rightarrow M^\vee_{-2-\lambda} \rightarrow 0.$$ 
The map $M_{-2-\lambda} \rightarrow M_\lambda$ sends the basis vector $\bar{X}^n$ of $M_{-2-\lambda}$ to $\bar{X}^{n+\lambda+2}$ in $M_\lambda$. The  map $M^\vee_\lambda \rightarrow M^\vee_{-2-\lambda} $ is dual to the map  $M_{-2-\lambda} \rightarrow M_{\lambda}$. 
\end{enumerate}
\end{cor}
\begin{proof} Giving a map  $M_\lambda \rightarrow M_\lambda^\vee$ of
  $U(\mathfrak{g})$-modules amounts to giving a map of $\mathfrak{b}$-modules,
  $E(\lambda) \rightarrow M_\lambda^\vee$. Since $M_\lambda^\vee$ has a unique
  vector of weight~ $\lambda$, namely $(\bar{X}^0)^\star$,  the space of maps is
  one dimensional, generated by the map $\bar{X}^0 \mapsto
  (\bar{X}^0)^\star$. Then we see by Lemma \ref{lem-actingVermaSL2} that
  $\bar{X}^n \mapsto \bar{X}^n. (\bar{X}^0)^\star = n! \lambda(\lambda-1)\cdots
  (\lambda-(n-1)) (\bar{X}^n)^\star$. If $\lambda \notin \ZZ_{\geq 0}$, this map
  is an isomorphism. Otherwise, let $L_\lambda$ be the finite 
  dimensional irreducible representation of highest weight~$\lambda$ (and dimension $\lambda+1$). 
There is a surjective map $M_\lambda \rightarrow L_\lambda \rightarrow 0$, fitting in an exact sequence:
$$  0 \rightarrow M_{-2-\lambda} \rightarrow M_{\lambda} \rightarrow L_\lambda \rightarrow 0$$
The dual in category $\ocal$ of this exact sequence gives:
$$0 \rightarrow L_\lambda \rightarrow M^\vee_\lambda \rightarrow M^\vee_{-2-\lambda} \rightarrow 0$$
and combining these exact sequences concludes the proof. 
\end{proof}

 \subsection{Intertwining maps}
 Let $\beta$ be a simple root.   Let us consider the corresponding parabolic
 ~$P_\beta$,   and  the partial flag variety $FL_\beta := P_\beta \backslash G$. We have a map $ \pi_\beta: FL \rightarrow FL_\beta$ which is a $\mathbb{P}^1$-fibration. 
 
 For each $w \in W$, we let $D_w = P_\beta \backslash P_\beta w B$ be
  the corresponding Bruhat cell. As in Example~\ref{ex:an-stack-str-on-Cw}, it is equipped with the following
  analytic stack structure. Its structure ring is $E \otimes_\Z
  \Z[T_\alpha, \alpha \in w^{-1} \Phi^{-,M_\beta} \cap
  \Phi^+]\llb T_\alpha, \alpha \in w^{-1} \Phi^{-,M_\beta} \cap
  \Phi^-\rrb$. Its modules are solid $\Z$-modules which are $E\otimes_\Z
  \Z[T_\alpha, \alpha \in w^{-1} \Phi^{-,M_\beta} \cap
  \Phi^+]\llb T_\alpha, \alpha \in w^{-1} \Phi^{-,M_\beta} \cap
  \Phi^-\rrb$-modules and are $(T_\alpha, \alpha \in w^{-1}
  \Phi^{-,M_\beta} \cap \Phi^-)$-derived complete. Assume from  now on that
  $\ell(s_\beta w) = \ell(w)+1$; then
we have $\pi_\beta^{-1} (D_w) = C_w \cup C_{s_\beta w}$. %
  Each $D_w$ is a $\hat{G} \rtimes B$-orbit. We let  $Stab(w)_\beta$ be the stabilizer of $w$. 
  We again have an equivalence between $Stab(w)_\beta$-representations and
  $\hat{G} \rtimes B$-equivariant sheaves on $D_w$. %

 \begin{lem}\label{lem:Stabw-beta} ${Stab(w)_\beta} = \widehat{B \cap P_{s_\beta w}} \backslash [ (\hat{P}_{s_\beta w} \times \hat{B}) \rtimes B \cap P_{s_\beta w}]$. 
 \end{lem}
  \begin{proof}The same as Lemma~\ref{lem:stab-w-analytic-stack}.
  \end{proof}
 
 \begin{lem}\label{lem-trivial-D_w} Let $\mathscr{F}$ be  an $\hat{G} \rtimes B$-equivariant sheaf on $D_w$.  There is a $\hat{T}$-equivariant isomorphism $\mathscr{F}\vert_w \otimes \oscr_{D_w} \rightarrow \mathscr{F}$.
 \end{lem}
 \begin{proof} This is the same as Lemma \ref{lem-triv-of-sheaf}.
 \end{proof}
 
 Given any pair of characters $(\lambda, \nu) \in X^\star(T)_E^2$ with $\lambda + \nu \in X^\star(T)$, one can construct  representations 
 $M_{\lambda}(\nu)$ and $M^\vee_\lambda(\nu)$ of $Stab(w)_\beta$ as follows: %
the underlying representation of $\hat{P}_{s_\beta w}$
  factors through $\hat{M}_{s_\beta w}$, and is respectively given
  by~$M(\mf{m}_{s_\beta w})_{\lambda}$ or~$M(\mf{m}_{s_\beta w})_{\lambda}^{\vee}$; and we let $\hat{B}$ act via
  $\nu$. The product of these actions integrates to an action of $B \cap
  P_{s_\beta w}$. %
 
 By Corollary \ref{coro-intertwining}, we see that we have intertwining maps:
 $I:  M_{\lambda}(\nu) \rightarrow M^\vee_\lambda(\nu).$
 
 \begin{lem} We have $$\mathrm{R}\Gamma_c(D_w, \oscr_{D_w}) =  \oplus_{ k_\alpha \geq 0, \alpha \in w^{-1} \Phi^{-,M_\beta} \cap \Phi^+, k_\alpha < 0, \alpha \in w^{-1} \Phi^{-,M_\beta} \cap \Phi^-} E\prod T_{\alpha}^{k_\alpha}.$$
 \end{lem}
 \begin{proof} This  follows from Lemma \ref{lemkey-comput}.
 \end{proof}
 
 \begin{prop}\label{prop-intertwining}  Assume that  $\ell(s_{\beta} w) = \ell(w) +1$. There is an intertwining map of $\mathfrak{g}$-modules: $$M_\lambda^{w} \rightarrow M_{\lambda}^{s_\beta w}.$$ This map is given (as $\mathfrak{h}$-modules)  by the map 
 
 $$ M_{\lambda + w^{-1} s_\beta \rho + \rho} \otimes_E \mathrm{R}\Gamma_c(D_w, \oscr_{D_w}) $$ $$ \stackrel{I \otimes Id}\rightarrow M^\vee_{\lambda + w^{-1} s_\beta \rho + \rho}  \otimes_E \mathrm{R}\Gamma_c(D_w, \oscr_{D_w}). $$
 
 \begin{enumerate}
 \item  If $\langle \lambda+\rho, w^{-1}\beta^\vee \rangle  \in \ZZ_{> 0}$, we have a long exact sequence:
 
 $$ 0 \rightarrow M^w_{ s_{ w^{-1} \beta}\cdot\lambda} \rightarrow M_\lambda^w \rightarrow M_{\lambda}^{s_\beta w} \rightarrow M^{s_\beta w}_{ s_{ w^{-1} \beta}\cdot\lambda} \rightarrow 0$$
  which is the tensor product of the long exact sequence of Corollary
  \ref{coro-intertwining} with  $\mathrm{R}\Gamma_c(D_w, \oscr_{D_w})$ (as
  $\mathfrak{h}$-modules). Furthermore  $M^w_{ s_{ w^{-1} \beta}\cdot\lambda}  \simeq M^{s_{\beta}w}_{ s_{ w^{-1} \beta}\cdot\lambda}$.
\item  Otherwise, the intertwining map is an isomorphism $M_\lambda^w \simeq M_{\lambda}^{s_\beta
    w}$. %
\end{enumerate}

 \end{prop}
 \begin{proof}  Consider the map $ \pi: FL \rightarrow FL_\beta$  and the maps $\pi_w: C_w \rightarrow D_w$ and $\pi_{s_\beta w}: C_{s_\beta w } \rightarrow D_w$. We construct a map:
 $\pi_{w, !} \mathcal{L}_{\lambda + w^{-1} \rho + \rho}(\nu) \rightarrow \pi_{s_\beta w, !} \mathcal{L}_{\lambda + w^{-1} s_\beta \rho + \rho}(\nu)$.
 For this, we observe that both are $\hat{G} \rtimes B$-equivariant sheaves on $D_w$. We compute the corresponding $Stab(w)_\beta$-representations. For this we can work over the fiber at $w$ by  proper base change. 
 By Proposition \ref{prop-computation-of-local-coho}, (3) and (4) we deduce that 
 the fiber  $\pi_{w, !} \mathcal{L}_{\lambda + w^{-1} \rho + \rho}\vert_w $
 corresponds to the representation $M_{\lambda + w^{-1} s_\beta \rho +
   \rho}(\nu)$, and the fiber $\pi_{s_\beta w, !} \mathcal{L}_{\lambda + w^{-1}
   s_\beta \rho + \rho}\vert_w$ corresponds to the representation
 $M^\vee_{\lambda + w^{-1} s_\beta \rho + \rho}(\nu)$.  As noted above, we get a map between these representations by using the intertwining map $I$ defined in Corollary \ref{coro-intertwining}.
 Moreover, by Lemma \ref{lem-trivial-D_w},   both sheaves are trivial, and are
 respectively $\widehat{T}$-equivariantly isomorphic  to $M_{\lambda + w^{-1} s_\beta \rho + \rho} \otimes_E \oscr_{D_w}$ and $M^\vee_{\lambda + w^{-1} s_\beta \rho + \rho} \otimes_E \oscr_{D_w}$. 
  We now take cohomology with compact support so that the cohomologies 
 $\mathrm{R}\Gamma_c (D_w, \pi_{w, !} \mathcal{L}_{\lambda + w^{-1} \rho + \rho}
 ) $ and $\mathrm{R}\Gamma_c (D_w, \pi_{s_\beta w, !} \mathcal{L}_{\lambda +
   w^{-1} s_\beta \rho + \rho} )$ are indeed given by the claimed formulas (use the projection formulas). 
 
  To see~(1) and~(2), observe that $\langle \lambda+\rho, w^{-1}\beta^\vee \rangle  \in \ZZ_{> 0}$ then again by Corollary \ref{coro-intertwining} (noting that $\lambda$ there is $\langle\lambda+w^{-1}s_\beta\rho+\rho,w^{-1}\beta^\vee\rangle=\langle\lambda+\rho,w^{-1}\beta^\vee\rangle-1$), we actually have a long exact sequence of sheaves: 
 $$0 \rightarrow \pi_{w, !} \mathcal{L}_{s_{w^{-1}\beta}\cdot\lambda + w^{-1} \rho + \rho}(\nu) \rightarrow \pi_{w, !} \mathcal{L}_{\lambda + w^{-1} \rho + \rho}(\nu) \rightarrow $$ $$ \pi_{s_\beta w, !} \mathcal{L}_{\lambda + w^{-1} s_\beta \rho + \rho}(\nu) \rightarrow  \pi_{s_\beta w, !} \mathcal{L}_{ s_{w^{-1} \beta}\cdot\lambda + w^{-1} s_\beta \rho + \rho}(\nu)\rightarrow 0$$
 inducing the expected long exact sequence on cohomology. 
 Otherwise, the intertwining map of sheaves is an isomorphism, inducing an isomorphism on cohomology. 
  \end{proof}

\subsection{Topology}

In this section we consider $E$-vector spaces $V$ equipped with a
  weight space  decomposition  $V  = \oplus_{\nu\in X^\star(T)} V_\nu$
  where each $V_\nu$ is finite dimensional and equipped with a norm
  $|-|_\nu$. We can define the norms~$|-|_{\nu}$ by choosing a basis
  for~$V_\nu$, and decreeing the basis vectors to have norm~$1$.

  Fix a basis $\{e_i\}$ of
  $X^\star(T)$; then we have a function $|-|: X^\star(T) \rightarrow
  \mathbb{N}$ measuring the size of $\nu$ as follows:  any $ \nu$ can be written as $\sum n_i e_i$ and we put $|
  \nu | = 1+\sum | n_i |$ (where in contrast to the rest of this section,
  $|n_i|$ is the archimedean norm of~$n_i$). 

Let $r \in \mathbb{R}_{>0}$. We define a norm $| - |_r$ on $V$ by letting $|
\sum v_\nu |_r = \sup_\nu | v_\nu |_\nu r^{\vert \nu \vert}$.  We write
$\hat{V}_r$ for  the Banach space completion;  concretely,
$$\hat{V}_r=\{(v_\nu)\in\prod V_\nu\mid \limsup_\nu|v_\nu|_\nu r^{|\nu|}=0\}.$$

We let $\mathcal{T}_{nat}$ be the natural topology on $V$ defined by the family
of norms $\{\vert - \vert_r\}_{r \geq 1}$, making $V$  a locally convex $E$-vector space. We let $\hat{V}_{nat}$  be $\lim_{r \rightarrow \infty} \hat{V}_r$, the completion of $V$ for $\mathcal{T}_{nat}$.   This is a Fr\'{e}chet space.

This applies in particular  to  $V=U(\mathfrak{g})$. Fixing a PBW basis gives a
decomposition of $U(\mathfrak{g})$ into weight spaces and  defines the natural
topology.  We have  $\hat{U}(\mathfrak{g}) = \hat{U}(\mathfrak{g})_{nat}$.  It
follows that any object $M$ of $\ocal$, being a finitely generated
$U(\mathfrak{g})$-module, inherits a canonical topology $\mathcal{T}_{can}$
which is a locally convex topology. For any such $M$, its  completion is
$\hat{M} = M \otimes_{U(\mathfrak{g})} \hat{U}(\mathfrak{g})$.

 \begin{lem}\label{lem:canonical-strict} Any map $M \rightarrow N$ in category $\ocal$ is strict for the canonical topology.
 \end{lem}%
 \begin{proof}
See for example~\cite[Prop.\ 3.1.1]{MR3072116}.
 \end{proof}
We can consider the twisted Verma module $M^w_\lambda$, which admits the basis 
 $$ \oplus_{ k_\alpha \geq 0, \alpha \in w^{-1} \Phi^- \cap \Phi^+, k_\alpha <
   0, \alpha \in w^{-1} \Phi^{-} \cap \Phi^-} E\prod T_{\alpha}^{k_\alpha}.$$ We
 use this basis to define the natural topology $\mathcal{T}_{nat}$ as above.  It follows that twisted Verma modules have two topologies $\mathcal{T}_{nat}$ and $\mathcal{T}_{can}$. 
 It is immediate  that in case $w = Id$, the canonical topology and the
 natural topology coincide. We will next show that  they coincide more
 generally if~$\lambda$ is $p$-adically non-Liouville (see
 Definitions~\ref{defn-number-non-Liouville} and~\ref{defn-weight-non-Liouville}).
\begin{lem}\label{lem-george-proved}

Suppose $x\in E$ is $p$-adically non-Liouville and $x\not\in\mathbb{Z}_{> 0}$.  Then there exists a constant $C>0$ such that
$\prod_{r=1}^n|x-r|\geq C^n$ for all $n$.%
\end{lem}%
\begin{proof}%
If $x\not\in\Z_p$, there is a constant $C>0$ such that $|x-r|\geq C$ for all $r$ (e.g. take $C=\min_{r\in\Z_p}|x-r|$) so the result is clear in this case.  So assume $x\in\Z_p$.  By assumption, there is a constant $B>0$ such that $|x-r|\geq {B}^r$ for all $r\in\Z_{>0}$ (see~\eqref{formulaL}).
For a given $n$, let $p^{m-1}\leq n< p^{m}$ and choose $0< r_0\leq p^m$ with $r_0\equiv x\pmod {p^m}$.  Then, for $0< r\leq n$ with $r\not=r_0$, $|x-r|=|r_0-r|$, and so we may estimate
\[
\begin{aligned} \prod_{r=1}^n|x-r| = & \ |x - r_0| \cdot  \prod_{r=1, r \ne r_0}^{n} |r_0 - r| \\
\geq & \ B^n  \cdot |(p^m)! | \ge  C^n \end{aligned}\] for some~$C>0$.\end{proof}

\begin{lem}\label{lem-strictnesssl2}  If $\lambda$ is a non-Liouville number, the maps of  Corollary \ref{coro-intertwining} are strict for the natural topology.
\end{lem}
\begin{proof} When $\lambda \in \Z_{\geq 0}$, the statement is obvious. We  now
  assume that $\lambda \notin \Z_{\geq 0}$, so we need to show that the isomorphism $I: M_\lambda \rightarrow M_\lambda^\vee$ which sends $\bar{X}^n$ to $ n! \lambda(\lambda-1)\cdots (\lambda-(n-1)) (\bar{X}^n)^\star$ is strict. 
By Lemma \ref{lem-george-proved} (and the trivial bound $|n!|\ge p^{-n}$), we see that for~ $n>0$, we have $\vert  n!
\lambda(\lambda-1)\cdots (\lambda-(n-1)) \vert \geq C^n$ %
for a positive constant $C$. This easily implies strictness. 
\end{proof}
 
 \begin{lem}\label{lem:beta-strict} If $\lambda$ is a non-Liouville weight, the sequence of Proposition \ref{prop-intertwining}~(1)
  $$0 \rightarrow M^w_{ s_{ w^{-1} \beta}\cdot\lambda} \rightarrow M_\lambda^w \rightarrow M_{\lambda}^{s_\beta w} \rightarrow M^{s_\beta w}_{ s_{ w^{-1} \beta}\cdot\lambda} \rightarrow 0$$
  or the isomorphism $ M_\lambda^w \rightarrow M_{\lambda}^{s_\beta
    w}$ of Proposition \ref{prop-intertwining}~(2)   are strict for the natural  topology. 
 \end{lem}
 \begin{proof} This follows from Lemma \ref{lem-strictnesssl2}.
 \end{proof}

 \begin{lem}\label{lem-topo-ses} Let $0 \rightarrow V_1 \rightarrow V_2 \rightarrow V_3 \rightarrow 0$ be a short exact sequence of $E$-vector spaces. Assume that $V_2$ has two locally convex topologies $\mathcal{T}$ and $\mathcal{T}'$. Assume that the induced topologies $\mathcal{T}_1$ and $\mathcal{T}'_1$ on $V_1$, as well as the induced  topologies $\mathcal{T}_2$ and $\mathcal{T}'_3$ on $V_3$ coincide. Then $\mathcal{T}$ and $\mathcal{T}'$ coincide.
 \end{lem}
 \begin{proof} The topologies $\mathcal{T}$ and $\mathcal{T}'$ are given by  families of lattices $\{L_i\}_{i \in I}$ and $\{L_{i'}\}_{i' \in I'}$ subject to certain conditions (in particular, for any $i,j \in I$, there is a $k \in I$ such that $L_k \subseteq L_i \cap L_j$). By symmetry,  it suffices to prove that for any $i' \in I'$, there is $i \in I$ such that $L_i \subseteq L_{i'}$.  Any lattice $L$ in $V_2$ sits in an exact sequence $0 \rightarrow L_{1} \rightarrow L \rightarrow L_{3}  \rightarrow 0$, with $L_1 = L \cap V_1$. 
 By assumption, there is $i_1\in I$ such that $L_{i_1,1} \subseteq L_{i', 1}$ and $i_3 \in I$ such that $L_{i_3, 3} \subseteq L_{i',3}$. Picking $i \in I$ such that $L_i \subseteq L_{i_1} \cap L_{i_3}$, we find that $L_i \subseteq L_{i'}$. 
 \end{proof}

\begin{prop}\label{prop-2topequal} If $\lambda$ is a non-Liouville weight, the canonical and natural topologies on $M^w_\lambda$ coincide.
\end{prop}

\begin{proof} We use induction on the length of $w$. We know this is
  true if $w=Id$. We assume that this is true for $w$  and all
  $\lambda'$, and want to prove it for $M_{\lambda}^{s_\beta w}$,
  where $\ell(s_{\beta} w) = \ell(w) +1$. 
If the intertwining map $M_\lambda^w \rightarrow M_\lambda^{s_\beta
  w}$ is an isomorphism, it is strict for both the canonical and natural
topologies (by Lemmas~\ref{lem:canonical-strict} and~\ref{lem:beta-strict}), so we are done.

Otherwise, we have a long exact sequence: 
\[0 \rightarrow M^w_{ s_{ w^{-1} \beta}\cdot\lambda} \rightarrow M_\lambda^w \rightarrow M_{\lambda}^{s_\beta w} \rightarrow M^{s_\beta w}_{ s_{ w^{-1} \beta}\cdot\lambda} \rightarrow 0.\]
This long exact sequence is again strict for the natural and canonical
topologies. Since the natural and canonical topologies agree on $M^w_{
  s_{ w^{-1} \beta}\cdot\lambda}$ and $M_\lambda^w $, and also on $M^{s_\beta w}_{
  s_{ w^{-1} \beta}\cdot\lambda} \simeq M^{w}_{
  s_{ w^{-1} \beta}\cdot\lambda}$, %
we are done by Lemma \ref{lem-topo-ses}. 
\end{proof}

\subsection{Proof of Theorem \ref{thm-STrict}}

 \begin{lem}\label{lem-on-coho-affine-line} Let $\alpha$ be a root. Consider the complexes $C_1:  E[T_\alpha] \otimes \mathfrak{u}_\alpha \rightarrow E[T_\alpha]$ and $C_2:  E[T_\alpha, T_\alpha^{-1}]/E[T_\alpha] \otimes \mathfrak{u}_\alpha \rightarrow E[T_\alpha,T_\alpha^{-1}]/E[T_\alpha]$. Let $\hat{C_1}_{nat}$ and $\hat{C_2}_{nat}$ be the completions of these complexes for the natural topology. For $i=1, 2$, the natural map  $C_i \rightarrow \hat{C_i}_{nat}$ is a quasi-isomorphism.  In particular, the differentials in $\hat{C_i}_{nat}$ are strict. 
 \end{lem}
 \begin{proof} This is a standard computation; for instance,
   $\hat{C_1}_{nat}$ computes the de Rham cohomology of the analytic
   affine line, but in any case it is a simple explicit calculation as in Lemma~\ref{lemkey-comput}. %
 \end{proof}

We now restate Theorem~\ref{thm-STrict}, for the reader's convenience.
\begin{thm}[Theorem~\ref{thm-STrict}]\label{thm-STrict2}  Let $M \in \ocal_{nL}(\mathfrak{g}, \mathfrak{b})$. The canonical map $$\hat{U}(\mathfrak{m}_w)\otimes_{U(\mathfrak{m}_w)}(E \otimes^L_{U(\mathfrak{u}_{\mathfrak{p}_w})} M)   \rightarrow E \otimes^L_{\hat{U}(\mathfrak{u}_{\mathfrak{p}_w})} ( \hat{U}(\mathfrak{g})\otimes_{U(\mathfrak{g})} M)$$
is a quasi-isomorphism.
\end{thm}
\begin{proof}
  A standard argument using the five lemma shows that we may replace~$M$ by a
  resolution, and thus we reduce to the case that~$M=M_{\lambda}$ for some
  non-Liouville~$\lambda$. In fact, it is more convenient to handle all of the twisted
  Verma modules
  $M=M^{w'}_{\lambda}$ by induction on $\ell(w'w^{-1}w_0)$.

  We begin with the base case $w'=w_0 w$. Lemma \ref{lem-on-coho-affine-line} implies easily that   the Chevalley--Eilenberg complex computing $E \otimes^L_{U(\mathfrak{u}_{\mathfrak{p}_w})} M^{w_0w}_\lambda$ is a strict complex for the natural topology, so that the formula of Proposition \ref{prop-computing-coho-twistedverma} $$E \otimes^L_{U(\mathfrak{u}_{\mathfrak{p}_w})} M_\lambda^{w_0w} = M(\mathfrak{m}_w)^{(w^M)^{-1}w_{0,M}w}_{(w^M)^{-1}(w^M\cdot\lambda + 2\rho^M)}[d-\ell(w^M)]$$
passes to completions for the natural topology.  The result now follows since the canonical
and natural topologies coincide (on $M^{w_0w}_\lambda$, and  on
$M(\mathfrak{m}_w)^{(w^M)^{-1}w_{0,M}w}_{(w^M)^{-1}(w^M\cdot\lambda +
  2\rho^M)}$) by Proposition \ref{prop-2topequal}.

For the inductive step, we can suppose that there is some $\beta$ such that
$\ell(s_\beta w'w^{-1}w_0)=\ell(w'w^{-1}w_0)-1$. Then 
$\ell(s_\beta w')=\ell(w')\pm 1$. Let us assume
that it is $\ell(w')+1$, so that we have the intertwining map  $M^{w'}_\lambda \rightarrow M^{s_\beta
  w'}_\lambda$ (the other case is almost identical, using  the intertwining map
$ M^{s_\beta
  w'}_\lambda\rightarrow M^{w'}_\lambda $, and we leave it to the reader). 

If the intertwining map is an isomorphism, we are done. Otherwise, we have a long exact sequence: 
$0 \rightarrow M^{w'}_{ s_{ (w')^{-1} \beta}\cdot\lambda} \rightarrow M_\lambda^{w'} \rightarrow M_{\lambda}^{s_\beta w'} \rightarrow M^{s_\beta w'}_{ s_{ (w')^{-1} \beta}\cdot\lambda} \rightarrow 0$.
By induction, the theorem holds for $M_{\lambda}^{s_\beta w'} $ and $ M^{s_\beta w'}_{ s_{ (w')^{-1} \beta}\cdot\lambda}$. Since  
$M^{w'}_{ s_{ (w')^{-1} \beta}\cdot\lambda}  \rightarrow M^{s_\beta w'}_{ s_{
    (w')^{-1} \beta}\cdot\lambda} $ is an isomorphism, it also holds for
$M^{w'}_{ s_{ (w')^{-1} \beta}\cdot\lambda}$, and %
thus (again by the five lemma) for $M_\lambda^{w'}$, as required.
\end{proof}

We will use the following result in Section~\ref{subsec:localization-Higher-Coleman-GSp4}.
\begin{prop}\label{prop:intertwining-for-homology-w2} Assume $\beta$ is a simple root with $w=w_0^Ms_\beta\in\WM$ and $\langle\lambda,\beta^\vee\rangle\not\in\Z_{\geq 0}$.  Then
$$E\otimes^L_{U(\mathfrak{u}_{\mathfrak{p}_{w}})}M_\lambda^{Id}=M(\m_{w})_{w^{-1}(w\cdot\lambda+2\rho^M)}^{Id}[1]$$
\end{prop}
\begin{proof}
By the hypothesis that $\langle\lambda,\beta^\vee\rangle\not\in\Z_{\geq 0}$, there is an isomorphism
$M_\lambda^{Id}\isoto M_\lambda^{s_\beta}$ by Proposition
\ref{prop-intertwining}~(2).  The result then follows from Proposition
\ref{prop-computing-coho-twistedverma} (taking $w$ there to be  $w_{0,M}w$).
\end{proof}

\begin{rem}
We note that the condition that $w_0^Ms_\beta\in\WM$ is equivalent to $\beta$ not being a root of
$\m_{w_0^M}$.%
\end{rem}

\section{Equivariant sheaves on the flag variety and
  localization}\label{sec: classicality}
\subsection{Introduction}%
  This entire section is concerned with geometric representation theory. We fix a split reductive group $G$, a Parabolic $P$ with Levi $M$,  and a Borel $B \subseteq P$.  We consider the partial flag variety $P \backslash G$ and its Bruhat decomposition $P \backslash G = \coprod_{w \in \WM} P \backslash PwB$ into $B$-orbits, indexed by the subset $\WM$ of Kostant representatives of the Weyl group $W$ of $G$.  
  
  In Section ~\ref{subsec:equivariant-sheaves} and Section~\ref{subsec:equivariant-sheaves-Bruhat-strata} we consider
equivariant sheaves on the partial flag variety as well as (dagger neighbourhoods of) Bruhat cells, for the action of $G$, its Lie algebra $\mathfrak{g}$, or a Borel subgroup $B$, depending on the context. 
We also establish the connection between these equivariant sheaves  and twisted $D$-modules and introduce the
horizontal Levi action.
We begin with some
generalities on equivariant sheaves on adic and dagger spaces, before
turning to the specific cases that we need. We repeatedly make use of
the standard equivalence (given by passage to the fibre at a
point~$x\in X$) between $H$-equivariant sheaves on a
space~$X$ on which the group~$H$ acts transitively, and the
representations of the stabilizer group~$\Stab_{H}(x)$; however,
since we are working with topological (or rather solid) sheaves, we
have to go to some lengths to make precise the categories that we are
working with, and their interactions with these equivalences. (The
particular categories that we work with are ultimately dictated by the
use of geometric Sen theory in
Section~\ref{sec:application-to-Shimura}.)

\begin{rem} All the sheaves we consider will be sheaves on topological spaces, valued in the category of solid $E$-vector spaces (where $E$ is a finite extension of $\qq_p$). These form an abelian category. Our topological spaces will usually be adic spaces or dagger spaces, and   our sheaves will  also be  ``quasi-coherent'' and  often be twisted $D$-modules. This means that the objects we manipulate  would naturally fit in the   formalism of quasi-coherent sheaves on adic spaces of \cite{andreychev2021pseudocoherent}, and the formalism of analytic geometry and the de Rham stack of \cite{camargo2024analyticrhamstackrigid}. The much simpler perspective we adopt is sufficient for our purposes. 
\end{rem}

This preliminary material is used in 
Section~\ref{subsec:algebraic-locally-analytic-representations} to  
produce, for any $w \in \WM$,  a functor~$HCS$ (for ``Higher Coleman sheaf'') from category~$\cO$ for~ $\mf{m}_w$ (the Lie algebra of $w^{-1} Mw$), to the category
of $(\mf{g},B)$-equivariant sheaves on the dagger neighborhood  of the Bruhat cell  $P\backslash PwB$. %
In section  \ref{sec:application-to-Shimura} we will use these sheaves to produce sheaves on (open subsets of) Shimura varieties whose cohomology with support  is  Higher Coleman theory
of~\cite{boxer2021higher}.

In Section~\ref{subsec:localization-partial-flag-variety} we define
our localization functor on the partial flag variety.  This functor goes from category $\ocal$ for $\mathfrak{g}$ to twisted $D$-modules on the flag variety. 
In 
Theorem~\ref{thm-localization}  we describe the localization in terms of Higher Coleman sheaves. 
Namely, in 
$p$-adically non-Liouville weight, the restrictions to
the Bruhat cells of the cohomology sheaves of the localization of a Verma
module~$M$ of~$\mf{g}$ are given by the Higher Coleman sheaves associated
to the $\mathfrak{u}_{P_w}$-homology of~$M$.  (It is here that we use
Theorem~\ref{thm-STrict}.) Furthermore we give an explicit filtration
on these sheaves in Corollary~\ref{coro-ESflag}.%

Finally in Section~\ref{subsec:localization-Higher-Coleman-GSp4} we
specialize to the case $G=\GSp_4$ and prove the crucial
Theorem~\ref{thm-extequal}, which describes the cohomology of the horizontal Cartan action on the localization in a special case of  interest to us.

\subsection{Equivariant sheaves on partial flag varieties}\label{subsec:equivariant-sheaves} In this section we discuss several kind of equivariant sheaves. 
\subsubsection{Equivariant sheaves over adic
  spaces}\label{subsec-equivariant-sheaves}
  Let $C$ be a rank one field extension of~ $E$. In applications, $C$ is either $E$ or $\C_p$. 
Let $X$ be an adic space  which is locally of finite type over  $\Spa(C, \ocal_C)$. Its structure sheaf $\oscr_X$ is naturally a
topological sheaf, whose value on a  quasi-compact open subset is a Banach space. It follows that we can think of $\oscr_X$ as taking values in the category $\Mod(E)$. 
All the sheaves we will encounter will be sheaves of solid $E$-vector spaces.  By a solid $\oscr_X$-module we mean a sheaf valued in the category $\Mod(E)$ equipped with an $\oscr_X$-module structure.  We emphasize that we do not impose any kind of quasi-coherence condition in the definition of solid $\oscr_X$-modules.

\begin{defn} \leavevmode \begin{enumerate}
\item A sheaf  $\mathscr{F}$ of solid 
$\oscr_{X}$-modules is an \emph{orthonormalizable Banach sheaf} if there exists a Banach space $V$ over $E$ such that $\mathscr{F} = \oscr_X \otimes_E V$. 
\item  A sheaf  $\mathscr{F}$ of solid 
$\oscr_{X}$-modules is a \emph{summand of orthonormalizable Banach sheaf} if it
is a direct summand %
of  an orthonormalizable Banach sheaf. 
\item A sheaf  $\mathscr{F}$ of solid $\oscr_{X}$-modules is a     \emph{Banach sheaf} if there is a  covering
$X = \cup_i \Spa(A_i,A_i^+) $ %
and a Banach space $V_i$ over $E$ such that $\mathscr{F}\vert_{\Spa(A_i,
  A_i^+)}$ is a  direct summand of the sheaf $\oscr_{\Spa(A_i,A_i^+)} {\otimes}_E V_i$.
  \item A sheaf  $\mathscr{F}$ of solid $\oscr_{X}$-modules is an $\LB$-sheaf if there is a  covering
$X = \cup_i \Spa(A_i,A_i^+) $ %
and $LB$-spaces $V_i$ over $E$ such that $\mathscr{F}\vert_{\Spa(A_i, A_i^+)}$
is a  direct summand  of the sheaf $\oscr_{\Spa(A_i,A_i^+)} {\otimes}_E V_i$.
  \end{enumerate}Banach sheaves define a category  $B(X)$ and $\LB$-sheaves define a category $\LB(X)$. 
  \end{defn}  
  
Let $G$ be an analytic group acting on $X$.  We have two maps $\act, p: G \times
X \rightarrow X$, which are respectively the action and projection maps.  We let $B_G(X)$  be the category  of $G$-equivariant Banach sheaves, whose objects are objects $\mathscr{F}$ of $B(X)$ together with  an isomorphism $\act^\star \mathscr{F} \rightarrow p^\star \mathscr{F}$ (in the category $B(G \times X)$) satisfying the usual cocycle condition.  We let $\{G_n\}_{n \in \ZZ_{\geq 0}}$ be a system of neighborhoods of the identity $e$ in $G$, given by quasi-compact open subgroups.

  \begin{defn}
    The category $\LB_G(X)$ is the category whose objects are objects
    $\mathscr{F}$ of $\LB(X)$ together with an isomorphism
    $\act^\star \mathscr{F} \rightarrow p^\star \mathscr{F}$ (in the
    category $\LB(G \times X)$) satisfying the usual cocycle
    condition, and further satisfying the following finiteness
    condition:

    \begin{enumerate}
    \item There exists a covering $X = \cup_j U_j$ such that
      $\mathscr{F}\vert_{U_j} = \colim_{r \geq 0} \mathscr{F}_{j,r}$
      is a filtered countable inductive limit of orthonormalizable
      Banach sheaves with injective transition maps. %
    \item For all $j$, there exists a quasi-compact open subgroup
      $G_{r(j)} \subseteq G$ which stabilizes $U_j$ and we can upgrade
      $\mathscr{F}_{j,r} $ to an object of $B_{G_{r(j)}}(U_j)$, in
      such a way that the inductive system $ \{\mathscr{F}_{j,r}\}$ is
      an inductive system in $B_{G_{r(j)}}(U_j)$.
    \item The two $G_{r(j)}$-actions on $\mathscr{F}\vert_{U_j}$ (the
      one induced by the inclusion $G_{r(j)} \hookrightarrow G$, and
      the one obtained by taking the colimit of the
      $\mathscr{F}_{j,r}$) are the same.
    \end{enumerate}
  \end{defn}

We let $\mathfrak{g}$ be the Lie-algebra of $G$.   The action of $G$
on $X$ induces an action of $\mathfrak{g}$ by derivations on
$\oscr_X$.  
We let $\Mod'_{\mathfrak{g}}(X)$ be the following category: its objects are solid $\oscr_X$-modules $\mathscr{F}$ together with a
  map $\mathfrak{g} \otimes_E \mathscr{F} \rightarrow \mathscr{F}$
  inducing an action of $\mathfrak{g}$ on $\mathscr{F}$ by derivations in the following sense: 
  
  \begin{enumerate}
  \item For any $X,Y \in \mathfrak{g}$, we have $[X,Y] = XY-YX$ in
    $\mathrm{End}(\mathscr{F})$.
  \item For any
    $(X,a,f) \in \mathfrak{g} \times \oscr_{X} \times
    \mathscr{F}$, we have $X(af) = X(a) f + a X(f)$. \end{enumerate}
\begin{defn}\label{defn-Modg} We let $\Mod_{\mathfrak{g}}(X)$ be the subcategory of $\Mod'_{\mathfrak{g}}(X)$ generated under colimits by objects $\mathscr{F}$ which have the following property:
 for any quasi-compact open subset $U$ of $X$, there exists $r$ such that the action $\mathfrak{g} \otimes \mathscr{F}(U) \rightarrow \mathscr{F}(U)$ can be integrated to an action $D(G_r) \otimes \mathscr{F}(U) \rightarrow \mathscr{F}(U)$. 
 \end{defn}
 \begin{lem} The categories $\Mod_{\mathfrak{g}}(X)$ and $\Mod'_{\mathfrak{g}}(X)$ are Grothendieck abelian
   category, and in particular have enough injectives. %
 \end{lem}
 \begin{proof}We begin with the case of $\Mod'_{\mathfrak{g}}(X)$, where the
   only non-obvious point is the existence of a set of generators. For this we
   may take the sheaves   $ j_! (U(\mf{g}) \otimes R \otimes \oscr_U)$ for $R$ a
   generator of the category of solid $E$-modules, $j : U \hookrightarrow X$ a
   quasi-compact open subset and $r \in \qq_{\geq 0}$.  We now turn to
   $\Mod_{\mathfrak{g}}(X)$, where we first make a comment on the condition that the action $\mathfrak{g} \otimes \mathscr{F}(U) \rightarrow \mathscr{F}(U)$ can be integrated to an action $D(G_r) \otimes \mathscr{F}(U) \rightarrow \mathscr{F}(U)$ for some $r$. Let $D(G_{r^+}) = \lim_{r'>r} D(G_{r'})$. Then $D(G_{r^+}) \otimes_{U(\mathfrak{g})} D(G_{r^+}) =  D(G_{r^+})$ by \cite[Lem. 5.13]{MR4475468}. As a result, the extension of the $\mathfrak{g}$-action to an action of $D(G_{r^+})$ for some $r$ is a property of $\mathscr{F}(U)$ and not some extra data: it means that $\mathscr{F}(U) \otimes_{U(\mathfrak{g})} D(G_{r^+}) = \mathscr{F}(U)$ for some $r$. 
 By construction $\Mod_{\mathfrak{g}}(X)$ is an abelian subcategory of
 $\Mod'_{\mathfrak{g}}(X)$ stable under colimits, and filtered colimits are
 exact since they are exact in $\Mod'_{\mathfrak{g}}(X)$. Then a set of
 generators is given by the sheaves $ j_! (D(G_r) \otimes R \otimes \oscr_U)$
 for~$R,U$ as above. 
 \end{proof}

We define the subcategory $\LB_{\mathfrak{g}}(X)$ of $\Mod_{\mathfrak{g}}(X)$ as follows.
\begin{defn}\label{defn:LBgX}
  The objects of $\LB_{\mathfrak{g}}(X)$ are $\LB$-sheaves $\mathscr{F}$ on $X$ together with a
  map $\mathfrak{g} \otimes_E \mathscr{F} \rightarrow \mathscr{F}$
  inducing an action of $\mathfrak{g}$ on $\mathscr{F}$ by derivations.
  We furthermore impose that the $\mathfrak{g}$-action can locally be
  integrated to a locally analytic action. Here is the precise
  condition:

  \begin{enumerate}
  \item There exists a covering $X = \cup_j U_j$ such that
    $\mathscr{F}\vert_{U_j} = \colim_{r \geq 0} \mathscr{F}_{j,r}$ is
    an inductive limit of orthonormalizable Banach sheaves with
    injective transition maps.%
  \item For $r$ large enough, we can upgrade $\mathscr{F}_{j,r} $ to
    an object of $B_{ G_{r}}(U_j)$, in such a way that the transition
    maps $\mathscr{F}_{j,r} \rightarrow \mathscr{F}_{j,r'} $ are
    equivariant for the maps $G_{r'} \rightarrow G_r$.
 
  \item The two $\mathfrak{g}$-actions on $\mathscr{F}\vert_{U_j}$
    (the one induced by differentiating the action of $ G_{r}$ and
    passing to the colimit, and the one which is part of the original
    data) are the same.

  \end{enumerate}
\end{defn}

\subsubsection{Equivariant sheaves over topological and ringed  spaces}
We also need to consider the situation where a locally profinite group $M$ acts continuously on a locally spectral topological space $X$.
\begin{lem} Let $V$ be a quasi-compact open subset of $X$. There is a compact open subgroup $N$ of $M$ such that $N\cdot V = V$.
\end{lem}
\begin{proof} The map $\act: M \times X \rightarrow X$ is
  continuous. It follows that $act^{-1}(V)$ is open, so that for any $x \in V$,
  there exists a compact open subgroup $N_x$ of $M$ and an open neighborhood
  $V_x \subseteq V$ of~$x$ such that $N_x\cdot V_x \subseteq V$. Since $V =
  \cup_{x \in V} V_x$, and~$V$ is quasi-compact, 
   there is a finite collection of elements $\{x_i\}_{i \in I}$ such that $V = \cup_i V_{x_i}$. We deduce that $N = \cap_i N_{x_i}$ works. 
\end{proof}

\begin{defn}
  We let $\Mod_M(X)$ be the category consisting of:
  \begin{enumerate}
  \item A $\Mod(E)$-valued sheaf $\mathscr{F}$ on $X$.
  \item An abstract action of $M$ on $\mathscr{F}$ (that is for any $m \in M$, an isomorphism
    $m^{-1} \mathscr{F} \rightarrow \mathscr{F}$ satisfying
    compatibility conditions for various $m$).
  \item For any quasi-compact open subset $V \subseteq X$,  for one
    (equivalently for any) compact open subgroup $N_V$ of $M$
    stabilizing $V$, the abstract action of $N_V$ on $\mathscr{F}(V)$
    extends to a $E_{\blacksquare}[N_V]$-module structure.
  \item For any $U \subseteq V$, and for one (equivalently any)
    quasi-compact open subgroup $N_{V,U}$ stabilizing $U$ and $V$, the
    restriction maps $\mathscr{F}(V) \rightarrow \mathscr{F}(U)$ are
    $E_{\blacksquare}[N_{V,U}]$-equivariant.
  \end{enumerate}
  If in $(3)$, the action of $N_V$ is smooth,  we say that $\mathscr{F}$ is smooth. We thus have a subcategory category $\Mod_M^{\sm}(X)$ of $\Mod_M(X)$.
 \end{defn}

  \begin{lem}\label{coro-enough-inj} The categories  $\Mod_M^{\sm}(X)$ and $\Mod_M(X)$ are Grothendieck
    abelian categories, and in particular they have enough injectives. %
  \end{lem}
  \begin{proof} All claims are obvious, except for the existence of
    generators. Let us first prove the existence of a generator in
    $\Mod_M(X)$. Let $U$ be a quasi-compact open  in $X$ and  let $N$ be a compact
    open subgroup stabilizing   $U$. Let $R$ be a generator of the
    category $\Mod_N(E)$. %
      We consider the sheaf $L(U) = \oplus_{m \in M/N} j_{mU, !} R$ where $j_{mU}: mU \rightarrow X$ is the open immersion. It is endowed with the obvious $M$-equivariant action. Let $\mathscr{F}$ be an object of $\Mod_M(X)$. A map $L(U) \rightarrow \mathscr{F}$ amounts to a map $R \rightarrow \mathscr{F}(U)$ in the category $\Mod_N(E)$. It follows that $\oplus_U L(U)$ is a generator of $\Mod_M(X)$.  We construct similarly a generator of $\Mod_M^{\sm}(X)$ by the same construction, but replacing $R$ by a generator of $\Mod_N^{\sm}(E)$. 
  \end{proof}

 Let us briefly indicate some possible variations. If $X$ is equipped with a sheaf of algebras in solid $E$-vector spaces $\oscr_X$, which belongs to $\Mod_{M}(X)$, one can consider the category $\Mod_M(\oscr_X)$ of $M$-equivariant $\oscr_X$-modules. If the $M$-equivariant sheaf $\oscr_X$ is smooth, we also have a category $\Mod_{M}^{\sm}(\oscr_X)$.  In this case, we can also introduce a twist by a character $\lambda: M \rightarrow E^\times$. We say that an object $\mathscr{F}$ of $\Mod_{M}(\oscr_X)$ is $\lambda$-smooth if $\mathscr{F}\otimes E(-\lambda)$ is smooth. The category of $\lambda$-smooth objects is denoted by $\Mod_M^{\lambda-\sm}(X)$. The categories
 $\Mod_M(\oscr_X)$, $\Mod^{\lambda-\sm}_M(\oscr_X)$ are again Grothendieck
 abelian categories, and %
  in particular they have enough injectives. %

\subsubsection{Dagger spaces}\label{subsubsec: dagger spaces}

 We need to enlarge the category of adic spaces and also consider certain limits of adic spaces (for example dagger spaces in the sense of \cite{MR1739729}). Let $Z$ be a locally closed subset of $X$. We let $Z^{\dag,X}$ be the locally ringed space $\lim_{Z \subseteq U} U$ where $U$ runs through the open subsets of $X$ containing $Z$. As a topological space, $Z^{\dag,X}=Z$. It carries the structure sheaf $\oscr_{Z^{\dag,X}} = i_Z^{-1} \oscr_X$ where $i_Z: Z \rightarrow X$ is the inclusion.   This sheaf also takes values in the category $\Mod(E)$. 

\begin{rem} The locally ringed space $Z^{\dag,X}$ depends on $Z$ and on the embedding $Z \hookrightarrow X$. If $X$ is clear from the context, we simply denote $Z^{\dag,X}$ by $Z^\dag$. 
\end{rem}
\subsubsection{Notations for the flag variety}

We let $G$ be a connected split reductive group  over $E$  and we let $P$ be a
parabolic subgroup. %
We let $\mathcal{FL} = P \backslash G$ be the
partial flag variety (an analytic adic space). 

\begin{rem} In section \ref{sec:application-to-Shimura}, we will use the notation $\mathcal{FL}^{rat}$ for this analytic space, and we will let $\mathcal{FL}$ be its base change to $\Spa~(\C_p, \ocal_{\C_p})$. 
\end{rem}

We let $U_P$ be the unipotent radical of $P$,
with Levi quotient $M = P/U_P$. We let $B \subseteq P$ be a
Borel with maximal torus $T$ and unipotent radical $U_B$. %
We let $B_M$
be the induced Borel on $M$  and $U_M$ be the unipotent radical of
$B_M$. %
We use gothic letters
for the Lie algebras of all groups introduced so far, so that for
example~$\mathfrak{g}$ is the Lie algebra of~$G$ and~$\mathfrak{b}$ is
the Lie algebra of~$B$; the one exception is that following standard conventions, the Lie algebra of $T$ is denoted by $\mathfrak{h}$. 
For any $x \in \mathcal{FL}$, we let $P_x$ be $x^{-1} P x$, and let $U_{P_x}$
be its unipotent radical.  We adopt similar notation for other groups
or Lie algebras, and in addition for~$x$ replaced by an element~$w$ of a
Weyl group; so for example for~$w\in W$ (the Weyl group of~$G$) we
have ~$P_w$, $U_{P_w}$ and so on.

 From now on $G$ is viewed as an analytic group over $\Spa(E, \ocal_E)$ (this is a Stein space and is not quasi-compact unless $G=\{1\}$). 

We fix a reductive model for $G$ over $\ocal_E$. Its analytification defines a quasi-compact open subgroup $G_0 \subseteq G$. 
For any $r \in
\qq_{\geq 0}$, %
we let $G_r$ be the quasi-compact analytic subgroup
of $G_0$  of elements reducing to the identity $e$  modulo
$p^r$. %
Here are some slightly non-standard conventions and constructions:  
\begin{itemize}
\item If $H$ is an analytic
subgroup of $G$, %
we let $H_r = G_r \cap H$.

  \item If $H$ is an analytic subgroup of $G$, we let $H_e:=
\lim_{r\geq 0} H_r$ where the limit is taken in the category of  locally ringed spaces. Thus $H_e = \{e\}^{\dag, H}$.   As a space, $H_e$ has only one
point (the identity  $e$ of $G$), but it carries the  structure sheaf
$\oscr_{H,e}$ whose dual is the distribution algebra $\hat{U}(\mathrm{Lie}(H))$.

\end{itemize}

\subsubsection{ $G$-equivariant sheaves}
In section \ref{subsec-equivariant-sheaves}, we  have introduced  the
categories $B_G(\mathcal{FL})$ and $\LB_G(\mathcal{FL})$. We also have
the categories of representations $B_{P_w}(E)=B_{P_w}(\Spa(E,\cO_E))$ and
$\LB_{P_w}(E)=\LB_{P_w}(\Spa(E,\cO_E))$ for each~$w\in W$.

\begin{prop}\label{prop: equiv of cats BG Fl}Taking the fiber
  at~$w=P\backslash Pw\in\mathcal{FL}$ gives equivalences of categories between  $B_G(\mathcal{FL})$ and $B_{P_w}(E)$ and between $\LB_G(\mathcal{FL})$ and $\LB_{P_w}(E)$. 
\end{prop}
\begin{proof}We have an isomorphism $P\backslash G\isoto P_w\backslash
  G$, $x\mapsto w^{-1}\cdot x$, which takes $w\in \mathcal{FL}$
  to~$e\in P_w\backslash
  G$; so we can and do reduce to the case~$w=e$.
  Given a $G$-equivariant sheaf $\mathcal{V}$, we take its
  fiber at $e = P \backslash P \in \mathcal{FL}$, which is a
  representation of $P$. Conversely,   let $\pi: G \rightarrow
  \mathcal{FL}$, $g \mapsto eg$ be the uniformization map. The sheaf
  $\pi_\star \oscr_G$ is $G$-equivariant (via the action of $G$ by
  right translation on itself) and carries a $P$-action (via the
  action of $P$ by left translation on $G$). Given an object  $V$ of
  $\LB_P(E)$, we consider $(\pi_\star \oscr_{G} {\otimes}
  V)^{P}$. These two functors define the equivalences of
  categories of the proposition (and in particular match the various
  finiteness conditions); we leave the details to the reader. %
\end{proof}

\begin{example} We have a filtration $\mathfrak{u}_P \subseteq \mathfrak{p} \subseteq \mathfrak{g}$  of finite dimensional $P$-representations. Via our equivalence of categories, this corresponds to a filtration of $G$-equivariant coherent sheaves:
$\mathfrak{u}_P^0 \subseteq \mathfrak{p}^0 \subseteq \mathfrak{g}^0 = \oscr_{\mathcal{FL}}  \otimes \mathfrak{g}$. 
The fibers of this filtration at a point $x \in \mathcal{FL} $ %
are $\mathfrak{u}_{P_x} = x^{-1} \mathfrak{u}_P x \subseteq
\mathfrak{p}_x = x^{-1} \mathfrak{p} x \subseteq %
E(x)\otimes \mathfrak{g}$. %
Moreover, we have an isomorphism $\mathfrak{g}^0/\mathfrak{p}^0 =
T_{\mathcal{FL}}$. %
\end{example}

\begin{example}\label{example-construction-classical-sheaves} Let
  $\lambda \in X^\star(T)^{M,+}$. There is an associated highest
  weight representation of $M$ and via our equivalence of categories, this corresponds to a $G$-equivariant coherent sheaf $\mathcal{L}_\lambda$. Here is an equivalent geometric construction of this sheaf. 
Let $ \pi: U_P \backslash G \rightarrow \mathcal{FL}$. This is a $G$-equivariant $M$-torsor. 
Then $ \mathcal{L}_\lambda=(\pi_\star \oscr_{U_P \backslash G})[B_M =
-w_{0,M} \lambda]$, with the right translation action of $G$. 
\end{example}

\begin{rem} In the Siegel case, the tautological exact sequence over $\mathcal{FL}$ is (for $St = E^{2g}$ the standard $2g$-dimensional representation of $G$):
\[0 \rightarrow \mathcal{L}_{(0, \cdots,0, -1; 1)} \rightarrow \oscr_{\mathcal{FL}} \otimes St  \rightarrow \mathcal{L}_{(1,0, \cdots, 0; 1)} \rightarrow 0.\]
\end{rem}

\subsubsection{ $\mathfrak{g}$-equivariant sheaves and the horizontal action} Let $U$ be an open subset of $\mathcal{FL}$. By Definition \ref{defn:LBgX}, we have a category of $\mathfrak{g}$-equivariant sheaves $\LB_{\mathfrak{g}}(U)$. 

On any object of $\LB_{\mathfrak{g}}(U)$, the $\mathfrak{g}$-action extends linearly to an $\oscr_{U} \otimes \mathfrak{g}$-action.  We recall that we have the moving parabolic Lie-algebra $\mathfrak{p}^0 \subseteq \oscr_{U}\otimes\mathfrak{g}$. 

\begin{lem} We have that  $\mathfrak{p}^0 \subseteq
  \oscr_{U}\otimes\mathfrak{g}$ acts
  $\oscr_{U}$-linearly %
  and $G$-equivariantly on  any object of  $\LB_{\mathfrak{g}}(U)$. 
\end{lem}%
\begin{proof} We have that  $\mathfrak{p}^0$ is a $G$-equivariant subsheaf of $\oscr_{\mathcal{FL}} \otimes \mathfrak{g}$. Moreover, it acts trivially on $\oscr_{\mathcal{FL}}$ since $T_{\mathcal{FL}} = \mathfrak{g}^0/\mathfrak{p}^0$. 
\end{proof}

                                               \begin{defn}\label{defn:LBguP0}
                                                 We let
                                                 $\LB_{\mathfrak{g}}(U)^{\mathfrak{u}_P^0}$
                                                 be the full
                                                 subcategory of
                                                 $\LB_{\mathfrak{g}}(U)$
                                                 of objects which are
                                                 annihilated by
                                                 $\mathfrak{u}_{P}^0$.
                                               \end{defn}
                                               
For any $\mathscr{F} \in \LB_{\mathfrak{g}}(U)^{\mathfrak{u}_P^0}$, we have a $G$-equivariant map 
$ \mathfrak{m}^0 = \mathfrak{p}^0/\mathfrak{u}_P^0 \rightarrow
\underline{\mathrm{End}}_{\oscr_{U}}( \mathscr{F})$ which
can be extended to an algebra map: \numequation\label{eqn: map from U
m 0} U(\mathfrak{m}^0) \rightarrow
\underline{\mathrm{End}}_{\oscr_{U}}(
\mathscr{F}).\end{equation}Let $Z(\mathfrak{m})$ be the centre of $U(\mathfrak{m})$. 

\begin{lem}\label{lem-horizontal} We have an injective algebra homomorphism $Z(\mathfrak{m}) \hookrightarrow \HH^0( \mathcal{FL}, U(\mathfrak{m}^0))$.
\end{lem}
\begin{proof} The $G$-equivariant sheaf $ U(\mathfrak{m}^0)$ is
  associated via Proposition~\ref{prop: equiv of cats BG Fl} %
  to the $P$-representation $U(\mathfrak{m})$ (the fiber at $e$). 
We have a natural inclusion  $Z(\mathfrak{m}) \into U(\mathfrak{m})$,
and $Z(\mathfrak{m})$ identifies with the  $P$-invariant subspace of
$U(\mathfrak{m})$. It follows that we get an injective map of sheaves
$\oscr_{\mathcal{FL}}\otimes Z(\mathfrak{m}) \rightarrow U(\mathfrak{m}^0)$,
inducing the expected map on global sections.%
\end{proof}

\begin{defn} We define the horizontal action as the map $\Theta_{\hor}
: Z(\mathfrak{m}) \rightarrow  \mathrm{End}_{\oscr_{U}}(
  \mathscr{F})$ obtained by composing the map of Lemma
  \ref{lem-horizontal} and the map $$ \HH^0(U,
  U(\mathfrak{m}^0)) \rightarrow
  {\mathrm{End}}_{\oscr_{U}}( \mathscr{F})$$ obtained from ~\eqref{eqn: map from U
m 0}.
\end{defn}

\subsubsection{$(\mathfrak{g},G)$-equivariant sheaves}
We now consider $(\mathfrak{g}, G)$-equivariant sheaves. We sometimes
find it helpful to interpret these as $G_e \rtimes G$-equivariant
sheaves, where the action of~$G$ on~$G_e$ is via conjugation, see
Remark~\ref{rem:fracg-versus-Ge} below. We remark there is a group homomorphism $G_e \rtimes G
  \rightarrow G$, $(g_e,g)\mapsto g_eg$.
\begin{defn}\label{defn-LBgG}
  The category $\LB_{(\mathfrak{g},G)}(\mathcal{FL})$ has objects
  consisting of a $G$-equivariant sheaf %
  $\mathscr{F} \in \LB_G(\mathcal{FL})$ together with a map
  $\mathfrak{g} \otimes \mathscr{F} \rightarrow \mathscr{F}$ of
  $G$-equivariant sheaves (where $\mathfrak{g} \otimes \mathscr{F}$
  carries the diagonal
  $G$-action), %
  giving a Lie algebra action on $\mathscr{F}$:
  \begin{enumerate}
  \item For any $X,Y \in \mathfrak{g}$, we have $[X,Y] = XY-YX$ in
    $\mathrm{End}(\mathscr{F})$.
  \item For any
    $(X,a,f) \in \mathfrak{g} \times \oscr_{\mathcal{FL}} \times
    \mathscr{F}$, we have $X(af) = X(a) f + a X(f)$\footnote{We have a
      map $\mathfrak{g} \rightarrow T_{\mathcal{FL}}$ and thus an
      action of $\mathfrak{g}$ by derivations on
      $\oscr_{\mathcal{FL}}$.}.
  \end{enumerate}
  We furthermore impose that the $\mathfrak{g}$-action can locally be
  integrated to a locally analytic action. Here is the precise
  condition:

  \begin{enumerate}
  \item There exists a covering $\mathcal{FL} = \cup_j U_j$ such that
    $\mathscr{F}\vert_{U_j} = \colim_{r \geq 0} \mathscr{F}_{j,r}$ is
    an inductive limit of Banach sheaves with injective transition maps.%
  \item For all $j$, there exists a quasi-compact open subgroup
    $G_{r(j)} \subseteq G$ which stabilizes $U_j$.
  \item For $r$ large enough, we can upgrade $\mathscr{F}_{j,r} $ to
    an object of $B_{G_r \rtimes G_{r(j)}}(U_j)$, in such a way that
    the inductive system $ \{\mathscr{F}_{j,r}\}$ is now in
    $B_{G_r \rtimes G_{r(j)}}(U_j)$.  %
  \item The two $G_{r(j)}$-actions on $\mathscr{F}\vert_{U_j}$ (the
    one induced by the inclusion $ G_{r(j)} \hookrightarrow G$, and
    the one obtained by taking the colimit of the $\mathscr{F}_{j,r}$)
    are the same.
  \item The two $\mathfrak{g}$-actions on $\mathscr{F}\vert_{U_j}$
    (the one induced by differentiating the action of $ G_{r}$ and
    passing to the colimit, and the one which is part of the original
    data) are the same.

  \end{enumerate}
\end{defn}

 \begin{rem}\label{rem:fracg-versus-Ge} In particular, the $\mathfrak{g}$-action on an
   object~$\cF$ of $\LB_{(\mathfrak{g},G)}(\mathcal{FL})$ can be upgraded to  a $\hat{U}(\mathfrak{g})$-action on $\mathscr{F}$. We can thus think of an object $\mathscr{F}$ of $\LB_{(\mathfrak{g},G)}(\mathcal{FL})$ as a $G_e \rtimes G$-equivariant sheaf satisfying certain finiteness conditions. 
 \end{rem}
 
 \begin{rem}\label{rem-defn-LBgG} In section \ref{sec:application-to-Shimura}, we will also consider the category $\LB_{(\mathfrak{g},G)}(\mathcal{FL}_{\C_p})$ for the flag variety over $\C_p$, whose definition is the obvious variant of the definition
 \ref{defn-LBgG}. There is a base change map $\LB_{(\mathfrak{g},G)}(\mathcal{FL}) \rightarrow \LB_{(\mathfrak{g},G)}(\mathcal{FL}_{\C_p})$. 
 \end{rem}
 \begin{defn}
   We define $\LB_{(\mathfrak{g}, P_w)}(E)$ to be the category whose
   objects consist of an object $V$ of $\LB_{P_w}(E)$, together with a $P_w$-equivariant morphism
   $\mathfrak{g} \otimes V \rightarrow V$ in the category
   $\LB_{P_w}(E)$, inducing a Lie algebra action of $\mathfrak{g}$ on
   $V$: %
 for any $X,Y \in \mathfrak{g}$, we have $[X,Y] = XY-YX$ in
   $\mathrm{End}(V)$.   We further impose the following finiteness condition:
   \begin{enumerate}
   \item $V = \colim_r V_r$ is an inductive limit of Banach spaces
     with injective transition maps.%
   \item There exists $s$ such that for all $r$ large enough, $V_r$
     can be upgraded to an object of $B_{G_r \rtimes P_{w,s}}(E)$ and the
     maps $V_r \rightarrow V_{r'}$ are equivariant for the map
     $G_{r'} \rtimes P_{w,s} \rightarrow G_{r} \rtimes P_{w,s}$.
   \item The action of $P_{w,s}$ on $V$ obtained on the limit is the one
     induced by restriction from $P_w$ to $P_{w,s}$.
   \item The action of $G_r$ on $V_r$ induces an action of
     $\mathfrak{g}$, and the action of $\mathfrak{g}$ on $V$ coincides
     with the original action of $\mathfrak{g}$.
   \end{enumerate}
 \end{defn}

 Similarly to Proposition~\ref{prop: equiv of cats BG Fl}, we have the
 following equivalence of categories. As usual, this equivalence is
 obtained by passage to a fiber; note in particular that we are not
 simply restricting the $\mathfrak{g}$ and~$G$-actions, and indeed the action
 of~$\mathfrak{g}$ is obtained as the difference between the given
 $\mathfrak{g}$-action and the
 derivative of the $G$-action.
\begin{prop}\label{prop: LB g G equiv of cats}Taking the fiber at~$w$
  induces an equivalence of categories between
  $\LB_{(\mathfrak{g},G)}(\mathcal{FL})$ and
  $\LB_{(\mathfrak{g},P_w)}(E)$. \end{prop}%
\begin{proof}As in the proof of Proposition~\ref{prop: equiv of cats
    BG Fl}, we can without loss of generality take~$w=e$. We consider the uniformization  map $m: G_e \rtimes G
  \rightarrow \mathcal{FL}$, $(g_e,g) \mapsto e g_eg$.  The stabilizer
  $\Stab(e)$ of $e$ is the subgroup of elements $( g_e, g) \in G_e
  \rtimes G$, such that  $g_eg \in P$. This is also the semi-direct
  product $ G_e \rtimes P$ with $G_e = \{ (g_e^{-1},g_e) \}$ and $P =
  \{(1,p)\}$. %
A $(\mathfrak{g}, G)$-equivariant sheaf $\mathcal{V}$ gives a
$(\mathfrak{g}, P)$-module by taking the fiber at $e$. Conversely,
given a $\Stab(e)$-representation $V$, we consider the sheaf
$\mathcal{V} = (m_\star (\oscr_{G_e \rtimes G}) {\otimes} V)^{G_e
  \rtimes P}$. It has an action of $G_e$ given  by $g_e f( g'_e, g') =
f( g_e' g'g_e (g')^{-1}, g' )$ and an action of $G$ given by $g. f(
g'_e, g') = f(g'_e, g'g)$.  %
We check that this induces an equivalence between
$\LB_{(\mathfrak{g},G)}(\mathcal{FL})$ and
$\LB_{(\mathfrak{g},P)}(E)$.  For example assume that $V$ is an object
of $\LB_{(\mathfrak{g},P)}(E)$, thus  $V = \colim_{r} V_r$ where each
$V_r$ carries an action of $G_r \rtimes P_s$ for $s$ and $r$ large
enough.     We can consider the map $G_r \rtimes G_s \rightarrow
\mathcal{FL}$, $(g,g') \mapsto eg'g$. The image is a neighborhood $U$
of $e$ and the stabilizer of $e$ is $G_r \rtimes P_s$.  We deduce that $\mathcal{V}\vert_U = \colim_r (m_\star(\oscr_{G_r \rtimes G_s}) \otimes V_r)^{G_r \rtimes P_s}$ which proves that $\mathcal{V}$ is indeed an object of  $\LB_{(\mathfrak{g},G)}(\mathcal{FL})$.  The reverse computation is left to the reader. 
 \end{proof}

\begin{rem}%
  \label{rem:G-equi-to-g-G-equi}
  We can use the group homomorphism $G_e \rtimes G
  \rightarrow G$, $(g_e,g)\mapsto g_eg$ to turn a $G$-equivariant sheaf into a $G_e \rtimes G$-equivariant sheaf. This defines   a natural fully faithful functor $
  \LB_G(\mathcal{FL}) \rightarrow \LB_{(\mathfrak{g},G)}(\mathcal{FL})$.
  We can also interpret this functor as saying   that a $G$-equivariant sheaf is naturally a
  $(\mathfrak{g}, G)$-equivariant sheaf by differentiating the
  $G$-action. %
    There is also a forgetful functor   $\LB_{(\mathfrak{g},G)}(\mathcal{FL}) \rightarrow \LB_G(\mathcal{FL})$. 

Via the equivalences of categories of Propositions~\ref{prop: equiv of cats BG Fl} and~\ref{prop: LB g G equiv of cats}, the functor from $G$-equivariant sheaves to $(\mathfrak{g}, G)$-equivariant sheaves amounts to associating to a $P$-representation the $(\mathfrak{g}, P)$-representation with trivial $\mathfrak{g}$-action. 
Indeed, the $(\mathfrak{g},G)$-equivariant sheaf corresponding to a $(\mathfrak{g},P)$-representation
$V$ is $\mathscr{F} = (m_{*}\oscr_{G_e \rtimes G} \otimes V)^{G_e \rtimes
  P}$.  If $V$ has trivial $\mathfrak{g}$-action, then $ f( g_e g'_e,
(g'_e)^{-1} g) = f(g_e, g)$, and we deduce that $f(g_e, gg'_e) = f( g_e g(g_e)' g^{-1}, g )$,
which means that the two $\mathfrak{g}$-actions (the obvious one and the one coming from the $G$-action) coincide. The converse implication is similar. 
The forgetful functor from $(\mathfrak{g}, G)$-equivariant sheaves to
$G$-equivariant sheaves corresponds to the forgetful functor associating to a $(\mathfrak{g}, P)$-representation the underlying $P$-representation. 
\end{rem}%

                                               \begin{defn}\label{defn:LBguP01}
                                                 We let
                                                 $\LB_{(\mathfrak{g},G)}(\mathcal{FL})^{\mathfrak{u}_P^0}$
                                                 be the full
                                                 subcategory of
                                                 $\LB_{(\mathfrak{g},G)}(\mathcal{FL})$
                                                 of objects which are
                                                 annihilated by
                                                 $\mathfrak{u}_{P}^0$.
                                               \end{defn}

We let $\LB_{(\mathfrak{g}, P)}(E)^{\mathfrak{u}_{P}}$ be the full subcategory of  $\LB_{(\mathfrak{g}, P)}(E)$ whose objects have the
property  that the $\mathfrak{u}_P$-action coming from differentiating
the $P$-action coincides with the $\mathfrak{u}_P$-action coming from
the $\mathfrak{g}$-action. (We will shortly see that this is
equivalent to the category $\LB_{(\mathfrak{g},G)}(\mathcal{FL})^{\mathfrak{u}_P^0}$.)

Any object of  $\LB_{(\mathfrak{g}, P)}(E)^{\mathfrak{u}_P}$ carries an action $\Theta$ of $Z({\mathfrak{m}})$ defined as follows. The differentiation of the $P$-action gives an action $d\rho_P$ of $\mathfrak{p}$. On the other hand the action $\rho_{\mathfrak{g}}$ of $\mathfrak{g}$ restricts to an action of $\mathfrak{p}$. 
We then let $z \in {\mathfrak{m}}$ act via the formula 
$$ z\cdot v = d\rho_P(\tilde{z})\cdot v - \rho_{\mathfrak{g}}(\tilde{z})\cdot v$$ 
for any lift $\tilde{z}$ of $z$ in $\mathfrak{p}$ (this is independent of the lift).  This induces an action of $U(\mathfrak{m})$ and restricts to an action of $Z(\mathfrak{m})$. This action commutes with the $(\mathfrak{g},P)$-action. 

\begin{prop}\label{prop:equiv-of-u-torsion-cats} The equivalence of
  categories between $\LB_{(\mathfrak{g},G)}(\mathcal{FL})$ and
  $\LB_{(\mathfrak{g}, P)}(E)$ of Proposition~\ref{prop: LB g G equiv of cats} induces an equivalence between 
 $\LB_{(\mathfrak{g},G)}(\mathcal{FL})^{\mathfrak{u}_P^0}$ and $\LB_{(\mathfrak{g}, P)}(E)^{\mathfrak{u}_{P}}$. The actions $\Theta_{\hor}$ and $\Theta$ of $Z({\mathfrak{m}})$ correspond to each other. 
 \end{prop}

\begin{proof} Let $V \in \LB_{(\mathfrak{g}, P)}(E)^{\mathfrak{u}_{P}}$. On the sheaf $(m_\star \oscr_{G_e \rtimes G} {\otimes} V)^{G_e \rtimes P}$, we want to see that the action of $\mathfrak{u}_{P}^0$ is trivial. 
We have $((g_eg)^{-1} u (g_eg)).f(g_e, g)  = f( u g_e , g)$ for $u \in
(U_P)_e$.   On the other hand, our assumption implies that the action
of $(u,1) \in P_e \rtimes G \subseteq G_e \rtimes G$ is trivial on $V$ (since the actions of
$(u, u^{-1})$ and $(1,u^{-1})$ coincide). This tells us that  $f( u
g_e , g) = f(g_e,g)$, as required. The converse implication follows similarly. The actions of  $\Theta_{\hor}$ and
$\Theta$ correspond by construction. \end{proof}

\begin{rem}Let~$\Rep(M)$ be the category of algebraic representations
  of~$M$. Then the natural functor $\Rep(M) \rightarrow
  \LB_{(\mathfrak{g}, P)}(E)$ (induced by inflation from~$M$ to~$P$,
  and letting~$\mathfrak{g}$ act trivially, see Remark~\ref{rem:G-equi-to-g-G-equi}) factors through $\LB_{(\mathfrak{g},
    P)}(E)^{\mathfrak{u}_{P}}$, and the obvious  action of
  $Z({\mathfrak{m}})$ on $\Rep(M)$ induces the action
  $\Theta$. 
\end{rem}

\begin{rem} We have a natural functor \[\HH^0(\mathfrak{u}_P^0,-):
  \LB_{(\mathfrak{g},G)}(\mathcal{FL}) \rightarrow
  \LB_{(\mathfrak{g},G)}(\mathcal{FL})^{\mathfrak{u}_P^0},\] which can
  be defined as follows.  Via our equivalence of categories, it corresponds to the natural functor  $\HH^0(\mathfrak{u}_P,-): \LB_{(\mathfrak{g}, P)}(E) \rightarrow \LB_{(\mathfrak{g}, P)}(E)^{\mathfrak{u}_P}$.  In this last formula, the $\mathfrak{u}_P$ action is the diagonal  one. More precisely, on any object $V$ of $\LB_{(\mathfrak{g}, P)}(E)$, we can differentiate the $P$-action to obtain a $\mathfrak{p}$-action. We therefore have an action of $\mathfrak{g} \rtimes \mathfrak{p}$ and $\mathfrak{u}_P$ embeds diagonally via $u \mapsto (-u,u)$ as a normal sub-Lie algebra.  \end{rem}

\subsubsection{Twisted differential operators and the sheaf $\mathcal{C}^{\la}$}\label{section-def-cla}

\begin{defn} Let $\tilde{\mathcal{D}}^{\la} = \oscr_{\mathcal{FL}} {\otimes} \hat{U}(\mathfrak{g})/\mathfrak{u}_P^0\oscr_{\mathcal{FL}} \otimes \hat{U}(\mathfrak{g})$ be the ring of universal twisted differential operators. 
\end{defn} 
\begin{rem} We have that $\mathcal{D}^{\la} =  \oscr_{\mathcal{FL}} {\otimes} \hat{U}(\mathfrak{g})/\mathfrak{p}^0\oscr_{\mathcal{FL}} \otimes \hat{U}(\mathfrak{g})$ is the usual ring of differential operators.
\end{rem}%

\begin{rem} One also has an ``algebraic'' version of
  $\tilde{\mathcal{D}}^{\la}$. Namely, we let
  $\tilde{\mathcal{D}}^{\alg} = \oscr_{\mathcal{FL}^{\alg}} {\otimes}
  {U}(\mathfrak{g})/\mathfrak{u}_P^{0,\alg}\oscr_{\mathcal{FL}^{\alg}}
  {\otimes} {U}(\mathfrak{g})$ be the ring of (algebraic) universal
  twisted differential operator on the $E$-scheme
  $\mathcal{FL}^{\alg}$. %
\end{rem}

We have three commuting actions of $G$ on $\oscr_{G} \otimes \oscr_{\mathcal{FL}}$:
\begin{enumerate}
\item $h \star_1 f(g,x) = f(h^{-1} g,x)$,
\item $h \star_2 f(g,x) = f( gh,x)$,
\item $h \star_3 f(g,x) = f( g,xh)$.
\end{enumerate}
We write $\star_{1,3}$  for the composition of the $\star_1$ and $\star_3$
action, and similarly for~$\star_{1,2,3}$ and so on.

\begin{defn}  We let $\mathcal{C}^{\la}  = ( \oscr_{G,e} \otimes \oscr_{\mathcal{FL}} )^{\mathfrak{u}_P^0}$ where the invariants are for the $\star_{1,3}$-action.
\end{defn}%
 Elements of $\mathcal{C}^{\la}$ are functions $f(g,x)$ with $g \in G_e$, $x \in \mathcal{FL}$, satisfying $f( u_x g, x) = f(g,x)$ for $u_x \in U_{P_x, e}$.

In particular $\mathcal{C}^{\la}$  has a $\star_{1,3}$-action
of $\mathfrak{g}$ and a $\star_{1,2,3}$
action of $G$, and it is easy to check that this gives
it the structure of an object of
$\LB_{(\mathfrak{g},G)}(\mathcal{FL})^{\mathfrak{u}_P^0}$. It has an
extra linear $\star_2$-action of $\mathfrak{g}$.  Its fiber at $e$ is the
module $\mathcal{C}^{\la}_e = \oscr_{ U_P \backslash G,e}$. %
  Under the equivalence of \ref{prop: LB g G equiv of cats}, the
$(\mathfrak{g}, P)$-module structure is the conjugation action of $P$ and the
right translation action of $\mathfrak{g}$.  The linear $\star_2$-action of
$\mathfrak{g}$ also induces the right translation action on the fiber. %

\begin{rem} The subsheaf $(\mathcal{C}^{\la})^{\mathfrak{p}^0}  \subseteq (\oscr_{G,e} \otimes \oscr_{\mathcal{FL}})^{\mathfrak{u}_P^0}$ is some kind of  infinite jet bundle over $\oscr_{\mathcal{FL}}$. 
\end{rem}

 We have a map $\Theta_{\hor}: Z(\mathfrak{m}) \rightarrow
 \mathrm{End}_{\oscr_{\mathcal{FL}}} (\mathcal{C}^{\la})$. We also
 have a map $\star_2: Z(\mathfrak{g}) \rightarrow
 \mathrm{End}_{\oscr_{\mathcal{FL}}}(\mathcal{C}^{\la})$ induced by
 the $\star_2$ action of~$\mf{g}$. %
 These maps are related as follows:

\begin{lem}\label{lem-all-compatibilities} We let $\iota: Z(\mathfrak{g}) \rightarrow Z(\mathfrak{g})$ be the map induced by the inverse  on $G$. For any $z \in Z(\mathfrak{g})$, we have $\star_{2}(\iota z) =
  \Theta_{\hor}(HC(z))$, where~$HC$ is the map~\eqref{eqn:HCmap}.
\end{lem} 
\begin{proof}  The endomorphisms of $\mathcal{C}^{\la}$,  $\star_{2}(\iota z)$ and 
  $\Theta_{\hor}(HC(z))$ are $G$-equivariant and $\oscr_{\mathcal{FL}}$-linear. Therefore it suffices to understand what happens on the fiber at $e$, namely $\oscr_{ U_P \backslash G,e}$. We first observe that on $\oscr_{G,e}$ the two actions $\star_1$ and $\star_2$ of $G$ induce two actions of $Z(\mathfrak{g})$ and these actions are related by $\star_2(z) = \star_{1}(\iota z)$. %
For the second point, by the definition of the map $HC$, for any $z \in Z(\mathfrak{g})$, we have that $HC(z) = z + z'$ where $z' \in U(\mathfrak{g}) \mathfrak{u}_{P}$. It follows that $\star_1(z)$ and $HC(z)$ act in the same way on $\oscr_{ U_P \backslash G,e} = \HH^0(\mathfrak{u}_P, \oscr_{G,e})$. 
\end{proof}

We have a   left action of $\mathfrak{g}$ on $\oscr_{G,e}$, defined by
$g.f(\cdot) = f'(\mathrm{exp}(-tg)\cdot)\vert_{t=0}$. This induces a natural pairing $\hat{U}(\mathfrak{g})  \otimes \oscr_{G,e} \rightarrow E$, $(\mathfrak{g}, f) \mapsto (g.f)(e)$. This pairing induces a pairing 
$(\oscr_{G,e} {\otimes} \oscr_{\mathcal{FL}}) \otimes (\hat{U}(\mathfrak{g}) {\otimes} \oscr_{\mathcal{FL}}) \rightarrow \oscr_{\mathcal{FL}}$.
It passes to a pairing on the quotient: 
$$ \tilde{\mathcal{D}}^{\la} \otimes \mathcal{C}^{\la} \rightarrow \oscr_{\mathcal{FL}}.$$

\begin{prop}\label{prop-dualCla} We have that $\underline{\mathrm{RHom}}_{\oscr_{\mathcal{FL}}}(\mathcal{C}^{\la}, \oscr_{\mathcal{FL}}) = \tilde{\mathcal{D}}^{\la}$. 
\end{prop}
\begin{proof} We will prove that $\underline{\mathrm{RHom}}_{\oscr_{\mathcal{FL}}}((\oscr_{G,e} \otimes \oscr_{\mathcal{FL}}) , \oscr_{\mathcal{FL}}) = \hat{U}(\mathfrak{g}) \otimes \oscr_{\mathcal{FL}}$. Since $\mathcal{C}^{\la}$ is locally a direct summand in $(\oscr_{G,e} {\otimes} \oscr_{\mathcal{FL}}) $ and  $\tilde{\mathcal{D}}^{\la}$ is locally a direct summand in  $\hat{U}(\mathfrak{g}) {\otimes} \oscr_{\mathcal{FL}}$, this implies the claim. 
We take a presentation $\oscr_{G,e} = \colim_r V_r$ where the $V_r$ are Smith spaces. We deduce that 
\begin{eqnarray*}
\underline{\mathrm{RHom}}_{\oscr_{\mathcal{FL}}}(\oscr_{G,e} \otimes \oscr_{\mathcal{FL}}, \oscr_{\mathcal{FL}})&=& \mathrm{Rlim}_r \underline{\mathrm{RHom}}_{E}(V_r, \oscr_{\mathcal{FL}}) \\
&=& \mathrm{Rlim}_r (V_r^\vee \otimes \oscr_{\mathcal{FL}})\\
&=& \mathrm{lim}_r (V_r^\vee \otimes \oscr_{\mathcal{FL}})\\
&=& \hat{U}(\mathfrak{g}) \otimes \oscr_{\mathcal{FL}}
\end{eqnarray*}
Here, the first equality is formal, the second equality is a
consequence of the nuclearity of Banach spaces \cite[Cor.\ 3.7]{MR4475468}, the third equality follows from Mittag-Leffler
\cite[Lem.\ 3.27]{MR4475468} and the last equality follows from
\cite[Lem.\ 3.28]{MR4475468}.
\end{proof}

\begin{rem}Using
  \cite[ Lem.\ 3.10]{MR4475468}, one can prove conversely that
  $\underline{\mathrm{Hom}}_{\oscr_{\mathcal{FL}}}(\tilde{\mathcal{D}}^{\la},
  \oscr_{\mathcal{FL}}) = {\mathcal{C}}^{\la}$. Conjecture 3.41 of  \cite{MR4475468}
  would imply that furthermore   $\underline{\mathrm{Hom}}_{\oscr_{\mathcal{FL}}}(\tilde{\mathcal{D}}^{\la}, \oscr_{\mathcal{FL}}) = \underline{\mathrm{RHom}}_{\oscr_{\mathcal{FL}}}(\tilde{\mathcal{D}}^{\la}, \oscr_{\mathcal{FL}})$. 
\end{rem} 

\subsubsection{The admissible objects}%
We let $\Adm_{(\mathfrak{g},P)}(E)$ be the subcategory of $\LB_{(\mathfrak{g},P)}(E)$ whose objects are admissible $\hat{U}(\mathfrak{g})$-modules. This is an abelian category.    We let $\Adm_{(\mathfrak{g}, G)}(\mathcal{FL})$ be the subcategory of $\LB_{(\mathfrak{g},G)}(\mathcal{FL})$ which corresponds to $\Adm_{(\mathfrak{g},P)}(E)$.

\begin{example} We see that $\mathcal{C}^{\la}$ is an object of $\Adm_{(\mathfrak{g},G)}(\mathcal{FL})$.
\end{example}%

\subsection{Equivariant sheaves on Bruhat
  cells}\label{subsec:equivariant-sheaves-Bruhat-strata}

We will now consider the stratification of $\mathcal{FL}$ into its
                                               $B$-orbits (i.e.\
                                               Bruhat cells).  %

\subsubsection{$(\mf{g},Q)$-equivariant sheaves}%
Recall
  that $\WM$ denotes the Kostant representatives of $W_M \backslash W$, and for $w \in \WM$ we let $C_w =
P\backslash P w B$ be the Bruhat cell. More
generally,  let~$Q$ be a standard parabolic, i.e.\ $B\subseteq
Q$, and write~$M_Q$ for the Levi quotient of~$Q$. We write
$\CwQ:=P\backslash PwQ \stackrel{j}\hookrightarrow \mathcal{FL}$ for the corresponding $Q$-orbit in~$\mathcal{FL}$. %
We write $\CwQdag = \lim_{ \CwQ \subseteq U} U$ where $U$ runs through the neighborhoods of $\CwQ$ in $\mathcal{FL}$. By definition $\CwQdag$ is the space $\CwQ$ equipped  with the sheaf $\oscr_{\CwQdag} =   j^{-1}\oscr_{\mathcal{FL}}$.

We consider the semi-direct product $G_e\rtimes Q$.   We have a  product map $G_e \rtimes Q \rightarrow  G $, $(g,q) \mapsto gq$. The group $G_e \rtimes Q$ acts on $\CwQdag$.

\begin{defn}
  We let $\LB_{(\mathfrak{g},Q)}(\CwQdag)$ be the category whose
  objects consist of the following list of data:
  \begin{enumerate}
  \item An $\LB$-sheaf $\mathscr{F}$ over $\CwQdag$: this is a sheaf
    of $\oscr_{\CwQdag}$-modules such that there is a covering
    $\CwQ = \cup_i U_i$ by quasi-compact opens, and for each $i$ a family $\{U_{i,j}\}_{j}$ of quasi-compact opens of $\mathcal{FL}$ with $U_{i,j}\cap \CwQ=U_i$ and
    $\cap_j U_{i,j} = U_i$, and Banach sheaves
    $ \mathscr{F}_{i,j}$ over $U_{i,j}$ such that
    $\{\mathscr{F}_{i,j} \vert_{U_i}\}_j$ form an inductive system and
    $\mathscr{F}\vert_{U_i} = \colim_j \mathscr{F}_{i,j}|_{U_i}$.
  \item We have a $Q$-equivariant map
    $\mathfrak{g} \otimes \mathscr{F} \rightarrow \mathscr{F}$ of
    $\LB$-sheaves providing a Lie algebra action of $\mathfrak{g}$ on
    $\mathscr{F}$.
  \item For each~$i$, there exists $s(i)$ such that each $U_{i,j}$ is stable under
    the action of $G_j \rtimes Q_{s(i)}$ and $\mathscr{F}_{i,j}$ is an
    object of $B_{G_j \rtimes Q_{s(i)}}(U_{i,j})$ and the maps
    $\mathscr{F}_{i,j} \rightarrow \mathscr{F}_{i,j'}$ are equivariant
    for the maps
    $G_{j'} \rtimes Q_{s(i)} \rightarrow G_{j} \rtimes Q_{s(i)}$.
  \item The induced action of $G_e \rtimes Q_{s(i)}$ on
    $\mathscr{F}\vert_{U_i}$ coincides with the restriction of the
    $(\mathfrak{g},Q)$-action.
  \end{enumerate}
\end{defn}

 \begin{rem}\label{rem-defn-LBgQ} Similar to remark \ref{rem-defn-LBgG}, we also can define  a category $\LB_{(\mathfrak{g},Q)}((\CwQdag)_{\C_p})$ with a base change functor $\LB_{(\mathfrak{g},Q)}(\CwQdag) \rightarrow \LB_{(\mathfrak{g},Q)}((\CwQdag)_{\C_p})$.
 \end{rem}

\subsubsection{An equivalence of categories}\label{sec-equ-cat-repLB-Q}

We consider the uniformization: 
\begin{eqnarray*}
m:  G_e \rtimes Q  &\rightarrow &\CwQdag \\
 (g,q) & \mapsto & w gq.
 \end{eqnarray*} 
 We let $\StabQw$ be the stabilizer  of $w$ for this action, so that
 \[\StabQw = \{ (g,q) \in G_e \rtimes Q, gq \in P_w\}.\]
 We have an injective homomorphism \[\StabQw\into P_w\times Q\]
 given by $(g,q)\mapsto(gq,q)$, which induces an isomorphism
 \numequation\label{eqn:StabQe-description}\StabQw_e\isoto
 P_{w,e}\times Q_{e}.\end{equation} From now on we will frequently
identify $P_{w,e}\times Q_{e}$ with $\StabQw_e$
via~\eqref{eqn:StabQe-description}, and in particular we will frequently
regard~$Q_e$ as a subgroup of $ G_e \rtimes Q$
via~\eqref{eqn:StabQe-description} (and the inclusion
$\StabQw\subset G_e \rtimes Q$), i.e.\ as the subgroup of elements
$(q^{-1},q)$ with $q\in Q_e$.

 \begin{lem}\label{lem-decompose-stab(w)-Q}\leavevmode
   \begin{enumerate}
   \item The group $\StabQw$ is generated by its  subgroups
     $\StabQw_e=P_{w,e}\times Q_{e}$ and~$P_w\cap Q$.
   \item   There is an
     isomorphism %
     $ (P_w \cap Q)_e \backslash \big((P_{w,e} \times Q_e) \rtimes (P_w
     \cap Q)\big) \isoto \StabQw$. %
   \end{enumerate}
 \end{lem}
 \begin{proof} Consider an element  $(g,q)\in\StabQw\subseteq G_e
   \rtimes Q$, so that  $gq
   \in P_w$. Then~$g\in G_e\cap P_wQ=P_{w,e}Q_e$, so we can write
   $g=pq'$ with~$p\in P_{w,e},q'\in Q_e$. Then~$q'q\in Q$, and
   since~$gq\in P_w$ and~$p\in P_w$, we in fact have $q'q\in P_w\cap
   Q$. Thus we can write \[(g,q)=(p,1)(q',q'^{-1})(1,q'q)\in G_e
   \rtimes Q\] with $(p,1)\in P_{w,e}$, $(q',q'^{-1})\in Q_e$  and $(1,q'q)\in P_w\cap Q$,
 completing the proof of the first part.

 It follows from the first part  that  we have a surjective
homomorphism \[\StabQw_e\rtimes ( P_w\cap Q)\to \StabQw\]given
by $(g,q)\mapsto gq$. The
kernel of this homomorphism is $\StabQw_e\cap (P_w\cap Q)=(P_w\cap Q)_e$, and the second
part is immediate.   
 \end{proof}

 \begin{cor} A representation $(V, \rho)$ of $\StabQw$ is  the data
   of a  representation $(V, \rho_1)$ of $P_{w,e}$, a representation
   $(V, \rho_2)$ of $Q_e$ and a representation $(V,\rho_3)$ of $P_w \cap Q$, satisfying:
 
\begin{enumerate}
\item $\rho_1$ and $\rho_2$ commute.
\item $\mathrm{Ad}\rho_3(a) ( \rho_1(b)\rho_2(c)) = \rho_1( a ba^{-1}) \rho_2(a c a^{-1})$. 
 \item  $\rho_3 = \rho_1  \rho_2$ on $(P_w \cap Q)_e$. 
\end{enumerate}
\end{cor}
\begin{proof}
  This is immediate from Lemma~\ref{lem-decompose-stab(w)-Q}.
\end{proof}
 
 \begin{example}\label{ex:stalk-of-Cla-at-w}%
   We consider the sheaf $\mathcal{C}^{\la}\vert_{\CwQdag}$. We see that the
   fiber $\mathcal{C}^{\la}_w$ is $\oscr_{U_{P_w} \backslash G, e}$, and the action of $\StabQw$ is given by $(g,q) f(g')= f( q^{-1} g^{-1} g' q)$. In other words, it has a $P_{w,e}$-action by left translation, a $Q_e$-action by right translation, and a $P_w \cap Q$-action by conjugation. 
 \end{example}

If $r, s \geq 1$, we can consider the semi-direct product $G_r \rtimes
Q_s$. We let $\StabQw_{r,s}$ be the stabilizer of $w$ for the
action of $G_r \rtimes Q_s$; again, this is the subgroup of $G_r
\rtimes Q_s$ of elements $(g,q)$ such that $gq \in P_w$. 

\begin{lem} If $r \geq s$, we have an isomorphism:  \[ (P_w \cap Q)_r \backslash \big((P_{w,r} \times Q_r) \rtimes (P_w
     \cap Q_{s})\big) \isoto \StabQw_{r,s}.\]
\end{lem}
\begin{proof} This follows exactly as in the proof of Lemma \ref{lem-decompose-stab(w)-Q}.
\end{proof}
We now define a  category $\LB_{\StabQw}(E)$ as follows.
\begin{defn}\label{defn-Stabw-rep-Q}
  The category $\LB_{\StabQw}(E)$ has objects
  consisting of the following list of data:
  \begin{enumerate}
  \item An $\LB$-space $V = \colim_r V_r$ over $E$,
  \item An action of $\StabQw$ on $V$.
  \item For each~$r$, there exists $s$ such that $V_r \in B_{\StabQw_{r,s}}(E)$ and
    the maps $V_r \rightarrow V_{r'}$ are equivariant with respect to
    the maps $\StabQw_{r',s} \rightarrow \StabQw_{r,s}$.
  \item This induces on $V$ an action of $\StabQw_e$ and
    $P_w\cap Q_s$ which coincides with the restriction of the
    $\StabQw$-action of $V$.
  \end{enumerate}
\end{defn}
 \begin{prop}\label{prop-gQeq} Taking the fiber at~$w$ gives an equivalence between the categories  $\LB_{(\mathfrak{g}, Q)}(\CwQdag)$ and $\LB_{\StabQw}(E)$. 
 \end{prop}
\begin{proof}  As usual, by taking the fiber at $w$, we obtain a $\StabQw$-module. Conversely, we attach to an object $V$ of $\LB_{\StabQw}(E)$ the sheaf $(m_\star (\oscr_{G_e \rtimes Q}) \hat{\otimes }V)^{\StabQw}$. 
 \end{proof}

 \begin{rem} We have a restriction map
   $\LB_{(\mathfrak{g},G)}(\mathcal{FL}) \rightarrow
   \LB_{(\mathfrak{g},Q)}(\CwQdag)$. The category
   $\LB_{(\mathfrak{g},G)}(\mathcal{FL}) $ is equivalent to
   $\LB_{(\mathfrak{g}, P_w)}(E)$ by Proposition~\ref{prop: LB g G
     equiv of cats}, %
 and  this
   restriction map corresponds to the map
 $\LB_{(\mathfrak{g}, P_w)}(E) \rightarrow \LB_{\StabQw}(E)$ which is induced by the inclusion $\StabQw \into G_e \rtimes P_w$ where 
 $ q\in Q_e \mapsto (q, 1)$, $p \in P_{w,e} \mapsto (p^{-1}, p)$, $r \in P_w \cap Q \mapsto (1,r)$. 
 \end{rem}
 
 \begin{rem}\label{rem-lambda-smooth-Q} For any object $\mathscr{F}$ of $\LB_{(\mathfrak{g}, Q)}(\CwQdag)$, one can differentiate the $Q$-equivariant structure and thus obtain a $Q$-equivariant map $\act_{\mf{q}}: \mf{q} \otimes \mathscr{F} \rightarrow \mathscr{F}$.  One may want to compare this map with the map $\act_{\mathfrak{g}}: \mathfrak{g} \otimes \mathscr{F} \rightarrow \mathscr{F}$ which is given by the $\mathfrak{g}$-equivariant action. The difference $\act_{\mf{q}} - \act_{\mathfrak{g}}\vert_\mf{q}: \mf{q} \otimes \mathscr{F} \rightarrow \mathscr{F}$ is a $Q$-equivariant linear map. It is therefore entirely determined by its fiber at $w$. If $\mathscr{F}$ corresponds to a $\StabQw$-representation $V$ via the equivalence  of Proposition \ref{prop-gQeq}, then the $\mf{q}$-action $\act_{\mf{q}} - \act_{\mathfrak{g}}\vert_\mf{q}$  is induced by the $Q_e$-action on $V$. 
 \end{rem} 
 
 We let $\LB_{(\mathfrak{g}, Q)}(\CwQdag)^{\mathfrak{u}_{P}^0}$ be
 the subcategory of objects which are killed by $\mathfrak{u}_P^0$. We
 let $\LB_{\StabQw}(E)^{\mathfrak{u}_{P_w}}$ be the full subcategory
 of objects with trivial action of the subgroup $(U_{P_w})_e
 \hookrightarrow P_{w,e}$. The objects of $\LB_{\StabQw}(E)^{\mathfrak{u}_{P_w}}$ carry an action $\Theta$ of
 $Z({\mathfrak{m}_w})$.  %

 \begin{prop}\label{prop:equiv-LB-Cw-Stab-u-Q} The equivalence  of categories between the categories $\LB_{(\mathfrak{g}, Q)}(\CwQdag)$ and $\LB_{\StabQw}(E)$ induces an equivalence between $\LB_{(\mathfrak{g}, Q)}(\CwQdag)^{\mathfrak{u}_{P}^0}$ and $\LB_{\StabQw}(E)^{\mathfrak{u}_{P_w}}$. Via this equivalence, the action of $\Theta_{\hor}$  of $Z({\mathfrak{m}})$ corresponds to the action $\Theta$ of $Z({\mathfrak{m}_w})$ via conjugation by $w^{-1}$. 
 \end{prop}
 \begin{proof}
This follows from Proposition~\ref{prop-gQeq}, exactly as in the proof
of Proposition~\ref{prop:equiv-of-u-torsion-cats}.
 \end{proof}
 
\begin{rem}\label{rem-twist-Q} Let $\mu \in X^\star(M_Q^{\ab})$ be an algebraic character. There is a functor $\LB_{(\mathfrak{g}, Q)}(\CwQdag) \rightarrow \LB_{(\mathfrak{g}, Q)}(\CwQdag)$, $\mathscr{F} \mapsto \mathscr{F} \otimes E(\mu)$, corresponding to twisting the $Q$-action by $\mu$. 
There is a map $\StabQw \hookrightarrow G_e \rtimes Q   \rightarrow Q
\rightarrow M_Q^{\ab}$, so that any character $\mu \in X^\star(M_Q^{\ab})$ induces a character of $\StabQw$. The operation of twisting the $Q$-action by $\mu$ corresponds to the operation of twisting a $\StabQw$-representation by $\mu$. 
\end{rem}

\subsection{Algebraic and locally analytic representations}\label{subsec:algebraic-locally-analytic-representations}  In this
section we will explain how to define a functor from a subcategory of the algebraic category $\ocal(\mathfrak{m}_w,\mathfrak{b}_{\mathfrak{m}_{w}})$ %
to representations of $\StabQ(w)$. We begin with some more general
considerations. %

\subsubsection{Completion of category~$\cO$}\label{subsubsec:completion-category-O}  In this subsection we consider a reductive group $G$ with Borel $B$ and maximal torus $T$. Its   Lie algebra is $\mathfrak{g}$ with Borel $\mathfrak{b}$ and Cartan $\mathfrak{h}$. We recall that $\ocal(\mathfrak{g}, \mathfrak{b})$ is the corresponding BGG subcategory of $U(\mathfrak{g})$-modules. Let $B \subseteq Q \subseteq G$ be a parabolic with Levi $M$. Let  $\mathfrak{q} \subseteq \mathfrak{g}$ be its Lie algebra, with Levi $\mathfrak{m}$. We let $\ocal(\mathfrak{g}, \mathfrak{q})$ be the parabolic BGG category, which is the subcategory of $\ocal(\mathfrak{g}, \mathfrak{b})$ of objects whose restriction to $\mathfrak{m}$ is a direct sum of finite dimensional representations. 

We start with the following definition which is nothing but Definition  \ref{defn:LBgX} specialized to $X = \Spa(E, \ocal_E)$. 
\begin{defn}
  We let $\LB_{\mathfrak{g}}(E)$ be the category of $\mathfrak{g}$-representations on $\LB$-spaces. More precisely, its objects are $\LB$ spaces over
  $E$, $V$ which are $\hat{U}(\mathfrak{g})$-modules and satisfy the
  following conditions:

  \begin{enumerate}
  \item We have $V= \colim_r V_r$ and for $r$ large enough
    $V_r \in B_{G_r}(E)$. Moreover, the transition maps
    $V_{r} \rightarrow V_{r'}$ for $r \leq r'$ are equivariant for the
    map $G_{r} \rightarrow G_{r'}$.
  \item The actions of $G_r$ on $V_r$ induce the action of
    $\mathfrak{g}$ on the limit.
  \end{enumerate}
\end{defn}

 \begin{prop} \label{prop:completion-functor-on-cat-O}

There is an exact contravariant functor: 
\begin{eqnarray*}
\ocal(\mathfrak{g},\mathfrak{b}) &\rightarrow& \LB_{\mathfrak{g}}(E) \\
M& \mapsto & [M \otimes_{U(\mathfrak{g})}\hat{U}(\mathfrak{g})]^\vee = \hat{M}^{\vee}
\end{eqnarray*}
\end{prop}

\begin{proof} This follows from Theorem \ref{thm-schmidt}. Indeed, our
  functor is the composition of the functor $\ocal \rightarrow
  \hat{\ocal}$ which is an equivalence of abelian categories, and then of the
  duality functor (which is an exact anti-equivalence of categories, and
  turns a coadmissible $\hat{U}(\mathfrak{g})$-module into an
  admissible $\hat{U}(\mathfrak{g})$-module), and finally the forgetful
  functor to the category of $\LB$-spaces equipped with a
  $\mathfrak{g}$-action. It remains to justify that    $[M
  \otimes_{U(\mathfrak{g})}\hat{U}(\mathfrak{g})]^\vee$ belongs to
  $\LB_{\mathfrak{g}}(E)$.  To see this, note that $\hat{U}(\mathfrak{g}) \otimes_{U(\mathfrak{g})} M= \lim_r M \otimes_{U(\mathfrak{g})} D(G_r)$, so that $\hat{M}^\vee  = \colim_r V_r$ where $V_r = (M \otimes_{U(\mathfrak{g})} D(G_r))^\vee$ is a  Banach space. Moreover,  there is  an action of $G_r$ on $V_r$. 
\end{proof}

\begin{rem}\label{rem-making-explicit} We  can explicate what this
  functor is doing on Verma modules. Let $\lambda \in X^\star(T)_{E}$
  be a character of $\mathfrak{b}$. Let $M = U(\mathfrak{g})
  \otimes_{U(\mathfrak{b})} \lambda$. We see that $\hat{M}^\vee$ is
  the submodule of $\oscr_{G,e}$ of functions $f$ which satisfy $f(gb)
  = \lambda(b)f(g)$ for $(g,b) \in G_e \times B_e$, with the action
  of~$\mf{g}$ being that given by the action of~$G_e$ as $(gf)(x)=f(g^{-1}x)$. More generally, let $V$ be a finite dimensional representation of $\mathfrak{b}$ with dual $V^\vee$. Let $M = U(\mathfrak{g}) \otimes_{U(\mathfrak{b})} V$. We see that $\hat{M}^\vee$ is the submodule of the space of functions $f: {G_e} \rightarrow V^\vee$ which satisfy $bf(gb) = f(g)$. 
\end{rem}

We let $X^\star(M)$ be the character space of $M$. 

\begin{defn}\label{defn:O-lambda-alg}%
  Let $\lambda \in X^\star(M)_{E}$. We let $\ocal(\mf{g},\mf{q})_{\lambda-\alg}$ be
  the full subcategory of $\cO(\mf{g},\mathfrak{q})$, whose
  objects are those~
  $V$ which have the property that in the weight
  decomposition $V = \oplus_{\nu \in X^\star(T)_{E}} V[\nu]$ for the
  action of $\mathfrak{h}$\footnote{ By definition,
    $V[\nu] = \{ m \in V, h.m = \nu(h)m\}$. In the direct sum, we
    suppose $V[\nu] \neq 0$.}, we have $\nu - \lambda \in X^\star(T)$.
  This is an abelian category.%
\end{defn}
\begin{lem}%
  \label{rem:M-twisted-minus-lambda}If~$V\in
  \ocal(\mf{g},\mf{q})_{\lambda-\alg}$, then the $\mathfrak{q}$-action on the twisted module ~$V(-\lambda)$ integrates to  an
  action of~$Q$. 
\end{lem}  
  \begin{proof} We observe that $V(-\lambda)$ is a union of finite dimensional
    representations of $\mathfrak{q}$. We claim that on any finite dimensional
    representation $(W, d\rho)$ of $\mathfrak{q}$,  the action integrates to an
    action $(W, \rho)$ of $Q$ as long as the action of $\mathfrak{h}$ integrates
    to an action of $T$ (which is the reason why we are introducing a twist).
    Indeed, the Lie algebra action gives a map $U(\mathfrak{q}) \otimes W
    \rightarrow W$ and since $W$ is finite dimensional, we can dualize this map
    to a map $W \rightarrow W \otimes \widehat{\oscr_{Q,e}}$ where $
    \widehat{\oscr_{Q,e}}$ is the completion of the local ring at $e$. We claim
    that this map factorizes through $W \otimes \oscr_{Q}$ and gives the
    coaction map; this establishes the lemma. 
  This claim must be well known but we could not find a reference so we sketch the argument.  Using the Levi decomposition $Q = M \ltimes U_{Q}$ it suffices to treat the case of $M$ and $U_{Q}$ separately. 
  
To see that   the action of the unipotent radical $U_Q$ of $Q$ integrates we can
just consider a root groups $U_\beta \simeq \mathbb{G}_a$ inside $U_Q$ with Lie
algebra $\mathfrak{u}_\beta$ generated by $u_\beta$, and the rule $V \rightarrow V \otimes E[T]$, $v \mapsto \exp(T u_\beta)v$ defines an action of $U_\beta$ (we use here that the action of $u_\beta$ is locally nilpotent). 
  On the other hand,   the action of $\mathfrak{h}$ on $V(-\lambda)$ integrates to an action of $T$. Similarly, the action of $M^{der}$ integrates (from the action of $\mathfrak{m}$). Indeed, for a semi-simple group, the categories of finite dimensional representations of the group and of its Lie algebra are equivalent. This actions of $M^{der}$ and $T$ combine together to an action of $M$.   
  \end{proof}

\subsubsection{A particular class of representations of $\StabQw$}
We now go back to our original setting. Let $Q_{M_w} = P_w \cap Q/U_{P_w \cap Q}$.

\begin{lem}\label{rem:Kostant-Borel}%
  If $w \in \WM$ then $Q_{M_w}$ is a parabolic subgroup of $M_w$ containing $B_{M_w}$.    %
\end{lem}
\begin{proof} Clearly, $Q_{M_w}$ contains  $P_w \cap B/U_{P_w \cap Q}$ which is a Borel subgroup of $M_w$. Any closed subgroup of a reductive group containing a Borel is a parabolic subgroup. 
\end{proof}

Let $\lambda \in X^\star(M_{Q})_E$.

\begin{defn}%
  We let $\Adm_{(\mathfrak{m}_w, Q_{M_w})}(E)_\lambda$ be the following
  category. Its objects are admissible
  $\hat{U}( \mathfrak{m}_{w})$-modules $V = \colim_r V_r$ which admit
  an action of $Q_{M_w}$, compatible with the action of $Q_{M_w}$ by
  conjugation on $D(M_{w,e})$. We further demand the following conditions
:
  \begin{enumerate}
  \item For $q \in Q_{M_w}$, $g \in \mathfrak{m}_w$, and $m \in V$, we
    have $\mathrm{Ad}(q)(g).m = q.g. q^{-1}.m$.
  \item there exists $s \in \ZZ_{>0}$ and an action of $Q_{M_w,s}$ on
    each $V_r$, inducing an action of $V = \colim_r V_r$ which
    coincides with the restriction of the action of $Q_{M_w}$ to
    $Q_{M_w,s}$.
  \item Let us denote by $\rho_{Q_{M_w}}$ the action of $Q_{M_w}$ on
    $V$. This action differentiates to an action
    $\mathrm{d} \rho_{Q_{M_w}}$ of $\mathfrak{q}_{M_w}$.  We let
    $\rho_{\mathfrak{m}_w}$ be the action of $\mathfrak{m}_w$.  Then
    we ask that
    $\mathrm{d} \rho_{Q_{M_w}} = \lambda +
    \rho_{\mathfrak{m}_w}\vert_{\mathfrak{q}_{M_w}} $. %
  \end{enumerate}
\end{defn}%

\begin{rem}%
  \label{rem:adm-is-full-in-LB}The category $\Adm_{(\mathfrak{m}_w,
    Q_{M_w})}(E)_\lambda$ is an abelian category by  general results on coadmissible and admissible modules over Fr\'echet--Stein algebras (see \cite{MR1990669}, sect. 3).

  The category $\Adm_{(\mathfrak{m}_w,
    Q_{M_w})}(E)_\lambda$ is a full subcategory of
  $\LB_{\mathfrak{m}_w}(E)$. In particular, since~$Q_{M_w}$ is
  connected, the action of~$Q_{M_w}$ in the third condition is
  uniquely determined by the $\mathfrak{m}_w$-action and $\lambda$.%
\end{rem}

\begin{rem} By Proposition~\ref{prop: LB g G equiv of cats} (applied with~$G$ replaced by~$M_{w}$), we see that objects of $\Adm_{(\mathfrak{m}_w,
    Q_{M_w})}(E)_\lambda$ define $M_w$-equivariant
  $D_{\lambda}$-modules on the flag variety $X = Q_{M_w} \backslash
  M_w$, where $D_\lambda:=\oscr_{X} \otimes \hat{U}(\mathfrak{m}_w)
  \otimes_{\mathfrak{q}_{M_w}^0} \lambda$ is a ring of twisted differential
  operators. %
\end{rem}

\begin{prop}\label{prop:full-faithful-Adm-to-LB-Stab}
There is a natural fully faithful functor: $\Adm_{(\mathfrak{m}_w, Q_{M_w})}(E)_\lambda \rightarrow \LB_{\StabQw}(E)^{\mathfrak{u}_{\mathfrak{p}_w}}$. 
\end{prop}

\begin{proof}It suffices to exhibit an equivalence of categories
  between  
  $\Adm_{(\mathfrak{m}_w, Q_{M_w})}(E)_\lambda $ and a full subcategory of
  $\LB_{\StabQw}(E)^{\mathfrak{u}_{\mathfrak{p}_w}}$. To this end, recall that  by Lemma \ref{lem-decompose-stab(w)-Q},
  $\StabQw = (P_w \cap Q)_e \backslash \big((P_{w,e} \times Q_e) \rtimes (P_w
     \cap Q)\big)$. %
 Accordingly,
     we may consider  the full subcategory  of representations of
     $\StabQw$ which factor through $(Q_{M_w})_e \backslash
     \big(((M_w)_e \times M_{Q,e}^{\ab} \big) \rtimes Q_{M_w}$ and have the property that $M_{Q,e}^{\ab}$ acts through the character $\lambda$.  
In other words, we consider the full subcategory of $\LB_{\StabQw}(E)$ whose objects $V$ satisfy the following  properties: 
\begin{enumerate}%
\item The action of $P_{w,e}$ factors through an action $\rho_1$ of
  $M_{w,e}$. Moreover $V$, viewed as an object of
  $\LB_{\mathfrak{m}_{w}}(E)$, is admissible. %
\item The action of $Q_e$ factors through %
  $M_{Q,e}^{\ab}$ %
  acting via $\lambda$. 
\item The action of $Q \cap P_w$ factors through an action $\rho_3$ of $Q_{M_{w}}$. 
\end{enumerate}
Clearly, this is equivalent to $\Adm_{(\mathfrak{m}_w, B_{M_w})}(E)_\lambda$. 
\end{proof}

We now consider the parabolic BGG category $\ocal(\mf{m}_w,\mf{q}_{M_w})$   for $\mathfrak{m}_w$ and the parabolic
$\mathfrak{q}\cap \mathfrak{m}_w = \mathfrak{q}_{M_w}$. %
We  think of
$\mathfrak{h} \hookrightarrow \mathfrak{b}$ as giving the Cartan of
$\mathfrak{m}_w$. We now apply the material of
Section~\ref{subsubsec:completion-category-O} %
to~$\cO(\mf{m}_w,\mathfrak{q}_{M_w})$. %

By Proposition~\ref{prop:completion-functor-on-cat-O} (applied
to~$\cO(\mathfrak{m}_w,\mf{q}_{M_w})$), we have a completion functor
$\cO\to \LB_{\m_w}(E)$, which restricts to a functor $\ocal(\mf{m}_w,\mf{q}_{M_w})_{\lambda-\alg}\to\LB_{\m_w}(E)$.

\begin{prop}\label{prop-completion-Stab} 
We can uniquely upgrade the  completion functor 
\begin{eqnarray*}
\ocal(\mf{m}_w,\mf{q}_{M_w})_{\lambda-\alg} &\rightarrow& \LB_{\mathfrak{m}_{w}}(E) \\
V& \mapsto & \hat{V}^\vee
\end{eqnarray*}
 to a fully faithful  functor
\begin{eqnarray*}
\ocal(\mf{m}_w,\mf{q}_{M_w})_{\lambda-\alg} &\rightarrow& \Adm_{(\mathfrak{m}_w, Q_{M_w})}(E)_\lambda   \\
V& \mapsto & \hat{V}^\vee(\lambda).
\end{eqnarray*} 
Its essential image is the subcategory of $\Adm_{(\mathfrak{m}_w,
  Q_{M_w})}(E)_\lambda$ of objects which are in the image of the
functor $\ocal(\mf{m}_w,\mf{q}_{M_w})_{\lambda-\alg} \rightarrow
\LB_{\mathfrak{m}_w}(E)$ when viewed as
$\mathfrak{m}_{w}$-representations. %
\end{prop}

\begin{rem}  As the notation $\hat{V}^\vee(\lambda)$ suggests, we make
  a twist of the action of $\mathfrak{q}_{M_w}$ on $\hat{V}^\vee$ by
  $\lambda$ so that it extends to an action of $Q_{M_w}$. (See also Lemma~\ref{rem:M-twisted-minus-lambda}.)
\end{rem}

\begin{proof}[Proof of Proposition~\ref{prop-completion-Stab}] %
  Let $V \in \ocal(\mf{m}_w,\mf{q}_{M_w})_{\lambda-\alg}$. We consider
  $\hat{V}^\vee$. This space carries an action of $P_{w,e}$ (factoring
  through $M_{w,e}$). We can also define an action of $(M_{Q})_e$ on~$\hat{V}^{\vee}$ via scalar
  multiplication by the character $\lambda$. Clearly these two actions commute. The
  product of the two actions defines an action of $(M_{Q})_e \cap P_{w,e}$
  factoring through $Q_{M_w,e}$. We claim that we can extend it to an
  action of $Q \cap P_w$ factoring through $Q_{M_w}$.
  It follows from Lemma \ref{rem:M-twisted-minus-lambda} that we have an action on $V(-\lambda)$. We thus get an action on $\hat{V}(-\lambda) = V(-\lambda) \otimes_{U(\mathfrak{m}_w)} \hat{U}(\mathfrak{m}_w)$, since $M_w$ acts on $\hat{U}(\mathfrak{m}_w)$ and $U(\mathfrak{m}_w)$ via the adjoint representation. 
  \end{proof} 

\begin{rem}\label{rem-on-twisting} Clearly the categories $\ocal(\mf{m}_w,\mf{q}_{M_w})_{\lambda-\alg}$ and $\ocal(\mf{m}_w,\mf{q}_{M_w})_{\lambda+\mu-\alg}$ are equivalent if $\mu \in X^\star(M_Q)$. However, the functor $\ocal(\mf{m}_w,\mf{q}_{M_w})_{\lambda-\alg} \rightarrow \Adm_{(\mathfrak{m}_w, Q_{M_w})}(E)_\lambda$ depends on the choice of $\lambda$ (as clearly the target category depends on $\lambda$). We have a functor $\Adm_{(\mathfrak{m}_w, Q_{M_w})}(E)_\lambda \rightarrow \Adm_{(\mathfrak{m}_w, Q_{M_w})}(E)_{\lambda + \mu}$, $V \mapsto V \otimes E(\mu)$.  We have the following commutative  diagram of functors (telling us that $\hat{V}^\vee(\lambda) \otimes E(\mu) = \hat{V}^\vee(\lambda + \mu)$):
\begin{eqnarray*} 
 \xymatrix{  \ocal(\mf{m}_w,\mf{q}_{M_w})_{\lambda-\alg} \ar[r] \ar[d]  &  \Adm_{(\mathfrak{m}_w, Q_{M_w})}(E)_\lambda \ar[d]\\
 \ocal(\mf{m}_w,\mf{q}_{M_w})_{\lambda+\mu-\alg} \ar[r] & \Adm_{(\mathfrak{m}_w, Q_{M_w})}(E)_{\lambda+\mu}}
 \end{eqnarray*}
\end{rem}

\subsubsection{Higher Coleman sheaves}
\begin{defn}\label{defn:HCS-functor} %
  We now define a contravariant exact functor
$$  HCS_{Q,w,\lambda}: \ocal(\mf{m}_w,\mf{q}_{M_w})_{\lambda-\alg} \rightarrow \LB_{(\mathfrak{g},Q)}(C_{w,Q}^\dag)^{\mathfrak{u}_{\mathfrak{p}}^0}$$(where
``HCS'' stands for ``higher Coleman sheaf'') as as the
composite
\[\ocal(\mf{m}_w,\mf{q}_{M_w})_{\lambda-\alg} \rightarrow \Adm_{(\mathfrak{m}_w,
    Q_{M_w})}(E)_{\lambda} \rightarrow
  \LB_{\StabQw}(E)^{\mathfrak{u}_{\mathfrak{p}_w}} \rightarrow
  \LB_{(\mathfrak{g},Q)}(C_{w,Q}^\dag)^{\mathfrak{u}_{\mathfrak{p}}^0},\]
where the first functor is the one defined in
Proposition~\ref{prop-completion-Stab}, the second is the fully
faithful functor of
Proposition~\ref{prop:full-faithful-Adm-to-LB-Stab}, and the third is
the equivalence of Proposition~\ref{prop:equiv-LB-Cw-Stab-u-Q}.

In the case~$Q=B$ we write $HCS_{w,\lambda}$ for $HCS_{B,w,\lambda}$.
\end{defn}
\begin{prop}\label{TwistingHCS} Let $M \in \ocal(\mf{m}_w,\mf{b}_{M_w})_{\lambda-\alg}$. Let $\mu \in X^\star(T)$. We have $HCS_{w, \lambda}(M) \otimes E(\mu) = HCS_{w, \lambda+\mu}(M)$. 
\end{prop}
\begin{proof} This follows from Remark \ref{rem-on-twisting}.
\end{proof}

We have an action of $Z(\mathfrak{m}_w)$ on
$\ocal(\mf{m}_w,\mf{b}_{M_w})_{\lambda-\alg}$. We also have an action
$Z(\mathfrak{m})$ on $
\LB_{(\mathfrak{g},B)}(C_w^\dag)^{\mathfrak{u}_{\mathfrak{p}}^0}$ via
$\Theta_{\hor}$. These two action are related by the following
lemma. We let $\iota: Z(\mathfrak{m}) \rightarrow Z(\mathfrak{m})$ be
the map induced by the inverse map on $M$, and let $w: \mathfrak{m}_w \rightarrow \mathfrak{m}$ be  conjugation by $w$. 

\begin{prop}\label{prop-action-centre-enveloping} Let us consider the map $w: Z(\mathfrak{m}_w) \rightarrow Z(\mathfrak{m})$. 
Then we have that $HCS_{w,\lambda}(z) = \Theta_{\hor}( \iota w z)$ for any $z \in Z(\mathfrak{m}_w)$.
\end{prop}
\begin{proof} This follows directly from the construction, bearing in mind Proposition~\ref{prop:equiv-LB-Cw-Stab-u-Q}.
\end{proof}

 \subsection{Localization on the partial flag variety}\label{subsec:localization-partial-flag-variety}
 \subsubsection{Statement of the localization problem}

Recall that  in Section \ref{section-def-cla} we  defined an object  $\mathcal{C}^{\la} \in \LB_{(\mathfrak{g},G)}(\mathcal{FL})$. This is a   $\tilde{\mathcal{D}}^{\la}$-module and it carries an action $\star_2$ of $\mathfrak{g}$ which commutes with the $\tilde{\mathcal{D}}^{\la}$-module structure. 

We define a localization  functor:   
\begin{eqnarray*}
\mathrm{Loc}: D^-(U(\mathfrak{g})) &\rightarrow &D( \tilde{\mathcal{D}}^{\la}) \\
M& \mapsto & \RHom_{\mathfrak{g}, \star_2}(M, \mathcal{C}^{\la})
\end{eqnarray*}%
where $D(\tilde{\mathcal{D}}^{\la})$ is the  derived category of solid
$\tilde{\mathcal{D}}^{\la}$-modules. We will sometimes drop the
subscript~$\star_2$ from the notation, and simply write $\RHom_{\mathfrak{g}}(M, \mathcal{C}^{\la})$.  %

Recall that if $M$ is an object of $\mathrm{Mod}^{\fing}(U(\mathfrak{g}))$
then we let $\hat{M} =
 M\otimes_{U(\mathfrak{g})}\hat{U}(\mathfrak{g} )$ %
and $\hat{M}^\vee = \underline{\mathrm{Hom}}( \hat{M},
E)$. The following lemma gives another description of~$\Loc(M)$ for $M \in \mathrm{Mod}^{\fing}(U(\mathfrak{g}))$.

\begin{lem} \label{lem-passingtotheotherside} Assume that $M \in
  \mathrm{Mod}^{\fing}(U(\mathfrak{g}))$. Then we have:
\begin{eqnarray*}
\RHom_{\mathfrak{g}, \star_2}(M, \mathcal{C}^{\la}) &=&  \RHom_{\mathfrak{g}, \star_2}(\hat{M}, \mathcal{C}^{\la}) \\
 &=&  \RHom_{\mathfrak{g}, \star_2}(E , \hat{M}^\vee \otimes\mathcal{C}^{\la}).
 \end{eqnarray*}
 \end{lem}%
 \begin{proof} For the first equality, we use that  $\mathcal{C}^{\la}$
   is a $\hat{U}(\mathfrak{g})$-module so that $\RHom_{U(\mathfrak{g}), \star_2}(M, \mathcal{C}^{\la}) =  \RHom_{\hat{U}(\mathfrak{g}), \star_2}({M}\otimes^L_{U(\mathfrak{g})} \hat{U}(\mathfrak{g}), \mathcal{C}^{\la})$.  It follows from Corollary \ref{cor-exact-U(g)-to-hat} that ${M}\otimes^L_{U(\mathfrak{g})} \hat{U}(\mathfrak{g}) = \hat{M}[0]$.

 For the second equality, we  have an obvious map 
$$  \RHom_{\mathfrak{g}, \star_2}(E , \hat{M}^\vee
\otimes\mathcal{C}^{\la}) \rightarrow  \RHom_{\mathfrak{g}, \star_2}(
\hat{M},\mathcal{C}^{\la})=\RHom_{\mathfrak{g}, \star_2}(M,
\mathcal{C}^{\la}). $$ By resolving $M$ by free modules it suffices to check that this map is a quasi-isomorphism for $M=U(\mathfrak{g})$.  But then we have
$$\RHom_{\mathfrak{g},\star_2}(U(\mathfrak{g}),\mathcal{C}^{\la})=\mathcal{C}^{\la}[0]=\RHom_{\mathfrak{g},\star_2}(E,\O_{G,e}\otimes\mathcal{C}^{\la})$$
where the first equality is obvious. %
For the other equality, we think of  $ \oscr_{G,e} \otimes \mathcal{C}^{\la} $ as a submodule of  $\oscr_{G,e} \otimes \oscr_{G,e} \otimes \oscr_{\mathcal{FL}}$ which is the   germs of functions $f( g',g,x)$ at $(e,e)$ in $G \times G\times \mathcal{FL}$. The  $\mathfrak{g}$-action is induced from $g''.f(g',g,x) = \lambda^{-1}(b) f( (g'')^{-1}g', gg'',x)$ for $g'' \in G_e$. We consider the automorphism of $ \oscr_{G,e} \otimes \mathcal{C}^{\la}$ given by the map $f(g',g,x)
\mapsto [(g',g) \mapsto f( g', gg',x)]$. Via this automorphism, the
$\mathfrak{g}$-action becomes $g''.f(g,g',x) = f((g'')^{-1}g', g,x)$ for $g''
\in G_e$, and is therefore only on the first factor. We can now  use the
flatness of $\mathcal{C}^{\la}$ over $E$ and Lemma \ref{lem:OGe-is-acyclic} to
conclude.
\end{proof}

We recall that  $Z(\mathfrak{g})$ lies in the centre of
$D^{-}(U(\mathfrak{g}))$. We also have defined a map $\Theta_{\hor} : Z(\mathfrak{m})  \rightarrow \mathrm{End}_{\oscr_{\mathcal{FL}}}(\mathcal{C}^{\la})$ in section \ref{section-def-cla}.

\begin{lem}\label{lem-homogeneous-Loc}  For any $z \in Z(\mathfrak{g})$, we have 
$\mathrm{Loc}(z) = \Theta_{\hor}(HC( \iota z))$. %
\end{lem}
\begin{proof} This follows from Lemma \ref{lem-all-compatibilities}. 
\end{proof}

\begin{cor}%
 Let $M \in \mathrm{Mod}(U(\mathfrak{g}))$ be a module with
 infinitesimal character $\lambda \in X^\star(T)_E$ (modulo dotted
 $W$-action). Then on $\mathrm{Loc}(M)$, the horizontal action of
 $Z(\mathfrak{m})$ factors through an action of $Z(\mathfrak{m})
 \otimes_{ HC, Z(\mathfrak{g})} (-w_0\lambda)$. %
\end{cor}
\begin{proof} We recall that for the map $HC_{\mathfrak{g}}: Z(\mathfrak{g}) \rightarrow U(\mathfrak{h})$, we have $HC_{\mathfrak{g}} \circ \iota = -w_0 HC_{\mathfrak{g}}$. The rest follows from Lemma \ref{lem-homogeneous-Loc}. 
\end{proof}

\begin{rem} Let us describe $\Spec~Z(\mathfrak{m}) \otimes_{ HC,
    Z(\mathfrak{g})} (-w_0\lambda)$, %
  or equivalently the idempotents in this finite $E$-algebra (note
  however that in singular weight %
  $Z(\mathfrak{m}) \otimes_{ HC, Z(\mathfrak{g})} (-w_0\lambda)$ is
  not reduced).  The possible characters of $Z(\mathfrak{m}) \otimes_{
    HC, Z(\mathfrak{g})} (-w_0\lambda)$  range through the set $\{
  w\cdot (-w_0\lambda), w \in \WM\}$. 
Since $w\cdot (-w_0\lambda) = -w_{0,M}( w_{0,M}w w_0\cdot\lambda + 2 \rho^M)$
we deduce that %
 $$ \Spec~Z(\mathfrak{m}) \otimes_{ HC, Z(\mathfrak{g})} (-w_0\lambda)
 = \{-w_{0,M}(w\cdot\lambda + 2 \rho^M), w \in \WM\}.$$ %
It follows  that if~
$M$ has infinitesimal character~ $\lambda$, then $$\mathrm{Loc}(M) = \oplus_{  w \in \WM}  \mathrm{Loc}(M)_{-w_{0,M}(w\cdot\lambda + 2 \rho^M)}$$ where $\mathrm{Loc}(M)_{-w_{0,M}(w\cdot\lambda + 2 \rho^M)}$ is the direct factor which corresponds to the idempotent in $Z(\mathfrak{m}) \otimes_{ HC, Z(\mathfrak{g})} (-w_0\lambda)$ given by $-w_{0,M}(w\cdot\lambda + 2 \rho^M)$. 
\end{rem}

 We conclude our generalities on our localization problem by showing
 that it is pre-dual
 to an obvious variant of the classical localization problem as in
 \cite{MR733805} (which is of course formulated in the algebraic
 context, and involves a fixed choice of (generalized) infinitesimal
 character). In order to do so, we recall the Chevalley--Eilenberg resolution of $E$ (with $d= \dim \mathfrak{g}$):
$$ 0 \rightarrow U(\mathfrak{g}) \otimes \Lambda^d\mathfrak{g} \rightarrow \cdots \rightarrow U(\mathfrak{g}) \rightarrow E \rightarrow 0.$$ 

By Lemma \ref{lem-passingtotheotherside}, for a finitely generated $U(\mathfrak{g})$-module $M$,
$\Loc(M)=\RHom_{\mathfrak{g}, \star_2}(M, \mathcal{C}^{\la})$ is computed by the following complex  of $\LB$-sheaves (in degree $[0,d]$):
\numequation\label{eqn:complex-of-LB-sheaves-for-loc}0 \rightarrow \mathcal{C}^{\la} \otimes \hat{M}^\vee \rightarrow \cdots \rightarrow \mathcal{C}^{\la} \otimes \hat{M}^\vee \otimes \Lambda^{d} \mathfrak{g}^\vee \rightarrow 0.\end{equation}

\begin{prop}%
  For $M \in
  \mathrm{Mod}^{\fing}(U(\mathfrak{g}))$, we have $\underline{\mathrm{RHom}}_{\oscr_{\mathcal{FL}}}(\RHom_{\mathfrak{g}, \star_2}(M, \mathcal{C}^{\la}), \oscr_{\mathcal{FL}}) = \tilde{\mathcal{D}}^{\la} \otimes^L_{{U}(\mathfrak{g})} {M}$. 
\end{prop}
\begin{proof}The same computation as in Proposition \ref{prop-dualCla} shows that   the derived $\oscr_{\mathcal{FL}}$-dual of the complex:
$$0 \rightarrow \mathcal{C}^{\la} \otimes \hat{M}^\vee \rightarrow \cdots
\rightarrow \mathcal{C}^{\la} \otimes \hat{M}^\vee \otimes \Lambda^{d}
\mathfrak{g}^\vee \rightarrow 0$$is the complex:
\[0 \rightarrow \tilde{\mathcal{D}}^{\la} \otimes \hat{M} \otimes \Lambda^{d} \mathfrak{g} \rightarrow \cdots \rightarrow \tilde{\mathcal{D}}^{\la} \otimes \hat{M} \rightarrow 0.\qedhere\]
\end{proof}

\subsubsection{$B$-action} We can also exploit the $B$-equivariant
structure on $\mathrm{Loc}(M)$. To this end, suppose that $M \in
\ocal(\mathfrak{g}, \mathfrak{b})_{\lambda\dalg}$, so that~$M(-\lambda)$ has an
action of~$B$ (see Remark~\ref{rem:M-twisted-minus-lambda}). Then by Proposition~\ref{prop:completion-functor-on-cat-O}, $\hat{M}^\vee(\lambda) \in \LB_{(\mathfrak{g},B)}(E)$. We deduce that the complex~\eqref{eqn:complex-of-LB-sheaves-for-loc} computing $\Loc((M(-\lambda))=\RHom_{\mathfrak{g}, \star_2}(M(-\lambda), \mathcal{C}^{\la})$ is a complex in $\LB_{(\mathfrak{g},B)}(\mathcal{FL})$. More precisely, $\mathcal{C}^{\la} \otimes \hat{M}^\vee \otimes \Lambda^i \mathfrak{g}^\vee$ carries  the induced $\mathfrak{g}$-action from the $*_{1,3}$-action on $\mathcal{C}^{\la}$ and the $B$-action which is the tensor product of the $B$-action on $\mathcal{C}^{\la}$ and the $B$-action on $\hat{M}^\vee \otimes \Lambda^i \mathfrak{g}^\vee$. There is another $\mathfrak{g}$-action which is the tensor product of the $\star_2$ action on $\mathcal{C}^{\la}$ and the $\mathfrak{g}$-action on $\hat{M}$ (and which is used to construct the differentials in the complex).

\subsubsection{Main theorem} For $M \in 
\ocal(\mathfrak{g},\mathfrak{b})_{\lambda\dalg}$, the
cohomology sheaves of $\Loc(M(-\lambda))$ are $(\mathfrak{g},
B)$-equivariant sheaves that we want to describe.  As a first step we
intend to  describe their restrictions to~ $C_{w}^\dag$ for each  $w  \in \WM$. %
\begin{rem} In principle the cohomology sheaves could be sheaves of
  solid $E$-vector spaces which need not arise from nice sheaves of
  topological spaces (in more classical language, the cohomology could
  be non-separated). However, under the assumption that $\lambda$ is
  non-Liouville, we  see as a consequence of the following theorem that they are actually separated objects and again belong to the category $\LB_{(\mathfrak{g},B)}(C_w^\dag)$. 
\end{rem}
\begin{thm}\label{thm-localization} Let $M \in \ocal(\mf{g},\mf{b})_{\lambda-\alg}$ and assume that $\lambda$ is non-Liouville. 
Then we have: $$\underline{\HH}^i(\mathrm{Loc}(M(-\lambda)))|_{C_{w}^\dag} = HCS_{w,\lambda}( \HH_i(\mathfrak{u}_{P_w}, M)).$$
\end{thm}

\begin{proof}%

By the definition of the functor~$HCS$ the sheaf $HCS_{w,\lambda}(
\HH_i(\mathfrak{u}_{P_w}, M))$ corresponds via the equivalence of categories of  Proposition
~\ref{prop-gQeq} to the
$\Stab(w)$-representation $\widehat{\HH_i(\mathfrak{u}_{\mathfrak{p}}, M)}^\vee(\lambda)$.
  By
 Corollary \ref{Main-cor-strictness}, we can identify this with
 $\mathrm{Ext}^i_{\mathfrak{u}_{P_w}}(E, \hat{M}^\vee (\lambda) )
 $, and by Proposition~\ref{prop-qis-g-coh} below, this can in turn be
 identified with \numequation\label{eqn:HCS-fibre-Ext-description}\Ext^i_{\mathfrak{g}, \star_2}(\hat{M}(-\lambda),
  \oscr_{U_{P_w} \backslash G,
    e}).\end{equation}
 By definition we
  have \numequation\label{eqn:restriction-to-Cw-is-Ext}\underline{\HH}^i(\mathrm{Loc}(M(-\lambda)))|_{C_{w}^\dag}=\Ext^i_{\mathfrak{g},
      \star_2}(M(- \lambda),
    \mathcal{C}^{\la}\vert_{C_{w}^\dag}).\end{equation}
  Morally, it remains to show that passage to the fiber at~$w$
  identifies the right hand side
  of~\eqref{eqn:restriction-to-Cw-is-Ext}
  with~\ref{eqn:HCS-fibre-Ext-description}. However, we have to be a
  little careful with this comparison, because we have not developed a
  theory which allows us to consider arbitrary sheaves of solid
  $E$-vector spaces.

To this end, we  consider the
 Chevalley--Eilenberg complex computing $\Ext^i_{\mathfrak{g},
      \star_2}(M(- \lambda), \mathcal{C}^{\la}\vert_{C_{w}^\dag})$. Under the equivalence of categories of  Proposition
~\ref{prop-gQeq} (given by taking the fiber at $w$), %
this complex corresponds to the following complex of $\Stab(w)$-representations, 
$$0 \rightarrow  \oscr_{U_{P_w} \backslash G, e} \otimes \hat{M}^\vee (\lambda) \rightarrow  \oscr_{U_{P_w} \backslash G, e} \otimes \hat{M}^\vee(\lambda) \otimes \mathfrak{g}^\vee \rightarrow  \oscr_{U_{P_w} \backslash G, e} \otimes \hat{M}^\vee (\lambda) \otimes  \Lambda^2\mathfrak{g}^\vee  \rightarrow \dots $$
which  computes $\RHom_{\mathfrak{g}, \star_2}(\hat{M}(-\lambda),
\oscr_{U_{P_w} \backslash G, e})$, as required.

(This cohomology is computed in the category of solid $E$-vector
 spaces. Again, the cohomology groups could be very pathological (from
 the classical perspective). For clarity, we can make explicit   the action of $\Stab(w)$ on $\oscr_{U_{P_w} \backslash G, e} \otimes \hat{M}^\vee (\lambda) \otimes  \Lambda^i\mathfrak{g}^\vee $. This action consists of:
\begin{itemize}
\item An action of $(P_w)_e$   induced by the action on  $\oscr_{U_{P_w} \backslash G, e}$ via $p.f(-) = f(p^{-1}-)$. 
\item An action of $B_e$, which is the tensor product of the action on $\oscr_{U_{P_w} \backslash G, e}$ via $b.f(-) = f(-b)$ and of the restriction to $B_e$ of the $B$-action on $\hat{M}^\vee (\lambda) \otimes  \Lambda^i\mathfrak{g}^\vee$.
\item An action of $B\cap P_w$ which is the tensor product of the action on $\oscr_{U_{P_w} \backslash G, e}$ via $b.f(-) = f(b^{-1}-b)$ and the action of $B\cap P_w$ on $\hat{M}^\vee (\lambda) \otimes  \Lambda^i\mathfrak{g}^\vee$.
\end{itemize}

The differentials in the complex involve the $\mathfrak{g}$-action which is the tensor product of the $\star_2$ action on $\oscr_{U_{P_w} \backslash G, e}$ and the $\mathfrak{g}$-action on $\hat{M}^\vee(\lambda) \otimes \mathfrak{g}^\vee$.)
\end{proof}

\begin{prop}
  \label{prop-qis-g-coh}If $M \in 
\ocal(\mathfrak{g},\mathfrak{b})_{\lambda\dalg}$, then we have a
  quasi-isomorphism $\RHom_{\mathfrak{g}, \star_2}({M}(-\lambda),
  \oscr_{U_{P_w} \backslash G, e}) = \RHom_{\mathfrak{u}_{P_w}}(E,
  \hat{M}^\vee (\lambda) )$.
\end{prop}
\begin{proof}
  Indeed we have \begin{eqnarray*}
                   \RHom_{\mathfrak{g}, \star_2}({M}(-\lambda),  \oscr_{U_{P_w} \backslash G, e}) &= & \RHom_{\mathfrak{g}, \star_2}(E, \hat{M}^\vee(\lambda) \otimes  \oscr_{U_{P_w} \backslash G, e}) \\
                                                                                                  &=&  \RHom_{\mathfrak{g} \oplus \mathfrak{u}_{P_w}}(E, \hat{M}^\vee(\lambda) \otimes  \oscr_{ G, e}) \\
                                                                                                  &=&
                                                                                                      \RHom_{\mathfrak{u}_{P_w}}(E,
                                                                                                      \hat{M}^\vee
                                                                                                      (\lambda)
                                                                                                      ).
                 \end{eqnarray*}
                 Here the first equality follows from the same
                 argument as in Lemma \ref{lem-passingtotheotherside},
                 the second equality uses that
                 $\RHom_{\mathfrak{u}_{P_w}}(E, \oscr_{ G, e}) =
                 \oscr_{U_{P_w} \backslash G, e}[0]$ and the flatness
                 of $\hat{M}^\vee(\lambda)$ (as it is a colimit of
                 Smith spaces), and the last equality follows from Lemma~\ref{lem:OGe-is-acyclic}.
               \end{proof}
\subsubsection{Localization of finite dimensional representations}
Let $\lambda \in X^\star(T)^{+}$. Let $V_\lambda$ be the irreducible finite dimensional representation of $G$ of highest weight $\lambda$ viewed as an object of $\mathcal{O}(\mathfrak{g}, \mathfrak{b})_{0-alg}$.
Let $\lambda \in X^\star(T)^{+,M}$. We let $L_\lambda$ be the irreducible finite dimensional representation of $M$ of highest weight $\lambda$. We also let $d = \mathrm{dim} (\mathfrak{u}_P)$. We recall the following theorem of Kostant:
\begin{thm}\label{thm: of Kostant} We have that $\HH_i(\mathfrak{u}_{P},V_{\lambda}) = \oplus_{w \in \WM, \ell(w)=d-i} L_{w\cdot\lambda + 2\rho^M}$. 
\end{thm} %
\begin{proof}
  See for example \cite[Thm.\ 4.2.1]{MR2513265} (together
  with~\eqref{eqn:Poincare-duality-Lie-algebra} to pass from
  cohomology to homology).%
\end{proof}

\begin{prop}\label{prop-loc-fd}  We have $\mathrm{Loc}(V_{\lambda}) =  \bigoplus_{w \in \WM} \mathrm{Loc}(V_{\lambda})_{-w_{0,M}(w\cdot\lambda + 2\rho^M)}$ and 
$ \mathrm{Loc}(V_{\lambda})_{-w_{0,M}(w\cdot\lambda + 2\rho^M)} = \mathcal{L}_{-w_{0,M}(w\cdot\lambda + 2\rho^M)}[-d+\ell(w)]$.
\end{prop}
\begin{proof} Since there is a $G$-action on $V_{\lambda}$, we see
  that $\mathrm{Loc}(V_{\lambda})$ is in fact computed  by a complex
  in  $\LB_{(\mathfrak{g},G)}(\mathcal{FL})$ (and not just in
  $\LB_{(\mathfrak{g},B)}(\mathcal{FL})$).  Via the equivalence of
  categories of Proposition \ref{prop: LB g G equiv of cats}, this
  complex corresponds to the complex $\RHom_{\mathfrak{g},
    \star_2}(V_\lambda,  \oscr_{U_{P_w} \backslash G, e})$ in
  $L_{(\mathfrak{g}, P)}(E)$,   and as in Proposition \ref{prop-qis-g-coh},
we find that this is quasi-isomorphic to
$\RHom_{\mathfrak{u}_{P_w}}(E, {V_\lambda}^\vee )$. It follows that the cohomology groups
are  simply  given by the representations of $M$:
$\Ext^i_{\mathfrak{u}_{P_w}}(E, {V_\lambda}^\vee ) = \HH_i(\mathfrak{u}_P,
V_\lambda)^\vee$. By Kostant's Theorem~\ref{thm: of Kostant}, these correspond
to $\oplus_{w, \ell(w) = d-i}\mathcal{L}_{-w_{0,M}(w\cdot\lambda +
  2\rho^M)}$, as required.
\end{proof}

\begin{rem} Proposition~\ref{prop-loc-fd} is clearly compatible with
  Theorem~\ref{thm-localization}, since %
\begin{eqnarray*}
\mathrm{Loc}(V_{\lambda}) \vert_{C_w^\dag} & =& \oplus_{w' \in \WM} \mathcal{L}_{-w_{0,M}(w'\cdot\lambda + 2\rho^M)}[\ell(w')-d]\vert_{C_w^\dag} \\
&=& \oplus_{i=0}^{d} HCS_{w,\lambda}( \HH_i(\mathfrak{u}_{P_w}, V_\lambda))
\end{eqnarray*}
but it also gives more information as it describes all the extensions between the sheaves on the Bruhat strata. 
\end{rem}

\subsubsection{Localization of Verma modules in the non-Liouville
  case} %
Let $\lambda \in X^\star(T)_{E}$ be non-Liouville. We %
view the Verma module~$M_{\lambda}$ of weight $\lambda$ as an object of
$\ocal(\mf{g},\mf{b})_{\lambda-\alg}$. Thanks to Theorem \ref{thm-localization},
we see that understanding $\Loc(M_\lambda(-\lambda))\vert_{C_w^\dag}$ boils down
to understanding the cohomology of some Verma modules, as in
Theorem~\ref{thm-algebraic-computationStrict}.  We also  assume for simplicity
that~$\mf{u}_P$ is abelian (this assumption holds in our applications to Shimura varieties).

\begin{cor}\label{coro-ESflag} Assume that $\lambda \in X^\star(T)_E$
  is non-Liouville and that $\mathfrak{u}_{P_w}$ is abelian. %
  Then the
  following hold: %
\begin{enumerate}
\item All the cohomology groups $\underline{\HH}^i(\mathrm{Loc}(M_\lambda(-\lambda)))|_{C_{w}^\dag}$ belong to the image of the functor $HCS_{w,\lambda}: \mathcal{O}(\mf{m}_w,\mf{b}_{M_w})_{\lambda\dalg} \rightarrow \LB_{(\mathfrak{g},B)}( C_w^\dag)$.
\item The cohomology groups are zero if $i > d-\ell(w)$. 
\item There is a surjective ``highest weight'' map: \[
  \underline{\HH}^{d-\ell(w)}(\mathrm{Loc}(M_\lambda(-\lambda)))|_{C_{w}^\dag} \rightarrow HCS_{w,
    \lambda}(\Mmw_{\lambda + w^{-1}w_{0,M}\rho +\rho}).\]%
\item The kernel of the highest weight map, and the cohomology groups
  $\underline{\HH}^{i}(\mathrm{Loc}(M_\lambda(-\lambda)))|_{C_{w}^\dag}$ for $i < d-\ell(w)$, admit
  finite filtrations with sub-quotients ranging among the sheaves
  $\mathcal{L}_{w, -w_{0,M}(w'\cdot\lambda + 2\rho^M)} \otimes E(
  \lambda- w^{-1}w'\cdot\lambda - w^{-1}w_{0,M} \rho- \rho)$, where
  $w' \in wW_{< \lambda}$. %
\end{enumerate}

\end{cor}%
\begin{proof} This is immediate from Theorem \ref{thm-localization}
  and Theorem \ref{thm-algebraic-computationStrict}, bearing in
  mind~\eqref{eqn:two-ways-to-write-the-weights} and Proposition~\ref{TwistingHCS}. %
\end{proof}

\begin{rem}\label{rem-antidom} In particular, if $\lambda$ is non-Liouville and
  antidominant 
  (i.e.\ that $W_{<\lambda} = \emptyset$), we see that the cohomology is
  concentrated in degree $d-\ell(w)$ and that the highest weight map
  is an isomorphism on this cohomology. %
\end{rem}

\subsubsection{Localization of Verma modules in general}
Let $\lambda \in X^\star(T)_{E}$. %
 For sake of completeness, in this section we give the  following general result (without any non-Liouville
assumption on $\lambda$) which is a weaker form of Corollary
\ref{coro-ESflag}. We still  assume for simplicity   that~$\mf{u}_P$ is
abelian. This result was obtained by Juan Esteban Rodriguez-Camargo in his PhD thesis. 

\begin{thm}[Rodriguez-Camargo]\label{thm-general-coho-verma}%
  Let $M_{\lambda} \in \ocal(\mf{g},\mf{b})_{\lambda-\alg}$ be the
  Verma module of weight $\lambda$. %
Let $w \in \WM$. The following is true:
\begin{enumerate}
\item $\underline{\HH}^{i}(\mathrm{Loc}(M_\lambda(-\lambda)))|_{C_{w}^\dag}$  vanishes unless $i \in [ 0, d-\ell(w)]$. 
\item We have a surjective ``highest weight''  map \[
 \underline{\HH}^{d-\ell(w)}(\mathrm{Loc}(M_\lambda(-\lambda)))|_{C_{w}^\dag} \rightarrow HCS_{w, \lambda}(\Mmw_{\lambda + w^{-1}w_{0,M}\rho +\rho}).\] %
\item If $w=w_0^M$, the above map is an isomorphism. 
\end{enumerate}
\end{thm}

\begin{proof} 
We remark that \numequation\label{eqn:Loc-of-Verma-is-b-Hom-lambda}\mathrm{Loc}(M_{\lambda}) = \RHom_{\mathfrak{g}}( M_{\lambda}, \mathcal{C}^{\la}) = \RHom_{\mathfrak{b}}(E(\lambda), \mathcal{C}^{\la}).\end{equation}
Thus, the fiber of $\mathrm{Loc}(M_\lambda)$ at $w$   is the
$Stab(w)$-representation $ \mathrm{Ext}_{\mathfrak{b} ,\star_2}
(E(\lambda), \oscr_{U_{P_w} \backslash G, e})$, and the result is
immediate from Proposition~\ref{prop:coho-of-Verma-for-Juan-theorem} below. 
 \end{proof}
 \begin{prop}\label{prop:coho-of-Verma-for-Juan-theorem}  The cohomology groups %
  $\mathrm{Ext}^i_{\mathfrak{b} ,\star_2} (E(\lambda), \oscr_{U_{P_w} \backslash G, e})$ vanish outside degrees $[0,d-\ell(w)]$. Moreover, there is a canonical surjective map $$\mathrm{Ext}^{d-\ell(w)}_{\mathfrak{b} ,\star_2} (E(\lambda), \oscr_{U_{P_w} \backslash G, e}) \rightarrow %
    \Mmw_{\lambda + w^{-1}w_{0,M}\rho +\rho}^\vee.$$  If $w=w_0^M$, this map is an isomorphism. 
\end{prop}
 \begin{proof} We recall that by Proposition \ref{prop-qis-g-coh} (and
   its proof), we have: 
 \begin{eqnarray*}
 \RHom_{\mathfrak{b} ,\star_2} (E(\lambda), \oscr_{U_{P_w} \backslash G, e}) &= &\RHom_{ \mathfrak{u}_{P_w} \oplus \mathfrak{b} }(E, \oscr_{G,e}(-\lambda)) \\ &=& \RHom_{ \mathfrak{u}_{P_w}  }(E, \hat{M}^\vee_\lambda).
 \end{eqnarray*}Now, the cohomology $\mathrm{R}\Gamma(\mathfrak{u}_{P_w} \cap
\bar{\mathfrak{b}} , \hat{M}_\lambda^\vee )$ is concentrated in degree
$0$. %
(Indeed, recall from
Remark~\ref{rem-making-explicit} that $\hat{M}_\lambda^\vee =
(\oscr_{  G/U, e}(-\lambda))^{\mathfrak{h}}$. As a  $\mathfrak{u}_{P_w}
\cap \bar{\mathfrak{b}}$-module, this module can be written in the
form $\oscr_{U_{P_w} \cap \bar{B},e} \otimes V$ where $V$ is an
$\LB$-space of compact type with trivial action. We observe that
$\HH^i(\mathfrak{u}_{P_w} \cap \bar{\mathfrak{b}}, \oscr_{U_{P_w} \cap
  \bar{B},e}) =0 $  if $i>0$.) %

 We therefore have:
\begin{eqnarray*}
\mathrm{R}\Gamma( \mathfrak{u}_{P_w} ,  \hat{M}_\lambda^\vee) & =&  \mathrm{R}\Gamma( \mathfrak{u}_{P_w} \cap \mathfrak{b}, \mathrm{R}\Gamma(\mathfrak{u}_{P_w} \cap \bar{\mathfrak{b}} , \hat{M}_\lambda^\vee ))\\
&=& \mathrm{R}\Gamma( \mathfrak{u}_{P_w} \cap \mathfrak{b}, \HH^0(\mathfrak{u}_{P_w}\cap\bar{\mathfrak{b}}, \hat{M}_\lambda^\vee)),
\end{eqnarray*}and we see in particular that the cohomology vanishes
above degree $d-\ell(w)=\dim \mf{u}_{P_w}\cap\mf{b}$.

We now consider the surjective ``restriction''  map
  $\oscr_{G,e} \rightarrow \oscr_{P_w(B\cap \bar{U}_{P_w}),e}$ induced
  by the inclusion $P_w(B\cap \bar{U}_{P_w}) \hookrightarrow
  G$.
We deduce a map:
$$\RHom_{ \mathfrak{u}_{P_w} \oplus \mathfrak{b} }(E(\lambda), \oscr_{G,e}) \rightarrow \RHom_{ \mathfrak{u}_{P_w} \oplus \mathfrak{b} }(E(\lambda), \oscr_{P_w(B\cap \bar{U}_{P_w}),e}).$$
We claim that this this map is surjective in degree $d-\ell(w)$, and
that it is an isomorphism if~$w=w_0 ^M$. 
To see this, we first take $\mathfrak{b}$ and $\mathfrak{u}_{P_w}\cap\mathfrak{b}$-cohomology 
which gives a  surjective map (the cohomology is still in degree $0$): 
$$(\oscr_{ U_{P_w} \cap \bar{B}\backslash G/B,e}(-\lambda))^{\mathfrak{h}} \rightarrow (\oscr_{U_{P_w \cap \bar{B}}\backslash P_w/B\cap P_w,e}(-\lambda))^{\mathfrak{h}}.$$
Taking the cohomology of $\mathfrak{u}_{P_w} \cap
{\mathfrak{b}}$,  we see that the above surjective map induces a
surjective map in top degree cohomology (which is an isomorphism  if~$w=w_0 ^M$). 

It remains to identify the target of the surjection
with~$\Mmw_{\lambda + w^{-1}w_{0,M}\rho +\rho}^\vee$. We cannot
immediately deduce this, because the
computation above is not $\mathfrak{m}_w$-equivariant (because the
decomposition of $\mathfrak{u}_{P_w} = \mathfrak{u}_{P_w} \cap
\mathfrak{b} \oplus \mathfrak{u}_{P_w} \cap \bar{\mathfrak{b}}$ is not
$\mathfrak{m}_w$-invariant). In order to identify the
$\mathfrak{m}_w$-module structure of $\RHom_{
  \mathfrak{u}_{P_w} \oplus \mathfrak{b} }(E(\lambda),
\oscr_{P_w(B\cap \bar{U}_{P_w}),e})$ we compute the cohomology in a
different way, by first considering $\mathfrak{u}_{P_w}$-cohomology,
and then $\mathfrak{b}$-cohomology.

Certainly
\[\RHom_{
  \mathfrak{u}_{P_w} }(E,
\oscr_{P_w(B\cap \bar{U}_{P_w}),e})=\oscr_{ U_{P_w}  \backslash
  P_w(B\cap \bar{U}_{P_w}), e}[0],\]so it remains to show that  $$\mathrm{Ext}^{d-\ell(w)}_{\mathfrak{b},\star_2}(E(\lambda),  \oscr_{ U_{P_w}  \backslash P_w(B\cap \bar{U}_{P_w}), e}) =  (\oscr_{M_w/U_{M_w},e}(-\lambda - w^{-1}w_{0,M} \rho - \rho))^{\mathfrak{h}}.$$ 
 We will then be done, because the right hand side is %
  $\hat{M}(\mathfrak{m}_w)^\vee_{\lambda + w^{-1}w_{0,M}\rho + \rho}$ by Remark \ref{rem-making-explicit} and the proof of Proposition \ref{prop-completion-Stab}. %

We now show the claim.   Note that %
 $\mathfrak{b} = \mathfrak{b} \cap \mathfrak{p}_w  \oplus \mathfrak{b} \cap \bar{\mathfrak{u}}_{P_w}$. 
 We first compute %
$$\mathrm{R}\Gamma_{\star_2}( \mathfrak{b} \cap
\bar{\mathfrak{u}}_{P_w},  \oscr_{ U_{P_w}  \backslash P_w(B\cap
  \bar{U}_{P_w}), e}) = \oscr_{ U_{P_w}  \backslash P_w, e}[0] =
\oscr_{M_w,e}[0].$$We now observe that $\oscr_{M_w,e}$ is a
$\mathfrak{b} \cap \mathfrak{p}_w$-representation, with $\mathfrak{b}
\cap \mathfrak{u}_{P_w}$ acting trivially. 
 We deduce that 
 \begin{align*}
 \mathrm{Ext}^{d-\ell(w)}_{ \mathfrak{b} \cap {\mathfrak{p}}_w, \star_2}( E(\lambda),  \oscr_{ M_w, e}) & =  \Hom_{ \mathfrak{b} \cap {\mathfrak{m}}_{w}, \star_2}( E(\lambda),  \oscr_{M_w, e} \otimes \wedge^{d-\ell(w)} (\mathfrak{u}_{P_w} \cap \mathfrak{b})^\vee) \\
 &= (\oscr_{M_w/U_{M_w},e}(-\lambda - w^{-1}w_{0,M} \rho - \rho))^{\mathfrak{h}}. \qedhere
 \end{align*}%
 \end{proof}

\subsection{Localization and higher Coleman sheaves at singular weight}\label{subsec:localization-Higher-Coleman-GSp4} %
In this section we study localization
at a singular weight for $G = \mathrm{GSp}_4$ and  $P$ is the Siegel parabolic
associated to the cocharacter $\mu =  (-1/2,-1/2;1/2) \in
X_{\star}(T)_E$. %
We freely use our notation for $\mathrm{GSp}_4$ (see \ref{notn:GSp4}). 
We consider the Klingen parabolic $Q \supseteq B $ attached to the  simple root $\beta$. We denote by $M_Q$ the associated Levi which is a group of semi-simple rank $1$. 
It is important for us that $w_{0}^M s_\beta \in \WM$; this implies that the stratum $C_{w_0^M,Q}$ is the union of two $B$-orbits  $C_{w_0^M}$ and $C_{w_0^Ms_\beta}$.
We wish to study the localization $\Loc(M(\mathfrak{g})_\lambda(-\lambda))\vert_{C_{w_0^M,Q}}$ for $\lambda = (1,1;w)$. We notice that in this singular weight the horizontal action is not semi-simple. We are going to study this  action and describe the semi-simple part.

\subsubsection{Geometry of the  strata}\label{sec-actionB}

We consider the $Q$-orbit $C_{w_0^M,Q}$ on~$\mathcal{FL}$, which is
the union of the two Bruhat strata $C_{w_0^M}$
and~$C_{w_0^Ms_{\beta}}$. 

Note that the map $Q \rightarrow C_{w_0^M,Q}$, $q \mapsto w_0^M q$ induces an isomorphism $P_{w_0^M} \cap Q \backslash Q \rightarrow C_{w_0^M,Q}$. 
The projection $Q \rightarrow M_Q$ induces a map $(P_{w_0^M} \cap Q) \backslash Q \rightarrow (P_{w_0^M} \cap M_Q) \backslash M_Q$. Since $M_Q$ has semi-simple rank one and $(P_{w_0^M} \cap M_Q)$ is a Borel subgroup, we can identify $(P_{w_0^M} \cap M_Q) \backslash M_Q$ with $\mathbb{P}^1$, with $\infty$ the image of $w_{0}^Ms_\beta$ and $0$ the image of $w_0^M$. 
We therefore have a  natural
map $\pi:C_{w_0^M,Q}\to\mathbf{P}^1$, with $\pi^{-1}(\{\infty\})=
C_{w_0^Ms_{\beta}}$ and $\pi^{-1}(\mathbb{A}^1) = C_{w_0^M}$; moreover, $B
\cap s_{\beta} B s_{\beta}^{-1}$ acts transitively on  $\pi^{-1}(\mathbb{G}_m) =
C_{w_0^M} \cap   C_{w_0^M}s_\beta$.

\subsubsection{Sheaves on the union of strata}

  The  inclusion  $G_e\rtimes B\into G_e\rtimes Q$ induces an
  inclusion $\Stab_B(w_0^M) \into \StabQ(w_0^M)$, which induces the
  restriction map $\LB_{(\mathfrak{g}, Q)}(C^{\dagger}_{w_0^M,Q}) \rightarrow
  \LB_{(\mathfrak{g}, B)}(C^{\dagger}_{w_0^M}) $. %
 Similarly, conjugation by $s_{\beta}$ gives an inclusion $\Stab_{B}(w_0^Ms_{\beta}) \into
\Stab_{Q}(w_0^M)$, %
which %
 induces the restriction map $\LB_{(\mathfrak{g}, Q)}(C^{\dagger}_{w_0^M,Q}) \rightarrow
  \LB_{(\mathfrak{g}, B)}(C^{\dagger}_{w_0^Ms_{\beta}}) $.

We also note that the image of  $Q \cap P_{w_{0}^M}$ in $M_{w_0^M}$ is the Borel
$B_{M_{w_0^M}}$ (and the image of  $Q \cap P_{w_{0}^Ms_\beta}$ in
$M_{w_0^Ms_\beta}$ is the Borel $B_{M_{w_0^Ms_\beta}}$). Therefore, the source
category  $\ocal(\mathfrak{m}_{w_0^M}, \mathfrak{q}_{M_{w_0^M}})$ for producing
sheaves on  $C_{w_0^M,Q}$ (see Definition~ \ref{defn:HCS-functor}) is
tautologically equal to $\ocal(\mathfrak{m}_{w_0^M}, \mathfrak{b}_{M_{w_0^M}})$. 

 \begin{prop}\label{lem-extendingunionstrata} Let $\nu \in
   X^{\star}(M_Q)_E$. %
 \begin{enumerate}
 \item Conjugation by $s_\beta$ induces an equivalence of categories $s_{\beta}: \ocal(\mathfrak{m}_{w_0^M}, \mathfrak{b}_{M_{w_0^M}})_{\nu\dalg} \rightarrow \ocal(\mathfrak{m}_{w_0^Ms_{\beta}}, \mathfrak{b}_{M_{w_0^Ms_{\beta}}})_{\nu\dalg}$. 
 \item For any~$M\in\cO(\mf{m}_{w_0^M},\mf{b}_{M_{w_0^M}})_{\nu\dalg}$, we have
 \[HCS_{Q,w_0^M, \nu}(M)\vert_{C^{\dagger}_{w_0^M}} =
HCS_{w_0^M,\nu}(M),\] \[HCS_{Q,w_0^M, \nu}(M)\vert_{C^{\dagger}_{w_0^Ms_{\beta}}} = HCS_{w_0^Ms_{\beta},\nu}(s_\beta M).\]
 \end{enumerate}
 \end{prop}
 
 \begin{proof} The first point is obvious, and the second is immediate
   from the definition of the functors~$HCS$ (i.e.\ from the
   construction in the proof of Proposition~\ref{prop:full-faithful-Adm-to-LB-Stab}).
 \end{proof}

  \subsubsection{Singular localization}
Write $j_{w_0^M}:C^{\dagger}_{w_0^M}\into
C^{\dagger}_{w_0^M,Q}$,  $j_{w_0^Ms_{\beta}}:C^{\dagger}_{w_0^Ms_{\beta}}\into
C^{\dagger}_{w_0^M,Q}$ for the inclusions. 
We take $\lambda = (1,1;w) \in X^\star(T)$. Recall  that $\rho = (-1,-2;0)$, so
that $\lambda + \rho = (0,-1;w)$ %
is invariant under $s_\beta$. However,  this character is not integral. Let us define $\eta = (0,-1;1)$ which  differs from $\lambda + \rho$ by a character of the centre, and is still invariant under $s_\beta$. 

Applying Proposition \ref{lem-extendingunionstrata} with~$\nu=\eta$, we have a short exact sequence  of $(\mathfrak{g}, Q)$-equivariant sheaves: 
\numequation\label{eqn:HCS-ses}
0 \rightarrow (j_{w_0^M})_! HCS_{w_0^M, \eta}( M(\mathfrak{m}_{w_0^M})_\lambda) \rightarrow HCS_{Q,w_0^M, \eta}(  M(\mathfrak{m}_{w_0^M})_\lambda) \rightarrow HCS_{w_0^Ms_\beta , \eta}( M(\mathfrak{m}_{w_0^Ms_\beta})_{s_\beta\lambda}) \rightarrow 0 \end{equation} 

We want to study $\Loc(M(\mathfrak{g})_\lambda(-\lambda))\vert_{C^{\dagger}_{w_0^M,Q}} =  \RHom_{\mathfrak{b}, \star_2} (\lambda, \mathcal{C}^{\la}\vert_{C^{\dagger}_{w_0^M,Q}}) $. 
\begin{prop}\label{prop:we-have-exact-triangle}We have   an exact triangle: %
$$ (j_{w_0^M})_! HCS_{w_0^M, \lambda}( M(\mathfrak{m}_{w_0^M})_\lambda) \rightarrow Loc(M(\mathfrak{g})_\lambda(-\lambda))\vert_{C^{\dagger}_{w_0^M,Q}} \rightarrow HCS_{w_0^Ms_\beta , \lambda}( M(\mathfrak{m}_{w_0^Ms_\beta})_{s_\beta\lambda}) [-1] \stackrel{+1}\rightarrow $$ 
\end{prop}
\begin{proof} This is immediate from consequence of Theorem
  \ref{thm-localization} together with Theorem~\ref{thm-algebraic-computationStrict}
and Proposition~\ref{prop:intertwining-for-homology-w2}.
\end{proof}

There is a horizontal action~$\Theta_{\hor}$ of $Z({\mathfrak{m}})$ on
$\Loc(M_\lambda(-\lambda))\vert_{C^{\dagger}_{w_0^M,Q}}$. 
By Proposition~\ref{prop-action-centre-enveloping} (see also
Remark~\ref{rem:horizontal-on-lambda-normalization} below), this action is via $\mu_0 :=\langle
-w_0 \lambda, \mu \rangle$ on both $HCS_{w_0^M, \lambda}( M(\mathfrak{m}_{w_0^M})_\lambda)$
and $HCS_{w_0^Ms_\beta, \lambda}( M(\mathfrak{m}_{w_0^M})_{s_\beta\lambda})$ (and in particular  the action of $\mu$ doesn't split the triangle). Taking the derived invariants for  $\mu-\mu_0$ yields a triangle:

\numequation\label{eqn:triangle-derived-mu-invariants}
\begin{aligned}
    & (j_{w_0^M})_! HCS_{w_0^M, \lambda}( M(\mathfrak{m}_{w_0^M})_\lambda) \oplus (j_{w_0^M})_! HCS_{w_0^M, \lambda}( M(\mathfrak{m}_{w_0^M})_\lambda)[-1] \rightarrow \\
    & \RHom_{\mathfrak{b}, \star_2; \mu} ((\lambda,\mu_0), \mathcal{C}^{\la}\vert_{C^{\dagger}_{w_0^M,Q}}) \rightarrow \\
    & HCS_{w_0^Ms_\beta, \lambda}( M(\mathfrak{m}_{w_0^M})_{s_\beta\lambda})[-1] \oplus HCS_{w_0^Ms_\beta, \lambda}( M(\mathfrak{m}_{w_0^M})_{s_\beta\lambda})[-2] \stackrel{+1}\rightarrow
\end{aligned}
\end{equation}

Taking the $\HH^1$ yields the following short exact sequence of $(\mathfrak{g},B)$-equivariant sheaves:

\numequation\label{eqn:localization-as-extension-on-2-strata}0
\rightarrow  (j_{w_0^M})_!  HCS_{w_0^M, \lambda}( M(\mathfrak{m}_{w_0^M})_\lambda)
\rightarrow  \mathrm{Ext}^1_{\mathfrak{b}, \star_2; \mu}
((\lambda,\mu_0), \mathcal{C}^{\la}\vert_{C^{\dagger}_{w_0^M,Q}}) 
\rightarrow  HCS_{w_0^Ms_\beta, \lambda}( M(\mathfrak{m}_{w_0^M})_{s_\beta\lambda}) %
\rightarrow 0 \end{equation}

The main result of this section is  the following:

\begin{thm}\label{thm-extequal} The two extensions of  
 $(\mathfrak{g},B)$-equivariant sheaves~\eqref{eqn:localization-as-extension-on-2-strata} and \eqref{eqn:HCS-ses} differ by a twist of the $B$-action by the character $\lambda-\eta$ and multiplication by a  scalar in $E^\times$. 
  \end{thm}

\begin{rem} This arguments in the remainder of this section admit a simpler analogue in the
  $\mathrm{GL}_2$-context; see \cite[\S 6]{PilloniVB}.%
\end{rem} 
\subsubsection{Preparations for the proof of Theorem
  \ref{thm-extequal}} We have commuting actions of $\mathfrak{b}$ (via
$\star_2$) and $\mu$ (via the horizontal action) on
$\mathcal{C}^{\la}$.    We begin by isolating a certain sub-Lie
algebra of   $\mathfrak{b} \oplus  E\mu$ whose cohomology on
$\mathcal{C}^{\la}\vert_{C^{\dagger}_{w_0^M,Q}}$ is in degree $0$.

We use the usual standard basis elements~$X_{\beta},X_{-\beta},H_{\beta}$ with
~$H_{\beta}=[X_{\beta},X_{-\beta}]$. %
We set $\mf{h}':=\ker(\beta)$, so that $\mf{b}$-cohomology can be
obtained by first taking $\mf{h}'\oplus \mf{u}_Q$-cohomology, and then
taking $EX_{\beta}\oplus E H_{\beta}$-cohomology.  We write $\lambda' =
  \lambda\vert\mathfrak{h}'$.

\begin{lem}\label{lem:coh-hprime-deg0} $\mathcal{C}^{\la, \lambda', \mu_0}:=\mathrm{Ext}_{\mathfrak{h}' \oplus \mf{u}_Q, \mu}( (\lambda' ; \mu_0), \mathcal{C}^{\la}\vert_{C^{\dagger}_{w_0^M,Q}})$ is concentrated in degree $0$.
\end{lem}%
\begin{proof}
We can do the computation separately on each of the
  strata $C^{\dagger}_{w_0^M}$ and $C^{\dagger}_{w_0^Ms_{\beta}}$, so we reduce
  to showing that  the cohomology  on the fibers at $w_0^M$ and $ w_0^M s_\beta
  $ of  $\mathcal{C}^{\la}$ vanishes in positive degrees. These fibers are
  respectively $\oscr_{U_{P_{w_0^M}}\backslash G,e}$ and
  $\oscr_{U_{P_{w_0^Ms_\beta}}\backslash G,e}$, and the required vanishing
  follows from a consideration of the actions of $\bar{\mathfrak{p}}_{w_0^M} \cap \mathfrak{b} = \mathfrak{b}$
and  $\bar{\mathfrak{p}}_{w_0^Ms_{\beta}} \cap \mathfrak{b} =
\mathfrak{h} \oplus
\mf{u}_Q$ respectively.
(Note that $\mu = \diag(0,0,1,1)$, so that $w_0^M \mu = \diag(1,1,0,0)$ and
$s_\beta w_0^M \mu = \diag (0,1,0,1)$, and therefore $\mathfrak{h}' \oplus E\cdot w_0^M \mu = \mathfrak{h}' \oplus  E\cdot s_\beta w_0^M \mu = \mathfrak{h}$.)
\end{proof}%
We thus see that \[ \RHom_{\mathfrak{b}, \star_2; \mu}
((\lambda,\mu_0), \mathcal{C}^{\la}\vert_{C^{\dagger}_{w_0^M,Q}})  =
\RHom_{ EX_{\beta}\oplus EH_{\beta}} ( \lambda, \mathcal{C}^{\la,
  \lambda', \mu_0})\] (where the restriction of~$\lambda$ to $EX_{\beta}\oplus E H_{\beta}$ takes $X_\beta \mapsto 0$). %
This cohomology is  represented by the following Chevalley--Eilenberg complex
$K^\bullet$ (in degrees $0,1$ and $2$): %
$$\mathcal{C}^{\la, \lambda', \mu_0} \stackrel{\begin{pmatrix} %
      X_\beta & H_{\beta}-\lambda(H_{\beta})
     \end{pmatrix}}\rightarrow \mathcal{C}^{\la, \lambda', \mu_0} \otimes E(-\beta) \oplus \mathcal{C}^{\la, \lambda', \mu_0} \stackrel{\begin{pmatrix} %
        H_{\beta}-\lambda(H_{\beta}) \\ -X_\beta
     \end{pmatrix}}\rightarrow \mathcal{C}^{\la, \lambda', \mu_0} \otimes E(-\beta).$$
     
  Since  \[HCS_{w_0^M, \lambda}( M(\mathfrak{m}_{w_0^M})_\lambda) = \mathrm{Ext}^0_{\mathfrak{b}}(
  \lambda, \mathcal{C}^{\la}|_{C^{\dagger}_{w_0^M}}),\] we deduce that
  $s_{\beta}^\star HCS_{w_0^M, \lambda}( M(\mathfrak{m}_{w_0^M})_\lambda) \otimes E(-\beta) =
  \mathrm{Ext}^0_{s_{\beta}\mathfrak{b}}( s_{\beta}\lambda, \mathcal{C}^{\la}
  \vert_{ C^{\dagger}_{w_0^M} s_{\beta}}) \otimes E(-\beta)$. On this sheaf,
  $\mathfrak{h}$ acts via $s_{\beta}\lambda - \beta =
  \lambda$  (this is a crucial place where we use that the weight is singular), %
    and $\mu$ still acts via $\mu_0$. Thus, we can consider the map 
     \begin{eqnarray*}
     \mathrm{Ext}^0_{s_{\beta}\mathfrak{b}}( s_{\beta}\lambda, \mathcal{C}^{\la} \vert_{ C^{\dagger}_{w_0^M} s_{\beta}}) \otimes E(-\beta) &\rightarrow& K^1\vert_{ C^{\dagger}_{w_0^M} s_{\beta}} = \mathcal{C}^{\la, \lambda', \mu_0}\vert_{ C^{\dagger}_{w_0^M} s_{\beta}} \otimes E(-\beta) \oplus \mathcal{C}^{\la, \lambda', \mu_0}\vert_{ C^{\dagger}_{w_0^M} s_{\beta}} \\
      s& \mapsto& (s,0),
      \end{eqnarray*}
      which induces a map
      
\numequation\label{eqn:Ext0-to-Ext1}\mathrm{Ext}^0_{s_{\beta}\mathfrak{b}}(
s_{\beta}\lambda, \mathcal{C}^{\la} \vert_{ C^{\dagger}_{w_0^M} s_{\beta}})
\otimes E(-\beta) \rightarrow \mathrm{Ext}^1_{\mathfrak{b}, \star_2;
  \mu} ((\lambda,\mu_0),
\mathcal{C}^{\la}\vert_{C^{\dagger}_{w_0^M}s_{\beta}}). \end{equation}%
We will show below that~\eqref{eqn:Ext0-to-Ext1} is an isomorphism. 
We begin by studying its restrictions to $C^{\dagger}_{w_0^M} \cap
    C^{\dagger}_{w_0^M}s_{\beta}$ and to $
  C^{\dagger}_{{w_0^M}s_{\beta}}$; it is immediate from the definitions that the
  former restriction is a map in  $\LB_{(\mathfrak{g},B \cap
      B_{s_{\beta}})}(C^{\dagger}_{w_0^M} \cap
    C^{\dagger}_{w_0^M}s_{\beta})$, and the latter is a map in  $\LB_{(\mathfrak{g},B
    \cap B_{s_{\beta}})}( C^{\dagger}_{{w_0^M}s_{\beta}})$.
\begin{lem} 
  \label{lem:lemmas-towards-Cartan-Sen-computation}\leavevmode
  \begin{enumerate}
  \item\label{item:restrictions-are-in-LB} The restrictions of the left and right hand sides of ~\eqref{eqn:Ext0-to-Ext1} to $C^{\dagger}_{w_0^M} \cap
    C^{\dagger}_{w_0^M}s_{\beta}$ are 
    isomorphic, and their endomorphism algebras in $\LB_{(\mathfrak{g},B \cap
      B_{s_{\beta}})}(C^{\dagger}_{w_0^M} \cap
    C^{\dagger}_{w_0^M}s_{\beta})$
 are  the scalars~$E$.
\item The restrictions of the left and right hand sides  of ~\eqref{eqn:Ext0-to-Ext1} to $C^{\dagger}_{{w_0^M}s_{\beta}}$ are isomorphic, and their endomorphisms in  $\LB_{(\mathfrak{g},B
    \cap B_{s_{\beta}})}(C^{\dagger}_{{w_0^M}s_{\beta}})$ are the scalars~$E$.
  \end{enumerate}
\end{lem}
\begin{proof} To begin, we note that it follows from the $Q$-equivariance of  $HCS_{Q,w_0^M, \eta}(  M(\mathfrak{m}_{w_0^M})_\lambda)$ that  $$s_\beta^\star (HCS_{Q,w_0^M, \eta}(  M(\mathfrak{m}_{w_0^M})_\lambda) \otimes E(\lambda-\eta +\beta)) = HCS_{Q,w_0^M, \eta}(  M(\mathfrak{m}_{w_0^M})_\lambda) \otimes E(\lambda-\eta).$$
It follows that %
\begin{eqnarray*}
\mathrm{Ext}^0_{s_{\beta}\mathfrak{b}}( s_{\beta}\lambda, \mathcal{C}^{\la} \vert_{ C^{\dagger}_{w_0^M}s_{\beta}} )\vert_{ C^{\dagger}_{w_0^M} \cap C^{\dagger}_{w_0^M} s_{\beta}} &=&  s_{\beta}^\star  HCS_{w_0^M, \lambda}( M(\mathfrak{m}_{w_0^M})_\lambda) \otimes E(-\beta) \ \vert_{ C^{\dagger}_{w_0^M} \cap C^{\dagger}_{w_0^M} s_{\beta}}   \\
& =& HCS_{w_0^M, \lambda}( M(\mathfrak{m}_{w_0^M})_\lambda)  \ \vert_{ C^{\dagger}_{w_0^M} \cap C^{\dagger}_{w_0^M} s_{\beta}}. %
\end{eqnarray*} 
On the other hand \[\mathrm{Ext}^1_{\mathfrak{b}, \star_2; \mu}
  ((\lambda,\mu_0),
  \mathcal{C}^{\la}\vert_{C^\dagger_{w_0^M}s_{\beta}})\vert_{
    C^\dagger_{w_0^M} \cap C^\dagger_{w_0^M} s_{\beta}} =HCS_{w_0^M,
    \lambda}( M(\mathfrak{m}_{w_0^M})_\lambda)  \ \vert_{
    C^\dagger_{w_0^M} \cap C^\dagger_{w_0^M} s_{\beta}}\](by
Proposition~\ref{lem-extendingunionstrata} and the proof of
Proposition~\ref{prop:we-have-exact-triangle}.)   %

We now check that the  endomorphisms of the sheaf $HCS_{w_0^M, \lambda}( M(\mathfrak{m}_{w_0^M})_\lambda)  \ \vert_{ C^\dagger_{w_0^M} \cap C^\dagger_{w_0^M} s_{\beta}}$ in the category $\LB_{(\mathfrak{g},B \cap B_{s_{\beta}})}(C^{\dagger}_{w_0^M} \cap C^{\dagger}_{w_0^M}s_{\beta})$ are scalars.  

Since $B \cap B_{s_{\beta}}$ acts transitively on $C_{w_0^M} \cap C_{w_0^M}s_{\beta}$  by section \ref{sec-actionB},
an endomorphism is determined by its behavior at any fiber. 
More precisely,  let $x \in C_{w_0^M} \cap C_{w_0^M}s_{\beta}$. We consider the uniformization map 
\begin{eqnarray*}
G_e \rtimes (B \cap B_{s_{\beta}}) &\rightarrow &C^{\dagger}_{w_0^M} \cap C^{\dagger}_{w_0^M}s_{\beta} \\
(g,b) \mapsto xgb.
\end{eqnarray*}

 Let $\Stab(x)$ be the stabilizer of $x$ for this action. We have
 $\Stab(x) = \{ (g,b), gb \in P_x\}$. Exactly as in the proof of Lemma~\ref{lem-decompose-stab(w)-Q}, this group is   
 $$ (P_x  \cap B \cap  B_{s_{\beta}})_e \backslash [ ((P_x)_e \times (B \cap B_{s_{\beta}})_e) \rtimes (P_x \cap B \cap  B_{s_{\beta}})].$$
 The fiber we consider is the completion of a dual Verma module for
 $\mathfrak{m}_x$ (by Definitions~\ref{defn:HCS-functor} and~
 \ref{defn-HigherColeman}), and the first part of the lemma follows from the
 property that the endomorphisms of a Verma module are the scalars
 (together with Theorem~\ref{thm-schmidt}).

  The second part is proved in the same way, as follows. We first observe that 
\begin{eqnarray*}
\mathrm{Ext}^0_{s_{\beta}\mathfrak{b}}( s_{\beta}\lambda, \mathcal{C}^{\la} \vert_{ C^{\dagger}_{w_0^M} s_{\beta}}) \otimes E(-\beta)\vert_{C^\dagger_{w_0^Ms_{\beta}} }& =& HCS_{w_0^Ms_\beta, \lambda}( M(\mathfrak{m}_{w_0^Ms_\beta})_{s_\beta\lambda})  %
  \\
&=& \mathrm{Ext}^1_{\mathfrak{b}, \star_2; \mu} ((\lambda,\mu_0), \mathcal{C}^{\la}\vert_{C^{\dagger}_{w_0^Ms_{\beta}}})
\end{eqnarray*}
We reduce to showing that the endomorphisms of $HCS_{w_0^Ms_\beta, \lambda}( M(\mathfrak{m}_{w_0^Ms_\beta})_{s_\beta\lambda})$ in the category $\LB_{(\mathfrak{g},B \cap
  B_{s_{\beta}})}(C^{\dagger}_{w_0^Ms_{\beta}})$ are scalars which again follows from the property that endomorphisms of Verma modules are scalar. 
\end{proof}

\begin{prop}
  \label{key-prop}The map~\eqref{eqn:Ext0-to-Ext1} is an isomorphism.
\end{prop}
\begin{proof}
  By Lemma~\ref{lem:lemmas-towards-Cartan-Sen-computation} %
   it
  suffices to show that ~\eqref{eqn:Ext0-to-Ext1} is nonzero on the
  fiber at  one point of $C_{w_0^Ms_{\beta}}$ and at one point of
  $C_{w_0^M} \cap C_{w_0^M}s_{\beta}$. It suffices in turn to prove that the map  \[
\begin{tikzcd}
\mathcal{C}^{\la, \lambda', \mu_0} \arrow[rr, "{H_\beta-\lambda(H_{\beta})}"] && \mathcal{C}^{\la, \lambda', \mu_0}
\end{tikzcd}
\]
   induces an injective map on
   the fibers at $w_0^Ms_{\beta}$ and at
   one point of  $C_{w_0^M} \cap
   C_{w_0^M}s_{\beta}$. %

  By definition, the kernel of this map on the fiber at a
  point~$x$ is
  \[\Hom_{\mf{h}\oplus\mf{u}_Q, \star_2; \mu}
((\lambda,\mu_0), \mathcal{C}^{\la}_x).\]  We first consider the fiber
at $x=w_0 ^Ms_{\beta}$, where $\mathcal{C}^{\la}_{w_0 ^Ms_{\beta}}  =
\oscr_{U_{P_{^2w}}\backslash G,e}$. We have  \[\mathfrak{u}_{P_{w_0
      ^Ms_{\beta}}} = E X_{\beta}  \oplus E X_{-\alpha} \oplus E X_{-\delta}, ~~\mathfrak{u}_Q = E X_\alpha \oplus X_{\delta} \oplus X_{\gamma}.\]

For any~$s \in \Phi$ we write $U_{s}$ for the corresponding
$1-$parameter subgroup.  We pick a coordinate   $x_s$ on $U_s$ with
the property that the corresponding vector field is ~ $X_s$. Then elements of $\Hom_{\mf{h}\oplus\mf{u}_Q, \star_2}
(\lambda, \mathcal{C}^{\la}_x)$
are germs of analytic functions on $U_{-\beta} U_{-\gamma}T$, which can
be written as
$$f(x_{-\beta},x_{-\gamma}, t) = \sum_{k_{-\beta}, k_{-\gamma} \in \ZZ_{\geq 0}} x_{-\beta}^{k_{-\beta}} x_{-\gamma}^{k_{-\gamma}} \lambda(t).$$

We need to show that if $(\Theta_{\hor}(\mu)-\mu_0)f=0$, then~$f=0$. By definition, $\Theta_{\hor}(\mu)$ acts on the left via the action of $-(w_0 ^Ms_{\beta})^{-1}\mu\in \mf{h}$.
 It follows that   
\begin{eqnarray*}
\Theta_{\hor}(\mu) & =& \langle (w_0 ^Ms_{\beta})^{-1}\mu,  \beta \rangle x_{-\beta}\partial_{x_{-\beta}} + \langle (w_0 ^Ms_{\beta})^{-1}\mu,  \gamma \rangle x_{-\gamma}\partial_{x_{-\gamma}} - \langle (w_0 ^Ms_{\beta})^{-1}\mu,  \lambda \rangle\\
&=& x_{-\beta}\partial_{x_{-\beta}} - \frac{w}{2}.
\end{eqnarray*}
(We are using here that  $(w_0 ^Ms_{\beta})^{-1}\mu = (-\frac{1}{2},
\frac{1}{2}; \frac{1}{2})$, $\beta = (-2,0;0)$, $\gamma = (-1,-1;0)$,
$\lambda = (1,1;w)$.) Thus $\Theta_{\hor}(\mu)$ acts by $k_{-\beta}-\frac{w}{2}$ on  $x_{-\beta}^{k_{-\beta}} x_{-\gamma}^{k_{-\gamma}} \lambda(t)$. Since $\mu_0 = -1- \frac{w}{2}$, we deduce that $\Theta_{\hor}(\mu)-\mu_0$ is injective.

We now consider the point $w_0 ^Ms_{\beta}g_{-\beta}\in C_{w_0 ^M}
  \cap C_{w_0 ^M}s_\beta $, where $g_{-\beta}=\exp(X_{-\beta})$. We see that \[(w_0 ^Ms_{\beta}g_{-\beta})^{-1}\mu = (w_0 ^Ms_{\beta})^{-1}\mu + \langle (w_0 ^Ms_{\beta})^{-1}\mu,\beta \rangle X_{-\beta}.\]
An element of $\mathcal{C}^{\la, \lambda', \mu_0}_{^2wg_{-\beta} }$
can still be expressed as a germ  of an analytic function on
$U_{-\beta} U_{-\gamma}T$, and can thus be written as
$$f(x_{-\beta},x_{-\gamma}, t) = \sum_{k_{-\beta}, k_{-\gamma} \in \ZZ_{\geq 0}} x_{-\beta}^{k_{-\beta}} x_{-\gamma}^{k_{-\gamma}} \lambda(t).$$

We now find that $\Theta_{\hor}(\mu) =  x_{-\beta}\partial_{x_{-\beta}} - \frac{w}{2} + \partial_{x_{-\beta}}$, so that 
\[\Theta_{\hor}(\mu) x_{-\beta}^{k_{-\beta}} x_{-\gamma}^{k_{-\gamma}}
\lambda(t) = (k_{-\beta}-\frac{w}{2}) x_{-\beta}^{k_{-\beta}}
x_{-\gamma}^{k_{-\gamma}} \lambda(t)  + k_{-\beta}
x_{-\beta}^{k_{-\beta}-1} x_{-\gamma}^{k_{-\gamma}} \lambda(t),\] and
we again deduce that $\Theta_{\hor}(\mu)-\mu_0$ is injective. 
\end{proof}
\begin{proof}[Proof of Theorem \ref{thm-extequal}]
  By Proposition \ref{key-prop}, the sheaf  $\mathrm{Ext}^1_{\mathfrak{b}, \star_2; \mu}
((\lambda,\mu_0), \mathcal{C}^{\la}\vert_{C^{\dagger}_{w_0^M,Q}})$ is
obtained by gluing   $HCS_{w_0^M,
      \lambda}(M(\mf{m}_{w_0^M})_{\lambda})$ and  $s_\beta^\star HCS_{w_0^M,
      \lambda}(M(\mf{m}_{w_0^M})_{\lambda})\otimes E(-\beta)$ along $C^\dagger_{w_0^M} \cap
    C^\dagger_{w_0^M}s_{\beta}$. The same is true of  $HCS_{Q,w_0^M,
      \eta}(M(\mf{m}_{w_0^M})_{\lambda}) \otimes E(\lambda-\eta)$ by construction.  The gluing data is that of an isomorphism of
  $(\mathfrak{g}, B \cap B_{s_\beta})$-equivariant sheaves, and by
  Lemma~\ref{lem:lemmas-towards-Cartan-Sen-computation} the space of
  such isomorphisms identifies with $E^\times$, so there is (up to
  isomorphism) a unique way to glue, as required.%
\end{proof}

\section{$p$-adic Eichler--Shimura theory}   \label{sec:application-to-Shimura}%

\subsection{Introduction}
\label{sec:eichlerintro}
The main goal of this section (as mentioned in~\S\ref{sec:babyintro})
is to relate higher Coleman
theory to completed cohomology, so that (ultimately) we can connect
the Galois-theoretic properties of a $p$-adic ordinary (overconvergent)
modular form (in terms of the action of the Sen operator) to
its classicality.

Before proceeding, we introduce some notation.
We fix a Hodge type Shimura datum $(G,X)$. We assume that $G_{\Q_p}$ is quasi-split. Let $P$ be the parabolic corresponding to $\mu$ with Levi $M$.  Let $B \subseteq G_{\Q_p}$ be a Borel subgroup. We pick a maximal torus $T \subseteq B$.  The relevant flag variety is  $FL = P \backslash G$ and we have the decomposition $ FL = \coprod_{w \in \WM} P \backslash P w B$ where $\WM \subseteq W$ is the set of Kostant representatives in the absolute Weyl group. We also fix a coefficient field $E$ which is a finite extension of $\qq_p$ and admits a map from the reflex field of the Shimura datum. 
\subsubsection{Higher Coleman theory} 
Higher Coleman theory~\cite{boxer2021higher} is a
theory of (higher) overconvergent modular forms. The different higher Coleman
theories are  parameterized by two parameters: an element $w \in \WM$, and a weight.
\begin{rem} We note that $\WM$ parametrizes chambers in the weight space which are $M$-dominant, and the $w$-theory will interpolate those classical cohomologies whose weights belong to the $w$-chamber. 
\end{rem}

To describe the weight parameter, we fix $w \in \WM$.  Let $\mathfrak{m} =
\mathrm{Lie}(M)$, and let $\mathfrak{m}_w = w^{-1} \mathfrak{m} w$. Let $\mathcal{O}(\mathfrak{m}_w, \mathfrak{b}_{\mathfrak{m}_w})$ be the BGG category $\ocal$ of $U(\mathfrak{m}_w)$-modules for the Borel $\mathfrak{b}_{\mathfrak{m}_w} = \mathrm{Lie}(B) \cap \mathfrak{m}_w$. In the same way that weights of modular forms are  finite dimensional representations of $M$, weights of (higher) overconvergent modular forms (parameterized by $w \in \WM$) are objects of $\mathcal{O}(\mathfrak{m}_w, \mathfrak{b}_{\mathfrak{m}_w})$.  For any $\lambda \in X^\star(T)_{E}$, we let $\mathcal{O}(\mathfrak{m}_w, \mathfrak{b}_{\mathfrak{m}_w})_{\lambda-\alg}$ be the subcategory of 
category $\mathcal{O}(\mathfrak{m}_w, \mathfrak{b}_{\mathfrak{m}_w})$ with weights in $\lambda + X^\star(T)$. 

We have higher Coleman functors (see Definition  \ref{defn:HCF}):
\numequation
\label{thefunctorsintro}
 HC_{w,\lambda}, HC_{\cusp, w, \lambda}:  \mathcal{O}(\mathfrak{m}_w,
\mathfrak{b}_{\mathfrak{m}_w})_{\lambda-\alg}^{\mathrm{op}} \rightarrow
D(\Mod^{\lambda-\sm}_{B(\qq_p)}(E))
\end{equation}
 where $\Mod^{\lambda-\sm}_{B(\qq_p)}(E)$ is the category of locally analytic representations of $B(\qq_p)$ with $\mathfrak{b}$ acting like $\lambda$. For example, if $\lambda=0$, this is just the category of smooth $B(\qq_p)$-representations.

\begin{rem} The parameter $\lambda$ is there to specify the $B(\qq_p)$-action. If $\eta \in X^\star(T)$, the categories 
$\mathcal{O}(\mathfrak{m}_w, \mathfrak{b}_{\mathfrak{m}_w})_{\lambda-\alg}$ and $\mathcal{O}(\mathfrak{m}_w, \mathfrak{b}_{\mathfrak{m}_w})_{\lambda+\eta-\alg}$ are canonically the same and one has $HC_{w,\lambda+\eta}(M) = HC_{w,\lambda}(M) \otimes E(\eta)$ where $-\otimes E(\eta)$ means a twist of the $B(\qq_p)$-action by $\eta$. 
\end{rem}

The functors~\eqref{thefunctorsintro} are defined by first attaching to every object $M$ of 
$\mathcal{O}(\mathfrak{m}_w, \mathfrak{b}_{\mathfrak{m}_w})_{\lambda-\alg}$ a ``quasi-coherent'' $B(\qq_p)$-equivariant sheaf over the pullback of the Bruhat stratum  $P \backslash PwB$ via the Hodge--Tate period map and taking its cohomology with suitable support condition.

 We also have a finite slope part functor $D(\Mod^{\lambda-\sm}_{B(\qq_p)}(E)) \rightarrow D(\Mod^{\lambda-\sm}_{T(\qq_p)}(E))$ and we can speak of the finite slope part of higher Coleman functors (Section \ref{sec-fspHC}): 
 
 $$HC^{\nfs}_{w,\lambda}, HC^{\nfs}_{\cusp, w, \lambda}:  \mathcal{O}(\mathfrak{m}_w, \mathfrak{b}_{\mathfrak{m}_w})_{\lambda-\alg}^{\mathrm{op}} \rightarrow D(\Mod^{\lambda-\sm}_{T(\qq_p)}(E)).$$

\begin{rem}
Let  $\kappa \in X^\star(T)^{+,M}$. Let $L(\mathfrak{m}_w)_{-w^{-1}w_{0,M}\kappa}$ be the finite dimensional representation of $\mathfrak{m}_w$ of   of highest weight $-w^{-1}w_{0,M} \kappa$. We show  that $$HC^{\nfs}_{w, 0}(L(\mathfrak{m}_w)_{-w^{-1}w_{0,M}\kappa})$$ is  the direct sum of the $\mathrm{R}\Gamma_w(K^p, \kappa, \chi)^{+,\nfs}$ of \cite{boxer2021higher} (Theorem \ref{thm-compare-Higher-Coleman-old-new}). These are higher Coleman theories with value in the classical sheaf of weight $\kappa$. 
   On the other hand, if  $M(\mathfrak{m}_w)_{-w^{-1}w_{0,M}\kappa}$ denotes the Verma of highest weight ${-w^{-1}w_{0,M}\kappa}$, then $$HC_{w, 0}^{\nfs}(M(\mathfrak{m}_w)_{-w^{-1}w_{0,M}\kappa})$$  corresponds to  higher Coleman theory  with value in the  big ``induction'' sheaf (Theorem \ref{thm:another-BP21-comparison}).  The surjective map $M(\mathfrak{m}_w)_{-w^{-1}w_{0,M}\kappa} \rightarrow L(\mathfrak{m}_w)_{-w^{-1}w_{0,M}\kappa}$ induces a 
 map $$HC_{w,0}( L(\mathfrak{m}_w)_{-w^{-1}w_{0,M} \kappa}) \rightarrow HC_{w,\lambda}(M(\mathfrak{m}_w)_{-w^{-1}w_{0,M}\kappa})$$ and similarly on the finite slope part. 
 In summary, the improvements on \cite{boxer2021higher} are the following:
 \begin{itemize}
 \item We extend the definitions to the infinite slope part (in \cite{boxer2021higher}, only the finite slope part was canonically defined). 
 \item We introduce a more functorial perspective on the weights. In
   \cite{boxer2021higher} we allowed weights to be either finite dimensional
   representations or Verma modules, which of course  generate the BGG category. 
 \end{itemize}
\end{rem}

Our main results on higher Coleman theory can be summarized as follows.

\begin{thm}[Theorems \ref{thm-coho-amplitudeHC},  \ref{coro-vanishingThm}, and \ref{thm-bounds-on-slopes}]
\leavevmode
\begin{enumerate}
\item $HC_{w,\lambda}$,  $HC^{\nfs}_{w,\lambda}$ have cohomological amplitude  $[\ell(w), d]$ and $HC_{\cusp, w,\lambda}$, $HC^{\nfs}_{\cusp, w,\lambda}$  have cohomological amplitude  $[0, \ell(w)]$. 
\item Let $M \in \ocal(\mathfrak{m}_w, \mathfrak{b}_{\mathfrak{m}_w})_{\lambda-\alg}$ be a module generated by a highest weight vector of weight $\nu$. Assume that the Shimura variety is proper or that we are in the Siegel case. 
The slopes appearing in $HC^{\nfs}_{\cusp, w, \lambda}(M)$ and $HC^{\nfs}_{w, \lambda}(M)$ are $\geq \lambda-\nu + w^{-1}w_{0,M} \rho + \rho$. 
\end{enumerate}

\end{thm}

\begin{rem} 
In the proper case, the functors $HC_{w,\lambda}$ and $HC^{\nfs}_{w, \lambda}$ are exact. 
\end{rem}

\subsubsection{Completed cohomology}
We let $\mathrm{R}\Gamma(\Sh_{K^p}^{\tor}, \qq_p)$ be completed cohomology and $\mathrm{R}\Gamma_c(\Sh_{K^p}, \qq_p)$ denote completed cohomology with compact support.

We let $\ocal(\mathfrak{g}, \mathfrak{b})$ be the BGG category $\ocal$ for $\mathfrak{g}$ and $\mathfrak{b}$. 
We define  functors (Section \ref{sect-general-padicES}):

\begin{eqnarray*}
CC_{\lambda}: \ocal(\mathfrak{g}, \mathfrak{b})_{\lambda-\alg} &\rightarrow& D( \Mod^{\lambda-\sm}_{B(\qq_p)}(E)) \\
M & \mapsto & \mathrm{RHom}_{\mathfrak{g}}( M, \mathrm{R}\Gamma(\Sh_{K^p}, \qq_p)^{\la} )
\end{eqnarray*}
\begin{eqnarray*}
CC_{\cusp, \lambda}:  \ocal(\mathfrak{g}, \mathfrak{b})_{\lambda-\alg} &\rightarrow& D( \Mod^{\lambda-\sm}_{B(\qq_p)}(E)) \\
M & \mapsto & \mathrm{RHom}_{\mathfrak{g}}( M, \mathrm{R}\Gamma_c(\Sh_{K^p}, \qq_p)^{\la} )
\end{eqnarray*}

\begin{rem}
 The first natural example is to apply these functors to  $M$  a finite dimensional representation of $G$, in which case we recover classical \'etale cohomology with weight $M^\vee$ (it is natural to take $\lambda=0$). 
For a  non-classical example, we can  take $M = U(\mathfrak{g})\otimes_{U(\mathfrak{b})} \lambda$ to be a  Verma module of weight $\lambda$ so that we are computing $\mathfrak{b}$-cohomology.
\end{rem}

\subsubsection{$p$-adic Eichler Shimura} We are now ready to state our main
result
comparing completed cohomology and higher Coleman theory.
 It holds under a non-Liouville condition on $\lambda$ (see Definition \ref{defn-number-non-Liouville}). We observe that if $\lambda \in X^\star(T)$ is algebraic, it is non-Liouville. 

\begin{thm}[Theorem \ref{thm-p-adic-ES}]  Assume $\lambda$ is non-Liouville and $M$ is an object of $\ocal(\mathfrak{g}, \mathfrak{b})_{\lambda-\alg}$. We have a spectral sequence 
$$E_1^{p,q} = \oplus_{w \in \WM, \ell(w)=p} H^{p+q}( HC_{w,\lambda}( M \otimes_{\mathfrak{u}_{\mathfrak{p}_w}}^L E))$$
converging to $H^{p+q}(CC_\lambda(M))\otimes \C_p$. Moreover, the Sen operator is    given by $w\mu \in Z({\mathfrak{m}_w})$ acting on  $H_\star(\mathfrak{u}_{\mathfrak{p}_w},M)$. 
\end{thm}

\begin{rem} The functor $ -  \otimes_{U(\mathfrak{u}_{\mathfrak{p}_w})} E:  D(\ocal(\mathfrak{g}, \mathfrak{b}))\rightarrow D(\ocal(\mathfrak{m}_w, \mathfrak{b}_{\mathfrak{m}_w}))$ is the Lie algebra homology of the unipotent radical $\mathfrak{u}_{\mathfrak{p}_w}$ of $\mathfrak{p}_w = w^{-1}\mathrm{Lie}(P) w$. 
It is computed by the Koszul complex  (in degree $-d$ to $0$):

$$0 \rightarrow M \otimes \Lambda^d \mathfrak{u}_{\mathfrak{p}_w} \rightarrow  \cdots \rightarrow M \rightarrow 0$$
\end{rem}
\begin{rem}  In Section \ref{subsec:localization-partial-flag-variety} we attached to $M \in \ocal(\mathfrak{g}, \mathfrak{b})$ a certain twisted $D$-module $\Loc(M)$ on the flag variety.  This is a version of Beilinson--Bernstein localization. This twisted $D$-module completely encodes  the $p$-adic Eichler--Shimura theory (see Theorem \ref{thm-p-adic-ES} for a precise statement). We observe that $\Loc(M)$ is  ``constant'' on each Bruhat stratum $C_w$ and  its restriction to each $C_w$ is determined (in the non-Liouville case) by the Lie algebra homology $M \otimes_{\mathfrak{u}_{\mathfrak{p}_w}}^L E$. 
\end{rem} 

Under favorable ``genericity'' assumptions, the spectral sequence simplifies a lot. Let us denote by $M(\mathfrak{g})_\lambda =  U(\mathfrak{g})\otimes_{U(\mathfrak{b})} \lambda$ the Verma module of weight $\lambda$. We adopt a similar notation to denote Verma modules for other reductive Lie algebras. 

 \begin{cor}[Corollary \ref{coro-simplification-anti}] Assume that $\lambda$ is non-Liouville and  antidominant in the sense of
  Remark~\ref{rem-antidom}, and that the Shimura variety is proper.  Then  $CC_\lambda(M(\mathfrak{g})_\lambda)$ is concentrated in the middle degree $d$ and moreover, it has a decreasing filtration $\mathrm{Fil}^i \mathrm{H}^d(CC_\lambda(M(\mathfrak{g})_\lambda))$ with 
\begin{itemize}
\item $ \mathrm{Fil}^{d+1} \mathrm{H}^d(CC_\lambda(M(\mathfrak{g})_\lambda)) = 0 $,
\item $  \mathrm{Fil}^0 \mathrm{H}^d(CC_\lambda(M(\mathfrak{g})_\lambda)) =  \mathrm{H}^d(CC_\lambda(M(\mathfrak{g})_\lambda))$,
\item  $\mathrm{Gr}^p \mathrm{H}^d(CC_\lambda(M(\mathfrak{g})_\lambda)) =
  \oplus_{ w \in \WM, \ell(w) = p}
  \HH^{p}(HC_{w,\lambda}(M(\mathfrak{m}_w)_{\lambda + w^{-1}w_{0,M}\rho +
    \rho}))$.%
\end{itemize}
\end{cor}
We also refer to Theorem \ref{thm-p-adic-ES3} for a similar result in the ordinary case.

\subsection{Perfectoid Shimura varieties}%
We consider a Hodge-type Shimura datum $(G,X)$.  We also fix a map from the reflex field of the Shimura datum to 
$E$. %
We let $\Sh_{K_pK^p}^{rat}$ be  the analytic space over $\Spa(E,\O_E)$ attached to the Shimura variety of level $K_pK^p$ over $\Spec~E$. We let $\Sh_{K_pK^p} = \Sh_{K_pK^p}^{rat} \times_{\Spa(E, \ocal_E)} \Spa (\C_p, \ocal_{\C_p}) $.  
We let $\Sh_{K_pK^p,\Sigma}^{rat, \tor}$ be a toroidal compactification of $\Sh_{K_pK^p}^{rat}$  over $\Spa(E, \ocal_E)$  %
 for a specific choice of cone decomposition $\Sigma$ (see \cite{MR1083353}, \cite{MR3186092}). For $K_p' \subseteq K_p$, we have a natural map $\Sh_{K'_pK^p,\Sigma}^{rat, \tor} \rightarrow \Sh_{K_pK^p,\Sigma}^{rat, \tor}$. 
We let  $\Sh_{K^p,\Sigma}^{rat,\tor} = \lim_{K'_p}\Sh_{K'_pK^p,\Sigma}^{rat,\tor}$.  By \cite{scholze-torsion} (and \cite{MR3512528}, \cite{MR4492226} for
the extension to the compactification), this is a perfectoid space. We note that the same cone
decomposition $\Sigma$ is used at each stage of the limit and that there is some restriction on the choice of cone decomposition if we are not in the Siegel case. 
Concretely, the underlying topological space of $\Sh_{K^p,\Sigma}^{rat,\tor}$ is the inverse limit of the topological spaces of the  $\Sh_{K'_pK^p,\Sigma}^{rat,\tor}$, and there is a basis of affinoid opens $\Spa(A, A^+)$ of $\Sh_{K^p,\Sigma}^{rat,\tor}$, which are pull backs of affinoid opens $\Spa (A_{K'_p}, A_{K'_p}^+)$ in $\Sh_{K'_pK^p,\Sigma}^{rat,\tor}$ for small enough $K'_p$, and   such that $A^+$ is the $p$-adic completion of $\colim_{K'_p} A_{K'_p}^+$. 

   We let $\Sh_{K^p}^{rat} \hookrightarrow \Sh_{K^p,\Sigma}^{rat,\tor}$ be the open subspace $ \lim_{K'_p}\Sh_{K'_pK^p}^{rat}$.
We let 
$\pi^{rat}_{HT,\Sigma}: \Sh_{K^p,\Sigma}^{rat,\tor} \rightarrow \mathcal{FL}^{rat}$ be
the Hodge--Tate period map, where $\mathcal{FL}^{rat}$ is the Flag
variety over $\Spa(E, \ocal_E)$. %

We let $\Sh_{K_pK^p,\Sigma}^{\tor}$ be the base change to $\Spa( \C_p,
\ocal_{\C_p})$ of $\Sh_{K_pK^p,\Sigma}^{rat,\tor}$ (which therefore
carries an action of $\mathrm{Gal}(\bar{E}/E)$). We similarly define
$\Sh_{K^p,\Sigma}^{\tor}$, $\Sh_{K^p}$  and $\mathcal{FL}$, and the period map
$\pi_{HT,\Sigma}: \Sh_{K^p,\Sigma}^{\tor} \rightarrow \mathcal{FL}$
which is $\mathrm{Gal}(\bar{E}/E)$-equivariant.  %
The period map is also $K_p$-equivariant.
The action of ~ $G(\Qp)$ on~  $\Sh_{K^p, \Sigma}$ does not extend to an action on
$\Sh_{K^p, \Sigma}^{\tor}$ but  for a general $g \in G(\qq_p)$ we still have diagrams: 

\begin{eqnarray*}
\xymatrix{   \Sh_{K^p,\Sigma}^{\tor} \ar[r]^g \ar[d]  & \Sh_{K^p, g\Sigma}^{\tor} \ar[d] \\
\mathcal{FL} \ar[r]^g & \mathcal{FL}}
\end{eqnarray*}

\begin{rem}\label{rem-Sigma} The choice of a specific $\Sigma$ does not usually play any role. If no confusion is likely to arise, we fix a  $\Sigma$ and  drop it from the notation. We will eventually allow ourselves to change $\Sigma$. It is also important to note that all the cohomologies we will consider  (coherent cohomology, completed cohomology) do not depend on $\Sigma$. 
\end{rem}

\subsection{Smooth and locally analytic vectors of the structure sheaf}
We let $\oscr_{\Sh_{K^p}^{\tor}}$ be the structure sheaf of the
perfectoid space. If $U \subseteq \Sh_{K^p}^{\tor}$ is a quasi-compact
open, then $U$ is stabilized by an open subgroup $K_{p}$ %
of $G(\qq_p)$ and we get a continuous action of $K_{p}$ on $ \oscr_{\Sh_{K^p}^{\tor}}(U)$.  One can speak of the smooth  and locally analytic vectors of $\oscr_{\Sh_{K^p}^{\tor}}(U)$. We thus obtain subsheaves of smooth and locally analytic vectors: 
$$\oscr^{\sm}_{\Sh_{K^p}^{\tor}} \subseteq \oscr^{\la}_{\Sh_{K^p}^{\tor}} \subseteq \oscr_{\Sh_{K^p}^{\tor}}.$$

\begin{prop} For any compact open subgroup $K'_p \subseteq G(\qq_p)$, let $\pi_{K'_p} : \Sh_{K^p}^{\tor}
\rightarrow \Sh_{K^pK'_p}^{\tor}$ be the natural map. The pullback map  $\colim_{K'_p} \pi_{K'_p}^{-1} \oscr_{\Sh_{K'_pK^p}^{\tor}} \rightarrow \oscr^{\sm}_{\Sh_{K^p}^{\tor}}$ is an isomorphism. 
\end{prop}
\begin{proof} We consider the map of sites $\nu :
  (\Sh_{K^pK_p}^{\tor})_{\proket} \rightarrow
  (\Sh_{K^pK_p}^{\tor})_{\mathrm{ket}}$ from the pro-Kummer-\'etale site to the
  Kummer-\'etale site. It follows from \cite[Coro. 6.19]{MR3090230} (which is
  easily extended to the Kummer-\'etale case via the machinery of
  \cite{diao2022logarithmic}) %
  that $\nu_\star \hat{\oscr}_{\Sh_{K^pK_p}^{\tor}} = \oscr_{\Sh_{K^pK_p}^{\tor}}$. We are going to see that the proposition follows directly from this statement. 
We consider the pro-Kummer \'etale cover $\pi_{K_p} : \Sh_{K^p}^{\tor}
\rightarrow \Sh_{K^pK_p}^{\tor}$. Since $\Sh_{K^p}^{\tor}$ is perfectoid, it
follows  from \cite[Thm.\ 5.4.3]{diao2022logarithmic} that for any open affine
$U\subseteq  \Sh_{K^p}^{\tor}$,   $\hat{\oscr}_{\Sh_{K^pK_p}^{\tor}}(U) =
\oscr_{\Sh_{K^p}^{\tor}}(U)$. Since %
$\Sh_{K^p}^{\tor} \times_{\Sh_{K^pK_p}^{\tor}} \Sh_{K^p}^{\tor} = K_p \times
\Sh_{K^p}^{\tor}$, we deduce that  ${\oscr}_{\Sh_{K^pK_p}^{\tor}}=\nu_\star \hat{\oscr}_{\Sh_{K^pK_p}^{\tor}}
= \HH^0(K_p, \pi_{K_p,\star} \oscr_{\Sh_{K^p}}^{\tor})$. The proposition 
follows by taking the colimit over~$K_p$.
\end{proof}

We let $\mathscr{I}_{\Sh_{K_pK^p}^{\tor}}$ be the ideal of the (reduced) boundary $D_{K_p}$ in $\Sh_{K_pK^p}^{\tor}$.   This is an invertible ideal in $\mathscr{O}_{\Sh_{K_pK^p}^{\tor}}$ that we also denote by $\oscr_{\Sh_{K_pK^p}^{\tor}}(-D_{K_p})$.   
We let $\mathscr{I}^{\sm}_{\Sh_{K^p}^{\tor}}  = \colim_{K'_p} \mathscr{I}_{\Sh_{K'_pK^p}^{\tor}}$. 

We let $\mathscr{I}_{\Sh_{K^p}^{\tor}}$ be the ideal of the boundary
in the structure sheaf $\oscr_{\Sh_{K^p}^{\tor}}$. Its subsheaf of
locally analytic vectors is denoted by
$\mathscr{I}_{\Sh_{K^p}^{\tor}}^{\la}$ and turns out to be equal to
$\oscr_{\Sh^{\tor}_{K^p}}^{\la}\otimes_{\oscr^{\sm}_{\Sh^{\tor}_{K^p}}} \mathscr{I}^{\sm}_{\Sh_{K^p}^{\tor}}$ by the following lemma. 

\begin{lem}\label{lem-deducing-ideal} The natural map: $\oscr_{\Sh^{\tor}_{K^p}}^{\la}\otimes_{\oscr^{\sm}_{\Sh^{\tor}_{K^p}}} \mathscr{I}^{\sm}_{\Sh_{K^p}^{\tor}} \rightarrow 
\mathscr{I}_{\Sh_{K^p}^{\tor}}^{\la}$ is an isomorphism. 
\end{lem}
\begin{proof} This is a consequence of \cite[Thm.\
  3.4.1.]{camargo2023geometric}.%
\end{proof} 

\subsection{Completed cohomology}\label{subsect-def-compcoh} Let $\Sh_{K^pK_p}^{alg}$ denote the Shimura variety, viewed as a scheme,  defined over its reflex field $E(G,X)$. Let $\Sh_{K^p}^{alg} = \lim_{K_p} \Sh_{K^pK_p}^{alg}$. The limit exists as a scheme since the transition maps are affine.  We define completed cohomology with $\qq_p$  coefficients to be $\mathrm{R}\Gamma_{\mathrm{proet}}(\Sh^{alg}_{K^p, \bar{\qq}}, \qq_p)$. The cohomology groups of 
$\mathrm{R}\Gamma_{\mathrm{proet}}(\Sh_{K^p, \bar{\qq}}^{alg}, \qq_p)$ identify with the usual completed
cohomology groups with~$\Qp$-coefficients (as defined
for example in \cite{MR2905536}). This cohomology has a $G(\qq_p)$-action, an action of the Hecke algebra away from $p$, and an action of $G_{E(G,X)}$, the absolute Galois group of $E(G,X)$.

Using  comparison theorems in \cite[page 30]{MR1734903} and \cite[Cor.\
6.1.7]{camargo2022locally},
 completed cohomology with $\qq_p$-coefficients identifies with  $\mathrm{R}\Gamma_{\proket}(\Sh_{K^p}^{\tor}, \qq_p)$.  Here the subscript $\proket$ is
 is short for ``pro-Kummer-\'etale'', in the sense
 of~\cite{diao2022logarithmic}; this modification of the
 pro-\'etale site is needed because our Shimura varieties are not
 compact. We usually omit this subscript from now on. %
Similarly, we identify $\mathrm{R}\Gamma_{c}(\Sh_{K^p}, \qq_p)$ with the completed
 cohomology with compact support. %

We let $\mathrm{R}\Gamma(\Sh_{K^p}^{\tor}, \qq_p)^{\la}$ be the
(derived) locally analytic vectors in completed cohomology. %
Similarly, we let $\mathrm{R}\Gamma_{c}(\Sh_{K^p}, \qq_p)^{\la}$ be the (derived) locally analytic vectors in completed cohomology with compact support.  
We remark that both $\mathrm{R}\Gamma(\Sh_{K^p}^{\tor}, \qq_p)$ and
$\mathrm{R}\Gamma_c(\Sh_{K^p}, \qq_p)$  have    admissible
cohomology groups, and that  the passage to locally analytic vectors
is an exact functor on admissible representations by \cite{MR1990669},
Thm.\ 7.1 (see also~\cite[Prop.\ 4.48]{MR4475468}). Therefore, the cohomology groups of  $\mathrm{R}\Gamma(\Sh_{K^p}^{\tor}, \qq_p)^{\la}$ and $\mathrm{R}\Gamma_c(\Sh_{K^p}, \qq_p)^{\la}$ are the locally analytic vectors in the  completed cohomology groups.

\begin{thm}%
  \label{thm-prim-compar} We have 
\begin{eqnarray*}
\mathrm{R}\Gamma(\Sh_{K^p}^{\tor}, \qq_p) {\otimes}_{\qq_p} \C_p &=& \mathrm{R}\Gamma_{\an}(\Sh_{K^p}^{\tor}, \oscr_{\Sh_{K^p}^{\tor}})\\
\mathrm{R}\Gamma_c(\Sh_{K^p}, \qq_p){\otimes}_{\qq_p} \C_p &=& \mathrm{R}\Gamma_{\an}(\Sh_{K^p}^{\tor}, \mathscr{I}_{\Sh_{K^p}^{\tor}})
\end{eqnarray*}
We have  
\begin{eqnarray*}
\mathrm{R}\Gamma(\Sh_{K^p}^{\tor}, \qq_p)^{\la}{\otimes}_{\qq_p} \C_p & =& \mathrm{R}\Gamma_{\an}(\Sh_{K^p}^{\tor}, \oscr^{\la}_{\Sh_{K^p}^{\tor}}) \\
\mathrm{R}\Gamma_c(\Sh_{K^p}, \qq_p)^{\la} {\otimes}_{\qq_p} \C_p & =& \mathrm{R}\Gamma_{\an}(\Sh_{K^p}^{\tor}, \mathscr{I}^{\la}_{\Sh_{K^p}^{\tor}})
\end{eqnarray*}
\end{thm}
\begin{proof}
The first part is immediate from  \cite[Thm. \ 6.2.1]{diao2022logarithmic}, by
passing to the  limit as in the proof of \cite[Thm.\ 4.2.1]{scholze-torsion}. %
The second part is \cite[Thm.\ 6.2.6]{camargo2022locally}.
\end{proof}

\begin{rem}\label{rem:ideas-of-proof-Scholze-Pan}The main ideas in the proof of
  Theorem~\ref{thm-prim-compar} are due to Scholze and Pan. More
  precisely, the statements regarding completed cohomology (before taking locally
  analytic vectors) are a consequence of Scholze's primitive
  comparison theorem; see e.g.\ \cite[Thm.\ 4.2.1]{scholze-torsion}. %
For the locally analytic vectors, in the case of usual (i.e.\ not compactly supported) cohomology it
  is a consequence of the fact that
  $\oscr^{\la}_{\Sh_{K^p}^{\tor}}= \oscr^{\Rla}_{\Sh_{K^p}^{\tor}}$
  where $\Rla$ are the derived locally analytic vectors (in the sense
  of \cite{MR4475468}). %
  This was proved for modular curves in
  \cite{MR4390302}, Thm.\ 4.4.6. %
In  \emph{loc. cit.}
it is also proved for modular curves that
  $\mathrm{R}\Gamma_{\an}(\Sh_{K^p}^{\tor},
  \oscr^{\la}_{\Sh_{K^p}^{\tor}})=
  \mathrm{R}\Gamma_{\an}(\mathcal{FL}, (\pi_{HT})_\star
  \oscr^{\la}_{\Sh_{K^p}^{\tor}})$. We do not need (and have not proved) this
  fact more generally.%
\end{rem}

\subsection{The functor $\VB$}  We now introduce a functor which turns equivariant sheaves on the flag variety into sheaves on the perfectoid Shimura variety. 
\subsubsection{Definition of the functor and main properties} 
Let us briefly reintroduce $\Sigma$ to the notation (see Remark \ref{rem-Sigma}).
Let %
$U^{rat}_{\mathcal{FL}}$ be a quasi-compact open subset of $\mathcal{FL}^{rat}$
and let $U_{\mathcal{FL}}$ be its base change to $\Spa(\C_p, \ocal_{\C_p})$.  %
We
write $U_{K^p,\Sigma}:=\pi_{HT,\Sigma}^{-1}(U_{\mathcal{FL}})  =
\lim_{K'_p}  U_{K'_pK^p,\Sigma}$ where $U_{K'_pK^p,\Sigma}$ is a
quasi-compact open subset of $\Sh_{K'_pK^p,\Sigma}$ for $K'_p$ small
enough.  In %
Definition \ref{defn:LBgX} (see also Remark~\ref{rem-defn-LBgG}) we have defined the categories $\LB_{\mathfrak{g}}(U^{rat}_{\mathcal{FL}})$ and
$\LB_{\mathfrak{g}}(U_{\mathcal{FL}})$ and there is a base change functor
$\LB_{\mathfrak{g}}(U^{rat}_{\mathcal{FL}}) \rightarrow
\LB_{\mathfrak{g}}(U_{\mathcal{FL}})$. %
We define a functor
\begin{eqnarray*}
\VBzero_{\Sigma}: \LB_{\mathfrak{g}}(U_\mathcal{FL})& \rightarrow& \mathrm{Mod}( \oscr^{\sm}_{U_{K^p,\Sigma}})  \\
\mathscr{F} & \mapsto & \colim_{K'_p}((\pi_{HT,\Sigma}^{-1}\mathscr{F}) {\otimes}_{\pi_{HT,\Sigma}^{-1}\oscr_{\mathcal{FL}}} \oscr_{U_{K^p,\Sigma}})^{K'_p}.
\end{eqnarray*}
\begin{rem}\label{rem:how-to-define-VBzero}
  Note that this functor is (well-) defined, because for any quasi-compact open
  subset $V_{\mathcal{FL}} \subseteq U_{\mathcal{FL}}$, we have
  $\mathscr{F}(V_{\mathcal{FL}}) = \colim_{K'_p} %
  \mathscr{F}(V_{\mathcal{FL}})_{K'_p}$ where
  $\mathscr{F}(V_{\mathcal{FL}})_{K'_p}$ is a submodule of
  $\mathscr{F}(V_{\mathcal{FL}})$ where the action of $\mathfrak{g}$ integrates
  to an action of $K'_p$ (compatibly for the transition maps).
\end{rem}

We also need a derived version of this functor.  We recall that in Definition
\ref{defn-Modg} we introduced abelian categories
$\Mod_{\mathfrak{g}}(U^{rat}_{\mathcal{FL}})$ and
$\Mod_{\mathfrak{g}}(U_{\mathcal{FL}})$ which respectively contain $\LB_{\mathfrak{g}}(U^{rat}_{\mathcal{FL}})$ and $\LB_{\mathfrak{g}}(U_{\mathcal{FL}})$.  Moreover there is a base change functor $\Mod_{\mathfrak{g}}(U^{rat}_{\mathcal{FL}}) \rightarrow \Mod_{\mathfrak{g}}(U_{\mathcal{FL}})$.  We define a functor  %
\begin{eqnarray*}
\VB_{\Sigma}: D(\Mod_{\mathfrak{g}}(U_{\mathcal{FL}})) & \rightarrow& D(\mathrm{Mod}( \oscr^{\sm}_{U_{K^p,\Sigma}}))  \\
\mathscr{F} & \mapsto & \colim_{K'_p}\mathrm{R}\Gamma(K_p',(\pi_{HT,\Sigma}^{-1}\mathscr{F}) {\otimes}^L_{\pi_{HT,\Sigma}^{-1}\oscr_{\mathcal{FL}}} \oscr_{U_{K^p,\Sigma}}).
\end{eqnarray*}%

\begin{rem}\label{rem:not-writing-piHT-*}
The pullback functor  \[\cF\mapsto(\pi_{HT,\Sigma}^{-1}\mathscr{F})
{\otimes}^{L}_{\pi_{HT,\Sigma}^{-1}\oscr_{\mathcal{FL}}} \oscr_{U_{K^p,\Sigma}}\]
in the definition of $\VB_\Sigma$ is exact on the category
$\LB_\mathfrak{g}(U_\mathcal{FL})$ (essentially by the definition of
$\LB$-sheaves).  Thus on $\LB_{\mathfrak{g}}(U_\mathcal{FL})$, $\VBzero_\Sigma$
is left exact, and we may think of $\VB_\Sigma$ as its right derived functor.

It might be more natural to denote this pullback by~$\pi^{*}_{\HT,\Sigma}$, but
we reserve this notation below for the underived pullback \[\cF\mapsto(\pi_{HT,\Sigma}^{-1}\mathscr{F})
{\otimes}_{\pi_{HT,\Sigma}^{-1}\oscr_{\mathcal{FL}}} \oscr_{U_{K^p,\Sigma}}.\] %
\end{rem}

Since the Hodge--Tate period map is $\mathrm{Gal}(\bar{E}/E)$-equivariant,   if $\mathscr{F}$ is an object of $\Mod_{\mathfrak{g}}(U_{\mathcal{FL}})$ which comes from $\Mod_{\mathfrak{g}}(U^{rat}_{\mathcal{FL}})$ by base change, then $\VB_{\Sigma}(\mathscr{F})$ carries a semi-linear
action of $\mathrm{Gal}(\bar{E}/E)$. 
One can  study this action as follows. We can consider the category $\mathrm{Mod}_{\mathrm{Gal}(\bar{E}/E)}(\oscr^{\sm}_{\Sh^{\tor}_{K^p}})$ of sheaves of $\oscr^{\sm}_{\Sh^{\tor}_{K^p}}$-modules, carrying a semi-linear continuous Galois action. 
We let $\oscr^{\sm, \mathrm{Gal}(\bar{E}/E)-\sm}_{\Sh^{\tor}_{K^p}}$ be the
subsheaf of smooth vectors for the action of  $\mathrm{Gal}(\bar{E}/E)$. %
We
remark that 
$\oscr^{\sm, \mathrm{Gal}(\bar{E}/E)-\sm}_{\Sh^{\tor}_{K^p}} = \oscr^{\sm}_{\Sh^{ rat, \tor}_{K^p}} \otimes_E \bar{E}$.

We define an arithmetic Sen  functor: %
\begin{eqnarray*}
S^{\arit} : \mathrm{Mod}_{\mathrm{Gal}(\bar{E}/E)}(\oscr^{\sm}_{\Sh^{\tor}_{K^p}})& \rightarrow & \mathrm{Mod}(\oscr^{\sm, \mathrm{Gal}(\bar{E}/E)-\sm}_{\Sh^{\tor}_{K^p}}) \\
\mathscr{F} & \rightarrow &  \colim_{E'} \mathscr{F}^{ \mathrm{Gal}(\bar{E}/E'_{\cycl}), \mathrm{Gal}(E'_{\cycl}/E')-an}
\end{eqnarray*}
where the colimit goes over all finite extensions $E'/E$ and the superscript $$(-)^{ \mathrm{Gal}(\bar{E}/E'_{\cycl}), \mathrm{Gal}(E'_{\cycl}/E')-an}$$ means the $\mathrm{Gal}(\bar{E}/E'_{\cycl})$-fixed  and $\mathrm{Gal}(E'_{\cycl}/E')$-analytic vectors (where $\mathrm{Gal}(E'_{\cycl}/E')$ is viewed as a subgroup of $\ZZ_p^\times$ via the cyclotomic character). 
We observe that $S^{\arit}(\mathscr{F})$ carries an $\oscr^{\sm,
  \mathrm{Gal}(\bar{E}/E)-\sm}_{\Sh^{\tor}_{K^p}}$-linear  arithmetic Sen
operator, %
obtained by differentiating the $\mathrm{Gal}(E'_{\cycl}/E')$-action on
$(\mathscr{F})^{ \mathrm{Gal}(\bar{E}/E'_{\cycl}),
  \mathrm{Gal}(E'_{\cycl}/E')-an}$ and passing to the colimit. %

The following theorem is implicit  in  \cite{MR4390302} in the modular curve
case,  see  \cite{PilloniVB} for a formulation in this spirit. In higher
dimension it is essentially a  direct consequence of the results
\cite{camargo2022locally,camargo2023geometric}, as we will see in
the course of the proof. 

\begin{thm}\label{thm-VB-primitive} \leavevmode
\begin{enumerate}
\item For any %
  $\mathscr{F} \in
  \LB_{\mathfrak{g}}(U_{\mathcal{FL}})^{\mathfrak{u}_{P}^0}$,
  there exists a covering by open affinoids $U_{K^p,\Sigma}= \cup_i V_{i}$
  with $V_i = \lim_{K_p} V_{i,K_p}$  (for $K_p$ small enough; here the~$V_{i,K_p}$ are open affinoid) %
  a
  sequence of compact open subgroups $\{K_{p,r}\}_{r\in \Z_{\geq0}}$ and summand of orthonormalizable Banach sheaves  $\VBzero_{\Sigma,K_{p,r},V_i}(\mathscr{F})$ over  $V_{i,K_{p,r}}$  such that  we have: 
$$\VBzero_{\Sigma}(\mathscr{F})\vert_{V_i} = \colim_{r}   \oscr_{V_i}^{\sm} \otimes_{\oscr_{V_{i,K_{p,r}}}} \VBzero_{\Sigma,K_{p,r},V_i}(\mathscr{F}).$$
Moreover,  we have  \numequation\label{eqn:admissibility-VBzero}\VBzero_{\Sigma}(\mathscr{F})
{\otimes}_{\oscr^{\sm}_{U_{K^p,\Sigma}}} \oscr_{U_{K^p,\Sigma}} =
(\pi_{HT,\Sigma}^{-1}\mathscr{F}) {\otimes}_{\pi_{HT,\Sigma}^{-1}(\oscr_{U_{\mathcal{FL}}})}
\oscr_{U_{K^p,\Sigma}}.\end{equation} %
  \item\label{item:VB-exact-uP-0} The  restriction of the functor  $\VBzero_{\Sigma}$ to the category
  $\LB_{\mathfrak{g}}(U_{\mathcal{FL}})^{\mathfrak{u}_{P}^0}$  is an exact functor, and  for an object $\mathscr{F}$ of  $\LB_{\mathfrak{g}}(U_{\mathcal{FL}})^{\mathfrak{u}_{P}^0}$ we have  $$\VB_\Sigma(\mathscr{F}) = \VB^0(\mathscr{F}) \otimes_{\oscr_{U_{K^p,\Sigma}}^{\sm}}( \colim_{K_p}
\oplus_{i=0}^{d}
\Omega^{i}_{U_{K_pK^p,\Sigma}}(\mathrm{log}(D_{K_p,\Sigma}))[-i]).$$

\item\label{item:VB-0-cohomology} Let $\mathscr{F} \in
  \LB_{\mathfrak{g}}(U_{\mathcal{FL}})$. Assume that for all~ $i$, we have $\HH^i(\mathfrak{u}_P^0, \mathscr{F}) \in \LB_{\mathfrak{g}}(U_{\mathcal{FL}})$. 
Then we have isomorphisms: 
$\VBzero_\Sigma(\HH^i(\mathfrak{u}_P^0, \mathscr{F})) \isoto
\HH^i(\VB_\Sigma(\mathscr{F}))$. %

\item Let $\mathscr{F} \in
  \LB_{\mathfrak{g}}(U_{\mathcal{FL}})^{\mathfrak{u}_{P}^0}$. Assume that $\mathscr{F}$ arises from an object of $\LB_{\mathfrak{g}}(U^{rat}_{\mathcal{FL}})$. %
  Then we have $$S^{\arit}(\VBzero_{\Sigma}(\mathscr{F})) \otimes_{\oscr^{\sm, \mathrm{Gal}(\bar{E}/E)-\sm}_{U_{K^p,\Sigma}}}\oscr^{sm }_{U_{K^p,\Sigma}} = \VBzero_\Sigma(\mathscr{F})$$ and the action of $\mu$ via $\Theta_{\hor}$ is an arithmetic Sen
  operator on $S^{\arit}(\VBzero_{\Sigma}(\mathscr{F}))$. %
 
 More precisely, in the notation of $(1)$, we can suppose that the covering $U_{K^p,\Sigma}= \cup_i V_{i}$ comes from a covering $U^{rat}_{K^p,\Sigma}= \cup_i V^{rat}_{i}$, and there exists a finite extension $E_{i,n}$ of $E$ such that 
  $$\VBzero_{\Sigma,K_{p,n},V_i}(\mathscr{F}) = \VBzero_{\Sigma,K_{p,n},V_i}(\mathscr{F})^{\mathrm{Gal}(\bar{E}/E_{i,n,cycl}), \mathrm{Gal}(E_{i,n,cycl}/E_{i,n})-an} \otimes_{E_{i,n}} \C_p.$$

\end{enumerate}
\end{thm}

\begin{proof}%
  We explain why the theorem follows from  the results of \cite{camargo2022locally} and \cite{camargo2023geometric}.  %
The statement is local so we can assume that $\mathscr{F} = \colim
\mathscr{F}_r$ is a colimit of orthonormalizable Banach sheaves with injective
transition maps and that $G_r$ acts on each $\mathscr{F}_r$.  We can assume that
$U_{\mathcal{FL}} = \Spa(C, C^+)$ is affinoid. We can also consider
$V_{K^p,\Sigma} =  \Spa(B, B^+)$ an open affinoid subset of
$U_{K^p,\Sigma}$. %
The open subset $V_{K^p,\Sigma}$ descends to $V_{K^pK_p,\Sigma} =  \Spa(B_{K_p}, B_{K_p}^+)$ for $K_p$ small enough.

For $K_p$ small enough (such that $K_p \subseteq G_r(\qq_p)$), the pull
back %
$\pi_{HT,\Sigma}^\star\mathscr{F}_r$ to $\Sh^{\tor}_{K^p,\Sigma}$ carries a $K_p$-action 
Let us put  $F_r:=\mathscr{F}_r(U_{\mathcal{FL}}) $.  We pick a $C^+$ lattice $F_r^+$ (that is, $F_r^+$ is the completion of a free $C^+$-module and $F_r^+ \otimes_{C^+} C = F_r$). The action of $G_r$ amounts  to a co-action map $ c: F_r \rightarrow F_r \otimes_{C} \oscr_{G_r}$. By continuity, there exists  $n$ such that $c(F_r^+)\subseteq p^{-n} F_r^+ \otimes \oscr^+_{G_r}$. 
 We claim that  for  $r' = r+n+1$  the  restriction of the co-action
 map   $c': F_r \rightarrow F_r \otimes_{C} \oscr_{G_{r'}}$ induces a
 map $F^+_r \rightarrow F^+_r \otimes_{C} \oscr^+_{G_{r'}}$ and
 moreover, this co-action map is trivial modulo $p$. %
To see this,  we may write (for example by using the exponential map) $\oscr_{G_r} = \C_p\langle X_1, \cdots, X_t \rangle$
so that $\oscr_{G_{r'}} = \C_p\langle p^{-n-1}X_1, \cdots, p^{-n-1}X_t
\rangle$. %
For any $f \in F_r^+$, we write $c(f) = \sum_{\underline{i}} f_{\underline{i}}
X^{\underline{i}}$, where $f_{\underline{0}} = f$ and $f_{\underline{i}}  \in p^{-n} F_r^+$ is tending to $0$. Our claim is thus clear. By shrinking $K_p$, we can assume that $K_p \subseteq G_{r'}(\qq_p)$. 
We remark that we have in particular checked that the action of $K_p$  is ``locally analytic''  in the sense of  \cite[Defn.\
1.0.1]{camargo2023geometric}; more precisely, the pro-Kummer-\'etale
$\widehat{\cO}_{V_{K^pK_p,\Sigma}}$-module corresponding to $\pi_{HT,\Sigma}^\star\mathscr{F}_r$ is
relatively analytic ON Banach in the sense of  \cite[Defn.\
1.0.1]{camargo2023geometric}. %
Note that this is the familiar smallness
condition in $p$-adic Simpson theory (see~\cite[Rem.\
1.0.1]{camargo2023geometric} ).

We let:
\begin{itemize}
\item $\mathbf{T} = \Spa(\C_p\langle T, T^{-1} \rangle, \C_p^+\langle T, T^{-1} \rangle)$,
\item $\mathbf{T}_n = \Spa(\C_p\langle T^{p^{-n}}, T^{-p^{-n}} \rangle, \C_p^+\langle T^{p^{-n}}, T^{-p^{-n}} \rangle)$ for any $n\geq 1$.
\item $\mathbf{D} =  \Spa(\C_p\langle T \rangle, \C_p^+\langle T \rangle)$,
\item $\mathbf{D}_n =  \Spa(\C_p\langle T^{p^{-n}} \rangle, \C_p^+\langle T^{p^{-n}} \rangle)$ for any $n\geq 1$. 
\end{itemize}
 
We  can also assume (after shrinking $V_{K^p,\Sigma}$ and  taking $K_p$ small
enough) that  we have a toric chart (in the sense of~\cite[Prop.\ 3.1.10]{diao2022logarithmic}) %
$V_{K_pK^p,\Sigma} \rightarrow  \mathbf{T}^{e} \times \mathbf{D}^{d-e}$.
We let $\Spa(B_{K_p,n}, B_{K_p, n}^+) = U_{K_pK^p,\Sigma} \times_{ \mathbf{T}^{e} \times \mathbf{D}^{d-e}}\mathbf{T}^{e}_n \times \mathbf{D}^{d-e}_n$. We let $\Spa(B_{K_p, \infty}, B_{K_p, \infty}^+) = \lim_n \Spa(B_{K_p,n}, B_{K_p, n}^+)$. We let $\Spa(B_{n}, B_{ n}^+) = U_{K^p,\Sigma} \times_{\mathbf{T}^{e} \times \mathbf{D}^{d-e}}\mathbf{T}^{e}_n \times \mathbf{D}^{d-e}_n$. We let $\Spa(B_\infty, B_\infty^+) = \lim_n \Spa(B_{n}, B_{ n}^+)$. 
After making a choice of compatible $p$-power roots of unity, We let  $\Gamma = \ZZ_p^{d}$, acting on $\mathbf{T}_n^{e} \times
\mathbf{D}_n^{d-e}$. %
We thus have an action of $K_p \times \Gamma$ on
$B_\infty$. By \cite[Prop.\ 3.2.3]{camargo2023geometric}, the triple $(B_\infty, K_p \times \Gamma, pr_2: K_p \times \Gamma \rightarrow \Gamma)$ is a strongly decomposable Sen theory in the sense of \cite[ Defn.\ 2.2.6]{camargo2023geometric}.

We consider the semi-linear representation of $K_p$,  $B \otimes_{C} F_r$. Our
goal is to compute $\colim_{K_p} \HH^i(K_p, B \otimes_{C} F_r)$ using Sen
theory. By almost purity, we have \numequation\label{eqn:almost-purity-swap-K-Gamma}\mathrm{R}\Gamma(K_p, B \otimes_{C} F_r) = \mathrm{R}\Gamma(K_p \times \Gamma, B_\infty \otimes_C F_r)= \mathrm{R}\Gamma( \Gamma, \HH^0(K_p, B_\infty \otimes_C F_r)).\end{equation}
By \cite[Thm.\ 2.4.3]{camargo2023geometric}, we have (after possibly shrinking
$K_p$ and for all $n$ large enough) the Sen module \numequation\label{eqn:Sen-module-defn}S_{K_p, n}(F_r):= (B\otimes_{C}F_r)^{K_p, p^n\Gamma-an}\end{equation} which is obtained by taking the $K_p$-invariants and
the $p^n\Gamma$-analytic vectors.  This is an orthonormalizable Banach
$B_{K_p,n}$-module with a locally analytic action of $\Gamma$,  and it satisfies
\numequation\label{eqn:Sen-Kp-n-decompletes}B_{\infty} \otimes_{B_{K_p,n}}
S_{K_p, n}(F_r) = B_\infty \otimes_{C} F_r.\end{equation}\ %
In addition
by~\eqref{eqn:almost-purity-swap-K-Gamma} we have \numequation\label{eqn:Gamma-invariants-in-Sen-finite-level} \HH^{0}(\Gamma,S_{K_p, n}(F_r))=\HH^0(K_p, B \otimes_{C} F_r).\end{equation}
We have an action of $\mathrm{Lie}(\Gamma)$ on $S_{K_p, n}(F_r)$, which are the
``geometric'' Sen operators, and by \cite[Thm.\ 1.7]{MR4475468} %
we have
\numequation\label{eqn:Gamma-cohomology-from-Sen-module}\mathrm{R}\Gamma( \Gamma, \HH^0(K_p, B_\infty \otimes_C F_r)) =  \HH^0(\Gamma, \mathrm{R}\Gamma(\mathrm{Lie}(\Gamma), S_{K_p, n}(F_r)))\end{equation}
where $\mathrm{R}\Gamma(\mathrm{Lie}(\Gamma), S_{K_p, n}(F_r))$ is a complex of
smooth $\Gamma$-modules and  $\HH^0(\Gamma,-)$ is  the exact functor of
$\Gamma$-invariants on smooth $\Gamma$-modules. 
It is a consequence of \cite[Thm.\ 1.1.5]{camargo2022locally}  that these Sen
operators  are induced by functoriality from the map $\mathfrak{u}_{P}^0 \otimes
\mathscr{F}_r \rightarrow \mathscr{F}_r$, using the identification
$\pi_{HT}^\star\mathfrak{u}_P^0  \simeq \mathrm{Lie}( \Gamma) \otimes
\oscr_{V_{K^p,\Sigma}}$. %
Since the transition maps in the colimit $\cF=\colim\cF_r$ are injective, it
follows in particular that   $\mathscr{F} \in \LB_{\mathfrak{g}}(U_{\mathcal{FL}})^{\mathfrak{u}_{P}^0}$ if and only if the geometric Sen operator  of each $\mathscr{F}_r$ is trivial. 
Let us assume that this is the case. 
Then the action of $\Gamma$ on~$S_{K_p, n}(F_r)$  is smooth, and this action
factors through $\Gamma/p^n\Gamma$. By finite \'etale descent, we find that
\numequation\label{eqn:finite-etale-descent} B_{K_p,n} \otimes_{B_{K_p}}S_{K_p, n}(F_r)^{\Gamma} = S_{K_p,n}(F_r)\end{equation} and that
$S_{K_p,n}(F_r)^\Gamma$ is a direct summand of the orthonormalizable Banach
$B_{K_p}$-module $S_{K_p,n}(F_r)$. Taking $p^{n}\Gamma$-invariants
in~\eqref{eqn:Sen-Kp-n-decompletes}, we obtain \[B_{n} \otimes_{B_{K_p,n}}
S_{K_p, n}(F_r) = B_n\otimes_{C} F_r\]and thus
(using~\eqref{eqn:finite-etale-descent})\[B_{n}  \otimes_{B_{K_p}}S_{K_p,
    n}(F_r)^{\Gamma} = B_n\otimes_{C} F_r.\]Taking $\Gamma$-invariants and
using~\eqref{eqn:Gamma-invariants-in-Sen-finite-level}, 
we deduce that \numequation\label{eqn:B-descent-H0KpFr}B\otimes_{B_{K_{p}}}\HH^0(K_p, B \otimes_{C} F_r) = B \otimes_{C} F_r.\end{equation} 
 Passing to the colimit over $r$ and $K_p$, we obtain~\eqref{eqn:admissibility-VBzero}. This completes the
 proof of~(1) (taking $\VBzero_{\Sigma,K_{p,r},V_i}(\mathscr{F})$ to be the
 sheaf associated to
$\HH^0(K_p, B \otimes_{C} F_r) $, and~$K_{p,r}$ to be~$K_p$). %

We now prove~(2).  We have  $$\mathrm{R}\Gamma(\mathrm{Lie}(\Gamma), S_{K_p, n}(F_r)) = \oplus S_{K_p, n}(F_r) \otimes \Lambda^i (\mathrm{Lie}(\Gamma))^\vee[-i]$$ from which we deduce that \numequation\label{eqn:VB-factors-into-VB0}\VB_\Sigma(\mathscr{F}) = \VB^0(\mathscr{F}) \otimes_{\oscr_{U_{K^p,\Sigma}}^{\sm}}( \colim_{K_p}
\oplus_{i=0}^{d}
\Omega^{i}_{U_{K_pK^p,\Sigma}}(\mathrm{log}(D_{K_p,\Sigma}))[-i]).\end{equation}

Let $0 \rightarrow \mathscr{F} \rightarrow \mathscr{L} \rightarrow \mathscr{H} \rightarrow 0$ be an exact sequence in $\LB_{\mathfrak{g}}(U_{\mathcal{FL}})^{\mathfrak{u}_{P}^0}$. Applying $\VB_\Sigma$ yields exact sequences:
\[0 \rightarrow \VBzero_\Sigma(\mathscr{F}) \rightarrow
\VBzero_\Sigma(\mathscr{L}) \rightarrow \VBzero_\Sigma(\mathscr{H})\] and
\[ \VBzero_\Sigma(\mathscr{F})\otimes_{\oscr_{U}^{\sm}} \VBzero_\Sigma(
  \Lambda^d (\mathfrak{u}_P^0)^\vee) \rightarrow
  \VBzero_\Sigma(\mathscr{L})\otimes_{\oscr_{U}^{\sm}} \VBzero_\Sigma(
  \Lambda^d(\mathfrak{u}_P^0)^\vee) \rightarrow \VBzero_\Sigma(\mathscr{H})
  \otimes_{\oscr_{U}^{\sm}} \VBzero_\Sigma( \Lambda^d(\mathfrak{u}_P^0)^\vee)
  \rightarrow 0.\] Since  $\VBzero_\Sigma( \Lambda^d (\mathfrak{u}_P^0)^\vee)$
is an invertible sheaf, we conclude that $0 \rightarrow
\VBzero_\Sigma(\mathscr{F}) \rightarrow \VBzero_\Sigma(\mathscr{L}) \rightarrow
\VBzero_\Sigma(\mathscr{H}) \rightarrow 0$ is exact, as required. %

We now turn to~(3), so we no longer assume that $\mathscr{F}$ is killed by
$\mathfrak{u}_P^0$.   We let $S_{K_p}(F_r) = \colim_{n} S_{K_p,n}(F_r)$, we let
$S(F_r) = \colim_{K_p} S_{K_p}(F_r)$ and finally we let $S(F) = \colim_{r}
S(F_r)$. 
We claim that
\numequation\label{eqn:Lie-Gamma-cohomology-uP0}\HH^i(\mathrm{Lie}(\Gamma),
S(F)) = S(\HH^i(\mathfrak{u}_P^0, F)). \end{equation}
Granting~\eqref{eqn:Lie-Gamma-cohomology-uP0}, we claim that taking
$\Gamma$-invariants  gives %
\numequation\label{eqn:Gamma-invariants-uP0}\colim_{r,K_{p}}\HH^i(K_p, B \otimes_C F_{r}) =
\colim_{r,K_{p}} \HH^0(K_p, B\otimes_{C}\HH^i(\mathfrak{u}_P^0,
F_{r})),\end{equation}which immediately gives~(3). Indeed,
by~\ref{eqn:Gamma-cohomology-from-Sen-module} and~\eqref{eqn:almost-purity-swap-K-Gamma} we have\[H^0(\Gamma,\HH^i(\mathrm{Lie}(\Gamma),
S_{K_p,n}(F_r)))=\HH^i(\Gamma, \HH^0(K_p, B_\infty \otimes_C
F_r))=\HH^{i}(K_p,B\otimes_CF_r).\]On the other hand, passing to  colimits
in~\eqref{eqn:Gamma-invariants-in-Sen-finite-level} we see that \numequation\label{eqn:Gamma-invariants-in-Sen-infinite-level}\colim_{K_p}
\HH^0(K_p, B \otimes_C F_r) = \HH^0(\Gamma, S(F_r)),\end{equation}
so that (replacing~$F_r$ by $\HH^i(\mathfrak{u}_P^0,
F_{r})$)
\[H^0(\Gamma,S(\HH^i(\mathfrak{u}_P^0, F)))=\colim_{r,K_{p}}\HH^0(K_p, B
  \otimes_C \HH^i(\mathfrak{u}_P^0,
F_{r})), \]as required.

We now establish~\eqref{eqn:Lie-Gamma-cohomology-uP0}. Firstly, we claim that
$B_{K_p, \infty}$ is an orthonormalizable $B_{K_p,n}$-module. Indeed the algebra
of   $\mathbf{T}_\infty = \lim_n \mathbf{T}_n $ has a topological basis $\{T^{i}\}_{i \in \qq_p/\ZZ_p}$ over $\C_p\langle T, T^{-1} \rangle$, and similarly the algebra of $\mathbb{D}_{\infty}\lim_n \mathbf{D}_n$ has a topological basis $\{T^{i}\}_{i \in \qq_p/\ZZ_p}$ over $\C_p \langle T \rangle$. 
We next claim that  $B_\infty$ is a direct summand of an orthonormalizable
$B_{K_p, \infty}$-module. %
To see this, 
let us fix a decreasing  sequence of compact open subgroups $\{K_{p,r}\}_{r \geq
  0}$  tending to $\{e\}$ with $K_{p,0}=K_p$. Since $B^+_{K_{p,r}, \infty}
\rightarrow B^+_{K_{p,r+1}, \infty}$ is almost \'etale, there exist  finite
$B^+_{K_{p,r}, \infty}$ modules  $X_{K_{p,r+1}}$ and $Y_{K_{p,r+1}}$ together
with:   \begin{itemize}%
\item an injective  map  $ B^+_{K_{p,r},   \infty} \oplus X_{K_{p,r+1}} \rightarrow B^+_{K_{p,r+1}, \infty} $ whose cokernel is annihilated by $p$, and
\item an integer $n_r \in \ZZ_{\geq 0}$ and  an injective map $X_{K_{p,r+1}} \oplus Y_{K_{p,r+1}} \rightarrow (B^+_{K_{p},r,  \infty})^{n_r}$ whose cokernel is annihilated by $p$. 
\end{itemize}
We deduce that $B_{\infty} \bigoplus (\widehat{\oplus_{r} Y_{K_{p,r+1}}}[1/p])$
is orthonormalizable over $B_{K_p, \infty}$, as required.

In particular, we have shown that $B_\infty$ is flat over $B_{K_p,n}$. 
We deduce that $\HH^i(\mathrm{Lie}(\Gamma),B_\infty \otimes_{B_{K_p,n}} S_{K_p,n}(F_r) ) = \HH^i(\mathrm{Lie}(\Gamma),B_\infty \otimes_{B_{K_p,n}} S_{K_p,n}(F_r))$, and passing to the colimit over $K_p,n,r$ we obtain that (for $B_\infty^{\sm}$ the subring of smooth vectors for the $K_p\times \Gamma$-action):  
$$\HH^i(\mathrm{Lie}(\Gamma),B_\infty\otimes_{B^{\sm}_{\infty}} S(F) ) = B_\infty\otimes_{B^{\sm}_{\infty}}\HH^i(\mathrm{Lie}(\Gamma), S(F) ).$$
On the other hand, $\HH^i(\mathfrak{u}_P^0, B_\infty \otimes_C F) =
B_{\infty}\otimes_{C}\HH^i(\mathfrak{u}_P^0, F) $ since  all the $\HH^j(\mathfrak{u}_P^0,
F)$ are  flat over~$C$ %
(here we use our assumption that  $\HH^i(\mathfrak{u}_P^0, \mathscr{F}) \in \LB_{\mathfrak{g}}(U_{\mathcal{FL}})$). 
Recalling~\eqref{eqn:Sen-Kp-n-decompletes}, we deduce that we have $K_p \times \Gamma$- equivariant isomorphisms: 
 $$B_\infty\otimes_{B^{\sm}_{\infty}}\HH^i(\mathrm{Lie}(\Gamma), S(F) )= B_{\infty}\otimes_{C}\HH^i(\mathfrak{u}_P^0, F) =B_\infty\otimes_{B^{\sm}_{\infty}} S(\HH^i(\mathfrak{u}_P^0, F) ) $$
 Taking $K_p \times \Gamma$ smooth vectors yields $\HH^i(\mathrm{Lie}(\Gamma),
 S(F)) = S(\HH^i(\mathfrak{u}_P^0, F))$, as required. 

We now prove the last point. We take $\mathscr{F} \in
\LB_{\mathfrak{g}}(U_{\mathcal{FL}})^{\mathfrak{u}_{P}^0}$, arising from an
object of~$\LB_{\mathfrak{g}}(U^{rat}_{\mathcal{FL}})$. %
We can therefore choose the~$\cF_r$ to be defined over $E$,  so that there is a semi-linear $\mathrm{Gal}(\bar{E}/E)$-action on $B \otimes_C
F_r$. %
Moreover, taking a topological  basis $\{v_i\}_{i\in I}$ of $F_r$ defined over $E$, we see that 
the Galois action is trivial in this basis. Assuming that $K_p$ acts trivially
modulo $p^2$ on $F_r^+$ (which we can always arrange after shrinking $K_p$), we
deduce that $S_{K_p,n}(F_r)$ has a topological  basis $\{v'_i\}_{i \in I}$ where
the change of basis matrix (for the isomorphism~\eqref{eqn:Sen-Kp-n-decompletes})  from $\{v_i\}$ to $\{v'_i\}$ is congruent to $1$
modulo $p$ %
(see \cite[Thm.\ 2.4.3, (1),
(b)]{camargo2023geometric}). %
As a result, the
matrix of the Galois action  on $S_{K_p,n}(F_r)$  is congruent to $1$ modulo
$p$. One can therefore apply  Sen theory to the extension $B^{rat}_{K_p,n}
\rightarrow B_{K_p,n} $ %
where $B_{K_p,n}^{rat} =
(B_{K_p,n})^{\mathrm{Gal}{\bar{E}/E(\zeta_{p^n})}}$. We let %
$S^{\arit, s}(S_{K_p,n}(F_r)) = S_{K_p,n}(F_r)^{G_{L^{\cycl}},
  G_{L(\zeta_{p^s})}-an}$. For $s$ large enough, $S^{\arit, s}(S_{K_p,n}(F_r))
\otimes_{B^{rat}_{K_p,n}(\zeta_{p^s})} B_{K_p,n} = S_{K_p,n}(F_r)$ and the
derivative of the $\mathrm{Gal}(E_{\cycl}/E(\zeta_{p^s}))$ action provides an
arithmetic Sen operator.

In order to prove~(4), it only remains to identify this Sen operator with the operator $\mu$ coming from the horizontal action. 
The orbit map provides an embedding $F_r \into\mathcal{C}^{an}(K_p, F_r)$, $f \mapsto [k \mapsto k.f]$. It intertwines the action of $K_p$ on $F_r$ with the action of $K_p$ on functions $h(-) \in \mathcal{C}^{an}(K_p, F_r)$ via $k\star_2 h(-) = h(-k)$ (therefore the $\star_2$ action does not depend on the $K_p$-action on $F_r$). Note that $\mathcal{C}^{an}(K_p, F_r)$ is a $C$-module as there is a orbit map $C \rightarrow \mathcal{C}^{an}(K_p, C)$ and $\mathcal{C}^{an}(K_p, C)$ acts on $\mathcal{C}^{an}(K_p, F_r)$ naturally.
Moreover, the embedding $F_r \into\mathcal{C}^{an}(K_p, F_r)$ factors through
$\mathcal{C}^{an}(K_p, F_r)^{\mathfrak{u}_{P}^0}$. It therefore suffices to identify the arithmetic Sen operator of  $\mathcal{C}^{an}(K_p, F_r)^{\mathfrak{u}_{P}^0}$. Since $F_r$ has a topological basis over $C$ we can reduce to understanding the Sen operator of $\mathcal{C}^{an}(K_p, C)^{\mathfrak{u}_{P}^0}$. 
The orbit map $C \mapsto \mathcal{C}^{an}(K_p, C)$ induces an isomorphism
$\mathcal{C}^{an}(K_p, \qq_p) \otimes_{\qq_p} C \rightarrow
\mathcal{C}^{an}(K_p, C)$. Moreover, the subspace of algebraic functions
$\mathcal{C}^{alg}(K_p,\qq_p) \hookrightarrow \mathcal{C}^{an}(K_p, \qq_p)$ is
dense, and it induces a dense map %
$(\mathcal{C}^{alg}(K_p, \qq_p) \otimes_{\qq_p} C)^{\mathfrak{u}_{P}^0} \hookrightarrow (\mathcal{C}^{an}(K_p, \qq_p) \otimes_{\qq_p} C)^{\mathfrak{u}_{P}^0}$. Viewing the arithmetic Sen operator as an endomorphism of $$B_\infty \otimes_{B_{K_p,n}}S_{K_p,n}(\mathcal{C}^{an}(K_p, \qq_p) \otimes_{\qq_p} C)^{\mathfrak{u}_{P}^0})  = B_\infty\otimes_{C}(\mathcal{C}^{an}(K_p, \qq_p) \otimes_{\qq_p} C)^{\mathfrak{u}_{P}^0}$$ we deduce that it suffices to prove it coincides with $\mu$ on $B_{\infty}\otimes_{C}(\mathcal{C}^{alg}(K_p, \qq_p) \otimes_{\qq_p} C)^{\mathfrak{u}_{P}^0}$. Since $\mathcal{C}^{alg}(K_p, \qq_p)  = \oplus_{\kappa \in X^\star(T)^+} V_\kappa \otimes V_\kappa^\vee$ and 
$(V_\kappa \otimes V_\kappa^\vee \otimes_{\qq_p}
\oscr_{\mathcal{FL}})^{\mathfrak{u}_P^0} = \mathcal{L}_\kappa \otimes_{\qq_p}
V^\vee_\kappa$ we deduce that it suffices to understand the Sen operator of  the
classical automorphic vector bundles. It therefore suffices to show that
$\VBzero_{\Sigma}(\mathcal{L}_\kappa) = \omega^{\kappa,sm}( \kappa(\mu))$, where
$\kappa(\mu)$ is the Tate twist and $\omega^{\kappa, \sm} = \colim_{K_p}
\omega^\kappa_{K_p}$ is the colimit of the automorphic vector bundles defined
over $\Sh^{\tor}_{K^pK_p}$ (equipped with their rational structure). %
This
follows from an inspection of the rationality properties of the Hodge--Tate map
and of the universal $M$-torsor, see
\cite[Thm.\ 4.2.1]{camargo2022locally}. %
\end{proof}
Observe that by~\eqref{eqn:VB-factors-into-VB0} %
$$\VB_\Sigma(\oscr_{U_\mathcal{FL}}) = \colim_{K_p}
\oplus_{i=0}^{d}
\Omega^{i}_{U_{K_pK^p,\Sigma}}(\mathrm{log}(D_{K_p,\Sigma}))[-i]$$ is a DG
algebra %
and admits an augmentation map to
$\oscr_{U_{K^p,\Sigma}}^{\sm}[0]$. We can therefore make the following definition.

\begin{defn}\label{defn:VBred}%
  We let $\VBred_{\Sigma}: D(\Mod_{\mathfrak{g}}(U_{\mathcal{FL}})) \to
  D(\mathrm{Mod}( \oscr^{\sm}_{U_{K^p,\Sigma}})) $ be the functor given by \[\VBred_{\Sigma}( \mathscr{F}) := \VB_\Sigma(\mathscr{F}) \otimes^{L}_{\VB_\Sigma(  \oscr_{U_\mathcal{FL}})} \oscr_{U}^{\sm}[0].\]
\end{defn}

\begin{rem}\label{rem:apply-Vbzero-termwise} %
By Theorem~\ref{thm-VB-primitive}~\eqref{item:VB-exact-uP-0}, if $\mathscr{F}
\in \LB_{\mathfrak{g}}(U_{\mathcal{FL}})^{\mathfrak{u}_{P}^0}$, then
$\VBred_\Sigma(\mathscr{F}) = \VBzero_{\Sigma}(\mathscr{F})[0]$. Consequently
given a complex of objects in
$\LB_{\mathfrak{g}}(U_{\mathcal{FL}})^{\mathfrak{u}_{P}^0}$, we
can evaluate~$\VBred_{\Sigma}$ by applying $\VBzero_{\Sigma}$ termwise.
\end{rem}

\subsubsection{Variants} We now introduce variants of the above
functor carrying extra structure. Let us define $\mathrm{Mod}_{G(\qq_p)} ( \oscr^{\sm}_{\Sh^{\tor}_{K^p}})$ to be the category of sheaves $(\mathscr{F}_{\Sigma})_{\Sigma}$ of $\oscr^{\sm}_{\Sh^{\tor}_{K^p,\Sigma}}$-modules with the properties:
\begin{enumerate}
\item For any refinement $\Sigma'$ of  $\Sigma$, inducing a map $\pi_{\Sigma', \Sigma} : \Sh^{\tor}_{K^p, \Sigma'} \rightarrow \Sh^{\tor}_{K^p, \Sigma}$, we have an  isomorphism  $\pi_{\Sigma', \Sigma}^\star \mathscr{F}_{\Sigma} \rightarrow  \mathscr{F}_{\Sigma'} $ of $\oscr^{\sm}_{\Sh^{\tor}_{K^p,\Sigma'}}$-modules (and these isomorphisms are compatible).
\item For any $g \in G(\qq_p)$, inducing an isomorphism $g : \Sh^{\tor}_{K^p, \Sigma} \rightarrow \Sh^{\tor}_{K^p, g\Sigma}$, there is an isomorphism  $g^\star \mathscr{F}_{g\Sigma} \rightarrow \mathscr{F}_{\Sigma}$ of $\oscr^{\sm}_{\Sh^{\tor}_{K_p, \Sigma}}$-modules (and they satisfy the usual cocycle condition). 
\end{enumerate}
Then we have a functor  \[\VBzero:
\LB_{(\mathfrak{g},G)}( \mathcal{FL})  \rightarrow
\mathrm{Mod}_{G(\qq_p)}( \oscr^{\sm}_{\Sh^{\tor}_{K^p}}),\]constructed
as follows. Composing the functor~$\VBzero_{\Sigma}$
with the forgetful functor $\LB_{(\mathfrak{g},G)}( \mathcal{FL})\to
\LB_{\mathfrak{g}}( \mathcal{FL})$ gives a functor $\VBzero_{\Sigma}:
\LB_{(\mathfrak{g},G)}( \mathcal{FL})\to\mathrm{Mod} ( \oscr^{\sm}_{\Sh^{\tor}_{K^p}})$.  %
For each $g \in G(\qq_p)$,  there is  a map $ g:\Sh^{\tor}_{K^p,\Sigma} \rightarrow \Sh^{\tor}_{K^p, g \Sigma}$ and a map $g^\star \VBzero_{g\Sigma}(\mathscr{F}) \rightarrow \VBzero_{\Sigma}(\mathscr{F})$ satisfying the usual cocycle  condition.  
The various $\VBzero_{\Sigma}$ thus define a functor $\VBzero:
\LB_{(\mathfrak{g},G)}( \mathcal{FL})  \rightarrow
\mathrm{Mod}_{G(\qq_p)}( \oscr^{\sm}_{\Sh^{\tor}_{K^p}})$ as claimed. Note that
in practice, we fix some $\Sigma$ and really work with the functor
$\VBzero_{\Sigma}$ (but see Remark \ref{rem-Sigma} for our notational convention).

We recall the stratification into $B$-orbits  $\mathcal{FL} = \coprod_{w \in \WM} C_w$, with $C_{w} = P \backslash P {w} B$. 
We let $j_{w}: C_{w} \hookrightarrow \mathcal{FL}$ be the locally closed immersion. It induces $j_{w, \Sh^{\tor}_{K^p,\Sigma}}: \pi_{HT,\Sigma}^{-1}(C_{w}) \rightarrow \Sh^{\tor}_{K^p,\Sigma}$.   
Instead of working on the whole Shimura variety, we can
also work over $\pi_{HT,\Sigma}^{-1} (C_w^\dag)$ for any $w \in \WM$. We recall that this is a ringed space, whose underlying topological space is $\pi_{HT,\Sigma}^{-1} (C_w)$ and whose structure sheaf is $j^{-1}_{w, \Sh^{\tor}_{K^p,\Sigma}}\oscr_{\Sh^{\tor}_{K^p,\Sigma}}$. In this case, we can  consider  functors: 
$$\VBzero_{\Sigma}: \LB_{(\mathfrak{g},B)}( C_w^\dag)
\rightarrow \mathrm{Mod}( \oscr^{\sm}_{\pi_{HT,\Sigma}^{-1}
  (C_w^\dag)}).$$

 If $b \in B(\qq_p)$,  there is  a map $ b:  \pi_{HT,\Sigma}^{-1} C_w^\dag \rightarrow \pi_{HT,b\Sigma}^{-1} C_w^\dag$ and a map $b^\star \VBzero_{b\Sigma}(\mathscr{F}) \rightarrow \VBzero_{\Sigma}(\mathscr{F})$ satisfying the usual cocycle  condition. This leads us to consider the category $\mathrm{Mod}_{B(\qq_p)}( \oscr^{\sm}_{\pi_{HT}^{-1} C_w^\dag})$ whose objects are collections of sheaves $(\mathscr{F}_{\Sigma})_{\Sigma}$ of $\oscr^{\sm}_{\pi_{HT,\Sigma}^{-1} C_w^\dag}$-modules, such that for any refinement $\Sigma'$ of  $\Sigma$, inducing a map $\pi_{\Sigma', \Sigma} : \Sh^{\tor}_{K^p, \Sigma'} \rightarrow \Sh^{\tor}_{K^p, \Sigma}$, we have an  isomorphism  $\pi_{\Sigma', \Sigma}^\star \mathscr{F}_{\Sigma} \rightarrow  \mathscr{F}_{\Sigma'} $ of $\oscr^{\sm}_{\pi_{HT,\Sigma'}^{-1} C_w^\dag}$-modules (and these isomorphisms are compatible),
  and such that for any  $b \in B(\qq_p)$, there is a map  $b^\star \mathscr{F}_{b\Sigma} \rightarrow \mathscr{F}_{\Sigma}$ satisfying the usual cocycle condition. 
The various $\VBzero_{\Sigma}$ thus define a functor $\VBzero:
\LB_{(\mathfrak{g},B)}( C_w^\dag)^{\mathfrak{u}_P^0} \rightarrow
\mathrm{Mod}_{B(\qq_p)}( \oscr^{\sm}_{\pi_{HT}^{-1}
  C_w^\dag})$. Again,  in practice, we fix some $\Sigma$ and really
work with the functor $\VBzero_{\Sigma}$ (but drop ~$\Sigma$ from the notation). %

We also remark that we have  $E$-rational structures $\mathcal{FL}^{rat}$ on $\mathcal{FL}$  and $C^{rat}_w$ on $C_w$, and we can consider the categories $\LB_{(\mathfrak{g},G)}( \mathcal{FL}^{rat})$ and $ \LB_{(\mathfrak{g},B)}( C_w^{rat,\dag})$ which admit base change functors to the categories $\LB_{(\mathfrak{g},G)}( \mathcal{FL})$ and $ \LB_{(\mathfrak{g},B)}( C_w^{\dag})$.

\begin{thm}\label{thm-VB} \leavevmode
\begin{enumerate}
\item The functor $$\VBzero: \LB_{(\mathfrak{g},G)}( \mathcal{FL})^{\mathfrak{u}_P^0} \rightarrow \mathrm{Mod}_{G(\qq_p)}( \oscr^{\sm}_{\Sh^{\tor}_{K^p}})$$
is an exact functor.
\item For any $\mathscr{F} \in \LB_{(\mathfrak{g},G)}( \mathcal{FL})^{\mathfrak{u}_P^0}$, we have an analytic covering $\Sh^{\tor}_{K^p} = \cup V_i$, a sequence of compact open subgroups $K_{p,n}$ and summand of orthonormalizable Banach sheaves  $\VBzero_{K_{p,n}, V_i}(\mathscr{F})$ over $V_{i,K_{p,n}}$ such that 
$$\VB^0(\mathscr{F}) \vert_{V_i} = \colim_n \VB^0_{K_{p,n}, V_i}(\mathscr{F}) \otimes_{\oscr_{V_{i,K_{p,n}}}} \oscr_{V_i}^{\sm}.$$
Moreover, there is a compact open subgroup $K_p$ fixing $V_i$ such that all sheaves $\VB^0_{K_{p,n}, V_i}(\mathscr{F}) \otimes_{\oscr_{V_{i,K_{p,n}}}} \oscr_{V_i}^{\sm}$ are $K_p$-equivariant (compatibly with $n$) and this induces in the limit the $K_p$-equivariant structure on $\VB^0(\mathscr{F})\vert_{V_i}$. 
\item For any $\mathscr{F} \in \LB_{(\mathfrak{g},G)}(
  \mathcal{FL})^{\mathfrak{u}_P^0}$,  we have $$\VB^0(\mathscr{F})
  {\otimes}_{\oscr^{\sm}_{\Sh^{\tor}_{K^p}}} \oscr_{\Sh^{\tor}_{K^p}}
  = \pi_{HT}^{-1} \mathscr{F} {\otimes}_{\pi_{HT}^{-1}\oscr_{\mathcal{FL}}}
  \oscr_{\Sh^{\tor}_{K^p}}.$$ %
\item\label{irem:VBofClaisOla} We have that $\VBzero(\mathcal{C}^{\la}) = \oscr_{\Sh^{\tor}_{K^p}}^{\la}$. 
\item The functor $$\VBzero: \LB_{(\mathfrak{g},B)}( C_w^\dag)^{\mathfrak{u}_P^0} \rightarrow \mathrm{Mod}_{B(\qq_p)}( \oscr^{\sm}_{\pi_{HT}^{-1} C_w^\dag})$$
is an exact functor.
\item For any $\mathscr{F} \in \LB_{(\mathfrak{g},B)}( C_w^\dag)^{\mathfrak{u}_P^0}$, there exists  an analytic covering by quasi-compact subsets $\pi_{HT}^{-1}(C_{w}) = \cup_i V_i$, a cofinal  decreasing  family of  quasi-compact  strict neighborhoods of $V_i$:  $V_{i,n} = \lim_{K_p} V_{i,n, K_p}$,  compact open subgroups $K_{p,n}$,  and summand of orthonormalizable Banach sheaves $\VBzero_{K_{p,n}, V_{i,n}}(\mathscr{F})$ over $V_{i,n, K_{p,n}}$  such that 
$$\VBzero(\mathscr{F})\vert_{V_i} = \colim_{n} \VBzero_{K_{p,n}, V_{i,n}}(\mathscr{F}) \otimes_{\oscr_{V_{i,n,K_{p,n}}}} \oscr_{V_{i,n}^{\sm}}$$
Moreover, for each $i$, there is a compact open subgroup $K_B \subseteq B(\qq_p)$ stabilizing $V_i$ and all $V_{i,n}$ and such that the sheaves  $\VBzero_{K_{p,n}, V_{i,n}}(\mathscr{F}) \otimes_{\oscr_{V_{i,n,K_{p,n}}}} \oscr_{V_{i,n}^{\sm}}$ are $K_B$-equivariant, compatibly in $n$, and this induces in the limit the  $K_B$-equivariant structure on $\VBzero(\mathscr{F})\vert_{V_i}$. 
\item   For any $\mathscr{F} \in \LB_{(\mathfrak{g},B)}( C_w^\dag)^{\mathfrak{u}_P^0}$, we have that $$\VBzero(\mathscr{F}) {\otimes}_{\oscr^{\sm}_{\pi_{HT}^{-1} C_w^\dag} } \oscr_{\pi_{HT}^{-1} C_w^\dag} = \pi_{\HT}^{-1}\mathscr{F} {\otimes}_{\pi_{\HT}^{-1}\oscr_{C_w^\dag}} \oscr_{\pi_{HT}^{-1} C_w^\dag}.$$
\item The following diagram of functors is commutative (where the horizontal
  functors are given by~$\VBzero$, and the vertical functors are the natural
  restriction functors): 
\begin{eqnarray*}
\xymatrix{ \LB_{(\mathfrak{g},G)}( \mathcal{FL})^{\mathfrak{u}_P^0} \ar[r] \ar[d] & \mathrm{Mod}_{G(\qq_p)}( \oscr^{\sm}_{\Sh^{\tor}_{K^p}}) \ar[d] \\ 
\LB_{(\mathfrak{g},B)}( C_w^\dag)^{\mathfrak{u}_P^0} \ar[r] & \mathrm{Mod}_{B(\qq_p)}( \oscr^{\sm}_{\pi_{HT}^{-1} C_w^\dag})}
\end{eqnarray*}
\item If $\mathscr{F}$ arises from $\LB_{(\mathfrak{g},G)}( \mathcal{FL}^{rat})^{\mathfrak{u}_P^0}$ or $ \LB_{(\mathfrak{g},B)}( C_w^{rat,\dag})^{\mathfrak{u}_P^0}$, the action of $\mu$ via $\Theta_{\hor}$ is an arithmetic Sen operator on $\VBzero(\mathscr{F})$. 
\end{enumerate}
\end{thm}
\begin{proof}Everything but part \eqref{irem:VBofClaisOla} is an
  immediate consequence of Theorem \ref{thm-VB-primitive}. To
  see~\eqref{irem:VBofClaisOla}, note firstly that since $\RGamma(\mathfrak{u}_P,\cO_{G,e})=\cO_{U_{P_{w}}\backslash G,e}[0]$
 (\emph{cf.}\ Lemma \ref{lem:OGe-is-acyclic}), it follows from Theorem
 \ref{thm-VB-primitive}~\eqref{item:VB-0-cohomology} that  $\VB(\oscr_{G,e}
 \otimes \oscr_{\mathcal{FL}}) = \VBzero( \mathcal{C}^{\la})[0]$. %
  On the other hand, bearing in mind Remark~\ref{rem:how-to-define-VBzero}, we see that \[
 \VB(\oscr_{G,e} \otimes
 \oscr_{\mathcal{FL}})=\oscr_{\Sh^{\tor}_{K^p}}^{\Rla},\]where~$\Rla$ is the
functor of derived locally analytic vectors defined in~\cite[Defn.\ 4.40]{MR4475468}. The result follows immediately.    %
\end{proof} 
\begin{rem}
  \label{rem:O-Rla} In particular the proof of Theorem~\ref{thm-VB} showed that 
 $\oscr^{\la}_{\Sh_{K^p}^{\tor}}= \oscr^{\Rla}_{\Sh_{K^p}^{\tor}}$, confirming
 Remark~\ref{rem:ideas-of-proof-Scholze-Pan}.
\end{rem}

\begin{prop}\label{prop-smooth-lambdasmooth}\leavevmode \begin{enumerate}
\item  Assume that $\mathscr{F} \in \LB_{(\mathfrak{g},G)}(
  \mathcal{FL})^{\mathfrak{u}_P^0}$ is such that the
  $\mathfrak{g}$-action is the derivative of the $G$-action. Then $\VBzero(\mathscr{F}) \in \Mod^{\sm}_{G(\qq_p)}(\oscr_{{\Sh}^{\tor}_{K^p}})$.
\item  Assume that $\mathscr{F} \in \LB_{(\mathfrak{g},B)}(
  C_w^\dag)^{\mathfrak{u}_P^0}$ is such that the restriction to $\mathfrak{b}$
  of the  $\mathfrak{g}$-action  is the derivative of the $B$-action. Then $\VBzero(\mathscr{F}) \in \Mod^{\sm}_{B(\qq_p)}(\oscr^{\sm}_{\pi_{HT}^{-1} C_w^\dag})$.
\end{enumerate}
\end{prop}
\begin{proof}This is immediate from the definition of~$\VB_{\Sigma}$
  (and of the $G(\Qp)$-action).
 \end{proof}
 
 We now address the existence of an arithmetic  Sen operator on the locally analytic vectors in completed cohomology. One can consider the category 
 $\Mod_{\mathrm{Gal}(\bar{E}/E)}(\C_p)$ of semi-linear $\C_p$-representations and define a Sen module functor
 \begin{eqnarray*}
S^{\arit} : \mathrm{Mod}_{\mathrm{Gal}(\bar{E}/E)}(\C_p)& \rightarrow & \mathrm{Mod}(\bar{E}) \\
V & \rightarrow &  \colim_{E'} (V)^{ \mathrm{Gal}(\bar{E}/E'_{\cycl}), \mathrm{Gal}(E'_{\cycl}/E')-an}
\end{eqnarray*}
 We still denote by  $S^{\arit}$ the right derived functor. The following is
 \cite[Cor.\ 6.3.6]{camargo2022locally}. %
\begin{thm}\label{thm-existence-sen-oncc}\leavevmode
\begin{enumerate}
\item  We have that 
\begin{eqnarray*}
\mathrm{R}\Gamma(\Sh_{K^p}^{\tor}, \qq_p)^{\la}{\otimes}_{\qq_p} \C_p & = &S^{\arit}(\mathrm{R}\Gamma(\Sh_{K^p}^{\tor}, \qq_p)^{\la}{\otimes}_{\qq_p} \C_p) \otimes_{\bar{E}} \C_p, \\ 
\mathrm{R}\Gamma_c(\Sh_{K^p}, \qq_p)^{\la}{\otimes}_{\qq_p} \C_p & = &S^{\arit}(\mathrm{R}\Gamma_c(\Sh_{K^p}, \qq_p)^{\la}{\otimes}_{\qq_p} \C_p) \otimes_{\bar{E}} \C_p.
\end{eqnarray*}
\item The action of $\mu$ via $\Theta_{\hor}$  on $\mathrm{R}\Gamma(\Sh_{K^p}^{\tor}, \qq_p)^{\la}{\otimes}_{\qq_p} \C_p =  \mathrm{R}\Gamma_{\an}(\Sh_{K^p}^{\tor}, \oscr^{\la}_{\Sh_{K^p}^{\tor}})$ and $\mathrm{R}\Gamma_c(\Sh_{K^p}^{\tor}, \qq_p)^{\la}{\otimes}_{\qq_p} \C_p =  \mathrm{R}\Gamma_{\an}(\Sh_{K^p}^{\tor}, \mathscr{I}^{\la}_{\Sh_{K^p}^{\tor}})$ is the  arithmetic Sen operator for the semilinear action of $\mathrm{Gal}(\bar{E}/E)$ (whose existence was guaranteed by (1)).
\end{enumerate}
\end{thm}
\begin{proof} We can take an affinoid  covering  $\mathcal{V} = \{ V_i\}$ of $\Sh^{\tor}_{K^p}$ with the property that $\oscr^{\la}_{\Sh^{\tor}_{K^p}}\vert_{V_i}$ is a colimit 
of acyclic %
sheaves $VB^0_{K_{p,n}, V_i}(\mathcal{C}^{\la})$.
It follows   that the \v{C}ech complex $\mathcal{C}(\mathcal{V}, \oscr^{\la}_{\Sh^{\tor}_{K^p}})$ represents $\mathrm{R}\Gamma(\Sh_{K^p}^{\tor}, \qq_p)^{\la}{\otimes}_{\qq_p} \C_p$ and carries a semi-linear  $\mathrm{Gal}(\bar{E}/E)$-action.  It follows from Theorem  \ref{thm-VB-primitive} (4), that we have  $$\mathcal{C}(\mathcal{U}, \oscr^{\la}_{\Sh^{\tor}_{K^p}}) = \mathcal{C}(\mathcal{U}, S^{\arit}(\oscr^{\la}_{\Sh^{\tor}_{K^p}})) \otimes_{\bar{E}} \C_p,$$ where $\mathcal{C}(\mathcal{U}, S^{\arit}(\oscr^{\la}_{\Sh^{\tor}_{K^p}}))$ is therefore the sub-complex of $\mathrm{Gal}(\bar{E}/E^{\cycl})$-smooth and $\mathrm{Gal}(\bar{E}/E)$-locally analytic vectors. We deduce that $\mathrm{R}\Gamma(\Sh_{K^p}^{\tor}, \qq_p)^{\la}{\otimes}_{\qq_p} \C_p$ admits a Sen operator and it is given by $\Theta_{\hor}(\mu)$ by Theorem \ref{thm-VB}, (9). 
\end{proof}

\subsection{Higher Coleman theory}

\subsubsection{Automorphic vector bundles}%
One can apply the  functor $\VBzero$ to the $G$-equivariant locally free sheaves of finite rank which are  parameterized by finite dimensional representations of $M$.
Let $\kappa \in X^\star(T)^{M,+}$. We let $\omega^{\kappa,\sm} = \VBzero(\mathcal{L}_\kappa)$ (where $\mathcal{L}_\kappa$ is constructed in Example \ref{example-construction-classical-sheaves}). The sheaf $\omega^{\kappa,\sm}$ descends to a sheaf $\omega^\kappa_{K_p}$ on the Shimura variety $\Sh_{K_pK^p}^{\tor}$ (the usual sheaf of modular forms of weight $\kappa$). 
By construction $\mathrm{R}\Gamma( \Sh_{K^p}^{\tor}, \omega^{\kappa,\sm})$ is a  complex of smooth admissible  $G(\qq_p)$-representations, equal to $\colim_{K_p} \mathrm{R}\Gamma(\Sh_{K_pK^p}^{\tor}, \omega^\kappa_{K_p})$. 
Recall that we have  denoted by $D_{K_p}$ the divisor of the boundary in $\Sh_{K_pK^p}^{\tor}$. We then consider the cuspidal subsheaf $\omega^\kappa_{K_p}(-D_{K_p})$. Passing to the limit, we get $$\omega^{\kappa,\sm}(-D) = \colim \pi_{K_p}^{-1}\omega^\kappa_{K_p}(-D_{K_p}) = \omega^{\kappa,\sm} \otimes_{\oscr_{\Sh_{K^p}^{\tor}}^{\sm}} \mathscr{I}^{\sm}_{\Sh_{K^p}^{\tor}}.$$ 
Similarly, $\mathrm{R}\Gamma( \Sh_{K^p}^{\tor}, \omega^{\kappa,\sm}(-D))$ is a  complex of smooth admissible  $G(\qq_p)$-representations, equal to $\colim_{K_p} \mathrm{R}\Gamma(\Sh_{K_pK^p}^{\tor}, \omega^\kappa_{K_p}(-D_{K_p}))$.

\begin{rem}For~$G=\GSp_4 $, the tautological exact sequence over $\mathcal{FL}$ is 
$$ 0 \rightarrow \mathcal{L}_{(0,  -1; 1)} \rightarrow St \otimes \oscr_{\mathcal{FL}} \rightarrow \mathcal{L}_{(1,  0; 1)} \rightarrow 0$$ which pulls back to $$0 \rightarrow Lie(A)_{K_p}(1) \otimes_{ \oscr_{\Sh_{K_pK^p}^{\tor}}} \oscr_{\Sh_{K^p}^{\tor}} \rightarrow T_pA \otimes_{\ZZ_p} \oscr_{\Sh_{K^p}^{\tor}} \rightarrow (\omega_{A^t})_{K_p} \otimes_{ \oscr_{\Sh_{K_pK^p}^{\tor}}} \oscr_{\Sh_{K^p}^{\tor}} \rightarrow 0$$ for any level  $K_p$. 
We deduce that $\omega^{(0, -1; 1),\sm} = Lie(A)_{K_p}(1) \otimes
\oscr^{\sm}_{\Sh_{K^p}^{\tor}} $ and thus that $$\omega^{(1, 0;
  -1),\sm} = (\omega_{A})_{K_p}(-1) \otimes_{
  \oscr_{\Sh_{K_pK^p}^{\tor}}}   \oscr^{\sm}_{\Sh_{K^p}^{\tor}} $$ by
duality; so our normalization of the weights of Siegel
modular forms is the standard one.%
\end{rem}

\subsubsection{Higher Coleman sheaves} We now fix $w \in \WM$. Let $\lambda \in X^\star(T)_{E}$. We consider the exact  functor:

$$ \VBzero \circ HCS_{w,\lambda}: \ocal(\mf{m}_w,\mf{b}_{M_w})_{\lambda-\alg} \rightarrow \LB_{(\mathfrak{g}, B)}(C_w^\dag) \rightarrow \mathrm{Mod}_{B(\qq_p)}( \oscr^{\sm}_{\pi_{HT}^{-1} C_w^\dag}). $$

\begin{lem} The functor $\VBzero \circ HCS_{w,\lambda}$ factors through the category $\mathrm{Mod}^{\lambda-\sm}_{B(\qq_p)}( \oscr^{\sm}_{\pi_{HT}^{-1} C_w^\dag})$. 
\end{lem}
\begin{proof} This is a combination of Proposition \ref{prop-smooth-lambdasmooth} and of Remark \ref{rem-lambda-smooth-Q}.
\end{proof}

\begin{lem}\label{lem-VBandclassicalsheaf} Let $\lambda \in
  X^\star(T)^{M,+}$. Let $L(\mf{m}_w)_{-w^{-1}w_{0,M}\lambda} \in \ocal(\mathfrak{m}_w,
  \mf{b}_{M_w})_{0-\alg}$ be the finite dimensional irreducible representation of highest weight $-w^{-1}w_{0,M}\lambda$. Then $ \VBzero \circ HCS_{w,0}(L(\mf{m}_w)_{-w^{-1}w_{0,M}\lambda} ) = \omega^{\lambda,sm}\vert_{ \pi_{HT}^{-1} C_w^\dag}$. 
\end{lem}
\begin{proof}.  Let $V_\lambda$ be the highest weight
  $\lambda$-representation of $M$. Then $V_\lambda$ identifies
  with~$\Lmw_{w^{-1}\lambda}$ if we conjugate $\mathfrak{m}_w$ to $\mathfrak{m}$. It follows from the definitions that $\mathcal{L}_\lambda\vert_{C_w^\dag} = HCS_{w,0}(L(\mf{m}_w)_{-w^{-1}w_{0,M}\lambda} )$ and the conclusion follows from applying $VB^0$. 
  \end{proof}

As in
Definition~\ref{defn:Verma-module-notation}, we have the Verma module
$\Mmw_{\lambda}$
of weight $\lambda$. %

\begin{defn}\label{defn-HigherColeman} We define the following object of $\mathrm{Mod}^{-w^{-1}w_{0,M}\lambda-\sm}_{B(\qq_p)}( \oscr^{\sm}_{\pi_{HT}^{-1} C_w^\dag})$: 
$$\omega^{\dag, \lambda}_{w}:= 
\VBzero \circ HCS_{w,-w^{-1}w_{0,M}\lambda}(\Mmw_{-w^{-1}w_{0,M} \lambda}).$$
\end{defn}

\begin{rem} This definition compares with \cite[\S 6.3]{boxer2021higher} as follows. 
 In that reference we (GB+VP) defined Banach sheaves $\mathcal{V}_{\nu}^{n-an}$ for  characters $\nu : T(\Z_p) \rightarrow \C_p^\times$, $n$ large enough, over certain quasi-compact  open subspaces $\pi^{-1}_{HT}(]C_{w,k}[_{n,n}K_p)$ of $\Sh^{\tor}_{K^pK_p}$ (for $K_p$ small enough), where $]C_{w,k}[_{n,n}K_p$ is a quasi-compact open subset of $C_w$. 
 The colimit over $K_p$,  over all $\nu$ with $\mathrm{d\nu} = \lambda$ and over
 all $n$ of the $\mathcal{V}_{\nu}^{n-an}$ identifies canonically with the germ
 of $\omega^{\dag, \lambda}_w$ at $\pi_{HT}^{-1}(\{w\})$. A slight change of
 perspective from \cite{boxer2021higher} is therefore this passage to the limit,
 and the fact that we define the sheaf $\omega^{\dag, \lambda}_w$ over the
 entire $\pi_{HT}^{-1}(C_w)$. For the definition and computation of the finite
 slope part of higher Coleman theory, the sheaves $\mathcal{V}_{\nu}^{n-an}$ are
 however sufficient (and seemed to us easier to define in the first place). See
 Section \ref{subsubsec:comparison-to-BP21} for further details.  
 \end{rem}

One reason for the twist from $\lambda$ to
 $-w^{-1}w_{0,M}\lambda$ is in order to obtain
 Remark~\ref{rem:horizontal-on-lambda-normalization} and
 Proposition~\ref{prop:integral-weight-simple-sheaf-computed}; more
 conceptually, the multiplication by $w^{-1}$  is justified by the change of base
 point in the flag variety, and the appearance of~ $-w_{0,M}$ is due to the usual duality
 involution on highest weights which comes from the contravariance of
 $HCS_{w,\lambda}$.

\begin{rem}\label{rem:horizontal-on-lambda-normalization} We see from Proposition \ref{prop-action-centre-enveloping} that the $\Theta_{\hor}$-action on $$HCS_{w,-w^{-1}w_{0,M}\lambda}(\Mmw_{-w^{-1}w_{0,M} \lambda})$$  is via $\lambda$. %
(We recall from
  Remark~\ref{rem:iota-on-HC} that $HC_{\mathfrak{m}}\iota = -w_{0,M} HC_{\mathfrak{m}}$.) 
  Thus, the arithmetic Sen operator acts via $\langle \mu , \lambda \rangle$ on $\omega^{\dag,\lambda}_w$. 
\end{rem}%

 \begin{prop}\label{prop:integral-weight-simple-sheaf-computed} Assume
   that $\lambda \in X^\star(T)^{M, +}$. We have an injective, $B(\qq_p)$-equivariant map:
   $\omega^{ \lambda, \sm}\vert_{\pi_{HT}^{-1}(C_w^\dag)} \otimes E(-w^{-1} w_{0,M} \lambda) \rightarrow \omega^{\dag, \lambda}_w.$
   \end{prop}

\begin{proof}%
  By Lemma \ref{lem-VBandclassicalsheaf},  $\omega^{\lambda, \sm}\vert_{\pi_{HT}^{-1}(C_w^\dag)}= \VBzero \circ HCS_{w,0}( \Lmw_{-w^{-1}w_{0,M} \lambda})$.  The map of the proposition  is thus given by applying $\VBzero \circ HCS_{w,-w^{-1}w_{0,M}\lambda}$ to the surjection $\Mmw_{-w^{-1}w_{0,M} \lambda} \rightarrow \Lmw_{-w^{-1}w_{0,M} \lambda}$, using  Proposition
\ref{TwistingHCS}.%
\end{proof}

\subsubsection{The Bruhat stratification and a Cousin spectral sequence}\label{sect-review-Cousin}

We recall the stratification into $B$-orbits  $\mathcal{FL} = \coprod_{w \in \WM} C_w$, with $C_{w} = P \backslash P {w} B$. 
We let $X_{w} = \cup_{ w' \leq w} C_{w'}$ be the Schubert variety. %
We let $j_{w}: C_{w} \hookrightarrow \mathcal{FL}$ be the locally closed immersion. It induces $j_{w, \Sh^{\tor}_{K^p}}: \pi_{HT}^{-1}(C_{w}) \rightarrow \Sh^{\tor}_{K^p}$.   
\begin{defn} If $\mathscr{F}$ is any sheaf  of solid  $E$-modules on $\pi_{HT}^{-1}(C_w)$,  we let $$ \mathrm{R}\Gamma_w( \Sh^{\tor}_{K^p}, \mathscr{F}) =  \mathrm{R}\Gamma( \Sh^{\tor}_{K^p}, (j_{w, \Sh^{\tor}_{K^p}})_! \mathscr{F}).$$
\end{defn}
In this definition, $(j_{w, \Sh^{\tor}_{K^p}})_! $ is the extension by zero functor on abelian sheaves of solid abelian groups.  By abuse of notation, if $\mathscr{F}$ is defined over any  subset of $\Sh^{\tor}_{K^p}$ containing $\pi_{HT}^{-1}(C_w)$, we let 
$ \mathrm{R}\Gamma_w( \Sh^{\tor}_{K^p}, \mathscr{F}) =  \mathrm{R}\Gamma( \Sh^{\tor}_{K^p}, (j_{w, \Sh^{\tor}_{K^p}})_! \mathscr{F}\vert_{\pi_{HT}^{-1}(C_w)})$. 

We now explain that if $\mathscr{F}$ carries a $B(\qq_p)$-equivariant structure, the cohomology is a $B(\qq_p)$-representation. 

\begin{lem}\label{lem:upgrade-to-B-action}The functor
\begin{eqnarray*}
 \Mod(\pi_{HT}^{-1}(C_w)) &\rightarrow& D(\Mod(E)) \\  
 \mathscr{F} &\mapsto& \mathrm{R}\Gamma_w(\Sh^{\tor}_{K^p}, \mathscr{F})
 \end{eqnarray*} can be upgraded to functors: 
\begin{eqnarray*}
\Mod_{B(\qq_p)}( \pi_{HT}^{-1}(C_w)) &\rightarrow& D( \Mod_{B(\qq_p)}(E)). \\
\Mod^{\lambda-\sm}_{B(\qq_p)}( \pi_{HT}^{-1}(C_w)) &\rightarrow& D( \Mod^{\lambda-\sm}_{B(\qq_p)}(E)).
\end{eqnarray*}
\end{lem}
\begin{proof} The existence of $\Mod_{B(\qq_p)}( \pi_{HT}^{-1}(C_w))
  \rightarrow D( \Mod_{B(\qq_p)}(E))$ follows from the fact that
  $\Mod_{B(\qq_p)}( \pi_{HT}^{-1}(C_w))$ has enough injectives (see
  Lemma~\ref{coro-enough-inj}). %

  We verify that this induces a functor $\Mod^{\lambda-\sm}_{B(\qq_p)}( \pi_{HT}^{-1}(C_w))
  \rightarrow D( \Mod^{\lambda-\sm}_{B(\qq_p)}(E))$. We can reduce to
  the case that $\lambda=0$. 
Then $ (j_{w, \Sh^{\tor}_{K^p}})_!
\mathscr{F}$ is a smooth $B(\qq_p)$-equivariant sheaf on
$\Sh^{\tor}_{K^p}$, which is quasi-compact, so by definition the
global sections  $\HH^0_w(\Sh^{\tor}_{K^p}, \mathscr{F}) \in
\Mod^{\sm}_{B(\qq_p)}(E)$ are a smooth $B(\qq_p)$-representation. The
result follows by applying this to an injective resolution.
\end{proof}

\begin{prop}\label{prop:abstract-Cousin-spectral-sequence} Let $\mathscr{F}$ be a solid abelian sheaf defined over $\Sh^{\tor}_{K^p}$.
We have a spectral sequence (Cousin spectral sequence):
$$ E_1^{p,q} =  \oplus_{w \in \WM, \ell(w)=p} \HH^{p+q}_{w}( \Sh^{\tor}_{K^p}, \mathscr{F}) \Rightarrow \HH^{p+q}( \Sh^{\tor}_{K^p}, \mathscr{F}).$$
\end{prop}
\begin{proof}
  See e.g.\ \cite[\S 2.3]{boxer2021higher}.
\end{proof}%
\begin{rem}If~$\cF$ in Proposition~\ref{prop:abstract-Cousin-spectral-sequence}
  is $B(\Qp)$-equivariant, then so (by construction) is the Cousin spectral sequence.
  \label{rem:Cousin-B-equivariant}
\end{rem}

\subsubsection{Tools for computing the cohomology}\label{subsubsec:tools-for-computing-cohomology}
In this section, we give a few basic tools for computing cohomology.  All adic
spaces are locally of finite type over $\Spa(E, \ocal_E)$ and are separated
unless specifically mentioned (thus they correspond to  separated rigid analytic
spaces in the sense of Tate). We will often consider the case of Stein or
quasi-Stein spaces, which we recalled in Definition~\ref{defn:quasi-Stein-space}.

Let $\mathcal{X}$ be an adic space over $\Spa(E, \ocal_E)$. Let $\mathscr{F}$ be a coherent sheaf over $\mathcal{X}$. We can think of $\mathscr{F}$ as an abelian sheaf of solid $E$-modules and compute $\mathrm{R}\Gamma( \mathcal{X}, \mathscr{F})$ in the derived category of solid $E$-modules, $D(\Mod(E))$.
\begin{thm}[Tate acyclicity]\label{thm-Tateacycli}  Let $\mathcal{X}$ be an affinoid adic space. Let $\mathscr{F}$ be a coherent sheaf over $\mathcal{X}$.
The cohomology $\mathrm{R}\Gamma( \mathcal{X}, \mathscr{F}) \in D(\Mod(E))$ is concentrated in degree $0$. 
\end{thm}
\begin{proof} Let $\{\mathcal{U}_i\}_{i \in I}$ be a finite affinoid covering of $\mathcal{X}$. We know by \cite[Thm. 2.5]{MR1306024} that the augmented Čech complex of $\{\mathcal{U}_i\}_{i \in I}$ is a long exact sequence of classical Banach spaces. Therefore, it is also an exact sequence of solid Banach spaces (see Remark \ref{rem-exact-Banach-classical-solid}). One concludes by \cite[3.8, coro. 4]{MR0102537}.
\end{proof}

\begin{cor}\label{cor-tateacycli} Let $\mathcal{X}$ be a finite type, separated adic space over $\Spa(E, \ocal_E)$. Let $\mathscr{F}$ be a coherent sheaf over $\mathcal{X}$. The cohomology of $\mathscr{F}$ agrees with Čech cohomology.
\end{cor}
\begin{proof} This is an application of \cite[3.8, coro. 4]{MR0102537} and Theorem \ref{thm-Tateacycli}.
\end{proof}

\begin{cor} Let $\mathcal{X}$ be a quasi-Stein affinoid adic space. Let $\mathscr{F}$ be a coherent sheaf over $\mathcal{X}$.
 The cohomology $\mathrm{R}\Gamma( \mathcal{X}, \mathscr{F}) \in D(\Mod(E))$ is concentrated in degree $0$. 
\end{cor}
\begin{proof} By Corollary \ref{cor-tateacycli}, the cohomology is given by
  $\mathrm{R}\lim_n \mathrm{H}^0( \mathcal{X}_n, \mathscr{F})$. One knows by
  \cite{MR0210949} that  the maps $ \mathrm{H}^0( \mathcal{X}_{n+1},
  \mathscr{F}) \rightarrow  \mathrm{H}^0( \mathcal{X}_n, \mathscr{F})$ of
  classical Banach spaces have dense image. Then  the topological Mittag-Leffler
  \cite[Lem.\ 3.27]{MR4475468} shows that $\mathrm{R}\lim_n \mathrm{H}^0(
  \mathcal{X}_n, \mathscr{F}) = \lim_n  \mathrm{H}^0( \mathcal{X}_n,
  \mathscr{F})$, as required.
\end{proof}

By \cite[Lem.\ 3.21]{MR4475468} a Banach or Smith space over $E$ is flat. 
\begin{cor}\label{coro-coho-quasi-steinV} Let $\mathcal{X}$ be a finite type, separated  adic space  over $\Spa(E,\ocal_E)$, and let $\mathscr{F}$ be a coherent sheaf.  Let $V$ be a Banach or Smith space over $E$. Then $\mathrm{R}\Gamma( \mathcal{X}, \mathscr{F} \otimes_EV) = \mathrm{R}\Gamma( \mathcal{X}, \mathscr{F})\otimes_EV$.
\end{cor}
\begin{proof} In the affinoid case, we check that $\mathrm{R}\Gamma(
  \mathcal{X}, \mathscr{F} \otimes_EV) = \HH^0( \mathcal{X},
  \mathscr{F})\otimes_E V[0]$. To see this, we simply check it on the
  Čech cohomology of arbitrary finite affinoid covers, and this follows from the
  flatness of $V$. For a general $\mathcal{X}$, we deduce (again using~\cite[3.8, coro. 4]{MR0102537}) that the cohomology of
  $\mathscr{F} \otimes V$ is computed by the Čech cohomology of a finite
  affinoid cover, and we use one more time that $V$ is flat. %
\end{proof}

\begin{cor}\label{cor-cohosteinV}  Let $\mathscr{F}$ be a coherent
  sheaf over $\mathcal{X}$, and let $V$ be a Banach or Smith space
  over $L$. Suppose that $\mathcal{X}$ is covered by finitely many
  quasi-Stein spaces. Then $\mathrm{R}\Gamma( \mathcal{X}, \mathscr{F}
  \otimes_EV) = \mathrm{R}\Gamma( \mathcal{X}, \mathscr{F})\otimes_E
  V$, and if  $\mathcal{X}$ is itself a quasi-Stein  space, then we have $\mathrm{R}\Gamma( \mathcal{X}, \mathscr{F} \otimes_E V) = \HH^0(\mathcal{X}, \mathscr{F}) \otimes_E V[0]$. 
\end{cor}
\begin{proof}If~$\cX$ is quasi-Stein, then by Corollary
  \ref{coro-coho-quasi-steinV}, the cohomology is given by
  $\mathrm{R}\lim_n (\mathrm{H}^0( \mathcal{X}_n, \mathscr{F})
  \otimes_E V)$. By the topological Mittag-Leffler and   \cite[Lem.\
  3.28]{MR4475468}, this limit is simply $(\lim_n \mathrm{H}^0(
  \mathcal{X}_n, \mathscr{F})) \otimes_E V$. From this, we deduce that
  if  $\mathcal{X}$ is covered by finitely many
  quasi-Stein spaces, then the cohomology of $\mathscr{F} \otimes_E V$
  is computed by the Čech cohomology of any finite quasi-Stein cover,
  and the claim  follows again from flatness of $V$. 
\end{proof}

We also consider duality and cohomology with compact support. Let $\mathcal{X}$ be an adic space. Let $\mathscr{F}$ be a solid sheaf of abelian groups. Following \cite[5.2]{MR1734903}, we
let $\HH^0_c( \mathcal{X}, \mathscr{F}) = \colim_{\mathcal{Z}} \HH^0_{\mathcal{Z}}(\mathscr{F})$ be the space of sections  with compact support where $\mathcal{Z}$ runs through the poset of  closed subsets of $\mathcal{X}$ which are proper over $\Spa (E, \ocal_E)$. 
We let  $\mathrm{R}\Gamma_{c}(\mathcal{X}, \mathscr{F}) = \colim_{\mathcal{Z}}
\mathrm{R}\Gamma_{\mathcal{Z}}(\mathcal{X}, \mathscr{F})$ be the cohomology
with compact
support. (Here $\mathrm{R}\Gamma_{\mathcal{Z}}(\mathcal{X},
  \mathscr{F})$ is by definition equal to
  $\mathrm{R}\Gamma(\mathcal{X},(i_{\cZ})_{*}i_{\cZ}^{!}\mathscr{F})$, where $i_{\cZ}:\cZ\into\cZ$.)

\begin{rem} For example,  $\mathcal{Z}$ is proper if it arises as the inverse image of a proper subset of a formal model of a quasi-compact open subset of $\mathcal{X}$.
\end{rem}

Consider a quasi-compact and separated adic space of finite type $\mathcal{X}$. Let $j: \mathcal{U} \hookrightarrow \mathcal{X}$ be an open adic subspace of $\mathcal{X}$. We assume that $\mathcal{U}$ admits an increasing covering by quasi-compact opens $\mathcal{U} = \cup_n \mathcal{U}_n$ with the property that $\mathcal{U}_n$ is relatively compact in $\mathcal{U}_{n+1}$. %

\begin{example} One can take $\mathcal{X} = \mathbb{P}^1$, $\mathcal{U} =
  \mathbb{A}^{1,an}$. More generally, one can take $\mathcal{X}$ a proper finite
  type analytic space and $\mathcal{U}$ be the complement of a Zariski closed
  subset  of $\mathcal{X}$ (see \cite[5.9]{MR1032938}). One can also take
  $\mathcal{X} = \Spa (L\langle T \rangle, \ocal_L \langle T \rangle )$ (the
  closed unit ball of radius~$1$), with $\mathcal{U}$  the open unit ball of radius $1$ inside $\mathcal{X}$. %
\end{example}

\begin{lem}\label{lemma-basic-compact-support} Let $\mathscr{F}$ be a
  solid abelian sheaf over $\mathcal{U}$ as above. Then $\mathrm{R}\Gamma(\mathcal{X}, j_! \mathscr{F}) = \mathrm{R}\Gamma_c(\mathcal{U}, \mathscr{F}) = \colim \mathrm{R}\Gamma_{\cZ}(\mathcal{U}, \mathscr{F})$ where $\cZ$ runs through all closed subsets of $\mathcal{X}$ contained in $\mathcal{U}$.  
Moreover, there exists an increasing family of closed subspaces $\cU_n\subseteq
\cZ_n\subseteq\cU_{n+1}$, each with quasi-compact complement in $\mathcal{X}$, such that $\mathrm{R}\Gamma_c(\mathcal{U}, \mathscr{F}) = \colim_n \mathrm{R}\Gamma_{\cZ_n}(\mathcal{U}, \mathscr{F})$.
\end{lem}  
\begin{proof} Firstly, we check that a subset $\cZ\subset\cU$ is closed
  in $\cX$ if and only if $\cZ$ is proper over~$\Spa(E,\cO_E)$.
  Suppose that $\cZ$ is a  closed subset of $\mathcal{X}$ contained in $\mathcal{U}$.  %
Then  $\cZ$ is quasi-compact (as it is a closed subset of
$\mathcal{X}$). Let $\mathcal{U} = \cup_n \mathcal{U}_n$. By
assumption, $\overline{\mathcal{U}_n}$ (the closure of $\mathcal{U}_n$
in $\mathcal{U}$) is proper over $\Spa(E, \ocal_E)$. Since $\cZ$ is
quasi-compact, $\cZ \subseteq \mathcal{U}_n$ for $n$ large
enough. Hence, $\cZ$ is closed in $\overline{\mathcal{U}_n}$, thus
partially proper. Conversely, if $\cZ$ is a  subset of $\mathcal{U}$,
proper over $\Spa(E, \ocal_E)$, then  ${\cZ} \rightarrow \mathcal{X}$ is
proper, as claimed.

For any closed subset $\cZ\subseteq\cU$, the counit of adjunction for $i_{\cZ}:\cZ\into\cX$ gives a map $(i_{\cZ})_\star i_{\cZ}^! j_{!}\mathscr{F}
\rightarrow j_! \mathscr{F}$ %
which induces a map $\colim_{\cZ} (i_{\cZ})_\star
i_{\cZ}^!j_{!} \mathscr{F} \rightarrow j_! \mathscr{F}$. This map is injective
(both sheaves are subsheaves of $j_\star \mathscr{F}$), and both sheaves
restrict to $\mathscr{F}$ on $\mathcal{U}$ and have zero stalk at points of
$\mathcal{X} \setminus \mathcal{U}$; thus the map is an isomorphism. Since
$\mathcal{X}$ is quasi-compact and separated, this implies that
$\mathrm{R}\Gamma(\mathcal{X}, j_! \mathscr{F}) = \colim_{{Z}} \mathrm{R}\Gamma_{\cZ}(\mathcal{X},j_{!} \mathscr{F})$. 
Moreover, $\colim_{{\cZ}} \mathrm{R}\Gamma_{\cZ}(\mathcal{X}, j_{!}\mathscr{F}) =
\mathrm{R}\Gamma_c(\mathcal{U}, \mathscr{F})$ by definition. %
We finally claim that for each~$n$ there exists a closed subspace $\mathcal{U}_n \subseteq \cZ_n \subseteq \mathcal{U}_{n+1}$ with quasi-compact complement in $\mathcal{X}$. %
We let $\cup_{i \in I} V_i$ be a covering by quasi-compact opens of $\overline{\mathcal{U}_n}^c$. We claim that there exists a finite subset $I'$ of $I$ such that $\cup_{i \in I'} V_i$ contains $\mathcal{U}_{n+1}^c$. This follows from endowing $\mathcal{X}$ with the constructible topology and noticing that $\mathcal{U}_{n+1}^c$ is compact in this topology. We can take $\cZ_n = (\cup_{i \in I'} V_i)^c$.  The $\cZ_n$ are clearly cofinal among all $\cZ$'s. 
\end{proof}

The basic duality statement is the following. 

\begin{thm}[\cite{MR1094850}]\label{thm-Chiarellotto}%
  Let $\mathscr{F}$ be a locally free coherent sheaf defined over a smooth Stein space $\mathcal{X}$ of dimension $d$.  The cohomology $\mathrm{R}\Gamma_{c}(\mathcal{X}, \mathscr{F})$ is concentrated in degree $d$ and is an $LB$-space of compact type.  Moreover, 
$$\HH^d_{c}(\mathcal{X}, \mathscr{F}) = \mathrm{Hom}_E( \HH^0(\mathcal{X}, D(\mathscr{F})), E)$$ where 
 $D(\mathscr{F}) =\underline{ \mathrm{Hom}}( \mathscr{F}, \Omega^{d}_{\mathcal{X}})$. 
\end{thm}

\begin{rem}\label{rem-using-Stein} Although we use the formalism of compactly supported cohomology for abelian sheaves, we get the correct ``coherent duality''. This is a favorable property of Stein spaces. 
\end{rem}

We now extend this result to orthonormalizable Banach sheaves. 

\begin{cor}\label{coro-dualitystein} Let $\mathscr{F}$ be a locally free
  coherent sheaf defined  over a Stein space  $\mathcal{X}$ of dimension~$d$. Let $V$ be a Smith space over $E$.  We have that  $\mathrm{R}\Gamma_{c}(\mathcal{X}, \mathscr{F} \otimes_E V) =  \mathrm{R}\Gamma_{c}( \mathcal{X}, \mathscr{F})\otimes_E V$ is concentrated in degree $d$ and moreover, 
$$\HH^d_{c}(\mathcal{X}, \mathscr{F} \otimes V) = \mathrm{Hom}_E( \HH^0(\mathcal{X}, D(\mathscr{F}) \otimes V^\vee), E).$$
\end{cor}
\begin{proof} We let $\mathcal{X} = \cup_n \mathcal{X}_n$ be an open cover by
  affinoids, with ${\mathcal{X}_n}$ relatively compact in
  $\mathcal{X}_{n+1}$. We also let $\mathcal{X}_n \subseteq \mathcal{Z}_n
  \subseteq \mathcal{X}_{n+1} $ be closed subspaces with quasi-compact
  complement as in the statement of Lemma \ref{lemma-basic-compact-support}
  (with~$\cU_n=\cX_n$). It follows from  Lemma \ref{lemma-basic-compact-support} that $\mathrm{R}\Gamma_{c}(\mathcal{X}, \mathscr{F} \otimes_E V)
= \colim_n \mathrm{R}\Gamma_{\mathcal{Z}_n}(\mathcal{X}_{n+1},
\mathscr{F} \otimes_E V)$. %
We observe that  $\mathrm{R}\Gamma_{\mathcal{Z}_n}(\mathcal{X}_{n+1}, \mathscr{F} \otimes V) =  \mathrm{R}\Gamma_{\mathcal{Z}_n}(\mathcal{X}_{n+1}, \mathscr{F} )\otimes V$ by an easy reduction to Corollary \ref{coro-coho-quasi-steinV}. 
We deduce that $\mathrm{R}\Gamma_{c}(\mathcal{X}, \mathscr{F} \otimes_E V) = \mathrm{R}\Gamma_{c}(\mathcal{X}, \mathscr{F} )\otimes_E V$. 
Next, we see that \begin{align*}
\mathrm{Hom}_E( \HH^0( \mathcal{X}, D(\mathscr{F}) \otimes V^\vee) , E) &= \mathrm{Hom}_E( \HH^0( \mathcal{X}, D(\mathscr{F})) \otimes V^\vee , E) ~ \textrm{by Corollary \ref{cor-cohosteinV}}\\
&= \mathrm{Hom}_E( \HH^0( \mathcal{X},D(\mathscr{F})), V)~\textrm{by adjunction and  Proposition \ref{prop: antieq Smith Banach}} \\
&= \mathrm{Hom}_E( \HH^0( \mathcal{X},D(\mathscr{F})), E) \otimes_E V~\textrm{by \cite[Thm.\ 3.40]{MR4475468}} \\
&= \HH^d_{c}(\mathcal{X}, \mathscr{F}) \otimes V~\textrm{by Theorem \ref{thm-Chiarellotto}}\\
&= \HH^d_{c}(\mathcal{X}, \mathscr{F}\otimes V)~\textrm{by the first point.}\qedhere
\end{align*}
\end{proof}%

\subsubsection{On the computation of the local cohomology}

We want to give some formulas for computing $ \mathrm{R}\Gamma_w( \Sh^{\tor}_{K^p}, \mathscr{F})$.

The Bruhat cell $C_w$ is an affine space of dimension $\ell(w)$ and
its closure $X_w$ is a compactification of this affine space. We
 recall (see \cite[Lemma 3.1.3]{boxer2021higher} for example) that $C_w = w \prod_{\alpha \in \Phi^+ \cap w^{-1} \Phi^{-,M}} U_{\alpha}$ where $U_{\alpha}$ is the $\alpha$-root space, isomorphic to $\mathbb{A}^{1,an}$. 
We also have a neighborhood of $C_w$, $U_w = w \prod_{ \alpha \in w^{-1} \Phi^{-,M}} U_{\alpha}$. Let us pick a coordinate $u_\alpha$ on each $U_\alpha$. 

\begin{defn}\label{example-chooseopen} 
We  now define certain subsets of  $Z_n, Z_{n,m}$, $U_{n}$ and $U_{n,m}$ of $C_w$, for $n,m \in \ZZ_{\geq 0}$. 
Let  $n_0 \geq 0$ be fixed. We take: 
\begin{itemize}
\item  $Z_n = \{ x \in C_w, \forall \alpha \in \Phi^+ \cap w^{-1} \Phi^{-,M}, \vert u_{\alpha} \vert_x < \vert p \vert^{n_0-n}_{x} \}$.
\item  $U_n =  \{ x \in C_w, \forall \alpha \in \Phi^+ \cap w^{-1} \Phi^{-,M}, \vert u_{\alpha} \vert_x \leq \vert p \vert^{n_0-n-1}_{x} \}$. 
\item $Z_{n,m} = \{ x \in U_w, \forall \alpha \in \Phi^+ \cap w^{-1} \Phi^{-,M}, \vert u_{\alpha} \vert_x < \vert p \vert^{n_0-n}_{x}, \forall \alpha \in \Phi^- \cap w^{-1} \Phi^{-,M}, \vert u_{\alpha} \vert_x < \vert p \vert^{m}_{x} \}$,
\item  $U_{n,m} =  \{ x \in U_w, \forall \alpha \in \Phi^+ \cap w^{-1} \Phi^{-,M}, \vert u_{\alpha} \vert_x \leq \vert p \vert^{n_0-n-1}_{x}, \forall \alpha \in \Phi^- \cap w^{-1} \Phi^{-,M}, \vert u_{\alpha} \vert_x \leq \vert p \vert^{m}_{x} \}$.
\end{itemize}
\end{defn}
 
 Here are some obvious properties of these sets. 
 
 \begin{itemize}
\item $\cup_n Z_n =  \cup_n U_n =C_w$,
\item The complement of $Z_n$ in $X_w$ is a quasi-compact open subset,
\item $U_n$ is a quasi-compact open subset of $X_w$, 
\item $U_n = \cap_m  U_{n,m}$, and $U_{n,m} \cap X_w = U_n$,
\item $U_{n,m}$ is a quasi-compact open of $\mathcal{FL}$,
\item $Z_{n,m}$  is a closed subset of $U_{n,m}$ with quasi-compact complement,
\item  $Z_n = U_n \cap Z_{n,m}$.
\end{itemize}

We let $F_n = \pi_{HT}^{-1}(Z_n)$, $V_n = \pi_{HT}^{-1}(U_n)$, $F_{n,m} = \pi_{HT}^{-1}(Z_{n,m})$ and $V_{n,m} = \pi_{HT}^{-1}(U_{n,m})$.

\begin{lem}\label{lem-computeasinductivelimit} Let $\mathscr{F}$ be a solid abelian sheaf over $\pi_{HT}^{-1}(C_w)$.  We have  $\mathrm{R}\Gamma_w( \Sh^{\tor}_{K^p}, \mathscr{F}) = \colim_n \mathrm{R}\Gamma_{ F_n}(V_n, \mathscr{F})$.
\end{lem}
\begin{proof}  This is a consequence of Lemma
  \ref{lemma-basic-compact-support} (since the~$V_n$ give a cofinal system of
  closed subsets of~$Z_w$). %
\end{proof}

\begin{lem}\label{lem-second-comp-inductivelimi} Let $m_0 \geq 0$ and let $\mathscr{F}'$ be a sheaf on  $V_{n,m_0}$ and let  $\mathscr{F} = i^{-1} \mathscr{F}'$ on $V_n$ for  $i: V_n \rightarrow V_{n,m_0}$. Then we have: 
$\mathrm{R}\Gamma_{F_n}(V_n, \mathscr{F}) = \colim_{m \geq m_0} \mathrm{R}\Gamma_{F_{n,m}}(V_{n,m}, \mathscr{F}')$. 
\end{lem}
\begin{proof} We consider the triangle:
$$\mathrm{R}\Gamma_{F_{n,m}}(V_{n,m}, \mathscr{F}') \rightarrow  \mathrm{R}\Gamma(V_{n,m}, \mathscr{F}') \rightarrow \mathrm{R}\Gamma(V_{n,m} \setminus F_{n,m}), \mathscr{F}') \stackrel{+1}\rightarrow $$
Passing to the colimit over $m$, we need to see that $\colim_{m}  \mathrm{R}\Gamma(V_{n,m}, \mathscr{F}') =  \mathrm{R}\Gamma(V_n, \mathscr{F})$ and 
$\colim_m \mathrm{R}\Gamma(V_{n,m} \setminus F_{n,m}, \mathscr{F}') =
\mathrm{R}\Gamma(V_n\setminus F_{n}, \mathscr{F})$. This follows from the 
fact that all spaces involved are quasi-compact.
\end{proof}

The combination of lemmas \ref{lem-computeasinductivelimit} and \ref{lem-second-comp-inductivelimi} allows us to compute the cohomology as a colimit of cohomologies of the shape $\mathrm{R}\Gamma_{F_{n,m}}(V_{n,m}, \mathscr{F}')$ where 
$\mathrm{R}\Gamma_{F_{n,m}}(V_{n,m}, \mathscr{F}')$ fits in a triangle: $$\mathrm{R}\Gamma_{F_{n,m}}(V_{n,m}, \mathscr{F}') \rightarrow  \mathrm{R}\Gamma(V_{n,m}, \mathscr{F}') \rightarrow \mathrm{R}\Gamma(V_{n,m} \setminus F_{n,m}, \mathscr{F}') \stackrel{+1}\rightarrow $$
and both $V_{n,m}$ and $V_{n,m} \setminus F_{n,m}$ are quasi-compact opens.

We now give another similar presentation of the cohomology which puts emphasis on the cohomology with compact support. We will be using Stein spaces because of remark \ref{rem-using-Stein}. 

\begin{defn}\label{example-open-coho-support} We define sets $X_n$, $X_{n,m}$ and $T_{n,m,r}$ for $n,m, r \in \ZZ_{\geq 0}$.  %

 \begin{itemize}
\item For any $\epsilon \in \Q_{>0}$, we put  $X_{n, \epsilon} =  \{ x
  \in C_w, \forall \alpha \in \Phi^+ \cap w^{-1} \Phi^{-,M}, \vert
  u_{\alpha} \vert_x \leq \vert p \vert^{n_0-n+\epsilon}_{x} \}$ and set $X_n = \cup_{\epsilon >0} X_{n, \epsilon}$. 
\item $X_{n,m, \epsilon} =  \{ x \in U_w, \forall \alpha \in \Phi^+
  \cap w^{-1} \Phi^{-,M}, \vert u_{\alpha} \vert_x \leq \vert p
  \vert^{n_0-n+\epsilon}_{x}, \forall \alpha \in \Phi^- \cap w^{-1}
  \Phi^{-,M}, \vert u_{\alpha} \vert_x \leq \vert p \vert^{m+
    \epsilon}_{x} \}$, and set $X_{n,m} = \cup_{\epsilon >0} X_{n,m, \epsilon}$. 
\item We also put $T_{n,m, r, \eta} =  \{ x \in X_{n,m }, \exists \alpha  \in \Phi^- \cap w^{-1} \Phi^{-,M}, \vert u_{\alpha} \vert_x \geq \vert p \vert^{r- \eta} \}$. We let  $T_{n,m,r}  = \cup_{\eta >0} T_{n,m,r,  \eta}$. 
\end{itemize}

\end{defn}

We see that $X_n$ is an increasing sequence of Stein subsets of $C_w$
with the property that $C_w = \cup_n X_n$ and that each $X_n$ is
included in a quasi-compact open subset~$Z_{n}$ of $C_w$. We see that
$X_{n,m}$ is a  decreasing (in~$m$) family of Stein open subsets  of
$\mathcal{FL}$ with $X_n = X_{n,m} \cap X_w$ and $\cap_m X_{n,m} =
X_n$. We have $\cup_{r>0} T_{n,m,r} = X_{n,m}\setminus X_n$. 

Let $Y_n = \pi_{HT}^{-1}(X_n)$.  We let $Y_{n,m} = \pi_{HT}^{-1}(X_{n,m})$ and $W_{n,m,r} = \pi_{HT}^{-1}(T_{m,n,r})$.

\begin{lem}\label{lem-computeasinductivelimitI} Let $\mathscr{F}$ be a solid abelian sheaf over $\pi_{HT}^{-1}(C_w)$.  We have  $\mathrm{R}\Gamma_w( \Sh^{\tor}_{K^p}, \mathscr{F}) = \colim_n \mathrm{R}\Gamma_c(Y_n, \mathscr{F})$.
\end{lem}
\begin{proof} This follows from Lemma
  \ref{lem-computeasinductivelimit}. Indeed the two inductive systems
  are equivalent as  $Y_n \subseteq F_n \subseteq Y_{n+1}$, and we
  have a  series of natural maps $\mathrm{R}\Gamma_{F_n}( V_n, \mathscr{F}) \rightarrow  \mathrm{R}\Gamma_c(Y_n,  \mathscr{F}) \rightarrow \mathrm{R}\Gamma_{F_{n+1}}( V_{n+1}, \mathscr{F})$. 
\end{proof}

\begin{lem}\label{lem-computeasinductivelimitII} Let $\mathscr{F}'$ be a sheaf on  $Y_{n,m}$ and let  $\mathscr{F} = i^{-1} \mathscr{F}'$ on $Y_n$.  Then we have a triangle:
$$\mathrm{R}\Gamma_c( Y_{n,m} \setminus Y_{n}, \mathscr{F}') \rightarrow \mathrm{R}\Gamma_c(Y_{n,m}, \mathscr{F}') \rightarrow \mathrm{R}\Gamma_c(Y_n, \mathscr{F}) \stackrel{+1}{\rightarrow} $$
Moreover, $\mathrm{R}\Gamma_c( Y_{n,m} \setminus Y_{n}, \mathscr{F}') = \colim_r \mathrm{R}\Gamma_c( W_{n,m,r}, \mathscr{F}')$. 
\end{lem}

\begin{proof}
Let $m \geq 0$. We consider the following commutative diagram:
\begin{eqnarray*}
\xymatrix{ Y_n \ar[r]^{j_n}  \ar[d]^{i} & \pi_{HT}^{-1}(X_w) \ar[d]^{i'} \\
Y_{n,m} \ar[r]^{j_{n,m}} &  \Sh^{\tor}_{K^p} \\
Y_{n,m} \setminus Y_n \ar[u]^j \ar[ru]^{i_{n,m}} &}
\end{eqnarray*}

 We have to prove that we have a triangle $$\mathrm{R}\Gamma(\Sh^{\tor}_{K^p},(i_{n,m})_! \mathscr{F}') \rightarrow \mathrm{R}\Gamma(\Sh^{\tor}_{K^p}, (j_{n,m})! \mathscr{F}')  \rightarrow \mathrm{R}\Gamma(\pi_{HT}^{-1}(X_w) , (j_n)_!\mathscr{F}) \stackrel{+1}\rightarrow  $$

We have the triangle   $j_! j^{-1} \mathscr{F}' \rightarrow \mathscr{F}' \rightarrow  i_\star  \mathscr{F} \stackrel{+1}\rightarrow$. We apply $(j_{n,m})_!$, and we get a triangle:
$(i_{n,m})_! j^{-1} \mathscr{F}' \rightarrow (j_{n,m})_! \mathscr{F}' \rightarrow   (j_{n,m})_! i_\star  \mathscr{F} \stackrel{+1}\rightarrow$.
Observe that  $(j_{n,m})_! i_\star  \mathscr{F} = i'_\star (j_n)_! \mathscr{F}$ (as $i_\star = i_!$ and $i'_\star = i'_!$).  We then apply $\mathrm{R}\Gamma(\Sh^{\tor}_{K^p},-)$.
\end{proof}

\subsubsection{Higher Coleman functors}

\begin{defn}\label{defn:HCF}
  Using Lemma~\ref{lem:upgrade-to-B-action}, we define the contravariant higher Coleman
  functor:
\begin{eqnarray*}
  HC_{w,\lambda}: \ocal(\mf{m}_w,\mf{b}_{M_w})_{\lambda-\alg} & \rightarrow & D( \Mod^{\lambda-\sm}_{B(\qq_p)}(E)) \\
  M & \mapsto &\mathrm{R}\Gamma_w( \Sh^{\tor}_{K^p}, \VBzero(HCS_{w,\lambda}(M)))%
\end{eqnarray*}
and the contravariant cuspidal higher Coleman functor:
\begin{eqnarray*}
  HC_{\cusp, w,\lambda}: \ocal(\mf{m}_w,\mf{b}_{M_w})_{\lambda-\alg} & \rightarrow & D(\Mod^{\lambda-\sm}_{B(\qq_p)}(E)) \\
  M & \mapsto &\mathrm{R}\Gamma_w( \Sh^{\tor}_{K^p}, \VBzero(HCS_{w,\lambda}(M)(-D))).%
\end{eqnarray*}
\end{defn}

We extend these functors  to the derived category $$ HC_{ w,\lambda}, HC_{\cusp, w,\lambda}: D^b(\ocal(\mf{m}_w,\mf{b}_{M_w})_{\lambda-\alg}) \rightarrow  D( \Mod^{\lambda-\sm}_{B(\qq_p)}(E))$$ by putting $HC_{w,\lambda}(M) = \mathrm{R}\Gamma_w( \Sh^{\tor}_{K^p}, \VB^{\red}(HCS_{w,\lambda}(M)))$ and similarly in the cuspidal case.
We note that $HCS_{w,\lambda}$ is an exact functor $\ocal(\mf{m}_w,\mf{b}_{M_w})_{\lambda-\alg} \rightarrow \LB_{(\mathfrak{g},B)}(C_{w}^\dag)^{\mathfrak{u}_{\mathfrak{p}}^0}$, and $\VB^0$ is also exact on $\LB_{(\mathfrak{g},B)}(C_{w}^\dag)^{\mathfrak{u}_{\mathfrak{p}}^0}$. Therefore, if $M$ is a complex, $\VB^{\red}(HCS_{w,\lambda}(M))$ is computed by applying $\VB^0(HCS_{w,\lambda}(-))$ to each term of a  complex representing $M$. Moreover,  its $i$-th cohomology sheaf $\underline{\HH}^i(\VB^{\red}(HCS_{w,\lambda}(M))$ is $\VB^0(HCS_{w,\lambda}(\HH^i(M)))$.

\subsubsection{Formal models and the Hodge--Tate period map}

Our goal is to compute the cohomological amplitude of the higher Coleman functors. The main source of vanishing is the affineness of the Hodge--Tate period map. To perform the argument,  we need to consider formal models for some of the spaces and sheaves introduced so far. 
 We can consider the following diagram  where $\Sh^{\star}_{K^p}$ is the minimal compactification. 
\begin{eqnarray*}
\xymatrix{ \Sh^{\tor}_{K^p} \ar[d]\ar[rd]^{\pi_{HT}}& \\
\Sh^{\star}_{K^p} \ar[r]^{\pi^\star_{HT}} & \mathcal{FL}}
\end{eqnarray*}%

 In \cite{boxer2021higher}, sect. 4.4.31, we constructed a formal model of this diagram:  

\begin{eqnarray*}
\xymatrix{ \mathfrak{Sh}^{\tor,mod}_{K^p} \ar[d]\ar[rd]& \\
\mathfrak{Sh}^{\star,mod}_{K^p} \ar[r] & \mathfrak{FL}}
\end{eqnarray*}
Moreover $\mathfrak{Sh}^{\tor,mod}_{K^p} = \lim_{K_p} \mathfrak{Sh}^{\tor,mod}_{K_pK^p}$ and $\mathfrak{Sh}^{\star,mod}_{K^p} = \lim_{K_p} \mathfrak{Sh}^{\star,mod}_{K_pK^p}$ for all $K_p$ small enough, where $\mathfrak{Sh}^{\tor,mod}_{K_pK^p}$ is a formal model of $\Sh^{\tor}_{K_pK^p}$ and $\mathfrak{Sh}^{\star,mod}_{K_pK^p}$ is a formal model of $\Sh^{\star}_{K_pK^p}$.

Let $\mathcal{U}$ be a quasi-compact  open of $\mathcal{FL}$. Let
$\mathfrak{U} \rightarrow \mathfrak{FL}$ be a formal model for the map
$\mathcal{U} \rightarrow \mathcal{FL}$. Using the notation of  \cite{boxer2021higher}, sect. 4.4.31, we have a formal model for $\pi_{HT}^{-1}(\mathcal{U})$, denoted by  $\mathfrak{Sh}^{\tor,mod}_{K^p,\mathfrak{U}}$  as well as a formal model for $(\pi^\star_{HT})^{-1}(\mathcal{U})$ denoted by  $\mathfrak{Sh}^{\star,mod}_{K^p,\mathfrak{U}}$. There is a diagram: 

\begin{eqnarray*}
\xymatrix{ \mathfrak{Sh}^{\tor,mod}_{K^p, \mathfrak{U}} \ar[d]\ar[rd]& \\
\mathfrak{Sh}^{\star,mod}_{K^p,\mathfrak{U}} \ar[r] & \mathfrak{U}}
\end{eqnarray*}
 
Moreover, we have $\pi_{HT}^{-1}(\mathcal{U}) = \lim_{K_p}
\pi_{HT}^{-1}(\mathcal{U})_{K_p}$ for $K_p$ small enough and
accordingly $\mathfrak{Sh}^{\tor,mod}_{K^p,\mathfrak{U}}= \lim_{K_p}
\mathfrak{Sh}^{\tor,mod}_{K_pK^p,\mathfrak{U}}$. Similarly, for $K_p$ small enough $(\pi^\star_{HT})^{-1}(\mathcal{U}) = \lim_{K_p}( \pi^\star_{HT})^{-1}(\mathcal{U})_{K_p}$  and  $\mathfrak{Sh}^{\star,mod}_{K^p,\mathfrak{U}}=\lim_{K_p} \mathfrak{Sh}^{\star,mod}_{K_pK^p,\mathfrak{U}}$.

\subsubsection{Formal Banach sheaves and formal Smith sheaves} In this
section we consider a flat $p$-adic formal scheme $\mathfrak{X}$ which locally of topologically finite
type  over $\Spf~\ocal_E$   %
 (in other words, Zariski locally, $\mathfrak{X}$ is $\Spf~A$ where $A$ is a
 quotient of  $\ocal_E\langle X_1, \cdots, X_n \rangle$). We write~$\cX$for the
 generic fiber of~ $\mathfrak{X}$.  The sheaf $\oscr_{\mathfrak{X}}$ is a sheaf
 of solid $\ocal_E$-modules, as are all the sheaves we will consider.   We define the scheme   $X_n = \mathfrak{X} \times_{\Spf~\ocal_E} \Spec~\ocal_E/p^n\ocal_E$.
\begin{defn} \leavevmode
\begin{enumerate}
\item  A locally trivial formal Banach sheaf over $\mathfrak{X}$ is a sheaf
  $\mathfrak{F}$ of $\oscr_{\mathfrak{X}}$-modules  which  is flat as
  an $\ocal_E$-module, and such that  $\mathfrak{F} = \lim_n \mathscr{F}_n$ for $\mathscr{F}_n = \mathfrak{F}/p^n
  \mathfrak{F}$,  and there exists a covering $\mathfrak{X} = \cup_i \mathfrak{V}_i$ and   sets
  $I_i$ and such that
  $\mathscr{F}_n\vert_{V_{i,n}} = \oscr_{V_{i,n}}
  \otimes_{\ocal_E/p^n\ocal_E}(\oplus_{ s\in I_i} \ocal_E/p^n
  \ocal_E)$ and the transition maps $\mathscr{F}_n\vert_{V_{i,n}} \rightarrow \mathscr{F}_{n-1}\vert_{V_{i,{n-1}}}$ are the obvious ones.

\item A  very  small formal Banach sheaf  over $\mathfrak{X}$ is a sheaf
  $\mathfrak{F}$ of $\oscr_{\mathfrak{X}}$-modules which  is flat as
  an $\ocal_E$-module, and such that  $\mathfrak{F} = \lim_n \mathscr{F}_n$ for $\mathscr{F}_n = \mathfrak{F}/p^n
  \mathfrak{F}$, and $\mathscr{F}_1 = \mathscr{G}_1 \otimes_{\ocal_E} (\oplus_{ s\in I} \ocal_E/p)$ for some coherent sheaf $\mathscr{G}_1$ and some set $I$. 
\end{enumerate}
\end{defn}

 \begin{defn} \leavevmode
\begin{enumerate}
\item A locally trivial formal Smith sheaf over $\mathfrak{X}$ a sheaf
  $\mathfrak{F}$ of $\oscr_{\mathfrak{X}}$-modules which  is flat as an
  $\ocal_E$-module,  such that  $\mathfrak{F} = \lim_n \mathscr{F}_n$ for $ \mathscr{F}_n =   \mathfrak{F}/p^n
  \mathfrak{F}$, and  there exists  a covering $\mathfrak{X} = \cup_{i} \mathfrak{V}_i$ and sets $I_i$ such that
  $\mathscr{F}_n = \oscr_{V_{i,n}}
  \otimes_{\ocal_E/p^n\ocal_E}(\ocal_E/p^n \ocal_E)^{I_i}$ with the obvious transition maps $\mathscr{F}_n\vert_{V_{i,n}}\rightarrow \mathscr{F}_{n-1}\vert_{V_{i,n-1}}$. %
  \item  A  very small formal Smith sheaf over $\mathfrak{X}$ is a sheaf
  $\mathscr{F}$ of $\oscr_{\mathfrak{X}}$-modules which  is flat as an
  $\ocal_E$-module,  such that  $\mathfrak{F} = \lim_n \mathscr{F}_n$ for $ \mathscr{F}_n =   \mathfrak{F}/p^n
  \mathfrak{F}$, and $\mathscr{F}_1 = \mathscr{G}_1 \otimes_E (\oscr_{E}/p \oscr_E)^I$ for a coherent sheaf $\mathscr{G}_1$ and some set $I$. 
\end{enumerate}
\end{defn}

 If $\mathfrak{F}$ is a locally trivial formal Smith sheaf, then it
 has a generic fiber $\mathscr{F}$ over $\mathcal{X}$ defined as
 follows: we take a covering $\mathcal{X} = \cup_i \mathfrak{V}_i$,
 such that $\mathfrak{F}\vert_{\mathfrak{V}_i} =
 \oscr_{\mathfrak{V}_i} \otimes \ocal_E^{I_i}$. %
 Over the generic fibre $\mathcal{V}_i$ of~$\mathfrak{V}_i$, we define $\mathscr{F}\vert_{\mathcal{V}_i} = \oscr_{\mathcal{V}_i} \otimes_E( \ocal_E^{I_i}[1/p])$ and we use the gluing data of $\mathfrak{F}$ to glue the $\mathscr{F}\vert_{\mathcal{V}_i}$. 

A similar construction applies to a locally trivial Banach sheaf. See
also \cite[Thm. 2.5.9]{boxer2021higher} for a more general
statement. %

\begin{prop}\label{prop:small-ample-vanishing-affinoid} If $\mathfrak{F}$ is a very small formal Banach sheaf or a very small formal Smith sheaf, $\mathfrak{X}$ admits an ample invertible sheaf and $\mathcal{X}$ is affinoid, then $\HH^i(\mathfrak{X}, \mathscr{F}) \otimes_{\ocal_E} E =0$ for all $i >0$. 
\end{prop}
\begin{proof} The formal Banach case is
  \cite[Thm. 2.5.8]{boxer2021higher}. The same proof with minor modification applies to the Smith case. 
\end{proof}

\subsubsection{Cohomological amplitude of the higher Coleman functors}

Let $M \in \mathcal{O}_{\lambda-\alg}$. Let us consider the sheaf
$\cF:=\VBzero(HCS_{w,\lambda}(M))$. Theorem \ref{thm-VB} addresses the local
structure of this sheaf. We also need to produce some integral structure, as in
the following lemma.

\begin{lem}\label{lem-exists-model}  There exists a quasi-compact open subset
  $U_0 \subseteq C_w$ containing $w$ such that if we set   $V_0 =
  \pi_{HT}^{-1}(U_0)$, then we can write  $\mathscr{F}\vert_{V_0} = \colim i_m^{-1} \mathscr{F}_m$ where:
\begin{enumerate} 
\item $U_0 = \cap_m U_{0,m}$ where $\{U_{0,m}\}_{m \in \ZZ_{\geq 0}}$ is a decreasing sequence of quasi-compact open affinoid subsets of $\mathcal{FL}$; 
\item $V_{0,m} = \pi_{HT}^{-1}(U_{0,m})$ is a quasi-compact open subset of $\Sh^{\tor}_{K^p}$, stable under  a compact open subgroup $K_{p,m} \subseteq G(\qq_p)$;
\item  $i_m : V_0 \rightarrow V_{0,m}$ is the natural inclusion;
\item $\mathscr{F}_m$ is a sheaf over $V_{0,m}$;
\item $V_{0,m,K_{p,m}} \hookrightarrow \Sh^{\tor}_{K^pK_{p,m}}$ is a quasi-compact open which  descends $V_{0,m}$ to finite level $K_{p,m}$, and $\pi_{K_{p,m}} : V_{0,m} \rightarrow V_{0,m,K_{p,m}}$ is the induced projection;
\item We have a Banach sheaf $\mathscr{G}_{m,K_{p,m}}$ on $V_{0,m,K_{p,m}}$ and $\mathscr{F}_m =   \pi_{K_{p,m}}^{-1} \mathscr{G}_{m,K_{p,m}} \otimes_{\pi_{K_{p,m}}^{-1}\oscr_{V_{0,m,K_{p,m}}}} \oscr^{\sm}_{V_{0,m}}$;
\item The sheaves $\mathscr{G}_{m,K_{p,m}}$ arise as the generic fibers of
locally trivial,  very small, formal Banach sheaves
$\mathfrak{G}_{m,K_{p,m}}$ over
$\mathfrak{Sh}^{\tor,mod}_{K^pK_{p,m},\mathfrak{U}_{0,m}}$.
\end{enumerate}\end{lem}%
\begin{proof}  The first $6$ points follow from Theorem \ref{thm-VB}. We  need to give a  more explicit construction of $\mathscr{G}_{m,K_{p,m}}$ in order to be able  to produce an integral structure. We will follow closely  the proof of 
  \cite[Lem.\ 6.6.2]{boxer2021higher}. As $M \in
  \ocal_{\lambda-\alg}$, we deduce that $\hat{M}^\vee(\lambda) =
  \colim_r M_r$ where each $M_r$ is a Banach space representation of
  the group denoted $\Stab(w)_{r,s}$ in Definition
  \ref{defn-Stabw-rep-Q} (with~$Q=B$). Unraveling the definition, we deduce in particular that $M_r$ is a Banach space representation of $M_{w,r}U_{M_w,s} =: M_{w,r,s} \hookrightarrow M_w$. Moreover, after possibly changing $r$ and $s$, we can also assume that there is a lattice $M_r^+ \subseteq M_r$  with the property that the co-action map $M_r \rightarrow M_r \otimes \oscr_{M_{w,r,s}}$ induces a map  $M^+_r \rightarrow M^+_r \otimes \oscr^+_{M_{w,r,s}}$, trivial on $M_{r}^+/pM_{r}^+$. 
We consider the $M_{w,r,s}$-torsor $ w (U_{P_w} \cap G_rU_s)\backslash
G_rU_s   \rightarrow w ({P_w} \cap G_rU_s)\backslash  G_rU_s$. This is
a reduction of structure of the standard $M$-torsor $ U_{P} \backslash
G \rightarrow P \backslash G$ over $w (P_w \cap G_rU_s)\backslash
G_rU_s$. Let us write $\mathcal{U}_r = w (P_w \cap G_rU_s)\backslash
G_rU_s$, which is a quasi-compact open subset of
$\mathcal{FL}$. %
 Over $\pi_{HT}^{-1}(\mathcal{U}_r)$, we pull back  to get a  $M_{w,r,s}$-torsor that we denote by $\mathcal{M}_{dR, r,s, \mathcal{U}_r}$. For $K_p$ small enough,  by \cite[prop. 4.6.12]{boxer2021higher}, it descends to a  
$M_{w,r,s}$-torsor $\mathcal{M}_{dR, r,s, \mathcal{U}_r, K_p}$ over
$\pi_{HT}^{-1}(\mathcal{U}_r)_{K_p}$. Moreover, by
\cite[prop. 4.6.15]{boxer2021higher} %
for any affinoid open subset
$V_{K_p} \subseteq \pi_{HT}^{-1}(\mathcal{U})_{K_p}$, we can find
$K_p' \subseteq K_p$ such that the torsor $\mathcal{M}_{dR, r,s,
  \mathcal{U}_r, K'_p}\vert_{V_{K_p'}}$ is trivial.  After
rescaling, %
we can assume that $U_{0,r} \subseteq \mathcal{U}_r$. %
We now attach to $M_r$ a locally projective small formal Banach sheaf
on the formal model %
$\mathfrak{Sh}^{\tor,mod}_{K^pK_{p,r},\mathfrak{U}_{0,r}}$  of $V_{0,r,K_{p,r}} = \pi_{HT}^{-1}({U}_{0,r})_{K_{p,r}}$ for a small enough $K_{p,r}$. Indeed, we pick a finite affine covering $\cup_i \mathfrak{V}_i = \Spf(A_i)$ of $\mathfrak{Sh}^{\tor,mod}_{K^pK_{p,r},\mathfrak{U}_{0,r}}$. After replacing $K_{p,r}$ by a smaller compact open, we can assume that $\mathcal{M}_{dR, r,s, {U}_{0,r}, K_{p,r}}$ is trivial over the generic fiber of  every $ \mathfrak{V}_i$. It follows that the torsor is described by a 1-cocycle $\{m_{i,j} \in M_{w,r,s}( (A_{i,j}[1/p], A_{i,j}))\}_{i,j}$ (where $\mathfrak{V}_{i,j} = \mathfrak{V}_i \cap \mathfrak{V}_j = \Spf(A_{i,j})$). We can use the elements  $m_{i,j}$ to glue the trivial formal Banach sheaf $\oscr_{\mathfrak{V}_i} \otimes_{\ocal_L} M_{r}^+$ to get the very small, locally trivial,  formal Banach sheaf $\mathfrak{G}_{r,K_{p,r}}$. The very smallness property comes from the fact that the $m_{i,j}$ act trivially on $M_{r}^+/p$. 
\end{proof}

We give a technical variant of this description, using formal Smith sheaves instead. 

\begin{lem}\label{lem-exists-model2} There exists a quasi-compact open subset
  $U_0 \subseteq C_w$ containing $w$ such that if we set   $V_0 =
  \pi_{HT}^{-1}(U_0)$, then we can write  $\mathscr{F}\vert_{V_0} = \colim i_m^{-1} \mathscr{F}_m$ where:
\begin{enumerate} 
\item $U_0 = \cap_m U_{0,m}$ where $\{U_{0,m}\}_{m \in \ZZ_{\geq 0}}$ is a decreasing sequence of quasi-compact open affinoid subsets of $\mathcal{FL}$, 
\item $V_{0,m} = \pi_{HT}^{-1}(U_{0,m})$ is a quasi-compact open subset of $\Sh^{\tor}_{K^p}$, stable under  a compact open subgroup $K_{p,m} \subseteq G(\qq_p)$,
\item  $i_m : V_0 \rightarrow V_{0,m}$ is the natural inclusion,
\item $\mathscr{F}_m$ is a sheaf over $V_{0,m}$, 
\item $V_{0,m,K_{p,m}} \hookrightarrow \Sh^{\tor}_{K^pK_{p,m}}$ is a quasi-compact open which  descends $V_{0,m}$ to finite level $K_{p,m}$, and $\pi_{K_{p,m}} : V_{0,m} \rightarrow V_{0,m,K_{p,m}}$ is the induced projection,
\item We have a  sheaf $\mathscr{G}'_{m,K_{p,m}}$ on $V_{0,m,K_{p,m}}$ and $\mathscr{F}_m =   \pi_{K_{p,m}}^{-1} \mathscr{G}_{m,K_{p,m}} \otimes_{\pi_{K_{p,m}}^{-1}\oscr_{V_{0,m,K_{p,m}}}} \oscr^{\sm}_{V_{0,m}}$. 
\item The sheaves $\mathscr{G}'_{m,K_{p,m}}$ arise as the generic fibers of
locally trivial,  very small, formal Smith sheaves
$\mathfrak{G}_{m,K_{p,m}}$ over
$\mathfrak{Sh}^{\tor,mod}_{K^pK_{p,m},\mathfrak{U}_{0,m}}$.
\end{enumerate}

\end{lem}

\begin{proof} The argument is almost identical to the proof of Lemma \ref{lem-exists-model}. We have that  $\hat{M}^\vee(\lambda) = \colim_r M_r$  where each $M_r$ is a Banach space representation of the group denoted $\Stab(w)_{r,s}$ in Definition \ref{defn-Stabw-rep-Q}. Since this is a $LB$-space of compact type, it also admits a presentation as a $LS$-space of compact type, $\hat{M}^\vee(\lambda) = \colim_r M'_r$  where each $M'_r$ is a Smith space representation of $\Stab(w)_{r,s}$ (see \cite[Cor. 3.38]{MR4475468}). We then pick lattices $(M'_r)^+$ in $M'_r$ and glue the sheaves  $\oscr_{\mathfrak{V}_i} \otimes_{\ocal_L} (M'_{r})^+$ to get the very small formal Smith sheaf $\mathfrak{G}'_{r,K_{p,r}}$.
\end{proof}

For any $K_p \subseteq K_{p,m}$, we have maps $\pi_{K_p, K_{p,m}} :
\mathfrak{Sh}^{\tor,mod}_{K^pK_{p},\mathfrak{U}_{0,m}} \rightarrow
\mathfrak{Sh}^{\tor,mod}_{K^pK_{p,m},\mathfrak{U}_{0,m}}$ and  we write
$\mathfrak{G}_{m,K_{p}} := \pi_{K_p, K_{p,m}}^\star \mathfrak{G}_{m, K_{p,m}}$
(a locally trivial, very small formal Banach sheaf) and $\mathfrak{G}'_{m,K_{p}}
:= \pi_{K_p, K_{p,m}}^\star \mathfrak{G}'_{m, K_{p,m}}$ (a locally trivial, very
small formal Smith sheaf). We denote their generic fibers by $\mathcal{G}_{m,K_{p}}$ and $\mathcal{G}'_{m,K_p}$.

\begin{lem}\label{lem-seed-vanishing}\leavevmode
\begin{enumerate}%
\item  For any affinoid $U' \subseteq U_{0,m}$, any compact open subgroup $K_{p} \subseteq K_{p,m'}$  fixing $U'$, and letting $V' = \pi_{HT}^{-1}(U')$,   we have $\HH^i(V'_{K_{p}}, \mathscr{G}_{m,K_{p}}(-D)) =0$ for all $i >0$.
\item For any Stein space $S' \subseteq U_{0,m}$ stable under a compact open subgroup $K_{p} \subseteq K_{p,m'}$, let $V'' = \pi_{HT}^{-1}(S')$. We have $\HH^i_c(V_{K_p}'', \mathscr{G}'_{m,K_{p}}) =0$ for all $i \neq d$.
\end{enumerate}
\end{lem}
\begin{proof} The first part follows as in lemma 6.6.2 of
  \cite{boxer2021higher}.  We briefly recall the argument. We take a
  formal model $\mathfrak{U}' \rightarrow \mathfrak{FL}$ of $\cU'\to\mathcal{FL}$. By Lemma \ref{lem-exists-model},  we have a locally trivial, very small, formal Banach sheaf $\mathfrak{G}_{m, K_p}$ over $\mathfrak{Sh}^{\tor}_{K_pK^p, \mathfrak{U}'}$.
Let $\pi: \mathfrak{Sh}^{\tor}_{K_pK^p, \mathfrak{U}'} \rightarrow
\mathfrak{Sh}^{\star}_{K_pK^p, \mathfrak{U}'}$.  We have that
$\mathrm{R}^i\pi_\star \mathfrak{G}_{m, K_p}(-D) = 0 $ for all $i
>0$. Indeed, using the very smallness,  this reduces to the vanishing
(\cite[Thm. 4.4.37]{boxer2021higher}) of $\mathrm{R}^i \pi_\star
\oscr_{\mathfrak{Sh}^{\tor}_{K_pK^p, \mathfrak{U}'}}(-D)$. We deduce
that $\pi_{\star} \mathfrak{G}_{m, K_p}(-D)$ is a very small  formal Banach
sheaf and thus Proposition \ref{prop:small-ample-vanishing-affinoid} implies that 
$\HH^i( \mathfrak{Sh}^{\tor}_{K_pK^p, \mathfrak{U}'} ,
\mathfrak{G}_{m, K_p}(-D))\otimes E =0$ for all $i >0$.

We claim that $\HH^i( \mathfrak{Sh}^{\tor}_{K_pK^p, \mathfrak{U}'} ,
\mathfrak{G}_{m, K_p}(-D))\otimes E = \HH^i(V'_{K_{p}},
\mathscr{G}_{m,K_{p}}(-D))$.  To see this, take a Zariski open affine cover $\{\mathfrak{W}_s\}_{s}$ of $ \mathfrak{Sh}^{\tor}_{K_pK^p, \mathfrak{U}'}$ with the property that $\mathfrak{G}_{m, K_p}(-D)\vert_{\mathfrak{W}_s}$ is a trivial formal Banach sheaf. 
By Corollary \ref{coro-coho-quasi-steinV} $\HH^i(V'_{K_{p}}, \mathscr{G}_{m,K_{p}}(-D))$ is computed by the \v{C}ech complex  associated to the generic fiber of this cover, which computes $\HH^i( \mathfrak{Sh}^{\tor}_{K_pK^p, \mathfrak{U}'} ,
\mathfrak{G}_{m, K_p}(-D))\otimes E$. We
deduce that %
$\HH^i(V'_{K_{p}}, \mathscr{G}_{m,K_{p}}(-D)) =0$ for all
$i >0$, as required.

The second part follows from a certain form of  duality. Let us define
the ``Serre  dual'' of $ \mathscr{G}'_{m,K_{p}}$, %
$$D(\mathscr{G}'_{m,K_{p}}) =
\underline{\mathrm{Hom}}_{\oscr_{V_{0,m,K_p}}}( \mathscr{G}'_{m,K_{p}},
\Omega^d_{V_{0,m,K_{p}}}),$$  by which we mean the following.
We take a finite Stein  covering $V''_{K_p} = \cup_i V''_{K_p,i} $ with the
property that  $\Omega^d_{V_{m,K_{p}}}\vert_{V''_{K_p,i}}$ is trivial and that
$\mathscr{G}'_{m,K_{p}} = \oscr_{V''_{K_p,i}} \otimes M_i$ for a Smith space
$M_i$. Then $D(\mathscr{G}'_{m,K_{p}}) \vert_{V''_{K_p,i}}=\oscr_{V''_{K_p,i}} \otimes M^\vee_i$ where $M^\vee_i$ is a Banach space. We consider the Čech complex of Fr\'echet spaces:
$$C:  \prod_i \HH^0(V''_{K_p,i}, D(\mathscr{G}'_{m,K_{p}})) \rightarrow \prod_{i,j} \HH^0(V''_{K_p,i} \cap V''_{K_p,i} , D(\mathscr{G}'_{m,K_{p}})) \rightarrow \cdots$$
as well as the Čech complex of $LS$-spaces:  $$D: \HH^d_c( \cap_i V''_{K_p,i}, \mathscr{G}'_{m,K_{p}}) \rightarrow \prod_{j} \HH^d_c( \cap_{i \neq j} V''_{K_p,i}, \mathscr{G}'_{m,K_{p}}) \rightarrow  \cdots$$
The two complexes are termwise dual of each other by the duality theory  (see Corollary \ref{coro-dualitystein}).  The complex $D$ computes $\mathrm{R}
\Gamma_c(V''_{K_p}, \mathscr{G}'_{m,K_{p}})$ and the complex $C$ computes $\mathrm{R}\Gamma(V''_{K_p}, D(\mathscr{G}'_{m,K_{p}}))$ by Corollary \ref{cor-cohosteinV} and \ref{coro-dualitystein}. 
Thus, it suffices to prove that the complex $C$ has cohomology concentrated in
degree $0$.   Write $S' = \cup S'_{r}$ as a countable increasing union of
affinoids, and let $V''_{r} = \pi_{HT}^{-1}(S'_{r})$ and $V''_{r} = \cup_i V''_{r,i}$. It follows from Lemma  \ref{lem-exists-model2} that   $D(\mathscr{G}'_{m,K_{p}})$ admits a formal model which is a very small formal Banach sheaf   over  a formal model of $V''_{r,K_p}$.  We deduce that  $\mathrm{R}\Gamma(V''_{r,K_p}, D(\mathscr{G}'_{m,K_{p}}))$ is concentrated in degree $0$ by the same argument as for the first point of the lemma. Thus, the Čech complex $C_r$ of $D(\mathscr{G}'_{m,K_{p}})$ with respect to $\{ V''_{r,i,K_p} \}$ has only cohomology in degree $0$.   We then  use topological Mittag-Leffler to deduce that $C = \lim_r C_r$ is concentrated in degree $0$. 
\end{proof}

\begin{thm}\label{thm-coho-amplitudeHC} For any $M \in \ocal_{\lambda-\alg}$, $HC_{\cusp, w,\lambda}(M)$ has amplitude $[0,\ell(w)]$ and $HC_{w,\lambda}(M)$ has amplitude $[\ell(w), d]$. 
\end{thm}

\begin{proof} We first prove the claim regarding $HC_{\cusp, w,\lambda}(M)$. Let us denote $\mathscr{F} = \VBzero(HCS_{w,\lambda}(M))(-D)$. 
Let us fix a closed  neighborhood $Z_0$ of $w$ in $C_w$ and an open
neighborhood $U_0$ of $Z_0$ in $C_w$ as in definition
\ref{example-chooseopen}. Let $t \in T^{++}(\qq_p)$. %
It is easy to see that $\cup_{n\in \Z_{\geq 0}} Z_0 t^n =
C_w$. As above, we put $V_0 = \pi_{HT}^{-1}(U_0)$ and $F_0 =
\pi_{HT}^{-1}Z_0$. By Lemma \ref{lem-computeasinductivelimit} (noting that $Z_0t^n$ and $Z_n$ are cofinal), we have
that $\mathrm{R}\Gamma_w( \Sh^{\tor}_{K^p}, \mathscr{F}) = \colim_n
\mathrm{R}\Gamma_{F_0t^n} (V_0t^n, \mathscr{F})$. %
It follows that each  $\HH^i_w( \Sh^{\tor}_{K^p}, \mathscr{F})$ is generated as a $B(\qq_p)$-module by the image of  $\HH^i_{F_0}(V_0, \mathscr{F})$ in $\HH^i_w( \Sh^{\tor}_{K^p}, \mathscr{F})$. Thus it suffices to see that 
$\mathrm{R}\Gamma_{F_0}(V_0, \mathscr{F})$ has amplitude $[0,
\ell(w)]$.

We are free to chose $n_0$ in definition \ref{example-chooseopen}, and
to make $F_0$ and $V_0$ arbitrarily small. By Lemma
\ref{lem-exists-model},  we can assume that $\mathscr{F}\vert_{V_0 } =
\colim \mathscr{F}_m$ where $\mathscr{F}_m = \mathscr{G}_{m,K_{p,m}}
\otimes_{\oscr_{V_{0,m,K_{p,m}}}} \oscr_{V_{0,m}}(-D)$ %
where $\mathscr{G}_{m,K_{p,m}}$ is a Banach sheaf  admitting a locally trivial, very small, formal model over $V_{0,m,K_{p,m}}$. %
By Lemma \ref{lem-second-comp-inductivelimi},  we have that $\mathrm{R}\Gamma_{F_0}(V_0, \mathscr{F})  = \colim_m \mathrm{R}\Gamma_{F_0}(V_0, \mathscr{F}_m)$ and 
$\mathrm{R}\Gamma_{F_0}(V_0, \mathscr{F}_m) = \colim_{s,n} \mathrm{R}\Gamma_{F_{0,s,K_{p,n}}}(V_{0,s, K_{p,n}}, \mathscr{G}_{m,K_{p,n}}(-D))$ where $K_{p,n}$ tends to $\{e\}$. Thus it suffices to prove that 
$\mathrm{R}\Gamma_{F_{0, s,K_{p,n}}}(V_{0,s, K_{p,n}},
\mathscr{G}_{m,K_{p,n}}(-D))$ has cohomological amplitude in $[0,
\ell(w)]$. Recall from definition  \ref{example-chooseopen} that  $F_{0,s} =
\pi_{HT}^{-1}(Z_{0,s})$ and $V_{0,s} = \pi_{HT}^{-1}(U_{0,s})$ for $s$ large
enough (and  a fixed choice of $n_0$ big enough). We observe that $U_{0,s}$ is
affinoid  while $U_{0,s} \setminus Z_{0,s}$ is covered by $\ell(w)$ affinoids. 
Using the triangle: 
$$\mathrm{R}\Gamma_{F_{0,s,K_{p,n}}}(V_{0,s, K_{p,n}},  \mathscr{G}_{m,K_{p,n}}(-D)) \rightarrow \mathrm{R}\Gamma(V_{0,s, K_{p,n}},  \mathscr{G}_{m,K_{p,n}}(-D)) $$ $$\rightarrow 
\mathrm{R}\Gamma(V_{0,s, K_{p,n}} \setminus F_{0,s,K_{p,n}},  \mathscr{G}_{m,K_{p,n}}(-D)) \stackrel{+1}\rightarrow $$
together with Lemma \ref{lem-seed-vanishing}, we arrive at the desired
conclusion.

We now turn to the case of usual cohomology, which follows along
similar lines to the above, using the presentation via cohomology with
compact support.   Let us denote  now $\mathscr{F} =
\VBzero(HCS_{w,\lambda}(M))$. We first see  by lemma
\ref{lem-computeasinductivelimitI} and a similar argument  using the
$B(\qq_p)$-action, that it is enough to check that
$\mathrm{R}\Gamma_c(Y_0, \mathscr{F})$ is concentrated in degrees
$[\ell(w), d]$ where $Y_0= \pi_{HT}^{-1} X_0$ and $X_0$ is a Stein open
neighborhood of $w$ in $C_w$ (see Definition
\ref{example-open-coho-support}). Next, we use Lemma
\ref{lem-exists-model2} to see  that $\mathscr{F}\vert_{Y_0} = \colim
\mathscr{F}'_m$ where $\mathscr{F}'_m  = \mathscr{G}'_{m,K_{p,m}}
\otimes_{\oscr_{Y_{0,m,K_{p,m}}}} \oscr_{Y_{0,m}}$ %
and $Y_{0,m} = \pi_{HT}^{-1}(X_{0,m})$  for
$X_{0,m}$  a Stein neighborhood of $X_0$ in $\mathcal{FL}$. We deduce
that $\mathrm{R}\Gamma_c(Y_0,  \mathscr{F}) = \colim
\mathrm{R}\Gamma_c(Y_0,  \mathscr{F}_m')$.  Let  $\mathscr{G}'_{m} =
\colim_{K_p} \mathscr{G}'_{m,K_p}$ where $\mathscr{G}'_{m,K_p}$ is the pull back of $\mathscr{G}'_{m,K_{p,m}}$ to $Y_{0,m,K_{p}}$ for $K_p \subseteq K_{p,m}$ and $\mathscr{G}'_{m}$ is viewed as a sheaf on $Y_{0,m}$.  %
We deduce from lemma  \ref{lem-computeasinductivelimitII} that we have
a triangle $$\mathrm{R}\Gamma_c(Y_{0,m}\setminus Y_0, \mathscr{G}'_m) \rightarrow \mathrm{R}\Gamma_c(Y_{0,m}, \mathscr{G}'_m) \rightarrow \mathrm{R}\Gamma_c( Y_0, \mathscr{F}'_m) \stackrel{+1}\rightarrow .$$
Recall that  $W_{0,m,r} = \pi_{HT}^{-1}(T_{0,m,r})$  (see definition \ref{example-open-coho-support}).
We have that $\mathrm{R}\Gamma_c(Y_{0,m}, \mathscr{G}'_m) = \colim_{K_p}
\mathrm{R}\Gamma_c(Y_{0,m,K_p}, \mathscr{G}'_{m,K_p})$ is concentrated in degree $d$ by Lemma \ref{lem-seed-vanishing}. 
Similarly, $\mathrm{R}\Gamma_c(Y_{0,m} \setminus Y_0, \mathscr{G}'_m)
= \colim_r \colim_{K_p} \mathrm{R}\Gamma_c( W_{0,m,r,K_p},
\mathscr{G}'_{m,K_p})$, and $\mathrm{R}\Gamma_c( W_{0,m,r,K_p},
\mathscr{G}'_{m,K_p})$ has cohomology  concentrated in degrees
$[\ell(w)+1, d]$  by Lemma \ref{lem-seed-vanishing} (as $T_{0,m,r}$  is the union of $\ell(w)$ Stein spaces).
\end{proof}

\subsubsection{Finite slope projector}\label{sect-finite-slope-proj}

We follow the notation introduced in Section~\ref{notation-torus}. Let $\mathcal{Z}$ be the character space of $T(\qq_p)$. Fixing an
isomorphism $T(\qq_p) = \ZZ^r \times T(\Z_p)$ (with $r$ the rank of $T^d$, the maximal split torus in $T$), we see that
$\mathcal{Z} = \mathcal{W} \times (\mathbb{G}_m^{an})^r$ with
$\mathcal{W} = \Spa (\Z_p\llb T(\Z_p)\rrb , \Z_p\llb T(\Z_p)\rrb)
\times_{\Spa(\Z_p)} \Spa(\qq_p, \Z_p)$. %

 We fix an affinoid increasing covering of the Stein space $\mathcal{Z} = \cup_n
 \mathcal{Z}_n$. As in \cite{andreychev2021pseudocoherent}, we let
 $D(\mathcal{Z}_n)$ denote the category $D( (\oscr_{\mathcal{Z}_n},
 \oscr_{\mathcal{Z}_n}^+)_{\blacksquare})$, %
and we let $D(\mathcal{Z})$ be the derived category of
 quasi-coherent sheaves over $\mathcal{Z}$ (see
 \cite[Thm. 1.6]{andreychev2021pseudocoherent}). 
 We have defined in Section \ref{notation-torus} the monoids $T^{++}(\qq_p) \subseteq T^+(\qq_p) \subseteq T(\qq_p)$.
Note that $T^+(\qq_p) = T(\Z_p) \times
\Z^{s}\times\ZZ_{\geq 0}^{r-s}$, where~$s$ is the rank of the maximal
split torus in~$Z(G)$. %
We let $(\qq_p)_{\blacksquare}[T^+(\qq_p)]$ be the solidification of
$\qq_p[T^+(\qq_p)]$,  and we let $(\qq_p)_{\blacksquare}[T(\qq_p)]$ the solidification of $\qq_p[T(\qq_p)]$. 
The categories of  solid modules over these rings  are denoted by $\Mod_{T^+(\qq_p)}(\qq_p)$ and $\Mod_{T(\qq_p)}(\qq_p)$ (this is consistent with definition \ref{def-locally-pro-rep}). 
Their derived categories are denoted  $D(\Mod_{T^+(\qq_p)}(\qq_p))$ and $D(\Mod_{T(\qq_p)}(\qq_p))$. 

We define a functor $f_n^\star:  D( \Mod_{T^+(\qq_p)}(\qq_p)) \rightarrow
D(\mathcal{Z}_n)$, $f_n^\star M= M \otimes_{ (\qq_p)_{\blacksquare}[T^+(\qq_p)]}
(\oscr_{\mathcal{Z}_n}, \oscr_{\mathcal{Z}_n}^+)_{\blacksquare}$. This functor
has a right adjoint given by the forgetful  functor $(f_n)_\star: D(\mathcal{Z}_n) \rightarrow D( \Mod_{T(\qq_p)}(\qq_p))$. 
These functors induce adjoint functors $f^\star: D( \Mod_{T^+(\qq_p)}(\qq_p)) \rightarrow  D(\mathcal{Z})$ and $f_\star:  D(\mathcal{Z}) \rightarrow D( \Mod_{T(\qq_p)}(\qq_p))$. 

We finally define the finite slope functor: 
$$(-)^{\nfs} = f_\star f^\star: D( \Mod_{T^+(\qq_p)}(\qq_p))  \rightarrow D( \Mod_{T(\qq_p)}(\qq_p)).$$ The unit of adjunction
gives a natural map $M \rightarrow M^{\nfs}$ in  $D( \Mod_{T^+(\qq_p)}(\qq_p))$. Note that $M^{\nfs} = \lim_n (f_n)_\star f_n^\star M$.
\begin{rem} The finite slope functor is a localization functor, which factors over the functor $- \otimes_{ (\qq_p)_{\blacksquare}[T^+(\qq_p)]}  (\qq_p)_{\blacksquare}[T(\qq_p)]$. 
However, it is a stronger form of localization. For example, let us
(abusively!) consider the case that $T^+(\qq_p) = \ZZ_{\geq 0}$ and $T(\qq_p) = \ZZ$. 
In this case, $(\qq_p)_{\blacksquare}[T^+(\qq_p)] = \qq_p[X]$ and $(\qq_p)_{\blacksquare}[T(\qq_p)] = \qq_p[X,X^{-1}]$.  We can also suppose that $\mathcal{Z}_n = \Spa (\qq_p\langle p^nX, p^n X^{-1} \rangle, \ZZ_p\langle p^nX, p^n X^{-1} \rangle)$.   We claim that the module $\qq_p((X))$  is a solid $\qq_p[X,X^{-1}]$-module, whose finite slope part is trivial. By definition we need  to show that $\qq_p((X)) \otimes^L_{\ZZ_p[X]}  \ZZ_p\langle p^nX, p^n X^{-1} \rangle =0$ for every $n$. 
To see this, we first show note that  $\qq_p((X))$ is a solid $\qq_p\langle
p^{-s}X \rangle$ for every $s \geq 0$.  Indeed  $\qq_p[X]/(X^\ell)$ is a
solid $\qq_p\langle p^{-s}X \rangle$-module for any $\ell,s \in \Z_{\geq 0}$
and so $\qq_p\llb X\rrb = \lim_\ell \qq_p[X]/(X^\ell)$ is also a $\qq_p\langle
p^{-s}X \rangle$-module, and thus $\qq_p((X))$ is a solid $\qq_p\langle
p^{-s}X \rangle$-module, as claimed.  It remains to observe that  if $s \geq
n$, then $\qq_p \langle p^{-s}X \rangle \otimes^L_{\qq_p[X]} \qq_p\langle p^nX,
p^n X^{-1} \rangle = 0$. To see this, note that this is represented by the following complex:
 $$[  \qq_p \langle U , p^nX, p^n X^{-1} \rangle \stackrel{ p^{s-n}U p^n X^{-1} -1}\longrightarrow  \qq_p \langle U, p^nX, p^n X^{-1} \rangle].$$ But  $1- p^{s-n} U (p^n X^{-1})$ is  invertible, with inverse $\sum_{\ell \geq 0} (p^{s-n} U (p^n X^{-1}))^\ell$. 
\end{rem}

Sometimes, one wants to consider not only the finite slope part, but to specify
the slope. For any rank one point $\Spa(C, \ocal_C) \rightarrow \mathcal{Z}$
corresponding to a character $\chi: T(\qq_p) \rightarrow C^\times$, we define
the slope of $\chi$ as follows. Let $v: C \rightarrow \mathbb{R}\cup \{ +
\infty\}$ be the valuation, normalized by $v(p)=1$. Composing the valuation and
$\chi$ we obtain a map $v(\chi): T(\qq_p) \rightarrow \mathbb{R}$. On the other
hand, via the  exact sequence~\eqref{eqn:slope-left-exact-tori} %
we can think of $v(\chi)$ as an element of $X^\star(T^d)_{\mathbb{R}}$ (see also \cite[Sect. 5.9]{boxer2021higher}). 
This defines a continuous  ``slope'' map $s: \mathcal{Z} \rightarrow
X^\star(T^d)_{\mathbb{R}}$ (which factors through the Berkovich space of~$\cZ$). %

\begin{rem} If (again abusively) we consider the case $T(\qq_p) =
  \ZZ$, then $\mathcal{Z} = \mathbb{G}_{m}^{an}$, and $s: \mathbb{G}_{m}^{an} \rightarrow \mathbb{R}$ extends the  $p$-adic valuation on classical rigid analytic points.
\end{rem} 

We have a partial  order relation on $X^\star(T^d)_{\mathbb{R}}$ where $\lambda \geq \lambda'$ if $\lambda(v(t))-\lambda'(v(t)) \geq 0$ for any $t \in T^+(\qq_p)$. %
For any $\lambda \in X^\star(T^d)_{\mathbb{\qq}}$, we can consider the cone $\{ \lambda' \in X^\star(T^d)_{\mathbb{R}}, \lambda' \leq \lambda\}$. 
The subset $s^{-1}( \{ \lambda' \in X^\star(T^d)_{\mathbb{R}}, \lambda' \leq
\lambda\})$ is the closure of a  unique rational open %
that we denote by  $\mathcal{Z}_{\leq \lambda} \stackrel{i_\lambda}\rightarrow \mathcal{Z}$. 
We have functors
 $f_{\leq \lambda}^\star:  D( \Mod_{T^+(\qq_p)}(\qq_p)) \rightarrow  D(\mathcal{Z}_{\leq \lambda})$, $f_{\leq \lambda}^\star M=  i_\lambda^\star f^\star M$. We also have a forgetful  functor $(f_{\leq \lambda})_\star: D(\mathcal{Z}_{\leq \lambda}) \rightarrow D( \Mod_{T(\qq_p)}(\qq_p))$. 
We finally define the  slope $\leq \lambda$ functor: 
$$(-)^{\leq \lambda} = (f_{\leq \lambda})_\star f_{\leq \lambda}^\star: D( \Mod_{T^+(\qq_p)}(\qq_p))  \rightarrow D( \Mod_{T(\qq_p)}(\qq_p)).$$
If $M \in D( \Mod_{T^+(\qq_p)}(\qq_p))$, then we have  $\lim_{\lambda} M^{\leq \lambda} = M^{\nfs}$.
We leave it to the reader to define the slope $\lambda$ functor $M \mapsto M^{=\lambda}$ and the slope $\geq \lambda$ functor, $M \mapsto M^{\geq \lambda}$. 

\begin{rem}\label{rem-finite-slope-Banach} Let us assume that $M$ is a Banach
  space equipped with a compact endomorphism $X$.  In other words, we are
  (again, abusively) in the situation  $T^+(\qq_p) = \ZZ_{\geq 0}$ and $T(\qq_p) = \ZZ$, and $1 \in  \ZZ_{\geq 0}$ acts like $X$. It follows from \cite[Prop. 12]{MR0144186}  that for any $\lambda \in \mathbb{\qq}$,  $M^{\leq \lambda}$ is a direct summand of $M$ and is finite dimensional. In particular, $f^\star M$ defines a coherent sheaf on $\mathcal{Z}$ and $M^{\nfs}$ is a pro-finite vector space.  
\end{rem}

\subsubsection{Hecke algebra action}

Let $K_U$ be a compact open subgroup of $U(\qq_p) \subseteq B(\qq_p)$ which
admits an Iwahori factorization. Let
$\lambda \in
X^\star(T)_{\qq_p}$. %
For any $M \in \Mod^{\lambda-\sm}_{B(\qq_p)}(E)$, we consider the
submodule %
$M^{K_U}$ of $K_{U}$-invariants. It is canonically a direct summand, since one
can define a trace $\mathrm{Tr}_{K_U}: M \rightarrow M^{K_U}$. Indeed, we have
$M = \colim_{K'} M^{K'}$ where $K'$ runs through the compact open subgroups of
$U(\qq_p)$.  Then for any $K_U' \subseteq K_U$ we can define a normalized
trace
\[\Tr_{K'_U/K_U}=\frac{1}{[K_U: K_U']}\left(\sum_{k \in K_U/K'_U} k\right): M^{K'_U}
  \rightarrow M^{K_U},\] and passing to the inductive limit
over~$\underline{K}_U$ yields the map $\mathrm{Tr}_{K_U}$.

Let $t \in T^+(\qq_p)$. We define an action of $t$ on $M^{K_U}$ as follows:
 $$M^{K_U} \stackrel{t}{\rightarrow} M^{ tK_U t^{-1}} \stackrel{\mathrm{Tr}_{tK_U t^{-1}/K_U}}\rightarrow M^{K_U}$$

 \begin{lem} The above rule defines an action of the commutative monoid $T^+(\qq_p)$ on  $M^{K_U}$. 
\end{lem}
\begin{proof} This follows from \cite[Lem. 4.1.5]{CasselmanNotes}.
\end{proof}

We deduce that we have an  exact  functor $$(-)^{K_U}:  D(\Mod^{\lambda-\sm}_{B(\qq_p)}(E))  \rightarrow D(\Mod^{\lambda-\sm}_{T^+(\qq_p)}(E)).$$

\begin{lem} Let $K'_U \subseteq K_U$. We have a natural transformation
  $(-)^{K'_U} \Rightarrow (-)^{K_U}$, induced by the trace. %
\end{lem}
\begin{proof} We consider the map $\mathrm{Tr}_{K'_U/K_U}: M^{K'_U} \rightarrow M^{K_U}$. 
It is elementary to check that this map commutes with the action of $t$. 
\end{proof}

We can now consider the composite functor: 
$$(-)^{K_U, \nfs}: D(\Mod^{\lambda-\sm}_{B(\qq_p)}(E))  \stackrel{(-)^{K_U}}\rightarrow D(\Mod^{\lambda-\sm}_{T^+(\qq_p)}(E)) \stackrel{(-)^{\nfs}}\rightarrow  D(\Mod^{\lambda-\sm}_{T(\qq_p)}(E))$$

\begin{lem} Let $K'_U \subseteq K_U$.  The natural transformation  $(-)^{K'_U, \nfs} \Rightarrow (-)^{K_U,\nfs}$ is an isomorphism.  
\end{lem}
\begin{proof} We pick $t \in T^{++}(\qq_p)$ such that $tK_Ut^{-1} \subseteq K'_U$. 
We have a commutative diagram (where the vertical maps are trace maps, and the
diagonal map is the composite $M^{K_U} \stackrel{t}{\rightarrow} M^{ tK_U t^{-1}} \stackrel{\mathrm{Tr}_{tK_U t^{-1}/K'_U}}\rightarrow M^{K_U}$): %
\begin{eqnarray*}
\xymatrix{ M^{K_U} \ar[r]^{t}\ar[rd] & M^{K_U} \\
M^{K'_U} \ar[u] \ar[r]^{t} & M^{K'_U} \ar[u]}
\end{eqnarray*}

On the finite slope quotient, the horizontal maps $t$ are isomorphisms, showing that the maps $M^{K'_U,\nfs} \rightarrow M^{K_U,\nfs}$ are isomorphisms.
\end{proof} 
In view of this lemma, we often simply write $M^{\nfs}$ instead of $M^{K_U,\nfs}$, as the choice of $K_U$ is irrelevant. 
\subsubsection{The finite slope part of the higher Coleman functor}\label{sec-fspHC}

We  define the finite slope part of the  higher Coleman functor: 

\begin{eqnarray*}
HC^{\nfs}_{w,\lambda}: \ocal_{\lambda-\alg} & \rightarrow & D(\Mod^{\lambda-\sm}_{T(\qq_p)}(E)) \\
M & \mapsto &\mathrm{R}\Gamma_w( \Sh^{\tor}_{K^p}, \VBzero(HCS_{w,\lambda}(M)))^{\nfs}%
\end{eqnarray*}

\begin{eqnarray*}
HC^{\nfs}_{\cusp, w,\lambda}: \ocal_{\lambda-\alg} & \rightarrow & D(\Mod^{\lambda-\sm}_{T(\qq_p)}(E)) \\
M & \mapsto &\mathrm{R}\Gamma_w( \Sh^{\tor}_{K^p}, \VBzero(HCS_{w,\lambda}(M)(-D)))^{\nfs}%
\end{eqnarray*}

\subsubsection{Comparison with higher Coleman theory \cite{boxer2021higher}}\label{subsubsec:comparison-to-BP21}
Let $\kappa \in X^\star(T)^{M,+}$ and $\chi: T(\ZZ_p)
\rightarrow \overline{\qq}_p^\times$ be a  finite order character. %
We have defined in  \cite[Sect.\ 5]{boxer2021higher} cohomology theories $\mathrm{R}\Gamma_{w}(K^p, \kappa,  \chi)^{+,\nfs}$ and $\mathrm{R}\Gamma_{w}(K^p, \kappa,  \chi, \cusp)^{+,\nfs}$.

\begin{thm}\label{thm-compare-Higher-Coleman-old-new} We have canonical isomorphisms of smooth $T(\qq_p)$-modules (where the decomposition on the right hand side corresponds to  the decomposition into isotypic parts for the action of $T(\ZZ_p)$):  
\begin{eqnarray*}
HC^{\nfs}_{\cusp,0, w}(L(\mf{m}_w)_{-w^{-1}w_{0,M}\kappa}) \otimes_{\qq_p} \Qpbar &=& \oplus_{\chi: T(\ZZ_p) \rightarrow \overline{\qq}_p^\times}  \mathrm{R}\Gamma_{w}(K^p, \kappa,  \chi, \cusp)^{+,\nfs} \\
HC^{\nfs}_{0, w}(L(\mf{m}_w)_{-w^{-1}w_{0,M}\kappa}) \otimes_{\qq_p} \Qpbar &=& \oplus_{\chi: T(\ZZ_p) \rightarrow \overline{\qq}_p^\times}  \mathrm{R}\Gamma_{w}(K^p, \kappa,  \chi)^{+,\nfs}
\end{eqnarray*}
\end{thm}

\begin{proof} By Lemma \ref{lem-VBandclassicalsheaf}, we have  $VB^0\circ
  HCS_{w,0}( L(\mf{m}_w)_{-w^{-1}w_{0,M}\kappa}) =
  \omega^{\kappa,sm}\vert_{\pi_{HT}^{-1}(C_w^\dag)}$, so that by definition we have $HC_{ w,0}(L(\mf{m}_w)_{-w^{-1}w_{0,M}\kappa}) = \mathrm{R}\Gamma_w( \Sh^{\tor}_{K^p}, \omega^{\kappa,\sm})$. 
 Let $Z_0$ be a closed neighborhood of $w$ in $C_w$ stable under  a compact open
 subgroup $K_U \subseteq U(\qq_p)$ (admitting an Iwahori factorization) and let $U_0$ be an open neighborhood of $Z_0$ in $C_w$. We let $F_0 = \pi_{HT}^{-1}(Z_0)$ and $V_0 =  \pi_{HT}^{-1}(U_0)$. We can define an action of $t \in T^+(\qq_p)$ on 
$\mathrm{R}\Gamma_{F_0}(V_0, \omega^{\kappa,\sm} )$ as follows (where the second
and third maps are respectively given by restriction and the trace): %

$$ \mathrm{R}\Gamma_{F_0}(V_0, \mathscr{F})^{K_U} \stackrel{t}\rightarrow \mathrm{R}\Gamma_{F_0t^{-1}}(V_0 t^{-1}, \omega^{\kappa,\sm} )^{tK_Ut^{-1}} \rightarrow $$ $$\mathrm{R}\Gamma_{F_0}(V_0, \mathscr{F})^{tK_Ut^{-1}} \rightarrow \mathrm{R}\Gamma_{F_0}(V_0, \omega^{\kappa,\sm})^{K_U}.$$

 The map $\mathrm{R}\Gamma_{F_0}(V_0, \omega^{\kappa,\sm}) \rightarrow \mathrm{R}\Gamma_w(\Sh^{\tor}_{K^p}, \omega^{\kappa,\sm})$ induces a 
$T^+(\qq_p)$-equivariant map $\mathrm{R}\Gamma_{F_0}(V_0,
\omega^{\kappa,\sm})^{K_U} \rightarrow
\mathrm{R}\Gamma_w(\Sh^{\tor}_{K^p}, \omega^{\kappa,\sm})^{K_U}$  which is
a quasi-isomorphism on the finite slope part. Indeed,  it is easy to see that we have factorizations: 
\begin{eqnarray*}
\xymatrix{ \mathrm{R}\Gamma_{F_0}(V_0,\omega^{\kappa,\sm})^{K_U} \ar[r]^{t}\ar[rd] & \mathrm{R}\Gamma_{F_0}(V_0, \omega^{\kappa,\sm})^{K_U}\\
\mathrm{R}\Gamma_{F_0t^{-1}}(V_0 t^{-1} , \omega^{\kappa,\sm})^{K_U} \ar[r]^{t}\ar[u]  & \mathrm{R}\Gamma_{F_0t^{-1}}(V_0t^{-1},  \omega^{\kappa,\sm})^{K_U} \ar[u]}
\end{eqnarray*}
On the finite slope part, the maps $t$ are quasi-isomorphisms, therefore the maps $\mathrm{R}\Gamma_{F_0t^{-1}}(V_0t^{-1}, \omega^{\kappa,\sm})^{K_U,\nfs}  \rightarrow \mathrm{R}\Gamma_{F_0}(V_0, \omega^{\kappa,\sm})^{K_U,\nfs} $ are also quasi-isomorphisms. Now we can take $t \in T^{++}(\qq_p)$. Since $\mathrm{R}\Gamma_w(\Sh^{\tor}_{K^p}, \omega^{\kappa,\sm})^{K_U} = \colim_n \mathrm{R}\Gamma_{F_0t^{n}}(V_0t^{-n}, \mathscr{F})^{K_U,\nfs} $,  the colimit is constant and we conclude.

We now consider compact open subgroups $K_{p,n}$ of $G(\qq_p)$ which
admits an Iwahori decomposition $K_{p,n} = K_U \times  K_{\bar{B},n}$
where $K_{\bar{B},n}$ is the principal level $p^n$ congruence
subgroup  in $\bar{B}(\qq_p)$. We let $K'_{p,n}$  be the compact open
subgroup of $G(\qq_p)$ which admit an Iwahori decomposition $K'_{p,n} = K_U \times T(\ZZ_p) \times  K_{\bar{U},n}$ where $K_{\bar{U},n}$ is the principal level $p^n$ congruence  subgroup  in $\bar{U}(\qq_p)$. Note that $K_{p,n} \subseteq K'_{p,n}$ is a normal subgroup and that $K'_{p,n}/K_{p,n} = T(\ZZ_p/p^n\ZZ_p)$. 

We have that  $\mathrm{R}\Gamma_{F_0}(V_0, \omega^{\kappa,\sm})^{K_U}  = \colim_{m, n} \mathrm{R}\Gamma_{F_{0,m, K_{p,n}}}(V_{0, m, K_{p,n}}, \omega^{\kappa}_{K_{p,n}})$.  Moreover, $$\mathrm{R}\Gamma_{F_{0,m, K_{p,n}}}(V_{0, m, K_{p,n}}, \omega^{\kappa}_{K_{p,n}})\otimes_{\qq_p} \Qpbar = \oplus_{\chi: T(\ZZ/p^n\ZZ_p) \rightarrow \Qpbar^\times} \mathrm{R}\Gamma_{F_{0,m, K'_{p,n}}}(V_{0, m, K'_{p,n}}, \omega^{\kappa}_{K'_{p,n}}(\chi)).$$

We claim that each $ \mathrm{R}\Gamma_{F_{0,m, K'_{p,n}}}(V_{0, m, K'_{p,n}}, \omega^{\kappa}_{K'_{p,n}}(\chi))$ can be equipped with an action of $T^+(\qq_p)$ and on the finite slope quotient the maps: 
$$\mathrm{R}\Gamma_{F_{0,m, K'_{p,n}}}(V_{0, m, K'_{p,n}}, \omega^{\kappa}_{K'_{p,n}}(\chi)) \rightarrow \mathrm{R}\Gamma_{F_{0,m, K'_{p,n'}}}(V_{0, m, K'_{p,n'}}, \omega^{\kappa}_{K'_{p,n}}(\chi))$$ for $n' \geq n$ are quasi-isomorphisms, and the maps: 
$$\mathrm{R}\Gamma_{F_{0,m, K'_{p,n}}}(V_{0, m, K'_{p,n}}, \omega^{\kappa}_{K'_{p,n}}(\chi)) \rightarrow \mathrm{R}\Gamma_{F_{0,m', K'_{p,n}}}(V_{0, m', K'_{p,n}}, \omega^{\kappa}_{K'_{p,n}}(\chi))$$ for $m' \geq m$ are quasi-isomorphism. 
This follows from the property that for $t \in T^{++}(\qq_p)$ sufficiently regular, we have factorizations:
 \begin{eqnarray*}
\xymatrix{ \mathrm{R}\Gamma_{F_{0,m, K'_{p,n'}}}(V_{0,m, K'_{p,n'}},\omega^{\kappa}_{K'_{p,n'}}(\chi)) \ar[r]^{t}\ar[rd] & \mathrm{R}\Gamma_{F_{0,m, K'_{p,n'}}}(V_{0,m, K'_{p,n}}, \omega^{\kappa}_{K'_{p,n'}}(\chi))\\
\mathrm{R}\Gamma_{F_{0,m, K'_{p,n}}}(V_{0,m, K'_{p,n}}  , \omega^{\kappa}_{K'_{p,n}}(\chi)) \ar[r]^{t}\ar[u]  & \mathrm{R}\Gamma_{F_{0,m, K'_{p,n}}}(V_{0,m, K'_{p,n}}, \omega^{\kappa}_{K'_{p,n}}(\chi)) \ar[u]}
\end{eqnarray*}
as well as factorizations:
 \begin{eqnarray*}
\xymatrix{ \mathrm{R}\Gamma_{F_{0,m', K'_{p,n}}}(V_{0,m', K'_{p,n'}},\omega^{\kappa}_{K'_{p,n}}(\chi)) \ar[r]^{t}\ar[rd] & \mathrm{R}\Gamma_{F_{0,m', K'_{p,n}}}(V_{0,m', K'_{p,n}}, \omega^{\kappa}_{K'_{p,n}}(\chi))\\
\mathrm{R}\Gamma_{F_{0,m, K'_{p,n}}}(V_{0,m, K'_{p,n}}  , \omega^{\kappa}_{K'_{p,n}}(\chi)) \ar[r]^{t}\ar[u]  & \mathrm{R}\Gamma_{F_{0,m, K'_{p,n}}}(V_{0,m, K'_{p,n}}, \omega^{\kappa}_{K'_{p,n}}(\chi)) \ar[u]}
\end{eqnarray*}
We conclude, since by definition (see just after Theorem 5.4.14
in~\cite{boxer2021higher}) we have $\mathrm{R}\Gamma_{F_{0,m', K'_{p,n}}}(V_{0, m', K'_{p,n}}, \omega^{\kappa}_{K'_{p,n}}(\chi))^{\nfs} = \mathrm{R}\Gamma_{w}(K^p, \kappa,  \chi)^{+,\nfs}$. 
\end{proof}

Let $\nu: T(\ZZ_p) \rightarrow \C_p^\times$ be an analytic character.    \cite[Sect. 5]{boxer2021higher} %
We have cohomology theories $\mathrm{R}\Gamma_{w,an}(K^p,
\nu)^{+,\nfs}$ and $\mathrm{R}\Gamma_{w,an}(K^p, \nu, \cusp)^{+,\nfs}$, defined in
\cite[Sect.\ 6.4]{boxer2021higher} (see just after Theorem
  6.4.10
there). Let $d\nu \in X^\star(T)_{\C_p}$ be the differential of $\nu$. Define $\kappa  \in X^\star(T)_{\C_p}$ by the formula 
$d\nu = -w^{-1}w_{0,M}(\kappa + \rho)- \rho$. 

\begin{thm}\label{thm:another-BP21-comparison} We have canonical isomorphisms of  $T(\qq_p)$-modules (where the decomposition on the right hand side corresponds to  the decomposition into isotypic parts for the action of $T(\ZZ_p)$):  
\begin{eqnarray*}
HC^{\nfs}_{\cusp,d\nu , w}(M(\mathfrak{m}_w)_{-w^{-1}w_{0,M}\kappa})  &=& \oplus_{\nu': T(\ZZ_p) \rightarrow \CC_p^\times,~d\nu' = d\nu}  \mathrm{R}\Gamma_{w}(K^p, \nu',\cusp)^{+,\nfs} \\
HC^{\nfs}_{d\nu , w}(M(\mathfrak{m}_w)_{-w^{-1}w_{0,M}\kappa})  &=& \oplus_{\nu': T(\ZZ_p) \rightarrow \CC_p^\times,~d\nu' = d\nu}  \mathrm{R}\Gamma_{w}(K^p, \nu')^{+,\nfs}
\end{eqnarray*}
\end{thm}
\begin{proof} This is similar to the proof of Theorem \ref{thm-compare-Higher-Coleman-old-new} and left to the reader.
\end{proof}

We have  the following theorem which slightly generalizes  \cite[Thm. 1.2.1, Thm. 6.7.3]{boxer2021higher}.
\begin{thm}\label{coro-vanishingThm} The functor $HC^{\nfs}_{w,\lambda}$ has
  cohomological amplitude $[\ell(w),d]$,  and $HC^{\nfs}_{\cusp, w,\lambda}$ has cohomological
  amplitude  $[0, \ell(w)]$. 
\end{thm}
\begin{proof}  Given Theorem \ref{thm-coho-amplitudeHC}, we simply need to see
  that the functor $(-)^{\nfs}$ is exact on higher Coleman theories. But one sees
  (see Theorem \ref{thm-compare-Higher-Coleman-old-new} and its proof)  that the
  finite slope part is obtained by taking  the finite slope part of  an
  inductive system of complexes of Banach spaces acted on by a compact operator,
  and $(-)^{\nfs}$ is exact in this case, see \cite[Proposition
  5.1.4]{boxer2021higher}. We remark that we could also deduce this theorem
  directly from \cite[Thm. 1.2.1, Thm. 6.7.3]{boxer2021higher} (which is the
  current theorem in the case of $HC^{\nfs}_{w,\lambda}$ and $HC^{\nfs}_{\cusp,
    w,\lambda}$ applied to Vermas) by using a diagram chase  similar to the
  proof of Theorem \ref{thm-STrict}.%
\end{proof}
\subsubsection{Bounds on slopes for higher Coleman theory}%
Using \cite{boxer2021higher} and \cite{boxer2023higher}, we can obtain bounds for the slopes. 

\begin{thm}\label{thm-bounds-on-slopes} Assume that either the Shimura variety is proper or that we are in the Siegel case. Let $w \in \WM$. Let $M \in \ocal(\mathfrak{m}_w, \mathfrak{b}_{M_w})_{\lambda-alg}$ be a module generated by  a highest weight vector of weight $\nu$. 
Then the slopes of $HC^{\nfs}_{\cusp,w, \lambda}(M)$ and $HC^{\nfs}_{w,\lambda}(M)$ are $\geq \lambda- \nu + w^{-1}w_{0,M} \rho + \rho$. 
\end{thm}%

\begin{rem} The finite slope projector and the notion of slope $\geq \lambda- \nu + w^{-1}w_{0,M} \rho + \rho$ were explained in section \ref{sect-finite-slope-proj}.
\end{rem}

\begin{rem} We conjecture (\cite[Conj.\ 6.8.1]{boxer2021higher}) that the  theorem should hold for any Shimura variety. The slightly weaker bound $\geq \lambda- \nu$ is currently available in full generality by \cite[Thm.\ 6.8.3]{boxer2021higher}. 
\end{rem}

\begin{proof}[Proof of Theorem~\ref{thm-bounds-on-slopes}] For $M = M(\mathfrak{m}_w)_{\nu}$ a Verma module with highest
  weight $\nu$, this follows from \cite[Thm.\ 6.10.1]{boxer2021higher}  for the proper case, and 
    \cite[Cor.\ 6.2.16]{boxer2023higher} in the Siegel case, %
  together with
  Theorem~\ref{thm:another-BP21-comparison} (the paper \cite{boxer2023higher}
  shows  that the strongly small slope condition which is sometimes  needed in
  \cite{boxer2021higher} can be weakened to the  small slope condition in the
  symplectic case). Since any 
   module $M$ as in the statement of the theorem admits a resolution  by Verma
   modules   $M(\mathfrak{m}_w)_{\nu'}$ with $\nu  \geq  \nu'$, we are done. %
\end{proof}
\subsection{$p$-adic Eichler Shimura theory}\label{sect-general-padicES}

 Let $\lambda \in
X^\star(T)_{E}$. %
We define functors:
\begin{eqnarray*}
CC_{\lambda}: \ocal(\mathfrak{g}, \mathfrak{b})_{\lambda-alg} &\rightarrow& D( \Mod^{\lambda-\sm}_{B(\qq_p)}(E)) \\
M & \mapsto & \mathrm{RHom}_{\mathfrak{g}}( M, \mathrm{R}\Gamma(\Sh_{K^p}, \qq_p)^{\la} )
\end{eqnarray*}
\begin{eqnarray*}
CC_{\cusp, \lambda}:  \ocal(\mathfrak{g}, \mathfrak{b})_{\lambda-alg} &\rightarrow& D( \Mod^{\lambda-\sm}_{B(\qq_p)}(E)) \\
M & \mapsto & \mathrm{RHom}_{\mathfrak{g}}( M, \mathrm{R}\Gamma_c(\Sh_{K^p}, \qq_p)^{\la} ).
\end{eqnarray*}
\begin{thm}\label{thm-p-adic-ES} %
For any $M \in  \ocal(\mathfrak{g}, \mathfrak{b})_{\lambda-alg}$, we have 
$$CC_\lambda(M) \otimes \C_p= \mathrm{R}\Gamma(\Sh^{\tor}_{K^p}, \VBred(\Loc(M)))$$
and
$$CC_{\cusp, \lambda}(M) \otimes \C_p = \mathrm{R}\Gamma(\Sh^{\tor}_{K^p}, \VBred(\Loc(M)) \otimes_{\mathscr{O}^{\sm}_{\Sh^{\tor}_{K^p}}} \mathscr{I}^{\sm}_{\Sh^{\tor}_{K^p}})$$ 
Moreover, the action of $\mu \in Z(\mathfrak{m})$  on $\Loc(M)$ via the
horizontal action induces an arithmetic Sen operator on the left hand side. %
\end{thm}
\begin{proof} By Theorem \ref{thm-prim-compar}, we have $\mathrm{R}\Gamma(\Sh_{K^p}^{\tor}, \qq_p)^{\la}{\otimes}_{\qq_p} \C_p  = \mathrm{R}\Gamma_{\an}(\Sh_{K^p}^{\tor}, \oscr^{\la}_{\Sh_{K^p}^{\tor}})$.
We consider the category %
$\Mod'_{\mathfrak{g}}(\oscr^{\sm}_{\Sh_{K^p}^{\tor}})$ of sheaves of
solid $\oscr^{\sm}_{\Sh_{K^p}^{\tor}}$-modules equipped with an action
of $\mathfrak{g}$, and the similarly-defined category
$\Mod'_{\mathfrak{g}}(\C_p)$.. These are abelian categories, with enough
injectives. %
 We consider the diagram of functors:
\begin{eqnarray*}
\xymatrix{\Mod'_{\mathfrak{g}}(\oscr^{\sm}_{\Sh_{K^p}^{\tor}}) \ar[r] \ar[d] & \Mod (\oscr^{\sm}_{\Sh_{K^p}^{\tor}}) \ar[d] \\
\Mod'_{\mathfrak{g}}(\C_p) \ar[r] & \Mod(\C_p)}
\end{eqnarray*}
where the horizontal arrows are given by taking $\mathfrak{g}$-invariants and the vertical arrows are given by taking global sections. 
This diagram is $2$-commutative,  and it induces a $2$-commutative diagram at the level of bounded derived categories. 
We deduce that  $CC_\lambda(M) \otimes \C_p =
\mathrm{R}\Gamma(\Sh^{\tor}_{K^p}, \mathrm{RHom}_{\mathfrak{g}}(M,
\oscr^{\la}_{\Sh_{K^p}^{\tor}}))$. We therefore need to show that $\VBred(\Loc(M)) = \mathrm{RHom}_{\mathfrak{g}}(M,
\oscr^{\la}_{\Sh_{K^p}^{\tor}})$. To this end, the  Chevalley--Eilenberg resolution
yields %
\[\mathrm{RHom}_{\mathfrak{g}}(M, \oscr^{\la}_{\Sh_{K^p}^{\tor}}) =
 [M^\vee \otimes \oscr^{\la}_{\Sh_{K^p}^{\tor}} \rightarrow  M^\vee \otimes \oscr^{\la}_{\Sh_{K^p}^{\tor}} \otimes \mathfrak{g}^\vee  \rightarrow \cdots \rightarrow M^\vee \otimes \oscr^{\la}_{\Sh_{K^p}^{\tor}} \otimes \Lambda^{r}\mathfrak{g}^\vee]\]
in degrees $0$ up to $ r = \mathrm{dim}\mathfrak{g}$. 
By   Theorem~\ref{thm-VB}(4) (and the flatness of $M^\vee \otimes\Lambda^{i}\mf{g}^{\vee}$), we have \[\VBzero(  M^\vee \otimes \mathcal{C}^{\la} \otimes
  \Lambda^i\mathfrak{g}^\vee)  = M^\vee \otimes \oscr^{\la}_{\Sh_{K^p}^{\tor}}
  \otimes \Lambda^i\mathfrak{g}^\vee,\]%
and we deduce (see Remark~\ref{rem:apply-Vbzero-termwise}) that 
$\VBred(\Loc(M)) = \mathrm{RHom}_{\mathfrak{g}}(M,
\oscr^{\la}_{\Sh_{K^p}^{\tor}})$. The cuspidal case is identical.
The final claim regarding the Sen operator follows from Theorem \ref{thm-existence-sen-oncc}.
\end{proof}

\begin{thm}\label{thm-p-adic-ES2} Assume that $\lambda$ is non-Liouville.
For any $M \in \ocal(\mf{g},\mf{b})_{\lambda-\alg}$,  we have that $\HH^i(\VBred(\Loc(M))) = \VBzero(\HH^i(\Loc(M))$.
 We have  a spectral sequence: %
$$E_1^{p,q} = \oplus_{w \in \WM, \ell(w)=p} \HH^{p+q}(HC_{w,\lambda}(E\otimes^L_{\mf{u}_{P_w}}M)) $$ converging to 
$\HH^{p+q}( CC_\lambda(M))\otimes \C_p$.   
Similarly, we have  a spectral sequence: 
$$E_1^{p,q} = \oplus_{w \in \WM, \ell(w)=p} \HH^{p+q}(HC_{\cusp,w,\lambda}(E\otimes^L_{\mf{u}_{P_w}}M)) $$ converging to 
$\HH^{p+q}( CC_{\cusp,\lambda}(M))\otimes \C_p$. 
In all cases $w\mu \in Z(\mathfrak{m}_w)$ acting on
$\HH_\star(\mathfrak{u}_{P_w},M)$ induces an arithmetic Sen operator.
\end{thm}
\begin{proof}To prove the first part, it suffices to show that the cohomology
  sheaves $\HH^i(\Loc(M))$ is acyclic for the functor $\VBred$. We can check
  this acyclicity locally, and in particular after restricting to Bruhat strata,
  where it follows from  Theorems~ \ref{thm-localization} and~\ref{thm-VB}(5).
  Finally, the spectral sequence is then a simple consequence of
  Proposition~\ref{prop:abstract-Cousin-spectral-sequence} together with Theorems~\ref{thm-p-adic-ES} and~\ref{thm-localization}.
\end{proof}

Let $\lambda \in X^\star(T)_E$ be a non-Liouville weight. We apply the above
theorem  to the Verma module  $ M(\mathfrak{g})_\lambda$, giving the following
results.

\begin{cor}\label{coro-simplification-anti} Assume that $\lambda$ is antidominant in the sense of
  Remark~\ref{rem-antidom}, and that the Shimura variety is proper. Then  $CC_\lambda(M(\mathfrak{g})_\lambda)$ is concentrated in the middle degree $d$ and moreover, it has a decreasing filtration $\mathrm{Fil}^i \mathrm{H}^d(CC_\lambda(M(\mathfrak{g})_\lambda))$ with 
\begin{itemize}
\item $ \mathrm{Fil}^{d+1} \mathrm{H}^d(CC_\lambda(M(\mathfrak{g})_\lambda)) = 0 $,
\item $  \mathrm{Fil}^0 \mathrm{H}^d(CC_\lambda(M(\mathfrak{g})_\lambda)) =  \mathrm{H}^d(CC_\lambda(M(\mathfrak{g})_\lambda))$,
\item  $\mathrm{Gr}^p \mathrm{H}^d(CC_\lambda(M(\mathfrak{g})_\lambda)) =
  \oplus_{ w \in \WM, \ell(w) = p}
  \HH^{p}(HC_{w,\lambda}(M(\mathfrak{m}_w)_{\lambda + w^{-1}w_{0,M}\rho +
    \rho}))$.%
\end{itemize}
\end{cor}
\begin{proof} This follows from Theorem \ref{thm-p-adic-ES2}, because the
  spectral sequence degenerates by a combination of Corollary  \ref{coro-ESflag}
  (noting Remark \ref{rem-antidom}) and  Theorem \ref{thm-coho-amplitudeHC}
  (noting that since the Shimura variety is proper, $HC_{\cusp, w,\lambda}=HC_{w,\lambda}$).
\end{proof}

\begin{rem} We see that in the antidominant case, the highest weights appearing in the $p$-adic Eichler--Shimura decomposition follow the exact   same pattern as the highest weights appearing in the classical  Eichler--Shimura theory.
\end{rem}

 We now consider the general case where~$\lambda$ need not be antidominant, and we take the ``ordinary'' part. 

\begin{thm}\label{thm-p-adic-ES3} Assume that we are in the Siegel case or that
  the Shimura variety is proper. Let $\lambda \in X^\star(T)_E$ be a non-Liouville weight. We have that
  $CC^{\nfs}_{\lambda}(M(\mathfrak{g})_\lambda)$ and $CC^{\nfs}_{\cusp,
    \lambda}(M(\mathfrak{g})_\lambda)$ have slope $\geq 0$. Moreover, we have a
  spectral sequence: %
$$E_1^{p,q} = \oplus_{w \in \WM, \ell(w)=p} \HH^{2p+q-d}(HC^{=0}_{w,\lambda}(M(\mathfrak{m}_w)_{\lambda + w^{-1}w_{0,M}\rho + \rho})) $$ converging to 
$\HH^{p+q}( CC^{=0}_\lambda(M(\mathfrak{g})_\lambda))\otimes \C_p$,
and similarly for cuspidal cohomology. 
\end{thm}
\begin{proof} This is a combination of Theorem \ref{thm-p-adic-ES2}, Corollary  \ref{coro-ESflag} %
  and Theorem \ref{thm-bounds-on-slopes}.
\end{proof}
Finally we have the following corollary.
\begin{cor}\label{coro-p-adic-ES-ord} Assume that the Shimura variety is
  proper. Let $\lambda \in X^\star(T)_E$ be a non-Liouville weight. Then  $CC^{=0}_\lambda(M(\mathfrak{g})_\lambda))$ is concentrated in the middle degree~ $d$ and moreover, it has a decreasing filtration $\mathrm{Fil}^i \mathrm{H}^d(CC^{=0}_\lambda(M(\mathfrak{g})_\lambda))$ with 
\begin{itemize}
\item $ \mathrm{Fil}^{d+1} \mathrm{H}^d(CC^{=0}_\lambda(M(\mathfrak{g})_\lambda)) = 0 $,
\item $  \mathrm{Fil}^0 \mathrm{H}^d(CC^{=0}_\lambda(M(\mathfrak{g})_\lambda)) =  \mathrm{H}^d(CC^{=0}_\lambda(M(\mathfrak{g})_\lambda))$,
\item  $\mathrm{Gr}^p \mathrm{H}^d(CC^{=0}_\lambda(M(\mathfrak{g})_\lambda)) = \oplus_{w \in \WM, \ell(w) = p}   \HH^{p}(HC^{=0}_{w,\lambda}(M(\mathfrak{m}_w)_{\lambda + w^{-1}w_{0,M}\rho + \rho}))$.
\end{itemize}
\end{cor}
\begin{proof} This is a consequence of Theorem \ref{thm-p-adic-ES3} and Theorem \ref{thm:another-BP21-comparison}. 
\end{proof}

\begin{rem} Thus, we see that on the ordinary part the highest weights appearing in the $p$-adic Eichler--Shimura decomposition follow  the same pattern as the highest weights appearing in the classical  Eichler--Shimura theory. 
\end{rem}

\subsection{The classical Hodge--Tate decomposition for
  $\mathrm{GSp}_4$} \label{sub:classical} %
We now specialize the theory to the group $\mathrm{GSp}_4$ and first review the
classical Hodge--Tate decomposition.  From now on we use the notation for
$\GSp_4 $ introduced in Section~\ref{notn:GSp4}.
Let $\kappa = (k_1, k_2; w) \in X^\star(T)^+$ be a  dominant weight for $\mathrm{GSp}_4$. The dominance condition is $0 \geq k_1 \geq k_2$. We let $V_\kappa$ be the corresponding highest weight representation.  We let  $\mathcal{V}_{\kappa,K_p}^\vee$ be the pro-Kummer \'etale local system on $\Sh^{\tor}_{K_pK^p}$, attached to $V_\kappa^\vee$. 

The coherent weights appearing in the Hodge--Tate decomposition of the local system $\mathcal{V}_{\kappa,K_p}^\vee$ are  the 
$\{ -w_{0,M}(w.\kappa + 2 \rho^M)\}$ where $w \in {\WM}$ (see \cite{MR1083353}, Thm.\ 6.2). 
We make this explicit: %
 the set $ \{ -w_{0,M}(w.\kappa + 2 \rho^M)\}$ consists exactly of $(3-k_2,3-k_1;-w), (3-k_2, k_1+1;-w), (2-k_1,k_2;-w),  (k_1,k_2;-w)$. 
We recall that our conventions are that the cyclotomic  character $\qq_p(1)$ has
Hodge--Tate weight $-1$ and that the Sen operator acts via $1$ on the Sen module
of $\qq_p(1)$, so that the (generalized)  Hodge--Tate weights are  the negatives
of the eigenvalues of the Sen operator. 
Given our choice of parabolic $P_\mu$, we have  $\mu = (-1/2,  -1/2; 1/2) \in
X_\star(T)_{E}$. %
By Theorem \ref{thm-VB},  $\mu$ is an arithmetic Sen operator.

\begin{rem}[reality check] This is consistent with the fact that the tautological exact sequence over $\mathcal{FL}$ is 
$$ 0 \rightarrow \mathcal{L}_{(0,  -1; 1)} \rightarrow St \otimes \oscr_{\mathcal{FL}} \rightarrow \mathcal{L}_{(1,  0; 1)} \rightarrow 0$$ which pulls back to \[0 \rightarrow Lie(A)_{K_p}(1) \otimes_{ \oscr_{\Sh_{K_pK^p}^{\tor}}} \oscr_{\Sh_{K^p}^{\tor}} \rightarrow T_pA \otimes_{\ZZ_p} \oscr_{\Sh_{K^p}^{\tor}} \rightarrow (\omega_{A^t})_{K_p} \otimes_{ \oscr_{\Sh_{K_pK^p}^{\tor}}} \oscr_{\Sh_{K^p}^{\tor}} \rightarrow 0.\] We see that $Lie(A)$ has Sen weight $$1 = \langle (0,  -1; 1), (-1/2,  -1/2; 1/2) \rangle  $$ (and Hodge--Tate weight $-1$), while $w_{A^t}$ has  weight $0 = \langle (1, 0; 1), (-1/2,  -1/2; 1/2) \rangle$. 
\end{rem}
The Hodge--Tate weight attached to the sheaf $\omega^{(l_1,l_2;w)}$ is  $\frac{l_1+l_2-w}{2}$. Thus, in the Hodge--Tate decomposition of the cohomology of $V_\kappa^\vee$, the Hodge--Tate weights are given by the formula:
$(k_1,k_2;-w) \mapsto \frac{k_1+k_2+w}{2}$, $(2-k_1,k_2;-w) \mapsto \frac{2-k_1+k_2+w}{2}$, $(3-k_2, k_1+1;-w) \mapsto \frac{4-k_2+ k_1+w}{2}$, $ (3-k_2,3-k_1;-w) \mapsto \frac{6-k_1-k_2+w}{2}$. 

\begin{thm}[\cite{MR1083353}, Thm.\ 6.2]\label{thm-FC} We have  the following $G_{\qq_p} \times \mathbb{T}_{K_pK^p}$-equivariant isomorphisms: 
$$\mathrm{H}^i(\Sh^{\tor}_{K_pK^p}, \mathcal{V}^\vee_{\kappa,K_p}) \otimes_{\qq_p} \C_p= $$ $$\HH^{i}(\Sh^{\tor}_{K_pK^p}, \omega_{K_p}^{(k_1,k_2;-w)})(\frac{-k_1-k_2-w}{2}) \oplus 
\HH^{i-1}(\Sh^{\tor}_{K_pK^p}, \omega_{K_p}^{(2-k_1,k_2;-w)})(\frac{-2+k_1-k_2-w}{2}) $$ $$\oplus \HH^{i-2}(\Sh^{\tor}_{K_pK^p}, \omega_{K_p}^{(3-k_2,k_1+1;-w)})(\frac{-4+k_2-k_1-w}{2})
\oplus \HH^{i-3}(\Sh^{\tor}_{K_pK^p}, \omega_{K_p}^{(3-k_2, 3-k_1;-w)})(\frac{-6+k_2+k_1-w}{2})$$

$$\mathrm{H}^i_c(\Sh^{\tor}_{K_pK^p}, \mathcal{V}^\vee_{\kappa,K_p}) \otimes_{\qq_p} \C_p= $$ $$\HH^{i}(\Sh^{\tor}_{K_pK^p}, \omega_{K_p}^{(k_1,k_2;-w)}(-D_{K_p}))(\frac{-k_1-k_2-w}{2}) \oplus $$ $$
\HH^{i-1}(\Sh^{\tor}_{K_pK^p}, \omega_{K_p}^{(2-k_1,k_2;-w)}(-D_{K_p}))(\frac{-2+k_1-k_2-w}{2}) \oplus $$ $$ \HH^{i-2}(\Sh^{\tor}_{K_pK^p}, \omega_{K_p}^{(3-k_2,k_1+1;-w)}(-D_{K_p}))(\frac{-4+k_2-k_1-w}{2}) \oplus $$ $$
 \HH^{i-3}(\Sh^{\tor}_{K_pK^p}, \omega_{K_p}^{(3-k_2, 3-k_1;-w)}(-D_{K_p}))(\frac{-6+k_2+k_1-w}{2})$$

\end{thm}

We can also  state a similar result using completed cohomology. %
We first recall the following theorem:

\begin{thm} Let $V_\kappa$ be a finite dimensional representation of $G$ of highest weight $\kappa$. Then  
\begin{eqnarray*}
CC_{0}(V_\kappa) = \mathrm{RHom}_{\mathfrak{g}}( V_\kappa, \mathrm{R}\Gamma( \Sh^{\tor}_{K^p}, \qq_p)^{\la}) &=& \colim_{K_p} \mathrm{R}\Gamma( \Sh^{\tor}_{K_pK^p}, \mathcal{V}^\vee_{\kappa,K_p}) \\
CC_{\cusp, 0}(V_\kappa) = \mathrm{RHom}_{\mathfrak{g}}( V_\kappa, \mathrm{R}\Gamma_c( \Sh^{\tor}_{K^p}, \qq_p)^{\la}) &=& \colim_{K_p} \mathrm{R}\Gamma_c( \Sh^{\tor}_{K_pK^p}, \mathcal{V}^\vee_{\kappa,K_p})
\end{eqnarray*}
\end{thm}
\begin{proof}See e.g.\ \cite[Cor.\
  2.2.18]{MR2207783}. %
\end{proof} 
We deduce the following: 
\begin{cor} We have that  $$\mathrm{H}^i(\mathrm{RHom}_{\mathfrak{g}}( V_\kappa, \mathrm{R}\Gamma( \Sh^{\tor}_{K^p}, \qq_p)^{\la})) {\otimes} \C_p= $$ $$  \HH^i(\Sh^{\tor}_{K^p}, \omega^{(k_1,k_2;-w),\sm})(\frac{-k_1-k_2-w}{2}) \oplus 
\HH^{i-1}(\Sh^{\tor}_{K^p}, \omega^{(2-k_1,k_2;-w),\sm})(\frac{-2+k_1-k_2-w}{2}) $$ $$\oplus \HH^{i-2}(\Sh^{\tor}_{K^p}, \omega^{(3-k_2,k_1+1;-w),\sm})(\frac{-4+k_2-k_1-w}{2})
\oplus \HH^{i-3}(\Sh^{\tor}_{K^p}, \omega^{(3-k_2, 3-k_1;-w),\sm})(\frac{-6+k_2+k_1-w}{2})$$
There is a similar statement for compactly supported cohomology and the cuspidal coherent cohomology. 
\end{cor}
\begin{proof} This is simply obtained by passing to the limit over $K_p$ in
  Theorem \ref{thm-FC}. But alternatively, this is a combination of Theorem
  \ref{thm-p-adic-ES}, Lemma~\ref{lem-VBandclassicalsheaf} and Proposition \ref{prop-loc-fd}.
\end{proof}

\subsection{The $p$-adic Eichler--Shimura theory for $\mathrm{GSp}_4$}

We continue to assume that $G=\GSp_4 $, and let $\lambda \in
X^\star(T)_{E}$. We are ultimately interested in the case that
$\lambda = (1,1;w)$. We will now specialize the results of Section
\ref{sect-general-padicES}, in order to compute the Hodge--Tate structure of the
ordinary part of
$CC_\lambda( M(\mathfrak{g})_\lambda) = \mathrm{RHom}_{\mathfrak{b}}(\lambda,
\mathrm{R}\Gamma( \Sh^{\tor}_{K^p}, \qq_p)^{\la})$.  We will use some more
suggestive notation for the higher Coleman sheaves now that we specialize to
$\mathrm{GSp}_4$. Let us briefly summarize who are the main players.
 \begin{itemize}
 \item For all $\kappa \in X^{\star}(T)^{M,+}$, we have the classical  modular sheaves $\omega^{\kappa,\sm}$, computing classical cohomology. 
 \item  Following Definition  \ref{defn-HigherColeman}, we have the ``big'' sheaves $\omega^{\dag,\kappa}_w$ on $\pi_{HT}^{-1}(C_w^\dag)$ for all $\kappa \in X^\star(T)_{\C_p}$ and all $w \in \WM$. If $\kappa \in X^{\star}(T)^{M,+}$  we have maps $\omega^{\kappa, \sm}\vert_{\pi_{HT}^{-1}(C_w^\dag)} \rightarrow \omega^{\dag,\kappa}_w \otimes E(-w^{-1} w_{0,M} \kappa)$ by Proposition \ref{prop:integral-weight-simple-sheaf-computed} (where the twist is a twist of the $B(\qq_p)$-action). 
 \item We have higher Coleman theories $\mathrm{R}\Gamma_w( \Sh^{\tor}_{K^p}, \omega^{\dag,\kappa}_w)$ for the big sheaves and $\mathrm{R}\Gamma_w( \Sh^{\tor}_{K^p}, \omega^{\kappa,sm})$ for   the classical sheaves.
 \end{itemize}

The superscript  $(-)^{\ord}$    means the ordinary part, which is the minimal
slope part; we caution the reader that precisely what the ``minimal slope'' part
occasionally depends on the context, but will always be spelled out. On $CC_\lambda( M(\mathfrak{g})_\lambda)$, the ordinary part is the slope $=0$ part by Theorem \ref{thm-p-adic-ES3}. 
 Theorem \ref{thm-p-adic-ES3}  specializes  as follows. %

\begin{thm}\label{thm:ss-ord-part-gsp4} There is a spectral sequence:
$$ E_1^{p,q}:  \big(\HH^{2p+q-d}_{^pw}( \Sh^{\tor}_{K^p}, \omega_{^pw}^{\dag, -w_{0,M}({^pw}\cdot\lambda + 2 \rho^M)})\otimes \C_p(-{^pw}^{-1}w_{0,M} \rho-\rho)\big)^{\ord} \Rightarrow $$ $$ \HH^{p+q}(\mathrm{RHom}_{\mathfrak{b}}( \lambda, \mathrm{R}\Gamma(\Sh^{\tor}_{K^p}, {\qq}_p)^{\la})^{\ord} \otimes \C_p)$$
and similarly: 
$$ E_1^{p,q}:  \big(\HH^{2p+q-d}_{^pw}( \Sh^{\tor}_{K^p}, \omega_{^pw}^{\dag, -w_{0,M}({^pw}\cdot\lambda + 2 \rho^M)}(-D))\otimes \C_p(-{^pw}^{-1}w_{0,M} \rho-\rho)\big)^{\ord} \Rightarrow $$ $$ \HH^{p+q}(\mathrm{RHom}_{\mathfrak{b}}( \lambda, \mathrm{R}\Gamma_c(\Sh^{\tor}_{K^p}, {\qq}_p)^{\la})^{\ord} \otimes \C_p).$$
\end{thm}

If the Shimura variety were proper, we could use Corollary \ref{coro-p-adic-ES-ord} to simplify the spectral sequence. In our case we will arrive to a similar conclusion after making a non-Eisenstein localization. We can give a first   analysis of the spectral sequence with the help of some vanishing theorems. The following lemma comes as a complement to Theorem \ref{coro-vanishingThm}.

\begin{lem}\label{lem-complement-vanishingTHM}\leavevmode \begin{enumerate}
\item For all $\kappa \in X^\star(T)_{E}$, we have that \[\HH^0_{w}( \Sh^{\tor}_{K^p}, \omega_w^{ \dag,
  \kappa}(-D))^{\nfs} = \HH^0_{w}( \Sh^{\tor}_{K^p}, \omega_w^{ \dag,
  \kappa})^{\nfs}=0\] for all $w \neq \Id$ and  \[\HH^3_{w}( \Sh^{\tor}_{K^p}, \omega_w^{ \dag, \kappa}(-D))^{\nfs} = \HH^3_{w}( \Sh^{\tor}_{K^p}, \omega_w^{ \dag, \kappa})^{\nfs}=0\] for all $w \ne w_0^M$.
  \item For all $\kappa \in X^\star(T)^{M,+}$, we have that \[\HH^0_{w}( \Sh^{\tor}_{K^p}, \omega^{
  \kappa,\sm}(-D))^{\nfs} = \HH^0_{w}( \Sh^{\tor}_{K^p}, \omega^{
  \kappa,\sm})^{\nfs}=0\] for all $w \neq \Id$ and  \[\HH^3_{w}( \Sh^{\tor}_{K^p}, \omega^{  \kappa, \sm}(-D))^{\nfs} = \HH^3_{w}( \Sh^{\tor}_{K^p}, \omega^{ \kappa,\sm})^{\nfs}=0\] for all $w \ne w_0^M$.
  \end{enumerate}
  \end{lem}
  \begin{proof} The  statements regarding $\HH^0$ and $\HH^3$ are equivalent under Serre duality (\cite[Thm. 6.7.2]{boxer2021higher}). The vanishing of $\HH^0_{w}( \Sh^{\tor}_{K^p}, \omega_w^{ \dag,
  \kappa})^{\nfs}$ for $w \neq \Id$ follows from Theorem \ref{coro-vanishingThm}.  The injective map of sheaves $\omega^{\dag, \kappa}_w(-D) \rightarrow \omega^{\dag, \kappa}_w$ induces an injective map $\HH^0_{w}( \Sh^{\tor}_{K^p}, \omega_w^{ \dag,
  \kappa}(-D)) \rightarrow \HH^0_{w}( \Sh^{\tor}_{K^p}, \omega_w^{ \dag,
  \kappa})$ which  implies the vanishing of $\HH^0_{w}( \Sh^{\tor}_{K^p}, \omega_w^{ \dag,
  \kappa}(-D))^{\nfs}$ for $w \neq \Id$. One argues similarly for the sheaf~ $\omega^{\kappa,\sm}$. 
  \end{proof}

\begin{prop}\label{prop-pre-HT}\leavevmode
  \begin{enumerate}\item
    $\mathrm{RHom}_{\mathfrak{b}}( \lambda, \mathrm{R}\Gamma_c(\Sh^{\tor}_{K^p},
    \Qp)^{\la})^{\ord}$ is supported in degrees in the range $[1,3]$. Moreover
    (with obvious notation) the graded pieces for the Hodge--Tate decompositions
    are:

    \begin{enumerate}%
    \item$ \HH^3 \otimes \C_p: \mathrm{coker} (H^1_{^2w} (-D)\rightarrow
      H^3_{^3w}(-D)), H^2_{^2w}(-D), H^1_{^1w}(-D), H^0_{^0w}(-D)$.
    \item
      $\HH^2 \otimes \C_p: H^2_{^3w}(-D), \mathrm{Ker}(H^1_{^2w}(-D) \rightarrow
      H^3_{^3w}(-D))$.
    \item $\HH^1 \otimes \C_p: H^1_{^3w}(-D)$.
    \end{enumerate}
  \item     $\mathrm{RHom}_{\mathfrak{b}}( \lambda, \mathrm{R}\Gamma(\Sh^{\tor}_{K^p},
    \Qp)^{\la})^{\ord}$ is supported in degrees in the range $[3,5]$. Moreover
    the graded for the Hodge--Tate decompositions are:

    \begin{enumerate}
    \item$ \HH^3 \otimes \C_p: H^3_{^3w}, H^2_{^2w}, H^1_{^1w},
      \mathrm{Ker}(H^0_{^0w} \rightarrow H^2_{^1w})$.
    \item
      $\HH^4 \otimes \C_p: \mathrm{Coker}(H^0_{^0w} \rightarrow H^2_{^1w}),
      H^1_{^0w}$.
    \item $\HH^5 \otimes \C_p: H^2_{^0w}$.
    \end{enumerate}
  \end{enumerate}
\end{prop} %
\begin{proof} From Theorem \ref{coro-vanishingThm} and Lemma \ref{lem-complement-vanishingTHM}, we have that    the cohomology $\mathrm{R}\Gamma_{w}( \Sh^{\tor}_{K^p}, \omega_w^{ \dag, \kappa})^{\nfs}$ is supported in the range:
\begin{itemize}
\item $[0,2]$ for $w= ~^0w$,
\item $[1,2]$ for $w =~ ^1w$,
\item $[2]$ for $w=~^2w$,
\item $[3]$ for $w=~^3w$. 
\end{itemize}
On the other hand, the cohomology $\mathrm{R}\Gamma_{w}( \Sh^{\tor}_{K^p}, \omega_w^{ \dag, \kappa}(-D))^{\nfs}$ is supported in the range:
\begin{itemize}
\item $[0]$ for $w= ~^0w$,
\item $[1]$ for $w =~ ^1w$,
\item $[1,2]$ for $w=~^2w$,
\item $[1,3]$ for $w=~^3w$. 
\end{itemize}
It follows that the spectral sequences of Theorem~\ref{thm:ss-ord-part-gsp4}
degenerate on the second page.
\end{proof}

We now fix an irreducible residual representation $\bar{\rho}: G_{\mathbb{Q}}
\rightarrow \mathrm{GSp}_{4}(\Fbar_p)$. We define a maximal ideal $\mathfrak{m}_{\bar{\rho}}$
 of the abstract spherical Hecke algebra of level prime to $S\cup\{p\}$ (where $S$ is
 the set of primes at which~$\rhobar$ is ramified or~$K^p$ is not hyperspecial) %
 by the formula:
  $$P_\ell(X) \mod \mathfrak{m}_{\bar{\rho}} = \mathrm{det}(X- \bar{\rho}(\mathrm{Frob}_{\ell})),~\ell \notin S\cup\{p\}$$
  where $P_\ell(X)$ is the Hecke polynomial defined in~\eqref{eq:P-ell}. %

\begin{thm}\label{thm-van-bound} The map $\mathrm{R}\Gamma_c(\Sh^{\tor}_{K^p}, \Qpbar)_{\mathfrak{m}_{\bar{\rho}}} \rightarrow \mathrm{R}\Gamma(\Sh^{\tor}_{K^p}, \Qpbar)_{\mathfrak{m}_{\bar{\rho}}}$ is a quasi-isomorphism.
\end{thm}%
\begin{proof} Indeed the cohomology of the boundary is Eisenstein by identical
  arguments to those of~\cite[\S 4]{new-tho}. 
\end{proof}%

\begin{cor}\label{cor:quasi-iso-of-usual-to-cuspidal}
 The maps $\mathrm{R}\Gamma_{w}(\Sh^{\tor}_{K^p}, \omega_{w}^{\dag, \kappa}(-D))^{\ord}_{\mathfrak{m}_{\bar{\rho}}} \rightarrow \mathrm{R}\Gamma_{w}(\Sh^{\tor}_{K^p}, \omega_{w}^{\dag, \kappa})^{\ord}_{\mathfrak{m}_{\bar{\rho}}}$ are quasi-isomorphisms. %
 \end{cor}
  \begin{proof}While this could be proved by analyzing the cohomology
    of the boundary, we argue as follows. By Proposition
    \ref{prop-pre-HT} and Theorem \ref{thm-van-bound}, it suffices to show that
    the maps  $$ \HH^0_{^0w}(\Sh^{\tor}_{K^p}, \omega_{^0w}^{\dag,
      -w_{0,M}({^0w}\cdot\lambda + 2 \rho^M)})^{\ord} \rightarrow
    \HH^2_{^1w}(\Sh^{\tor}_{K^p}, \omega_{^1w}^{\dag, -w_{0,M}({^1w}\cdot\lambda
      + 2 \rho^M)})^{\ord},$$   $$ \HH^1_{^2w}(\Sh^{\tor}_{K^p},
    \omega_{^2w}^{\dag, -w_{0,M}({^2w}\cdot\lambda + 2 \rho^M)}(-D))^{\ord}
    \rightarrow  \HH^3_{^3w}(\Sh^{\tor}_{K^p}, \omega_{^3w}^{\dag,
      -w_{0,M}({^3w}\cdot\lambda + 2 \rho^M)}(-D))^{\ord}$$  are~$0$ after
 localizing at $\mathfrak{m}_{\bar{\rho}}$. The second statement follows
    from the first by duality, and the first statement follows from the fact
    that the natural map \[\HH^0_{^0w}(\Sh^{\tor}_{K^p}, \omega_{^0w}^{\dag,
      \kappa}(-D))^{\ord}_{\mathfrak{m}_{\bar{\rho}}} \rightarrow  \HH^0_{^0w}(\Sh^{\tor}_{K^p}, \omega_{^0w}^{\dag,
      \kappa})^{\ord}_{\mathfrak{m}_{\bar{\rho}}} \]is an isomorphism, %
  which is Lemma~ \ref{lem-boundaryannoyance} below.
  \end{proof}
  
  \begin{lem}\label{lem-boundaryannoyance}
  The map $\HH^0_{^0w}(\Sh^{\tor}_{K^p}, \omega_{^0w}^{\dag,
      \kappa}(-D))^{\ord}_{\mathfrak{m}_{\bar{\rho}}} \rightarrow  \HH^0_{^0w}(\Sh^{\tor}_{K^p}, \omega_{^0w}^{\dag,
      \kappa})^{\ord}_{\mathfrak{m}_{\bar{\rho}}} $ is an isomorphism. 
  \end{lem}
  
  \begin{proof}
  This is similar to \cite[Cor.\ 15.2.3.1]{pilloniHidacomplexes}  and \cite[Lem.\ 3.10.7]{BCGP}, except that we work with ordinary $p$-adic modular forms rather than classical forms.  We translate the statement to a result about Hida complexes which can then be proved as in these previous results but working mod $p$ and with the structure sheaf.
  
  The Hida complexes we consider are constructed in \cite{boxer2023higher} and
  also recalled below in Section~ \ref{higherhidaRT}; they are perfect complexes
  $M^\bullet_{^0w,\cusp}$ and $  M^\bullet_{^0w}$
  of $\Lambda=\Z_p\llb T(\Z_p)\rrb$-modules, and there is a natural morphism $M^\bullet_{^0w,\cusp}\to
  M^\bullet_{^0w}$. The complex  $M^\bullet_{^0w,\cusp}$ is projective in degree
  0 (and in fact it is the classical object constructed by Hida
  \cite{MR1954939}), while $M^\bullet_{^0w}$ has amplitude $[0,2]$.  We can also
  consider the boundary Hida complex
  $M^\bullet_{^0w,\partial}=\mathrm{cone}(M^\bullet_{^0w,\cusp}\to
  M^\bullet_{^0w})$.  A priori this has amplitude $[-1,2]$ but we will recall
  below the simple geometric reason that it has amplitude $[0,2]$. 
  
  Below we shall prove that after non Eisenstein localization
  $M^\bullet_{^0w,\partial,\mathfrak{m}_{\rhobar}}$ has amplitude $[1,2]$ (in
  fact it actually vanishes, but since we do not need this, we do not prove it).
  This statement implies
  that for any continuous homomorphism $\nu:T(\Q_p)\to\Qpbar^\times$ the morphism
\numequation\label{eq-boundarynak}\HH^0(M^\bullet_{^0w,\cusp,\mathfrak{m}_{\rhobar}}\otimes^L_{\Lambda,\nu}\Qpbar)\to \HH^0(M^\bullet_{^0w,\mathfrak{m}_{\rhobar}}\otimes^L_{\Lambda,\nu}\Qpbar)\end{equation}
  is an isomorphism, and by the comparison between higher Hida and Coleman theory \cite[Thm 6.2.15]{boxer2023higher} and Theorem \ref{thm-compare-Higher-Coleman-old-new} this is exactly the statement of the lemma.
  
  To prove the claim, by Nakayama's lemma it suffices to prove that
  $$\HH^0(M^\bullet_{^0w,\cusp,\mathfrak{m}_{\rhobar}}\otimes^L_{\Lambda}\FF_p[T(\FF_p)])\to\HH^0(M^\bullet_{^0w,\mathfrak{m}_{\rhobar}}\otimes^L_{\Lambda}\FF_p[T(\FF_p)])$$
  is an isomorphism.  We now translate this back into a statement about mod $p$ modular forms on the ordinary locus which we prove by analyzing the boundary.
  
  We consider $IG/\FF_p$, the special fiber of the (ordinary) Igusa variety
  corresponding to the pro-$p$ Iwahori subgroup of $P'(\Qp)$, in the notation of
  \cite[\S 3.4.5]{boxer2023higher}.  We let $\pi:IG^{\tor}\to IG^{\star}$ be its
  (partial) toroidal and minimal compactifications.  By the very construction of
  the Hida complexes, the map \eqref{eq-boundarynak} is nothing but the natural map
  \numequation\label{eq-modpordinary}\HH^0(IG^{\tor},\oscr(-D))^{\ord}_{\mathfrak{m}_{\rhobar}}\to\HH^0(IG^{\tor},\oscr)^{\ord}_{\mathfrak{m}_{\rhobar}}\end{equation}
  (For this and the meaning of the ordinary part, see \cite[\S
  5.2]{boxer2023higher}, but note that the setup is substantially simplified
  because $w={}^0w$.)  We remark that even without the  non-Eisenstein localization this map is
  always injective, which justifies the assertion made above that the boundary cohomology always has amplitude $[0,2]$.
  
  We now follow the strategy of \cite[Cor.\ 15.2.3.1]{pilloniHidacomplexes}  and
  \cite[Lem.\ 3.10.7]{BCGP}  to prove that \eqref{eq-modpordinary} is surjective
  after non-Eisenstein localization.  If not then there is a non-Eisenstein
  Hecke eigenvector occurring in $\HH^0(IG^{\tor},\O_D)$.  We write
  $D^\star\subseteq IG^*$ for the (reduced) boundary.  We have
  $\pi_\star\O_D=\O_{D^\star}$ so that
  $\HH^0(D,\O_D)=\HH^0(D^\star,\O_{D^\star})$. The boundary  $D^\star$ is a union of (ordinary) Igusa curves crossing at cusps.  We write $\tilde{D}^\star$ for the normalization, which is a disjoint union of ordinary Igusa curves.  We have an injective pullback map
  $$\HH^0(D^\star,\O_{D^\star})\to\HH^0(\tilde{D}^\star,\O_{\tilde{D}^\star})$$
  There is a compatibility between the $\GSp_4$ Hecke action and the $\GL_2$
  Hecke action at primes away from $p$ and the tame level (see \cite[Lem.\ 3.10.7]{BCGP} for a precise statement).  This implies that the systems of Hecke
  eigenvalues in $\HH^0(\tilde{D}^\star,\O_{\tilde{D}^\star})$ are Eisenstein,
  as required. \end{proof}

 \begin{thm}\label{thm-ECpadic}For any irreducible  representation $\bar{\rho}: G_{\mathbb{Q}}
\rightarrow \mathrm{GSp}_{4}(\Fbar_p)$, the localization \[V := \mathrm{RHom}_{\mathfrak{b}}( \lambda, \mathrm{R}\Gamma(\Sh^{\tor}_{K^p}, \Qpbar)^{\la})_{ \mathfrak{m}_{\bar{\rho}}}^{\ord}=CC_{\lambda}(M(\mf{g})_{\lambda})_{ \mathfrak{m}_{\bar{\rho}}}^{\ord}\] is concentrated in degree $3$. 
 Moreover, there is a $G_{\qq_p} \times \mathbf{T}_{K^p} \times T(\qq_p)$-equivariant  filtration $\{Fil^iV_{\C_p}\}_{i=0,1,2, 3}$ on $V \otimes \C_p$  and $$\mathrm{Gr}^i V_{\C_p}= \big( \HH^{i}_{^iw}( \Sh^{\tor}_{K^p}, \omega_{^iw}^{\dag, -w_{0,M}({^iw}\cdot\lambda + 2 \rho^M)})\otimes \C_p(-{^iw}^{-1}w_{0,M} \rho-\rho)\big)_{\mathfrak{m}_{\bar{\rho}}}^{\ord}.$$
 The Sen operator is scalar on  $\mathrm{Gr}^i V _{\C_p}$ and acts via $\frac{-\lambda_1 - \lambda_2 -w}{2}$, $\frac{-2+\lambda_1 - \lambda_2 -w}{2}$, $\frac{-4 + \lambda_2 - \lambda_1 -w}{2}$, $\frac{-6 +\lambda_1+\lambda_2  -w}{2}$ for $i=3,2,1,0$ respectively.  
\end{thm}
\begin{proof}
  This is immediate from Proposition~\ref{prop-pre-HT}, Theorem~\ref{thm-van-bound} and Corollary~\ref{cor:quasi-iso-of-usual-to-cuspidal}.
\end{proof}

\subsection{Sen and Cousin}\label{sectsencous} We now specialize to $\lambda = (1,1; w)$. 

\subsubsection{The Cousin map for the classical sheaves}%
  All the action will be happening on the union of Bruhat strata (in fact a $Q$-orbit) $C_{^3w,Q} = C_{^3w} \cup C_{^2w}$. 
We have an extension over $\pi_{HT}^{-1}(C_{^3w,Q})$, corresponding to the stratification of $C_{^3w,Q}$ into $B$-orbits, with $ j_{^3w, \Sh^{\tor}_{K^p}} :  \pi_{HT}^{-1}( C_{^3w}) \hookrightarrow \pi_{HT}^{-1}(C_{^3w,Q})$:
\numequation\label{cous-eq}
 0 \rightarrow  (j_{^3w, \Sh^{\tor}_{K^p}})_! \omega^ {(1,1;-w),\sm}\vert_{\pi_{HT}^{-1}(C_{^3w})} \rightarrow  \omega^ {(1,1;-w),\sm}\vert_{\pi_{HT}^{-1}(C_{^3w,Q})}   \rightarrow \omega^ {(1,1;-w),\sm}\vert_{\pi_{HT}^{-1}(C_{^2w})} \rightarrow 0  
 \end{equation}
\begin{prop}\label{prop-cousinsmall}\leavevmode
\begin{enumerate}
\item  The natural map: 
$$\mathrm{R}\Gamma_c( \pi_{HT}^{-1}(C_{^3w,Q}), \omega^{(1,1;-w),\sm})    \rightarrow \mathrm{R}\Gamma( \Sh_{K^p}^{\tor}, \omega^{(1,1;-w),\sm}) $$ induces a quasi-isomorphism on the ordinary part (the slope $= -(1,1;w)$-part). 
\item Moreover, $\mathrm{R}\Gamma_c( \pi_{HT}^{-1}(C_{^3w,Q}), \omega^{(1,1;-w),\sm})^{\ord}$ is computed by the following complex in degrees $2,3$ where $Cous$ is induced by the class of the extension \eqref{cous-eq}
$$ \HH^2_{^2w}(\Sh^{\tor}_{K^p}, \omega^{(1,1;-w),\sm})^{\ord} \stackrel{Cous}\rightarrow \HH^3_{^3w}(\Sh^{\tor}_{K^p}, \omega^{ (1,1;-w),\sm})^{\ord}.$$
\end{enumerate}
\end{prop}
\begin{proof} By Proposition~\ref{prop:abstract-Cousin-spectral-sequence},  %
  we have a
  spectral sequence (the Cousin spectral sequence) from local
  cohomologies converging to classical cohomology. The first statement is
  equivalent to the vanishing of the ordinary part of the higher Coleman
  theories for the elements $^0w$ and $^1w$. By Theorem
  \ref{thm-bounds-on-slopes}, we find that the slopes on $\HH^i_{^0w}(\Sh^{\tor}_{K^p}, \omega^{(1,1;-w),\sm})$ are $\geq -(2,2;w)$ and the slopes on 
  $\HH^i_{^1w}(\Sh^{\tor}_{K^p}, \omega^{(1,1;-w),\sm})$ are $\geq (0,-2;-w)$. Since $-(1,1;w) + \gamma = -(2,2;w)$ and $-(1,1;w) + \alpha = (0,-2;-w)$ we conclude that the ordinary part vanishes. 
  The second statement is a consequence of  Theorem \ref{coro-vanishingThm} and Lemma \ref{lem-complement-vanishingTHM}.   
\end{proof}

\subsubsection{The Cousin map for the big sheaves}
Applying the functor $VB^0$ to the sheaf $HCS_{Q,^3w, \eta}(  M(\mathfrak{m}_{^3w})_\lambda)$ of ~\eqref{eqn:HCS-ses} %
and twisting the $B(\qq_p)$-action by $\lambda-\eta$ yields an extension %
over
$\pi_{HT}^{-1}(C_{^3w,Q})$: %
\numequation\label{extension-cousin-equation} 0 \rightarrow  j_{^3w, \Sh^{\tor}_{K^p}} \omega_{^3w}^{\dag, (1,1;-w)}  \rightarrow  \VBzero(HCS_{Q,^3w, \eta}(  M(\mathfrak{m}_{^3w})_\lambda))   \otimes \C_p(\lambda-\eta) \rightarrow \omega_{^2w}^{\dag, (1,1;-w)}  \otimes \C_p ((2,0;0)) \rightarrow 0
\end{equation}
The natural map in $\ocal(\mathfrak{m}_{^3w}, \mathfrak{b}_{M_{^3w}})$ (mapping
a Verma of dominant weight to its finite dimensional quotient)
$M(\mathfrak{m}_{^3w})_\lambda \rightarrow L(\mathfrak{m}_{^3w})_\lambda$ yields
a map $$HCS_{Q,^3w, \eta}(  L(\mathfrak{m}_{^3w})_\lambda) \rightarrow
HCS_{Q,^3w, \eta}(  M(\mathfrak{m}_{^3w})_\lambda).$$ As in 
  Lemma~\ref{lem-VBandclassicalsheaf}, applying $VB^0$ gives a
map   $$ \omega^ {(1,1;-w),\sm} \otimes \C_p(\eta)\vert_{\pi_{HT}^{-1}(C_{^3w,Q})} \rightarrow  \VBzero(HCS_{Q,^3w, \eta}(  M(\mathfrak{m}_{^3w})_\lambda)).$$ 
We deduce that there is a  map of  extensions from \eqref{cous-eq} to  \eqref{extension-cousin-equation}  (for clarity we drop the twist of the
$B(\qq_p)$-action in this diagram): %

\begin{eqnarray*}
  \resizebox{\textwidth}{!}{
\xymatrix{0\ar[r]&   j_{^3w, \Sh^{\tor}_{K^p}} \omega_{^3w}^{\dag, (1,1;-w)}  \ar[r]  & \VBzero(HCS_{Q,^3w, \eta}(  M(\mathfrak{m}_{^3w})_\lambda))  \ar[r] & \omega_{^2w}^{\dag, (1,1;-w)} \ar[r]&0  \\
0\ar[r]& (j_{^3w, \Sh^{\tor}_{K^p}})_! \omega^ {(1,1;-w),\sm}\vert_{\pi_{HT}^{-1}(C_{^3w})} \ar[u] \ar[r] &   \omega^ {(1,1;-w),\sm}\vert_{\pi_{HT}^{-1}(C_{^3w,Q})} \ar[u] \ar[r] & \omega^ {(1,1;-w),\sm}\vert_{\pi_{HT}^{-1}(C_{^2w})} \ar[u]\ar[r]&0 }}
\end{eqnarray*}

\begin{prop}\label{prop-COusin}The maps 
\begin{eqnarray*}
\mathrm{R}\Gamma_{^2w}(\Sh^{\tor}_{K^p}, \omega^{(1,1;-w),\sm})  \otimes \C_p ((1,1;w))^{\ord} 
&\rightarrow &\mathrm{R}\Gamma_{^2w}(\Sh^{\tor}_{K^p}, \omega_{^2w}^{\dag, (1,1;-w)}) \otimes \C_p ((2,0;0))^{\ord} \\
  \mathrm{R}\Gamma_{^3w}(\Sh^{\tor}_{K^p}, \omega^{(1,1;-w),\sm})  \otimes \C_p ((1,1;w))^{\ord} 
&\rightarrow& \mathrm{R}\Gamma_{^3w}(\Sh^{\tor}_{K^p}, \omega_{^3w}^{\dag, (1,1;-w)})^{\ord}
\end{eqnarray*}  are quasi-isomorphisms.

  Consequently we have a quasi-isomorphism $$\mathrm{R}\Gamma( \Sh_{K^p}^{\tor}, \omega^{(1,1;-w),\sm}) \otimes \C_p ((1,1;w))^{\ord} =$$
$$ \HH^2_{^2w}(\Sh^{\tor}_{K^p}, \omega_{^2w}^{\dag, (1,1;-w)}) \otimes \C_p ((2,0;0))^{\ord} \stackrel{Cous}\rightarrow \HH^3_{^3w}(\Sh^{\tor}_{K^p}, \omega_{^3w}^{\dag, (1,1;-w)})^{\ord} $$ where the complex is in degree $[2,3]$ and the map $Cous$ is induced by the class of the extension \eqref{extension-cousin-equation}.
\end{prop}
\begin{proof} We have the BGG short exact sequence  $$0 \rightarrow M(\mathfrak{m}_{^3w})_{(0,2;w)} \rightarrow M(\mathfrak{m}_{^3w})_{(1,1;w)} \rightarrow L(\mathfrak{m}_{^3w})_{(1,1;w)} \rightarrow 0.$$
Applying $HC_{^3w,\lambda}$ gives a triangle 
$$HC_{^3w,\lambda}(L(\mathfrak{m}_{^3w})_{(1,1;w)}) \rightarrow HC_{^3w,\lambda}(M(\mathfrak{m}_{^3w})_{(1,1;w)}) \rightarrow HC_{^3w,\lambda}(  M(\mathfrak{m}_{^3w})_{(0,2;w)}) \stackrel{+1}\rightarrow$$
By Theorem \ref{thm-bounds-on-slopes}, the ordinary part of $HC_{^3w,\lambda}(  M(\mathfrak{m}_{^3w})_{(0,2;w)})$ is trivial, so that we get a quasi-isomorphism $HC_{^3w,\lambda}(L(\mathfrak{m}_{^3w})_{(1,1;w)}))^{\ord} \rightarrow HC_{^3w,\lambda}(M(\mathfrak{m}_{^3w})_{(1,1;w)})^{\ord}$. 
This translates into the quasi-isomorphism $\mathrm{R}\Gamma_{^3w}(\Sh^{\tor}_{K^p}, \omega^{(1,1;-w),\sm})  \otimes \C_p ((1,1;w))^{\ord} 
\rightarrow \mathrm{R}\Gamma_{^3w}(\Sh^{\tor}_{K^p}, \omega_{^3w}^{\dag, (1,1;-w)})^{\ord}$. 
The quasi-isomorphism $\mathrm{R}\Gamma_{^2w}(\Sh^{\tor}_{K^p}, \omega^{(1,1;-w),\sm})  \otimes \C_p ((1,1;w))^{\ord} 
\rightarrow \mathrm{R}\Gamma_{^2w}(\Sh^{\tor}_{K^p}, \omega_{^2w}^{\dag, (1,1;-w)}) \otimes \C_p ((2,0;0))^{\ord}$ follows by similar considerations. 
The second part is then
immediate from Proposition~\ref{prop-cousinsmall}.%
\end{proof}

\begin{rem} We also have a quasi-isomorphism $$\mathrm{R}\Gamma( \Sh_{K^p}^{\tor}, \omega^{(2,2;-w),\sm}(-D)) \otimes \C_p ((2,2;w))^{\ord} =$$ 
$$ \HH^0_{^0w}(\Sh^{\tor}_{K^p}, \omega_{^0w}^{\dag, (1,1;-w)}(-D))\otimes \C_p((3,3; 0))^{\ord}  \stackrel{Cous}\rightarrow \HH^1_{^1w}(\Sh^{\tor}_{K^p}, \omega_{^1w}^{\dag, (2,2;-w)}(-D)) \otimes \C_p((-1,3; 0))^{\ord} $$ where the complex is in degrees $[0,1]$. 
This statement is  Serre dual to  Proposition \ref{prop-cousinsmall}.
\end{rem}

\subsubsection{The Sen map}
Let   $V = \HH^3(\mathrm{RHom}_{\mathfrak{b}}( \lambda, \mathrm{R}\Gamma(\Sh^{\tor}_{K^p}, \Qpbar)^{\la})_{ \mathfrak{m}_{\bar{\rho}}}^{\ord})$. By Theorem \ref{thm-ECpadic},  $V_{\C_p}$ carries a filtration where 
\begin{itemize}
\item $\mathrm{Gr}^3 V_{\C_p}  = H^3_{^3w}( \Sh^{\tor}_{K^p}, \omega_{^3w}^{\dag, (1,1;-w)}) \otimes \C_p(0,0; 0)^{\ord}_{ \mathfrak{m}_{\bar{\rho}}}$,
\item $\mathrm{Gr}^2 V_{\C_p}  = H^2_{^2w}( \Sh^{\tor}_{K^p}, \omega_{^2w}^{\dag, (1,1;-w)}) \otimes \C_p(2,0; 0)^{\ord}_{ \mathfrak{m}_{\bar{\rho}}}$,
\item $\mathrm{Gr}^1 V_{\C_p}  = H^1_{^1w}( \Sh^{\tor}_{K^p}, \omega_{^1w}^{\dag, (2,2;-w)}) \otimes \C_p(-1,3; 0)^{\ord}_{ \mathfrak{m}_{\bar{\rho}}}$,
\item $\mathrm{Gr}^0 V_{\C_p}  = H^0_{^0w}( \Sh^{\tor}_{K^p}, \omega_{^0w}^{\dag, (2,2;-w)}) \otimes \C_p(3,3; 0)^{\ord}_{ \mathfrak{m}_{\bar{\rho}}}$.
\end{itemize}
Acting on $V_{\C_p}$ we have   a Sen operator whose eigenvalues are $-1 - \frac{w}{2}, -1 - \frac{w}{2}, -2 - \frac{w}{2}, -2 - \frac{w}{2}$. 
The generalized Hodge--Tate weight $1 + \frac{w}{2}$-part of  $V_{\C_p}$ fits in  the following short exact sequence (where $\Theta$ is the Sen operator and we drop the twist of the $B(\qq_p)$-action to lighten the notation): 

\numequation\label{eqn:Theta-ses} 0 \rightarrow  \HH^3_{^3w}(\Sh^{\tor}_{K^p}, \omega_{^3w}^{\dag, (1,1;-w)})^{\ord}_{ \mathfrak{m}_{\bar{\rho}}} \rightarrow V_{{\C_p}}[(\Theta +1+ \frac{w}{2})^2] \rightarrow 
 \HH^2_{^2w}(\Sh^{\tor}_{K^p}, \omega_{^2w}^{\dag, (1,1;-w)})^{\ord}_{ \mathfrak{m}_{\bar{\rho}}}\rightarrow 0.\end{equation}
 Since $ Sen:=\Theta +1+ \frac{w}{2} $ is nilpotent, it induces a  map: 
 \numequation\label{eqn:the-Sen-operator-we-will-use}\HH^2_{^2w}(\Sh^{\tor}_{K^p}, \omega_{^2w}^{\dag, (1,1;-w)}) \otimes \C_p(2,0;0)^{\ord}_{ \mathfrak{m}_{\bar{\rho}}} \stackrel{Sen}\rightarrow \HH^3_{^3w}(\Sh^{\tor}_{K^p}, \omega_{^3w}^{\dag, (1,1;-w)})^{\ord}_{ \mathfrak{m}_{\bar{\rho}}}.\end{equation}
 
 Similarly, by looking at the weight $2 + \frac{w}{2}$-part of  $V$ we obtain the following map: 
 
 $$\HH^0_{^0w}(\Sh^{\tor}_{K^p}, \omega_{^0w}^{\dag, (2,2;-w)}) \otimes  \C_p(3,3; 0)^{\ord}_{ \mathfrak{m}_{\bar{\rho}}} \stackrel{Sen}\rightarrow \HH^1_{^1w}(\Sh^{\tor}_{K^p}, \omega_{^1w}^{\dag, (2,2;-w)}) \otimes \C_p(-1,3; 0)^{\ord}_{ \mathfrak{m}_{\bar{\rho}}}$$

 \subsubsection{Comparison between the Sen and Cousin map}
 The following theorem is one of the main results in this section. It is a
 generalization of \cite[Thm.\ 5.3.18]{MR4390302}, in the modular curve case. We
 will  follow the method of proof of   \cite[Thm.\ 6.1]{PilloniVB}.
 
 \begin{thm}\label{thm-Sen-Cousin} The two maps
   \[Cous, Sen:\HH^2_{^2w}(\Sh^{\tor}_{K^p}, \omega_{^2w}^{\dag, (1,1;-w)})
     \otimes \C_p(2,0;0)^{\ord}_{ \mathfrak{m}_{\bar{\rho}}} \rightarrow
     \HH^3_{^3w}(\Sh^{\tor}_{K^p}, \omega_{^3w}^{\dag, (1,1;-w)})^{\ord}_{
       \mathfrak{m}_{\bar{\rho}}}\](coming respectively from
   Proposition~\ref{prop-COusin} and~\eqref{eqn:the-Sen-operator-we-will-use})agree up to a non-zero scalar.%
 \end{thm}%
 \begin{proof}
We consider $%
 \RHom_{\mathfrak{b}, \star_2}(\lambda,
 \oscr^{\la}_{\Sh^{\tor}_{K^p}})\vert_{\pi_{HT}^{-1}C_{^3w,Q}} =
 VB^{\red}(\Loc(M(\mathfrak{g})_\lambda)\vert_{\pi_{HT}^{-1}C_{^3w,Q}})$, which  fits in the triangle (obtained by applying $VB^{\red}$ to Proposition \ref{prop:we-have-exact-triangle}): 

$$ \Ext^0_{\mathfrak{b}, \star_2}(\lambda, \oscr^{\la}_{\Sh^{\tor}_{K^p}}\vert_{\pi_{HT}^{-1}(C_{^3w,Q})}) \rightarrow \RHom_{\mathfrak{b}, \star_2}(\lambda, \oscr^{\la}_{\Sh^{\tor}_{K^p}}\vert_{\pi_{HT}^{-1}(C_{^3w,Q})})  $$ $$\rightarrow \mathrm{Ext}^1_{\mathfrak{b}, \star_2}(\lambda, \oscr^{\la}_{\Sh^{\tor}_{K^p}}\vert_{\pi_{HT}^{-1}(C_{^3w,Q})})[-1] \stackrel{+1}{\rightarrow}$$
where 
$ \mathrm{Ext}^0_{\mathfrak{b}, \star_2}(\lambda, \oscr^{\la}_{\Sh^{\tor}_{K^p}}\vert_{\pi_{HT}^{-1}C_{^3w,Q}}) = (j_{^3w, \Sh^{\tor}_{K^p}})_! \omega_{^3w}^{\dag, (1,1;-w)}$ and $\mathrm{Ext}^1_{\mathfrak{b}, \star_2}(\lambda, \oscr^{\la}_{\Sh^{\tor}_{K^p}}\vert_{\pi_{HT}^{-1}C_{^3w,Q}}) = \omega_{^2w}^{\dag, (1,1;-w)} \otimes \C_p((2,0,0))$.
 
 On $\RHom_{\mathfrak{b}, \star_2}(\lambda, \oscr^{\la}_{\Sh^{\tor}_{K^p}}\vert_{\pi_{HT}^{-1}(C_{^3w,Q})})$, we
 have  $(\Theta+1 + \frac{w}{2})^2 =0$. 
Let us introduce some simplifying notations and denote by:
\begin{itemize}
\item $X = (j_{^3w, \Sh^{\tor}_{K^p}})_! \omega_{^3w}^{\dag, (1,1;-w)}$, 
\item $Y = \RHom_{\mathfrak{b}, \star_2}(\lambda, \oscr^{\la}_{\Sh^{\tor}_{K^p}}\vert_{\pi_{HT}^{-1}(C_{^3w,Q})})$,
\item $Z = \omega_{^2w}^{\dag, (1,1;-w)} \otimes \C_p((2,0,0))[-1]$,
\item $W = \RHom_{\mathfrak{b}, \star_2, \Theta}(\lambda, -1 - \frac{w}{2} , \oscr^{\la}_{\Sh^{\tor}_{K^p}}\vert_{\pi_{HT}^{-1}(C_{^3w,Q})})$. 
\end{itemize}Applying~\cite[\href{https://stacks.math.columbia.edu/tag/05R0}{Tag 05R0}]{stacks-project}  to the commutative diagram \begin{eqnarray*}
\xymatrix{ 
X \ar[r] \ar[d]^{0} & Y \ar[d]^{(\Theta+1 + \frac{w}{2})}  \\
X \ar[r]  & Y }
\end{eqnarray*}
we obtain the following  commutative diagram, where all lines and columns are part of  distinguished triangles:
\begin{eqnarray*}
\xymatrix{ X \oplus X[-1] \ar[r] \ar[d]  & W \ar[r] \ar[d] & Z \oplus Z[-1] \ar[r]  \ar[d] & X[1] \oplus X \ar[d] \\
X \ar[r] \ar[d]^{0} & Y \ar[r] \ar[d]^{(\Theta+1 + \frac{w}{2})} & Z \ar[r] \ar[d]^{0} & X[1] \ar[d]^{0}  \\
X \ar[r]  & Y \ar[r]  & Z \ar[r]  & X[1]}
\end{eqnarray*}
The top horizontal triangle can be written as
$$   (j_{^3w, \Sh^{\tor}_{K^p}})_! \omega_{^3w}^{\dag, (1,1;-w)} \oplus  (j_{^3w, \Sh^{\tor}_{K^p}})_! \omega_{^3w}^{\dag, (1,1;-w)}[-1] \rightarrow  \RHom_{\mathfrak{b}, \star_2, \Theta}(\lambda, -1 - \frac{w}{2} , \oscr^{\la}_{\Sh^{\tor}_{K^p}})  $$ $$\rightarrow  \omega_{^2w}^{\dag, (1,1;-w)}[-1] \oplus \omega_{^2w}^{\dag, (1,1;-w)}[-2] \stackrel{+1}\rightarrow$$(compare~\eqref{eqn:triangle-derived-mu-invariants}; we again drop the twist of the $B(\qq_p)$-action for simplicity and we note that taking the $\Theta+1 + \frac{w}{2}$-cohomology localizes over $\pi_{HT}^{-1}(C_{^3w,Q})$ so we also drop it from the notation.)

Taking cohomology yields a long exact sequence: 
$$\HH^3(\RHom_{\mathfrak{b}, \star_2, \Theta}(\lambda, -1 - \frac{w}{2} , \oscr^{\la}_{\Sh^{\tor}_{K^p}})) \rightarrow  \HH^2_{^2w}(\Sh^{\tor}_{K^p}, \omega^{\dag, (1,1;-w)}) \stackrel{\delta}\rightarrow \HH^3_{^3w}(\Sh^{\tor}_{K^p}, \omega^{\dag, (1,1;-w)})$$
where the map $\delta$ is induced by the class of the extension (extracted from the top horizontal triangle):
\numequation\label{ext-delta} 0 \rightarrow  j_{^3w, \Sh^{\tor}_{K^p}} \omega_{^3w}^{\dag, (1,1;-w)} \rightarrow   \mathrm{Ext}^1_{\mathfrak{b}, \star_2, \Theta}(\lambda, -1 - \frac{w}{2} , \oscr^{\la}_{\Sh^{\tor}_{K^p}}) \rightarrow \omega_{^2w}^{\dag, (1,1;-w)} \otimes \C_p((2,0;0)) \rightarrow 0.
\end{equation}By definition, the map~$Sen$ is equal to~$\delta$.
Now, by Theorem \ref{thm-extequal} and Theorem \ref{thm-VB}, the extensions
\eqref{extension-cousin-equation} and \eqref{ext-delta} agree up to a non-zero
scalar. Thus the maps $Cous,Sen$ agree up to a non-zero
scalar, as claimed.
\end{proof}

\subsection{The Eichler--Shimura relation and
  semi-simplicity}\label{subsec:Eichler--Shimura-Nekovar}
The following is a special case of a result of Nekov\'{a}\v{r},
\cite{MR3942040}. If~$r$ is a finite-dimensional representation of a
group~$\Gamma$, and $g\in\Gamma$, then we write $\cha_{r(g)}$ for the
characteristic polynomial of~$r(g)$.

\begin{prop}
  \label{prop:Eichler--Shimura-semi-simple}Let $\rho:G_{\Q}\to\GSp_4
  (\Qpbar)$ and $s:G_{\Q}\to\GL_n(\Qpbar)$ be continuous representations \emph{(}for some
  $n\ge 1$\emph{)} and assume that
  \begin{enumerate}\item the Zariski closure of $\rho(G_{\Q})$
    contains~$\Sp_4 $, and 
  \item \label{primescondition} for a density one set of primes~$l$, we have
    $\cha_{\rho(\Frob_{l})}\bigl(s(\Frob_{l})\bigr)=0$.
  \end{enumerate}
  Then we have $s\cong \rho^{\oplus m}$ for some
  integer~$m\ge 1$.
\end{prop}
\begin{proof}We claim that the result is an immediate application
  of~\cite[Prop.\ 3.10]{MR3942040} (with $\ell$ replaced by~$p$),
   taking $\Gamma=\Gamma'=G_{\Q}$, $a=r=1$, and
   the representations~$\rho$, $\rho_1$ of~\cite[Prop.\ 3.10]{MR3942040} to
   be~$s$ and~$\rho$ respectively.

  The hypothesis~(C') of~\cite[Prop.\ 3.10]{MR3942040} is immediate
  from our hypotheses, taking~$\Sigma$ to be the set of
  the~$\Frob_{l}$ for primes~$l$ satisfying Condition~(\ref{primescondition}). For hypothesis~(A'), since the Zariski closure
  of $\rho(G_{\Q})$ contains~$\Sp_4$ by hypothesis (and is contained
  in~$\GSp_4$),
  the Lie algebra $\Qpbar\cdot\Lie(\rho(G_{\Q}))$ is
  equal to~$\franksp_{4}$ or ~$\frankgsp_{4}$, 
     and the representation of
  this Lie algebra induced by~$\rho$ is the standard $4$-dimensional
  representation, which is minuscule.  We are thus in the situation of
  part~(3) of~\cite[Prop.\ 3.10]{MR3942040}, and the proposition follows.
\end{proof}

We have the following variant of Proposition~\ref{subsec:Eichler--Shimura-Nekovar}
in the induced case:

\begin{prop} \label{inducedversion}
Let $\rho:G_{\Q}\to\GSp_4
  (\Qpbar)$ and $s:G_{\Q}\to\GL_n(\Qpbar)$ be continuous representations and assume that
  \begin{enumerate}
  \item the Zariski closure of $\rho(G_{\Q})$
    contains~$\SL_2 \times \SL_2 $,
    \item $\rho$ is absolutely irreducible but becomes reducible
  on some index two subgroup~$G_{E}$, and
  \item \label{assumptionfrob} for a density one set of primes~$l$, we have
    $\cha_{\rho(\Frob_{l})}\bigl(s(\Frob_{l})\bigr)=0$.
  \end{enumerate}
  Then we have  $s\cong \rho^{\oplus m}$ for some
  integer~$m\ge 1$.
\end{prop}

\begin{proof}
We may write~$\rho |_{G_E}= \varrho \oplus \varrho^c$, where~$\Gal(E/\Q)$ permutes the factors. 
The assumption that~$\rho$ is symplectic (together with the assumptions on the
image of~$\rho$)
implies that~$\det(\varrho)$ and~$\det(\varrho^c)$
are both the restriction to~$G_E$ of the similitude character of~$\rho$. %

Let~$t$ be an irreducible subquotient of~$s|_{G_E}$.
Assumption~(\ref{assumptionfrob}) implies that we have~$\cha_{\rho(h)}\bigl(t(h)\bigr)=0$
for a dense set of elements~$h \in G_E$.
We may assume that~$\varrho$, $\varrho^c$, and~$t$
all have models over the ring of integers~$\cO$ of some finite extension of~$\Q_p$.
Let~$T$ and~$P$ denote the
 Zariski closures  of~$t$ and~$\varrho \oplus \varrho^c$ respectively. 
By assumption, $T$ is reductive, and 
the Zariski closure of~$t \oplus \varrho \oplus \varrho^c$ inside~$T \oplus P$ is also reductive,
and by Goursat's lemma is the graph of some projections~$\pi_1: T \rightarrow G$,
$\pi_2: P \rightarrow G$ onto a common quotient~$G$.
By the Chebotarev density theorem and continuity,
Assumption~(\ref{assumptionfrob}) implies that the minimal polynomial of any element in~$K = \ker(\pi_1)$ %
divides~$(X-1)^4$.  This is because elements in the image of~$s$
which lie in~$K$ are limits of~$s(\Frob_l)$ for Frobenius elements~$\Frob_l$, and
by assumption these will satisfy~$(s(\Frob_l)-1)^4 \OL^n \subset \pi^m \OL^n$ for larger
and larger~$m$.
This implies that~$K$ is
unipotent, which --- since~$T$ is reductive --- implies that~$K$  is trivial.
Hence~$G=T$ and thus~$T$ is a quotient of~$P$. But now from the fact that~$P$
contains~$\SL_2 \times \SL_2$, we see that the only possibilities for~$t$
up to twist are~$\Sym^i \varrho \otimes \Sym^j \varrho^c$, from which one easily sees that~$t$
must either be~$\varrho$ or~$\varrho^c$, and thus
 any irreducible subquotient of~$s |_{G_E}$
is either~$\varrho$ or~$\varrho^c$.

We claim that~$s |_{G_E}$ cannot contain any non-trivial extensions of~$\varrho$ by~$\varrho$ (equally,
of~$\varrho^c$ by~$\varrho^c$). To see this, note that a generic element in the image of~$\varrho^c$
acts invertibly on~$\varrho$.  Hence the assumptions imply that  for of a dense
set of~$h \in G_E$, the characteristic polynomial of~$\varrho(h)$
annihilates this extension of~$\varrho$ by~$\varrho$.
But then the results follow from the~$\GL_2$-version 
of~\cite[Prop.\ 3.10]{MR3942040} (first proved in~\cite{BLR}).

 Now return to representations of~$G_{\Q}$.
By what we have shown for~$s|_{G_E}$, we deduce that every irreducible
subquotient of~$s$ is isomorphic to~$\rho$.
Hence it suffices to rule out the case that~$s$ is of the form
$$0 \rightarrow \rho \rightarrow W \rightarrow \rho \rightarrow 0$$
for some non-split extension~$W$.
Note that the restriction ~$W|_{G_{E}}$ is an extension: %
$$0 \rightarrow \varrho \oplus \varrho^c \rightarrow W |_{G_E} \rightarrow \varrho \oplus \varrho^c \rightarrow 0.$$
In particular, $W|_{G_E}$  corresponds to a~$\Gal(E/\Q)$-invariant class in
\[\Ext^1_{G_{E}}(\varrho,\varrho^c) \oplus
\Ext^1_{G_{E}}(\varrho,\varrho) \oplus
\Ext^1_{G_{E}}(\varrho^c,\varrho^c) \oplus
\Ext^1_{G_{E}}(\varrho^c,\varrho).\]
As already shown, the corresponding extensions of~$\varrho$ by~$\varrho$ and~$\varrho^c$ by~$\varrho^c$
both split, so 
the projection of this class to both $\Ext^1_{G_{E}}(\varrho^c,\varrho)$
and $\Ext^1_{G_{E}}(\varrho,\varrho^c)$  (they are permuted by the Galois
action) must be non-trivial.
In particular,
we may %
assume that~$W |_{G_E}$ 
has a subquotient which is
 a genuine extension of~$\varrho$ by~$\varrho^c$. %
Our assumptions on the Zariski closure of the image of~$\rho$ imply that a generic element of~$\rho$
is regular semi-simple. Thus, by Chebotarev and assumption~(\ref{assumptionfrob}), it follows
that there is a dense set of elements of~$G_{\Q}$ which act semi-simply on~$W$.
Moreover, since the tensor product of two semi-simple matrices is semi-simple, the same holds
for the tensor product~$W \otimes \rho^{\vee}$ and so consequently also for any subquotient of this representation.
We have an exact sequence:
\numequation \label{tensored}
0 \rightarrow \rho \otimes \rho^{\vee} \rightarrow W \otimes \rho^{\vee}
 \rightarrow \rho \otimes \rho^{\vee} \rightarrow 0.\end{equation}
Write~$\chi$ for the quadratic character of~$\Gal(E/\Q)$, 
 and let~$\eta$ be the similitude character of~$\rho$. %
 We have (compare~\cite[\S 7.5.16]{BCGP}) a decomposition 
$$\rho \otimes \rho^{\vee}
\simeq \Qpbar \oplus \Qpbar(\chi) \oplus \Ind^{G_{\Q}}_{G_E}
\ad^0(\varrho) \oplus \As(\varrho) \otimes \eta^{-1} \oplus \As(\varrho) \otimes
\eta^{-1} \otimes \chi,$$ where $\As(\varrho)$ is the Asai representation (i.e.\
the tensor induction).
Hence taking  suitable subquotients of~(\ref{tensored}), we arrive at  a pair of extensions
\numequation \label{oneisnonsplit}
\begin{aligned}
0 \rightarrow \As(\varrho) \otimes \eta^{-1}  \rightarrow U \rightarrow \Qpbar \rightarrow 0, \\
0 \rightarrow \As(\varrho) \otimes \eta^{-1} \otimes \chi \rightarrow V \rightarrow \Qpbar \rightarrow 0.
\end{aligned}
\end{equation}

We claim that at least one of these sequences must be non-split. The point is that
$$\Ext^1_{G_{E}}(\varrho,\varrho^c) =  H^1(E,\Hom(\varrho,\varrho^c)),$$
but~$\Hom(\varrho,\varrho^c)$ is the restriction of~$\As(\varrho) \otimes \eta^{-1}$ to~$G_E$, and thus,
by Shapiro's lemma, %
 we have
$$\begin{aligned}
\Ext^1_{G_E}(\varrho,\varrho^c)  & = H^1(E,\left(\As(\varrho) \otimes \eta^{-1}\right) |_{G_E}) \\
 & = H^1(\Q,\left(\As(\varrho) \otimes \eta^{-1}\right) \otimes \Ind^{G_{\Q}}_{G_{E}} \Qbar_p) \\
& = 
H^1(\Q,\As(\varrho) \otimes \eta^{-1} \otimes \chi) \oplus
H^1(\Q,\As(\varrho) \otimes \eta^{-1}); \end{aligned}$$
 and by construction, the elements of
the right hand side corresponding to our non-split extension of ~$\varrho$
by~$\varrho^c$ coming from~$W|_{G_E}$ are the extensions~$U,V$ of ~(\ref{oneisnonsplit}). 
We consider the case that~$U$ is non-split, the case of~$V$ being entirely similar.

To complete the proof, it suffices to show that for any such non-split extension, there cannot
be a dense set of~$g \in G_{\Q}$ which act semi-simply.
Write~$A = \As(\varrho) \otimes \eta^{-1}$, %
and let~$G$ and~$H$ be the Zariski closures of the images of~$G_{\Q}$
and~$G_{E}$ on~$U$ respectively.
By constriction, the Zariski closure of the image of~$G_{\Q}$ on~$A$ is the orthogonal
group~$\mathrm{O}_4$, and the Zariski closure of the image of~$G_E$ on~$A$
is the index two subgroup~$\mathrm{SO}_4$ which is isomorphic to  the image of~$\SL_2 \times \SL_2$.
Any generic~$h \in G_E$ \emph{will} act semi-simply on~$U$ because it will
have distinct eigenvalues. However, any~$g \in G_{\Q} \setminus G_E$  has~$1$ as an eigenvalue on
both~$A = \As(\varrho) \otimes \eta^{-1}$ and~$A \otimes \chi$
(the eigenvalues in either case take the form~$1,-1,\lambda,-\lambda$ for
some~$\lambda$), so that in particular the eigenvalues of~$g$ on~$U$  %
are contained within the eigenvalues of~$g$ on~$A$.
By our semi-simplicity hypothesis, there is therefore a dense set  of such~$g$
with the property that the image of~$g$ in~$\End(U)$ is annihilated
by the characteristic polynomial of~$g$ on~$A$. By continuity, this extends to \emph{all} elements
of~$G_{\Q} \setminus G_E$ and also to all elements of~$G\setminus H$.

Now, $G$ is a subgroup of the semi-direct product ~$N \rtimes \mathrm{O}_4$,
where~$N$ is the standard representation of~$\mathrm{O}_4$ of dimension~$4$;
and the projection $G\to \mathrm{O}_4$ is surjective. Since~$N$ is abelian, the
conjugation action of~$G$ on~$N$ factors through this surjection, and since~$N$
is an irreducible representation of ~$\mathrm{O}_4$, we see that
either~$G=\mathrm{O}_4$ (in which case the extension~\eqref{oneisnonsplit}
splits, and we are done), or $G=N \rtimes \mathrm{O}_4$, in which case there are
elements $g\in G\setminus H$ whose minimal polynomial has degree~$5$, a contradiction.
\end{proof}

\subsection{A classicality theorem}\label{subsec:the-classicity-thm}

Let~$S$ be the finite set of primes at which~$K^p$
is not hyperspecial, together with the prime~$p$. In this section we will consider an ordinary overconvergent modular form
$f\in\HH^0_{^0w}( \Sh_{K^p}^{\tor}, \omega^{(2,2;-w),\sm})^{\ord}$, which is an
eigenform for $T(\Qp)$ as well as for the spherical Hecke operators at the
places not contained in~$S$.  We write $\chi_f^{sm} : T(\qq_p) \rightarrow
\Qpbar^\times$  for the smooth character corresponding to $f$. Via the
identification of dual groups, $\chi_f^{sm}$ induces a  cocharacter
$\qq_p^\times \rightarrow T(\Qpbar)$ %
that we denote by  $t \mapsto \mathrm{diag}(\chi_1(t),\chi_2(t), \chi_3(t), \chi_4(t))$.
We write $$\m_f\subseteq
\Qpbar[T(\Qp)]\otimes\bigotimes'_{l\not\in
  S}\Qpbar[\GSp_4(\Q_l)/\kern-0.2em{/}{\GSp_4(\Z_l)}]$$ %
for the maximal ideal corresponding to $f$. 

\begin{lem}\label{lem-Galoisrepattachedp-adic} We have a continuous  semi-simple
  Galois representation  $\rho_f:G_\Q\to\GSp_4(\Qpbar)$, which satisfies the
  following properties (where $P_{\ell}(X)$ is as in~\eqref{eq:P-ell}):
\begin{enumerate}%
\item $\rho_f$ is unramified at primes $\ell \notin S$, and $$P_\ell(X) \mod \mathfrak{m}_f = \det(X- \rho_f(\mathrm{Frob}_\ell)).$$ 
\item $\rho_f\vert_{G_{\qq_p}}$    can be conjugated to a representation 
$\rho_f \vert_{G_{\qq_p}} : G_{\qq_p} \rightarrow B(\Qpbar)$ where  the diagonal is given, via class field theory, by the
cocharacter  
$$z \mapsto \mathrm{diag}( \chi^{-1}_4(z) z^{-1- \frac{w}{2}}, \chi^{-1}_3(z) z^{-1- \frac{w}{2}}, \chi^{-1}_2(z)^{-1} z^{-2- \frac{w}{2}},  \chi_1^{-1}(z)^{-1} z^{-2- \frac{w}{2}})$$ for $z \in \Z_p^\times$,
 and $$p \mapsto \mathrm{diag}( \chi^{-1}_4(p) p^{-1- \frac{w}{2}}, \chi^{-1}_3(p) p^{-1- \frac{w}{2}}, \chi^{-1}_2(p)^{-1} p^{-2- \frac{w}{2}},  \chi^{-1}_1(p)^{-1} p^{-2- \frac{w}{2}}).$$
 
\end{enumerate}
\end{lem} 

\begin{proof} Note that any cohomological, $C$-algebraic automorphic  representation $\pi$ has an associated Galois representation
$\rho_{\pi,p}$ (see Section \ref{subsec-GalAAR}). This Galois representation is
furthermore ordinary if $\pi_p$ is ordinary. %
By a standard
argument using $p$-adic families, %
we can interpolate the Galois representations~ $\rho_{\pi,p}$ associated
to~$\pi$ of regular weight (see Section~\ref{subsec-GalAAR}),
and we define~$\rho_f$ to be the representation corresponding to the
interpolation of $\rho_{\pi,p}^\vee \otimes \varepsilon^{-3}$. %
\end{proof}

\begin{rem} The reason for considering $\rho_f$ rather than $\rho_f^\vee \otimes \varepsilon^{-3}$ is that $\rho_f$ is the Galois representation we are likely to realize in the completed cohomology, in view of Theorem \ref{thm-ESRFC}. 
\end{rem}

We
will also assume that $\rhobar_f$ is irreducible and write
$\mathfrak{m}_{\rhobar_f}$ for the corresponding maximal ideal of the spherical
Hecke algebra with $\Fpbar$ coefficients, as in the previous sections. %
We let  $V = \HH^3(\mathrm{RHom}_{\mathfrak{b}} ( \lambda,
\mathrm{R}\Gamma(\Sh_{K^p}^{\tor},
\Qpbar)^{\la})^{\ord}_{\mathfrak{m}_{\rhobar_f}})$, but we think of this space
as $CC_0(M(\mathfrak{g})_\lambda)^{\ord}_{\mathfrak{m}_{\rhobar_f}}$ and not as
$CC_\lambda(M(\mathfrak{g})_\lambda)^{\ord}_{\mathfrak{m}_{\rhobar_f}}$. %
In other words we now twist the $B(\qq_p)$ action to make it  smooth (and not $\lambda$-smooth as in section \ref{sectsencous}). 

\begin{lem}\label{lem-Galois-rep-cc} Let $f \in \HH^0_{^0w}( \Sh_{K^p}^{\tor}, \omega^{(2,2;-w),\sm})^{\ord}$ be an ordinary overconvergent modular eigenform with Galois representation
$\rho_f: G_{\qq} \rightarrow \mathrm{GSp}_4(\Qpbar)$.  Let $\mathfrak{m}_f$ be the corresponding maximal ideal. 
We assume that $\rhobar_f$ is irreducible and either %
\begin{enumerate}
\item the Zariski closure of $\rho_f(G_{\qq})$ contains $\mathrm{Sp}_4$; or
\item the Zariski closure of $\rho_{f}(G_{\Q})$
    contains~$\SL_2 \times \SL_2 $, and $\rho$ is irreducible but becomes reducible
  on some index two subgroup~$G_{E}$.
\end{enumerate}
Let $V = \HH^3(\mathrm{RHom}_{\mathfrak{b}} ( \lambda, \mathrm{R}\Gamma(\Sh_{K^p}^{\tor}, \Qpbar)^{\la})^{\ord}_{\mathfrak{m}_{\rhobar_f}})$. 
Then $V[\mathfrak{m}_f] = \rho_f \otimes_{\Qpbar} W$ for some finite-dimensional
vector space $W\neq 0$. %
\end{lem}
\begin{proof} Note that $V$ carries a global Galois action since by results recalled in  \ref{subsect-def-compcoh},    $\mathrm{R}\Gamma(\Sh_{K^p}^{\tor}, \Qpbar) =  \mathrm{R}\Gamma(\Sh^{alg}_{K^p}, \Qpbar)$. By Theorem \ref{thm-ECpadic},
 $V_{\Cp}:=V \otimes_{\Qpbar} \C_p$ has a decreasing filtration with graded pieces
$\{\mathrm{Gr}^iV_{\C_p}\}_{i=0,1,2, 3}$ which are  certain explicit spaces of
ordinary higher Coleman theories localized at $\mathfrak{m}_{\rhobar_f}$. In
particular (bearing in mind Proposition~\ref{prop-COusin}),
$$ \mathrm{Gr}^0V_{\C_p}= (\HH^0_{^0w}(\Sh^{\tor}_{K^p}, \omega^{(2,2;-w),\sm}))^{\ord}_{\mathfrak{m}_{\rhobar_f}}.$$
Let $\chi$ be the Nebentypus of $f$ (the finite order character giving the action of $T(\ZZ_p)$ on ~$f$). Then each $\mathrm{Gr}^iV_{\C_p}[\chi]$ is finite dimensional (since we are
fixing the slope to be ordinary, and also the action of $T(\Z_p)$). Taking the $\chi$-isotypic part is an exact operation for smooth $T(\ZZ_p)$-modules in characteristic $0$. We deduce that $V_{\Cp}[\chi]$ is finite-dimensional and has a filtration with graded pieces the 
$\{\mathrm{Gr}^iV_{\C_p}[\chi]\}_{i=0,1,2, 3}$. 
Since $\HH^0_{^0w}(\Sh^{\tor}_{K^p}, \omega^{(2,2;-w),\sm})[\mathfrak{m}_f] \neq 0$, this implies that
$V[\mathfrak{m}_f] \neq 0$ is finite dimensional. %
The result follows from the Eichler--Shimura relation (Theorem \ref{thm-ESRFC})
together with Proposition~ \ref{prop:Eichler--Shimura-semi-simple} and
Corollary~\ref{inducedversion}.%
\end{proof}

By Lemma \ref{lem-Galoisrepattachedp-adic}, the representation  $\rho_f\vert_{G_{\qq_p}}$ is ordinary,
i.e.\ $\rho_f\vert_{G_{\qq_p}}$ preserves a Borel. The representation $\rho_f\vert_{G_{\qq_p}}$ is de Rham  if it fits  in an extension:
$$ 0 \rightarrow \rho_f^{(1 + \frac{w}{2})}\vert_{G_{\qq_p}}\rightarrow \rho_f\vert_{G_{\qq_p}} \rightarrow \rho^{(2+\frac{w}{2})}_f\vert_{G_{\qq_p}} \rightarrow 0$$
where 
$\rho_f^{(1 + \frac{w}{2})}\vert_{G_{\qq_p}}(1+\frac{w}{2})$ and
$\rho^{(2+\frac{w}{2})}_f\vert_{G_{\qq_p}}(2+\frac{w}{2})$ are potentially
unramified $2$-dimensional representations. Equivalently, this means that the
Sen operator of $D_{Sen}(\rho_f\vert_{G_{\qq_p}})$ is semi-simple with
eigenvalues $-1-\frac{w}{2}$ and $-2-\frac{w}{2}$. %

\begin{thm}\label{thm:multiplicity-one-implies-classical}  Let $f \in \HH^0_{^0w}( \Sh_{K^p}, \omega^{(2,2;-w),\sm})^{\ord}$ be an ordinary overconvergent modular eigenform with Galois representation
$\rho_f: G_{\qq} \rightarrow \mathrm{GSp}_4(\Qpbar)$.  Let $\mathfrak{m}_f$
be the corresponding maximal ideal. Let $V = \HH^3(\mathrm{RHom}_{\mathfrak{b}} ( \lambda, \mathrm{R}\Gamma(\Sh_{K^p}^{\tor}, \Qpbar)^{\la})^{\ord}_{\mathfrak{m}_{\rhobar_f}})$. We assume that:  
\begin{enumerate}
\item Either
\begin{enumerate}[(a)]
\item the Zariski closure of $\rho_f(G_{\qq})$ contains $\mathrm{Sp}_4$; or
\item the Zariski closure of $\rho_f(G_{\Q})$
    contains~$\SL_2 \times \SL_2 $, and $\rho_f$ is irreducible but becomes reducible
  on some index two subgroup~$G_{E}$.
\end{enumerate}
\item The representation $\rho_f\vert_{G_{\qq_p}}$ is de Rham.
\item\label{ass:dimension-graded-add-up} There exists an integer $n$ such that $\dim_{\C_p}
  \mathrm{Gr}^iV_{\C_p}[\mathfrak{m}_f]  = n $  for each~$0\le i\le 3$, and $\dim_{\Qpbar} V[\mathfrak{m}_f] = 4n$. 
  \item The representation $\rhobar_f$ is irreducible.
\end{enumerate}
Then $f$ is a classical modular form. 
\end{thm}
\begin{rem}\label{rem:mult-1-suffices-for-classicality} The simplest way to verify the assumption $(3)$ is to
  prove the multiplicity one statement that  $\mathrm{dim}_{\C_p}
  \HH^i_{^iw}(\Sh^{\tor}_{K^p}, \omega^{(2,2;-w),\sm})[\mathfrak{m}_f]
  =1$ for each~$i$ (when this statement is true). This forces $\dim_{\Qpbar} V[\mathfrak{m}_f] = 4$. Indeed, $\dim_{\Qpbar} V[\mathfrak{m}_f]$ is a multiple of $4$ by Lemma \ref{lem-Galois-rep-cc}. 
\end{rem}

\begin{proof}By Proposition~\ref{prop-COusin}
  and~\eqref{eqn:Theta-ses}, we have an exact sequence (where $\Theta$ is the arithmetic Sen
  operator): %
  $$ 0 \rightarrow  \HH^3_{^3w}(\Sh^{\tor}_{K^p}, \omega^{ (1,1;-w),\sm})^{\ord}_{ \mathfrak{m}_{\rhobar_f}} \rightarrow V_{\C_p}[(\Theta +1+ \frac{w}{2})^2] \rightarrow 
 \HH^2_{^2w}(\Sh^{\tor}_{K^p}, \omega^{(1,1;-w),\sm})^{\ord}_{ \mathfrak{m}_{\rhobar_f}}\rightarrow 0 $$
 Since $\Theta +1+ \frac{w}{2} $ is nilpotent, it induces a  map: %
 $$\HH^2_{^2w}(\Sh^{\tor}_{K^p}, \omega^{ (1,1;-w),\sm})^{\ord}_{ \mathfrak{m}_{\rhobar_f}} \stackrel{\Theta +1+ \frac{w}{2}}\rightarrow \HH^3_{^3w}(\Sh^{\tor}_{K^p}, \omega^{ (1,1;-w),\sm})^{\ord}_{ \mathfrak{m}_{\rhobar_f}}.$$
We pass to $\mathfrak{m}_{f}$-isotypic components. 
Assumption~\eqref{ass:dimension-graded-add-up} implies that the sequence:
$$ 0 \rightarrow  \HH^3_{^3w}(\Sh^{\tor}_{K^p}, \omega^{ (1,1;-w),\sm})^{\ord}[\mathfrak{m}_f] \rightarrow V_{\C_p}[(\Theta +1+ \frac{w}{2})^2][\mathfrak{m}_f]  \rightarrow 
 \HH^2_{^2w}(\Sh^{\tor}_{K^p}, \omega^{(1,1;-w),\sm})^{\ord}[\mathfrak{m}_f] \rightarrow 0 $$ remains exact. %
 Therefore, we get a map 
 $$\HH^2_{^2w}(\Sh^{\tor}_{K^p}, \omega^{ (1,1;-w),\sm})^{\ord}[\mathfrak{m}_f] \stackrel{\Theta +1+ \frac{w}{2}}\rightarrow \HH^3_{^3w}(\Sh^{\tor}_{K^p}, \omega^{ (1,1;-w),\sm})^{\ord}[\mathfrak{m}_f].$$
 The kernel of this map is the space of classical forms $\HH^2(\Sh^{\tor}_{K^p},
 \omega^{ (1,1;-w),\sm})[\mathfrak{m}_f]$ because of Theorem
 \ref{thm-Sen-Cousin} and Proposition \ref{prop-COusin}. %
 The map is zero since   $V_{\C_p}[\mathfrak{m}_f] =
 D_{Sen}(\rho_f)^{\oplus n}$ by Lemma \ref{lem-Galois-rep-cc}, and  $\Theta$
 identifies with the Sen operator of $\rho_f$ which is semi-simple by assumption. 
 It follows that $\HH^2(\Sh^{\tor}_{K^p}, \omega^{
   (1,1;-w),\sm})^{\ord})[\mathfrak{m}_f]  = \HH^2_{^2w}(\Sh^{\tor}_{K^p},
 \omega^{ (1,1;-w),\sm})^{\ord}[\mathfrak{m}_f]$. This implies that $\dim
 \HH^0(\Sh^{\tor}_{K^p}, \omega^{ (2,2;-w),\sm}))[\mathfrak{m}_f] = n$ (by the
 stability of the  $L$-packet of automorphic forms corresponding to $f$) and
 therefore $\HH^0(\Sh^{\tor}_{K^p}, \omega^{ (2,2;-w),\sm}))[\mathfrak{m}_f]  =
 \HH^0_{^0w}(\Sh^{\tor}_{K^p}, \omega^{ (2,2;-w),\sm})^{\ord}[\mathfrak{m}_f]$,
 and we are done. %
 \end{proof}
 
 \begin{rem} The way we have explained the argument, we naturally proved the
   classicality of the relevant space of  degree $2$ cohomology classes, and used
   Arthur's classification of automorphic forms to deduce the classicality of
   $f$. This may seem a bit strange. We should first explain  why we focused on
   degree $2$ and degree $3$ cohomology classes in the argument, instead of
   degree $0$ and degree $1$ cohomology classes. The reason is that  the
   localization $\Loc(M(\mathfrak{g})_{\lambda})$ is simpler on $C_{^3w,Q}$,  but
   has some richer structure on its complement. This rich structure is
   irrelevant for the ordinary case, but would cause minor technical problems in
   Section \ref{sectsencous}.
   
Still focusing on degree $2$ and $3$ cohomology as we did, we could also have
proven the classicality of the relevant space of degree $3$ cohomology, and then used
Serre duality instead of Arthur's classification to deduce the classicality of
$f$. In order to do this, one would consider $V/\mathfrak{m}_f V$ instead of $V[\mathfrak{m}_f]$, and the same argument would  apply with minor modifications. 
\end{rem}

\section{An ordinary modularity lifting theorem for unitary groups with $p=2$}\label{sec:
  R=T unitary}%

Our goal in this section is to prove Theorem~\ref{thm: U(n) ordinary
  automorphy lifting including Ihara avoidance}, which combines an
ordinary $2$-adic automorphy lifting theorem with a finiteness theorem
for a universal deformation ring. This result is a slight variant on
the $2$-adic automorphy lifting theorems for unitary groups proved by
Thorne in~\cite{MR3598803}; we work with ordinary representations and
use Ihara avoidance, and we use a slighter weaker definition of
adequacy (see Definition~\ref{defn: nearly adequate}). We emphasize
that there are no significant innovations here, and indeed the
arguments of~\cite{MR3598803} go through verbatim using this weaker
notion of adequacy. The one minor novelty in our arguments is an
argument using base change to allow us to use an auxiliary prime 
in order to work at neat level; the usual choices of such primes for
$p>2$ rely on choosing a prime at which all Galois deformations are
unramified, which is impossible when~$p=2$.

We have endeavored to write out the arguments in enough detail to
make it easy for the reader familiar with automorphy lifting theorems
for unitary groups which assume that $p>2$ (e.g.\ \cite{jack, BLGGT})
but unfamiliar with~\cite{MR3598803} to check the details (although
where these details are literally identical to those
of~\cite{MR3598803} we do not repeat the proofs). We follow the
notation of~\cite{MR3598803} closely, although we assume throughout
that~$p=2$, as the analogue for $p>2$ of our results is already known
(see e.g.\ \cite[Cor.\ 7.3]{MR3598803}), and in any case the only use
of automorphy lifting theorems for unitary groups in this paper is in
the case~$p=2$.

In~\S\ref{defofcG}, we recall the notion
of polarized representations
and clarify the relationship between
essentially self-dual representations~$G_{F^{+}}
\rightarrow \GSp_{2n}(R)$ and their associated
polarized representations~$G_{F} \rightarrow \cG_{2n}(R)$.
In~\S\ref{subsec: oddness in Un}, we recall the notion
of strong residual oddness defined in~\cite{MR3598803} and
establish some
basic facts concerning what this entails
for polarized representations
associated to essentially self-dual
representations over totally real fields with image in~$\GSp_{2n}(k)$.
(In practice, we only use the special
case corresponding to~$\GSp_4(\F_2)$.)
In~\S\ref{subsec: nearly adequate},  we discuss variants of 
the notion of adequateness in characteristic~$2$
as introduced in~\cite{MR3598803}.
Finally, in sections~\S\ref{subsec: Galois
  deformations Un}, \S\ref{subsec: TW systems}, \S\ref{subsec: local
  deformation problems Un},  and~\S\ref{subsec: auto forms Un}, we adapt the
arguments of~\cite{MR3598803} to our precise setting.

\subsection{Polarized representations}\label{subsec: polarized representations}     \label{defofcG}

       Let $\cG_n$ denote the semi-direct product of $\cG_n^0=\GL_n \times \GL_1$ by the group $\{1 , \jmath\}$
where 
\[ \jmath (g,a) \jmath^{-1}=(ag^{-t},a). \]
We let $\nu:\cG_n \to \GL_1$ be the character which sends $(g,a)$ to
$a$ and sends $\jmath$ to $-1$. (This group, in the context of modularity lifting,
was first introduced in~\cite[\S2.1]{cht}.)

Let $F$ be an imaginary CM field with maximal totally real
   subfield~$F^+$.  
 For each infinite place~$v$ we let~$c_v \in G_{F^{+}}$ denote complex
 conjugation. Let~$c$ denote a fixed arbitrary choice of element in~$G_{F^{+}} \setminus G_{F}$  with~$c^2=e$
   (so for example one could take~$c = c_v$ for any~$v|\infty$).
   \begin{defn}
     \label{defn: polarized}Let  $k$ be a perfect field,
     $\rho:G_F\to\GL_n(k)$ an absolutely irreducible representation,
     and $\mu:G_{F^+}\to k^\times$ a character. We say that
     $(\rho,\mu)$ is \emph{polarized} if there is a perfect pairing %
     $\langle\cdot,\cdot\rangle:k^n\times k^n\to k$ such
     that \[\langle x,y\rangle=-\mu(c)\langle y,x\rangle,\] 
    and for all $g\in G_F$ we have
      \[\langle \rho(g)x,\rho(c g c^{-1})y\rangle=\mu(g)\langle x,y\rangle.\]
   \end{defn}
 We have (\cite[Lem.\ 2.2]{MR3598803}):
   \begin{lemma} \label{extendtoG}
   If $(\rho,\mu)$ is polarized, then  %
   we may extend~$\rho$ to $r:G_{F^+}\to\cG_n(k)$ with $\nu\circ
   r=\mu$ and $r^{-1}(\cG_n^0(k))=G_F$; and this extension is unique
   up to $\cG_n^0(k)$-conjugacy. If~$R$ is any ring, then
   $\cG_n^0(R)$ acts via conjugation on the set of homomorphisms
   $r:G_{F^+}\to\cG_n(R)$ with $r^{-1}(\cG_n^0(k))=G_F$.
   \end{lemma}

   \begin{rem}\label{rem: cht conjugacy versus jt}
     Note that the~$\GL_1(R)$ factor of~$\cG_n^0(R)$ is not in general
     acting trivially via conjugation (while it is in the centre
     of~$\cG_n^0(R)$, it is not in the centre of~$\cG_n(R)$). Lemmas
     2.1--2.5 of~\cite{MR3598803} are analogues for $\cG_n^0$-conjugacy of lemmas in~\cite[\S
     2.1]{cht} which work instead with~$\GL_n$-conjugacy. As well as
     giving cleaner statements (for example, the extension of~$r$
     to~$\rho$ above is unique up to $\cG_n^0(k)$-conjugacy, but the
     $\GL_n(k)$-conjugacy classes of~$\rho$ are in bijection with
     $k^\times/(k^\times)^2$), the versions of~\cite{MR3598803} hold
     in the case of residue characteristic~$2$, unlike their analogues
     in~\cite{cht}.
   \end{rem}

In Section~\ref{subsec: deducing R=T GSp4 from base change}, we will
need to relate representations $G_{F^+}\to\GSp_4(R)$ to
representations $G_{F^+}\to \cG_4(R)$. We now explain how to do this
following the construction of~\cite[Lem.\ 2.1]{cht}. Write~$\GSp_{2n}$
for the generalized symplectic group defined by an antisymmetric
matrix~$J_{2n}$ (in particular, we can take~$J_4$ to be the matrix~$J$
that we use to define~$\GSp_4$).

\begin{defn}\label{defn: multiplier extended to GSp times pm 1}
  If~$R$ is any ring, then we extend the multiplier
  $\nu:\GSp_{2n}(R)\to R^{\times}$ to a
  homomorphism \[\nu:\GSp_{2n}(R)\times\{\pm 1\}\to R^{\times}\] via
  projection to the~$\GSp_{2n}$-factor, i.e.\ we set
  $\nu(g,a):=\nu(g)$.
\end{defn}
\begin{lem} %
  \label{lem: symplectic to cG}%
There is an injective
  homomorphism \[r:\GSp_{2n}(R)\times
    \{\pm 1\}\to\cG_{2n}(R)\] defined as follows:
  \begin{enumerate}
\item $r((g,1)) = (g,\nu(g))$.
\item $r((g,-1)) = (g,\nu(g))\cdot(J_{2n}^{-1},-1) \jmath$.
\end{enumerate}
This homomorphism is compatible with~$\nu$ \emph{(}defined on the
source in Definition~\ref{defn: multiplier extended to GSp times pm
  1}\emph{)}.
\end{lem}
\begin{proof} The only non-trivial thing to check is that
$r((g,-1))r((h,-1)) = r((gh,1))$, which amounts to the claim
that \[(J_{2n}^{-1},-1) \jmath (h,\nu(h))(J_{2n}^{-1},-1) \jmath=(h,\nu(h)).\]
The left hand side is (noting~$\jmath = \jmath^{-1}$):
\begin{align*}
(J_{2n}^{-1},-1) \jmath   (h ,\nu(h)) (J_{2n}^{-1},-1) \jmath & = (J_{2n}^{-1},-1) \jmath   (hJ_{2n}^{-1} ,-\nu(h)) \jmath^{-1}\\
   &= (J_{2n}^{-1},-1)    (-\nu(h)(hJ_{2n}^{-1})^{-t} ,-\nu(h)) 
 \\ &= (-\nu(h) J_{2n}^{-1}(hJ_{2n}^{-1})^{-t},\nu(h)),
\end{align*}
so we need to show that $-\nu(h)
J_{2n}^{-1}(hJ_{2n}^{-1})^{-t}=h$. This can be rearranged
to \[h^tJ_{2n}h=-\nu(h)J_{2n}^t,\] and since $J_{2n}^t=-J_{2n}$ and~$h\in\GSp_{2n}(R)$, we are done.
\end{proof}

\begin{cor}
  \label{cor: another symplectic to unitary} If
  $\aaaa:G_{F^{+}}\to\GSp_{2n}(R)$ is a homomorphism, then there is a
  homomorphism~ $r_{\aaaa}:G_{F^+}\to\cG_{2n}(R)$ defined by \[r_{\aaaa}(g)=
    \begin{cases}
    ({\aaaa}(g),\nu\circ {\aaaa}(g))  & \text{if } g\in G_F\\
    ({\aaaa}(g),\nu\circ {\aaaa}(g))\cdot(J_{2n}^{-1},-1)\jmath    & \text{if } g\in G_{F^+}\setminus G_{F}.
    \end{cases}
  \]   Furthermore we have $r_{\aaaa}^{-1}(\cG_{2n}^{\circ}(R))=G_F$, and $\nu\circ r_{\aaaa}=\nu\circ {\aaaa}$.
\end{cor}
\begin{proof}There is obviously a homomorphism
  \[\taaaa:G_{F^+}\to \GSp_{2n}(R)\times\{\pm 1\} \] given by the product
  of~$\psi$ and the
  projection onto $G_{F^{+}}/G_F\cong\{\pm 1\}$, i.e.\
  \numequation\label{eqn: tilde a}\taaaa(g)=
    \begin{cases}
      ({\aaaa}(g),1)  & \text{if } g\in G_F\\
      ({\aaaa}(g),-1)    & \text{if } g\in G_{F^+}\setminus G_{F}.
    \end{cases}
  \end{equation} Furthermore we have $\nu\circ\taaaa=\nu\circ {\aaaa}$. By definition
  we have $r_{\aaaa}=r\circ \taaaa$, where~$r$ is as in Lemma~\ref{lem:
    symplectic to cG}, and the result follows immediately.  
\end{proof}
\begin{rem}
  \label{rem: how did we uniquely extend a symplectic representation
    to cG?}If~$R$ is a (perfect) field, then we may apply Lemma~\ref{extendtoG} to the
  representation~${\aaaa}|_{G_{F}}$ of  Corollary~\ref{cor: another symplectic
    to unitary}, and we see that the extension of~${\aaaa}|_{G_{F}}$ to a
  homomorphism $G_{F^+}\to\cG_{2n}(R)$ is 
  well-defined up to~$\cG_{2n}^{\circ}(R)$-conjugacy. Corollary~\ref{cor: another symplectic
    to unitary} provides a particular choice of extension (depending
  of course on our choice of symplectic form, i.e.\ on ~$J_{2n}$). %
  The reader may find it helpful to compare to the discussion at the
  end of~\cite[\S 1.1]{BLGGT}, which in the case that 
  $R$ is a field and ${ \aaaa}|_{G_F}$ is absolutely irreducible
  shows that the choice of a specific element in the
  ~$\cG_{2n}^{\circ}(R)$-conjugacy class %
  amounts to
  choosing~$\bbbb_c\in\GL_{2n}(R)$ with \[{\aaaa}(cgc^{-1})\cdot \bbbb_c
    {\aaaa}(g)^t=\nu(g)\bbbb_c\] for all~$g\in G_F$. The implicit choice of such
  an element in Corollary~\ref{cor: another symplectic
    to unitary} is~$\bbbb_c:=\aaaa(c)J_{2n}^{-1}$. 
\end{rem}

\subsection{Oddness}\label{subsec: oddness in Un}
We recall the following definition~\cite[Defn.\ 3.3]{MR3598803}.
      \begin{defn}
     \label{defn: strongly residually odd}Suppose that $(\rho,\mu)$ is
     polarized,  that~$k$ has characteristic~$2$, and that $n$ is
     even. If~$v$ is an infinite place of~$F^+$, then we say that
     $(\rho,\mu)$ is \emph{strongly residually odd at~$v$} if $r(c_v)$
     is $\GL_n(k)$-conjugate to $(1_n,1)\jmath$.     
   \end{defn}

    \begin{remark}%
The idea behind Definition~\ref{defn: strongly residually odd} is as follows. A representation
      $r:G_{F^+}\to\cG_{n}(\Qpbar)$ is defined to be (totally) odd if~$\nu\circ r(c_v)=-1$
      for all infinite places of~$F^+$ (see e.g.\ \cite[\S
      2.1]{BLGGT}). If~$p>2$, then a representation
      $\rbar:G_{F^+}\to\cG_{n}(\Fpbar)$ is odd if ~$\nu\circ \rbar(c_v)=-1$, and
      any lifting of~$\rbar$ will automatically be odd. Furthermore, in
      either case there is a single $\GL_{n}$-conjugacy class of
      elements~$(x,1)\jmath$ of order~$2$ (because all symmetric matrices
      are equivalent), so the analogue of
      Definition~\ref{defn: strongly residually odd} holds
      automatically.
\      If~$p=2$, in contrast, then the condition that $\nu\circ \rbar(c_v)=-1$ is
      automatic. However, if $n$ is even, then  there are
      two conjugacy classes of elements of the form~$(x,1)\jmath$ of order~$2$ (see~\cite[Lem.\
      2.16]{MR3598803}). In the situation of Definition~\ref{defn:
        strongly residually odd}, any lift of~$(\rho,\mu)$ 
      is automatically
      odd at~$v$, i.e.\ the lift of~$\mu$ is odd; this is one
      motivation for the terminology ``strongly residually
      odd''.
\end{remark}

In the remainder of this section, we examine when the
representations~$r_{\aaaa}$ of Corollary~\ref{cor: another symplectic to
  unitary} are strongly residually odd at some~$v$. We will ultimately
only need a single example for~$\GSp_4$, but as it is straightforward
to give a general treatment, we do so.
Assume for the rest of this subsection that~$k$ has
characteristic~$2$, so that in particular~$J_{2n}=-J_{2n}=J_{2n}^{-1}$. 
 Let~$G = \GL_{2n}(k) \times k^{\times}$, and let~$G.2 = G \rtimes  \Z/2\Z = \cG_{2n}(k)$ where the action by
the order~$2$ element~$\jmath$ is by the outer
 automorphism~$\jmath (g,a) \jmath^{-1}=(ag^{-t},a)$. Let~$\chi: G.2 \rightarrow \Z/2\Z$ be the canonical projection. By Lemma~\ref{lem: symplectic to cG}, 
there is a natural map~$\Sp_{2n}(k) \times \Z/2\Z \rightarrow G.2$ given by sending~$-1$
to~$J_{2n}\jmath$. %
If~$A \in \Sp_{2n}(k)$ satisfies~$A^2 = I$,
then~$(AJ_{2n},1) \cdot \jmath \in G.2 \setminus G$ has order~$2$.
Recall (e.g.\ from the proof of~\cite[Lem.\ 2.16]{MR3598803}) that any
invertible symmetric matrix~$M$ in~$\GL_{2n}(k)$ is either equivalent
(under $M\mapsto g^tMg$) to~$J_{2n}$ or to~$I$;
the former if and only if all diagonal entries of~$M$ are zero.
In the former case we say that~$M$ is alternating
 (since the corresponding non-degenerate pairing satisfies~$B(x,x)=0$ for all~$x$)
and otherwise we say that~$M$ is non-alternating.

\begin{lemma} \label{oddity} If~$A \in \Sp_{2n}(k)$ satisfies~$A^2 = I$, then~$AJ_{2n}$ is a symmetric matrix and~$AJ_{2n} \cdot \jmath \in G.2$
has order~$2$. This induces a map
$$
\begin{aligned}
 \left\{\begin{array}{c}
\text{Conjugacy classes of } \Sp_{2n}(k) \\
\text{of order dividing~$2$}
\end{array}\right\} &
 \to  \left\{\begin{array}{c}
\text{Conjugacy classes of } G.2 \setminus G \\
\text{of order~$2$}
\end{array}\right\}
\end{aligned}
$$
The target has order~$2$ and consists of
the conjugacy classes of~$\jmath$ and~$J_{2n} \cdot \jmath$.
The fibre over~$J_{2n} \cdot \jmath$ consists of~$A$ for which~$AJ_{2n}$
is alternating, and the fibre over~$\jmath$ consists of~$A$ for which~$AJ_{2n}$
is not alternating. The set on the right remains unchanged if we only
consider order~$2$ elements in~$G.2 \setminus G$ up to conjugation
by~$\GL_{2n}(k) \subset G \subset G.2$.
\end{lemma}

\begin{proof} For any element~$\gamma \in G.2 \setminus G$, the~$G.2$-conjugacy
class of~$\gamma$ coincides with the~$G$-conjugacy class of~$\gamma$, since
any element in~$G.2$ either has the form~$g \in G$ or~$g \gamma$ with~$g \in G$,
and~$(g \gamma) \gamma (g \gamma)^{-1} = g \gamma g^{-1}$.
If~$\gamma = (B,b) \jmath$ and~$g = (I_{2n},\lambda) \in G$, then
$$g \gamma g^{-1} =   (I_{2n},\lambda) (B,b) \jmath  (I_{2n},\lambda^{-1}) = (\lambda^{-1} \cdot B, b) \jmath,$$
but writing~$\lambda^{-1} = \mu^{2}$ (possible since~$k$ is finite of characteristic~$2$),
and taking~$h = (\mu I_{2n},1) \in \GL_{2n}(k) \subset G$, 
we also have
$$h \gamma h^{-1} = (\mu I_{2n},1)  (B,b) \jmath  (\mu^{-1} I_{2n},1)
=
 (\mu^2 \cdot B,b) \jmath = (\lambda^{-1} \cdot B, b) \jmath,$$
 and so the~$G$ and~$\GL_{2n}(k)$-conjugacy classes of~$\gamma \in G.2 \setminus G$ also coincide.

 If~$A^2 = I$, then~$(A,-1) \in \Sp_{2n}(k) \times \Z/2\Z$ has order~$2$ so the image~$AJ_{2n} \cdot \jmath
\in G.2$ certainly has order~$2$ and does not lie in~$G$. Moreover, this map certainly induces a map on conjugacy
classes because if~$A$ is conjugate to~$A'$ in~$\Sp_{2n}(k)$, then~$(A,-1)$ is conjugate to~$(A',-1)$
in~$\Sp_{2n}(k) \times \Z/2\Z$. The condition that~$(B,b) \cdot \jmath \in G.2$ has order~$2$ is equivalent to the equation
$$(B,b) \jmath (B,b) \jmath = (b \cdot B(B^t)^{-1},b^2) = (I,1),$$
which implies that~$b^2=1$ and so~$b = 1$ and~$B = B^t$ is symmetric. Let~$\Sigma_{2n} \subset \GL_{2n}(k)$ denote the set of symmetric matrices.
Conjugation by~$(M,1) \in \GL_{2n}(k) \subset G$ replaces~$(B,1)\jmath$ 
by~$(M B M^t,1) \jmath$.
Hence the order~$2$ conjugacy classes in~$G.2\setminus G$ are the orbits of~$\GL_{2n}(k)$
acting via $MBM^t$ on~$\Sigma_{2n}$.
But the orbits of~$\GL_{2n}(k)$ on this space are none other than the equivalence class
of perfect pairings on~$k^{2n}$, and as recalled above, there are two such orbits, corresponding to~$B=I$ and~$B=J_{2n}$.
\end{proof}

\subsubsection{Involutions in~$\Sp_{2n}(k)$}\label{subsubsec:involutions-in-Sp2n} (See the proof
of~\cite[Lem.\ 4.3]{bobthebuilder}).
An involution~$A$ in~$\Sp_{2n}(k)$ acting
on the natural representation~$V$ preserves a flag
\numequation\label{eqn: involution flag}0 \subset (A-1)V  \subset   V^A \subset V,\end{equation}
where we write~$r:=\dim((A-1)V)  \le n$ for the rank of~$A-1$. With respect to this flag, one can
(by~\cite[Lem.\ 4.4]{bobthebuilder}) write~$A$ in the form
$$A = \left( \begin{matrix} I_r & 0 & S_r\\ 0 & I_{2n-2r} & 0 \\ 0 & 0 & I_r \end{matrix} \right),$$
where~$S_r$ has rank~$r$ and~$S_r J_r$ is symmetric. The parabolic
stabilizing the flag~\eqref{eqn: involution flag} acts on the matrices
                                              of this form, and the
                                               orbit corresponds to 
                                               all symmetric matrices
                                               equivalent to~$S_rJ_r$.
                                               Accordingly, if~$r$ is odd, the conjugacy class of~$A$ is determined by~$r$.
If~$r > 0$ is even, there are two conjugacy classes corresponding to~$S_r=I_r$ and~$S_r=J_r$.
In total  there are~$n+1+ \lfloor n/2 \rfloor$ conjugacy classes.
We see that
\numequation\label{eqn: AJ2n}AJ_{2n} =   \left( \begin{matrix} S_r J_r & 0 & J_r \\ 0 & J_{2n-2r} & 0 \\ J_r & 0 & 0 \end{matrix} \right)\end{equation}
This is non-alternating  if  either~$r$ is odd or~$r$ is even and~$S_r = J_r$.
Any involution~$A \in \GSp_{2n}(k)$ must have~$\nu(A)^2 = 1$ and thus~$\nu(A)=1$,
and hence must lie in~$\Sp_{2n}(k)$.
In particular,  from the discussion above and Lemma~\ref{oddity}, we have the following:

\begin{lemma} \label{checkodd} Let~$\aaaabar: G_{F^+} \rightarrow
  \GSp_{2n}(k)$ with~$k$ of characteristic~$2$, and
  let~$r_{\aaaabar}:G_{F^+}\to\cG_{2n}(k)$  be as in Lemma~\ref{cor:
    another symplectic to unitary}. Let~$v$ be an infinite place
  of~$F^+$. Then the polarized pair~$(\aaaabar|_{G_F},\nu\circ\aaaabar)$ is strongly residually odd at~$v$ if and only if either
\begin{enumerate}
\item  $(\aaaabar(c_v)-I)$ has odd rank, or
\item $(\aaaabar(c_v)-I)$ has even rank~$r > 0$,
and the matrix~$S_r J_r$ obtained from~\eqref{eqn: AJ2n} with~$A:=\aaaabar(c_v)$ is non-alternating.
\end{enumerate}
Equivalently~$(\aaaabar|_{G_F},\nu\circ\aaaabar)$ is strongly residually odd at~$v$ if and only if 
the quadratic form associated~$\aaaabar(c_v) J_{2n}$ is equivalent to the one associated with~$I_{2n}$, which occurs if and only if~$\aaaabar(c_v) J_{2n}$ has at least one non-zero diagonal entry.
\end{lemma}

\begin{rem}\label{rem: GL2 oddness is particularly simple}                                               
If~$n=1$, then either~$A$ is trivial or~$A-I$ has rank~$r=1$. Hence 
in this case strong residual oddness is equivalent to~$A \ne I$
(cf.~\cite[Lem.\ 3.5(ii)]{MR3598803}).
This is no longer true for~$n>1$;
 there are~$n$ conjugacy classes of involutions
giving rise to strongly residually odd representations and~$1 + \lfloor n/2 \rfloor$ conjugacy classes
of involutions which do not. In particular, for~$2n = 4$, there are two conjugacy
classes of involutions giving rise to strongly residually odd representations and two
(one of which is the identity) which do not. 
Explicit representatives for the odd
classes can be given as follows:
$$\left(\begin{matrix} 1 & 0 & 0 & 1 \\ 0 & 1 & 0 & 0 \\ 0 & 0 & 1 & 0 \\ 0 & 0 & 0 & 1 \end{matrix}\right), \quad
\left( \begin{matrix} 1 & 0 & 0 & 1 \\ 0 & 1 & 1 & 0 \\ 0 & 0 & 1 & 0 \\ 0 & 0 & 0 & 1 \end{matrix} \right),$$
where the latter is conjugate to~$J_4$,
and an explicit representative for the non-trivial non-odd involution is given by
$$\left( \begin{matrix} 1 & 1 & 0 & 0 \\ 0 & 1 & 0 & 0 \\ 0 & 0 & 1 & 1 \\ 0 & 0 & 0 & 1 \end{matrix} \right).$$
When we later fix (in Lemma~\ref{explicit}) an explicit isomorphism~$S_6 \simeq \Sp_4(\F_2)$,
the first two elements can be identified with the images of~$(1,2)$, and~$(1,2)(3,4)$ respectively, whereas the latter can be identified with~$(1, 2)(3, 5)(4, 6)$
(See also the proof of Lemma~\ref{toquotetwo}.)
\end{rem}

   \subsection{Adequacy}\label{subsec: nearly adequate}

Let~$\ad = \Hom(V,V)$ denote the adjoint representation and~$\ad^0 \subset \Hom(V,V)$
the submodule of trace zero endomorphisms.
We begin with the following lemma just to clarify that the definition
of weakly adequate used in~\cite[Defn.\ 2.20]{MR3598803}
(which is case~(\ref{jackcase}) of Lemma~\ref{itsclear} below)
agrees with other definitions in the literature. 

\begin{lemma} \label{itsclear}Let~$V$ be a finite-dimensional vector space over a finite field~$k$,
and let~$H \subseteq \GL(V)$.
The following conditions are equivalent:
\begin{enumerate}
\item  \label{jackcase}  For each simple~$\kbar[H]$-submodule~$W \subset \ad \otimes \kbar$, there exists a semi-simple 
element~$\sigma \in  H$ with an eigenvalue~$\alpha \in \kbar$ such that $\mathrm{tr}(e_{\sigma,\alpha} W)
\ne 0$. (Here $e_{g,\alpha}$ is the $g$-equivariant projection onto the generalized
$\alpha$-eigenspace of~$g$.)
\item  \label{jackcasetwo} For each simple~$\kbar[H]$-submodule~$W \subset \ad^0 \otimes \kbar$, there exists a semi-simple 
element~$\sigma \in  H$ with an eigenvalue~$\alpha \in \kbar$ such that $\mathrm{tr}(e_{\sigma,\alpha} W)
\ne 0$.
\item  \label{spansirreducible}  $\End(V)$
is spanned by the set~$H^{\semis}$ of semi-simple elements of~$H$. %
\end{enumerate}
\end{lemma}

\begin{proof}  This follows directly from the proof
  of~\cite[Lemma~A.1]{jack}. %
More precisely, it is shown there that we have an equality 
$$\begin{aligned} U: = & \  \{w \in \ad\otimes\kbar: \mathrm{tr}(g w) = 0 \quad \forall g \in H^{\mathrm{ss}}\} \\
= & \  \{w \in \ad\otimes\kbar: \mathrm{tr}(e_{g,\alpha} w) = 0 \quad \forall g \in H^{\mathrm{ss}}, \alpha \in \kbar\}.
\end{aligned}$$
Note that~$U$ is an~$H$-submodule of~$ \ad\otimes\kbar$; suppose that~$w \in U$,
$g \in H^{\mathrm{ss}}$, and~$h \in H$. The element~$h$ via the natural action sends~$w$
to~$h w h^{-1}$, and then
\[\mathrm{tr}(g h w h^{-1}) = \mathrm{tr}(h^{-1} g h w) = 0,\]
since~$h^{-1} g h \in H^{\mathrm{ss}}$ and~$w \in U$.
Condition~(\ref{spansirreducible}) is equivalent to~$U = 0$, whereas conditions~(\ref{jackcase}) and~(\ref{jackcasetwo})
 above
are equivalent to the intersection of $U$ with the socle of~$\ad\otimes\kbar$
(respectively, the socle of~$\ad^0\otimes\kbar$) being trivial. %
Since~$U \subset \ad^0\otimes\kbar$ (take~$g=1$), the result follows.
\end{proof}

   \begin{df}\label{defn: weakly adequate}We say that $H \subseteq \GL(V)$ is
     \emph{weakly adequate} if the equivalent conditions of
     Lemma~\ref{itsclear} hold for~$H$. A
     representation~$\rho: G \rightarrow \GL(V)$ is \emph{weakly
       adequate} if $\im(\rho)$ is weakly adequate.
   \end{df}
   As remarked in~\cite{MR3598803} (after Definition~2.20), if~$H$ is
   weakly adequate, then~$H$ acts absolutely irreducibly on~$V$ as a
   consequence of condition~\ref{spansirreducible} of
   Lemma~\ref{itsclear}. 

   \begin{defn}
     \label{defn: nearly adequate}Let $k$ be a subfield of~$\Ftwobar$. We say that a finite subgroup
     $H\subset\GL_n(k)$ is \emph{nearly adequate} if:
     \begin{enumerate}
     \item \label{item: weakly adequate} $H$ is weakly adequate.
     \item \label{item: H1 scalars vanishes} $H^1(H,k)=0$.
     \item \label{item: H1 gl vanishes} $H^1(H,\ad)=0$. %
   \end{enumerate}
       We say that a representation~$\rho: G \rightarrow \GL(V)$ is
     \emph{nearly adequate} if $\im(\rho)$ is nearly adequate.
   \end{defn}

   \begin{rem}
     \label{rem: adequate implies nearly adequate}The definition of a
     nearly adequate subgroup is almost the same as the definition of an
     adequate subgroup~\cite[Defn.\ 2.20]{MR3598803}. Indeed the only
     difference is that we are assuming that $H^1(H,\ad)=0$, rather
     than the stronger assumption that $H^1(H,\ad/k)=0$ (which implies
     the vanishing of $H^1(H,\ad)$ in conjunction with the assumption
     that $H^1(H,k)=0$).

     The point of our definition is that, as we explain below, the
     arguments of~\cite{MR3598803} apply unchanged with ``adequate''
     relaxed to ``nearly adequate'', and in our applications we will need to
     work with representations which are nearly adequate, but not
     adequate. (See Lemma~\ref{toquoteone} and Remark~\ref{rem: A5 is only nearly adequate}.)
   \end{rem}

\subsection{Galois deformation theory}\label{subsec: Galois
  deformations Un}We now recall some facts about Galois deformation
theory when $p=2$. The results we
need are essentially identical to those of~\cite[\S 2.1]{MR3598803},
except that we need to work relative to a larger coefficient ring (the
weight space~$\Lambda$), which we do following~\cite[\S 4]{KT}.

We continue to assume that~$F$ is an imaginary CM field with maximal totally real
subfield~$F^+$, and we assume that $F/F^+$ is everywhere unramified,
and that all places of~$F^+$ dividing~$2$ split in~$F$. 
We write~$S_2$ for the set of
places of~$F^+$ dividing~$2$, $S_\infty$ for the set of places
of~$F^+$ dividing~$\infty$, and~$S$ for a finite set of places
of~$F^+$ containing $S_2\cup S_\infty$. Let~$F(S)$ be the maximal
Galois extension of~$F$ unramified outside of~$S$, and write
$G_{F^+,S}:=\Gal(F(S)/F^+)$, $G_{F,S}:=\Gal(F(S)/F)$. Throughout this section we
will use the notation established in Section~\ref{notn}, specialized to the case
$p=2$, so that for example we have our field of coefficients~ $E/\Q_2$ with
ring of integers $\cO$, uniformizer $\varpi$, and residue field $k$. %

We fix a representation $\rbar:G_{F^+,S}\to \cG_n(k)$ such that
$\rbar^{-1}(\cG_n^0(k))=G_{F,S}$, together with a character
$\chi:G_{F^+,S}\to\cO^\times$ with $\chibar=\nu\circ\rbar$
and~$\chi(c_v)=-1$ for all~$v\in S_\infty$. We abusively write
$\rbar|_{G_{F,S}}$ for the representation given by restriction
to~$G_{F,S}$ and projection to the~ $\GL_n(k)$ factor in~$\cG_n^0(k)$.
We assume that $\rbar|_{G_{F,S}}$ is absolutely irreducible. We 
often write~$\rhobar$ for~$\rbar|_{G_{F,S}}$.

For each $v\in S$, we fix $\Lambda_v \in \CNL_{\cO}$, and set $\Lambda = \widehat{\otimes}_{v \in S} \Lambda_v$, where the completed 
tensor product is taken over $\cO$. 
For each~$v\in S$, the canonical map $\Lambda_v \rightarrow \Lambda$
induces the forgetful functor $\CNL_\Lambda \to \CNL_{\Lambda_v}$.

As in~\cite[Defn.\ 2.6]{MR3598803}, a \emph{lifting} of
$\rbar|_{G_{F^+_v}}$ to a $ \CNL_{\Lambda_v}$-algebra~$A$ is a
continuous homomorphism $r_v: G_{F^+_v} \rightarrow \cG_n(A)$ such
that $r_v \bmod \frakm_A = \rbar|_{G_{F^+_v}}$
and~$\nu\circ r_v=\chi|_{G_{F^+_v}}$. We let
$\Lift_v^\square:\CNL_{\Lambda_v}\to\Sets$ be the functor sending~$A$
to the set of liftings of~$\rbar|_{G_{F^+_v}}$. The functor
$\Lift_v^\square$ is representable by an
object~$R_v^\square\in\CNL_{\Lambda_v}$. If~$v$ splits in~$F$,
and~$r_v$ is a lifting of~$\rbar|_{G_{F^+_v}}$, then
$r_v(G_{F^+_v})\subseteq\cG_n^0(A)$, and we sometimes
write~$\rho_v:G_{F^{+}_v} \rightarrow \GL_n(A)$ for the projection
of~$r_v$ to the~$\GL_n$ factor.

 A \emph{local deformation problem} for $\rbar|_{G_{F^+_v}}$ is a
 representable subfunctor $\cD_v\subseteq\Lift_v^\square$ such that
 for all~$A\in\CNL_{\Lambda_v}$, the set $\cD_v(A)$ is invariant under
 the conjugation action of $\cGhat_n(A)$ on~$\Lift_v(A)$.%

 A \emph{global deformation problem} is a tuple %
  \[
   \cS = (F,\rbar, \cO, \chi,S, \{\Lambda_v\}_{v\in S}, \{\cD_v\}_{v\in S}),
  \]
 where:
 \begin{itemize}
  \item  $F$, $\rbar$,   $\cO$, $\chi$, $S$, and $\{\Lambda_v\}_{v\in S}$ are as above.
  \item For each $v\in S$, $\cD_v$ is a local deformation problem for $\rhobar|_{G_{F^+_v}}$.
 \end{itemize}

As in the local case, 
a \emph{lift} (or \emph{lifting}) of $\rbar$ is a continuous homomorphism $r: G_F \rightarrow \cG_n(A)$ to a $\CNL_\Lambda$-algebra $A$, 
such that $r \bmod \frakm_A = \rbar$. 
We say that two lifts $r_1,r_2: G_F \rightarrow \cG_n(A)$ are \emph{strictly equivalent} if there is an
$a\in \cGhat_n(A)$ such that $r_2 = ar_1 a^{-1}$. 
A \emph{deformation of $\rbar$} is a strict equivalence class of lifts of $\rbar$.

For a global deformation problem 
  \[
   \cS = (F, \rbar, \cO, \chi, S, \{\Lambda_v\}_{v\in S}, \{\cD_v\}_{v\in S}),
  \]
we say that a lift $r: G_F \rightarrow \cG_n(A)$ is of \emph{type $\cS$} if $r|_{G_{F_v}} \in \cD_v(A)$ for each $v\in S$. 
Note that if $r_1$ and $r_2$ are strictly equivalent lifts of $\rbar$, and $r_1$ is of type $\cS$, then so is $r_2$. 
A \emph{deformation of type $\cS$} is a strict equivalence class
of lifts of type $\cS$, and we denote by $\Def_{\cS}$ the set-valued functor that 
takes a $\CNL_\Lambda$-algebra $A$ to the set of deformations $r: G_F \rightarrow \cG_n(A)$ of type $\cS$.

Given a subset $T\subseteq S$, a \emph{$T$-framed lift of type $\cS$} is a 
tuple $(r,\{\alpha_v\}_{v\in T})$, where $r: G_F \rightarrow 
\cG_n(A)$ is a lift of $\rbar$ of type $\cS$ and $\alpha_v \in \cGhat_n(A)$ for each $v\in T$. 
We say that two $T$-framed lifts $(r_1,\{\alpha_v\}_{v\in T})$ and $(r_2,\{\beta_v\}_{v\in T})$ to a $\CNL_\Lambda$-algebra $A$ are strictly 
equivalent if there is an $a\in \cGhat_n(A)$ such that $r_2 = a r_1 a^{-1}$ and $\beta_v = a\alpha_v$ for each $v\in T$.
A strict equivalence class of $T$-framed lifts of type $\cS$ is called a \emph{$T$-framed deformation of type $\cS$}. 
We denote by $\Def_{\cS}^T$ the  functor  $\CNL_\Lambda\to\Sets$ taking~ $A$ to the set of $T$-framed deformations to $A$ of type $\cS$.

 Let $\cS = (F, \rbar, \cO, \chi, S, \{\Lambda_v\}_{v\in S}, \{\cD_v\}_{v\in S})$ be a global deformation problem, and let $T$ be a subset of $S$.
 The functors $\Def_{\cS}$ and $\Def_{\cS}^T$ are representable
 (see~\cite[Lem.\ 2.8, Lem.\ 2.10]{MR3598803}); we denote their
 representing objects by $R_{\cS}$ and $R_{\cS}^T$, respectively. If
 $T = \emptyset$, then tautologically $R_{\cS} = R_{\cS}^T$, while if $T$ is nonempty, $R_{\cS}^T$ is a formally smooth
$R_{\cS}$-algebra of  relative dimension $(n^2+1)\#T-1$.

Let $T$ be a (possibly empty) subset of $S$ such that $\Lambda_v =
\cO$ for all $v\in S \smallsetminus T$. Write~$R_v$ for the
representing object of~$\cD_v$,
and define $R_{\cS,T}^{\loc} = \widehat{\otimes}_{v\in T} R_v$, with the completed tensor product being taken over $\cO$. 
It is canonically a $\Lambda$-algebra, via the canonical isomorphism $\widehat{\otimes}_{v\in T} \Lambda_v \cong \widehat{\otimes}_{v\in S} \Lambda_v$. 
For each $v\in T$, the natural transformation $\Def_{\cS}^T \rightarrow \cD_v$ given by $(\rho,\{\alpha_v\}_{v\in T}) \mapsto \alpha_v^{-1}\rho|_{G_{F_v}} \alpha_v$ 
induces a  morphism $R_v \rightarrow R_{\cS}^T$ in $\CNL_{\Lambda_v}$. 
We thus have a morphism $R_{\cS,T}^{\loc} \rightarrow R_{\cS}^T$ in $\CNL_{\Lambda}$.  

In ~\cite[\S 2.2]{MR3598803}, the relative tangent space to this
morphism is computed via Galois cohomology. (Strictly speaking this
reference has $\Lambda_v=\cO$ for all~$v$, but the $\Lambda$-algebra
structure does not intervene in the calculation.) More precisely,
there is an explicit chain complex of~ $k$-vector spaces $C^i_{\cS,T}$ with cohomology
groups $H^i_{\cS,T}$ of $k$-dimensions~$h^i_{\cS,T}$, and by~\cite[Lem.\
2.12]{MR3598803}, we have $\dim_k
\m_{R_{\cS}^T}/(\m_{R_{\cS,T}^{\loc}},\m_{R_{\cS}^T}^2)=h^1_{\cS,T}$,
so that there is a surjection of $R_{\cS,T}^{\loc}$-algebras
$R_{\cS,T}^{\loc}\llb X_1,\dots,X_{h^1_{\cS,T}}\rrb \onto
R_{\cS}^T$. Since we do not need any properties of~$C^i_{\cS,T}$ and
its cohomology groups beyond those proved in~\cite[\S 2.2]{MR3598803},
we do not recall their definitions here. (It may however be helpful to
point out that there is a typo in the definition of~$C^i_{\cS,T}$
for~$i\ge 3$: the sum over places~$v\in S$ should be of
$C^{i-1}(F^+_v,\ad\rbar)$, not $C^{i-1}(F(S)/F^+,\ad\rbar)$ as written in~\cite{MR3598803}.) We do however need to
consider a certain dual Selmer group~$H^1_{\cS^\perp,T}$ of
$k$-dimension $h^1_{\cS^\perp,T}$, and we now turn to its definition
(see~\eqref{eqn: defn of dual Selmer group} below).

We identify ~$\ad\rbar$ with
$\widehat{\GL}_n(k[\epsilon])$, and we write~$\fg_n\rbar$ for the
adjoint representation on~$\widehat{\cG}_n(k[\epsilon])$, so that we
have an exact sequence
of~$G_{F^+,S}$-modules \numequation\label{eqn: ses ad and fgn}0\to\ad\rbar\to\fg_n\rbar\to k\to 0.\end{equation}

For each~$v\in S$, we as usual identify $\Lift_v^\square(k[\epsilon])$
with the cocycles $Z^1(F_v^+,\ad\rbar)$, so that two liftings to
$k[\epsilon]$ are $\cG_n(k[\epsilon])$-conjugate if and only if the
images of the corresponding cocycles in $H^1(F_v^+,\fg_n\rbar)$
coincide. (If~$v$ splits in~$F$, this is equivalent to their
images coinciding in $H^1(F_v^+,\ad\rbar)$.) We write
$\cL^1_v\subseteq Z^1(F_v^+,\ad\rbar)$ for the %
cocycles corresponding to liftings in $\cD_v(k[\epsilon])$, and
~$\cL_v$ for the image of~$\cL^1_v$ in $H^1(F_v^+,\ad\rbar)$. We write
$l_v^1,l_v$ for the dimensions of the $k$-vector spaces
$\cL_v^1,\cL_v$ respectively. We define
\[\mu_v:=\ker\left(H^1(F_v^+,\ad\rbar)\to
    H^1(F_v^+,\fg_n\rbar)\right).\] 
   From the long exact sequence in cohomology associated to~(\ref{eqn: ses ad and fgn}),
    there is a natural identification
\numequation\label{eqn: muv via boundary map}\mu_v \simeq \im \left(H^0(F_v^+,k) \to
    H^1(F_v^+,\ad\rbar)\right),\end{equation}
    and so~$\dim(\mu_v) \le 1$. 
     Note that since (by definition)
$\cD_v$ is stable under conjugation by~$\widehat{\cG}_n$, we have
$\mu_v\subseteq \cL_v$, and~$\mu_v$ is trivial if~$v$ splits in~$F$
(since in this case~\eqref{eqn: ses ad and fgn} splits as a sequence
of~$G_{F_v^+}$-modules). (Note however that the places in~$S_\infty$
do not split in~$F$, and indeed if~$(\rho,\mu)$ is strongly residually
odd at~$v\in S_\infty$ in the sense of Definition~\ref{defn: strongly
  residually odd}, then ~$\mu_v$ is 1-dimensional by~\cite[Lem.\
2.17(ii)]{MR3598803}.)

The trace pairing $(X,Y) \mapsto \tr(XY)$ on $\ad\rbar$ is perfect and $G_{F^+,S}$-equivariant, so $\ad\rbar(1)$ is isomorphic to the 
Tate dual of $\ad\rbar$. (Of course, since~$p=2$, the Tate twist is
trivial, and~$\ad\rbar$ is self-dual, but to avoid confusing the
reader who is used to the case~$p>2$, we preserve the Tate twist in
our notation below.) For each~$v\in S$ we let~$\cL_v^\perp\subseteq
H^1(F_v^+,\ad\rbar)$ be the annihilator of~$\cL_v$ under this pairing,
and write~$\mu_v^\perp\supseteq\cL_v^\perp$ for the annihilator
of~$\mu_v$.

For any $T\subseteq S$ as above (i.e.\
for any~$T$ such that $\Lambda_v = \cO$ for all
$v\in S \smallsetminus T$) we define the dual Selmer
group \nummultline\label{eqn: defn of dual Selmer group}H^1_{\cS^\perp,T}:=\\
  \ker\left(H^1(F(S)/F^+,\ad\rbar(1))\to\prod_{v\in
      T}H^1(F_v^+,\ad\rbar(1))/\mu_v^\perp\times \prod_{v\in S\setminus
      T}H^1(F_v^+,\ad\rbar(1))/\cL_v^\perp \right).\end{multline}

As usual, we write
$h^1_{\cS^\perp,T}:=\dim_kH^1_{\cS^\perp,T}$. Assume that~$T$ is
nonempty; then by~\cite[Lem.\ 2.15]{MR3598803} and our assumption
that~$\rbar|_{G_{F,S}}$ is absolutely irreducible, the~$h^i_{\cS,T}$
vanish for~$i\ne 1,2$, and we have
$h^2_{\cS,T}=h^1_{\cS^\perp,T}$. (The assumption that~$T$ is nonempty
guarantees the vanishing of~$h^0_{\cS,T}$, and the vanishing
of~$h^i_{\cS,T}$ for~$i\ge 4$ is automatic, as in the proof of~\cite[Lem.\ 2.13]{MR3598803}.) %

\begin{remark}[Remarks on~$\mu_v$ and our deformation problem.]
  We now try to explain where the terms~$\mu_v^\perp$ in~\eqref{eqn:
    defn of dual Selmer group} come from. A possibly unilluminating
  answer is that they are necessary in order to prove that
  $h^2_{\cS,T}=h^1_{\cS^\perp,T}$. Indeed, the proof ~\cite[Lem.\
  2.15]{MR3598803} is as usual via the Poitou--Tate sequence, and
  the~$\mu_v^\perp$ arise because of the appearance of the cohomology
  groups $H^1(F_v^+,\ad\rbar)^\eta$ in the long exact sequence
  in~\cite{MR3598803} computing the~$H^i_{\cS,T}$. Here
  $H^1(F_v^+,\ad\rbar)^\eta$ is by definition the image of the map \[H^1(F_v^+,\ad\rbar)\to
    H^1(F_v^+,\fg_n\rbar),\] and~$\mu_v$ is its kernel, whence the
  appearance of~$\mu_v^\perp$ as the dual Selmer condition.

  A possibly more helpful explanation is as follows.
  In the usual Kisin modification of the Taylor--Wiles
method, when one presents a global deformation ring over a (completed 
tensor product of) local deformation rings at primes~$v \in T$, the corresponding
Selmer condition~$\cL_v$ at~$v$ is trivial. In our setting
(exactly as in~\cite{MR3598803}) we are considering 
deformations~$\rho: G_{F^+_v} \rightarrow \cG_n(A)$ 
which are equivalent under~$\cGhat_n(A)$. As in Remark~\ref{rem: cht conjugacy versus jt},
conjugation by the~$\widehat{\GL}_1(A)$ factor of~$\cGhat(A)$
does not in general act trivially on deformations. %
The group~$\mu_v$ exactly
measures deformations~$\rho$ which
are not equivalent under conjugation by~$\widehat{\GL}_n(A)$ but become equivalent under
conjugation by~$\cGhat(A)$.

An alternative way to view these general deformation problems would be to work purely with conjugate
self-dual~$G_{F}$-representations. In this setting, the Selmer groups are naturally subgroups
of~$H^1(F(S)/F,\ad \rbar)^{\Gal(F/F^{+})}$. To compare these approaches, note that,
(since~$\rbar$ is irreducible so~$(\ad \rbar)^{G_F} = k$)  there is a natural
inflation--restriction sequence:
\numequation \label{infres}
\begin{aligned}
 0 \rightarrow H^1(F/F^{+},k)  & \rightarrow H^1(F(S)/F^{+},\ad \rbar)
\xrightarrow{\mathrm{res}}
H^1(F(S)/F,\ad \rbar)^{\Gal(F/F^{+})} \\
 \rightarrow H^2(F/F^{+},k) &
\rightarrow H^2(F(S)/F^{+},k).\end{aligned}
\end{equation}
The first group is~$1$-dimensional, and (for example by an explicit
cocycle computation as in the proof of Lemma~\ref{lem: reality check
  can't kill all of dual Selmer} below) its localization at any
prime~$v\in S$ agrees with~$\mu_v$.
On the other hand, since~$F/F^{+}$ is~CM, for any real place~$v$, the composite map
$$H^2(F/F^{+},k) 
\rightarrow H^2(F(S)/F^{+},k) \rightarrow H^2(F^{+}_{v},k)$$
 is injective (and indeed an isomorphism), so the restriction map in~(\ref{infres}) is surjective.
\end{remark}

\subsection{Taylor--Wiles systems}\label{subsec: TW systems}We briefly
recall the deformation condition that we use at Taylor--Wiles primes,
and the notion of a Taylor--Wiles system, following~\cite[\S\S
2.3.2, 2.4]{MR3598803}. A Taylor--Wiles prime is a finite place~$v$
of~$F^+$ which splits in~$F$ and is such
that~$\rbar$ is unramified at~$v$ with~$\rbar(\Frob_v)$ semi-simple. At
such a~$v$ we choose an eigenvalue~$\alpha_v\in k$ of
multiplicity~$n_1$, and decompose %
\numequation\label{eqn: TW prime mod p}
\rhobar|_{G_{F^+_v}}=\overline{A}_v\oplus\overline{B}_v
\end{equation}
with
$A_v(\Frob_v)=\alpha_v\cdot 1_{n_1}$.
The local deformation problem $\cD_v^{\TW}$ is given for each
$R\in\CNL_\Lambda$ by declaring that $r\in\cD_v^{\TW}(R)$ if there is
a decomposition \[\rho=A_v\oplus B_v\] lifting~\eqref{eqn: TW prime
  mod p}, with~$B_v$ unramified and ~$A_v|_{I_{F_v^+}}=\psi_v\cdot
1_{n_1}$ for some $\psi_v:I_{F_v^+}\to R^\times$. Note that
while~$\cD_v^{\TW}$ depends on the choice of~$\alpha_v$, it is omitted
from the notation. Note also that~$\cD_v^{\TW}$ is indeed a local
deformation problem, by~\cite[Lem.\ 4.2]{jack}.

We write~$\Delta_v=k(v)^\times(2)$ for the $2$-part
of~$k(v)^\times$. For any~$\rho\in\cD_v^{\TW}(R)$, the character $\psi_v\circ\Art_{F_v^+}$ gives a
canonical homomorphism $\Delta_v\to R^\times$.

\begin{defn}
  \label{defn: TW system}Let $\cS = (F,\rbar, \cO, \chi,S, \{\Lambda_v\}_{v\in S}, \{\cD_v\}_{v\in S})$ be a global deformation
  problem, and set $T=S\setminus S_\infty$. For each~$N\ge 1$, a
  \emph{Taylor--Wiles datum of level~$N$} is a pair
  $(\cQ,(\alpha_v)_{v\in \cQ})$ such that
  \begin{enumerate}[(i)]
  \item $\cQ$ is a finite set of places of~$F^+$.
  \item For each~$v\in \cQ$, we have~$v\notin S$, and~$v$ splits
    completely in~$F(\zeta_{2^N})$. 
  \item For each~$v\in \cQ$, $\rhobar(\Frob_v)$ is semi-simple,
    and~$\alpha_v\in k$  is an eigenvalue of $\rhobar(\Frob_v)$.
  \end{enumerate}
\end{defn}
For each Taylor--Wiles datum $(\cQ,(\alpha_v)_{v\in \cQ})$, we define
the corresponding 
\emph{augmented global deformation problem} \[\cS_{\cQ} = (F,\rbar, \cO,
  \chi,S \cup \cQ, \{\Lambda_v\}_{v\in S}, \{\cD_v\}_{v\in S}\cup
  \{\cD_v^{\TW}\}_{v\in \cQ}), \]where for each $v\in\cQ$ the local
deformation problem~$\cD_v^{\TW}$ is defined using the choice of
eigenvalue~$\alpha_v\in
k$. Write~$\Delta_{\cQ}:=\prod_{v\in\cQ}k(v)^\times(2)$. Then the
canonical homomorphisms $\Delta_v\to R^\times$ give a canonical
homomorphism $\cO[\Delta_\cQ]\to R_{\cS_{\cQ}}$, and we have a
canonical identification $R_{\cS_\cQ}\otimes_{\cO[\Delta_{\cQ}]}\cO=R_{\cS}$.

The following lemma shows that it is not possible to kill all the
classes in the dual Selmer group by adding Taylor--Wiles primes. We
will shortly see that it is however possible to kill all but one class
(more precisely, all but a one-dimensional space of classes), and that
                                    this is enough to patch.%
\begin{lem}
  \label{lem: reality check can't kill all of dual Selmer}Let $\cS = (F,\rbar, \cO, \chi,S, \{\Lambda_v\}_{v\in S}, \{\cD_v\}_{v\in S})$ be a global deformation
  problem. Set~$T=S\setminus S_\infty$. Then for~$N$ sufficiently large with respect to~$F^{+}$, and any Taylor--Wiles
  datum~$\cQ$ of level~$N$, we have $h^1_{\cS_{\cQ}^\perp,T}\ge 1$. In particular,
  taking $\cQ=\emptyset$, we have $h^1_{\cS^\perp,T}\ge 1$.
\end{lem}
\begin{proof} 
Let~$K^{+}_m$ be the maximal totally real subfield of~$\Q(\zeta_{2^m})$. 
Assume without loss of generality that~$K^{+}_{m-1} \subseteq F^{+} \subsetneq K^{+}_{m}$.
Let~$L^+/F^+$ be the totally real quadratic extension
  of~$F^+$ given by~$L^{+} =  F^{+}.K^{+}_{m}$, so~$L^{+} \subset F^{+}(\zeta_{2^m}) \subseteq  F(\zeta_{2^m})$, and assume that~$N \ge m$.
   We  claim that the class~$\psi$ in $H^1(F(S)/F^+,\ad\rbar(1))$ which is
  inflated from~$H^1(L^+/F^+,k)$ (via the inclusion of the scalar
  matrices into $\ad\rbar=\ad\rbar(1)$) is necessarily contained in
  $H^1_{\cS_{\cQ}^\perp,T}$ for all~$\cQ$ of level~$N\ge m$.

  Since~$L^{+}$ is totally real, $\psi$ is
  trivial at all of the infinite places of~$F^+$. In addition, each
  prime in~$\cQ$  splits completely in~$F(\zeta_{2^N})$ and hence also in~$F(\zeta_{2^m})$ and thus in~$L^{+}$.
  Hence~$\psi$ is also trivial at all of
  the places in~$\cQ$.
  
    It remains to show that (the restriction of)
  $\psi$ is contained in~$\mu_v^\perp$ for each finite place~$v\in S$. 
Let~$\Delta_v \subset H^1(F^{+}_v,\ad \rbar)$ denote the the image
  of the map~$H^1(F^{+}_v,k)\to H^1(F^{+}_v,\ad \rbar)$ induced by the map~$k \rightarrow \ad \rbar$.
  Certainly~$\psi_v \in \Delta_v$; we now show
  that~$\mu_v \subset \Delta_v$ and then analyze the pairing~$\Delta_{v} \times \Delta_v \rightarrow k$. %

From~\eqref{eqn: muv via boundary
  map} we have~$\mu_v = \mathrm{im}(\delta_v)$, where \[\delta_v: k \rightarrow H^1(F^{+}_v,\ad \rbar)\]
is the boundary map in the long exact sequence in
cohomology obtained from the action of~$G_{F^+_v}$ on~\eqref{eqn: ses
  ad and fgn}. The short exact sequence~\eqref{eqn: ses ad and fgn}
has a ~$G_F$-equivariant splitting~$\fg_n \rbar
\simeq \ad \rbar \oplus k$.
Choose a lift~$(0,1)$ of~$1 \in k$ compatible with this splitting. The
corresponding cocycle~$\phi = \delta_v(1)$  vanishes
on~$G_{F}$ and sends~$c \in G_{F^{+}} \setminus G_{F}$
to~$c (0,1) - (0,1)= (I_n,0)  \in \ad \rbar$%
, so~$\phi \in
\Delta_v$, i.e.\ $\mu_v\subset \Delta_v$.

Since the images of  cocycles in~$\Delta_v$ are contained in $k\subset\ad\rbar$
and since the self-duality on~$\ad \rbar$ is given by~$(X,Y) \mapsto \Tr(XY)$,
the Tate pairing~$\Delta_v \times \Delta_v \rightarrow k$ may be computed by first evaluating the pairing on~$k \subset \ad \rbar$
and then multiplying by~$n = \Tr(I_n)$. If~$n$ is even, it follows that
~$\Delta_v \subset \Delta^{\perp}_v$, so that $\Delta_v\subset
\mu_v^\perp$, as required. This completes the proof when~$n$ is even (which
is the case which we ultimately use). When~$n$ is odd, we must investigate more closely
the pairing on~$k \subset \ad \rbar$.

The Tate pairing
\numequation \label{tatepairing}
H^1(F^{+}_v,k) \times H^1(F^{+}_v,k) \rightarrow H^2(F^{+}_v,k) = k
\end{equation}
for any~$v$ is given by the (local) Hilbert symbol:
\numequation \label{hilbertpairing}
\langle \star,\star \rangle_v: F^{+,\times}/F^{+,\times 2} \times F^{+,\times}/F^{+,\times 2} \rightarrow \mu_2.
\end{equation}
More precisely,  the relation between~(\ref{tatepairing}) and~(\ref{hilbertpairing}) is given by
first identifying~$k$ with~$\mu_2 \otimes k$ and~$H^1(F^{+}_v,\mu_2)$ with~$F^{+,\times}/F^{+,\times 2}$,
and then tensoring~(\ref{hilbertpairing}) with~$k$.
The class associated to~$\mu_v \subset \Delta_v$ is the class~$c_{F,v}$  coming by localization
from the extension~$F/F^{+}$, and the class associated to~$\psi \subset
\Delta_v$ is the class~$c_{L^{+},v}$ coming from~$L/L^{+}$.
 It therefore suffices to show that %
  \[ \langle c_{F,v}, c_{L^{+},v} \rangle_v = 1.\]
By assumption, $F/F^{+}$ is unramified at all finite primes and~$L^{+}/L$ is unramified outside~$v|2$,
so the pairing vanishes at all finite primes away from~$v|2$. (For  primes of odd residue characteristic, the Hilbert symbol vanishes when restricted to units).
Since~$L^{+}/F^{+}$ is totally real, the Hilbert symbol also vanishes at~$v|\infty$.
Finally, we are assuming that the primes above~$2$ in~$F^{+}$ are totally split in~$F/F^{+}$, so~$c_{F,v}$ is trivial for~$v|2$ and the pairing also vanishes for~$v|2$.
   \end{proof}

The following is identical to~\cite[Prop.\
2.21]{MR3598803} 
except with  ``adequate'' replaced by  
  ``nearly adequate.''
    The proof is identical, but we go through it
in detail in order to show exactly where each hypothesis in
Definition~\ref{defn: nearly adequate} is used (or more precisely, to
show that the hypotheses in
Definition~\ref{defn: nearly adequate} are the only ones used in the
proof of~\cite[Prop.\
2.21]{MR3598803}, and the stronger assumption made there that $\rhobar(G_F)$ is
adequate is in fact never used).
\begin{prop}
  \label{prop: existence of TW systems}Let $\cS = (F,\rbar, \cO, \chi,S, \{\Lambda_v\}_{v\in S}, \{\cD_v\}_{v\in S})$ be a global deformation
  problem, and let $T=S\setminus S_\infty$. Assume that:
  \begin{enumerate}[(i)]
  \item \label{item: oddness and unrestricted at infinity} For each~$v\in S_\infty$, $\mu(c_v)=-1$ and
    ~$\cD_v=\Lift_v^\square$.
  \item \label{item: F is given by i} $F=F^+(\sqrt{-1})$.
  \item \label{item: strongly residually odd} If~$n$ is even, there exists an
       infinite place~$v$ of~$F^+$ such that $(\rhobar,\mu)$ is
       strongly residually odd at~$v$. %
     \item \label{item: near adequacy of rhobar} The group $\rhobar(G_F)\subseteq\GL_n(k)$ is nearly adequate.
     \end{enumerate}
Write~$q=h^1_{\cS^\perp,T}-1$      and $g=q+\#
T-1-[F^+:\Q]n(n-1)/2$. Then for each~$N\ge 1$ there are infinitely
many Taylor--Wiles data $(\cQ,(\alpha_v)_{v\in\cQ})$ of level~$N$ such
that $\# Q=q$ and the map $R^{\loc}_{\cS,T}\to R^T_{\cS_\cQ}$ can be
extended to a surjection $R^{\loc}_{\cS,T}\llb X_1,\dots,X_g\rrb\onto R^T_{\cS_\cQ}$.
\end{prop}
\begin{proof}
  We follow the proof of~\cite[Prop.\ 2.21]{MR3598803} very closely,
  assuming throughout that $p=2$. As usual in arguments constructing
  Taylor--Wiles data, the proof begins by using the material on Galois
  cohomology recalled above to reduce to showing that for each~$N\ge
  1$, there are infinitely many Taylor--Wiles data $(\cQ,(\alpha_v)_{v\in\cQ})$ of level~$N$ such
that  $h^1_{\cS_{\cQ}^\perp,T}=1=h^1_{\cS^\perp,T}-\#\cQ$. (This second equation is of course
equivalent to $\# Q=q$ by the definition of~$q$.) This reduction
uses assumption~\eqref{item: oddness and unrestricted at infinity} in
the statement of the proposition, but
makes no use of adequacy, so goes over unchanged under our assumptions.
Since a Taylor--Wiles datum of level~$N$ is also a Taylor--Wiles
datum of level~$M$ for any~$M\le N$, we can and do assume
that~$N$ is sufficiently large to ensure
that $F_N:=F(\zeta_{2^N})$ strictly contains~$F$. %
Let
$[\psi]\in H^1_{\cS^\perp,T}\subseteq H^1(F(S)/F^+,\ad\rbar(1))$ be a
cohomology class with nonzero image in $H^1(F(S)/F_N,\ad\rbar(1))$. We
claim that:

\begin{claim} 
\label{theclaim}
There are infinitely many Taylor--Wiles data
$(\{w\},\alpha_w)$ of level~$N$ with $[\psi]\notin H^1_{\cS^\perp_{\{w\}},T}$.
\end{claim}

Admitting the claim for now, the proof of the proposition is as
follows. Write~$s$ for the dimension of the image
of~$H^1_{\cS^\perp,T}$ in~$H^1(F(S)/F_N,\ad\rbar(1))$. Applying the
claim repeatedly, we see that there are infinitely many Taylor--Wiles
data $(\cQ,(\alpha_v)_{v\in\cQ})$ of level~$N$ such that $\#Q=s$,
$h^1_{\cS_{\cQ}^\perp,T}=h^1_{\cS^\perp,T}-\#Q$, and the morphism
$H^1_{\cS_{\cQ}^\perp,T}\to H^1(F(S)/F_N,\ad\rbar(1))$ is zero.  We
therefore have \numequation\label{eqn: H1 contained in kernel of
  restriction}H^1_{\cS_{\cQ}^\perp,T}\subseteq
H^1(F_N/F^+,\ad\rbar(1)^{G_{F_N}}) \end{equation} (because by
inflation-restriction, $H^1(F_N/F^+,\ad\rbar(1)^{G_{F_N}})$ is the
kernel of the restriction map
$H^1(F(S)/F^+,\ad\rbar(1))\to H^1(F(S)/F_N,\ad\rbar(1))$).

It only remains to show that $h^1_{\cS_{\cQ}^\perp,T}=1$. It is now
time to use that~$\rhobar(G_F)$ is nearly adequate. We begin by using
points~\eqref{item: weakly adequate} and~\ref{item: H1 scalars
  vanishes} of Definition~\ref{defn: nearly adequate}. The latter
implies that (indeed, is equivalent to) $\rhobar(G_F)$ has no normal
subgroups of index~$2$, so that~$\rhobar(G_{F_N})=\rhobar(G_F)$.  Then
the former implies that $\rhobar(G_{F_N})$ acts absolutely
irreducibly, so that $\ad\rbar(1)^{G_{F_N}}=\ad\rbar^{G_{F_N}}=k$, the
scalar matrices. In particular we have
$H^1(F_N/F^+,\ad\rbar(1)^{G_{F_N}})=H^1(F_N/F^+,k)$. 
We are assuming that $F = F^{+}(\sqrt{-1})$ 
(Assumption~\eqref{item: F is given by i} in the statement of the
proposition) and that~$F_N = F(\zeta_{2^N})$ is non-trivial over~$F$.
Together, these imply (since $p=2$) 
that $H^1(F_N/F^+,k)$ is two-dimensional. By assumption~\eqref{item:
  strongly residually odd}, together with \cite[Lem.\
2.17(ii)]{MR3598803} (and~\cite[Lem.\ 2.16]{MR3598803} in the case
that~$n$ is odd), there is a place $v\in S_\infty$ such that the
morphism $H^1(F_v^+,k)\to H^1(F_v^+,\ad\rbar)$ is injective (this
morphism being the one induced by the inclusion of the scalar matrices
in $\ad\rbar$). In particular, for such a place the composite
$H^1(F_N/F^+,k)\to H^1(F_v^+,k)\to H^1(F_v^+,\ad\rbar(1))$ is nonzero
(because the first map is nonzero, for example because $F/F^+$ is an
imaginary CM extension contained in~$F_N/F^+$ and~$F^{+}_v$ is real). (Note that
in~\cite{MR3598803} there is a typo, asserting that this composite is
injective, but being nonzero is all that is needed.) Now, by
definition (i.e.\ by~\eqref{eqn: defn of dual Selmer group}) the
restriction to $H^1(F_v^+,\ad\rbar(1))$ of any class
in~$H^1_{\cS^\perp_{\cQ},T}$ vanishes; indeed, our choice of~ $T$
gives $S_\infty=S\setminus T$, and by assumption~\eqref{item: oddness
  and unrestricted at infinity} in the statement of the proposition,
we have~$\cL_v^\perp=0$ for all~$v\in\infty$. Going back
to~\eqref{eqn: H1 contained in kernel of restriction} (and recalling
Lemma~\ref{lem: reality check can't kill all of dual Selmer}) we see
that $h^1_{\cS_Q^\perp,T}=1$, and we are done.

It remains to prove Claim~\ref{theclaim}. Accordingly, we let
$$[\psi]\in H^1_{\cS^\perp,T}\subseteq H^1(F(S)/F^+,\ad\rbar(1))$$ 
be a
cohomology class with nonzero image in
$H^1(F(S)/F_N,\ad\rbar(1))$. By~\cite[Lem.\ 2.19]{MR3598803}, finding
a Taylor--Wiles datum
$(\{w\},\alpha_w)$ of level~$N$ with $[\psi]\notin
H^1_{\cS^\perp_{\{w\}},T}$ amounts to choosing $w,\alpha_w$ such that
\begin{itemize}
\item  $w$ splits completely in~$F_N$, 
and~$\rhobar(\Frob_w)$ is semi-simple; and
\item $\alpha_w\in k$ is an eigenvalue of~$\rhobar(\Frob_w)$ such that
  $\tr e_{\Frob_w,\alpha_w}\psi(\Frob_w)\ne 0$, where
  $e_{\Frob_w,\alpha_w}$ is the unique idempotent
  in~$k[\rhobar(\Frob_w)]$ whose image is the $\alpha_w$-eigenspace
  of~$\rhobar(\Frob_w)$. 
\end{itemize}
By Chebotarev, it therefore suffices to find $\sigma\in G_{F_N}$ and $\alpha\in
k$ such that~$\rhobar(\sigma)$ is semi-simple, and $\alpha$ is an
eigenvalue of~$\rhobar(\sigma)$ with~$\tr
e_{\sigma,\alpha}\psi(\sigma)\ne 0$.

Let $K/F$ be the extension cut out by~$\ad\rhobar$, and write
$K_N=K\cdot F_N$. Let $f$ denote the image of~$[\psi]$ under the
restriction map \numequation\label{eqn: restricting psi to
  f}H^1(F(S)/F^+,\ad\rbar(1))\to
H^1(F(S)/K_N,\ad\rbar(1))^{G_{F^+}}.\end{equation} By the definition
of~ $K$, the action of $G_{K_N}$ on $\ad\rbar(1)$ is trivial so this
image is a homomorphism $f:\Gal(F(S)/K_N)\to\ad\rbar$. We claim that
$f\ne 0$. %
          To see this, note that the restriction
map~\eqref{eqn: restricting psi to f} factors through the restriction
map
\[H^1(F(S)/F^+,\ad\rbar(1))\to H^1(F(S)/F_N,\ad\rbar(1)),\]
and by assumption, the image of $[\psi]$ in
$H^1(F(S)/F_N,\ad\rbar(1))$ is nonzero. It therefore suffices to show that the kernel of the restriction
map
          \[H^1(F(S)/F_N,\ad\rbar(1))\to H^1(F(S)/K_N,\ad\rbar(1))\] vanishes. By
          inflation-restriction, this kernel is $H^1(K_N/F_N,\ad\rbar(1))$.  As we
          saw above, by Definition~\ref{defn: nearly
          adequate}~\eqref{item: H1 scalars vanishes} we have
          $\rhobar(G_{F_N})=\rhobar(G_F)$, so that
          $H^1(K_N/F_N,\ad\rbar(1))=H^1(\rhobar(G_F),\ad\rbar)$. This
          vanishes by Definition~\ref{defn: nearly
          adequate}~\eqref{item: H1 gl vanishes}, as required.

          Let $V\subseteq\ad\rhobar$ be the $k$-vector space generated by the
          image of~$f$. Since~$f$ is $G_{F_N}$-equivariant (as it is
          restricted from~$[\psi]$), $V$ is a $k[G_{F_N}]$-module, and
          we let~$W$ be a simple $k[G_{F_N}]$-submodule of~$V$. By Definition~\ref{defn: nearly
          adequate}~\eqref{item: weakly adequate} (and again using
          that $\rhobar(G_{F_N})=\rhobar(G_F)$), we may find $\sigma_0\in G_{F_N}$ and $\alpha_0\in
k$ such that~$\rhobar(\sigma_0)$ is semi-simple, and $\alpha_0$ is an
eigenvalue of~$\rhobar(\sigma_0)$ with~$\tr
          e_{\sigma_0,\alpha_0}W\ne 0$.

          If   $\tr
          e_{\sigma_0,\alpha}\psi(\sigma_0)\ne 0$, then we are done, taking
          $\sigma=\sigma_0$ and $\alpha=\alpha_0$. Suppose instead
          that $\tr
          e_{\sigma_0,\alpha_0}\psi(\sigma_0)=0$, and choose any $\tau\in K_N$ such
          that $\tr
          e_{\sigma_0,\alpha_0}f(\tau)\ne 0$. (Such a~$\tau$ exists, because $\tr
          e_{\sigma_0,\alpha_0}W\ne 0$, and by definition~$V$ is spanned
          as a $k$-vector space
          by the varying $f(\tau)$.) We set~$\sigma=\tau\sigma_0$, so
          that~$\rhobar(\sigma)$ is a scalar multiple
          of~$\rhobar(\sigma_0)$, and we let~$\alpha$ be the
          corresponding scalar multiple of~$\alpha_0$, so that
          $e_{\sigma,\alpha}=e_{\sigma_0,\alpha_0}$. We have
          $\psi(\sigma)=\psi(\sigma_0)+\psi(\tau)$, so that $\tr
          e_{\sigma,\alpha}\psi(\sigma)= \tr
          e_{\sigma_0,\alpha_0}f(\tau)\ne 0$, as required.
\end{proof}

\subsection{Local deformation problems}\label{subsec: local
  deformation problems Un}%
We now assume that all finite
places $v\in S$ of~$F^+$ split in~$F$, and choose a place~$\tv$ of~$F$
above each~$v\in S$. We write~$\tS$ for the set of places~$\tv$
with~$v\in S$ finite, and $\tS_2\subset\tS$ for the places lying over~$2$.
For each~$v\in S$ we can and do identify liftings
of~$\rbar|_{G_{F^+_v}}$ with liftings
of~$\rhobar|_{G_{F_\tv}}$. %
\subsubsection{Local deformation problems for~$v\nmid 2$}

The following two lemmas are presumably well known, but for lack of a
reference we give a proof.
\begin{lem}
  \label{lem: uniformly potentially unramified GLn}Suppose
  that~$v\nmid 2$.
    Then there is a finite extension $F'_\tv/F_\tv$ such that any
    lifting of $\rhobar|_{G_{F_\tv}}$ becomes unipotently ramified
    after restriction to~$G_{F'_\tv}$.
\end{lem}
\begin{proof}Since the universal lifting ring
  of~$\rhobar|_{G_{F_\tv}}$ is $\cO$-flat by~ \cite[Thm.\
  2.5]{shottonGLn}, it suffices to prove this for closed points of its
  generic fibre. Since this generic fibre has finitely many connected
  components, it suffices to prove the result for the closed points of
  any single connected component. For each connected
  component, it suffices to prove the result for a single point
  on that component by~\cite[Lem.\
  1.3.4(1)]{BLGGT} (a theorem of Choi), and the result is immediate.
  \end{proof}

\begin{lem}
  \label{lem: deformations of semi-simple are sums of
    characters}Suppose that $v\nmid 2$, that $\rhobar|_{G_{F_\tv}}$ is
  unramified, and that~$\rhobar(\Frob_\tv)$ is regular semi-simple. Then
  any lifting of~$\rbar|_{G_{F_\tv}}$ is strictly equivalent to a direct sum of
  characters. In particular, there is a finite extension $F'_\tv/F_\tv$
  such that any lifting of $\rhobar|_{G_{F_\tv}}$ becomes unramified after
  restriction to~$G_{F'_\tv}$.
\end{lem}
\begin{proof}
The second statement follows from the first by Lemma~\ref{lem:
  uniformly potentially unramified GLn}. Since~$q_v\equiv 1\pmod{2}$,
the first statement is standard, and may for example be proved by an
identical argument to the proof of Lemma~\ref{lem:anyliftisTW}.  %
\end{proof}

Now let
$\chi_{v,1}, \dots, \chi_{v,n}: \cO_{F_\tv}^\times \rightarrow
\cO^\times$ be finite order characters, which are trivial modulo
$\varpi$.  Suppose that $\rhobar|_{G_{F_\tv}}$ is trivial. We write
$\cD_v^{\chi_v}$ for the set of liftings $\rho$ of $\rhobar|_{G_{F_\tv}}$ to
objects of $\CNL_\cO$ such that for all $\sigma \in I_{F_\tv}$, we have
\[\cha_{\rho(\sigma)}(X) = \prod_{i=1}^n\left(X-\chi_{v,j}(\Art_{F_\tv}(\sigma)\right)^{-1}).\]
 Write $R_v^{\chi_v}$ for the corresponding local lifting ring. The
 following is~\cite[Prop.\ 3.16]{jack}.
 \begin{prop}\label{levelraising}\leavevmode \begin{enumerate}\item
     Suppose that $\chi_{v,j} = 1$ for each $j$. Then each irreducible
     component of $R_v^1$ has dimension $n^2 + 1$, and every prime of
     $R_v^1$ minimal over $\varpi$ contains a unique minimal
     prime. Every generic point of~$R_v^1$ is of characteristic zero.
   \item Suppose that the $\chi_{v,j}$ are pairwise distinct. Then
     $\Spec R_v^{\chi_v}$ is irreducible of dimension $n^2 + 1$, and
     its generic point is of characteristic zero.
\end{enumerate}
\end{prop}

\subsubsection{Local deformation rings for~$v|2$:  ordinary deformation rings for~$v|2$}\label{subsec:
  unitary group ordinary local def rings} We now  recall the ordinary deformation problems introduced
in~\cite[\S 3]{ger}, and studied there and in~\cite[\S
3.3.2]{MR3327536}. Suppose that~ $v|2$ and that~$\rhobar|_{G_{F_\tv}}$ can be
conjugated to an upper-triangular representation whose diagonal
characters are $\chibar_1,\dots,\chibar_n:G_{F_\tv}\to k^{\times}$
 (in
that order). Let $\Lambda_{\GL_n,v}$ be the completed group ring of the group
$I_{F_\tv}^{\ab}(2)^n$, where $(2)$ denotes pro-$2$
completion. %
Let
$(\chi_1,\dots,\chi_n)$ denote the universal $n$-tuple of characters
$I_{F_\tv}\to\Lambda_{\GL_n,v}^\times$
lifting~$(\chibar_1|_{I_{F_\tv}},\dots,\chibar_n|_{I_{F_\tv}})$.
  
\begin{prop}\label{prop: ord local def rings Un}  There is a local
  deformation problem~$\cD_{v}^{\triangle}$, represented by a
  $\CNL_{\Lambda_{\GL_n,v}}$-algebra $R^{\triangle}_v$ with the following
  properties.
  
  \begin{enumerate} %
\item $R^{\triangle}_v$ is  reduced and $2$-torsion free.
\item Let~$E'/E$ be a finite extension with ring of
  integers~$\cO_{E'}$, and fix a morphism of local $\cO$-algebras
  $\Lambda_{\GL_n,v}\to\cO_{E'}$. Then a morphism of local
  $\Lambda_{\GL_n,v}$-algebras $R_v^\square\to\cO_{E'}$ factors through
  $R^{\triangle}_v$ if and only if the corresponding representation
  $\rho:G_{F_\tv}\to\GL_n(\cO_{E'})$ is $\GL_n(\cO_{E'})$-conjugate to
  an upper-triangular representation whose ordered diagonal
  characters~$(\psi_1,\dots,\psi_n)$ are such that for each~ $i$,
  $\psi_i|_{I_{F_\tv}}$ is equal to the pushforward of
  $(\chi_1,\dots,\chi_n)$ along $\Lambda_{\GL_n,v}\to\cO_{E'}$.
\item Suppose that $\rhobar|_{G_{F_\tv}}$ is  trivial, and that
  $[F^+_v:\Q_2]>n(n-1)/2+1$.  Let $Q$ be a minimal prime of $\Lambda_{\GL_n,v}$. Then $\Spec
  R^{\triangle}_v/Q$ is geometrically irreducible of dimension
  $[F^+_v:\Q_2]n(n+1)/2+n^2+1$, and $R^{\triangle}_v/(Q,\varpi)$ is
  generically reduced.
\end{enumerate}

\end{prop}
\begin{proof}
The lifting ring~$R^{\triangle}_v$ is defined in~\cite[\S 3]{ger}; see
~\cite[\S 3.3.2]{MR3327536} for a summary of its definition.   It
is reduced and flat over~$\cO$ by construction. 
The remaining points are ~\cite[Lem.\ 3.3]{ger} and~\cite[Prop.\ 3.14(3)]{MR3327536}.
\end{proof}%
We will use the following remark in the proof of Proposition~\ref{prop: the modularity result for p equals 2 or 3}.
\begin{rem}
  \label{rem:Jack-shows-don't-need-flat-closure}Note that
~\cite[\S 3]{ger} defines $\Spec R^{\triangle}_v$ as the  flat closure of the
scheme-theoretic image of a projective morphism~$\pi:\cG_v\to \Spec
R^{\square}$, it is shown in~\cite[Lem.\ 3.11]{MR3327536} that if $\rhobar|_{G_{F_\tv}}$ is  trivial and
  $[F^+_v:\Q_2]>n(n-1)/2+1$, then ~$\cG_v$ is already flat over~$\cO$, so that
  $\Spec R^{\triangle}_v$ is equal to the
scheme-theoretic image $\pi:\cG_v\to \Spec
R^{\square}$.
\end{rem}
\begin{defn}\label{ord} Let $\lambda \in (\Z_+^n)^{\Hom(F_\tv,
    K)}$. We say that a continuous representation
  $\rho: G_{F_\tv} \rightarrow \GL_n(\cO)$ is \emph{ordinary of
    weight $\lambda$} if:
  \begin{enumerate} \item There exists a increasing invariant
    filtration $\Fil^i$ of $\cO^n$, with each $\gr^i \cO^n$ an
    $\cO$-module of rank one.
\item Write $\chi_i$ for the character
  $G_{F_\tv} \rightarrow \cO^\times$ giving the action on
  $\gr^i \cO^n$. Then for every $\alpha \in F_\tv^\times$ sufficiently
  close to 1, we have
\[(\chi_i \circ \Art_{F_\tv}(\alpha)) = \prod_\tau (\tau(\alpha))^{-(\lambda_{\tau, n - i + 1} + i - 1)}.\]
\end{enumerate}
\end{defn}%

\subsection{Automorphic forms on definite unitary groups}\label{subsec: auto forms Un}We now
introduce the spaces of automorphic forms that we work with,
following~\cite[\S 4]{MR3598803} and~\cite[\S 2]{ger}. We suppose throughout
Subsection~\ref{subsec: auto forms Un} that the following hypothesis
holds.

\begin{hypothesis}
  \label{hyp: existence of quasisplit unitary group}\leavevmode
  \begin{itemize}
  \item $F/F^+$ is everywhere unramified, and each place $v|2$
    of~$F^+$ splits in~$F$.
  \item $n[F^+:\Q]\equiv 0\pmod{4}$.
  \end{itemize}
\end{hypothesis}
Let $c$ denote the non-trivial element of $\Gal(F/F^+)$. By
Hypothesis~\ref{hyp: existence of quasisplit unitary group}, we can find a unitary
group~$G/F^+$ which splits over~$F$, and is such that:
\begin{itemize}
\item $G(F^+_v)$ is quasi-split at all finite places~$v$ of~$F^+$.
\item $G(F^+\otimes_{\Q}\R)$ is compact.
\end{itemize}

We can and do choose an integral model of~$G$ over~$\cO_{F^+}$ (which
we continue to denote by~$G$) in such a way that if
 $v$ is a finite place of $F^+$ which splits as $v = \tv \tv^c$ in
 $F$, then there is an isomorphism
\[\iota_\tv: G(\cO_{F_v^+}) \isoto \GL_n(\cO_{F_\tv}).\]

For each place~$v|2$ of~$F^{+}$ we choose a place $\tv|v$ of~$F$, and 
let~$\widetilde{S}_2$ be the set of~$\tv$ for~$v|2$. Let~ $\widetilde{I}_2$ denote the set of embeddings $F \hookrightarrow
E$ inducing a place in $\widetilde{S}_2$. To each $\lambda \in
(\Z^n_+)^{\widetilde{I}_2}$ %
there is an associated finite free
$\cO$-module $M_\lambda$ with a continuous action of $\prod_{v\in
  S_2}G(\cO_{F_v^+}) \isoto \prod_{v\in S_2}\GL_n(\cO_{F_\tv})$,
constructed as the tensor product over~$\tau\in\tS_2$ of the  the algebraic representations of $\GL_n/\cO$ with highest weight~$\lambda_\tau$.

We now write $S=S_\infty\disjointunion T$, and let~$R\subset T$ be a
(possibly empty)
set of places disjoint from~$S_2$. For each place~$v\in R$ we fix a choice of
~$\tv$ (a place of~$F$ dividing~$v$). Suppose that $U = \prod_v U_v$ is
an open compact subgroup of $G(\A^\infty_{F^+})$ such that
$U_v \subset \iota_\tv^{-1} \Iw(\tv)$ for $v \in R$.  (Here~$\Iw(\tv)$
is the Iwahori subgroup of $\GL_n(\cO_{F_\tv})$ consisting of matrices
which are upper-triangular modulo~$\tv$, with pro- $v$ Iwahori
subgroup $\Iw_1(\tv) \subset \Iw(\tv)$.)

For each $v \in R$, we choose a character 
\[\chi_v = \chi_{v, 1} \times \dots \times \chi_{v, n}: \Iw(\tv)/\Iw_1(\tv) \rightarrow \cO^\times,\]
the decomposition being with respect to the natural isomorphism 
\[\Iw(\tv)/\Iw_1(\tv) \cong (k(\tv)^\times)^n.\]
 We set 
\[M_{\lambda, \{\chi_v\}} = M_\lambda \otimes_\cO \left(\bigotimes_{v
      \in R} \cO(\chi_v)\right),\] %
a representation of  $G(\cO_{F^+,2}) \times
\prod_{v \in R} \Iw(\tv)$. %

If $A$ is an $\cO$-module, and $U_v \subset G(\cO_{F^+_v})$ for $v|2$,
then we write $S_{\lambda, \{\chi_v\}}(U, A)$ for the set of functions 
\[f: G(F^+)\backslash G(\A^\infty_{F^+}) \rightarrow M_{\lambda,\{\chi_v\}} \otimes_\cO A\]
such that for every $u \in U$, we have $f(gu) = u_{S_2 \cup R}^{-1} f(g)$, where $u_{S_2 \cup R}$ denotes the projection to $\prod_{v \in S_2 \cup R} U_v$.
If $R$ is empty then we write $S_{\lambda, \{\chi_v\}}(U, A) = S_\lambda(U, A)$.

We will sometimes assume that~$U$ is sufficiently small in the
following sense.

\begin{defn}
  \label{defn: suff small Un}We say that~$U$ is \emph{sufficiently
    small} if for some finite place~$v$ of~$F^+$, the projection
  of~$U$ to~$G(F^+_v)$ contains no element of finite order other than
  the identity.
\end{defn}

Let $w$ be a place of $F$ split over $F^+$ and not contained in $S$,
and let $\varpi_w$ be a uniformizer of $F_w$. Write \[\alpha^j_{\varpi_w}:= \diag(\underbrace{\varpi_w, \dots, \varpi_w}_j, 1, \dots, 1).\] The spaces $S_{\lambda,\{\chi_v\}}(U, A)$ receive an action of the Hecke operators
\[T_w^j:= \iota_w^{-1} \left( \left[ \GL_n(\cO_{F_w}) \alpha^j_{\varpi_w}\GL_n(\cO_{F_w})\right]\right).\]

For integers $0 \leq b \leq c$, and $v \in S_2$, we consider the
subgroup $\Iw(\tv^{b, c}) \subset \GL_n(\cO_{F_\tv})$ defined as those
matrices which are congruent to an upper-triangular matrix modulo
$\tv^c$ and congruent to a unipotent upper-triangular matrix modulo
$\tv^b$. We set
$U(\mathfrak{l}^{b, c}) = \prod_{v\notin S_2}U_v \times \prod_{v \in S_2} \Iw(\tv^{b,
  c})$. (Our use of~$\widetilde{l}$ is in order to follow the notation
of~\cite{ger}.)

We now recall from~\cite[Defn.\ 2.8]{ger} some additional Hecke operators at the places dividing $2$.
For each
$v\in S_2$ %
we let $\varpi_\tv$ be a
uniformizer of $F_\tv$. As above we write
\[\alpha^j_{\varpi_\tv} = \diag(\underbrace{\varpi_\tv, \dots,
    \varpi_\tv}_j, 1, \dots, 1),\] and we set
\[U^j_{\lambda, \varpi_\tv} = (w_0
  \lambda)(\alpha^j_{\varpi_\tv})^{-1} \left[U(\mathfrak{l}^{b,c})
    \iota_\tv^{-1} (\alpha^j_{\varpi_\tv})
    U(\mathfrak{l}^{b,c})\right],\]where as usual~$w_0$ is the longest
element of the Weyl group. If $u \in T(\cO_{F_\tv})$ then we write 
\[\langle u \rangle = \left[ U(\mathfrak{l}^{b,c}) \iota_\tv^{-1} (u)
    U(\mathfrak{l}^{b,c})\right].\]

By~\cite[Lem.\ 2.10]{ger}, these
operators commute with each other and act on the spaces
$S_{\lambda, \{\chi\}}(U(\mathfrak{l}^{b,c}), A)$, compatibly  with the
inclusions
\[S_{\lambda, \{\chi\}}(U(\mathfrak{l}^{b,c}), \cO) \subset S_{\lambda, \{\chi\}}(U(\mathfrak{l}^{b',c'}), \cO),\]
where $b \leq b'$ and $c \leq c'$.

We write $\T^T_{\lambda, \{\chi_v\}}(U(\mathfrak{l}^{b, c}), A)$ for
the $\cO$-subalgebra of $\End_\cO(S_{\lambda, \{\chi_v\}}(U, A))$
generated by the operators $T^j_w$ and $(T^n_w)^{-1}$ as above and all
the operators $\langle u \rangle:=\prod_{v\in S_2}\langle u_v\rangle$ for
\[u=(u_v)_{v\in S_2 } \in T(\cO_{F^+,2}) = \prod_{v \in S_2} T(\cO_{F^+_v}),\]where~$T$ denotes
the usual diagonal torus in~$\GL_n$. With
these identifications, the operators $\langle u \rangle$ endow each
Hecke algebra $\T^T_{\lambda, \{\chi_v\}}(U(\mathfrak{l}^{b, c}), A)$
with the structure of an algebra for the completed group ring
\numequation\label{eqn:Lambda-for-Un}\Lambda = \cO\llb T(\mathfrak{l})\rrb,\end{equation}
where $T(\mathfrak{l})$ is defined by the exact sequence
\[\xymatrix{ 0 \ar[r] & T(\mathfrak{l}) \ar[r] & \prod_{v \in S_2} T(\cO_{F^+_v}) \ar[r] & \prod_{v \in S_2} k(v)^\times \ar[r] & 0.}\]

We have the ordinary idempotent
$e = \lim_{r\rightarrow \infty} U(\mathfrak{l})^{r!}$, where we set
\[U(\mathfrak{l}) = \prod_{v \in S_2} \prod_{j=1}^n U^j_{\lambda, \varpi_\tv}.\]
We define the ordinary Hecke algebra 
\[\T^{T,\ord}_{\lambda, \{\chi_v\}}(U(\mathfrak{l}^{b,c}), A) = e\T^T_{\lambda, \{\chi_v\}}(U(\mathfrak{l}^{b,c}), A);\]
equivalently,
$\T^{T,\ord}_{\lambda, \{\chi_v\}}(U(\mathfrak{l}^{b,c}), A) $
is the image of the Hecke algebra $\T^T_{\lambda, \{\chi_v\}}(U(\mathfrak{l}^{b,c}), A)$
in
$\End_{\cO}\left(S^{\ord}_{\lambda,
    \{\chi_v\}}(U(\mathfrak{l}^{b,c}), A)\right)$.

We set 
\[S_{\lambda, \{ \chi_v \}}(U(\mathfrak{l}^\infty), E/\cO) = \ilim_c S_{\lambda, \{\chi_v\}}(U(\mathfrak{l}^{c, c}), E/\cO),\]
which receives a faithful action of the algebra
\[\T^T_{\lambda, \{\chi_v\}}(U(\mathfrak{l}^\infty), E/\cO) = \plim_c \T^T_{\lambda, \{\chi_v\}}(U(\mathfrak{l}^{c,c}), E/\cO).\]
By~\cite[Lem.\ 2.17]{ger},  this algebra is naturally isomorphic to
\[\T^T_{\lambda, \{\chi_v\}}(U(\mathfrak{l}^\infty), \cO) = \plim_c
  \T^T_{\lambda, \{\chi_v\}}(U(\mathfrak{l}^{c,c}), \cO).\]
We can again apply the idempotent $e$ to these spaces and rings, in which case
we again decorate them with `ord' superscripts.

Specializing to the case~$\lambda=0$, we define a homomorphism $T(\mathfrak{l}) \rightarrow \T^{T, \ord}_{0, \{\chi_v\}}(U(\mathfrak{l}^\infty), \cO)^\times$
 by 
 \[u \mapsto \left( \prod_{\tau \in \widetilde{I}_2} \prod_{i=1}^n
     \tau(u_i)^{1-i}\right) \langle u \rangle\](where the~$u_i$ are the
 coordinate entries of~$u$, recalling that~$T$ is the usual diagonal maximal
 torus in~$\GL_n$).
  This gives rise to an $\cO$-algebra homomorphism $\Lambda \rightarrow \T^{T, \ord}_{0, \{\chi_v\}}(U(\mathfrak{l}^\infty), \cO)$, and we write
\[\T^{T, \ord}_{\{\chi_v\}}(U(\mathfrak{l}^\infty), \cO)\]
for $\T^{T, \ord}_{0, \{\chi_v\}}(U(\mathfrak{l}^\infty), \cO)$
endowed with this $\Lambda$-algebra structure.  This is the
\emph{universal ordinary Hecke algebra of level $U$}. It is a finite
$\Lambda$-algebra by~\cite[Cor.\ 2.21]{ger}. Along with all of the
other Hecke algebras considered above, it is reduced (by~\cite[Lem.\
2.14]{ger}). %

We can pass back from the universal ordinary Hecke algebra to the
finite level Hecke algebras in the following way. Corresponding to
each~$\lambda$ is a prime ideal~$\wp_\lambda$ %
of~$\Spec \Lambda$ as defined in~\cite[Defn.\ 2.24(1)]{ger}; the prime
ideals~$\wp_\lambda$ are dense in~$\Spec \Lambda$ (see the
proof of~\cite[Cor.\ 3.4]{ger}). By~\cite[Lem.\ 2.25]{ger}, we have a
natural identification \numequation\label{eqn: Hida control
  Un}\Hom_{\cO}(\T^{T, \ord}_{0, \{\chi_v\}}(U(\mathfrak{l}^\infty),
\cO)/\wp_\lambda,\Qtwobar)=\Hom_{\cO}(\T^{T,\ord}_{\lambda,
  \{\chi_v\}}(U(\mathfrak{l}^{1,1}), \cO),\Qtwobar).\end{equation}

We say that a RACSDC automorphic representation~$\pi$ for~$\GL_n/F$ has weight~$\lambda\in (\Z^n_+)^{\widetilde{I}_2}$ if its
infinitesimal character agrees (after composing with our fixed isomorphism
$\imath:\Qtwobar\isoto\C$) with that of the algebraic representation of weight
$\lambda$. We say that~$\pi$ is ordinary if
$(\imath^{-1}(\otimes_{v|2}\pi_v))^{\ord}\ne 0$. (See \cite[\S 2]{BLGGT}.)
The relationship between the spaces of automorphic forms considered
above and ordinary RACSDC representations is as follows. For
each~$\lambda$, write
\[\cA_\lambda:=\varinjlim_{U} S_\lambda(U,\Qtwobar).\] This is a
semi-simple admissible $\Qtwobar[G(\A^\infty_{F^+})]$-module, and by~\cite[Cor.\ 5.3, Thm.\ 5.9]{labesse} the
irreducible submodules of $\imath\cA_\lambda$ are the finite parts of
automorphic representations of~$G/F^{+}$
which arise
as the descents of automorphic representations~$\pi$ of~$\GL_n/F$ of
weight~$\lambda$. These automorphic representations~$\pi$ are isobaric
direct sums of self dual representations, and in particular, they
include the  RACSDC representations of weight~$\lambda$; and after
localizing at a non-Eisenstein maximal ideal of an appropriate Hecke
algebra (as we will always do below), the RACSDC representations are
the only ones that contribute. Furthermore the irreducible submodules of
$\imath\cA_\lambda$ which have nonzero intersection with some
$S^{\ord}_\lambda(U(\mathfrak{l}^{b,b}),\Qtwobar)$ are precisely those
which correspond to those~ $\pi$ which are ordinary.

Now let~ $\pi$ be an ordinary RACSDC automorphic
representation~$\pi$ of $\GL_n/F$, and assume that
$\rhobar=\rbar_{\pi,2}:G_F\to\GL_n(\Ftwobar)$ is
irreducible. Assuming as always that our coefficient field~$E$ is
large enough, we fix an extension of~$\rhobar$
to~$\rbar:G_{F^+}\to\cG_n(k)$.
 As above, we let~$T\supset S_2$ be a
finite set of finite places of~$F^+$ which split in~$F$. Again we
consider a subset~$R\subset (T \setminus S_2)$, and for each~$v\in R$
we fix characters $\chi_v:\Iw(\tv)/\Iw_1(\tv)\to\cO^\times$. We assume
furthermore that:
\begin{hypothesis}\label{hyp: T is big enough}\leavevmode
  \begin{itemize}
  \item $T$ contains all finite places~ lying under a place~$w$ of~$F$
    at which~$\pi_w$ is ramified, and
  \item if~ $v\in R$ then $\rhobar|_{G_{F_\tv}}$ is trivial and
    ~$\pi_v^{\Iw(\tv)}\ne 0$.
  \end{itemize}
\end{hypothesis}

Set~$S=T\cup S_\infty$. If~$v\notin S_2$ then we set~$\Lambda_v=\cO$,
while if~$v\in S_2$ we take~$\Lambda_v=\Lambda_{\GL_n,v}$
where~$\Lambda_{\GL_n,v}$ is as in Section~\ref{subsec:
  unitary group ordinary local def rings}. We define the global deformation
problem

\nummultline
  \cS_{\{\chi_v\}}=\biggl(F, \rbar,\cO, \varepsilon^{1-n}
    \delta^{n}_{F/F^+}, S, \{\Lambda_v\}_{v\in S},
    \{\cD_v^{\chi_v}\}_{v \in R} \\ \cup \{\cD_{v}^{\triangle}\}_{v \in
      S_2} \cup \{ \Lift^{\square}_v\}_{v\in S\setminus (R\cup
      S_2)}\biggr).\label{eq: defn of S chi Un}
  \end{multline}

Using the natural isomorphism
$\widehat{\otimes}_{v \in S_2} \Lambda_v \cong \Lambda$ provided by
local class field theory, we see that $R^\text{loc}_{\{\chi_v\}}$ is
naturally a $\Lambda$-algebra.

\begin{lem}
  \label{lem: Un deformation ring dimension lower bound}Every
  irreducible component of~$\Spec  R_{\cS_{\{\chi_v\}}}$ has dimension
  at least~$\dim \Lambda=1+[F:\Q]n$.
\end{lem}
\begin{proof}
By Propositions~\ref{levelraising} and~\ref{prop: ord local def rings
  Un}, together with~\cite[Thm.\ 2.5]{shottonGLn}, the ring
$R_{\cS_{\{\chi_v\}},T}^{\loc}$ is
equidimensional of dimension \numequation\label{eqn: dim of local Un deformation ring}\dim R_{\cS_{\{\chi_v\}},T}^{\loc}=1+n^2\#T+[F^+:\Q]n(n+1)/2.\end{equation}
It therefore suffices to show that there is a presentation of the
form
\[R^\text{loc}_{\{\chi_v\}}\llb x_1,\dots,x_r\rrb /(f_1,\dots,f_{r+s})
  \isoto R_{\cS_{\{\chi_v\}}}\] for some~$r,s\ge 0$ with~$s\le n^2\#T+[F^+:\Q]n(n-1)/2$.
This follows from a standard deformation-obstruction argument, and can
for example be proved exactly as in~\cite[Cor.\ 2.2.12]{cht}, using
the complex~$C^i_{\cS,T}$. Alternatively, the existence of such a
presentation is a consequence of~\cite[Prop.\ 4.2.5]{MR3152673}. (As noted
in~\cite[\S 4.1]{Bellovin}, it is assumed in~\cite[\S 4.2]{MR3152673}
that the reductive group~$G$ there is connected, but in the proof of ~\cite[Prop.\ 4.2.5]{MR3152673} this assumption
is only used in order to cite results of~\cite{MR1643682} which do not
use this assumption.)
\end{proof}

\begin{prop}\label{prop: Galois rep to T in Un case}
  Suppose that for each~$v\notin T$, the compact open subgroup~ $U_v$
  is hyperspecial; and that for each $v \in R$, we have
  $U_v \subset \iota_\tv^{-1} \Iw(\tv)$.

  Then there is a %
  maximal ideal~$\m$ of $\T^{T,
    \ord}_{\{\chi_v\}}(U(\mathfrak{l}^\infty), \cO)$ such that
  there is a surjection of $\Lambda$-algebras  \numequation\label{eqn:
  R to T map Un}R_{\cS_{\{\chi_v\}}}\to \T^{T,
      \ord}_{\{\chi_v\}}(U(\mathfrak{l}^\infty), \cO)_\m.\end{equation} The
  corresponding \emph{(}unique up to strict equivalence\emph{)}
  representation \[r_\m: G_{F^+} \rightarrow \cG_n(\T^{T,
      \ord}_{\{\chi_v\}}(U(\mathfrak{l}^\infty), \cO)_\m)\] is
  characterized by the following property: if $v \not\in T$ is a finite place of $F^+$ which splits as $w w^c$ in $F$, then $r_\m$ is unramified at $w$ and $w^c$, and $r_\m(\Frob_w)$ has characteristic polynomial
\[X^n + \dots + (-1)^j(q_w)^{j(j-1)/2} T_w^j X^{n-j} + \dots + (-1)^n (q_w)^{n(n-1)/2} T_w^{n}.\]
\end{prop}
\begin{proof}
  This is proved in exactly the same way as~\cite[Prop.\ 2.29]{ger}
  (using \cite[Cor.\ 3.4]{ger} for the compatibility at the
  places~$v\in S_2$),
  using~\cite[Lem.\ 2.4]{MR3598803} in place of~\cite[Lem.\
  2.1.12]{cht} (which is used in the proof of~\cite[Lem.\
  3.4.4]{cht}, to which the proof of~\cite[Prop.\ 2.29]{ger}
  refers). %
  (See also~\cite[Thm.\
  4.1]{MR3598803} for a detailed proof of a very similar result.). 
\end{proof}

\begin{defn}\label{defn:R-and-T-for-Un-theorem}
  If~$R$ is empty then we write $\cS^{T,\ord}$ for~$\cS_{\{\chi_v\}}$
  and~$R^{T,\ord}$ for~$R_{\cS_{\{\chi_v\}}}$, and we write
  $\T^{T, \ord}(U(\mathfrak{l}^\infty))_\m$ for
  $\T^{T, \ord}_{\{\chi_v\}}(U(\mathfrak{l}^\infty), \cO)_\m$.
\end{defn}

Before stating and proving the main result of this section, we make
a  definition, using the following (presumably well known) lemma.
    \begin{lem}
      \label{lem: deformations have uniform fixed vectors Un}If~$w\nmid 2$ is
      a finite place of~$F$, then there is a compact open
      subgroup~$U_w$  of~$\GL_n(\cO_{F_w})$, depending only
      on~$\rbar_{\pi,2}|_{G_{F_w}}$, such that if~$\pi'$ is a
      RACSDC automorphic representation of $\GL_n/F$ with
      $\rbar_{\pi,2}\cong r_{\pi',2}$
      then~$(\pi'_w)^{U_w}\ne 0$.     
    \end{lem}
    \begin{proof}
      This follows from Lemma~\ref{lem: uniformly potentially
        unramified GLn} and local-global compatibility, together with
      the compatibility of the local Langlands correspondence with
      conductors. 
    \end{proof}
\begin{defn}
  \label{defn: sufficiently deep subgroup Un}Suppose
  that~$v\in T\setminus S_2$. Then we say that a compact open subgroup
  of~$\cG_n(F_v^+)\isoto\GL_n(F_{\tv})$ is \emph{sufficiently deep} if
  it satisfies the conclusion of Lemma~\ref{lem: deformations have
    uniform fixed vectors Un} (for $w=\tv$).
\end{defn}

   \begin{thm}
     \label{thm: U(n) ordinary automorphy lifting including Ihara
       avoidance}Let~$F$ be a CM field, and let~$n\ge 2$. %
     Fix  %
     a continuous
     representation \[\rhobar:G_F\to\GL_n(\Qtwobar)\]satisfying the
     following hypotheses.
     \begin{enumerate}
     \item\label{Un residually automorphic} There is an ordinary RACSDC automorphic
       representation~$\pi$ of 
       $\GL_n/F$ such that
       $\rbar_{\pi,2}\cong\rhobar$.%
\item\label{primes split in F}$F/F^+$ is everywhere unramified. All of
  the places~$v|2$ of~$F^+$ split in~$F$, as do  all places lying
  under a place at which~$\pi$ is ramified.
\item \label{existence of unitary group} $n[F^+:\Q]\equiv 0\pmod{4}$.
     \item\label{Un nearly adequate} $\rhobar(G_F)$ is nearly adequate.
     \item\label{Un sacrificial prime} $\rhobar(G_F)$ contains a regular semi-simple element.
     \item\label{Un strongly residually odd} If $n$ is even, then there exists an infinite place~$v$
       of~$F^+$ such that the polarized pair %
       $(\rhobar,\varepsilonbar^{1-n}\delta_{F/F^+}^n)$ is strongly
       residually odd at~$v$.
     \end{enumerate}%
     Let~$T$ be any finite set
     of finite places of~$F^+$ which split in~$F$, which contains all finite places~ lying under a place~$w$ of~$F$
     at which~$\pi_w$ is ramified, and all places dividing~$2$.
     
     Then ~ $R^{T,\ord}$ is a finite $\Lambda$-algebra. Furthermore,
     if for each~$v\in T\setminus S_2$ the group~$U_v$ is sufficiently
deep in the sense of Definition~\ref{defn: sufficiently deep subgroup Un},
     then the morphism
$R^{T,\ord}\to\T^{T,
      \ord}(U(\mathfrak{l}^\infty))_\m$ given by~\eqref{eqn:
  R to T map Un} has nilpotent kernel, i.e.\
\[(R^{T,\ord})^{\red}\isoto \T^{T,
      \ord}(U(\mathfrak{l}^\infty))_\m.\]
   \end{thm}
    \begin{rem}
     \label{rem: regular semi-simple} The assumption in Theorem~\ref{thm: U(n) ordinary automorphy lifting including Ihara avoidance}~\eqref{Un
       sacrificial prime}  that
     $\rhobar(G_F)$ contains a regular semi-simple element is used in
     order to ensure that the automorphic forms that we consider are
     of neat level.%
   \end{rem}

   \begin{proof}[Proof of Theorem~\ref{thm: U(n) ordinary automorphy
       lifting including Ihara avoidance}]
     We will begin by making a succession of solvable extensions of
     CM fields to put ourselves into a situation where we can apply
     the Taylor--Wiles patching method. In order to keep the notation
     compatible with that above we will continue to denote our CM
     field by~$F$ until the end of the argument, where we will descend
     to our original~$F$.

     We can
     and do replace~$F$ with a solvable extension (and replace~$T$ with
     the set of places lying over places in~$T$) and enlarge our coefficient field~$E$ so that in addition
     to the hypotheses of the theorem, we have:
     \begin{itemize}
     \item $F=F^+(\sqrt{-1})$.
     \item if~$v\in T$ then $\rhobar|_{G_{F_\tv}}$ is trivial.
     \item if~$v\in S_2$ then $[F_v^+:\Qtwo]>n(n-1)/2+1$.
     \item if~$v\in T\setminus S_2$ then:
       \begin{itemize}
       \item ~$\pi_\tv^{\Iw(\tv)}\ne 0$.
       \item  if $2^N \| (q_v - 1)$ then $2^N > n$ and $\cO$ contains
         a $2^N$th root of unity.
              \end{itemize}
     \end{itemize}
     (Note that if ~$(F')^+/F^+$ is a solvable extension of totally
     real fields then~$(F')^+F/F$ is a solvable extension of CM fields, so we
     can choose a solvable CM extension to realize any finite set of
     local extensions. All of these conditions are local except for
     the first condition that $F=F^+(\sqrt{-1})$. Since arranging this
     only involves a quadratic extension, and $\rhobar(G_F)$ has no
     normal subgroups of index~$2$ by the assumption that it is nearly
     adequate (which requires in particular that
     $H^1(\rhobar(G_F),k)=0$), this quadratic extension
     leaves~$\rhobar(G_F)$ unchanged.)  Choose a finite
     place~$v_1\notin T$ of~$F^+$ which splits in~$F$
     as~$v=\tv_1\tv_1^c$, for which
     $\rhobar(\Frob_{\tv_1})$ is regular semi-simple. (There are
     infinitely many such~$v_1$ by our assumption~\eqref{Un
       sacrificial prime}.)

     We replace~$T$ by $T\cup\{v_1\}$, and
     write $T=S_2\disjointunion \{v_1\}\disjointunion R$. Note in particular that
     Hypotheses~\ref{hyp: existence of quasisplit unitary group}
     and~\ref{hyp: T is big enough} hold. For each~$v\in R$ we choose pairwise distinct characters
     $\chi_{v,1}, \dots, \chi_{v,n}: \cO^\times_{F_\tv} \rightarrow
     \cO^\times$ which become trivial after reduction modulo
     $\varpi$. (We can do this by the conditions arranged in our
     initial base change.) We have the global deformation problem
     $\cS_{\{\chi_v\}}$ defined in~\eqref{eq: defn of S chi Un}, and
     we write~  $\cS_{\{1\}}$ for the global deformation problem
    defined in the same way but with all of the
     characters~$\chi_{v,i}$ replaced by the trivial character. By the
     definitions of the local deformation problems~$\cD_v^{\chi_v}$
     for~ $v\in R$, we can fix compatible isomorphisms $R_{\cS_{\{\chi_v\}},T}^{\loc} / \varpi \cong R_{\cS_{\{1\}},T}^{\loc} / \varpi$ and  $R_{\cS_{\{\chi_v\}}} / \varpi \cong R_{\cS_{\{1\}}} / \varpi$.

We now specify  open compact subgroups $U_v\subset G(F_v^+)$ as follows:
\begin{enumerate}%
 \item $U_v = G(\cO_{F^+_v})$ if $v \not\in T$ is split in $F$.
\item $U_v$ is a hyperspecial maximal compact  subgroup of $G(F^+_v)$ if $v$ is inert in $F$.
\item $U_v = \Iw(\tv)$ for $v \in R$.
\item $U_{v_1}$ is any torsion-free  compact open subgroup
  of~$G(F^+_{v_1})$. 
\end{enumerate}
(Note in particular that the choice of~$U_{v_1}$ means that for
any~$b,c$, the group~$U(\mathfrak{l}^{b,c})$ is sufficiently small in
the sense of Definition~\ref{defn: suff small Un}. While we will not
explicitly use this below, it is implicitly used multiple times,
ultimately in the form of~\cite[Lem.\ 2.6]{ger}.)

Since by assumption~$\pi$ is ordinary and unramified outside
of~$R\cup S_2$, and since $\pi_\tv^{\Iw(\tv)}\ne 0$ for all~$v\in R$,
there is a maximal ideal~$\m_{1}$ of
$\T^{T, \ord}_{\{1\}}(U(\mathfrak{l}^\infty), \cO)$ with residue field~$k$
such that~$\rhobar\cong \rbar_{\m_{1}}|_{G_F}$. (As ever, we feel free to
enlarge~$\cO$ if necessary.) Since the~
$\chi_v$ are trivial modulo~$\varpi$, we have \numequation\label{eqn:
  Un ord spaces identified mod varpi}S^{\ord}_{0,
    \{\chi_v\}}(U(\mathfrak{l}^{\infty}), k)=S^{\ord}_{0,
    \{1\}}(U(\mathfrak{l}^{\infty}), k),\end{equation}
    so~ $\m_{1}$ induces a
unique maximal ideal~ $\m_{\chi}$ of~$\T^{T,
  \ord}_{\{\chi_v\}}(U(\mathfrak{l}^\infty), \cO)$. After conjugating
we can and do assume that $\rbar_{\m_{\chi}}=\rbar_{\m_{1}}=\rhobar$.

Write \[H_\chi:=S^{\ord}_{0,
    \{\chi_v\}}(U(\mathfrak{l}^{\infty}),
  E/\cO)^\vee_{\m_{\chi}},\quad H_1:=S^{\ord}_{0,
    \{1\}}(U(\mathfrak{l}^{\infty}),
  E/\cO)^\vee_{\m_{1}}.\] By~\cite[Cor.\ 2.21]{ger}, $H_1$
is a faithful $\T^{T, \ord}_{\{1\}}(U(\mathfrak{l}^\infty),
\cO)_{\m_{1}}$-module, and is in particular
an~$R_{\cS_{\{1\}}}$-module via~\eqref{eqn:
  R to T map Un}. 
Similarly $H_\chi$ is an~$R_{\cS_{\{\chi_v\}}}$-module. By~\eqref{eqn:
  Un ord spaces identified mod varpi} we have a natural
isomorphism \[H_\chi/\varpi\cong H_1/\varpi,\] which is
compatible with the isomorphism  $R_{\cS_{\{\chi_v\}}} / \varpi \cong
R_{\cS_{\{1\}}} / \varpi$.

Write~$q=h^1_{\cS^\perp,T}-1$      and $g=q+\#
     T-1-[F^+:\Q]n(n-1)/2$. Write  \[S_\infty:=\Lambda\llb
     X_1,\cdots,X_{q+(n^2+1)\# T-1}\rrb\] with augmentation ideal
   $\mathfrak{a}_\infty = (X_1, \dots, X_{q+(n^2+1)\# T-1})$. (The
   number of formal variables here is given by the number of
   Taylor--Wiles primes plus the relative dimension of~$R_{\cS}^T$
   over~$R_{\cS}$.) We
   set \[R_{\chi,\infty}:=R_{\cS_{\{\chi_v\}},T}^{\loc}\llb
         Y_1,\dots,Y_g\rrb,\quad
         R_{1,\infty}:=R_{\cS_{\{1\}},T}^{\loc}\llb
         Y_1,\dots,Y_g\rrb.\]

By~\eqref{eqn: dim of local Un deformation ring} we have \numalign\label{eqn: dimensions add up Un}\dim
R_{\chi,\infty}=\dim R_{1,\infty}&= 1+ n^2\#T+[F^+:\Q]n(n+1)/2+g \\
\nonumber %
     &= [F^+:\Q]n+(n^2+1)\#T +q \\ \nonumber &=
\dim S_\infty.\end{align}
       
     Using Proposition~\ref{prop:
  existence of TW systems} in place of~\cite[Prop.\ 4.4]{jack}, a
standard patching argument exactly as in the proof of~\cite[Thm.\
8.6]{jack} provides us with the following:

\begin{itemize}
	\item $\CNL_\Lambda$-homomorphisms $S_\infty\to R_{1,\infty}$,
          $S_\infty\to R_{\chi,\infty}$.
        \item An $R_{1,\infty}$-module $H_{1,\infty}$, and an
          $R_{\chi,\infty}$-module $H_{\chi,\infty}$, each of which is
          free of finite rank over~$S_\infty$.
        \item A surjection of $R_{\cS_{\{1\}},T}^{\loc}$-algebras
          $R_{1,\infty}\onto R_{\cS_{\{1\}}}$, which factors through
          a $\Lambda$-algebra map
          $R_{1,\infty}/\mathfrak{a}_{\infty}\to R_{\cS_{\{1\}}}$; and
          similarly, a surjection of
          $R_{\cS_{\{\chi_v\}},T}^{\loc}$-algebras
          $R_{\chi,\infty}\onto R_{\cS_{\{\chi_v\}}}$, which factors
          through a $\Lambda$-algebra map
          $R_{\chi,\infty}/\mathfrak{a}_{\infty}\to
          R_{\cS_{\{\chi_v\}}}$.
        \item Isomorphisms $H_{1,\infty}/\mathfrak{a}_\infty\cong
          H_1$, $H_{\chi,\infty}/\mathfrak{a}_\infty\cong
          H_\chi$ compatible with the surjections $R_{1,\infty}\onto
          R_{\cS_{\{1\}}}$, $R_{\chi,\infty}\onto R_{\cS_{\{\chi_v\}}}$.
		\item Compatible identifications of all the above data for~$1$ and
          for~$\chi$ after reducing modulo~$\varpi$.
\end{itemize}

In particular since $H_{1,\infty}$ is a finite free $S_\infty$-module,
we deduce from~\eqref{eqn: dimensions add up Un}
that \[\depth_{R_{1,\infty}}H_{1,\infty}\ge
  \depth_{S_\infty}H_{1,\infty}=\dim S_\infty=\dim R_{1,\infty},\]
whence $\depth_{R_{1,\infty}}H_{1,\infty}=\dim R_{1,\infty}$, and the
support of~$H_{1,\infty}$ in $\Spec R_{1,\infty}$ is a union of
irreducible components (see~\cite[Lem.\ 2.3]{tay}). Similarly, the support of $H_{\chi,\infty}$ in $\Spec R_{\chi,\infty}$ is a union of
irreducible components.

We now examine the irreducible components of $\Spec R_{\chi,\infty}$
and~$\Spec R_{1,\infty}$. Bearing in mind
Propositions~\ref{levelraising} and~\ref{prop: ord local def rings
  Un}, an identical argument to the proof of~\cite[Lem.\
3.2.4]{BCGNT-Bianchi} shows that for each minimal prime~$Q$
of~$\Lambda$, we have the following properties.

\begin{enumerate}
  \item The generic points of $R_{\chi,\infty}/Q$ and $R_{1,\infty}/Q$ all
  have characteristic~$0$.
\item The irreducible components~$\cC$ of $\Spec R_{\chi,\infty}/Q$ and
  $\Spec R_{1,\infty}/Q$ biject with the products of the
  corresponding sets of irreducible components~$\cC_v$ of the local
  deformation rings for~$v\in R\cup\{v_1\}$.
\item The irreducible components~$\overline{\cC}$ of $\Spec R_{\chi,\infty}/(Q,\varpi)=\Spec R_{1,\infty}/(Q,\varpi)$ biject with the products
  over~$T$ of the corresponding sets of the irreducible components~$\overline{\cC}_v$ of
  the special fibres of the deformation rings for~$v\in R\cup\{v_1\}$.
\end{enumerate}
In view of these statements we will use the notation $\cC=\otimes_{v\in R\cup\{v_1\}}\cC_v$ and $\cCbar=\otimes_{v\in R\cup\{v_1\}}\cCbar_v$.
  
\begin{enumerate}[resume]
\item\label{item: R chi irred Un components} The irreducible components of
  $\Spec R_{\chi,\infty}/Q$  biject with the irreducible components of~$\Spec R_{v_1}^{\square}$.
\item \label{item: R 1 Un components lifting} For each irreducible
  component $\cCbar=\otimes_{v\in R\cup\{v_1\}}\cCbar_v$ of
  $\Spec R_{1,\infty}/(Q,\varpi)$, there are irreducible
  components~$\cC_v$ of~$\Spec R_v^1$ for~$v\in R$, and irreducible
  components ~$\cC_{v_1}^1,\dots\cC_{v_1}^s$ of
  ~$\Spec R_{v_1}^\square$ (for some~$s\ge 1$), such that the
  irreducible components of $\Spec R_{1,\infty}/Q$
  generalizing~$\cCbar$ are precisely the~$s$ components
  $\cC_{v_1}^i \otimes_{v\in R}\cC_v$.
\end{enumerate}

Fix for the moment a minimal prime~$Q$ of~$\Lambda$. The existence of~$\pi$
implies that  the support of $H_{1,\infty}$ in $\Spec R_{1,\infty}$ is
nonempty. Using the comparison modulo~$\varpi$, the same is true of
the support of $H_{\chi,\infty}$ in $\Spec R_{\chi,\infty}$. By points
~\eqref{item: R chi irred Un components} and ~\eqref{item: R 1 Un components
  lifting}, we conclude that for each set~$X$  of irreducible
components~$\cC_v$ of~$\Spec R_v^1$ for~$v\in R$, we can choose an
irreducible component $\cC_{X,v_1}$ of ~$\Spec R_{v_1}^\square$ such
that  $\cC_{X,v_1} \otimes_{v\in R, \cC_v\in X}\cC_v$ is in the
support of~ $H_{1,\infty}$. This choice of irreducible components
corresponds to a quotient~$R^{\loc}_{X}$  of~$R_{\cS_{\{1\}},T}^{\loc}$, and if we
set \[R_{X}:=R_{\cS_{\{1\}}}\otimes_{R^{\loc}_{\cS_{\{1\}},T}}R^{\loc}_{X},\]
then $\Spec R_{X}\subset \Spec R_{\cS_{\{1\}}}$ is contained in
the support of~$H_{1,\infty}'$; equivalently, $\Spec R_{X}$ is contained
in the support of~$H_1$ in  $\Spec R_{\cS_{\{1\}}}$.

In particular~$(R_X)^{\red}$ is a quotient of a Hecke algebra
$\T^{T, \ord}_{\{1\}}(U(\mathfrak{l}^\infty), \cO)_{\m_{1}}$, so it is
a finite~$\Lambda$-algebra, so that~$R_X$ itself is a
finite~$\Lambda$-algebra. (To see this, note by the topological
Nakayama lemma it suffices to observe $R_X/\m_{\Lambda}$ is Noetherian
and zero-dimensional, thus finite.) Since~$R_X$ has dimension at
least~$\dim\Lambda$ by Lemma~\ref{lem: Un deformation ring dimension
  lower bound}, we see that the morphism $\Spec R_X\to\Spec \Lambda/Q$
is dominant. In particular, we can choose a weight~$\lambda$ such that
$R_X/\wp_\lambda R_X$ is a nonzero finite~$\cO$-algebra of
dimension at least~$1$, and thus has a~$\Qtwobar$-point. Since
~$(R_X)^{\red}$ is a quotient of
$\T^{T, \ord}_{\{1\}}(U(\mathfrak{l}^\infty), \cO)_{\m_{1}}$, it
follows from ~\eqref{eqn: Hida control Un} that the Galois representation
corresponding to this~$\Qtwobar$-point comes from an ordinary
RACSDC representation of weight~$\lambda$.

Repeating this construction for all choices of~$Q$, we conclude that
for each choice set~$Y$ of irreducible components~$\cC_v$ for~$v\in R\cup
S_2$, there is an ordinary RACSDC automorphic
representation~$\pi_Y$ of~$\GL_n/F$ such that:
\begin{itemize}
\item $\rbar_{\pi_Y,2}\cong\rhobar$,
\item $r_{\pi_Y,2}|_{G_{F_w}}$ is unramified for all places~$w$ not
  lying over a place in~$T$,
\item and for each~$v\in R\cup S_2$, the  representation
  $r_{\pi_Y,2}|_{G_{F_\tv}}$ lies on ~$\cC_v$ and on no other
  irreducible component (by the genericity of~$\pi_{Y,v}$,
  see~\cite[Lem.\ 1.3.2(1)]{BLGGT}).
\end{itemize}

By Lemma~\ref{lem: deformations of semi-simple are sums of characters},
we can and do choose a solvable CM extension $L/F$, linearly disjoint
from~$\overline{F}^{\ker\rhobar}$ over~$F$, with the following
property: for any~$Y$ as above, and any place~$w_1$ of~$L^+$ lying
over~$v_1$,  the representation
$r_{\pi_Y,2}|_{G_{F_{\tw_1}}}$ is
    unramified (where we write~$w_1=\tw_1\tw_1^c$).

We now repeat the patching argument above with $F$ replaced
by~$L$. More precisely, we:
\begin{itemize}
\item replace $R$ by the set~$R'$ of places of~$L$ lying over places
  in~$R$;
\item choose a place~$v_1'\notin R'\cup S_2$ of~$L^+$ splitting in~$L$
  as $\tv_1'(\tv_1')^c$,  with $\rhobar(\Frob_{\tv_1'})$ being regular
  semi-simple; 
\item and replace~$T$ by~$T'=R'\cup\{v_1'\}$. 
\end{itemize}
Writing~$H_1'$, $R_{1,\infty}'$ for the corresponding objects
over~$L$, we find in particular that we have the patched
module~$H_{1,\infty}'$, whose support in $\Spec R_{1,\infty}'$ is a
union of irreducible components. Again, we write these irreducible
components as $\cC'=\otimes_{v'\in T'}\cC'_{v'}$, and for each set~$Y$
as above we let~$\cC'_Y$ denote the irreducible component determined
by letting~$\cC'_{v_1'}$ be the (unique) unramified component of
$\Spec R_{v_1'}^\square$, and letting $\cC'_{v'}$ for
$v'|v\in R\cup S_2$ be the image of the component $\cC_v$ for~$Y$ (via
the natural morphism $\Spec R_v^\square\to\Spec R_{v'}^\square$).

By considering the base changes to~$L$ of the~$\pi_Y$, we see that
each component~$\cC'_Y$ is in the support of~$H_{1,\infty}'$.  The
union of the irreducible components~$\cC'_Y$ corresponds to a
quotient~$R^{\loc}_{L/F,T'}$ of~$R^{\loc}_{\cS_{\{1\}},T'}$, and as
above, if we
set
\[R'_{L/F,T'}:=R'_{\cS_{\{1\}}}\otimes_{R^{\loc}_{\cS'_{\{1\}},T'}}R^{\loc}_{L/F,T'},\]
then $\Spec R'_{L/F,T'}\subset \Spec R'_{\cS'_{\{1\}}}$ is contained
in the support of~$H_{1,\infty}'$; equivalently, $\Spec R'_{L/F,T'}$
is contained in the support of~$H'_1$ in $\Spec
R'_{\cS'_{\{1\}}}$. Thus~$(R'_{L/F,T'})^{\red}$ is a quotient of a
Hecke algebra
$\T^{T', \ord}_{\{1\}}(U(\mathfrak{l}^\infty), \cO)_{\m_{1}}$, and in
particular every homomorphism $R'_{L/F,T'}\to\cO$ factors
through~$\T^{T', \ord}_{\{1\}}(U(\mathfrak{l}^\infty),
\cO)_{\m_{1}}$. Furthermore, it follows as above that~$R'_{L/F,T'}$
 is a finite~$\Lambda'$-algebra.

We now return to the original situation of the statement of the
theorem (so~$F$ is now the CM field that we started with, before we
made any base changes, and~$T$ is as in the statement of the
theorem).  By the choice of~$L$, we have a commutative
diagram \[\xymatrix{\Lambda'\ar[r]\ar[d] & R'_{L/F,T'} \ar[d]  \\
    \Lambda \ar[r] & R^{T,\ord} \ar[r] & \T^{T,
      \ord}(U(\mathfrak{l}^\infty))_\m }\] The
morphism~$R'_{L/F,T'}\to R^{T,\ord}$ is finite, by the obvious
generalization of \cite[Lem.\ 1.2.3(1)]{BLGGT} to the case~$p=2$,
which has an identical proof up to replacing the appeal to~\cite[Lem.\
2.1.12]{cht} with a citation of~\cite[Lem.\ 2.4]{MR3598803}. It
follows that $R^{T,\ord}$ is a finite~$\Lambda$-algebra, as
claimed.

Suppose now that the~$U_v$ for~$v\in T\setminus S_2$ are sufficiently
deep. Since every irreducible component of~$\Spec R^{T,\ord}$ has dimension
at least that of~$\Spec \Lambda$ by Lemma~\ref{lem: Un deformation ring
  dimension lower bound}, it follows that each irreducible component
dominates an irreducible component of~$\Spec \Lambda$. It follows that
the set of points $\Spec \Qtwobar\to \Spec R^{T,\ord}$ which lie over
points of~$\Lambda$ given by the primes~$\wp_\lambda$ is dense
in~$\Spec R^{T,\ord}$. It remains to show that each such point is in
$\Spec \T^{T, \ord}(U(\mathfrak{l}^\infty))_\m$. By~\eqref{eqn: Hida
  control Un} and the choice of~$U$, it is enough to check that the
corresponding Galois representations are automorphic. By solvable base
change, it is enough to check this after restriction to~$G_L$, where
it follows from another application of~\eqref{eqn: Hida control Un}
(and the observation above that ~$(R'_{L/F,T'})^{\red}$ is a quotient
of $\T^{T', \ord}_{\{1\}}(U(\mathfrak{l}^\infty), \cO)_{\m_{1}}$).
   \end{proof}

\section{Ordinary modularity lifting theorems for $\GSp_4$:
  preliminaries}\label{sec: R=T GSp4 preliminaries all p}
  
  The goal of this section is --- in part ---  to prove
an ordinary modularity lifting
theorem for~$\GSp_4$ for~$p \ge 3$ over totally real fields
in which~$p$ splits completely. Under suitable Taylor--Wiles hypotheses,
this can be used to show that $p$-adic Galois representations
coming from ordinary abelian surfaces give rise to quotients
of a certain $p$-adic Hecke algebra, but does not yet show that such classes are classical.
The modularity lifting theorem we prove in this section is 
 (in the language of~\cite{CG}) of~$l_0 = 0$
type rather than~$l_0 > 0$ type, and so is precisely
amenable to the usual Taylor--Wiles method.
Under a stronger hypothesis (that~$\rhobar$ is residually
$p$-distinguished) our results are actually directly contained
in~\cite{BCGP} (although that paper is generally concerned with the more subtle~$l_0 = 1$ situation),
and versions of this theorem go back as far as~\cite{MR2920881}. The methods we use here follow along generally
similar grounds, with some important technical improvements
due in several cases to Whitmore~\cite{Whitmore}.

In~\S\ref{subsec: notation and defns
  Galois}, we recall some general constructions and notation
for~$\GSp_4$-deformation problems. 
In~\S\ref{subsec: ord def
  rings}, we introduce the corresponding ordinary local
deformation rings and study their local properties. %
In~\S\ref{sec: Gsp4
  modularity lifting},  we carry out the Taylor--Wiles argument
(in part following~\cite{BCGP} and~\cite{Whitmore}). Finally, in~\S\ref{sec:subgroupsofGSp},
we explicitly analyze the subgroups of~$\GSp_4(\F_3)$ which
satisfy our running collection of ``big--image'' conditions.

\subsection{Notation and definitions}\label{subsec: notation and defns
  Galois}We now turn to modularity lifting theorems for~$\GSp_4$. Our
arguments have relatively little \emph{direct} overlap with those of our earlier paper
~\cite{BCGP}, although we will occasionally make references to
it. In particular, in order to avoid confusing clashes of notation with our results
for unitary groups in Section~\ref{sec: R=T unitary}, we continue to write~$F^+$
for a totally real field (whereas totally real fields were denoted~$F$
in~\cite{BCGP}).

Accordingly we let~$F^+$ denote a totally real field, and write $S_p$ for the set of places of
$F^+$ above $p$. %
We fix a continuous absolutely
irreducible homomorphism $\rhobar: G_{F^+} \to \GSp_4(k)$ with
similitude ~$\varepsilonbar^{-1}$.
When~$\rhobar$ is explicitly considered as a symplectic representation 
(as in this section), we
denote by~$\ad \rhobar$ and~$\ad^0 \rhobar$ the adjoint~$G_{F^{+}}$
action with respect to~$\GSp_4$ and~$\Sp_4$ respectively (so~$\dim \ad \rhobar = 11$
and~$\dim \ad^0 \rhobar = 10$). We warn the reader that  there is some tension in
this definition
with the notation of~\S\ref{sec: R=T unitary} where~$\ad \rbar$ denotes the adjoint action
with respect to~$\GL_n$; we hope  the precise meaning will  always be  clear from context.

Let $S$ be a finite set of finite places of $F^+$ containing $S_p$
and
all places at which $\rhobar$ is ramified. We write $F^+_S$ for the
maximal subextension of $\overline{F^+} / F^+$ which is unramified outside
$S$, and write~$G_{F^+,S}$ for~$\Gal(F^+_S/F^+)$.  For each $v\in S$, we fix
$\Lambda_v \in \CNL_{\cO}$, and set
$\Lambda = \widehat{\otimes}_{v \in S} \Lambda_v$, where the completed
tensor product is taken over $\cO$.  Then $\CNL_\Lambda$ is a
subcategory of $\CNL_{\Lambda_v}$ for each $v\in S$, via the canonical
map $\Lambda_v \to \Lambda$.

\begin{defn}\label{defn: local lift}
  A \emph{lift}, also called a \emph{lifting}, of $\rhobar|_{G_{F^+_v}}$
  is a continuous homomorphism
  $\rho: G_{F^+_v} \to \GSp_4(A)$ to a
  $\CNL_{\Lambda_v}$-algebra $A$ such that
  $\rho \bmod \frakm_A = \rhobar|_{G_{F^+_v}}$ and
  $\nu\circ\rho=\varepsilon^{-1}$.
\end{defn}

We let $\cD_v^\square$ %
denote the set-valued functor on $\CNL_{\Lambda_v}$ that sends $A$ to
the set of lifts of $\rhobar|_{G_{F^+_v}}$ to $A$.  This functor is
representable, and we denote the representing object by
$R_v^\square$. We can identify $\cD_v^{\square}(k[\epsilon])$ with the group of 
1-cocycles $Z^1(F_v,\ad^0 \rhobar)$ 
by associating a cocycle $\phi$ to
the lifting given
by \[\rho(\sigma)=(1+\epsilon\phi(\sigma))\rhobar(\sigma).\] Note that
two such liftings are $\widehat{\GSp_4}(k[\epsilon])$-conjugate if and only if the
images of the corresponding 1-cocycles in $H^1(F_v,\ad \rhobar)$ are equal.

\begin{defn}\label{def:locdefprob}
 A \emph{local deformation problem} for $\rhobar|_{G_{F^+_v}}$ is a subfunctor $\cD_v$ of $\cD_v^\square$ satisfying the following:
  \begin{itemize}
   \item $\cD_v$ is represented by a quotient $R_v$ of $R_v^\square$.
   \item For all $A \in \CNL_{\Lambda_v}$, $\rho \in \cD_v(A)$, and $a \in \widehat{\GSp}_4(A)$, 
   we have $a\rho a^{-1} \in \cD_v(A)$. 
 \end{itemize}
\end{defn}

\begin{defn}\label{def:globdefprob}%
 A \emph{global deformation problem} is a tuple
  \[
   \cS = (\rhobar, S, \{\Lambda_v\}_{v\in S},
   \{\cD_v\}_{v\in S})
  \]
 where:
 \begin{itemize}
  \item $\rhobar$, $S$, $\{\Lambda_v\}_{v\in S}$  are as
    above. %
  \item For each $v\in S$, $\cD_v$ is a local deformation problem  for $\rhobar|_{G_{F^+_v}}$.
  \end{itemize}
\end{defn}

As in the local case, a \emph{lift} (or \emph{lifting}) of $\rhobar$
is a continuous homomorphism $\rho: G_{F^+,S} \to \GSp_4(A)$
to a $\CNL_\Lambda$-algebra $A$, such that
$\rho \bmod \frakm_A = \rhobar$ and $\nu\circ\rho=\varepsilon^{-1}$.   We say that
two lifts $\rho_1,\rho_2: G_{F^+,S} \to \GSp_4(A)$ are
\emph{strictly equivalent} if there is an
$a\in \widehat{\GSp}_4(A)$ such that
$\rho_2 = a\rho_1 a^{-1}$.  A \emph{deformation of $\rhobar$} is a
strict equivalence class of lifts of $\rhobar$. %

For a global deformation problem 
  \[    \cS = (\rhobar, S, \{\Lambda_v\}_{v\in S},
   \{\cD_v\}_{v\in S})\]%
  we say that a lift $\rho: G_{F^+,S} \to \GSp_4(A)$ is of \emph{type $\cS$}
  if $\rho|_{G_{F^+_v}} \in \cD_v(A)$ for each $v\in S$. %
  If $\rho_1$ and $\rho_2$ are strictly equivalent
  lifts of $\rhobar$, and $\rho_1$ is of type $\cS$, then so is
  $\rho_2$.  A \emph{deformation of type $\cS$} is a strict
  equivalence class of lifts of type $\cS$, and we denote by
  $\cD_{\cS}$ the set-valued functor that takes a
  $\CNL_\Lambda$-algebra $A$ to the set of lifts
  $\rho: G_{F^+} \to \GSp_4(A)$ of type $\cS$.

Given a subset $T\subseteq S$, %
  a \emph{$T$-framed lift of type
    $\cS$} is a tuple $(\rho,\{\gamma_v\}_{v\in T})$, where $\rho$ is
  a lift of type $\cS$, and
  $\gamma_v \in \widehat{\GSp}_4(A)$ for each $v\in T$.  We
  say that two $T$-framed lifts $(\rho_1,\{\gamma_v\}_{v\in T})$ and
  $(\rho_2,\{\gamma'_v\}_{v\in T})$ to a $\CNL_\Lambda$-algebra $A$ are
  strictly equivalent if there is an $a\in \widehat{\GSp}_4(A)$
  such that $\rho_2 = a \rho_1 a^{-1}$, and $\gamma'_v = a\gamma_v$ for
  each $v\in T$.  A strict equivalence class of $T$-framed lifts of
  type $\cS$ is called a \emph{$T$-framed deformation of type $\cS$}.
  We denote by $\cD_{\cS}^T$ the set valued functor that sends a
  $\CNL_\Lambda$-algebra $A$ to the set of $T$-framed deformations to
  $A$ of type $\cS$. 

  The functors $\cD_{\cS}$,  $\cD_{\cS}^T$  are representable (as we are assuming
  that~$\rhobar$ is absolutely irreducible), and we denote their
  representing objects by $R_{\cS}$ and $R_{\cS}^{T}$
  respectively. Assume now that~$T$ is chosen so that $\Lambda_v =
\cO$ for all $v\in S \smallsetminus T$. Write~$R_v$ for the
representing object of~$\cD_v$,
and define $R_{\cS,T}^{\loc} = \widehat{\otimes}_{v\in T} R_v$, with the completed tensor product being taken over $\cO$. 
It is canonically a $\Lambda$-algebra, via the canonical isomorphism $\widehat{\otimes}_{v\in T} \Lambda_v \cong \widehat{\otimes}_{v\in S} \Lambda_v$. 
For each $v\in T$, the natural transformation $\Def_{\cS}^T \rightarrow \cD_v$ given by $(\rho,\{\alpha_v\}_{v\in T}) \mapsto \alpha_v^{-1}\rho|_{G_{F_v}} \alpha_v$ 
induces a  morphism $R_v \rightarrow R_{\cS}^T$ in $\CNL_{\Lambda_v}$. 
We thus have a morphism $R_{\cS,T}^{\loc} \rightarrow R_{\cS}^T$ in $\CNL_{\Lambda}$.

  If~$T$ is
  empty, then $R_{\cS}=R_{\cS}^T$, and otherwise the natural
  map~$R_{\cS}\to R_{\cS}^T$ is formally smooth of relative
  dimension~$11\#T-1$.  Indeed $\cD_{\cS}^T \rightarrow \cD_{\cS}$ is
  a torsor under
  $(\prod_{v \in T} \widehat{\GSp}_4)/\widehat{\mathbb{G}}_m$.
  \begin{defn}\label{defn: cT}
    Let 
    \[\cT:=\Lambda\llb x_1,\dots,x_{11\#T-1}\rrb\]
     be the coordinate ring of
    $(\prod_{v \in T} \widehat{\GSp}_4)/\widehat{\mathbb{G}}_m$ over
    $\Lambda$. %
  \end{defn}
  The choice of a representative
    $\rho_{\cS} \colon G_F \rightarrow \GSp_4(R_{\cS})$ for the
    universal type $\cS$ deformation determines  a splitting of the torsor $\cD_{\cS}^T \rightarrow \cD_{\cS}$  and   a canonical
    isomorphism \numequation\label{eqn: framed deformation over unframed}R_{\cS}^T \cong R_{\cS} \widehat{\otimes}_\Lambda
    \cT.\end{equation}

The following lemma and its proof are standard, but we include them in
order to reassure the reader that they remain valid for~$p=2$.
\begin{lem}\label{lem:anyliftisTW}Suppose that $q_v \equiv 1 \bmod p$, and that 
$\rhobar|_{G_{F^+_v}}$ is unramified, with  $\rhobar(\Frob_v)$ being regular
semi-simple with \emph{(}ordered\emph{)} eigenvalues $\alphabar_{v,1}$, $\alphabar_{v,2}$,
$\alphabar_{v,2}^{-1}$, $\alphabar_4=\alphabar_{v,1}^{-1}$. Let $\rho: G_{F^+_v} \rightarrow \GSp_4(A)$ be any lift of
$\rhobar$.

Then there  are unique continuous characters
  $\gamma_i: G_{F^+_v} \rightarrow A^\times$ for~$i=1,2$, such that
  $\rho$ is $\GSp_4(A)$-conjugate to a lift of the form
\[\gamma_1 \oplus \gamma_2 \oplus \gamma_2^{-1}\varepsilon^{-1}
  \oplus\gamma_1^{-1}\varepsilon^{-1},\] where
  $(\gamma_i \bmod \frakm_A)(\Frob_{v}) = \alphabar_{v,i}$ for each
  $i=1,2$.
\end{lem}
\begin{proof}Let~$\phi$ be a lifting of~$\Frob_v$
  to~$G_{F^+_v}$. Then~$\rho(\phi)$ is regular semi-simple, so is
  contained in~$T(A)$ for a unique torus~$T$, and we need to show that
  for each~$\sigma\in I_{F^+_v}$, we have $\rho(\sigma)\in T(A)$. To
  do this we will prove by induction that for each~$n\ge 1$ we
  have \[(\rho \bmod \frakm_A^n)(\sigma)\in T(A/\m_A^n), \]the
  case~$n=1$ being true by our hypotheses.

For the inductive step, we assume the result holds for~$n$ and deduce
it for~$n+1$; replacing~$A$ by~$A/\m_A^{n+1}$, we may furthermore
assume that~$\m_A^{n+1}=0$. By the inductive hypothesis (and the
smoothness of~$T$) we can write~$\rho(\sigma)=tu$ where~$t\in T(A)$
and~$u\equiv 1\pmod{\m^{n}}$. Since~$\rhobar(\sigma)=1$, we see
that~$u$ and~$t$ necessarily commute (as do~$\rho(\phi)$ and~$t$, as
both are contained in~$T(A)$).

Now, since~$\rho$ is tamely ramified, we have
$\rho(\phi)\rho(\sigma)\rho(\phi)^{-1}=\rho(\sigma)^{q_v}$. Since~$q_{v}\equiv
1\pmod{p}$, and~$u\equiv 1\pmod{\m^n}$, we see that $u^{q_v}=u$, and
thus that
$\rho(\phi)u\rho(\phi)^{-1}=t^{q_v-1}u$. Using again that~$u\equiv
1\pmod{\m^n}$ and that~$\rhobar(\phi)$ is
regular semi-simple, it follows that~$u=1$, as required.
\end{proof}
  
\subsection{Ordinary deformation rings}\label{subsec: ord def
  rings}

In this section we study some ordinary deformation rings for $\GSp_4$.  We assume that $v$ is a place of $F^+$ lying over $p$ such that $F^+_v=\Q_p$. 
Write
$\Lambda_{\GSp_4,v}=\cO\llb(\cO_{F^+_v}^\times(p))^2\rrb$, where
$\cO_{F^+_v}^\times(p)$ denotes the pro-$p$
completion of $\cO_{F^+_v}^\times$.  There is a canonical character
$I_{F^+_v}\to\cO_{F^+_v}^\times(p)$ given by~$\Art_{F^+_v}^{-1}$, and we
define a pair of characters~$\theta_{i}:I_{F^+_v}\to \Lambda_{\GSp_4,v}$,
$i=1,2$ by letting~$\theta_{i}$ correspond to the embedding
~$\cO_{F^+_v}^\times(p) \rightarrow (\cO_{F^+_v}^\times(p))^2$ given by the $i$th
copy.   When $p>2$, $\Lambda_{\GSp_4,v}=\cO\llb x_1,x_2\rrb$ is formally smooth, while when $p=2$, $\Spec\Lambda_{\GSp_4,v}$ has 4 irreducible components but the generic fiber is regular.

Assume that~$\rhobar|_{G_{F^{+_v}}}$ is ordinary, and fix a $p$-stabilization
~$(\chibar_1,\chibar_2)$ of~$\rhobar|_{G_{F^{+_v}}}$, so that $(\chibar_1 ,\chibar_2 )$ is an ordered pair of %
characters
$G_{F^+_v}\to k^\times$.  Then
for any $A\in  \CNL_{\cO}$ there is an obvious bijection  between homomorphisms
$\Lambda_{\GSp_4,v}\to A$ %
and ordered pairs of characters
$(\chi_1,\chi_2):I_{F^+_v}\to A^\times$ lifting~$(\chibar_1,\chibar_2)$,
given by multiplying the characters $(\theta_1,\theta_2)$ by the
Teichm\"uller lifts of $(\chibar_1,\chibar_2)$.

Similarly, write $\Lambdat_{\GSp_4,v}=\cO\llb \Gal(F^{+,\ab}_v/F^+_v)(p)^2\rrb$,
 where 
$\Gal(F^{+,\ab}_v/F^+_v)(p)$  is the
pro-$p$ completion of $\Gal(F^{+,\ab}_v/F^+_v)$. Then we have a universal
pair of characters $(\chit_1,\chit_2):G_{F^+_v}\to\Lambdat_{\GSp_4,v}$
lifting~$(\chibar_1,\chibar_2)$.

We now introduce the ordinary deformation ring we consider, following \cite[\S 3]{ger}.  Let~$\cF$ denote the flag variety for~$\GSp_4$ over~$\cO$, i.e.\ the
variety whose $S$-points, for any $S/\Spec\cO$, parameterize full
flags \[
     0 = \Fil_0 \subset \Fil_1 \subset \cdots \subset \Fil_4 = \cO_S^4
    \]  with~$\Fil_i$ being locally free of rank~$i$ and locally a
    direct summand, with the further property
    that~$(\Fil_i)^\perp=\Fil_{4-i}$ for each~$i$ (where ~${}^\perp$
    is with respect to our usual symplectic form
    on~$\cO_S^4$). %

    Write~$R_v^{\square}$ for the ring
    denoted~$R_v^\square$ in Section~\ref{subsec: notation and defns
      Galois} when~$\Lambda_v=\Lambdat_{\GSp_4,v}$. One shows as in the proof
    of~\cite[Lem.\ 3.2]{ger} that there is a closed subscheme
    $\cG_v$ of
    $\cF\times_{\cO}\Spec R_{v}^{\square}$, such
    that for any $\cO$-algebra~$A$, the $A$-points of
    $\cG_v$ are exactly the pairs
    $(\Fil_\bullet,R_{v}^{\square}\to A)$ consisting
    of a symplectic flag~$\Fil_\bullet$ on~$A^4$ and an $\cO$-algebra morphism
    $R_{v}^{\square}\to A$, such that the
    pushforward of the universal lifting over
    $R_{v}^{\square}\to A$ preserves~$\Fil_\bullet$,
    and for %
        $i=1,2,3,4$ the action of~$G_{F^+_v}$ on
    $\Fil_i/\Fil_{i-1}$ is via respectively (the pushforwards to~$A$
    of) the characters $\chit_1$, $\chit_2$,
    $\varepsilon^{-1}\chit_2^{-1}$, $\varepsilon^{-1}\chit_1^{-1}$.

We
write \[R_{v}^{\triangle}:=\im\left(R_{v}^{\square}\to \cO_{\cG_{v}}(\cG_{v})
  \right),\] so that $\Spec R_{v}^{\triangle}$ is the
scheme-theoretic image of the morphism $\cG_{v}
\to \Spec R_{v}^{\square}$. (Note that here we differ from~\cite{ger} by not
passing to the $p$-torsion-free quotient.)  We denote by $\cD_v^{\triangle}$ the corresponding deformation problem.  %

Exactly as in the proof of~\cite[Lem.\ 3.3]{ger}, it follows
immediately from the properness of~$\cF$ that 
  if~$E'/E$ is a finite extension with ring of integers~$\cO_{E'}$
  and $\Lambdat_{\GSp_4,v}\to\cO_{E'}$ is a morphism of $\cO$-algebras, then the
  $\cO_{E'}$-points of $\Spf R_{v}^{\triangle}$ are exactly those
  lifts $\rho$ of~$\rhobar$ having the property that 
 there is a symplectic flag  \[
     0 = \Fil_0 \subset \Fil_1 \subset \cdots \subset \Fil_4 = \cO_{E'}^4
    \] as above %
such that for $i=1,2,3,4$  the action of~$G_{F^+_v}$ on
$\Fil_i/\Fil_{i-1}$ is via respectively the characters $\chit_1$,
$\chit_2$, $\varepsilon^{-1}\chit_2^{-1}$,
$\varepsilon^{-1}\chit_1^{-1}$. Equivalently, these are the lifts~$\rho$ which
are ordinary with $p$-stabilization $(\chit_1 ,\chit_2 )$ in the sense of Definition~\ref{defn:ordinary-for-Galois-reps}.

\begin{remark}As in Definition~\ref{defn:ordinary-for-Galois-reps}, we say that~$\rhobar$ is \emph{residually $p$-distinguished} if the 4 characters
$\chibar_1,\chibar_2,\varepsilon^{-1}\chibar_2^{-1},\varepsilon^{-1}\chibar_1^{-1}$
are pairwise distinct (if $p\not=2$ this amounts to
$\chibar_1\not=\chibar_2$). In this case the filtration $\Fil_i$ in the
definition of $\cG_v$ is uniquely determined by the Galois representation, and
it follows that  the map $\cG_v\to\Spec R_v^\square$ is a closed immersion. In \cite[\S 7.3]{BCGP}, we made this %
assumption and assumed $p>2$ and studied $\cG_v$ under the name $R_v^{B,\alphabetabar_v}$.
\end{remark}

    Let $E'/E$ be a finite extension, and let
    $x:\Spec E'\to \cG_{v}[1/p]$ be a closed
    point. Let~$\rho_x:G_{F^+_v}\to\GSp_4(E')$ be the pushforward of
    the universal lift coming from the composite
    $\Spec E'\to \cG_{v}[1/p]\to \Spec
    R_{v}^{\square}$. Let~ $\ad^0\rho_x$ denote the
    adjoint representation with respect to~$\Sp_4$, and define a decreasing
    filtration~$\Fil^i\ad^0\rho_x$ on~$\rho_x$
    by
    \[\Fil^i\ad^0\rho_x:=\{A\in\ad^0\rho_x |
      A\Fil_j\rho_x\subseteq\Fil_{j-i}\rho_x \forall j\};\]
      in
    particular, $\Fil^0\ad^0
    \rho_x$ is the subspace of $\ad^0\rho_x$ preserving the flag $\Fil_\bullet\rho_x$. 
    We have~$\dim \ad^0 \rho_x = 10$, and~$\dim \Fil^i \ad^0 \rho_x = 6,4,2,1,0$, for $i = 0,\ldots,4$ respectively.
  \begin{lem}     \leavevmode     \label{lem: criterion for H2 vanishing}
\begin{enumerate}
\item

 \label{item: H2 zero implies regular}If $H^2(G_{F^+_v},\Fil^0\ad^0
 \rho_x)=0$ then  $x$ is a regular point
  of $\cG_{v}[1/p]$; and $x$ is contained
  in a unique irreducible component of
  $\cG_{v}[1/p]$, and this component has dimension $16$.
\label{item: H2 vanishes in cG} We have \[H^2(G_{F^+_v},\Fil^0\ad^0
    \rho_x)= 0\] if and only if
     \[H^0(G_{F^+_v},\bigl(\ad^0
     \rho_x/\Fil^1\ad^0\rho_x\bigr)(1))= 0.\] 
 \item The equivalent conditions of part~\eqref{item: H2 vanishes in cG} hold under any of the following circumstances:
  \begin{enumerate}
  \item \label{item:easy vanishes}  none of the specializations at~$x$ of the
   characters $\chit_1^2\varepsilon$, $\chit_2^2\varepsilon$,
   $\chit_1\chit_2\varepsilon$, $\chit_1\chit_2^{-1}$ are equal
   to~$\varepsilon$.
 \item \label{item: pure p dist vanishes} $\rho_x$ is pure and
   $p$-distinguished. %
 \item\label{item: pure crystalline vanishes} $\rho_x$ is pure and
   potentially crystalline.
  \end{enumerate}  
  \item \label{item: def ring regular}If $\rho_x$ is $p$-distinguished and the equivalent conditions of part~\eqref{item: H2 vanishes in cG} hold (in
    particular, if $\rho_x$ is $p$-distinguished and pure) then the image of $x$ in $\Spec R_v^\triangle$ is a regular point which is contained in a unique irreducible component of $\Spec R_v^\triangle$, which has relative (over $\cO$) dimension 16.
\end{enumerate} 
    \end{lem}
    \begin{proof}  %
 By a
    standard tangent-obstruction calculation exactly as in the proof of~\cite[Lem.\
    3.7]{ger} %
    (see also~\cite[\S 5.1]{MR1643682} for the case of general
    algebraic groups), the tangent space to
    $\cG_{\Lambda_{\GSp_4,v}}[1/p]$ at~$x$ has dimension
    \numequation\label{eqn: tangent spaces on Geraghty resolution
      Lambda}16+\dim_{E'}H^2(G_{F^+_v},\Fil^0\ad^0
    \rho_x),\end{equation}
    and there is an obstruction class in
  $H^2(G_{F^+_v},\Fil^0\ad^0 \rho_x)$ whose vanishing implies that~ $x$ is a regular point
  of $\cG_{\Lambda_{\GSp_4,v}}[1/p]$.

For the remaining claim in~\eqref{item: H2 zero implies
      regular}, note that the space~$\ad^0
      \rho_x$ is self-dual under the trace pairing~$(A,B) \rightarrow
      \Tr(AB)$ (since we are in characteristic zero).
    If~$A \in \Fil^0\ad^0
    \rho_x$ and~$B \in \Fil^1\ad^0
    \rho_x$ then~$AB \in \Fil^1\ad^0
    \rho_x$ and so~$\Tr(AB) = 0$. It follows that~$\Fil^1\ad^0
    \rho_x \subset (\Fil^0\ad^0
    \rho_x)^{\perp}$. By considering the dimensions of these spaces ($6$ and~$10-4 = 6$ respectively), we deduce that
    this is an equality, and hence the claim  follows by Tate local
    duality.

    For part~\eqref{item:easy vanishes}, note that $\Fil^0\ad^0
    \rho_x$
        has a filtration with graded pieces of rank~$1$,
    and the characters %
    through which~$G_{F^+_v}$ acts on these graded
    pieces are  $1$, $1$, $\chit_1^2\varepsilon$,
    $\chit_2^2\varepsilon$, $\chit_1\chit_2\varepsilon$, $\chit_1\chit_2^{-1}$. 
Part~\ref{item: pure crystalline vanishes} follows from  part~\eqref{item:easy vanishes}
 because a potentially crystalline pure representation cannot contain
 two Jordan--H\"older factors differing by a cyclotomic twist.
 
We now turn to part ~\eqref{item: pure
   p dist vanishes}, so that $\rho_x$ is pure and $p$-distinguished, and (since
 we have just established part~\ref{item: pure crystalline vanishes}) we can furthermore assume that we are not potentially
 crystalline. Assume
for the sake of contradiction that $H^0(G_{F^+_v},\bigl(\ad^0
 \rho_x/\Fil^1\ad^0\rho_x\bigr)(1))\ne 0$. We now argue as in the proof
 of~\cite[Lem. 3.7(3)]{ger}. By our assumption, there is some $A\in \ad^0
 \rho_x- \Fil^1\ad^0\rho_x$ such that for all $\sigma\in
 G_{F^+_v}$, we
 have~$\rho_x(\sigma) A \rho_x(\sigma)^{-1} = \varepsilon(\sigma)^{-1}A \bmod
 \Fil^1\ad^0\rho_x$; equivalently,
  \numequation\label{eqn: explicit matrix condition for non smoothness}\varepsilon(\sigma)\rho_x(\sigma)A-A\rho_x(\sigma)\in
   \Fil^1\ad^0\rho_x.  \end{equation} We can and do conjugate $\rho_{x}$ so that
 it is contained in the usual upper triangular Borel subgroup, so that each
 $\Fil_{x,i}$ is generated by~$e_1,\dots,e_i$. Write \[\rho_x=\begin{pmatrix}
    \chit_{x,1}&*&*&*\\
0 &\chit_{x,2}&*&*\\
0&0&\chit_{x,3}&*\\
0&0&0&\chit_{x,4}\end{pmatrix}\] (so
$\chit_{x,3}=\varepsilon^{-1}\chit_{x,2}^{-1}$ and
$\chit_{x,4}=\varepsilon^{-1}\chit_{x,1}^{-1}$). Let $1\le s\le 4$ be
minimal with $Ae_s\notin\Fil_{x,s-1}$ (such an~$s$ exists by the
hypothesis that $A\notin \Fil^1\ad^0\rho_x$), and let $r\ge s$ be the
unique integer with $Ae_s\in\Fil_{x,r}-\Fil_{x,r-1}$. By~\eqref{eqn:
  explicit matrix condition for non smoothness} and the assumption
on~$s$, we see that for all $\sigma\in
 G_{F^+_v}$, we
 have \[\rho_x(\sigma)(Ae_s) \equiv
   \varepsilon(\sigma)^{-1}A\rho_x(\sigma)e_s \equiv
   (\varepsilon^{-1}\chit_{x,s})(\sigma)(Ae_s)\pmod{\Fil_{x,s-1}}.\]Since~$r\ge
 s$ this congruence in particular holds modulo $\Fil_{x,r-1}$, whence
 $\chit_{x,r}=\varepsilon^{-1}\chit_{x,s}$, and consequently~$r>s$;
 and furthermore we see that  $E\cdot (Ae_s)+\Fil_{x,s-1}$ is
 $G_{F^+_v}$-stable. More precisely, we see that the 2-dimensional subquotient
 $(E\cdot (Ae_s)+\Fil_{x,s})/\Fil_{x,s-1}$ of~$\rho_{x}$ is isomorphic
 to $\varepsilon^{-1}\chit_{x,s}\oplus \chit_{x,s}$. 

  It is immediate from the
 definition of purity that no twist of~$\varepsilon\oplus 1$ can be a
 subrepresentation of~$\rho_{x}$, so we must have~$s=2$
 or~$s=3$. Since~$\Fil_{\bullet}$ is symplectic, the possibility
 $(s,r)=(3,4)$ is also ruled out (because we already saw that $(s,r)=(1,2)$
 is impossible), so we must have~$s=2$ and~$r=3$ or~$4$. However, since we are
 pure and not potentially crystalline, we see in either case that we have
 $\chit_{x,1}=\chit_{x,2}$, which contradicts our
 assumption that~$\rho_{x}$ is $p$-distinguished.  So $H^0(G_{F^+_v},\bigl(\ad^0
 \rho_x/\Fil^1\ad^0\rho_x\bigr)(1))=0$ after all,
 as claimed.
 
 Part \eqref{item: def ring regular} follows immediately as when $\rho_x$ is $p$-distinguished, the map $\cG_v\to\Spec R_v^\triangle$ is an isomorphism in a neighborhood of $x$.
    \end{proof}

In the rest of this section, we assume $p>2$.  In \cite[Prop. 7.3.4]{BCGP}, we
showed (in a somewhat hands-on manner) that if~ $\rhobar$ is residually
$p$-distinguished, then  $\Spec R_v^\triangle[1/p]$ is irreducible.  This was
used in the proof of our modularity lifting theorem.  We expect that the same
holds in general, %
but we don't prove this.  Instead we explain a softer way to proceed. We prove a series of Lemmas which will be used in our modularity lifting theorems.

\begin{lem}
  \label{lem: good ordinary reduction points are on a single component}Assume that $(\rhobar\otimes\varepsilonbar)|_{G_{F^+_v}}$ is finite flat.  Then any ordinary pure weight 2 crystalline lift lies on
  a unique irreducible component of $\Spec R_v^\triangle$, which is moreover
  independent of the lift.  This component has relative dimension 16
  over~$\cO$.%
\end{lem}
\begin{proof}
We first show that a point $\rho$ of $\Spec R_v^\triangle$ corresponding to an ordinary pure weight 2 crystalline lift lies on a unique irreducible component (if $\rho$ were $p$-distinguished this was already part of Lemma \ref{lem: criterion for H2 vanishing} \eqref{item: def ring regular}).  Consider the fiber in $\cG_v$ over $\rho$ in $\Spec R_v^\triangle$, or in other words consider the space of $G_{F_v^+}$-stable symplectic filtrations $\{\Fil_i\}$ on $\rho$ on which $G_{F_v^+}$ acts on $\Fil_i/\Fil_{i-1}$ by $\chi_{1},\chi_{2},\varepsilon^{-1}\chi_{2}^{-1},\varepsilon^{-1}\chi_{1}^{-1}$ for $i=1,2,3,4$.  By assumption $\chi_{1}$, $\chi_{2}$ are unramified and $\rho$ has two dimensional inertia invariants, which hence must be $\Fil_2$.  Then either $G_{F_v^+}$ has scalar action on $\Fil_2$, in which case the fiber is the $\Pone$ of possible $\Fil_1$'s, or there is a unique line on which $G_{F_v^+}$ acts through $\chi_{1,x}$, and the fiber is a point.  In particular either way this fiber is connected.
 
 By Lemma \ref{lem: criterion for H2 vanishing} \eqref{item: pure crystalline
   vanishes} each point of this fiber is contained in a unique irreducible
 component of $\cG_v$, and as the fiber is connected, the entire fiber is
 contained in this component.  It follows that the image of this component in
 $\Spec R_v^{\triangle}$ is  the unique irreducible component containing $\rho$.

Now we prove that all such $\rho$ lie on the same irreducible component.
Consider the closed subscheme $\cG_v^{\mathrm{flat}}\subseteq\cG_v$ whose points
for any $\cO$-algebra $A$ are pairs $(\Fil_\bullet,\rho)$ where
$\rho\otimes\varepsilon$ is finite flat and $\Fil_\bullet$ is a filtration with
$\Fil_2$ unramified.  We claim that the formal completions of
$\cG_v^{\mathrm{flat}}$ at $k'$-points for $k'/k$ finite are formally smooth.
By a standard tangent-obstruction calculation as in Lemma \ref{lem: criterion
  for H2 vanishing} this amounts to the vanishing of
$H^2_{\mathrm{flat}}(\Fil^0\ad^0\rhobar)$ (\emph{cf.}\ the proof of \cite[Prop.\
2.4.4]{MR2551765}).%

We now consider $R_v^{\triangle,\mathrm{flat}}$, the scheme-theoretic image of $\cG_v^{\mathrm{flat}}$ in $R_v^\triangle$.  It is irreducible by the same argument as above, as the fiber over $\rhobar$ is either a point or $\Pone$.  As every pure weight 2 crystalline point lies on this irreducible locus, they all lie on the same unique irreducible component of $R_v^\triangle$.

Finally the dimension can be computed at any $p$-distinguished point using Lemma \ref{lem: criterion for H2 vanishing}
\end{proof}

We finally prove a lemma which will help with ``Ihara avoidance''.  Let $Q\subset R_v^\triangle$ be a minimal prime, corresponding to an irreducible component of $\Spec R_v^\triangle$.

\begin{lem}\label{lem: for Ihara}
Suppose that $\Spec R_v^{\triangle}/Q \to\Spec \Lambda_{\GSp_4,v}$ is surjective.  Then there exists a minimal prime of $R_v^{\triangle}/(p)$ which contains $Q$ and no other minimal prime of $R_v^\triangle$.  Moreover $R_v^\triangle/Q$ has relative dimension 16 over $\cO$.
\end{lem}
\begin{proof}
By the hypothesis we can take an $\FF_q((t))$ valued point $x$ of $\Spec R_v^{\triangle}/Q$ so that no ratio of the characters $\chi_{1,x},\chi_{2,x},\chi_{2,x}^{-1}\varepsilon^{-1},\chi_{1,x}^{-1}\varepsilon^{-1}$ is $1$ or $\varepsilon$ (even on inertia).  By the same argument as in Lemma \ref{lem: criterion for H2 vanishing} \eqref{item:easy vanishes} and \eqref{item: def ring regular} we have that the local ring $R_{v,x}^{\triangle}$ is regular.  We now take any irreducible component of $\Spec R_v^\triangle/(p)$ containing $x$.
\end{proof}

We expect that the hypothesis in Lemma \ref{lem: for Ihara} is always satisfied.  Rather than attempt a direct local proof of this fact we will check it by global means in the application in the next section.

\subsection{An ordinary modularity lifting theorem for $\GSp_4$, $p>2$}\label{sec: Gsp4  modularity lifting}

We explain how to prove a modularity lifting theorem for a $p$-adic
Hida family of Hilbert--Siegel modular forms over a totally real
field~$F^+$ in which the odd prime~$p$ splits completely; this is a
slight adaptation of the arguments of our earlier paper
\cite{BCGP} and their improvements by
Whitmore~\cite{Whitmore}. Indeed,
under a residually $p$-distinguished hypothesis, our theorem is a very
special case of the theorems proved in those papers. (The entire
difficulty in~ \cite{BCGP} was about proving modularity lifting
theorems in the case $l_0>0$, but the $l_0=0$ case that we needed here
is completely routine.) It would of course be more natural not to
include the assumption that~$p$ splits completely in~$F^{+}$, but as
we do not know a reference for the relevant Hida families beyond this
case, %
we leave such
results for a future paper. %
We let $F^+$ be a totally real field in which $p>2$ splits completely and let $\pi$ be an ordinary cuspidal automorphic
representation for~$\GSp_4/F^{+}$ of central character $|\cdot|^2$
and weight~$((k_v, l_v;2))_{v\mid \infty}$. Using our fixed isomorphism
$\imath:\C\cong\Qpbar$ we identify the places~$v\mid
\infty$ and $v\mid p$ without further comment. %
Recall from Theorem~\ref{thm:irregular-ordinary-pi-Galois-rep}
that there is a Galois representation $\rho_{\pi,p}:G_{F^+}\to\GSp_4(\Qpbar)$
associated to~$\pi$.  We let $\rhobar=\rhobar_{\pi,p}$.  We fix $R$, a finite set of finite, prime to $p$ places of $F^+$ containing all the prime to $p$ places where $\pi$ is ramified.

We make the following assumptions:

\begin{hypothesis}\leavevmode
  \label{hypothesis: Q conditions we need for GSp4 modularity lifting}
     \begin{enumerate}
        \item $\rhobar$ is $\GSp_4$-reasonable, in the sense
          of~\cite[Defn.\
          3.19]{Whitmore}. In
          particular, $\rhobar$ is absolutely irreducible.
         \item $\rhobar$ is tidy, in the sense of~\cite[Defn.\
           7.5.11]{BCGP}. 
           \item\label{item:ordinary-of-weight-2-hypothesis} For each $v|p$,
$\rhobar|_{G_{F^+_v}}$ is ordinary of weight 2, with a fixed $p$-stabilization
$(\alphabar_v,\betabar_{v})$ which is compatible with a fixed choice of ordinary
$p$-stabilization of~$\pi_v$. %
           \item\label{item: less ludicrous lazy assumption} For each $v\in
             S_p$, the representation $\rho_{\pi,p}|_{G_{F_v^+}}$ lies on a
             unique irreducible component of $R_v^\triangle$
             (where~$R_v^\triangle$ is defined via the $p$-stabilization
             $(\alphabar_v,\betabar_v)$, and $\rho_{\pi,p}|_{G_{F_v^+}}$ is a
             point of~$R_v^{\triangle}$ via the chosen $p$-stabilization on~$\pi_v$).
\item\label{item: for Ihara} For each place $v\in R$ we have:
\begin{itemize}
\item  $\rhobar|_{G_{F_v^+}}$ is trivial.
\item  $q_v\equiv 1\pmod p$, and if $p=3$, then $q_v\equiv 1\pmod 9$.
\item $\pi_v^{\Iw(v)}\not=0$.
\end{itemize}
     \end{enumerate}
   \end{hypothesis}
Note that we do not assume that $\pi$ is ordinary of weight 2 (it will be in the main
application but we also allow $\pi$ to have regular weight in order to prove
Lemma~\ref{higherserre}). %

\begin{rem}\label{rem: remark clarifying hypotheses}
The main result of this section will be a ``minimal at $p$'' modularity lifting theorem for ordinary $p$-adic modular forms.  In particular hypothesis \eqref{item: for Ihara} will be used for Taylor's Ihara avoidance argument, in order to have no minimality hypotheses away from $p$.  In the application, $F^+/\Q$ will be a solvable extension chosen to ensure this, and the main modularity lifting theorem for $\Q$ will be deduced using base change.

We also note that the theorem is only ``minimal at $p$'' due to our failure to completely analyze the deformation rings $R_v^\triangle$ in the previous section.  However we emphasize that when $\rhobar|_{G_{F_v^+}}$ is residually $p$-distinguished then $\Spec R_v^\triangle[1/p]$ is irreducible and hypothesis \eqref{item: less ludicrous lazy assumption} is automatic.
\end{rem}

   \begin{rem}
     \label{rem: remark on Whitmore conditions}The most important applications of the
     result of this section are in the case~$p=3$ and $\rhobar(G_{F^+})\subseteq\GSp_4(\F_3)$. In this case the
     condition that $\rhobar$ is $\GSp_4$-reasonable can be made completely
     explicit, see~\S\ref{sec:subgroupsofGSp}.
   \end{rem}

By the assumption that~$\rhobar(G_{F^+})$ is tidy, we can choose an unramified
place~$w_0$ of~$F^+$ of residue characteristic greater than~$5$, with
$q_{w_0}\not\equiv 1\pmod{p}$, and such that no two eigenvalues
of~$\rhobar(\Frob_{w_0})$ have ratio~$q_{w_0}$. We set~$S=S_p\cup
R\cup\{w_0\}$. %

By the above hypotheses, we are in the situation of~\cite[Hyp.\
7.8.1]{BCGP}, except that we have not assumed that we are
residually $p$-distinguished. In the notation of~\cite[\S 7.8]{BCGP}
and~\cite[\S 7]{Whitmore}, we
take~$I=\emptyset$, so that by definition the ring~$\Lambda_I$ is
equal to
$\Lambda_{\GSp_4,F^+}:=\cotimes_{v|p}\Lambda_{\GSp_4,v}$. The only change that we make to
the setup of~\cite{BCGP,
  Whitmore} is that for $v|p$ we use the deformation problem~$\cD_v^{\triangle}$, taking
$\chibar_1,\chibar_2$ to be unramified with
$\chibar_1(\Frob_v)=\alphabar_p$, $\chibar_2(\Frob_v)=\betabar_p$. Note that in
the residually $p$-distinguished case that $\alphabar_p\ne\betabar_p$, this
agrees with the deformation problem denoted $\cD_v^{B,\overline{\alpha}_p}$ in
\cite{BCGP,Whitmore}. %

We can then carry out all of the constructions made in \cite[\S
7.8]{BCGP} and~\cite[\S
7]{Whitmore}, which for the most part
make no use of the hypothesis that $\alphabar_p\ne\betabar_p$: the proof of
\cite[Thm.\ 7.8]{Whitmore} generalizing \cite[Thm.\ 7.9.4]{BCGP} only uses the
values of $\cD_v^{\triangle}$ on $\cO_{E'}$ for $E'/E$ a finite extension
(recall that in the non residually $p$-distinguished case we don't necessarily
understand the values of $\cD_v^{\triangle}$ on general complete Noetherian
local rings due to its definition as a scheme-theoretic image). We now recall
the key points, allowing ourselves to simplify the notation slightly in
comparison to that of~\cite{BCGP}, by dropping the symbols ``$I$'' and
``$\alphabeta$'' appearing there. %

In particular, we have the global deformation problem %
\[ \cS_{1} = (\rhobar, S, \{\Lambda_{\GSp_4,v}\}_{v \in S_p}\cup\{ \cO
  \}_{v \in S\setminus S_p},\{\cD_v^{\triangle}
  \}_{v \in S_p} \cup \{ \cD_v^1 \}_{v \in R} \cup \{
  \cD_{w_0}^\square \}),\] where ~$\cD_v^1$ is defined in~\cite[\S
7.4.5]{BCGP} (it corresponds to unipotently ramified
liftings). %
There is a
surjection of $\Lambda_{\GSp_4,F^+}$-algebras
$R_{\cS_{1}}\to\T_{\cS_{1}}$, where~$\T_{\cS_{1}}$ is the Hida Hecke
algebra considered %
in~\cite[\S
7.9]{BCGP}: it acts faithfully on a finite free
~$\Lambda_{\GSp_4,F^+}$-module $M^{1}%
$, which is
obtained from the ordinary part of the coherent $\HH^0$ of
Hilbert--Siegel Shimura varieties.

By Hypothesis~\ref{hypothesis: Q conditions we need for GSp4 modularity
    lifting}~\eqref{item: less ludicrous lazy assumption}, for each $v|p$ the
  representation $\rho_{\pi,p}|_{G_{F^+_v}}$ lies on a unique irreducible component of $\Spec R_v^\triangle$, which we denote~$\Spec R_{v}^{\triangle,\pi}$. We let $\cD_v^{\triangle,\pi}$ be the
  deformation problem determined by this irreducible component, and
  write 
\[ \cS_{1,\pi} = (\rhobar, S, \{\Lambda_{\GSp_4,v}\}_{v \in S_p}\cup\{ \cO
  \}_{v \in S\setminus S_p}, \{\cD_v^{\triangle,\pi}
  \}_{v \in S_p} \cup \{ \cD_v^1 \}_{v \in R} \cup \{
  \cD_{w_0}^\square \}).\] We let $\T_{\cS_{1,\pi}}$ denote 
$R_{\cS_{1,\pi}}\otimes_{R_{\cS_{1}}}\T_{\cS_{1}}$.  These should be thought of as ``$p$-minimal'' deformation rings and Hecke algebras, see also Remark \ref{rem: remark clarifying hypotheses}.

   \begin{thm}
     \label{thm: GSp4 p bigger than 2 automorphy lifting}Assume that
     we are in the above situation, so that in particular Hypothesis~\ref{hypothesis: Q conditions we need for GSp4 modularity
    lifting} holds. Then~$R_{\cS_{1,\pi}}$ is a finite $\Lambda_{\GSp_4,F^+}$-algebra, and the
  morphism $R_{\cS_{1,\pi}}\to\T_{\cS_{1,\pi}}$ has nilpotent kernel, i.e.\ $(R_{\cS_{1,\pi}})^{\red}\isoto\T_{\cS_{1,\pi}}$.
\end{thm}

\begin{proof}
We first verify that for $v\mid p$, the irreducible components $R_v^{\triangle,\pi}$ satisfy the hypothesis of Lemma \ref{lem: for Ihara}.  For this consider a minimal prime $Q_\pi\subset \T_{\cS_1}$ contained in the prime ideal corresponding to $\pi$ and the chosen $p$-stabilizations, and consider the composition
$$R_v^\triangle\to R_{\cS_1}\to\T_{\cS_1}\to\T_{\cS_1}/Q_\pi.$$
By Hypothesis \ref{hypothesis: Q conditions we need for GSp4 modularity lifting} \eqref{item: less ludicrous lazy assumption}, the composite must factor through the component $R_v^{\triangle,\pi}$.  We now claim that the composite
$$\Spec \T_{\cS_1}/Q_\pi\to \Spec R_v^{\triangle,\pi}\to\Spec \Lambda_{\GSp_4,v}$$
is surjective, and hence the second map is surjective, which is what we are trying to prove.  To see this note that the composite is also
$$\Spec \T_{\cS_1}/Q_\pi\to\Spec\Lambda_{\GSp_4,F^+}\to\Spec\Lambda_{\GSp_4,v}.$$
Here the second map is clearly surjective, while the first map is surjective because $\Spec \T_{\cS_1}/Q_\pi$ is finite and torsion free as a $\Lambda_{\GSp_4,F^+}$-module, since $\T_{\cS_1}$ acts faithfully on a finite free $\Lambda_{\GSp_4,F^+}$-module.

 Now we proceed with the proof of the theorem.  It is enough to prove that the
  support of $M^{1}$ in $\Spec R_{\cS_1}$
  contains ~$\Spec R_{\cS_{1,\pi}}$. As in~\cite[\S
  7.3]{Whitmore}, we have a power
  series ring~ $S_\infty$ over~$\Lambda_{\GSp_4,F^+}$, and we write~$\mathfrak{a}_\infty$ for the augmentation ideal
  $\ker(S_\infty\to \Lambda_{\GSp_4,F^+})$.  We also have a power
  series ring $R^{1}_\infty$ over
  $R_{\cS_1,S}^{\loc}$. Set
  \[R_{\infty,\pi}:=R^{1}_\infty\otimes_{R_{\cS_1,S}^{\loc}}R_{\cS_{1,\pi},S}^{\loc}.\]

     The patching argument
     of~\cite[Prop.\ 7.11]{Whitmore}
     provides us in particular with:
     \begin{itemize}
     \item $\Lambda_{\GSp_4,F^+}$-algebra morphisms $S_\infty \to
       R^{1}_\infty\to  R_{\cS_1}$;
     \item a $R^{1}_\infty$ module $M^{1}_\infty$ which is free as an $S_\infty$ module (and hence has depth as an $R^{1}_\infty$ equal to the dimension of $S_\infty$).
     \item an isomorphism
       $M^{1}_\infty/\mathfrak{a}_\infty\cong
       M^{1}$;
     \item a commutative diagram of $S_\infty$-algebras  \[\xymatrix{R^{1}_\infty\ar[d]\ar[r]&\End_{S_\infty}(M^{1}_\infty)\ar[d]^{-\otimes_{S_\infty}\Lambda_{\GSp_4,F^+}}\\
          R_{\cS_1}\ar[r]&\End_{\Lambda_{\GSp_4,F^+}}(M^{1})}\]
    \end{itemize}

It thus suffices to prove that the support of $M^{1}_\infty$ contains every irreducible component of $R_{\infty,\pi}$.  Exactly as in
    the proof of Theorem~\ref{thm: U(n) ordinary automorphy lifting including Ihara
       avoidance}, we know that the support of ~$M^{1}_\infty/\varpi$
    contains $\Spec R_{\infty,\pi}/\varpi$. (This is Taylor's ``Ihara
     avoidance'' argument, using the data
     $M^{\chi}_\infty$ etc. from \cite[Prop.\
     7.11]{Whitmore} which we have
     not recalled here.)  Every irreducible component of the support of
     $M_\infty$ has dimension equal to that of $S_\infty$, and as every
     irreducible component of $\Spec R_{\infty,\pi}$ has dimension equal to that
     of~$S_\infty$, $M_\infty$ will be supported on an irreducible component of $\Spec R_{\infty,\pi}$ as soon as it is supported on some point which is only contained in that irreducible component.
     
     It thus suffices to show that for each irreducible component of $\Spec
     R_{\infty,\pi}$ there is a point of $\Spec R_{\infty,\pi}/p$ contained in
     it and in no other component.  Since $R_{\infty,\pi}$ is a power series
     ring over a completed tensor product of local deformation rings, it suffices to check the
     same property for each factor, i.e.\ to check that the same holds for $R_{v}^{\triangle,\pi}$ (for $v\in
S_p$), $\Spec R_v^1$ (for $v\in R$), and $\Spec R_{w_0}^\square$. The
first of these follows from Lemma~\ref{lem: for Ihara} (noting that we have verified the hypothesis above), the second from \cite[Prop.\
7.4.7]{BCGP}, and the last from our choice of~$w_0$, which
guarantees that~$R_{w_0}^\square$ is formally smooth over~$\cO$.
\end{proof}

\subsection{Subgroups of~$\GSp_4(\F_3)$} \label{sec:subgroupsofGSp}
In this section, we shall identify the precise subgroups of~$\GSp_4(\F_3)$
we are allowing for our modularity lifting theorems.
We first identify the regular semi-simple elements in~$\GSp_4(\F_3)$.
We have:
\begin{lemma}
There are three conjugacy classes of elements~$g \in \GSp_4(\F_3)$
such that~$\nu(g) = 1$ and~$g$ is regular semi-simple, namely:
\begin{enumerate}
\item The unique conjugacy class of elements of order~$5$,
\item The unique conjugacy class of elements of order~$10$,
\item The unique conjugacy class of elements of order~$8$ which lie in~$\Sp_4(\F_3)$.
\end{enumerate}
There are five conjugacy classes of elements~$g \in \GSp_4(\F_3)$
such that~$\nu(g) = -1$ and~$g$ is regular semi-simple, namely:
\begin{enumerate}
\item Both conjugacy classes of elements of order~$20$,
\item Three of the five conjugacy classes of elements of order~$8$,
namely those whose images in~$\PSp_4(\F_3)$ lie in the conjugacy
classes~$4D$ or~$8A$ but not~$4C$ in the notation
of Lemma~\ref{atlas}.
\end{enumerate}
\end{lemma}

We now turn to reasonableness~\cite[Defn.\
          3.19]{Whitmore}. Although this definition does not \emph{a priori}
          depend only on the image of~$\rhobar$, it shall turn out that that under our
          running assumptions this will be true.
          
\begin{lemma} \label{imageonly}
Suppose that~$A/\Q$ is an abelian surface, and that~$A$ has good
reduction at some~$p > 2$. Then 
the image of~$\rhobar_{A,p} |_{G_{\Q(\zeta_p)}}$ coincides with the image
of~$\rhobar_{A,p} |_{G_{\Q(\zeta_{p^n})}}$ for all~$n \ge 1$.
\end{lemma}

\begin{proof}  Let~$K/\Q_p$ be the fixed field of the kernel of~$\rhobar_{A,p}|_{G_{\Qp}}$, which
certainly contains~$\Q_p(\zeta_p)$. Since~$\Gal(\Q_p(\zeta_{p^n})/\Q_p(\zeta_p))$ is cyclic,
the lemma holds unless there exists an inclusion~$\Q_p(\zeta_{p^2}) \subseteq K$. Assume such
an inclusion exists. It follows that
 the root discriminant~$\delta_K$ is divisible by the root discriminant of~$\Q(\zeta_{p^2})$,
which is~$p^{(2p^2 - 3p)/p(p-1)} = p^{(2p-3)/(p-1)}$. The assumption that~$A$ has good reduction implies that~$A[p]/\Z_p$
is a finite flat group scheme, which by~\cite[2.1 Thm.\ 1]{MR807070} implies that the root discriminant 
of~$K$ satisfies
$$v_p(\delta_K) < 1 + \frac{1}{p - 1}.$$
Since~$p \ge 3$, this contradicts our lower bound:
\[v_p(\delta_K) \ge v_p(\delta_{\Q(\zeta_{p^2})})
= 1 + \frac{1}{p-1} + \frac{p-3}{p-1}. \qedhere \]
\end{proof}
We deduce:
\begin{lemma} \label{listofsubgroups} Let~$A/\Q$ be an abelian surface
  with a prime to~$3$ polarization and good reduction at~$3$, and let
$$\rhobar = \rhobar_{A,3}: G_{\Q} \rightarrow \GSp_4(\F_3)$$
denote the corresponding mod~$3$ representation. Then the following
hypotheses:
\begin{enumerate}
\item \label{reasonablefortable} $\rhobar$ is~$\GSp_4$-reasonable in the sense of~\cite[Defn.\
          3.19]{Whitmore},
          \item \label{tidyfortable} $\rhobar$ is tidy in the sense of~\cite[Defn.\
           7.5.11]{BCGP},
 \item \label{regfortable} $\rhobar(G_{\Q(\zeta_3)})$ contains a regular
 semi-simple element,
   \item \label{regular-semi-simple-not-splitting} $\rhobar(G_{\Q})\setminus \rhobar(G_{\Q(\zeta_3)})$ contains a regular
   semi-simple element,
 \end{enumerate}
 are satisfied precisely if~$\Gamma' = \rhobar(G_{\Q})$ 
 in Table~\ref{tabletwo} 
 has a tick, where otherwise
 the cross 
 indicates the corresponding obstruction to condition~\eqref{reasonablefortable},
 \eqref{tidyfortable}, \eqref{regfortable}, or~\eqref{regular-semi-simple-not-splitting}. 
 In particular, these
 conditions are all satisfied if $\rhobar(G_{\Q})=\GSp_4(\F_3)$.
 \end{lemma}

\numtable 
\begin{center}
\begin{tabular}{*3l|*4c|c}
\toprule
\multicolumn{3}{c|}{$\Gamma' \subset \GSp_4(\F_3)$ and $\Gamma = \Gamma' \cap \Sp_4(\F_3)$} & \multicolumn{4}{c|}{Conditions} &  \\
\midrule
LMFDB label & \multicolumn{2}{l|}{small group labels for $\Gamma'$, $\Gamma$} & \eqref{reasonablefortable} & \eqref{tidyfortable} & \eqref{regfortable} & \eqref{regular-semi-simple-not-splitting} &  \\  
\midrule
\texttt{3.1620.1} & $ \texttt{<64,258>}$ & $ \texttt{<32,44>}$ &     &    &   & & \TTTT  \\ 
\texttt{3.1620.2} & $ \texttt{<64,258>}$ &$ \texttt{<32,50>}$&     &    &  \FFFFF &   &  \\ 
\texttt{3.1620.5} & $ \texttt{<64,152>}$ & $ \texttt{<32,44>}$ &     &    &   & & \TTTT  \\ 
\texttt{3.1620.10} & $ \texttt{<64,152>}$ &$ \texttt{<32,8>}$ &     &    &   & & \TTTT  \\ 
\texttt{3.1296.1} & $ \texttt{<80,29>}$ & $\texttt{<40,3>}$ &     &    &   & & \TTTT  \\ 
\texttt{3.810.1} & $ \texttt{<128,137>}$& $ \texttt{<64,137>}$ &     &  \FFFFF  &  &   &  \\ 
\texttt{3.810.2} & $ \texttt{<128,2023>}$ & $ \texttt{<64,137>}$ &     &    &    & & \TTTT  \\ 
\texttt{3.810.5} & $ \texttt{<128,137>}$ & $ \texttt{<64,37>}$ &     &  \FFFFF  &  & &    \\ 
\texttt{3.810.6} & $ \texttt{<128,142>}$ &$ \texttt{<64,37>}$ &     &    &   & & \TTTT  \\ 
\texttt{3.540.1} & $ \texttt{<192,1485>}$ & $ \texttt{<96,190>}$ &     &    &   & & \TTTT   \\ 
\texttt{3.540.2} & $ \texttt{<192,1483>}$ & $ \texttt{<96,191>}$ &  \FFFFF  &    &  & &   \\ 
\texttt{3.540.3} & $ \texttt{<192,1485>}$&$ \texttt{<96,202>}$ &     &    &  \FFFFF  & &  \\ 
\texttt{3.540.5} & $ \texttt{<192,1018>}$&$ \texttt{<96,202>}$&     &    &  \FFFFF  & &   \\  
\texttt{3.540.7} & $ \texttt{<192,965>}$ & $ \texttt{<96,191>}$ &   \FFFFF  &    &  &  &  \\ 
\texttt{3.405.1} & $ \texttt{<256,6671>}$  & $ \texttt{<128,937>}$  &    &    &  & & \TTTT   \\ 
\texttt{3.270.1} & $ \texttt{<384,18045>}$ & $ \texttt{<192,989>}$ &    &    &  & & \TTTT   \\ 
\texttt{3.216.1} & $ \texttt{<480,948>}$ & $ \texttt{<240,90>}$  &    &    &   & & \TTTT  \\ 
\texttt{3.216.2} & $ \texttt{<480,947>}$ &$ \texttt{<240,89>}$ &  \FFFFF  &    &  & &    \\  
\texttt{3.162.1} & $ \texttt{<640,21454>}$  & $ \texttt{<320,1581>}$ &    &  \FFFFF  &  & &    \\ 
\texttt{3.135.1} & $\texttt{<768,1086054>}$ & $ \texttt{<384,618>}$ &    &    &   & & \TTTT  \\ 
\texttt{3.135.2} & $\texttt{<768,1086054>}$ & $ \texttt{<384,18130>}$&    &    &   & & \TTTT  \\ 
\texttt{3.45.1} & $ 2304$ & $1152$ &    &    &   & & \TTTT  \\ 
\texttt{3.36.1} & $ 2880$ & $1440$ &  \FFFFF  &    &  &    &  \\ 
\texttt{3.27.1} & $ 3840$ & $1920$ &    &    &   & & \TTTT  \\ 
\texttt{1.1.1} & $ 103680$ & $51840$ &    &    &   & & \TTTT   \\ 
\bottomrule
\end{tabular}
\caption{Conjugacy classes of subgroups~$\Gamma' \subset \GSp_4(\F_3)$
with~$\nu(\Gamma') \ne 1$ and~$\Gamma = \Gamma' \cap \Sp_4(\F_3)$ 
absolutely irreducible. LMFDB labels determine the conjugacy class
of~$\Gamma'$, the small group labels~\cite{MR1826989}
determine $\Gamma$, $\Gamma'$ up to abstract isomorphism.}  
\label{tabletwo}
\end{center}
\end{table}

\begin{proof}Note that if $\Gamma' = \rhobar(G_{\Q})$,
  and~$\Gamma=\rhobar(G_{\Q(\zeta_3)})$, then
  $\Gamma=\Gamma'\cap\GSp_4(\F_3)$. Furthermore,  %
reasonableness (which a priori depends on the image of~$\rhobar|_{G(\Q(\zeta_{3^n}))}$
for all~$n$) only depends on the image of~$\rhobar|_{G(\Q(\zeta_{3}))}$ by
Lemma~\ref{imageonly}.
As noted in~\cite[\S4.3]{Whitmore}, the spanning condition of reasonableness is satisfied
for all of these subgroups. 
We have  listed the abstract isomorphism types of~$\Gamma'$ and~$\Gamma$
according to the small groups database~\cite{MR1826989} when they are of small order.
The LMFDB subgroup labels~\cite{LMFDB}
(which for proper subgroups are of the form~$\texttt{3.i.n}$ where~$\texttt{i} = [\GSp_4(\F_3):\Gamma']$) determine~$\Gamma'$
up to conjugacy in~$\GSp_4(\F_3)$.
The groups of larger order are described more explicitly in~\cite[Lemma~7.5.21]{BCGP}.
The abstract isomorphism type  of~$\Gamma'$ is already enough to determine~$\Gamma'$
up to conjugation in~$\GSp_4(\F_3)$ (under the assumption that~$\Gamma = \Gamma' \cap \Sp_4(\F_3)$
acts absolutely irreducibly) except for two pairs with~$|\Gamma'| = 64$, one pair with~$|\Gamma'| = 128$, and the pair of groups
with~$|\Gamma'| = 768$, and in all such cases  they can be distinguished by the abstract isomorphism type
of~$\Gamma = \Gamma' \cap \Sp_4(\F_3)$.
\end{proof}

\begin{rem}[$\Sp_4(\F_3)$ is not~$\GL_4$-adequate]  \label{rem:curious} 
There is a natural inclusion
\numequation
\rhobar: \Sp_4(\F_3) \hookrightarrow \GL_4(\F_3).
\end{equation}

It turns out that the image~$G$ of~$\rhobar$ is \emph{not} adequate
(in the sense of ~\cite[Defn.\ 2.20]{MR3598803}), which is the reason
why, when~$p=3$, we need to use~$\GSp_4$ modularity lifting theorems
rather than~$U(4)$ automorphy lifting theorems (in contrast to our treatment of
the case~$p=2$ in Section~\ref{sec:
  R=T unitary}). The failure of adequacy can be seen directly as
follows. The group~$G \simeq \Sp_4(\F_3)$ has exactly two irreducible
representations~$V$, $V^{\sigma}$ of dimension~$4$ over~$\C$. The
representations are defined over the ring~$\Z[\zeta_3]$ and are
conjugate under the action of~$\Gal(\Q(\zeta_3)/\Q)$~\cite{Atlas}.
The mod~$\pi = (1 - \zeta_3)$ reduction of this representation
is~$\rhobar$. The corresponding mod~$\pi^2 = (3)$ reduction:
$$\rho: G \rightarrow \SL_4(\Z[\zeta]/3) = \SL_4(\F_3[\epsilon]/\epsilon^2)$$
gives a non-trivial deformation of~$\rhobar$ and a non-zero class
 in~$H^1(G,\sl_{4}) \subset H^1(G,\gl_{4})$. The deformation~$\rho$ is not, however,
valued in~$\Sp_4(\F_3[\epsilon]/\epsilon^2)$; this reflects the fact that~$V$ is not self-dual in characteristic zero;
we have~$V^{\vee} \simeq V^{\sigma}$. 
In particular, identifying~$\gl_4$ with~$\Hom(\rhobar,\rhobar) \simeq \rhobar
\otimes \rhobar  \simeq \Sym^2(\rhobar) \oplus \wedge^2(\rhobar)$, 
 this cohomology class lives in~$H^1(G,\wedge^2 \rhobar)$. %
In contrast,  for~$\rhobar$ to be adequate as a symplectic representation,  it suffices that
$H^1(G,\franksp_4)=0$ where (in this example)~$\franksp_4$ may be identified with~$\Sym^2(\rhobar)$.

 This example is similar to the failure of the image of the map
\numequation \label{noproblemforandrew}
\rhobar: \SL_2(\F_5) \hookrightarrow \GL_2(\F_5)
\end{equation}
to be adequate. This failure of adequacy for~(\ref{noproblemforandrew}) 
does not cause an issue  in~\cite{MR1333035};  one  exploits the fact that the fixed field of
the kernel of the adjoint representation of~$\rhobar_{E,5}$ for an elliptic curve~$E/\Q$
 does not contain~$\zeta_5$ (see  the proof 
of~\cite[Prop~1.11]{MR1333035}). On the other hand, for an abelian surface~$A$, the fixed of the kernel of
the adjoint representation of~$\rhobar_{A,3}$ always  contains~$\zeta_3$, so there is no
way to avoid this cohomological obstruction. Hence this situation is more analogous to
the problem of proving 
modularity lifting for elliptic curves over~$\Q(\sqrt{5})$
using~$5$-adic modularity lifting theorems; see the introduction
to~\cite{MR3648503} for an exposition of this case, and an explanation
of why there are classes in the dual Selmer group (so-called ``Lie
classes'') that cannot be killed by  Taylor--Wiles primes. 
\end{rem}

\section{Multiplicity one
  theorems}\label{sec: SW type argument}%
Our classicality theorems in low weight (in particular Theorem~\ref{thm:multiplicity-one-implies-classical})
require as input a multiplicity one theorem in characteristic zero.
The main goal of this section is to prove
 such a theorem. Note that multiplicity one really consists
of two separate statements --- firstly that the multiplicity is at least one
(a $p$-adic modularity statement), and secondly that the multiplicity
is at most one.

One approach to proving multiplicity one (following
Diamond~\cite{MR1440309}) would be to prove an~$R=\T$ theorem
for the corresponding ordinary Hida family and then, assuming
the local deformation ring at~$p$ is formally smooth, deduce
that the corresponding module~$M$ of modular forms is free,
and moreover free of rank one by specialization at classical points.
Such an argument would work 
if we made the additional
hypotheses that~ $p>2$ and that $\rhobar$ was $p$-distinguished (as is
done in~\cite{BCGP}). %
Since we are not making such assumptions,
a further argument is required.
The first point to note is that we are working in characteristic
zero and hence we only need  prove that  $M$ is free (and non-zero) after localizing
at a height one prime $\q$ corresponding to our characteristic zero
representation. To show that~$M$ is non-zero when~$p=2$ we are able to
appeal to the $R=\T$ theorem for unitary groups that we proved in Section~\ref{sec:
  R=T unitary}, while for~$p>2$ we use the $R=\T$ theorems
for~$\GSp_{4}$ proved in Section~\ref{sec: Gsp4
  modularity lifting}. In either case,  Diamond's argument applies (at
least in principle) providing
that the formal completion of the local deformation ring at $\q$
is regular, something that is ultimately true under our hypotheses. 

More precisely, what is ultimately required for our arguments is the following. First, we need 
an  $R^{\red}[1/p]= \T[1/p]$ theorem in our higher Hida theory (not yet classical) situation.
 (In truth, when~$p>2$,
we get away with a weaker version of
such a theorem in a neighbourhood of the prime $\q$, at the cost of some further local
complications already considered in~\S\ref{sec: R=T GSp4 preliminaries all p}.) This proves
that $M_{\q}$ is non-zero. Second, we want to control the relative
tangent space of $R$ at the prime~$\q$. This is closely related
to establishing the vanishing of the adjoint Bloch-Selmer group (in
characteristic zero) of our characteristic zero representation.
Theorems of this kind were proved by Newton 
and Thorne~\cite{MR4592862} in some generality for Galois
representations associated to  automorphic
representations of~$\mathrm{GL}_n$ of unitary type,
and we follow their arguments closely. 
 In fact our situation
 is for the most part simpler than theirs, since as we are
assuming that $\rhobar$ is absolutely irreducible we
do not have to use pseudorepresentations.
Finally, we need enough local properties of the local deformation
ring at~$p$ in characteristic zero at~$\q$, and this is what ultimately
requires the $p$-distinguished hypothesis in characteristic zero.

A summary of this section is as follows.
 In~\S\ref{subsec: TW primes for mult one}, we adapt the Galois-theoretic arguments of Newton--Thorne~\cite{MR4592862} to
our setting. 
In~\S\ref{subsec: abstract TW patching for mult one} we set up the basic patching formalism required
for our argument, 
and in~\S\ref{subsec: defn of GSp4 TW
  system for mult one}, we show how this can be applied
in the setting of ~$\GSp_4$ over~$\Q$, by patching modules coming from higher
Hida theory, as recalled in~\S\ref{higherhidaRT}. Note that there is quite a lot
of overlap between the arguments of~\S\ref{subsec: defn of GSp4 TW
  system for mult one}
and of similar ones in~\cite[\S 7.8, 7.9]{BCGP}
 --- the difference being that the
latter worked under a more restrictive hypothesis on $\rhobar$
but also proved strong integral statements. Finally, 
in~\S\ref{subsec: deducing R=T GSp4
     from base change}, we  prove
the desired multiplicity one theorem.
Note that our arguments certainly require understanding
the multiplicities of certain automorphic representations
in cohomology, which ultimately uses
Arthur's
   classification of discrete automorphic representations of~$\GSp_4$.

\subsection{Taylor--Wiles primes}\label{subsec: TW primes for mult one}
Let~$H$ be a  compact subgroup of~ $\GSp_4(\cO)$. After replacing~$E$ by a finite extension, we can assume that
for each element~$h\in H$, the characteristic polynomial of~$h$ is
already split over~$E$. We will assume this without comment from now
on. We write~$\Hbar$ for the image of~$H$ in~$\GSp_4(k)$. Throughout
this section, we write~$\ad$, $\ad^0$ for the Lie algebras associated
to~$\GSp_4$ and~$\Sp_4$ over~$\cO$.

\begin{defn}
  \label{defn: integrally enormous}A  compact subgroup~$H$ of
  $\GSp_4(\cO)$ is \emph{integrally enormous} if it acts absolutely
  irreducibly on~$E^4$, and if for all
  simple $E[H]$-submodules
  $W \subset E \otimes \ad^0$, there exists an  
  element
  $h \in H$ such that
  \begin{itemize}
 \item $1$ is an eigenvalue for the action of~$h$ on~$W$, and   \item
   the image~$\overline{h}$ of $h$ in~$\GSp_4(k)$ is regular semi-simple (i.e.\ has $4$ distinct
   eigenvalues).\end{itemize}
\end{defn}%

\begin{lem}[Examples of integrally enormous representations, I]
  \label{lem: sp4 is integrally enormous} \leavevmode
  \begin{enumerate}
  \item \label{sp4case}
      If~$\Hbar$ contains a regular semi-simple element, and the
    Zariski closure of~$H$ contains~$\Sp_4$, then~$H$ is integrally
    enormous.
    \item \label{inducedcase} %
    Suppose that the action of~$H \subset \Sp_4(\OL)$ is absolutely irreducible but becomes
    reducible after restriction to an index~$2$ subgroup~$G$. 
    Suppose that~$(\Hbar \setminus \Gbar)$ contains 
    a regular semi-simple
    element.
        Suppose
    that the Zariski closure of~$G$ is~$\SL_2 \times \SL_2$. Then~$H$
    is integrally enormous.
     \end{enumerate}
\end{lem}
\begin{proof} We first consider case~(\ref{sp4case}).
Since ~$\Sp_4$ acts irreducibly on both its standard
  representation and on~ $\ad^0$, $H$ acts
  irreducibly on both~$E^4$ and~$E\otimes\ad^0$. (Indeed if~$H$ preserves a subspace,
  we obtain a partial flag which is stabilized by the Zariski closure
  of~$H$.) Furthermore every element $h\in \GSp_4(\cO)$ has~$1$ as an
  eigenvalue on~$E\otimes\ad^0$; so if~$h\in H$  is any element whose
  reduction~$\overline{h}$ is regular semi-simple, then
  $h$ satisfies the
  conditions of Definition~\ref{defn: integrally enormous}.
  
  We now consider case~(\ref{inducedcase}).
  Our argument is essentially a characteristic zero version of the proof of~\cite[Lemma~7.5.17]{BCGP} (though note
  that the roles of~$G$ and~$H$ are reversed). 
  If~$h \in (H \setminus G) \cap \Sp_4$, then the eigenvalues of~$h$ 
  are of the form~$(\alpha,\alpha^{-1},-\alpha,-\alpha^{-1})$, and by assumption
  we may assume that the image of~$h$ in~$\Sp_4(k)$
  lands in~$(\Hbar \setminus \Gbar)$ and
  is regular semi-simple --- the latter condition being equivalent to the condition
   that~$\alphabar^4 \ne 1$.
 Now, following the proof of~\cite[Lemma~7.5.17]{BCGP},
  the representation~$E\otimes\ad^0$ decomposes over the algebraic closure of~$E$ into two irreducible
   representations
  of dimension~$6$ and~$4$ 
  on which~$h$ has eigenvalues~$(1,-1,\alpha^2,-\alpha^2,\alpha^{-2},-\alpha^{-2})$
  and~$(-1,1,\alpha^2,\alpha^{-2})$ respectively, both of which contain~$1$ as an eigenvalue.
  \end{proof}

We now let~$F^+$ be a totally real number field in which~$p$ splits completely,
and we fix a continuous representation $\rho:G_{F^+}\to\GSp_4(\cO)$ satisfying the
following hypothesis. 

\begin{hypothesis}\label{hyp: assumptions on rho for mult one}Assume that:
  \begin{enumerate}
  \item $\rho$ is unramified at all but finitely many places.
  \item $\nu\circ\rho=\varepsilon^{-1}$.
  \item $\rhobar:G_{F^+}\to\GSp_4(k)$ is absolutely irreducible.
  \item\label{item: mult one enormous hyp} $\rho(G_{F^+(\zeta_{p^\infty})})$ is integrally enormous.
  \item $\rho$ is  pure. %
  \item\label{item: weight 2 pure p dist} for all~ $v\in S_p$, $\rho|_{G_{F^+_v}}$ is 
    ordinary, semistable of weight~$2$, pure, and $p$-distinguished. %
    We  choose a $p$-stabilization $(\alpha_p,\beta_p)$ of
    $\rho|_{G_{F^+_v}}$ (and thus of $\rhobar|_{G_{F^+_v}}$,  so %
    that by definition ~$\rho|_{G_{F^+_v}}$
     corresponds to a 
    point of $\Spec
           R_{p}^{\triangle}$). %
  \item \label{item:p-2-Q-i-image}If~$p=2$ then
    $\rhobar(G_{F^+(i)})=\rhobar(G_{F^+})$.
  \item \label{item:p>2-regular-semi-simple}If~$p>2$ then
    $\rhobar(G_{F^+})\setminus \Sp_4 (\Fp)$ contains a regular semi-simple element.
  \end{enumerate}
\end{hypothesis}
We will use the following result to verify condition~\eqref{item: mult one
  enormous hyp} in Hypothesis~\ref{hyp: assumptions on rho for mult one}.
\begin{cor}[Examples of integrally enormous representations, II]
 \label{cor: still integrally enormous up the pro p tower}  
Suppose that  $\rho:G_{F^+}\to\GSp_4(\cO)$ satisfies either of the following two sets of conditions:
  \begin{enumerate}[label=(A\arabic*)]
  \item the
    Zariski closure of~$\rho(G_{F^+})$ contains~$\Sp_4$, and
\item $\rhobar(G_{F^+(\zeta_p)})$ contains a regular semi-simple
    element. 
     \end{enumerate}
  Or alternatively:
  \begin{enumerate}[label=(B\arabic*)]
       \item $\rho$ is induced from a quadratic extension~$K/F^{+}$ disjoint from
       the compositum~$F_{\infty}$ of~$F^{+}(\zeta_{p^{\infty}})$ with the fixed field of the similitude character.
The
  Zariski closure of~$\rho(G_{F^{+}})$ contains~$\SL_2 \times \SL_2$, and
  \item $\rhobar(G_{F_{\infty}}) \setminus \rhobar(G_{K.F_{\infty}})$ contains a regular semi-simple
    element.
  \end{enumerate}
 Then $\rho(G_{F^+(\zeta_{p^\infty})})$ is integrally enormous.
\end{cor}
\begin{proof}
We assume the first set of conditions.
  Since the regular semi-simple elements have order prime to~$p$, and
 since $F^+(\zeta_{p^\infty})/F^+(\zeta_p)$ is a pro-$p$ extension,
  $\rhobar(G_{F^+(\zeta_{p^\infty})})$ contains a regular semi-simple
  element. Since the Zariski closure of $\rho(G_{F^+})$
  contains~$\Sp_4$ by assumption, so does the Zariski closure of
  $\rho(G_{F^+(\zeta_{p^\infty})})$ (note that taking the derived
  subgroup is compatible with taking the Zariski closure,
  by~\cite[2.1(e)]{MR1102012}). The result follows from Lemma~\ref{lem: sp4 is integrally enormous}(\ref{sp4case}).
  
  Now we assume the second set of conditions. The property of being integrally enormous is inherited from
  subgroups so it suffices to show that~$H := \rho(G_{F_{\infty}})$ is
  integrally enormous. By construction~$H \subset \Sp_4(\OL)$
  is absolutely irreducible but becomes reducible after restriction to~$G = \rho(G_{K.F_{\infty}})$.
  Moreover, as in the first case, we deduce that the Zariski closure of~$H$ contains~$\SL_2 \times \SL_2$.
 Hence the result follows from Lemma~\ref{lem: sp4 is integrally enormous}(\ref{inducedcase}).
\end{proof}

We now use the notation for deformation rings that was introduced in  Section~\ref{subsec: notation and defns
  Galois}. Let~$S$ be a finite set
of places of~$F^+$ containing~$S_p$ and the places where~$\rho$ is ramified. %
 We define a global deformation problem~$\cS$ by
\numequation\label{defn: global deformation problem for char 0 mult one} \cS = (\rhobar, S, \{\Lambda_v\}_{v\in S_p}\cup\{\cO\}_{v\in
    S\setminus S_p},\{\cD_v^{\triangle}\}_{v\in S_p}\cup
  \{\cD_v^\square\}_{v\in S\setminus S_p}),\end{equation}where $\Lambda_v=\Lambda_{\GSp_4,v}$
and~$\cD_v^{\triangle}$  is as in Section~\ref{subsec: ord
  def rings}. %
    As in Section~\ref{subsec: notation and defns
  Galois}, we write $\Lambda_{\GSp_4,F^+} = \widehat{\otimes}_{v \in
  S_p} \Lambda_v$. Given a nonempty subset $T\subseteq S$, which we
assume contains~$S_p$, we fix an extension of~ $\rho$
to a $T$-framed lifting $(\rho,\{\gamma_v\}_{v\in T})$
of~$\rhobar$ (i.e.\ fix choices of $\gamma_v\in\widehat{\GSp}_4(\cO)$
for each $v\in T$).  Write~$\mathfrak{q}$ (resp.\ $\mathfrak{q}_{S,T}$) for the kernel of the
homomorphism $R_{\cS}\to\cO$ (resp.\ $R_{\cS}^T\to\cO$) corresponding to~$\rho$, and write
$\mathfrak{q}^{\loc}_T$ for the kernel of the composite
$R_{\cS,T}^{\loc}\to R_{\cS}^T\to\cO$.

Write~$W=\ad\rho$, $W_E=W\otimes_{\cO}E$, $W_{E/\cO}=W_E/W$, and for
each $m\ge 1$ write $W_m=W\otimes_{\cO}\cO/\varpi^m$; and we similarly
write~$W^0=\ad^0\rho$, $W^0_{E/\cO}=W^0\otimes_\cO{E/\cO}$,
$W^0_m=W^0\otimes_{\cO}\cO/\varpi^m$.  We write
\[H^i(F^+_S/F^+,W)':=\im\biggl(H^i(F^+_S/F^+,W^0)\to
  H^i(F^+_S/F^+,W)\biggr),\]and similarly we write $H^i(F^+_S/F^+,W_E)'$
and so on. For any place~$v$ of~$F^+$ we also have cohomology groups
$H^i(F_v^+,W)'$ etc.\ defined in the analogous way. We write
$h^i(F_S^+,W_m)'$ for the length of the finite $\cO$-module
$H^i(F_S^+,W_m)'$, and similarly for $h^i(F_v^+,W_{m})'$.

\begin{rem}
  \label{rem: it doesn't matter which lattices we choose}Write~$\mathfrak{z}$
  for the centre of~$\ad$. Since $\ad^0$
  and $\ad/\mathfrak{z}$ are isogenous,
   the quantities
  $h^i(F_S^+,W_m)'-h^i(F_S^+,W^0_m)$ are bounded independently of~$m$
  (and similarly for $h^i(F_v,W_m)'-h^i(F_v,W^0_m)$, and for the
  various Selmer groups introduced below).
\end{rem}%
We define Selmer groups
$H^1_{\cS,T}(F^+,W_m)$  by
\[ H_{\cS,T}^1(F^+,W_m):= \ker \left(H^i(F^+_S/F^+,W_m)'\to
    \prod_{v\in T} H^i(F_v^+,W_m)' \right) .\] We write
$h_{\cS,T}^1(F^+,W_m)$ for the length of $H_{\cS,T}^1(F^+,W_m)$. In
the same way we define \[ H_{\cS,T}^1(F^+,W):= \ker \left(H^i(F^+_S/F^+,W)'\to
    \prod_{v\in T} H^i(F_v^+,W)' \right),\]  \[ H_{\cS,T}^1(F^+,W_E):= \ker \left(H^i(F^+_S/F^+,W_E)'\to
    \prod_{v\in T} H^i(F_v^+,W_E)' \right),\] \[ H_{\cS,T}^1(F^+,W_{E/\cO}):= \ker \left(H^i(F^+_S/F^+,W_{E/\cO})'\to
    \prod_{v\in T} H^i(F_v^+,W_{E/\cO})' \right).\] Note that if we
take direct limits via the
injections $W_m \cong \varpi W_{m+1} \subset W_{m+1}$ we obtain \[H^1_{\cS,T}(F^+,W_{E/\cO}) = 
\varinjlim_m H^1_{\cS,T}(F^+, W_m), \]
while the Mittag-Leffler property means that if we take inverse limits 
with respect to the projection maps $W_{m+1}\to W_m$ we have\[H^1_{\cS,T}(F^+,W) = 
\varprojlim_m H^1_{\cS,T}(F^+, W_m)\] and
thus \[H^1_{\cS,T}(F^+,W_E) = 
  \left(\varprojlim_m H^1_{\cS,T}(F^+, W_m)\right)\otimes_{\cO}E.\]

\begin{lem}%
  \label{lem: length of the tangent space}Assume that~$S_p\subseteq
  T$. For each~$m\ge 1$, the
  length of the
  $\cO/\varpi^m$-module
  \numequation\label{eqn: tangent space at q mod pi m}\mathfrak{q}_{S,T}/(\q_{S,T}^2,\q_T^{\loc}\cdot
    R_{\cS}^T,\varpi^m)\end{equation} %
  is %
  \[h_{\cS,T}^1(F^+,W_m)+\sum_{v\in T}h^0(G_{F^+_v},W_m)-h^0(F^+_S,W_m).\]%
\end{lem}
\begin{proof}%
The length of~\eqref{eqn: tangent space at q mod pi m} is equal to the
length of its $\cO/\varpi^m$-dual, which equals\[\Hom_{\cO}(\mathfrak{q}_{S,T}/(\q_{S,T}^2,\q_T^{\loc}\cdot
    R_{\cS}^T),\cO/\varpi^m),\] and  elements of this latter group
  correspond to strict equivalence classes of $T$-framed liftings ~$(\rho',\{%
  \gamma'_v\}_{v\in T})$ of type~$\cS$
    of~$\rhobar$ to the ring 
~$\cO\oplus \epsilon \cO/\varpi^m$, which furthermore satisfy:
  \begin{enumerate}
  \item\label{item: rho' lifts rho} $\rho'\pmod{\epsilon}$ is strictly equivalent to $(\rho,\{\gamma_v\}_{v\in T})$, and
  \item\label{item: gamma' lifts gamma} for each~ $v\in T$,  %
    $(\gamma'_v)^{-1}(\rho'|_{G_{F^+_v}})\gamma'_v=\gamma_v^{-1}(\rho|_{G_{F^+_v}})\gamma_v$.
  \end{enumerate}

  (Recall that since we have assumed that~$S_{p}\subseteq T$, we are
  not imposing any conditions at the places~$v\in S\setminus T$.)
  Suppose that $[\phi]\in H_{\cS,T}^1(F^+,W_m)$, and let~$\phi$ be a
  lift of~$[\phi]$ to $Z^1(F_S^+,W_m)$. By definition, for each
  $v\in T$ we can write $\phi|_{G_{F_v^+}}=d \psi_v$ for
  some~$\psi_v\in Z^0(F^+_v,W_m)=W_m$. Then
  $((1+\epsilon\phi)\rho,\{(1-\epsilon \psi_v)\gamma_v\}_{v\in T})$
  defines a $T$-framed lifting of~$\rhobar$, and we claim that the
  strict equivalence classes satisfying the two conditions above are
  exactly those containing
  $((1+\epsilon\phi)\rho,\{(1+\epsilon (a_v-\psi_v))\gamma_v\}_{v\in
    T})$ for some $[\phi]\in H_{\cS,T}^1(F^+,W_m)$ and
  $a_v\in H^0(G_{F^+_v},W_m)$. (Implicit in the claim is that this set
  of strict equivalence classes does not depend on the choices of
  liftings~$\phi$ 
  and elements~$\psi_v$.) Given this claim, the
  lemma follows immediately (because the only strict equivalences
  between such $T$-framed liftings are given by replacing all
  the~$a_v$ with $a_v+a$ for some $a\in H^0(F_S^+,W_m)$).

To establish the claim, we firstly check that
$((1+\epsilon\phi)\rho,\{(1-\epsilon \psi_v)\gamma_v\}_{v\in T})$ is a lift of the
required form. Condition~\eqref{item: rho' lifts rho} is obvious, while~\eqref{item: gamma' lifts gamma} is
equivalent to asking that $\phi|_{G_{F^+_v}}(g)=\psi_v-g\psi_v$ for all~$v\in
T$, $g\in G_{F^+_v}$, which is true by the choice of
the~$\psi_v$.

Conversely, if $(\rho',\{ \gamma'_v\}_{v\in T})$ satisfies~\eqref{item: rho' lifts rho},
 then after replacing~$\rho'$ by a strictly equivalent
representation, we can and do assume that
~$(\rho',\{ \gamma'_v\}_{v\in
  T})\pmod{\epsilon}=(\rho,\{\gamma_v\}_{v\in T})$. We may write~
$\rho'=(1+\epsilon\phi)\rho$ for some~$\phi \in Z^1(F_S^+,W_m)'$, and
write~$\gamma'_v=(1+\epsilon a_v)\gamma_v$ with~$a_v\in W_m$. Then~\eqref{item: gamma' lifts gamma}
 is equivalent to
asking that for each~$v\in T$ and
we have $\phi|_{G_{F^+_v}}=-da_v$. Thus we require that $[\phi]\in
H_{\cS,T}^1(F^+,W_m)$,  and given this, the possible~$a_v$ differ by
elements of $H^0(G_{F^+_v},W_m)$, as required.
\end{proof}

We define dual Selmer groups as follows. We let~$W^{0*}$ be the
$\cO$-module dual of~$W^0$, so that $W^{0*}_m:=W^{0*}/\varpi^m$ is the
$\cO/\varpi^m$-module dual of~$W_m$ (and similarly~$W^{0*}_{E/\cO}
:= %
W^{0*}\otimes_{\cO}E/\cO \simeq \Hom(W^0,E/\cO)$).
Then we set \[ H_{\cS^\perp,T}^1(F^+,W^{0*}_m(1)):= \ker \left(H^i(F^+_S/F^+,W^{0*}_m(1))\to
    \prod_{v\in S\setminus T} H^i(F_v^+,W^{0*}_m(1)) \right),\] and similarly \[ H_{\cS^\perp,T}^1(F^+,W^{0*}(1)):= \ker \left(H^i(F^+_S/F^+,W^{0*}(1))\to
    \prod_{v\in S\setminus T} H^i(F_v^+,W^{0*}(1)) \right),\]
\[ H_{\cS^\perp,T}^1(F^+,W_E^{0*}(1)):= \ker \left(H^i(F^+_S/F^+,W_E^{0*}(1))\to
    \prod_{v\in S\setminus T} H^i(F_v^+,W_E^{0*}(1)) \right),\]
\[ H_{\cS^\perp,T}^1(F^+,W_{E/\cO}^{0*}(1)):= \ker \left(H^i(F^+_S/F^+,W_{E/\cO}^{0*}(1))\to
    \prod_{v\in S\setminus T} H^i(F_v^+,W^{0*}_{E/\cO}(1)) \right).\] Just as for the
Selmer groups, these satisfy the obvious compatibilities with direct and
inverse limits.

We now show the existence of appropriate sets of Taylor--Wiles primes. We closely follow
the proofs of~\cite[Cor.\ 2.21, Lem.\ 2.26, Cor.\ 2.27]{MR4592862},
beginning with the following lemma.
\begin{lem}\label{lem: killing a dual Selmer class rationally}Suppose
  that~$\rho:G_{F^+}\to\GSp_4(\cO)$ satisfies Hypothesis~\ref{hyp: assumptions on rho for mult one}.
  Let
  $q:=\corank_{\cO} H_{\cS^\perp,S}^1(F^+,W_{E/\cO}^{0*}(1))$. %

  Then there exist 
  $\sigma_1,\ldots,\sigma_q \in 
        G_{F^+(\zeta_{p^\infty})}$ such 
        that
        \begin{enumerate}
        	\item[(a)] for each $i = 1, \dots, q$,
                  $\rhobar(\sigma_i)$ is regular semi-simple, and
        	\item[(b)] the kernel of the map %
                  \begin{align*}H_{\cS^\perp,S}^1(F^+,W_{E/\cO}^{0*}(1))&\to \bigoplus_{i=1}^q 
        	H^1(\widehat{\Z}, W^{0*}_{E/\cO}(1))\\ &\cong
                \bigoplus_{i=1}^q W^{0*}_{E/\cO}(1) / (\sigma_i - 1) W^{0*}_{E/\cO}(1)
                  \end{align*} \emph{(}the product of the restriction
                  maps $\Res^{G_{F^+, S}}_{\langle \sigma_i \rangle}$ associated to the
                  homomorphisms $\widehat{\Z} \to G_{F^+, S}$, the
                  $i^\text{th}$ such homomorphism sending $1$ to
                  $\sigma_i$\emph{)} is a finite length $\cO$-module.%

  \end{enumerate}%
\end{lem}
\begin{proof}Since $H_{\cS^\perp,S}^1(F^+,W_{E/\cO}^{0*}(1))$ is cofinitely
  generated, it suffices (by an obvious inductive construction) to show that for any non-zero homomorphism
  $f: E / \cO \to H_{\cS^\perp,S}^1(F^+,W_{E/\cO}^{0*}(1))$, we can find
  $\sigma \in G_{F^{+}(\zeta_{p^\infty})}$ such that $\rhobar(\sigma)$ is
  regular semi-simple, and the restriction
  $\Res^{G_{F^+, S}}_{\langle \sigma \rangle} \circ f: E / \cO \to
  W^{0*}_{E / \cO}(1) / (\sigma - 1)W^{0*}_{E / \cO}(1)$ is non-zero.

Let  $L_\infty' / F^+$ be the 
extension cut out by $W^{0*}_E(1)$, and let $L_\infty =
L_\infty'(\zeta_{p^\infty})$. We claim that $H^1(L_\infty / F^+,
W^{0*}_E(1)) = 0$. To see this, note that the extension cut out by $W^*_E(1)=W^{0*}_E(1)\oplus
E(1)$ is~$L_\infty$, and since~$W^*_E(1)$ is pure, it follows from~ \cite[Lemma 
6.2]{MR2076924} that $H^1(L_\infty / F^+,
W^{*}_E(1)) = 0$, and thus $H^1(L_\infty / F^+,
W^{0*}_E(1)) = 0$ as claimed. 
Thus $H^1(L_\infty / F^+, W^{0*}_{E / \cO}(1))$ is killed by a power of 
$p$ (since it injects into the finitely generated $\cO$-module $H^2(L_\infty / F^+, W^{0*}_{\cO}(1))$), 
and hence the homomorphism 
\begin{align*} \Res^{G_{F^+, S}}_{G_{L_\infty,S_{L_\infty}}} \circ f: E / 
\cO &\to H^1(F_S/L_\infty, W^{0*}_{E / \cO}(1))^{G_{F^+, S}} \\&\cong \Hom_{G_{F^+, 
S}}(G_{L_\infty,S_{L_\infty}}, W^{0*}_{E / \cO}(1)) \end{align*}is
still non-zero (here~$S_{L_\infty}$ denotes the set of places
of~$L_\infty$ lying over places in~$S$).

Let $M \subset W^{0*}_{E / \cO}(1)$ be the $\cO$-submodule 
generated by the elements $f(x)(\sigma)$, $x \in E/\cO$, $\sigma \in 
G_{L_\infty}$; it is a non-zero divisible
$\cO[G_{F^+(\zeta_{p^\infty})}]$-submodule of 
$W^{0*}_{E / \cO}(1)$. By the assumption
that~$\rho(G_{F^+(\zeta_{p^\infty})})$ is integrally enormous, we deduce
that there 
exists $\sigma \in G_{F^+(\zeta_{p^\infty})}$ such that
$\rhobar(\sigma)$ is regular semi-simple and $M \not\subset (\sigma-1)W^{0*}_{E / \cO}(1)$. Consequently, 
there exists $m \geq 0$ and $\tau \in G_{L_\infty}$ such that $f(1 / 
\varpi^m)(\tau) \not\in (\sigma - 1)W^{0*}_{E / \cO}(1)$.

If $f(1 / \varpi^m)(\sigma) \not\in (\sigma - 1)W^{0*}_{E / \cO}(1)$,
then $\Res^{G_{F^+, S}}_{\langle \sigma \rangle} \circ f$ is non-zero,
as required.  On the other hand if
$f(1 / \varpi^m)(\sigma) \in (\sigma - 1)W^{0*}_{E / \cO}(1)$ then
$\Res^{G_{F^+, S}}_{\langle \tau \sigma \rangle} \circ f$ is non-zero
(because $f(1 / \varpi^m)(\tau\sigma) =f(1 / \varpi^m)(\tau)+f(1 /
\varpi^m)(\sigma)$ and $(\tau\sigma - 1)W^{0*}_{E / \cO}(1)=(\sigma -
1)W^{0*}_{E / \cO}(1)$).
By construction, 
$\tau \in G_{L_{\infty}} \subset  G_{F^{+}(\zeta_{p^{\infty}})}$
so~$\tau \sigma \in G_{F^{+}(\zeta_{p^{\infty}})}$.
Finally, since~$\tau$ lies in~$G_{L_{\infty}}$, $\rhobar(\tau)$ is scalar and hence~$\rhobar(\tau \sigma)$
is regular semi-simple,
so  we are done.
\end{proof}
\begin{defn}\label{defn: TW datum enormous image}
  A  \emph{set of Taylor--Wiles primes of
    level~$N$} is %
a finite set of finite places $Q$ of $F^+$, disjoint from $S$,
    such that for each $v\in Q$, we have $q_v \equiv 1 \bmod p^N$,
    and $\rhobar(\Frob_v)$ is regular semi-simple.
\end{defn}
Given a set of Taylor--Wiles primes~
$Q$, we define the
augmented global deformation problem
 \[ \cS_Q = (\rhobar, S\cup Q, \{\Lambda_v\}_{v\in S_p}\cup\{\cO\}_{v\in
    (S\cup Q)\setminus S_p},\{\cD_v^{\triangle}\}_{v\in S_p}\cup
  \{\cD_v^\square\}_{v\in (S\cup Q)\setminus S_p}).\]

\begin{lem}
  \label{lem: existence of TW primes rationally} Suppose
  that~$\rho:G_{F^+}\to\GSp_4(\cO)$ satisfies Hypothesis~\ref{hyp: assumptions on rho for mult one}.
  Let
  $q:=\corank_{\cO} H_{\cS^\perp,S}^1(F^+,W_{E/\cO}^{0*}(1))$.

  Then there exist constants ~$d_1,d_2\ge 0$ such that for each~$N\ge 1$ we
  can find a set of Taylor--Wiles primes~
  $Q_N$  of level
  ~$N$ such that 
 \begin{enumerate}
  \item\label{item: right number of TW primes for uniform bounds} $\# Q_N = q$.
\item\label{item: dual Selmer group uniformly bounded}   $h_{\cS_{Q_N}^\perp,S}^1(F^+,W_{N}^{0*}(1))\le d_1$.%
 \item\label{item: Selmer group uniform linear growth} for each $m\le
   N$, the length of \[\mathfrak{q}_{S\cup Q_N,S}/(\q_{S\cup Q_N,S}^2,\q_S^{\loc}\cdot
    R_{\cS_{Q_N}}^S,\varpi^m)\] is at most
  $d_2+m\bigl(2q-4[F^+:\Q]+\#S-1\bigr)$.%
 \end{enumerate}
\end{lem}
\begin{proof}%

Let~$Q$ be any set of Taylor--Wiles primes.
  We have by definition
   the exact sequence \[0\to
    H_{\cS_{Q}^\perp,S}^1(F^+,W_{N}^{0*}(1))\to
    H_{\cS^\perp,S}^1(F^+,W_{N}^{0*}(1))\to\oplus_{v\in
      Q}H^1(k(v),W_{N}^{0*}(1)).\] 
      More generally, for a set~$\QQQ$ of elements~$\sigma_i \in G_{F^{+}(\zeta_{p^N})}$
 with~$\rhobar(\sigma_i)$
      semi-simple,
     let us denote 
      by~$H_{\cS_{\QQQ}^\perp,S}^1(F^+,W_{N}^{0*}(1))$
      the kernel of the map
   \[ H_{\cS^\perp,S}^1(F^+,W_{N}^{0*}(1))\to\oplus_{\QQQ}
   H^1(\langle \sigma_i \rangle,W_{N}^{0*}(1)).\] 
      Accordingly,
      if~$Q = \{v_i\}$ is a set of Taylor--Wiles primes,
      and~$\QQQ = \{\sigma_i\}$ with~$\sigma_i = \Frob_{v_i}$, then
      \[H_{\cS_{\QQQ}^\perp,S}^1(F^+,W_{N}^{0*}(1))=
      H_{\cS_{Q}^\perp,S}^1(F^+,W_{N}^{0*}(1)).\]
Moreover, by the Chebotarev density theorem, for any
      such set~$\QQQ = \{\sigma_i\}$, there exists a set of
      places~$v_{i}$ with~$q_{v_i} \equiv 1 \bmod p^N$
      such that the action of~$\Frob_{v_i}$ on the finite
       module~$W_{N}^{0*}(1)$ coincides with the action
       of~$\sigma_i$, and such that for any~$m\le N$ and any~$f\in
       H_{\cS^\perp,S}^1(F^+,W_m^{0*}(1))$, we have
       $f(\sigma_i)=f(\Frob_{v_i})$. (Here we use that the
       $\cO$-modules $H_{\cS^\perp,S}^1(F^+,W_m^{0*}(1))$ are finitely
       generated and the modules~$W_m^{0*}(1)$ are  finite.) %
       This implies the equality %
           \numequation \label{itsanequality}
         H_{\cS_{\QQQ}^\perp,S}^1(F^+,W_{m}^{0*}(1)) =
      H_{\cS_{Q}^\perp,S}^1(F^+,W_{m}^{0*}(1)).\end{equation}
      for any~$m \le N$.
      Comparing to Lemma~\ref{lem: killing a dual Selmer class
  rationally}, we see that 
  we can and do choose a set of elements~$\QQQ = \{\sigma_i\}$ so that the groups
  $H_{\cS_{\QQQ}^\perp,S}^1(F^+,W_{E/\cO}^{0*}(1))$ are finite length
  $\cO$-modules of uniformly bounded length (indeed they are all
  isomorphic, but we do not need this),
 and let~$Q_N$ denote a corresponding set of Taylor--Wiles primes of level~$N$ so
 that  equality~\ref{itsanequality} holds.    We will now show that these
  sets in fact satisfy properties~\eqref{item: dual Selmer group
    uniformly bounded} and~\eqref{item: Selmer group uniform linear growth}.
  Considering the long exact sequences in Galois cohomology
  associated to the short exact sequence
  \[0\to W_N^{0*}(1)\to
    W_{E/\cO}^{0*}(1)\stackrel{\varpi^N}{\to}W_{E/\cO}^{0*}(1)\to 0\]
  we see that we have a
  morphism
  \[
  H_{\cS_{Q_N}^\perp,S}^1(F^+,W_{N}^{0*}(1)) = 
  H_{\cS_{\QQQ}^\perp,S}^1(F^+,W_{N}^{0*}(1))\to
    H_{\cS_{\QQQ}^\perp,S}^1(F^+,W_{E/\cO}^{0*}(1))[\varpi^N]\] 
    whose
  kernel has order bounded by that of $H^0(F^+,W_{E/\cO}^{0*}(1))/\varpi^N$. Since
  as explained above
  $h_{\cS_{\QQQ}^\perp,S}^1(F^+,W_{E/\cO}^{0*}(1))$ is uniformly
  bounded, in order to prove~\eqref{item: dual Selmer group uniformly
    bounded}, it remains to show that $H^0(F^+,W_{E/\cO}^{0*}(1))$ has
  finite length. It in turn suffices to check that
  $H^0(F^+,W_{E}^{0*}(1))=0$ (because then
  $H^0(F^+,W_{E/\cO}^{0*}(1))$ is cofinitely generated and injects
  into the finitely generated $\cO$-module $H^1(F^+,W^{0*}(1))$). But
  $H^0(F^+,W_{E}^{0*}(1))=H^0(F^+,W_{E}^0(1))$, and we even have
  $H^0(F^+,W_E(1))=0$, because $\rho$ is absolutely irreducible and
  not isomorphic to~$\rho(1)$ (e.g.\ because it has different
  Hodge--Tate weights).

We now turn to~\eqref{item: Selmer group uniform linear
  growth}. Write~$l$ for the length of a finite $\cO$-module. By
the Greenberg--Wiles formula %
together
with Remark~\ref{rem: it doesn't matter which lattices we choose}, the
quantity
\begin{multline*}h_{\cS_{Q_N},S}^1(F^+,W_{m})-h_{\cS_{Q_N}^\perp,S}^1(F^+,W_{m}^{0*}(1))+h^0(F^+,W_{m}^{0*}(1))-h^0(F^+,W_m^0)\\+\sum_{v\in
  S}h^0(F_v^+,W_m^0)+\sum_{v\in Q_N}\left(h^0(F_v^+,W_m^0)-h^1(F_v^+,W_m^0)\right)+\sum_{v|\infty}l((1+c_v)W_m^0)\end{multline*}is uniformly bounded independently of~$N$ and~$m\le N$ (and of our choice
of~$Q_N$).

Comparing to Lemma~\ref{lem: length of the tangent space},  we see that in order to
establish ~\eqref{item: Selmer group uniform linear growth}, it
suffices to show that the quantity %
\begin{multline*}
-h_{\cS_{Q_N}^\perp,S}^1(F^+,W_{m}^{0*}(1)) + h^0(F^+,W_{m}^{0*}(1)) \\+\biggl(h^0(F^+,W_m)-h^0(F^+,W_m^0)-m\biggr)+\sum_{v\in
    S}\biggl(h^0(F_v^+,W_m^0)-h^0(F_v^+,W_m)+m\biggr)\\+\sum_{v\in
    Q_N}\biggl(h^0(F_v^+,W_m^0)-h^1(F_v^+,W_m^0)+2m\biggr)+\sum_{v|\infty}\biggl(l\bigl((1+c_v\bigr)W_m^0)-4m\biggr)
\end{multline*}is uniformly bounded independently of~$N$ and of~$m\le N$. We will do
this by showing that each of the terms is uniformly bounded.

We begin with
$h_{\cS_{Q_N}^\perp,S}^1(F^+,W_{m}^{0*}(1))$. Considering the
morphisms \[W_{m}^{0*}(1)\to W_{N}^{0*}(1)\to W_{E/\cO}^{0*}(1)\] we
obtain a morphism
\[H_{\cS_{Q_N}^\perp,S}^1(F^+,W_{m}^{0*}(1))
= H_{\cS_{\QQQ}^\perp,S}^1(F^+,W_{m}^{0*}(1))
\to
  H_{\cS_{\QQQ}^\perp,S}^1(F^+,W_{N}^{0*}(1))\] whose kernel is
contained in the kernel of the
morphism
\[H^1(F^+_S/F^+,W_{m}^{0*}(1))\to
  H^1(F^+_S/F^+,W_{E/\cO}^{0*}(1)).\] This latter kernel is isomorphic
to a
subquotient of $H^0(F^+,W_{E/\cO}^{0*}(1))$, which we showed above is
a finite $\cO$-module.  The uniform boundedness of
$h_{\cS_{Q_N}^\perp,S}^1(F^+,W_{m}^{0*}(1))$ for~$m \le N$ then follows from that of
$h_{\cS_{Q_N}^\perp,S}^1(F^+,W_{N}^{0*}(1))$, i.e.\ from
~\eqref{item: dual Selmer group uniformly bounded}.

To show that the term $h^0(F^+,W_{m}^{0*}(1))$ is uniformly bounded,
we recall from above that $H^0(F^+,W_E^{0*}(1))=0$. It follows that we
have an injective map:
\[H^0(F^+,W_m^{0*}(1))\into H^1(F^+,W^{0*})[\varpi^m], \]
 and we are done
because $H^1(F^+,W^{0*})[\varpi^m]$ is uniformly bounded (since
$H^1(F^+,W^{0*}) $ is a finitely generated $\cO$-module). The
term \[h^0(F^+,W_m)-h^0(F^+,W_m^0)-m\] and the terms for~$v\in S$ can be handled similarly, using
that $H^0(F^+,W_E)=H^0(F^+,W^0_E)\oplus E$ (respectively
$H^0(F^+_v,W_E)=H^0(F^+_v,W^0_E)\oplus E$).

If $v\in Q_N$, then since $v$ splits in~$F^+(\zeta_{p^N})$ and~$m\le N$, the local
Euler characteristic formula and Tate local duality
give \[h^1(F_v^+,W_m^0)-h^0(F_v^+,W_m^0)=h^2(F_v^+,W_m^0)=h^0(F_v^+,W_m^{0*}(1))=h^0(F_v^+,W_m^{0*}),\]and
since~$\rhobar({\Frob_v})$ is regular semi-simple, we have
$h^0(F_v^+,W_m^{0*})=2m$, so these terms vanish identically.

Finally the claim for places~$v|\infty$ follows easily from
$\dim_EH^0(F_v^+,W_E^0)=4$ (which in turn follows from the assumption
that $\nu\circ\rho=\varepsilon^{-1}$).
\end{proof}

\begin{prop}
  \label{prop: actual result of TW systems to use rationally} 
  Suppose that~$\rho:G_{F^+}\to\GSp_4(\cO)$ satisfies
  Hypothesis~\ref{hyp: assumptions on rho for mult one}.  Let
  $q:=\corank_{\cO} H_{\cS^\perp,S}^1(F^+,W_{E/\cO}^{0*}(1))$.

  Then there exists an integer ~$d\ge 0$ such that for each~$N\ge 1$ we
  can find a set of Taylor--Wiles primes~
  $Q_N$  of level
  ~$N$, together with a morphism in~$\CNL_{\Lambda_{\GSp_4,F^+}}$
  \numequation\label{eqn: almost presenting deformation rings}R_{\cS,S}^{\loc}\llb X_1,\ldots,X_g \rrb \rightarrow
  R_{\cS_{Q_N}}^S,\end{equation} such that:
 \begin{enumerate}
  \item $\# Q_N = q$.
  \item $g=2q-4[F^+:\Q]+\#S-1$.
  \item\label{item: cokernel of almost presentation is almost nothing} Let $\q_{S,g}^{\loc}$ be the kernel of the composite morphism \[R_{\cS,S}^{\loc}\llb X_1,\ldots,X_g \rrb \rightarrow
  R_{\cS_{Q_N}}^S\to\cO\] determined by~$\rho$ and~\eqref{eqn: almost presenting deformation rings}.
    Then the finite $\cO$-module  \[\mathfrak{q}_{S\cup Q_N,S}/(\q_{S\cup Q_N,S}^2,\q_{S,g}^{\loc}\cdot
    R_{\cS_{Q_N}}^S)\] %
  is killed by~$\varpi^d$.
 \end{enumerate}
\end{prop}
\begin{proof}
Choose the Taylor--Wiles primes~$Q_N$ as in Lemma~\ref{lem: existence
  of TW primes rationally}, and  let~$d$ be the constant ~$d_2$ in the
statement of that lemma. We can without loss of generality assume
that~$N>d$ (because a set of Taylor--Wiles primes of level~$N+1$ is
also a set of Taylor--Wiles primes of level~$N$). 
If~$M = \mathfrak{q}_{S\cup Q_N,S}/(\q_{S\cup Q_N,S}^2,\q_S^{\loc}\cdot
    R_{\cS_{Q_N}}^S)$, 
    the length of~$M/\varpi^m$
    is at most~$gm+d$ for all~$m \le N$  by 
    Lemma~\ref{lem: existence of TW primes rationally} \eqref{item: Selmer group uniform linear growth}.
    Thus,
     by~\cite[Lem.\ 2.20]{MR4592862}, 
we can
find a map \[\cO^g\to M/\varpi^N = \mathfrak{q}_{S\cup Q_N,S}/(\q_{S\cup Q_N,S}^2,\q_S^{\loc}\cdot
    R_{\cS_{Q_N}}^S,\varpi^N)\] whose cokernel has length at
  most~$d$, and is in particular killed
by~$\varpi^{d}$. By the topological
version of Nakayama's lemma, we can find a presentation~\eqref{eqn:
  almost presenting deformation rings} such that \[\mathfrak{q}_{S\cup Q_N,S}/(\q_{S\cup Q_N,S}^2,\q_{S,g}^{\loc}\cdot
    R_{\cS_{Q_N}}^S,\varpi^N)\] is killed by~$\varpi^d$. Since we are
  assuming that~$N>d$, this implies that \[\mathfrak{q}_{S\cup Q_N,S}/(\q_{S\cup Q_N,S}^2,\q_{S,g}^{\loc}\cdot
    R_{\cS_{Q_N}}^S)\] is killed by~$\varpi^d$, as required.
\end{proof}

\subsection{An abstract patching argument}\label{subsec: abstract TW patching for mult
  one}%

Set~$\Delta_{\infty}:=\Zp^{2q}$. 
Suppose that we have the following data.
\begin{hypothesis}\label{hypothesis: data for the SW NT patching argument}\leavevmode
   \begin{enumerate}
 \item \label{item: SW module M free over Lambda} An $R_{\cS}$-module $M$, which is finite free as a
   $\Lambda_{\GSp_4,F^+}$-module. 
 \item\label{item: DeltaN not too small} For each integer~$N\ge 1$, a finite quotient~$\Delta_N$
   of~$\Delta_{\infty}$, such that
   $\Delta_{\infty}\onto\Delta_{\infty}/p^N$ factors through~$\Delta_{N}$.
  \item\label{item: SN to RN map} For each~$N\ge 1$, a homomorphism $\Lambda_{\GSp_4,F^+}[\Delta_N]\to R_{\cS_{Q_N}}$, and a finite
   $R_{\cS_{Q_N}}$-module $M_N$ which is finite free as a $\Lambda_{\GSp_4,F^+}[\Delta_N]$-module.
 \item \label{item: SW killing DeltaN} Isomorphisms of $\Lambda_{\GSp_4,F^+}$-modules
   \[M_N\otimes_{\Lambda_{\GSp_4,F^+}[\Delta_N]}\Lambda_{\GSp_4,F^+}\isoto
     M,\](where the homomorphism
   $\Lambda_{\GSp_4,F^+}[\Delta_N]\to\Lambda_{\GSp_4,F^+}$ is the
   augmentation map), compatible with the actions of~$R_{\cS_{Q_N}}$
   and~$R_{\cS}$ and the natural homomorphism $R_{\cS_{Q_N}}\to R_{\cS}$.
 \end{enumerate}
\end{hypothesis}

As in Definition~\ref{defn: cT}, we let  $\cT=\Lambda_{\GSp_4,F^+}\llb x_1,\dots,x_{11\#S-1}\rrb$ be  the coordinate ring of
  $(\prod_{v \in T} \widehat{\GSp}_4)/\widehat{\mathbb{G}_m}$ over
  $\Lambda_{\GSp_4,F^+}$.
Let \[S_\infty:=\Lambda_{\GSp_4,F^+}\llb\Delta_{\infty}\rrb\cotimes_{\Lambda_{\GSp_4,F^+}}\cT,\] %
and let  %
$\mathbf{a}_\infty$ be the kernel of the map
$S_\infty\to\Lambda_{\GSp_4,F^+}$ given by sending each element
of~$\Delta_{\infty}$ to~$1$ and each~$x_i$ to~$0$. Write
$S_N:=\Lambda_{\GSp_4,F^+}[\Delta_{N}]\cotimes_{\Lambda_{\GSp_4,F^+}}\cT=\cT[\Delta_N]$,
a quotient of~$S_{\infty}$. Let
\[R_\infty:=R_{\cS,S}^{\loc}\llb X_1,\ldots,X_g \rrb.\]
For each~$N\ge 1$ we write~$R_N:=R_{\cS_{Q_N}}^S$ and $\q_N=\q_{S\cup
  Q_N,S}$, so that~$R_N$ is an $R_\infty$-algebra
via~\eqref{eqn: almost presenting deformation
  rings} and an~$S_N$-algebra via~\eqref{eqn: framed deformation over
  unframed} and Hypothesis~\ref{hypothesis: data for the SW NT patching argument}\eqref{item: SN to RN map}. We
set \[M_N^{\square}:=M_N\otimes_{R_{\cS_{Q_N}}}R_N=M_N\otimes_{\Lambda_{\GSp_4,F^+}}\cT,\]
where the equality follows from~\eqref{eqn: framed deformation over unframed}. Write~$R=R_{\cS}$, which is naturally an
$R_N$-algebra for each~$N$. 

Then Hypothesis~\ref{hypothesis: data for the SW NT patching argument}
implies that:

\begin{itemize}
 \item The $R$-module $M$ is finite free as a
   $\Lambda_{\GSp_4,F^+}$-module. 
  \item For each~$N\ge 1$, we have a homomorphism $S_N\to R_{N}$,
    and~$M^{\square}_{N}$ is an $R_N$-module
  which is finite free as an $S_N$-module.
 \item We have isomorphisms of $\Lambda_{\GSp_4,F^+}$-modules
   \[M_N^{\square}/\mathbf{a}_{\infty}\isoto
     M,\] compatible with the actions of~$R_{N}$
   and~$R$ and the homomorphism $R_N\to R$.
 \end{itemize}

Fix a non-principal ultrafilter~$\cF$ on~$\N$, and
write $\mathbf{\Lambda}:=\prod_{N\in\N}\Lambda_{\GSp_4,F^+}$, and write
$\mathbf{\Lambda}_{x}$ for the localization of~$\mathbf{\Lambda}$
at the prime ideal \[x:=\{(x_N)_{N\in\N}|\exists I\in\cF\textrm{
  s.t. }\forall N\in I, x_N\in\m_{\Lambda_{\GSp_4,F^+}}\}.\] Then we
set \[M_\infty:=\varprojlim_n\biggl(\mathbf{\Lambda}_{x}\otimes_{\mathbf{\Lambda}}\prod_{N\gg
  0}M_N^{\square}/\mathbf{a}_\infty^n\biggr),\]  where the product is over the
cofinite (by Hypothesis~\ref{hypothesis: data for the SW NT patching argument}\eqref{item: DeltaN not too small}) set of~$N$ 
for which $\mathbf{a}_\infty^n\supseteq \ker(S_\infty\to
S_N)$, and we set
\[R^{\textrm{patch}}:=\varprojlim_n\biggl(\mathbf{\Lambda}_{x}\otimes_{\mathbf{\Lambda}}\prod_{N\ge
  1}R_N/\m_{R_N}^n\biggr).\] By for example~\cite[Prop.\ 3.4.16]{geenew}, we
have in particular produced the following structures.

    \begin{itemize}
    \item $\Lambda_{\GSp_4,F^+}$-algebra homomorphisms 
      $S_\infty\to R_\infty\to R^{\textrm{patch}}\onto R$.
    \item A finite free $S_\infty$-module $M_\infty$, together with an
      isomorphism of $\Lambda_{\GSp_4,F^+}$-modules \[M_\infty\otimes_{S_\infty}\Lambda_{\GSp_4,F^+}\isoto M.\]
    \item A commutative diagram of
      $S_\infty$-algebras \[\xymatrix{R_\infty\ar[d]\ar[r]&\End_{S_\infty}(M_\infty)\ar[d]^{-\otimes_{S_\infty}\Lambda_{\GSp_4,F^+}}\\
          R\ar[r]&\End_{\Lambda_{\GSp_4,F^+}}(M)} \]
    \end{itemize}
    Write~$\q^{\mathrm{patch}}\subset R^{\mathrm{patch}}$ and
    $q_\infty\subset R_\infty$ for the inverse images of~$\q\subset R$
    (i.e.\ the kernels of the composite morphisms
    $R_\infty\to R^{\textrm{patch}}\onto R\to\cO$ corresponding
    to~$\rho$).

\begin{lem}
  \label{lem: patched q is killed pi d} The $\cO$-module
  $\q^{\mathrm{patch}}/((\q^{\mathrm{patch}})^2,\q_\infty)$ is killed
  by~$\varpi^d$, where~$d$ is as in Proposition~\ref{prop: actual
    result of TW systems to use rationally}.
\end{lem}
\begin{proof}
  This is proved in exactly the same way as~\cite[Prop.\
  4.18]{MR4592862}. By Proposition~\ref{prop: actual result
    of TW systems to use rationally}, for each~$N\ge 1$ the
  $\cO$-module $\q_N/((\q_N)^2,\q_\infty)$ is killed by~$\varpi^d$, so
  the cokernel of
  \[\prod_{N\ge 1}\q_\infty/(\q_\infty)^2 \to \prod_{N\ge
      1}\q_N/(\q_N)^2\] is killed by~$\varpi^d$ (here the map
  $\q_{\infty}\to\q_N$ is the one induced by the morphism
  $R_{\infty}\to R_N$). By an identical
  argument to the proof of~\cite[Lem.\ 4.16, 4.17]{MR4592862}, the image of
  $\prod_{N\ge 1}\q_N$ (resp.\ $\prod_{N\ge 1}\q_N^2$) in~$R^{\mathrm{patch}}$
  is~$\q^{\mathrm{patch}}$ (resp.\ $(\q^{\mathrm{patch}})^2$).

  It remains to show that the image of
  $\prod_{N\ge 1}\q_\infty/(\q_\infty^2)$ in
  $\q^{\mathrm{patch}}/(\q^{\mathrm{patch}})^2$ agrees with the image
  of $\q_\infty/(\q_\infty)^2$. This follows from~\cite[Lem.\
  4.19]{MR4592862}, exactly as in the proof of~\cite[Prop.\
  4.18]{MR4592862}.
\end{proof}

Our abstract freeness result is the following, where for the
convenience of the reader we have recalled the running hypotheses and
notation in the statement.
\begin{prop}\label{prop: abstract mult one char 0}
  Assume that $\rho$ satisfies Hypothesis~\ref{hyp:
    assumptions on rho for mult one}, and that~$M$ satisfies
  Hypothesis~\ref{hypothesis: data for the SW NT patching argument}. Write~$\mathfrak{q}$ for the kernel of the
homomorphism $R_{\cS}\to\cO$  corresponding to~$\rho$. Then $M_{\q}$ is a finite free
  $(R_{\cS})_{\q}$-module.
\end{prop}
\begin{proof}By Lemma~\ref{lem: criterion for H2 vanishing}\eqref{item: def ring
    regular} and our assumption that~$\rho|_{G_{F_v^+}}$ is ordinary, pure and $p$-distinguished for all~$v|p$, together with~\cite[Lem.\
  7.1.3]{BCGP} and the assumption that~$\rho$ is pure, the
  local ring~$(R_\infty)_{\q_\infty}$ is regular. Write
  $(R_\infty)^{\wedge}_{\q_\infty}$ for its $\q_\infty$-adic
  completion, and similarly $(R^{\patch})^{\wedge}_{\q^{\patch}}$ for
  the $\q^{\patch}$-adic completion of $(R^{\patch})_{\q^{\patch}}$. Then
  the morphism \[(R_\infty)^{\wedge}_{\q_\infty}\to
    (R^{\patch})^{\wedge}_{\q^{\patch}}\] is surjective, because by
  Lemma~\ref{lem: patched q is killed pi d} the relative cotangent
  space  $\q^{\mathrm{patch}}/((\q^{\mathrm{patch}})^2,\q_\infty)$
  vanishes (note that localization at~$\q^{\patch}$ in particular
  inverts~$p$).

  Accordingly we can and do lift the morphism
  $S_\infty\to (R^{\patch})^{\wedge}_{\q^{\patch}}$ to a morphism
  $S_\infty\to(R_\infty)^{\wedge}_{\q_\infty}$.  Since~$M_\infty$ is a
  finite free $S_\infty$-module, we see that the
  $(R_\infty)^{\wedge}_{\q_\infty}$-module
  $(M_\infty)^{\wedge}_{\q_\infty}$ has depth
  $\dim S_\infty-1$. By the definition of~$g$ together with Lemma~\ref{lem: criterion for H2 vanishing}\eqref{item: def ring
    regular}%
   and~\cite[Thm.\
  3.3.3]{Bellovin},  we have $\dim
  (R_\infty)^{\wedge}_{\q_\infty}=\dim S_\infty-1$. Since
  $(R_\infty)^{\wedge}_{\q_\infty}$ is regular, it follows from
  Auslander--Buchsbaum that $(M_\infty)^{\wedge}_{\q_\infty}$ is a
  finite free $(R_\infty)^{\wedge}_{\q_\infty}$-module. Quotienting
  by~$\mathbf{a}_\infty$, we conclude that $M^{\wedge}_{\q}$ is a finite free
  $R^{\wedge}_\q$-module, and equivalently that $M_\q$ is a  finite free
  $R_\q$-module, as required.
\end{proof}

\subsection{Higher Hida theory}\label{higherhidaRT}  
We now specialize to the case~$F^+=\Q$ (but still allow $p\ge 2$ to
be arbitrary), and continue to assume that $\rho:G_{\Q}\to\GSp_4(\cO)$
satisfies Hypothesis~\ref{hyp: assumptions on rho for mult one}. In
particular since~$\rho|_{\Qp}$ is ordinary, semistable of weight~$2$, and
$p$-distinguished, and we have chosen a $p$-stabilization, we have
determined an ordered pair $(\alpha_p,\beta_p)$ of distinct elements of~$\cO$
such that \[\rho\cong \begin{pmatrix}
    \lambda_{\alpha_p}&0&*&*\\
0 &\lambda_{\beta_p}&*&*\\
0&0&\varepsilon^{-1}\lambda_{\beta_p}^{-1}&0\\
0&0&0&\varepsilon^{-1}\lambda_{\alpha_p}^{-1}
  \end{pmatrix}\] where~$\lambda_x$ is the unramified character with~$\lambda_x(\Frob_p)=x$.

Our next goal is to explain the construction of the data as in
Hypothesis~\ref{hypothesis: data for the SW NT patching argument}
which we will use to prove our multiplicity one theorems. This amounts
to constructing Taylor--Wiles systems out of higher Hida theory modules, which we
already did for usual (i.e. $\HH^0$) Hida modules in~\cite[\S 7.8, 7.9]{BCGP}, and
we will follow the account there where possible. We do make some changes to the setup however: we specialize to the case~$F^+=\Q$, allow~$p=2$, work
with different level structures at primes where~$\rho$ ramifies, and use a different argument to ensure that we can work at neat level.

We begin by recalling some of the main results of~\cite{boxer2023higher},
specialized to the case of~$\GSp_4 $. For each $w \in \WM$, and neat tame level
$K^p$, in \cite[\S 1.4, 5.4]{boxer2023higher}, we have defined perfect complexes
of $\ZZ_p\llb T(\ZZ_p)\rrb $-modules, which we denote by $M^\bullet_{w,\cusp}$
and $M^\bullet_{w}$. By~\cite[Prop.\ 5.6.3]{boxer2023higher} the complexes
$M^\bullet_{w,\cusp}$ have amplitude in the range $[0, \ell(w)]$,  and the complexes
$M^\bullet_{w}$ have amplitude in the range $[\ell(w), 3]$.  These
complexes have an action of $T(\qq_p)$ (extending the $T(\ZZ_p)$-action) and an
action of the Hecke algebra $\mathbb{T}^p$ of prime to $p$
level. %
We now set %
\[
M_w :=
\begin{cases}
    \HH^{\ell(w)}(M^\bullet_{w, \cusp}) \otimes_{\Zp} \cO, & \text{if } \ell(w) = 0, 1, \\[5pt]
    \HH^{\ell(w)}(M^\bullet_w) \otimes_{\Zp} \cO, & \text{if } \ell(w) = 2, 3.
\end{cases}
\]
  When we want to stress the dependence on the
  tame level~$K^p$, we write
  $M^\bullet_{w, \cusp,K^p}$, $M^\bullet_{w,K^p}$,
  and~$M_{w,K^p}$-respectively. 
We now summarize the main properties of these modules, writing
$\Iw(p)\subset\GSp_4(\Zp)$ for the Iwahori subgroup.
\begin{thm}\label{thm:properties-of-higher-Hida-theory} The following properties hold: %
\begin{enumerate}
\item The modules $M_w$ are finite projective $\cO\llb T(\ZZ_p)\rrb $-modules. 
\item\label{item-higher-Hida-gives-ordinary-cohomology}%
    For any dominant algebraic character  $\lambda \in
  X^\star(T)^+$, let 
  $$\kappa  = -w_{0,M}w(\lambda + \rho)-\rho.$$
   There are Hecke  equivariant isomorphisms \[
M_w \otimes_{\cO\llb T(\ZZ_p)\rrb, \lambda} E(-\lambda) =
\begin{cases}
    \HH^{\ell(w)}(\Sh_{\Iw(p)K^p}^{\tor}, \omega^\kappa(-D))^{\ord}, & \text{if } \ell(w) = 0, 1, \\[5pt]
    \HH^{\ell(w)}(\Sh_{\Iw(p)K^p}^{\tor}, \omega^\kappa)^{\ord}, & \text{if } \ell(w) = 2, 3.
\end{cases}
\]

\item\label{item-higher-Hida-higher-Coleman} For any algebraic character $\lambda\in X^\star(T)$ with
$$\kappa=-w_{0,M}w(\lambda+\rho)-\rho\in X^\star(T)^{M,+}$$
there are Hecke equivariant isomorphisms
\[
M_w \otimes_{\cO\llb T(\ZZ_p)\rrb, \lambda} E(-\lambda) =
\begin{cases}
    \HH_w^{\ell(w)}(\Sh_{K^p}^{\tor}, \omega^{\kappa, \sm}(-D))^{\ord, T(\Z_p)}, & \text{if } \ell(w) = 0, 1, \\[5pt]
    \HH_w^{\ell(w)}(\Sh_{K^p}^{\tor}, \omega^{\kappa, \sm})^{\ord, T(\Z_p)}, & \text{if } \ell(w) = 2, 3.
\end{cases}
\]

\item We have a perfect duality pairing:  \[M_w \otimes_{\cO\llb T(\ZZ_p)\rrb} M_{w_{0,M}ww_0} \rightarrow \cO\llb T(\ZZ_p)\rrb\] interpolating the classical Serre duality.
\item\label{item:projective-over-group-ring} If we have a normal subgroup
  $K_1^p\subseteq K_2^p$ then $M_{w,K_1^p}$ is a finite projective %
  $\cO\llb T(\ZZ_p)\rrb[K_2^p/K_1^p]$-module and the pullback and trace induce isomorphisms
$$M_{w,K_2^p}\isoto M_{w,K_1^p}^{K_2^p/K_1^p}\qquad (M_{w,K_1^p})_{K_2^p/K_1^p}\isoto M_{w,K_2^p}.$$
\end{enumerate}
\end{thm}

\begin{proof}
 By duality (\cite[Thm 5.5.2]{boxer2023higher}), point (1) follows from the vanishing theorem \cite[Prop.\
5.6.3]{boxer2023higher} combined with the fact that when $\ell(w)=1$ the $\HH^0$
with support vanishes.  %
Point (2) is the classicality theorem \cite[Cor 4.5.5]{boxer2023higher} combined
with the control theorem \cite[Thm 5.3.5]{boxer2023higher}.  Point (3) is the
comparison between higher Hida and Coleman theories \cite[Thm
6.2.9]{boxer2023higher} together with Theorem
\ref{thm-compare-Higher-Coleman-old-new}.  As already remarked, the duality in point (4) is \cite[Thm 5.5.2]{boxer2023higher}.

We give some justification for point (5) as this is not explained in
\cite{boxer2023higher}. The projectivity is a consequence of the other
statements and Lemma~\ref{lem:equivalent-conditions-to-be-projective-over-group-ring} below.  By the vanishing results just recalled it suffices to show that the corresponding statements for the higher Hida complexes.  In other words, pullback and trace induce quasi-isomorphisms
$$M_{w,K_2^p}^\bullet\to R\Gamma(K_2^p/K_1^p,M_{w,K_1^p}^\bullet),\quad M_{w,K_1^p}\otimes^L_{\cO\llb T(\ZZ_p)\rrb[K_2^p/K_1^p]}\cO\llb T(\ZZ_p)\rrb\to M_{w,K_2^p}^\bullet,$$
and the same statement holds for the cuspidal complexes.  For this, the key
point is to show that for a suitable choice of cone decomposition $\Sigma$, for
the integral models $\pi:\Sh_{K_pK_1^p,\Sigma}^{\tor}\to
\Sh_{K_pK_2^{p},\Sigma}^{\tor}$ considered in \cite{boxer2023higher} there is an
action of $K_2^p/K_1^p$ on $\Sh_{K_pK_1^p,\Sigma}^{\tor}$, and the natural map
$$\oscr_{\Sh_{K_pK_2,\Sigma}^{\tor}}\to R\Gamma(K_2^p/K_1^p,R\pi_\star\oscr_{\Sh_{K_pK_1^p,\Sigma}^{\tor}})$$
is an isomorphism, and the analogous statements for the sheaves
$\oscr_{\Sh_{K_pK_i,\Sigma}^{\tor}}(-D)$, and for traces.  For this, the
facts that
$R\pi_\star\oscr_{\Sh_{K_pK_1^p,\Sigma}^{\tor}}=\pi_\star\oscr_{\Sh_{K_pK_1^p,\Sigma}^{\tor}}$
and
$\pi_\star\oscr_{\Sh_{K_pK_1^p,\Sigma}^{\tor}}^{K_2^p/K_1^p}=\oscr_{\Sh_{K_pK_2^p,\Sigma}^{\tor}}$
are standard and follow from the local description of the toroidal
boundary (see \cite[Prop.\ 7.5]{Lan2016} for example).  We just need to explain why 
there is no higher $K_2^p/K_1^p$ cohomology (note that the action of
$K_2^p/K_1^p$ may not be free and $p$ may divide the order of
$K_2^p/K_1^p$); however it again follows from the local description of the toroidal boundary that the stabilizers have order prime to $p$.
\end{proof}

We used the following (presumably standard) lemma above.
\begin{lem}
  \label{lem:equivalent-conditions-to-be-projective-over-group-ring}Let~$R$ be a
  Noetherian local ring with residue field~$k$ of characteristic~$p$. Let~$G$ be
  a finite group with Sylow $p$-subgroup~$H$, and let~$M$ be a finitely
  generated $R[G]$-module. Then the following conditions are equivalent:
  \begin{enumerate}
  \item $M$ is a projective $R[G]$-module.
  \item $M$ is a projective (equivalently free) $R[H]$-module.
  \item  $M\otimes^L_{R[H]}R$ is a free $R$-module concentrated in degree~$0$.
  \item $M\otimes^L_{R[H]}k$ is concentrated in degree~$0$.
\end{enumerate}
\end{lem}
\begin{proof}The equivalence of the first two conditions follows from the usual
  averaging argument to promote a splitting as $R[H]$-modules to a splitting as
  $R[G]$-modules. Since~$R[H]$ is a local ring, the equivalence of the second,
  third and fourth conditions is immediate from the existence of minimal free
  resolutions, and in particular from~\cite[Lem.\ 2.1.7, Prop.\ 2.1.9]{geenew}.  
\end{proof}

\begin{rem}\label{rem:Hecke-action-at-p-on-higher-Hida}
  Let us spell out the Hecke action at $p$.  On
  $M_w \otimes_{\cO\llb T(\ZZ_p)\rrb, \lambda} E(-\lambda)$, $T(\qq_p)$ acts via
  a smooth character, trivial on $T(\ZZ_p)$ (this is the reason for the twist by
  $\lambda$). The isomorphism
  $M_w \otimes_{\cO\llb T(\ZZ_p)\rrb, \lambda} E(-\lambda) =
  \HH^{\ell(w)}(\Sh_{\Iw(p)K^p}^{\tor}, \omega^\kappa)^{\ord}$ matches the action
  of $t \in T^+(\qq_p)$ with the action of the double class $[\Iw(p) t \Iw(p)]$ (note
  that on the left hand side, the unitary action of $T(\qq_p)$ is twisted by
  $-\lambda$).
\end{rem}
\begin{rem}\label{rem-onGalrepHH}Let  $\kappa$ be such that $\kappa  =
  -w_{0,M}w(\lambda + \rho)-\rho$ with~$\lambda+\rho$ being $G$-dominant. Let  $c$ be  an eigenclass in $\HH^{\ell(w)}(\Sh_{\Iw(p)K^p}^{\tor},
  \omega^\kappa)^{\ord}$, corresponding to an automorphic
  representation~$\pi$. The torus $T(\qq_p)$ acts on the Jacquet module of $\pi$ 
  via a smooth
  character $\chi^{sm}_c$. 
 Assume furthermore that~$\pi_{\infty}$ is a discrete
  series representation (which is automatic if~$\lambda$ is sufficiently
  regular), so that we have a Galois
  representation  $\rho_{\pi,p} : G_{\qq} \rightarrow
  \mathrm{GSp}_4(\Qpbar)$ associated to~$\pi$ (see Theorem~\ref{thm:regular-weight-Galois-rep-recalled}).

Then $\rho_{\pi,p} \vert_{G_{\qq_p}}$ is conjugate to a $B(\Qpbar)$-valued
representation, which we can describe explicitly as follows
(see Remark~\ref{rem-recipe-diag}). 
If $\lambda = (\lambda_1, \lambda_2;w)$, the Hodge--Tate weights  of $\rho_{\pi}$  are given in increasing order by  $$\frac{\lambda_1 + \lambda_2-w}{2}, \frac{2-\lambda_1 + \lambda_2-w}{2}, \frac{ 4  + \lambda_1-\lambda_2 -w}{2},  \frac{6-\lambda_1 -\lambda_2- w}{2}.$$ %
If we use the upper triangular Borel in $B(\qq_p)$ in $\mathrm{GSp}_4(\qq_p)$, and use the local class field theory map $ \qq_p^\times \rightarrow G_{\qq_p}^{\ab}$, the  character on the diagonal is given by 
 $$t \mapsto \mathrm{diag}(t^{ \frac{ - \lambda_1-\lambda_2 +w}{2}}, t^{\frac{ -2+ \lambda_1-\lambda_2 +w}{2}}, t^{ \frac{ -4  - \lambda_1+\lambda_2 +w}{2}}, t^{\frac{\lambda_1+ \lambda_2 +w -6}{2}}),~t \in \Zptimes$$
while $p\in\Qptimes$ goes to
$$\mathrm{diag}(\chi^{sm}_c((1,1,p,p))p^{\frac{-\lambda_1 -
      \lambda_2 +w}{2}},  \chi^{sm}_c((p,1,p,1)) p^{\frac{\lambda_1 -
      \lambda_2 +w}{2}}, $$
      $$ \chi^{sm}_c((1,p,1,p)) p^{\frac{-\lambda_1 +
      \lambda_2 +w}{2}}, \chi^{sm}_c((p,p,1,1))p^{\frac{\lambda_1 +
      \lambda_2 +w}{2}}).$$
\end{rem}

In view of Remark~ \ref{rem-onGalrepHH}, we now give a more Galois-theoretic parametrization of the Higher Hida theories. 
 Recall that in Section~\ref{subsec: ord def
  rings} we  defined a complete local Noetherian  $\cO$-algebra $\Lambda_{\GSp_4,\Q} = \cO\llb
(\ZZ_p^\times(p))^2\rrb$. We have two maps $\ZZ_p^\times \rightarrow
\ZZ_p^\times(p)\times \ZZ_p^\times(p)$ (given by the projections to each factor),
corresponding to $\chi_1, \chi_2: \ZZ_p^\times \rightarrow
\Lambda_{\GSp_4,\Q}^\times$, and there is an associated homomorphism $\psi :
\ZZ_p^\times \rightarrow T(\Lambda_{\GSp_4,\Q})$ which is given by $$\psi(z)=
\mathrm{diag}(\chi_1(z), \chi_2(z), z^{-1} \chi_2^{-1}(z), z^{-1}
\chi_1(z)^{-1}).$$

On the other hand, consider the universal character \[\chi^{un} : T(\ZZ_p)
  \rightarrow  (\ZZ_p\llb T(\ZZ_p)\rrb)^\times.\] Using
Lemma~\ref{lem:cochar-T-to-Lambda} and the identification of $T$ and
$\widehat{T}$, we can view   $\chi^{un} \rho
\nu^{-\frac{3}{2}} : T(\ZZ_p)  \rightarrow  (\ZZ_p\llb T(\ZZ_p)\rrb)^\times$ as a
homomorphism \[\chi^{un} \rho
\nu^{-\frac{3}{2}} :\ZZ_p^\times \rightarrow T( \ZZ_p\llb T(\ZZ_p)\rrb).\]  

\begin{lem}\label{lem:cochar-T-to-Lambda} There is a unique algebra map $f:  \ZZ_p\llb T(\ZZ_p)\rrb \rightarrow \Lambda_{\GSp_4,\Q}$
such that $f \circ \chi^{un} \rho  \nu^{-\frac{3}{2}} = \psi$. 
\end{lem}
\begin{proof}Using
  Lemma~\ref{lem:converting-characters-to-cocharacters-split-torus} and  the
  identification of $T$ and $\widehat{T}$,    $\psi$ corresponds to a character $T(\ZZ_p) \rightarrow
  (\Lambda_{\GSp_4,\Q})^{\times}$, and the lemma follows from the universal
  properties  of
  $\chi^{un}$ and of~$\ZZ_p\llb T(\ZZ_p)\rrb$. 
\end{proof}

We used the following (presumably standard) lemma above.
\begin{lem}
  \label{lem:converting-characters-to-cocharacters-split-torus}If~$T/\Z$  is a
  split torus with dual~$\widehat{T}$, then for any two commutative rings~$R,S$,
  there is a natural bijection between group homomorphisms \[T(R)\to
    S^{\times}\] and group homomorphisms \[R^{\times}\to \widehat{T}(S).\]
\end{lem}
\begin{proof}
  This follows from the case~$T=\Gm$, which is obvious. More precisely,
  since~$X_{*}(\widehat{T})=X^{*}(T)$ is a free $\Z$-module, we
  have \[
\begin{aligned}
    \Hom(T(R), S^{\times}) 
    &= \Hom(X_{*}(T) \otimes R^{\times}, S^{\times}) \\
    &= \Hom(X_{*}(T), \Hom(R^{\times}, S^{\times})) \\
    &= X^{*}(T) \otimes \Hom(R^{\times}, S^{\times}) \\
    &= X_{*}(\widehat{T}) \otimes \Hom(R^{\times}, S^{\times}) \\
    &= \Hom(R^{\times}, X_{*}(\widehat{T}) \otimes S^{\times}) \\
    &= \Hom(R^{\times}, \widehat{T}(S)),
\end{aligned} 
\] as required.
\end{proof}

We now consider the action of the centre of the group $\mathrm{GSp}_4(\mathbb{A}_f)$. The action of $Z(\mathrm{GSp}_4(\mathbb{A}_f))$ on the toroidal compactification of our Shimura varieties factors into an action of $Z(\qq) \backslash Z(\mathrm{GSp}_4(\mathbb{A}_f)) = \{\pm 1 \} \backslash \prod_\ell \ZZ_\ell^\times$ (as can be seen by considering complex uniformization). From a modular perspective, $\ZZ_\ell^\times$ acts on the $\ell$-adic Tate module of an abelian surface by scalar multiplication. The fact that we  quotient by $\{\pm 1\}$ witnesses the fact that   $-1$  is an automorphism of any abelian surface. This action extends to an action on the $M_{w,K^p}$ (as part of the Hecke action on these modules). 

 When we have a $\ZZ_p$-module $M$, equipped with an action of $Z(\qq) \backslash Z(\mathrm{GSp}_4(\mathbb{A}_f))$, we let $M^{\vert.\vert^2}$ be the submodule where the centre acts via the character $(z_\ell)_\ell \mapsto z_p^{-2}$. More generally, if we have a finite set of primes $S$, we let $M^{\vert.\vert^2,S}$ be the submodule where the group $\prod_{\ell \notin S}\ZZ_\ell^\times$ acts via the character $(z_\ell)_{\ell \notin S} \mapsto z_p^{-2}$. 
 \begin{rem} This definition is motivated by Theorem \ref{thm:properties-of-higher-Hida-theory} (\ref{item-higher-Hida-gives-ordinary-cohomology}).  Under the classicality theorem, this condition corresponds to fixing the central character of automorphic forms contributing to coherent cohomology to be $\vert.\vert^2$. 

 By abuse of language, we say that the centre acts by $\vert. \vert^2$ on $M^{\vert.\vert^2}$.  
 \end{rem}

On the module $M_{w,K^p} \otimes_{\cO\llb T(\ZZ_p)\rrb} \Lambda_{\GSp_4,\Q}$, the subgroup $\ZZ_p^\times$ of $Z(\mathrm{GSp}_4(\mathbb{A}_f))$ acts via $z \mapsto z^{-2}$. Therefore, fixing the central action to be $\vert.\vert^2$ amounts to asking that the group $\prod_{\ell \neq p} \ZZ_\ell^\times$ acts trivially. 
We would therefore morally have that $(M_{w,K^p} \otimes_{\cO\llb T(\ZZ_p)\rrb} \Lambda_{\GSp_4,\Q})^{\vert . \vert^2} = M_{w,ZK^p} \otimes_{\cO\llb T(\ZZ_p)\rrb} \Lambda_{\GSp_4,\Q}$ where $ZK^p = \prod_{\ell \neq p } Z(\ZZ_\ell) K_\ell$. Note however that the group $ZK^p$ is not neat (a condition we have imposed so far on our tame level). We now explain a construction which addresses this issue. 

We choose a prime $\vzero >5$ such that $\vzero\not\equiv 1\pmod{p}$ and $\vzero \not\equiv 1 \pmod{4}$ if \(p=2\). We let $\Iw_1
(\vzero)$ denote the subgroup of matrices which are upper-triangular and unipotent modulo~$\vzero$. If $K^p$ is a compact open subgroup with $K_{\vzero} = \Iw_1(\vzero)$, then $K^p$ is neat   by~\cite[Lem.\
 7.8.3]{BCGP} (applied to~$v=\vzero$).

For any finite set of primes $S$ (possibly empty) not containing~$p$, we let
$Z^SK^p =   \prod_{\ell \notin S\cup\{ p\} } Z(\ZZ_\ell) K_\ell \prod_{\ell \in
  S} K_\ell$. %

\begin{lemma}\label{lem:freeness-over-the-neat-prime} Let  $K^p_1 \subseteq K^p_2$ be compact open
  subgroups,  $K^p_1$ normal in $K^p_2$, and the   $r$-components
  of~$K_1^p,K_2^p$ both being $\Iw_1(\vzero)$. Let $S$ be a finite set
  of primes, with $\vzero \in S$, and  $p \notin S$.
  \begin{enumerate}
  \item   \label{neata}
 The module $$(M_{w,K^p_1} \otimes_{\cO\llb T(\ZZ_p)\rrb} \Lambda_{\GSp_4,\Q})^{\vert . \vert^2,S}$$ is a finite projective $ \Lambda_{\GSp_4,\Q}[Z^SK^p_2/Z^SK^p_1]$ module.
 \item  \label{neatb} The module $$(M_{w,K^p_1} \otimes_{\cO\llb T(\ZZ_p)\rrb} \Lambda_{\GSp_4,\Q})^{\vert . \vert^2}$$ is a finite projective $ \Lambda_{\GSp_4,\Q}[ZK^p_2/ZK^p_1]$ module.  %
 \end{enumerate}
\end{lemma}
\begin{proof} Since $\vzero  \in S $,  it follows that the groups $Z^{S} K^p_i$ are neat
for~$i=1$, $2$. On the other hand, we have an identification
$$(M_{w,K^p_1} \otimes_{\cO\llb T(\ZZ_p)\rrb} \Lambda_{\GSp_4,\Q})^{\vert . \vert^2,S}
=  M_{w,Z^S K^p_1} \otimes_{\cO\llb T(\ZZ_p)\rrb} \Lambda_{\GSp_4,\Q},$$
  so~\eqref{neata} follows from Theorem \ref{thm:properties-of-higher-Hida-theory}. 
Next, we can take $S= \{\vzero\}$. Then $(M_{w,K^p_1} \otimes_{\cO\llb T(\ZZ_p)\rrb} \Lambda_{\GSp_4,\Q})^{\vert . \vert^2}$  is obtained from $(M_{w,K^p_1} \otimes_{\cO\llb T(\ZZ_p)\rrb} \Lambda_{\GSp_4,\Q})^{\vert . \vert^2,S }$ by considering the invariants for  $(\ZZ/\vzero\ZZ)^\times$. 
Observe that  $Z^SK^p_2/Z^SK^p_1= ZK^p_2/ZK^p_1$, so 
\[(M_{w,K^p_1} \otimes_{\cO\llb T(\ZZ_p)\rrb} \Lambda_{\GSp_4,\Q})^{\vert . \vert^2,S }\] 
is a finite 
projective $ \Lambda_{\GSp_4,\Q}[Z^SK^p_1/Z^SK^p_2]$ module. Note that $(\ZZ/\vzero\ZZ)^\times$ acts via 
$(\ZZ/\vzero\ZZ)^\times/\{ \pm 1\}$ and this group has order prime
to~$p$. Therefore the invariants are a direct factor, as required. %
\end{proof}

\subsection{Taylor--Wiles systems}\label{subsec: defn of GSp4 TW
  system for mult one}We continue to fix a continuous representation $\rho:G_{\Q}\to\GSp_4(\cO)$ satisfying Hypothesis~\ref{hyp: assumptions on rho for mult one}.
\begin{defn}\label{defn: neat prime}%
  A \emph{neat prime for~$\rhobar$} is a  prime~$\vzero>5$ such that
  \begin{itemize}
  \item $\vzero\ne p$,
     \item $\vzero\not\equiv 1\pmod{p}$ and $\vzero \not\equiv 1 \pmod{4}$ if \(p=2\).
  \item $\rho|_{G_{\Q_{\vzero}}}$ is unramified,
   \item $\rhobar(\Frob_{\vzero})$ is regular semi-simple,
  \end{itemize}together with a fixed ordering
$(\alphabar_{\vzero,1},\alphabar_{\vzero,2},\vzero\alphabar^{-1}_{\vzero,2},\vzero\alphabar_{\vzero,1}^{-1})$
  of the eigenvalues of~$\rhobar(\Frob_{\vzero})$.%
\end{defn}
 By Hypothesis~\ref{hyp: assumptions on rho for mult one} parts~\eqref{item:
  mult one enormous hyp}, \eqref{item:p-2-Q-i-image}, and~\eqref{item:p>2-regular-semi-simple} %
we can and do choose a neat prime~$\vzero$. %

\begin{defn}\label{defn: choice of S for mult one char 0}
  We let~$S$ be the union of~$\{r\}$ and the set of primes at
  which~$\rho$ is ramified (so in particular $p\in S$, because
  ~$\varepsilon$ is ramified at~$p$), and choose sets of Taylor--Wiles primes~$Q_N$ as
  in Proposition~\ref{prop: actual result of TW systems to use
    rationally}.
\end{defn}
For any prime~$l$, we let $\Iw(l)$ denote the subgroup of~$\GSp_4 (\Z_l)$
consisting of matrices which are upper-triangular modulo~$l$, and we let $\Iw_1
(l)$ denote the subgroup of matrices which are upper-triangular and unipotent modulo~$l$.
\begin{df} \label{df:levelstructure}
We define an open compact subgroup $K^p = \prod_l K_l$ of $\GSp_4(\A^{\infty,p})$ as follows:
\begin{itemize}
\item If $l \not\in S$, then $K_l = \GSp_4(\Z_l)$.
  \item  $K_{\vzero} = \Iw_1(\vzero)$.
\item If $l \in S\setminus \{\vzero\}$, then we allow any choice of
  open compact
  $K_l\subseteq \GSp_4(\Z_l)$. %
\end{itemize}
We have compact
subgroups $K^p_0(Q_N)$, $K^p_1(Q_N)$ of~$K^p$ given by
\begin{itemize}
\item If $l\not\in Q_N$, then $K^p_0(Q_N)_l=K^p_1(Q_N)_l=K^p_l$.
\item If $l\in Q_N$, then $K^p_0(Q_N)_l=\Iw(l)$, $K^p_1(Q_N)_l=\Iw_1(l)$.
\end{itemize}
 \end{df} 
 These groups are neat by  ~\cite[Lem.\
 7.8.3]{BCGP}.

We let \[\TT^{S}=\bigotimes_{l\not\in 
S}\cO[\GSp_4(\Q_l) \doubleslash \GSp_4(\Z_l)]\]
be the ring of spherical Hecke operators away from the bad places, and similarly we set \[\TT^{S\cup Q_N}=\bigotimes_{l\not\in 
S\cup Q_N}\cO[\GSp_4(\Q_l) \doubleslash \GSp_4(\Z_l)].\] %

We will also make use of some Hecke operators at Iwahori level, which
we recall from \cite[\S 2.4]{BCGP}. Let $l$ be a prime.  Assume
that~$E$ is large enough to contain a square root of~$l$; we fix
such a choice~$l^{1/2}$. The reader can easily check that
nothing before \cite[Lem.\ 2.4.3]{BCGP} makes any use of the
running assumption made there that~$p\ne 2$; indeed, these results are for the most part over~$E$, and
use only that it is a field of characteristic zero
containing~$l^{1/2}$. We define %
\[\cH=\cH(l)=\cO[G(\Q_l) \doubleslash \Iw(l)],\]
\[\cH_1=\cH_1(l)=\cO[G(\Q_l) \doubleslash \Iw_1(l)].\] With~$T$
denoting our usual maximal torus in~$\GSp_4$, we set  \[T(\Z_l)_1:= (\ker T(\Z_l)\to T(\F_l)),\] and exactly
as in~\cite[Prop.\ 2.4.2]{BCGP} we have an injective
homomorphism \numequation\label{eqn: torus to pro v Iwahori} T(\Q_l)/(\ker T(\Z_l)\to T(\F_l))\to
(\cH_1)^{\times}. \end{equation}%

 The injection~\eqref{eqn: torus to pro v
  Iwahori} induces an injective
homomorphism $\cO[T(\Q_l)/T(\Z_l)_1]\to \cH_1$, and we identify
$\cO[T(\Q_l)/T(\Z_l)_1]$ with its image in~$\cH_1$.
\begin{defn}\label{defn: maximal ideal to localize TW primes at}
  Assume that $l\equiv 1\pmod{p}$.  Given elements
  $\alphabar_{l,1},\alphabar_2\in k^{\times}$, we let
  $\mathfrak{m}_{\alphabar_1,\alphabar_2}$ denote the kernel of the
  homomorphism $\cO[T(\Q_l)/T(\Z_l)_1]\to k$ induced by the
  character $T(\Q_l)/T(\Z_l)_1\to k^\times$ sending
  $T(\Z_l)\mapsto 1$, $\beta_{l,0} \mapsto 1$,
  $\beta_{l,1} \mapsto \alphabar_{l,1}$,
  $\beta_{l,2} \mapsto \alphabar_{l,1}\alphabar_{l,2}$.
\end{defn}

 We define the following elements of
$T(\qq_p)$ which act on $M_w$: %
\begin{eqnarray*}
 U_{p,1} &= &\mathrm{diag}(1,1,p,p) \\
 U_{p,2} &= &\mathrm{diag}(p,1,p^2,p) \\
U_{p,0} &= &\mathrm{diag}(p,p,p,p) 
\end{eqnarray*}

We let 
$\m^{\mathrm{an}}\subset\TT^{S}$ be the maximal ideal
corresponding to $\overline{\rho}$; so by definition~$\m^{\mathrm{an}}$ contains~$\varpi$,
and the polynomials $\det(X-\rhobar(\Frob_l))$ and~$Q_l(X)$ are
congruent modulo~$\m^{\mathrm{an}}$
for each $l\not\in 
S$.
 We define
$\m^{\mathrm{an},Q_N}\subset\TT^{S \cup Q_N}$ in the same way. %
We let %
\begin{equation*}
\TT=\TT^{S}[U_{p,0},U_{p,1},U_{p,2}]
\end{equation*}
and
\begin{equation*}
\TT^{Q_N}=\TT^{S \cup Q_N}[U_{p,0},U_{p,1},U_{p,2}],
\end{equation*}
and additionally %
we let $\m\subset
\TT$ be the maximal ideal 
\[
\m=(\m^{\mathrm{an}},U_{p,0}-1,U_{p,1}-\alpha_p,U_{p,2}-\alpha_p\beta_p)\]
and we let $\m^{Q_N}\subset \TT^{Q_N}$ be the maximal ideal
\begin{equation*}
\m^{Q_N}=(\m^{\mathrm{an},Q_N},U_{p,0}-1,U_{p,1}-\alpha_p,U_{p,2}-\alpha_p\beta_p).
\end{equation*}

For each $q\in Q_N$, fix an ordering    $\alphabar_{q,1},\alphabar_{q,2},\alphabar_{q,2}^{-1},\alphabar_{q,1}^{-1}$
of the eigenvalues of $\rhobar(\Frob_q)$. Set
$\Delta_N = \Delta_{Q_N} = \prod_{q\in Q_N} \F_q^{\times}(p)^2$.  We now fix a choice
of~$w \in \WM$, and consider the finite free $\Lambda_{\GSp_4,\Q}$-module %
\numequation\label{eqn: definition of M}
M:=\Hom_{ \Lambda_{\GSp_4,\Q}}((M_{K^p} \otimes_{\cO\llb T(\ZZ_p)\rrb} \Lambda_{\GSp_4,\Q})^{|\cdot|^2},\Lambda_{\GSp_4,\Q})_{\m,\m_{\vzero}},
\end{equation}
and the finite $\Lambda_{\GSp_4,\Q}[\Delta_{Q_N}]$-modules %
\begin{equation*}
M_N:=\Hom_{ \Lambda_{\GSp_4,\Q}}((M_{K^p_1(Q_N)} \otimes_{\cO\llb T(\ZZ_p)\rrb} \Lambda_{\GSp_4,\Q})^{|\cdot|^2},\Lambda_{\GSp_4,\Q})_{\m^{Q_N},\m_{\vzero},\m_{Q_N}},
\end{equation*}
where:
\begin{itemize}
\item $M_{K^p}$ and~$M_{K^p_1(Q_N)}$ denote the modules~$M_w$ of
  Theorem~\ref{thm:properties-of-higher-Hida-theory}, taking~$K^p$ there to be
  respectively our~$K^p$ and~$K^p_1(Q_N)$.
 \item $\Lambda_{\GSp_4,\Q}$ is an $\cO\llb T(\ZZ_p)\rrb$-algebra via Lemma~\ref{lem:cochar-T-to-Lambda}.
\item The localizations
  $\m,{\m}^{Q_N}$ are
  defined above.  
\item The localization $\m_{\vzero}$ and the localization $\m_{Q_N}$ are with respect to the
  maximal ideals~$\m_l$ of the subalgebras $\cO[T(\Q_l)/T(\Z_l)_1]$ of the pro-$l$ Iwahori Hecke
  algebras~$\cH_1(l)$ for~$l\in Q_N\cup\{\vzero\}$ as in Definition~\ref{defn:
    maximal ideal to localize TW primes at}. %
\item The action of~$\Delta_{Q_N}$ on~$M_{K^p_1(Q_N)}$ is induced by
  the actions of $\cO[T(\Q_q)/T(\Z_q)_1]$ for~$q\in Q_{N}$,  by regarding
  $\F_q^{\times}(p)^2$ as the maximal $p$-power quotient of
  $T(\F_q)/Z(\F_{q})$, where~$Z$ denotes the centre
  of~$\GSp_4$. %
\item The superscript~$|\cdot|^2$ denotes that we are fixing the central
  character. %
\end{itemize}
\begin{lem}
  \label{lem: MN is free over DeltaN}
$M$ is a free~$\Lambda_{\GSp_4,\Q}$-module, and $M_N$ is a free $\Lambda_{\GSp_4,\Q}[\Delta_{Q_N}]$-module.
\end{lem}
\begin{proof}
 Since $\Lambda_{\GSp_4,\Q}[\Delta_{Q_N}]$ is a local ring, 
 this follows from Lemma~\ref{lem:freeness-over-the-neat-prime}(\ref{neatb}),
 since $\Iw(q) Z(\Z_q)/\Iw_1(q) Z(\Z_q) \simeq  T(\F_q)/Z(\F_{q})$.
\end{proof}%

The following is essentially~\cite[Lem.\ 7.1.1]{gg}, adapted slightly
to allow~$p=2$; we will use it in the proof of Proposition~\ref{prop:
  map R to T for SW argument}.%
\begin{lem}\label{lem: symplectic representation valued over trace}Let
  $\Gamma$ be a profinite group, and let $S\subset R$ be complete local
  Noetherian rings with $\m_R\cap S=\m_S$ and common residue
  field $k$. Let $\rho:\Gamma\to\GSp_4(R)$ be a
  continuous representation. Suppose that $\rho\text{ mod }\m_R$ is
  absolutely irreducible,
  that $\tr\rho(\Gamma)\subset S$, and that
  $\nu\circ\rho(\Gamma)\subset S^\times$. Then there is a
  $\widehat{\GSp}_4(R)$-conjugate of $\rho$ whose image is
  contained in $\GSp_4(S)$.
  \end{lem}
\begin{proof} By  \cite[Lem.\ 2.1.10]{cht}, there is some~$B\in
  \widehat{\GL}_4(R)$ such that
  $\rho':=B\rho B^{-1}$ is valued in $\GL_4(S)$. Since $J^{-1}\rho
  J=(\nu\circ\rho)\rho^{-t}$, we have \[(BJB^t)^{-1}\rho'
  (BJB^t)=(\nu\circ\rho')(\rho')^{-t}. \] %
By choosing a symplectic basis for the alternating form determined by
$BJB^t$, it follows that $\rho'$ is $\widehat{\GSp}_4(S)$-conjugate to
a representation~$\rho''$ valued in $\GSp_4(S)$. %
By Schur's lemma~\cite[Lem.\ 2.1.8]{cht},
we see the element in~$\widehat{\GL}_4(R)$ conjugating~$\rho$
to~$\rho'$ is necessarily contained in $\widehat{\GSp}_4(R)$, as
required.
\end{proof}

\begin{defn}\label{defn:Hecke-algebra-for-NT-argument}
  Let~$\TT_{\cS}$ (resp.\ $\TT_{\cS_{Q_N}}$) denote the image
  of~$\TT$ (resp.\ of~$\TT^{Q_N}$)
  in~$\End_{\Lambda_{\GSp_4,\Q}}(M)$ (resp.\ in
  $\End_{\Lambda_{\GSp_4,\Q}}(M_N))$. (These are objects
  of~$\CNL_{\Lambda_{\GSp_4,\Q}}$, but under the present hypotheses we
  do not know that these algebras are nonzero (because we do not know
  that~$\rhobar$ is modular).)
\end{defn}

Recall that we have defined the global deformation problems
\[ \cS = (\rhobar, S, \{\Lambda_{\GSp_4,p}\}\cup\{\cO\}_{l\in
    S\setminus\{p\} },\{{\cD}_p^{\triangle}\}\cup
  \{\cD_l^\square\}_{l\in S\setminus \{p\}}),\]
 \[ \cS_{Q_N} = (\rhobar, S\cup Q_N,\{\Lambda_{\GSp_4,p}\}\cup\{\cO\}_{l\in
    (S\cup Q_N)\setminus\{p\}},\{{\cD}_p^{\triangle}\}\cup
  \{\cD_l^\square\}_{l\in (S\cup Q_N)\setminus \{p\}}).\] The
deformation ring $R_{\cS_{Q_N}}$ is a
$\Lambda_{\GSp_4,\Q}[\Delta_{Q_N}]$-algebra via
Lemma~\ref{lem:anyliftisTW} (i.e.\ $\Delta_{Q_N}$ acts via the
characters $\gamma_{q,i}\circ \Art_{\Q_q}$).
For ease of notation, we sometimes (e.g.\ in the statement of the
following proposition) adopt the convention that $Q_0=\emptyset$, so
that for example $\TT_{\cS_{Q_N}}=\TT_{\cS}$. %

\begin{prop}
  \label{prop: map R to T for SW argument}For each~$N\ge 0$, the
  action of~$\Delta_{Q_N}$ on~$M_N$ makes $\TT_{\cS_{Q_N}}$ a
  $\Lambda_{\GSp_4,\Q}[\Delta_{Q_N}]$-algebra, and  there is
  a $\Lambda_{\GSp_4,\Q}[\Delta_{Q_N}]$-algebra homomorphism
  $R_{\cS_{Q_N}}\to\TT_{\cS_{Q_N}}$ with corresponding representation
  $\rho_{\cS_{Q_N}}:G_{\Q}\to\GSp_4(\TT_{\cS_{Q_N}})$ determined by
  the property that
 $\det(X-\rho_{\cS_{Q_N}}(\Frob_l))=Q_l(X)$ for all   $l\notin S\cup Q_N$.
\end{prop}
\begin{proof}Since $M_N$ is a finite free
  $\Lambda_{\GSp_4,\Q}$-module, this follows from
  local-global compatibility for the Galois representations at a dense
  set of points of regular weight. More precisely, it can be proved in exactly the same
  way as the case $I=\emptyset$ of~\cite[Thm.\ 7.9.4]{BCGP}, %
  using Lemma~\ref{lem: symplectic
    representation valued over trace} in place of~\cite[Lem.\
  7.1.1]{gg}. %
  (The fact that $\TT_{\cS_{Q_N}}$ is automatically a
  $\Lambda_{\GSp_4,\Q}[\Delta_{Q_N}]$-algebra was not recorded in
  \cite{BCGP}; it follows from local-global compatibility at the
  places in~$Q_N$.)  %
\end{proof}
In particular, Proposition~\ref{prop: map R to T for SW argument}
makes $M$ into an $R_{\cS}$-module, and each $M_N$ into an
$R_{\cS_{Q_N}}$-module. We can also regard~$M$ as an
$R_{\cS_{Q_N}}$-module via the natural map $R_{\cS_{Q_N}}\to R_{\cS}$. For each~$N \ge 1$ we fix a surjection
$\Delta_\infty \onto \Delta_{Q_N}$, and write~$\Delta_N$ for the
corresponding quotient of~$\Delta_{\infty}$.  The
kernel of this surjection is contained in $(p^N\Zp)^{2q}$, since each
$v\in Q_N$ satisfies $q_v \equiv 1 \bmod p^N$. At this point we have
established points \eqref{item: SW module M free over
  Lambda}--\eqref{item: SN to RN map} of Hypothesis~\ref{hypothesis: data for the
  SW NT patching argument}, so it only remains to check~\eqref{item:
  SW killing DeltaN}, which is the content of Lemma~\ref{lem: TW
  system going from Q1 to Q0} below.%

Before proving it, we recall some standard facts about Iwahori Hecke
algebras that were explained in~\cite[\S 2.4]{BCGP} under the
unnecessary assumption that~$p\ne 2$. Indeed the only place in~\cite[\S
2.4]{BCGP}  that relies on the assumption~$p\ne 2$ is  \cite[Lem.\
2.4.34]{BCGP}, i.e.\ the statement that the spherical invariants
of a module for the Iwahori Hecke algebra are a direct summand. This
is no longer valid for~$p=2$, but we will avoid this problem by making
use of Lemma~\ref{lem:freeness-over-the-neat-prime}. We do need to
make use of the (proofs of) \cite[Lem.\
2.4.36, 2.4.37]{BCGP}, but for the convenience of the reader we
will recall the necessary arguments as we use them (making it clear as
we do so that they remain valid for~$p=2$).

Suppose now that~$q$ is a prime with $q\equiv 1\pmod{p}$, and let
~$\cH=\cH(q)$ be the corresponding Iwahori Hecke algebra. The Bernstein
presentation of~$\cH$ is valid for all~$p$, so we can write
\[\cH= \cO[X_*(T)]\widetilde{\otimes}_{\cO}\cO[\Iw(q)\backslash
  \GSp_4(\Z_q)/\Iw(q)],\]where 
  the twisted tensor product is
determined by the relations~\cite[(2.4.32)]{BCGP}.  The centre
of~$\cH$ is $\cO[X_*(T)]^W$ (where as usual~$W\cong D_8$ is the Weyl group of~$\GSp_4$),
\[\cO[X_*(T)]^W= \cO[\GSp_4(\Q_q)\doubleslash \GSp_4(\Z_q)]\]
given by $x\mapsto [\GSp_4(\Z_q)]x$ (where we are regarding~$x$
as an element of~$\cH$); this isomorphism agrees with the usual Satake
isomorphism. (Indeed this presentation, and the compatibility with the
Satake isomorphism, are valid over any ring containing an invertible
square root of~$q$; see for example \cite{MR2122539} or the very
general results of~\cite{boumasmoud2021tale}.)

Since we are assuming that~$q=1$
in~$k$, we deduce exactly as in~\cite[Lem.\ 2.4.33]{BCGP} that
reduction modulo~$\varpi$ induces a natural
isomorphism \[\cH\otimes_{\cO}k\cong k[X_*(T)\rtimes W].\] Since we are assuming
that~$q\equiv 1\pmod{p}$, 
and in any case in our applications of these results in the global setting there
is a twist which makes all of the powers of~$q$ integral, we will ignore all powers of~$q^{1/2}$
from now on.

Exactly as
in~\cite{BCGP}, %
we let~$x_0$,
$x_1$, and~$x_2$ denote the following three cocharacters:
$$x_0: t \rightarrow \diag(t,t,1,1),$$
$$x_1: t \rightarrow \diag(1/t,1,1,t),$$
$$x_2: t \rightarrow \diag(1,1/t,t,1).$$
Then~$x^2_0 x_1 x_2$ is the cocharacter~$t \mapsto
\diag(t,t,t,t)$ and
$$\cO[X_*(T)]=\cO[x_0,x_1,x_2,(x^2_0 x_1 x_2)^{-1}] = \cO[x_0,x_1,x_2,(x_0 x_1 x_2)^{-1}].$$
The action of~$W$ preserves~$(x_0,x_0
x_1, x_0 x_2, x_0 x_1 x_2)$ considered as an  unordered
quadruple. Recalling that we are ignoring powers of~$q^{1/2}$, under
the identification of~$k[X_{*}(T)]^W$ with the spherical Hecke
algebra we have
\begin{align*}
  Q_q(X)&=X^4-T_{q,1}X^3+(T_{q,2}+2T_{q,0})X^2-T_{q,0}T_{q,1}X+T_{q,0}^2\\
        &=(X - x_0)(X - x_0 x_1)(X - x_0 x_2)(X - x_0 x_1 x_2).
\end{align*}Then $k[X_{*}(T)]^W=k[e_1,e_2,(e_3/e_1)^{\pm 1}]$ %
where
\[\sum  e_iX^i=(X - x_0)(X - x_0 x_1)(X - x_0 x_2)(X -
  x_0 x_1 x_2) .\]

Now suppose that $\alphabar_{q,1},\alphabar_{q,2}\in k^{\times}$ are
such that
$\alphabar_{q,1},\alphabar_{q,2},\alphabar_{q,2}^{-1},\alphabar_{q,1}^{-1}$
are pairwise distinct, and set $\gamma_0:=\alphabar_{v,1}$,
$\gamma_1:=\alphabar_{v,2}\alphabar_{v,1}^{-1}$,
$\gamma_2:=(\alphabar_{v,1}\alphabar_{v,2})^{-1}$. %

We let~$\mathfrak{n} \subset \cO[X_*(T)]^W$ be the maximal ideal
generated by~$\varpi$ and the $e_i-e_i(\gamma_0 ,\gamma_1 ,\gamma_2 )$, and for
each~$w\in W$ we let $\m_w$ be the maximal ideal of~$\cO[X_*(T)]$
generated by~$\varpi$ and the $w\cdot x_i-\gamma_i$. By our assumption on $\alphabar_{q,1},\alphabar_{q,2}$, the 8 %
ideals~$\m_w$ are pairwise
distinct, and exactly as in the proof of~\cite[Lem.\ 2.4.36]{BCGP},
we see that $\cO[X_{*}(T)]_{\mathfrak{n}}$ is a semi-local ring whose maximal
ideals are the~$\m_w$. In particular if~$M$ is an
$\cO[X_{*}(T)]$-module,  we can
write \[M_{\mathfrak{n}}=\oplus_{w\in W}M_{\m_w}.\]

\begin{lem}%
  \label{lem: basic properties of group algebra}Suppose that~$R$ is a commutative
  ring, and that~$G$ is a finite group. If~$H$ is a subgroup of~$G$ then we write
  $e_H:=\sum_{h\in H}h\in R[G]$, and let~$[G]_H=\sum_{g\in G/H}g\in
  R[G]$ \emph{(}the coset representatives being chosen arbitrarily\emph{)}. Then
  if~$M$ is a finite projective left
  $R[G]$-module, we have
  \begin{enumerate}
  \item \label{item: H fixed idempotent} $M^H=e_HM$.
      \item\label{item: G fixed from H Hecke operator} $M^G=[G]_HM^H$.
  \item \label{item: H invariants and base change}If $S$ is any $R$-algebra, then
    $(M\otimes_RS)^H=M^H\otimes_RS$.
  \item \label{item: H coinvariants and duality}The natural maps
    $M^H\to M$ and $M\to M_H$ induce isomorphisms of
    $R$-modules \[\Hom_R(M_H,R)\isoto\Hom_R(M,R)^H,\] \[\Hom_R(M,R)_H\isoto\Hom_R(M^H,R).\]
  \end{enumerate}
\end{lem}
\begin{proof}
Writing~$M$ as a direct summand of a free $R[G]$-module, we reduce to
the case $M=R[G]$. Then the first part is an easy calculation, while
the second part follows from the first and the relation $[G]_He_H=e_G$, which
is immediate from the definitions. The third part is immediate from
the first.

Turning to the final part,  it is easy to see that the isomorphism
$\Hom_R(M_H,R)\isoto\Hom_R(M,R)^H$ holds for any finite left
$R[G]$-module, projective or otherwise. It remains to show that if~$M$ is projective then
$\Hom_R(M,R)_H\isoto\Hom_R(M^H,R)$. We may again assume
that~$M=R[G]$, and since~$R[G]$ is a free~$R[H]$-module, we can
furthermore assume that~$H=G$.  We can identify
$R[G]$ with~$\Hom_R(R[G],R)$ by sending~$[1]$ to the 
map~$\phi$ such that~$\phi\bigl(\sum_gx_gg\bigr) = x_1$.
Combining this with the usual identification of
$R[G]_G$ with~$R$ via the trace map, we see that
$\Hom_R(M,R)_G$ is a free $R$-module of rank one generated by the
image of~$\phi$. 
By part~\eqref{item: H fixed
  idempotent}, we have $R[G]^G=R\cdot e_G$. Since~$\phi(e_G)=1$, we are
done.%
\end{proof}

\begin{lem}\label{lem: projections at q level are isomorphisms}Suppose that $\alphabar_{q,1},\alphabar_{q,2}\in k^{\times}$ are
such that
$\alphabar_{q,1},\alphabar_{q,2},\alphabar_{q,2}^{-1},\alphabar_{q,1}^{-1}$
are pairwise distinct, and define ideals~$\mathfrak{m}_w$,
$\mathfrak{n}$ as above. 

Let~$R$ be an object of~$\CNL_{\cO}$, and let~$N$ be an $R$-module with a smooth
action of $\GSp_4(\Q_q)$, with the property that if
$K_2\trianglelefteq K_1  $ are compact open subgroups of $\GSp_4
(\Q_q)$,  then $N^{K_2 }$ is a finite projective $R[K_1 /K_2
]$-module.

Then for each $w\in W$ the
projection \numequation\label{eqn: abstract projection Iwahori Hecke}\pr_w:(N^{\GSp_4(\Z_q)})_{\mathfrak{n}}\to
  (N^{\Iw(q)})_{\m_w} \end{equation}is an isomorphism.
\end{lem}
\begin{proof}Take~$K_1=\GSp_4 (\Z_q)$, $K_2=\ker\bigl(\GSp_4
  (\Z_q)\to\GSp_4 (\F_{q})\bigr)$. Set $G=K_1/K_2=
  \GSp_4(\F_q)$, and let $H=\Iw(q)/K_2=B(\F_q)<G$. Write~$M=N^{K_2}$,
  so that by assumption~$M$ is a finite projective $R[G]$-module, and
  we have $M^G=N^{\GSp_4(\Z_q)}$ and $M^H=N^{\Iw(q)}$.

  By Nakayama's lemma and  Lemma~\ref{lem: basic properties of group algebra}\eqref{item: H invariants and base
    change}, we
  can and do assume from now on that~$R=k$. 
  Let~$[\GSp_4(\Z_q)]\in\cH$ be the Hecke operator which is the indicator function
  on~$\GSp_4(\Z_q)$.
By Lemma~\ref{lem: basic
    properties of group algebra}\eqref{item: G fixed from H Hecke
    operator} we have 
    $$N^{\GSp_4(\Z_q)}=[\GSp_4(\F_q)]_{B(\F_q)} (N^{\Iw(q)}) =[\GSp_4(\Z_q)] (N^{\Iw(q)}),$$
    where the second equality follows from the very definition of the Hecke operators.
    We claim that (under our assumption that $R=k$)
  the Hecke operator $[\GSp_4(\Z_q)]$ is an inverse to~$\pr_w$.

To see this, note firstly that by the Bruhat decomposition we have
\[[\GSp_4(\Z_q)]=\sum_{w\in W}w\] (where we are using the
identification $\cH\otimes_{\cO}k=k[X_*(T)\rtimes W]$). We have the
decomposition \[ (N^{\Iw(q)})_{\mathfrak{n}}=\oplus_{w\in W}
  (N^{\Iw(q)})_{\m_w},\] so we can write any~$x\in
(N^{\Iw(q)})_{\mathfrak{n}}$  as $x=\sum_{w\in W}x_w$ for unique
elements~$x_w\in (N^{\Iw(q)})_{\m_w}$, and in particular if $x\in
(N^{\GSp_4(\Z_q)})_{\mathfrak{n}}$ then $\pr_w(x)=x_w$. By the definition of
the~$\m_w$, we see that the action of~$W$ on
$(N^{\Iw(q)})_{\mathfrak{n}}$ is via $w_1x_{w_2}=x_{w_1 w_2 }$ for all
$w_1 ,w_2 \in W$. 

We can therefore compute that if $x\in
(N^{\GSp_4(\Z_q)})_{\mathfrak{n}}$ then
$[\GSp_4(\Z_q)]\circ\pr_w(x)=[\GSp_4(\Z_q)]x_w=\sum_{w'\in
  W}w'x_w=\sum_{w'\in W}x_{w'w}=x$, while if $y\in
(N^{\Iw(q)})_{\m_w}$ then
$\pr_w\circ[\GSp_4(\Z_q)](y)=\pr_w\bigl(\sum_{w'\in W}w'y\bigr)=y$, as required.
  \end{proof}

  \begin{rem}
    \label{rem: compare to Whitmore}The second half of the proof of
    Lemma~\ref{lem: projections at q level are isomorphisms} is
    essentially identical to the proof of~\cite[Prop.\ 5.10]{Whitmore}
    in the special case that~$G=\GSp_4$ and~$P=B$. Indeed note that
    while it is assumed there that $p\nmid\# W$, all that is used in the
    proof is that $p\nmid \#W_L$, where~$L$ is the Levi factor
    of~$P$; and for~$P=B$ the group~$W_L$ is trivial.
  \end{rem}

We now establish Hypothesis~\ref{hypothesis: data for the
  SW NT patching argument}~\eqref{item:
  SW killing DeltaN}.
\begin{lem}
  \label{lem: TW system going from Q1 to Q0}The natural pullback map $M_{K^p}\to
  M_{K_1^p(Q_N)}$ %
induces an isomorphism
  of~$R_{\cS_{Q_N}}$-modules 
  $(M_N)_{\Delta_{Q_N}}\to M$.
\end{lem}
\begin{proof} First, we have that $$(M_N)_{\Delta_{Q_N}} = \Hom_{\Lambda_{\GSp_4,\Q}}((M_{K^p_0(Q_N)}\otimes_{\cO\llb T(\ZZ_p)\rrb} \Lambda_{\GSp_4,\Q})^{|\cdot|^2},\Lambda_{\GSp_4,\Q})_{\m^{Q_N},\m_{\vzero},\m_{Q_N}}$$ by
  part~\eqref{item: H coinvariants and duality} of Lemma~\ref{lem: basic properties of group algebra} and Lemma \ref{lem: MN is free over DeltaN}. 
 We claim that the map \[(M_{K^{p}} \otimes_{\cO\llb T(\ZZ_p)\rrb} \Lambda_{\GSp_4,\Q})^{|\cdot|^2}_{\m, \m_{\vzero}} \rightarrow (M_{K^p_0(Q_N)} \otimes_{\cO\llb T(\ZZ_p)\rrb} \Lambda_{\GSp_4,\Q})^{|\cdot|^2}_{\m^{Q_N},\m_{\vzero},\m_{Q_N}}\] is an isomorphism. It suffices to see (imposing that the central character acts by  $|\cdot|^2$)  that the map \[(M_{K^{p}} \otimes_{\cO\llb T(\ZZ_p)\rrb} \Lambda_{\GSp_4,\Q})^{|\cdot|^2, Q_N \cup \{r\}}_{\m, \m_{\vzero}}\to
  (M_{K^p_0(Q_N)} \otimes_{\cO\llb T(\ZZ_p)\rrb} \Lambda_{\GSp_4,\Q})_{\m^{Q_N},\m_{\vzero},\m_{Q_N}}^{|\cdot|^2, Q_N \cup \{r\}}\] is an isomorphism. Here the subscript~$(|\cdot|^2,Q_N \cup \{r\})$ indicates that the central character
acts by $|\cdot|^2$ up to primes in~$Q_N \cup \{r\}$ (as in
Lemma~\ref{lem:freeness-over-the-neat-prime}).

Let~$K(Q_N)=\prod_lK(Q_N)_l\subset\GSp_4(\A^{\infty,p})$ be defined by
  \[
    \begin{cases}
      K^p(Q_N)_l=K^p_l & \text{if } l\not\in Q_N,\\
 K(Q_N)_l=\ker\bigl(\GSp_4(\Z_l)\to\GSp_4(\F_l)\bigr) & \text{if } l\in Q_N
    \end{cases}
  \]
By Lemma~\ref{lem:freeness-over-the-neat-prime}(\ref{neata}), $(M_{K(Q_N)} \otimes_{\cO\llb T(\ZZ_p)\rrb} \Lambda_{\GSp_4,\Q})^{|\cdot|^2,Q_N \cup \{r\}}_{\m^{Q_N},\m_{\vzero}}$  is a finite projective
  $\Lambda_{\GSp_4,\Q}[\prod_{q\in Q_N}\GSp_4(\F_q)]$-module. We conclude by Lemma~\ref{lem: projections at q level are isomorphisms}.
  \end{proof}
We summarize our results so far in the following proposition; we
remind the reader that at this point  we do not know that~$M_{\q}$ is
nonzero. Indeed, even if we knew that~$\rho$ was modular, it could be
that our choice of subgroups~$K_l$ for~$l\in S$ forces~$M_{\q}$ to be
zero. In Proposition~\ref{prop: the modularity result for p equals 2
  or 3} we will establish sufficient conditions under which~$M_{\q}$
is free of rank~$1$.
\begin{prop}
  \label{prop: abstract freeness result without yet having modularity}
  Assume that $\rho$ satisfies Hypothesis~\ref{hyp: assumptions on rho
    for mult one}, and let~$\cS$ be the deformation
  problem~\eqref{defn: global deformation problem for char 0 mult one}
  with~$S$ as in Definition~\ref{defn: choice of S for mult one char
    0}.  Write~$\mathfrak{q}$ for the kernel of the homomorphism
  $R_{\cS}\to\cO$ corresponding to~$\rho$. Define~$M$ as
  in~\eqref{eqn: definition of M}. Then $M_{\q}$ is a finite free
  $(R_{\cS})_{\q}$-module.
\end{prop}
\begin{proof} By Lemmas~\ref{lem:freeness-over-the-neat-prime} and~\ref{lem: TW system going from Q1 to Q0}, our
  construction of the modules~$M_N$ above gives the data of
  Hypothesis~\ref{hypothesis: data for the SW NT patching argument}.
  The result is then immediate from Proposition~\ref{prop: abstract
    mult one char 0}.
\end{proof}

   \subsection{Multiplicity one}\label{subsec: deducing R=T GSp4
     from base change}Before establishing our main multiplicity one
   results we begin with some background material and preliminary lemmas. We
   refer to Section~\ref{sec:Arthur} for the relationship between cuspidal automorphic
   representations~$\pi$ of $\GSp_4 /\Q$ which are of general type, and cuspidal
   automorphic representations~$\Pi$ of~$\GL_4 /\Q$ of symplectic type.
   
      Assume from now on that~$F=\Q$, and that furthermore~$\pi$ has central
   character $\omega_{\pi}=|\cdot|^2$. We now recall some consequences of the theory of newforms due to
   Roberts and Schmidt~\cite{MR2344630}. (This theory assumes that we
   are working with representations of trivial central character, but
   this is harmless, as we can reduce to this case by twisting~$\pi$
   by the everywhere unramified character~$|\cdot|$.) Recall that for
   each prime~$l$ and each ~$n\ge 0$, the paramodular group of
   level~$n$ is 
   \[\Par(l^n):=\{g\in\GSp_{4}(\Q_l)\mid
     \nu(g)\in\Z_l^{\times}\}\cap
     \begin{pmatrix}
       \Z_l& \Z_l& \Z_l&  l^{-n}\Z_l\\
       l^n\Z_l& \Z_l& \Z_l&  \Z_l\\
      l^n \Z_l& \Z_l& \Z_l&  \Z_l\\
             l^n\Z_l& l^n\Z_l& l^n\Z_l&  \Z_l\\
     \end{pmatrix}
   \]We say that~$\pi_l$ is \emph{paramodular} if
   $(\pi_l)^{\Par(l^n)}\ne 0$ for some~$n$. The minimal such~$n$ is
   the \emph{paramodular level} $N_{\pi_l}$ of~$\pi_l$.

   The following result summarizes the facts that we need about
   paramodular vectors in cuspidal automorphic representations~$\pi$
   of~$\GSp_4/\Q$. %
   \begin{prop}
     \label{prop:multiplicity-one-paramodular-general-type}Suppose
     that~$\pi$ is of general type.
     \begin{enumerate}
     \item For each prime~$l$, there is a unique paramodular
   representation in the $L$-packet containing~$\pi_l$, namely the
   unique generic representation.
     \item If~$\pi_l$ is generic, then 
   $(\pi_l)^{\Par(l^{N_{\pi_l}})}$ is one-dimensional.
 \item The paramodular level~$N_{\pi_l}$ coincides
   with the conductor of the corresponding $L$-parameter
   $\recGT(\pi_l)$.
 \item If~$\pi$ is regular algebraic and ~$\rho_{\pi,p}$ is
   irreducible, then ~$N_{\pi_l}$ coincides with the conductor of
   $\rho_{\pi,p}|_{G_{\Q_l}}$.
     \end{enumerate}
   \end{prop}
   \begin{proof}The first part is~\cite[Thm.\ 1.1]{SchmidtParamodularPacket}, the
   second is ~\cite[Thm.\ 7.5.1]{MR2344630}, and the third is ~\cite[Thm.\
   2.3.5]{johnsonleung2022stable}. The last claim is~\cite[Thm.\
   2.7.1(2)]{BCGP} (see e.g.\ \cite{MR3477046} for the various
   equivalent definitions of the conductor of a Galois
   representation).     
   \end{proof}

 For the following lemma we return 
 to the setting of Section~\ref{subsec: notation and defns
  Galois}.
   \begin{lem}
     \label{lem: same component pure means same conductor}Suppose
     that~$v\nmid p$ and that~$x_1,x_2$ are two closed points
     of~$\Spec R_v^{\square}[1/p]$ which lie on a common irreducible
     component, and are such that the corresponding lifts ~$\rho_{x_1},\rho_{x_2}$
      of~$\rhobar|_{G_{F^+_v}}$ are both pure. Then we have the
     equality of conductors~$a(\rho_{x_1})=a(\rho_{x_2})$.
   \end{lem}
   \begin{proof}
     Let~$x'_i$ be the image of~$x_i$ in the spectrum of 
     the~$\GL_4$ lifting ring for~$\rhobar|_{G_{F^+_v}}$. These points lie on a common irreducible
     component by the assumption on~$x_1,x_2$, and the purity
     of the~$\rho_{x_i}$ ensures that this is the unique irreducible component
     that either lies on (e.g.\ by \cite[Cor.\ 3.3.4]{Bellovin} and
     the definition of purity),  The result follows  immediately from
     ~\cite[Lem.\ 1.3.4(2)]{BLGGT} (a lemma of Choi).
   \end{proof}

   We now establish some instances of solvable descent for~$\GSp_4$.
   \begin{prop}\leavevmode
     \label{prop:descent to symplectic plus base change results}
     \begin{enumerate}
     \item\label{item:1}If~$\Pi$ is a $|\cdot|^2$-self dual regular algebraic cuspidal
       automorphic representation of~$\GL_4/\Q$, then~$(\Pi,|\cdot|^2)$ is
       of symplectic type.
    
     \item\label{item:2}Let~$F/\Q$ be a solvable Galois extension
       with~$F$ CM. Suppose that 
       \[\rho:G_{\Q}\to\GSp_4(\Qpbar)\]
       has
       multiplier~$\varepsilon^{-1}$, that~$\rho|_{G_{F}}$ is
       irreducible, and that there is a RACSDC automorphic
       representation~$\pi_{F}$ of~$\GL_4/F$ such
       that~$\rho|_{G_{F}}\cong \rho_{\pi_{F},p}\otimes\varepsilon$.
 Then there is a regular algebraic cuspidal automorphic
       representation~$\pi$ of~$\GSp_4/\Q$ with central
       character~$|\cdot|^{2}$, such that $\rho\cong \rho_{\pi,p}$.
     \item\label{item:3}Let~$F^+/\Q$ be a solvable Galois extension
       with~$F^+$ totally real. Suppose that
       $\rho:G_{\Q}\to\GSp_4(\Qpbar)$ has
       multiplier~$\varepsilon^{-1}$, that $\rho|_{G_{F^+}}$ is
       irreducible, and that there is a regular algebraic cuspidal
       automorphic representation~$\pi_{F^+}$ of~$\GSp_4/F^{+}$,
       with central character~$|\cdot|^{2}$, and such
       that~$\rho|_{G_{F^+}}\cong \rho_{\pi_{F^+},p}$.

       Then there is a regular algebraic cuspidal automorphic
       representation~$\pi$ of~$\GSp_4/\Q$ with central
       character~$|\cdot|^{2}$, such that $\rho\cong \rho_{\pi,p}$.
            \end{enumerate}
   \end{prop}
   \begin{proof}We begin with part~\eqref{item:1}. If the pair
     $(\Pi,|\cdot|^2)$ is not of symplectic type, then it is of
     orthogonal type, so it descends to an automorphic
     representation~$\pi^{\alpha}$ of some~$\GSpin_4^{\alpha}/\Q$,
     with central character~$|\cdot|^2$. The central character
     of~$\pi^{\alpha}_{\infty}$ can be read off from its
     $L$-parameter.
     Under our assumption
     that~$\Pi$ is regular algebraic (i.e.\ $C$-algebraic), 
     it follows from~\cite[Lem.\
     3.2(3), 3.4]{MR3343873} that the central character must be odd,
     which means in particular that it
     cannot equal ~$|\cdot|^2$.
     This
     contradiction implies that $(\Pi,|\cdot|^2)$ is of symplectic type,
     as claimed.

By part~\eqref{item:1}, in each of parts~\eqref{item:2}
and~\eqref{item:3} it suffices to show that there is a
$|\cdot|^2$-self dual regular algebraic cuspidal
       automorphic representation ~$\Pi$ of~$\GL_4/\Q$ with
       $\rho\cong \rho_{\Pi,p}$. Indeed by ~\eqref{item:1} such
       a~$\Pi$ is of symplectic type, and we can take~$\pi$ to be a
       descent of~$\Pi$. Then part~\eqref{item:2} is a standard
       consequence of solvable descent for~$\GL_4$, and in particular
       is a special case of~\cite[Lem.\ 2.2.2]{BLGGT} (bearing in mind
       ~\cite[Lem.\ 2.2.1]{BLGGT}, which takes care of the twist
       by~$\varepsilon$).  Finally for part~\eqref{item:3}, since
       $\rho_{\pi_{F^+},p}$ is irreducible, we see that~$\pi_{F^+}$ is
       of general type. Its transfer~$\Pi_{{F^+}}$ is $|\cdot|^{2}$-self
       dual, and the result follows from another application of ~\cite[Lem.\ 2.2.2]{BLGGT}.     
   \end{proof}
We now return to our running hypotheses, so that~$\rho:G_{\Q}\to\GSp_4(\cO)$ is a continuous representation satisfying
 Hypothesis~\ref{hyp: assumptions on rho
    for mult one}, $\vzero$ is a fixed neat prime for~$\rhobar$, and ~$\cS$ is the deformation
  problem~\eqref{defn: global deformation problem for char 0 mult one}
  with~$S$ as in Definition~\ref{defn: choice of S for mult one char
    0}.  Write~$\mathfrak{q}$ for the kernel of the homomorphism
  $R_{\cS}\to\cO$ corresponding to~$\rho$.
\begin{lem}\label{lem: lower bound for dimension of generic fibre
       of GSp4 deformation ring}%
     \label{lem: lower bound irreducible components of deformation rings}Every
  irreducible component of~$\Spec  R_{\cS}$ containing $\q$ has dimension 
  at least~$\dim \Lambda_{\GSp_4,\Q}$.
\end{lem}
\begin{proof}
We claim there is a presentation
\[R_p^\triangle\llb x_1,\dots,x_r\rrb /(f_1,\dots,f_{r+14})
  \isoto R_{\cS}\] for some~$r\ge 0$.
 Indeed since ~$\{p\}\subsetneq S$, the map
~\cite[(4.2.1)]{MR3152673} is injective, and the existence of such a
presentation is a consequence of~\cite[Prop.\ 4.2.5, Rem.\ 4.2.6]{MR3152673}.  Now by Hypothesis \ref{hyp: assumptions on rho for mult one} \eqref{item: weight 2 pure p dist} and Lemma~\ref{lem: criterion for H2 vanishing} \eqref{item: def ring regular}, $\rho|_{G_{\Q_p}}$ lies on a unique irreducible component of $\Spec R_p^\triangle$ which has dimension $\dim\Lambda_{\GSp_4,\Q}+14$ and the result follows.

\end{proof}

We now put ourselves in the setting of Section~\ref{subsec: defn of GSp4 TW
  system for mult one}, and if $l \in S\setminus \{\vzero\}$, then we take \numequation\label{eqn:
  paramodular level of conductor}K_l
= \Par(l^{a(\rho|_{G_{\Q_l}})}).\end{equation}
  Define~$M$ as
  in~\eqref{eqn: definition of M}. Let~$\q'$ be the kernel of a
  homomorphism $R_{\cS}\to \cO$ corresponding to a lift
  $\rho':G_{\Q}\to\GSp_4(\cO)$  of~$\rhobar$, and assume that~$\q'$
  and~$\q$ lie on a common irreducible component of~$\Spec R_{\cS}[1/p]$ (in particular, we could
  take~$\q'=\q$, but we will also consider other choices below). Then
  by the definition of~$M$, we have

 \numequation\label{eqn:dim-fibre-of-M}\dim_{E}(M/\q'M)[1/p]=\dim_{E} \left( (M_{K^p})_{\m^{p},\m_{\vzero},|\cdot|^2}\otimes_{\cO}E \right)[\q'].\end{equation}%

\begin{defn}%
  \label{defn:classical weight} If~$R\in\Z_{\ge 1}$, we say that a $\Qpbar$-point of
  $\Spec\Lambda_{\GSp_4,\Q}$ is of \emph{$R$-regular classical weight} if the
  corresponding characters $\theta_1 ,\theta_2:I_{\Qp}\to\Qpbartimes$
  are algebraic with respective Hodge--Tate weights $h_1 ,h_2 $ satisfying $h_2 -h_1 ,
  1-2h_2 \ge R$.
\end{defn}%

\begin{lem}
  \label{lem:mult-one-at-regular-weight}For all sufficiently large~$R$, if the point
  of~$\Spec\Lambda_{\GSp_4,\Q}$ determined by~$\q'$ is of $R$-regular classical
  weight, %
  and~$M_{\q'}$ is nonzero, then $\dim_{E}(M/\q'M)[1/p]=1$.
\end{lem}
\begin{proof}
We claim that $\dim_{E}(M/\q'M)[1/p]$ is equal
to \numequation\label{eq:sum-giving-rank-of-M-at-q'}\sum_{\pi}\dim_{E}(\pi^{\infty})^{K^p\Iw(p),U_{p,1}=\alpha'_p,U_{p,2}=\alpha'_p\beta'_p,\beta_{\vzero,1}=\alpha'_{\vzero,1},\beta_{\vzero,2}=\alpha'_{\vzero,1}\alpha'_{\vzero,2}}\end{equation}where the sum is over the cuspidal automorphic representations~$\pi$ of
weight determined by~$\q'$ with central character~$|\cdot|^2$,
with~$\pi_{\infty}$ respectively holomorphic if~$l(w)=0$ or~$3$, and generic
if~$l(w)=1$ or~$l(w)=2$,
and which satisfy $\rho_{\pi,p}\cong\rho'$; and $\alpha'_p,\beta'_p$
are the lifts of~$\alphabar_p,\betabar_p$ determined by~$\q'$, and
similarly for~$\alpha'_{\vzero,1},\alpha'_{\vzero,2}$ (where the
$\beta_{\vzero,i}$ act as in~\eqref{eqn: torus to pro v Iwahori}). %
Indeed ~\eqref{eqn:dim-fibre-of-M} and
Theorem~\ref{thm:properties-of-higher-Hida-theory}~\eqref{item-higher-Hida-gives-ordinary-cohomology}
reduce this to the corresponding assertion about the (cuspidal or otherwise)
coherent cohomology of Shimura varieties, which holds by a standard argument
using~\cite{harris-ann-arb,blasius-harris-ramak, HZIII}. %
More precisely, the results of
Harris--Zucker allow us to reduce to the case of interior cohomology, and the
argument is then identical to that of the proof of~\cite[Thm.\ 3.10.1]{BCGP}. 

Since~$\rhobar$ is absolutely irreducible, so is~$\rho'$, so any
such~$\pi$ is of general type. We need to show that there is a
unique~$\pi$ with a nonzero contribution to~\eqref{eq:sum-giving-rank-of-M-at-q'}, and that this
contribution is~$1$. By strong multiplicity one, it suffices to show
that for each prime~$l$, there is a unique~$\pi_l$ in the $L$-packet
corresponding to~$\rho'|_{G_{\Q_l}}$ which contributes, and that it contributes with
multiplicity one. For $l\ne p,\vzero$, this follows from our choice of~$K^p$,
together with
Proposition~\ref{prop:multiplicity-one-paramodular-general-type} and
Lemma~\ref{lem: same component pure means same conductor}. For~$l=p$
it follows from ordinarity and the assumption that~$\q'$ is in $R$-regular classical
weight that~$\pi_p$ is an irreducible unramified principal series
representation, and that the simultaneous eigenspaces for the $U_{p,i}$-operators
are 1-dimensional (see~\cite[Prop.\ 2.4.24,
2.4.26]{BCGP}). Finally for~$l=\vzero$, since~$\rho'$ and~$\rho$ lie on a common irreducible
component of~$\Spec R_{\cS}[1/p]$, we know that~$\rho|_{G_{\Q_{\vzero}}}$ has
unipotent ramification, so~$\pi_{\vzero}$ has~$\Iw(\vzero)$-fixed
vectors. Since~$\rhobar(G_{\Q_{\vzero}})$ is regular semi-simple, we
again conclude that the simultaneous
$\beta_{\vzero,1},\beta_{\vzero,2}$-eigenspaces are 1-dimensional, as required.
\end{proof}

We now prove our multiplicity one criterion. For the convenience of the reader,
we incorporate our running hypotheses into the statement of the result.
   \begin{prop}%
     \label{prop: the modularity result for p equals 2 or 3}Suppose
     that  $\rho:G_{\Q}\to\GSp_4(\cO)$ satisfies the following conditions.
  \begin{enumerate}
  \item  $\rho$ is unramified at all but finitely many primes.
  \item $\nu\circ\rho=\varepsilon^{-1}$.
  \item $\rho(G_{\Q(\zeta_{p^\infty})})$ is integrally enormous.
  \item $\rho$ is pure.  %
  \item $\rho|_{G_{\Qp}}$ is ordinary, semistable of weight~$2$, and $p$-distinguished.
  \end{enumerate}
  Suppose furthermore that either $p=2$, and there is a solvable CM extension~$F/\Q$
  and an ordinary RACSDC automorphic
       representation~$\pi$ of 
       $\GL_4/F$ such that: %
  \begin{enumerate}[label=(A\arabic*)]%
     \item $\rbar_{\pi,2}\cong\rhobar|_{G_{F}}$. 
     \item \label{item: nearly adequate p 2 assumption for mult 1} $\rhobar(G_F)$ is nearly adequate.
     \item $\rhobar(G_F)$ contains a regular semi-simple element.
     \item There exists an infinite place~$v$
       of~$F^+$ such that the polarized pair 
       $(\rhobar|_{G_F},1)$ is strongly
       residually odd at~$v$.
       \item \label{item:easy-way-to-avoid-centre-causing-trouble}$\rhobar(G_{\Q(i)})=\rhobar(G_{\Q})$. %
     \end{enumerate}
  Or alternatively, suppose that~$p>2$, and there exists an ordinary cuspidal
  automorphic representation $\pi$ of $\GSp_4/\Q$ with central character
  $|\cdot|^2$ such that: %
  \begin{enumerate}[label=(B\arabic*)]%
  \item\label{item: this needs a label} $\rhobar_{\pi,p}\simeq\rhobar$.
  \item\label{item: reasonable p not 2 assumption for mult 1}  $\rhobar$ is $\GSp_4$-reasonable, in the sense
          of~\cite[Defn.\
          3.19]{Whitmore}. 
         \item  \label{tidylabel} $\rhobar$ is tidy, in the sense of~\cite[Defn.\
           7.5.11]{BCGP}. 
         \item \label{extra-regular-semi-simple-for-central-char} 
    $\rhobar(G_{\Q})\setminus \Sp_4 (\Fp)$ contains a regular semi-simple element.
    \item\label{item: same component}%
      There is a compatible choice of $p$-stabilizations
      of~ $\pi_p$ and~ $\rho|_{G_{\Qp}}$ such that $\rho_{\pi,p}|_{G_{\Q_p}}$ lies on a unique component of $\Spec R_p^\triangle$ and $\rho|_{G_{\Q_p}}$ lies on the same component.
  \end{enumerate}Let~$\vzero$ be a neat prime for~$\rhobar$, and
  define ~$S,\cS, M$  as in Definition~\ref{defn: choice of S for mult one char
    0} and  ~\eqref{defn: global deformation problem for char 0 mult one},
 \eqref{eqn: definition of M} respectively, where we make the choice of level structure~\eqref{eqn:
  paramodular level of conductor}. Let~$\q$ be the prime
 of~$R_{\cS}$ determined by~$\rho$. Then~$M_{\q}$ is a
  free~$(R_{\cS})_{\q}$-module of rank~$1$.
\end{prop}%
\begin{proof}
  Note firstly that in either case~$\rho$  satisfies
  Hypothesis~\ref{hyp: assumptions on rho for mult one} (with~$F^+$ there equal to~$\Q$), because
  ~$\rhobar$ is absolutely irreducible by whichever of~\ref{item: nearly adequate p 2 assumption for mult 1}
 and~\ref{item: reasonable p not 2 assumption for mult 1} applies,
  and if~$p=2$ then
  $\rhobar(G_{\Q(i)})=\rhobar(G_{\Q})$ by~\ref{item:easy-way-to-avoid-centre-causing-trouble}. We showed in Section~\ref{subsec: defn of GSp4 TW
  system for mult one} that the $R_{\cS}$-module~$M$ is part of a set
of data satisfying Hypothesis~\ref{hypothesis: data for the SW NT
  patching argument}, so~$M_\q$ is a
  free~$(R_{\cS})_{\q}$-module by Proposition~\ref{prop: abstract mult one
    char 0}. %

  Recall that we have a surjection $R_{\cS}\to\T_{\cS}$, where the
  Hecke algebra~$\T_{\cS}$ is defined in
  Definition~\ref{defn:Hecke-algebra-for-NT-argument}, and acts
  faithfully on~$M$ by definition. We next show in the case~$p=2$ that the induced map
  $R_{\cS}^{\red}[1/p]\to\T_{\cS}[1/p]$ is an isomorphism.  %
  In the case~$p>2$ we prove
  the weaker %
  statement that the image of the
  corresponding map on spectra contains an irreducible component
  containing~$\q$. In either case it follows that~$M_{\q}$ is nonzero, and we
  will conclude by showing that its rank is one.

  We begin with the case~$p=2$. After possibly replacing~$F^+$ with a
  quadratic extension, we can and do assume that all places of~$F^+$
  lying over~$S$ split in~$F$, as well as any place lying under a place at which $\pi$ is ramified. After making a further solvable extension of
  totally real fields, we can furthermore assume that $F/F^+$ is everywhere unramified and for each place~$v|2$
  of~$F$, $\rhobar|_{G_{F_v}}$ is trivial and we have $[F_v:\Q_2]>7$; so by
  Remark~\ref{rem:Jack-shows-don't-need-flat-closure} restriction induces a map from the deformation problem~$\cD_2^{\triangle}$ of Section~\ref{subsec: ord def
  rings} to the deformation problems~$\cD_v^{\triangle}$ of Proposition~\ref{prop: ord local def rings Un}. Let $T$ be the set of places of~$F^+$
  which lie over places in~$S$, and let ~$R^{T,\ord}$ be the global
  deformation problem for~$\rhobar|_{G_{F^+}}$ defined
  in Definition~\ref{defn:R-and-T-for-Un-theorem}. This is by
  definition a ~$\Lambda$-algebra, where $\Lambda$ is as
  in~\eqref{eqn:Lambda-for-Un}.   Then we have a
  morphism of~$\Lambda$-algebras
  $R^{T,\ord}\to R_{\cS}$ (the ~$\Lambda$-algebra
  structure on $R_{\cS}$ comes from the natural map
  $\Lambda_{\GSp_{4},F^+}\to\Lambda_{\GSp_4,\Q}$), defined
  by %
  applying  the
  construction of Corollary~\ref{cor: another symplectic to
    unitary} to~$(\rho\otimes\varepsilon^{-1})|_{G_{F^+}}$. (Note that
  this construction indeed gives a morphism of
  deformation problems rather than just framed deformation problems,
  because conjugating a lift of~$\rhobar$ by a
  matrix~$B\in\widehat{\GSp}_{4}(R)$ corresponds to conjugating the
  corresponding lift of $r_{\rhobar|_{G_{F^+}}}$
  by~$(B,\nu\circ(B))\in\widehat{\cG}^{\circ}_4(R)$.) %
By our assumptions, we can apply
Theorem~\ref{thm: U(n) ordinary automorphy lifting including Ihara
  avoidance}, and  %
conclude in particular that %
$R^{T,\ord}$ is a
finite~$\Lambda$-algebra. We have the  commutative diagram 
     \[\xymatrix{\Lambda\ar[r]\ar[d] & R^{T,\ord} \ar[d]  \\
    \Lambda_{\GSp_4,\Q} \ar[r] & R_{\cS}  \ar[r] & \T_{\cS},}\] and since $R^{T,\ord}\to R_{\cS}$ is a finite
morphism (by a standard argument exactly 
     as in the proof of~\cite[Lem.\ 1.2.3]{BLGGT}),  we deduce that $R_{\cS}$ is a finite $\Lambda$-algebra,
and thus a finite $\Lambda_{\GSp_4,\Q}$-algebra.
Combining this finiteness with Lemma~\ref{lem:
  lower bound irreducible components of deformation rings},
  it follows %
   that every
irreducible component of ~$\Spec R_{\cS}[1/2]$ dominates an irreducible component of $\Spec
\Lambda_{\GSp_4,\Q}[1/2]$, and thus that there is a dense set of closed points of
~$\Spec R_{\cS}[1/2]$ of $R$-regular classical weight. %
Let~$\rho':
G_{\Q}\to\GSp_4(\Qtwobar)$ be the lift of~$\rhobar$ corresponding to such
a point; then Theorem~\ref{thm: U(n) ordinary
  automorphy lifting including Ihara avoidance} shows
that~$\rho'|_{G_F}$ is automorphic. 
By solvable descent
(see Proposition~\ref{prop:descent to symplectic plus base change results}~\eqref{item:2}),
we deduce   that~$\rho'$ 
corresponds to a point of~$\Spec\T_{\cS}$. Since this applies to
a dense set of points, we see that the surjection
$R_{\cS}^{\red}[1/2]\to\T_{\cS}[1/2]$ is an isomorphism, as claimed.

We now turn to the case~$p>2$. %
 We choose a solvable totally real extension~$F^+/\Q$ disjoint from $\Qbar^{\ker\rhobar}$, in which $p$ splits completely, and so that for every prime $l\not=p$ for which either $\rho|_{G_{\Q_l}}$ or $\rho_{\pi,p}|_{G_{\Q_l}}$ ramifies, if $v\mid l$ is a prime of $F^+$, then:
 \begin{itemize}
 \item $q_v\equiv 1\pmod p$, and if $p=3$, then $q_v\equiv 1\pmod 9$
 \item $\rhobar|_{G_{F_v^+}}$ is trivial.
 \item If $\rho'$ is any lift of $\rhobar|_{G_{\Q_q}}$ then $\rho'|_{G_{F_v^+}}$ is unipotently ramified.  (In particular $\rho|_{G_{F_v^+}}$ and $\rho_{\pi,p}|_{G_{F_v^+}}$ are unipotently ramified.)
 \end{itemize}
  For the third point, see Lemma~\ref{lem: uniformly potentially unramified GLn}
  which we stated for~$p=2$, but whose proof makes no use of this assumption.
  We put ourselves in the setting of Section \ref{sec: Gsp4  modularity lifting}
  with the automorphic representation~ $\pi$ there being the base change
  ~$\pi_{F^+}$ of our $\pi$ to $F^+$, and $R$ the set of places lying above the primes $l\not=p$ where either $\rho|_{G_{\Q_l}}$ or $\rho_{\pi,p}|_{G_{\Q_l}}$ ramifies.  Then Hypothesis~\ref{hypothesis: Q conditions we need for GSp4 modularity lifting} is satisfied.

We have a diagram 
 \[\xymatrix{\Lambda_{\GSp_4,F^+}\ar[r]\ar[d] & R_{\cS_{1}} \ar[d]  \\
    \Lambda_{\GSp_4,\Q} \ar[r] & R_{\cS}  \ar[r] & \T_{\cS}}\]
where $R_{\cS}$ is a finite
$R_{\cS_{1}}$-algebra (again, this follows exactly
     as in~\cite[Lem.\ 1.2.3]{BLGGT}, using Lemma~\ref{lem: symplectic
       representation valued over trace} in place of~\cite[Lem.\
     2.1.12]{cht}).
     
     Consider a minimal prime $Q$ of $R_{\cS}$ contained in the prime corresponding to $\rho$.  By assumption \ref{item: same component}, the map $R_{\cS_1}\to R_{\cS}/Q$ factors through the quotient $R_{\cS_1,\pi_{F^+}}$ defined in section \ref{sec: Gsp4  modularity lifting}, which is finite over $\Lambda_{\GSp_4,F^+}$ by Theorem \ref{thm: GSp4 p bigger than 2 automorphy lifting}.  It follows that $R_{\cS_1}/Q$ is finite over $\Lambda_{\GSp_4,\Q}$.  As it also has dimension at least that of $\Lambda_{\GSp_4,\Q}$ by Lemma \ref{lem: lower bound for dimension of generic fibre
       of GSp4 deformation ring}, it follows that $\Spec R_{\cS_1}/Q\to \Spec \Lambda_{\GSp_4,\Q}$ is surjective.  Arguing exactly in the previous case we use Theorem~\ref{thm: GSp4 p bigger than 2 automorphy  lifting} and Proposition~\ref{prop:descent to symplectic plus base  change results}~\eqref{item:3} at a dense set of closed points of
$R$-regular classical weight to show that $Q$ contains the kernel of $R_{\cS}\to\T_{\cS}$,
and consequently~$M_{\q}$ is nonzero, as claimed.

It remains to show that~$M_{\q}$ is of rank~$1$. %
Let~$\eta$ be the generic point of any irreducible
component of ~$\Spec R_{\cS}[1/p]$ containing~$\rho$. Since~$M_{\q}$
is a nonzero free~$(R_{\cS})_{\q}$-module, it suffices to show
that~$M_{\eta}$ has rank~$1$ over $(R_{\cS})_{\eta}$. Since the rank
can only increase under specialization, it suffices to show that there
is some other $\rho':G_{\Q}\to\GSp_4(\Qpbar)$ corresponding to a
closed point of the component determined by~$\eta$, given by an
ideal~$\q'$ such that $M_{\q'}$ is free of rank~$1$ over
$(R_{\cS})_{\q'}$.

To see this, note that we have seen above that $\Spec R_{\cS}[1/p] $
has a dense set of closed points of $R$-regular classical weight. Let~$\q'$
correspond to such a point (on the component~$\eta$); then we are done by %
Lemma~\ref{lem:mult-one-at-regular-weight}.
\end{proof}

\begin{thm}
  \label{thm:rho-is-modular-from-mult-one-and-classicity}Suppose
  that~$\rho$ satisfies the hypotheses of Proposition~\ref{prop: the
    modularity result for p equals 2 or 3}, and that in addition either: 
\begin{enumerate}
\item the Zariski closure of $\rho(G_{\qq})$ contains $\mathrm{Sp}_4$; or
\item the Zariski closure of $\rho(G_{\Q})$
    contains~$\SL_2 \times \SL_2 $, and $\rho$ is irreducible but becomes reducible
  on some index two subgroup~$G_{E}$.
\end{enumerate} Then~$\rho$ is modular.
\end{thm}%
\begin{proof}
This is immediate from Proposition~\ref{prop: the modularity result for p equals
  2 or 3} and~\eqref{eqn:dim-fibre-of-M} (with $\mf{q}'=\mf{q}$), Theorem~\ref{thm:properties-of-higher-Hida-theory}~\eqref{item-higher-Hida-higher-Coleman} and Theorem~\ref{thm:multiplicity-one-implies-classical}, bearing in
mind Remark~\ref{rem:mult-1-suffices-for-classicality}.%
\end{proof}

\section{A 2-adic modularity theorem for abelian surfaces}\label{sec:
  modularity abelian surfaces}

In this section, we prove a %
modularity theorem (Theorem~\ref{thm: residually A5 implies modular}) for
abelian surfaces $A/\Q$ which are ordinary at~$2$ and whose mod-$2$ representation~$\rhobar: G_{\Q} \rightarrow \GSp_4(\F_2) \simeq S_6$ has a very particular form. Specifically, we demand that the image~$\Gamma$ of~$G_{\Q}$
contains a copy of~$A_5 \subset \GSp_4(\F_2)$ with index at most two which acts
absolutely irreducibly on~$\F^4_2$, and additionally require that  the image of complex conjugation is non-trivial and lands in~$A_5 \subseteq \Gamma$. 
After passing to the totally real field~$F^{+}$ (at most quadratic) such that~$\rhobar|_{G_{F^+}}$ has image~$A_5$,
we may identify~$\rhobar$ with the symmetric cube of a~$2$-dimensional representation with image~$A_5$,
and this allows us to deduce  that~$\rhobar$ is residually
modular (in regular weight) using known cases of the Artin conjecture for totally
real fields together with symmetric cube functoriality.
 In light of our previous modularity lifting theorems (in particular Theorem~\ref{thm:rho-is-modular-from-mult-one-and-classicity}), the  remaining work required to show that~$A$
 is modular is to show that the representation~$\rhobar$ is nearly adequate
in the sense of Definition~\ref{defn: nearly adequate} and strongly
residually odd in the sense of Definition~\ref{defn: strongly residually odd}. Using our results
from~\S\ref{subsec: oddness in Un} and~\S\ref{subsec: nearly adequate},
this reduces to some facts concerning the modular representation
theory of~$A_5$ in characteristic~$2$.

In~\S\ref{standardfactstorsion}, we recall some standard facts about $2$-torsion of abelian surfaces and fix once and for all a choice of isomorphism~$S_6 \rightarrow \GSp_4(\F_2)$.
In~\S\ref{A5modularmod2}, we carry out the necessary group-theoretic arguments concerning 
the mod-$2$ representation theory of~$A_5$. 
 Finally, in~\S\ref{subsec: modularity theorem},  we prove that the
 representations~$\rhobar$ we are considering are residually modular (although
 not \emph{a priori} in singular weight), and then use Theorem~\ref{thm:rho-is-modular-from-mult-one-and-classicity} to prove
the desired modularity theorem.

\subsection{The $2$-torsion of an abelian surface} \label{standardfactstorsion}
  We begin by recalling some standard facts  concerning the relationship between
Weierstrass points on a genus two curve and the~$2$-torsion on its Jacobian.
One source for the facts cited in this section is the introduction
to~\cite{VDG}.

Let~$A$ be the Jacobian of a genus two curve~$X$ over a field of
characteristic~$\ne 2$. %
   There is a Weil pairing on~$A[2]$ which defines a symplectic pairing~$\langle,\rangle$.
If one denotes the Weierstrass points (over the algebraic closure) by~$r_i$ for~$i = 1,\ldots,6$,
then for $i\ne j$ the element $r_i - r_j$
has order~$2$, and is thus a  non-trivial element of~$A[2]$.
The~$2$-torsion points~$r_i - r_j$ for $i<j$ are distinct, and they are precisely the nonzero elements of~$A[2]$. Moreover,
with respect to the Weil pairing, one has:
\numequation \label{wpairing}
\langle r_i - r_j, r_k - r_l \rangle =\#\{i,j\}\cap\{k,l\}\pmod 2
\end{equation}
for $i\not=j$, $k\not=l$.

 In~{\cite[5.1]{BPVY}, the following identification~$\iota: S_6 \isoto \Sp_4(\F_2) = \GSp_4(\F_2)$ is given:
let~$U:=\mathbf{F}^6_2$ with the bilinear form~$\langle x,y \rangle = \sum_{i=1}^{6} x_i y_i$,
let~$U^0 \subset U$ denote the trace zero subspace, and let~$L$ be the span
of~$(1,1,\ldots,1) \in U^0$. Let~$S_{6}$ act  in the obvious way
on~$=\mathbf{F}^6_2$. 
Then~$A[2] \simeq U^0/L$ where
the Weil pairing is the pairing inherited from~$U$.
To see that this isomorphism is compatible with the action on the Weierstrass points, it suffices to identify~$U$ with the~$\F_2$-space
generated by~$r_i$ for~$i = 1,\ldots, 6$. Certainly the~$r_i - r_j$ land in~$U^0$, and so it suffices to show that the divisor~$\sum (r_i)$
is congruent modulo~$2$ to a principal divisor. If one writes an affine model for~$X$ as~$y^2 = (x-r_1)(x - r_2) \cdots (x - r_6)$
and~$\{1,2,3,4,5,6\} = \{i,j,k\} \cup \{i',j',k'\}$ is any partition, then
$$r_i - r_{i'} + r_j - r_{j'} + r_k - r_{k'} = (x - r_i) + (x - r_j) + (x - r_k) - (y)$$
is principal.  Finally, the compatibility of the Weil pairing is a
consequence of equation~(\ref{wpairing}).

Under this identification, there are two conjugacy classes of subgroup~$S_5 \subset S_6$,
which one can denote~$S_5(a)$ and~$S_5(b)$, where~$S_5(b)$ is the subgroup which has a fixed 
point (we use here the same notation as~\cite[\S5.1]{BPVY}), that is, $S_5(b)$
is the standard copy of~$S_5$ in~$S_6$ (and below $A_5(b)$ denotes the copy
of~$A_5$ in ~$S_5(b)$).
It follows that~$X$ has a rational Weierstrass point,
if and only if~$\rhobar_{A,2}$ factors through a conjugate of~$S_5(b)$. 

The ~$S_6$-representation $U$ is  the
natural permutation representation.
Hence~$U$ as an~$A_5(b)$  representation is also the direct sum of the trivial
representation and the standard representation. The Brauer character of~$U$ satisfies
$$\quad \chi(1) = 5 + 1, \quad \chi((1,2,3,4,5)) = \chi((1,3,5,2,4)) =  0 + 1, \quad \chi((1,2,3)) = 2 + 1.$$
Consequently, if we let~$V = A[2] $ as an~$A_5$-representation, the Brauer character of~$V$ is
$$\chi(1) = 4, \quad \chi((1,2,3,4,5)) = \chi((1,3,5,2,4)) = -1, \quad \chi((1,2,3)) = 1.$$

\begin{lemma} The representation~$V\otimes\Fbar_2$
is the unique irreducible modular representation of~$A_5$ over~$\Fbar_2$ of dimension~$4$.
\end{lemma}

\begin{proof} This follows directly from the Brauer character table of~$A_5$;
see Lemma~\ref{irreps modulo 2}.
\end{proof}

We can make the identification~$\iota: S_6 \isoto \Sp_4(\F_2) = \GSp_4(\F_2)$ completely explicit:

\begin{lemma} \label{explicit}
An explicit isomorphism~$S_6 \rightarrow \GSp_4(\F_2)$ is given by:
$$(12)(34)(56) \mapsto \left( \begin{matrix} 1 & 0 & 1 & 0 \\ 0 & 1 & 0 & 1 \\ 0 & 0 & 1 & 0 \\ 0 & 0 & 0 & 1 \end{matrix} \right), (12) \mapsto \left( \begin{matrix} 1 & 0 & 0 & 1 \\ 0 & 1 & 0 & 0 \\ 0 & 0 & 1 & 0 \\ 0 & 0 & 0 & 1 \end{matrix} \right),
$$
$$
(12345) \mapsto \left( \begin{matrix} 0 & 0 & 1 & 1 \\ 1 & 1 & 0 & 0 \\ 1& 1 & 1 & 0 \\ 1 & 0 & 1 & 1 \end{matrix} \right),
$$
\end{lemma}

\begin{proof}
Let~$e_1 = r_1 - r_2$, $e_2 = r_3 - r_4$, $e_3 = r_3 - r_5$, and~$e_4 = r_1 - r_6$.
Then the~$e_i \in U^0$ span~$U^0/L \simeq A[2]$, and the corresponding Weil pairing
agrees with our usual choice of symplectic form~$J$. 
\end{proof}

\subsection{The modular representations of~$A_5$} \label{A5modularmod2}
 We now establish some easy group-theoretic lemmas and
  also prove some facts concerning mod-$2$ representations
  of~$A_5$. Everything here is elementary, but is included for
  completeness. We begin by describing the
irreducible modular representations of~$A_5$ in
characteristic~$2$. 

\begin{lemma}\label{irreps modulo 2}Let~$k$ be a subfield of $\Ftwobar$ which
contains~$\F_4$. Then the  irreducible representations
  of~$A_5$ over~$k$ are as follows; moreover, these representations
  are all absolutely irreducible, and in particular all absolutely
  irreducible representations of $A_5$ over $\Ftwobar$ are defined
  over~$\F_4$.
\begin{enumerate}
\item The trivial representation~$k$.
\item A two-dimensional  representation~$U$ obtained by choosing an identification~$A_5 \simeq \SL_2(\F_4)$
and then taking the tautological representation of~$\SL_2(\F_4)$ over~$k$.
\item The conjugate~$U^{\sigma}$ of~$U$ by~$\Gal(\F_4/\F_2)$ acting on~$k$.
\item \label{four} A four-dimensional representation~$V$ which is defined
  over~$\F_2$, which may be identified with~$U \otimes U^{\sigma}$,
  and also with $\Sym^3(U)$.
This lifts to the unique irreducible representation of~$A_5$ in characteristic zero of dimension four.
The representation~$V$ has a regular semi-simple element of order~$5$.
\end{enumerate}
Furthermore: there are exactly two blocks of~$A_5$, consisting of the trivial block and a block of defect zero consisting only of~$V$.
In particular, $V$ defines a projective module for~$k[A_5]$.
The Brauer character table of~$A_5$ is given as follows (where~$\zeta$
is a~$5$th root of unity):
\begin{center}
\begin{tabular}{*5c}
 \multicolumn{5}{c}{\emph{The Brauer character table of~$A_5$}}\\
 \toprule
$A_5$ & dim & $(1,2,3,4,5)$ & $(1,3,5,2,4)$ & $(1,2,3)$ \\
\midrule
$k$ & $1$ & $1$ & $1$ & $1$ \\
$U$ & $2$ & $\zeta+\zeta^{-1}$& $\zeta^2 + \zeta^{-2}$ & $-1$ \\
$U^{\sigma}$ & $2$ & $\zeta^2+\zeta^{-2}$& $\zeta + \zeta^{-1}$ & $-1$ \\
$V$ & $4$ & $-1$ & $-1$ & $1$ \\
\bottomrule
\end{tabular}
\end{center}
\end{lemma}

\begin{proof}
The group~$A_5$ has~$4$ conjugacy classes of order prime to~$2$,  and thus
has~$4$ distinct irreducible representations over~$k$. The trace of~$(1,2,3,4,5)$ on~$U$
is~$\zeta+\zeta^{-1}$. 
For the Brauer character table, see~\cite[Ch 4,
Example~8.5]{BrauerChar} or~\cite[Example~18.6]{MR0450380}.
From this table,
 the identifications~$V = U \otimes U^{\sigma}$ and
$V=\Sym^3(U)$ follow. (One can also deduce these identifications
from the Steinberg
tensor product theorem.) 
The facts concerning the blocks
can be read off from the decomposition matrix and Cartan matrix given
in~\cite[Ch 4,
Example~8.5]{BrauerChar}. The projectivity of~$V$ is also immediate
from~\cite[Prop.\ 46]{MR0450380}, since $4 = \dim(V)$ is the largest
power of~$2$ dividing~$|A_5| = 4 \cdot 15$.
 The fact that the order~$5$ elements act with distinct
eigenvalues on~$V$ is also apparent from the character table.
\end{proof}

\begin{lemma} \label{toquoteone}%
The representation~$V$ of~$G=A_5$ is nearly adequate \emph{(}in the sense of Definition~\ref{defn: nearly adequate}\emph{)}.
\end{lemma}
\begin{proof}By definition, we need to show that  
\begin{enumerate}
\item $V$ is weakly adequate.
\item $H^{1}(G,k)=0$.
\item $H^1(G,\ad V) = 0$.
\end{enumerate}

  The first claim follows directly from~\cite[Prop~9.1]{MR3626555} since~$A_5 = \SL_2(\F_4)$
and~$4 = 2^2 > 3$. (We also give a simple direct proof in
Lemma~\ref{lem:A5-is-weakly-adequate} below.) %
The second claim is immediate from the fact that~$A_{5}$ is perfect.
For the third claim, recall that~$V$ is projective, and thus~$\ad V
\simeq V \otimes V$ is also projective (e.g.\ using that a projective
module is a direct summand of a free module, and tensor products
commute with direct sums and preserve freeness.) 
Hence we deduce that~$H^n(G,\ad V)=0$
for~$n > 0$.
\end{proof}

\begin{remark}\label{rem: A5 is only nearly adequate} From the exact sequence~$0 \rightarrow k \rightarrow \ad \rightarrow \ad/k \rightarrow 0$, we deduce
that~$H^1(G,\ad/k) \simeq H^2(G,k) \simeq k$ since the Schur
multiplier of~$A_5$ is~$\Z/2\Z$. Thus~$V$ is not adequate in the sense
of~\cite[Defn.\ 2.20]{MR3598803}.
\end{remark}

\begin{lemma} \label{toquotetwo} Suppose that~$F^{+}$ is totally real and that~$\rhobar: G_{F^+}
  \rightarrow \GSp_4(\F_2)$ is an absolutely irreducible representation
with image~$G = A_5(b)$.  
 If~$v$ is an infinite place such that~$\rhobar(c_v)$ is non-trivial,
 then~$(\rhobar|_{G_F},1)$ is strongly residually odd at~$v$ in the
 sense of Definition~\ref{defn: strongly residually odd}. \label{A5odd}
\end{lemma}

\begin{proof}  
non-trivial involutions~$A$ in~$A_5(b) \subset \GSp_4(\F_2) \simeq S_6$ are
characterized by being squares of order~$4$ elements. This is most
obvious by thinking about conjugacy classes in~$S_6$ and noting that non-trivial
involutions in~$A_5$ have the cycle shape~$({*}{*})({*}{*})$;  this conjugacy
class is also preserved by the outer automorphism so this description does not
depend on any choice of isomorphism from~$\GSp_4(\F_2)$ to~$S_6$.
Judiciously choosing a suitable order~$4$ element~$\sigma$ of~$\GSp_4(\F_2)$, we find
 that any such~$A$ is conjugate in~$\GSp_4(\F_2)$ to
$$A \sim \sigma^2 = \left( \begin{matrix} 0 & 1 & 0 & 0 \\ 1 & 0 & 0 & 1 \\ 0 & 0 & 0 & 1 \\ 0 & 0 & 1 & 0 \end{matrix} \right)^2
= \left( \begin{matrix} 1 & 0 & 0 & 1 \\ 0 & 1 & 1 & 0 \\ 0 & 0 & 1 & 0 \\ 0 & 0 & 0 & 1 \end{matrix} \right)
= \left( \begin{matrix} I_2 & J_2 \\ 0 & I_2 \end{matrix}\right),$$
so that in the notation of Section~\ref{subsubsec:involutions-in-Sp2n}, $S_2 J_2 = J_2 \cdot J_2 = I_2$ which
is manifestly not alternating. (In the explicit isomorphism~$S_6 \simeq \GSp_4(\F_2)$
of Lemma~\ref{explicit}, we have~$\sigma = (1423)(56)$ and~$\sigma^2 = (12)(34)$.)
Hence the result follows from Lemma~\ref{checkodd}.
(See also Remark~\ref{rem: GL2 oddness is particularly simple}.)
\end{proof}

\begin{prop}
  \label{prop: A5 2 torsion tolerable}Suppose that~ $F^+$ is a
  totally real field, and that $\rhobar:G_{F^+}\to\GSp_4(\F_2)$ has
  image $A_5(b)$. Suppose that there is an infinite place $v$ of~$F^+$
  such that $\rhobar(c_v)\ne 1$. %
  Then for any imaginary CM quadratic
  extension $F/F^+$, 
   the polarized
  pair $(\rhobar|_{G_F},1)$ determined by~$\rhobar$ 
  is nearly adequate and strongly residually odd at~$v$. Furthermore,
  $\rhobar(G_F)$ contains a regular semi-simple element.
\end{prop}
\begin{proof} 
The representation of~$G = A_5(b)$ in~$\GSp_4(\Fbar_2)$
is irreducible and so coincides with the representation~$V$ of
dimension~$4$ in Lemma~\ref{irreps modulo 2} part~(\ref{four}),
and so in particular has a regular semi-simple element of
order~$5$. 
By Lemma~\ref{toquoteone}, %
$V$ is nearly adequate. Since~$G$ is perfect, we have
$\rhobar(G_F)=G$, so the polarized pair $(\rhobar|_{G_F},1)$ is indeed
nearly adequate. Finally, it is strongly residually odd at~$v$ by  Lemma~\ref{toquotetwo}.
\end{proof}

We end this section with our promised direct proof that~$V$ is weakly
adequate. Recall firstly that if~$W$ is a representation of a finite
group~$G$ over a field~$k$ of characteristic~$2$, then in addition to
the usual short exact sequence
  \numequation\label{eqn:WW-decomposition}
  0\to\wedge^2W\to W\otimes W\to
  \Sym^2W \rightarrow 0,
  \end{equation}
   there is a short exact sequence 
   \numequation\label{eqn:Sym2W-decomposition}
   0\to W(1)\to\Sym^2W\to\wedge^2W\to 0,
  \end{equation}
   where the first map
is the inclusion of the subspace spanned by the $x\otimes x$ for~$x\in
W$, and the second map is the one induced by $x\otimes y\mapsto
x\wedge y$. We can and do identify~$W(1)$ with the Frobenius twist
of~$W$.

\begin{lem}\label{lem:soc-V-otimes-V}
  The socle of~$V\otimes V$ is $k\oplus V$.
\end{lem}
\begin{proof}
Consider  ~\eqref{eqn:WW-decomposition}
and ~\eqref{eqn:Sym2W-decomposition} 
with~$G=A_5$
and~$W=V$. Then~$V=V(1)$, and since~$V$ is projective, we see that
$V\otimes V$ splits as a direct sum of~$V$ and an extension
of~$\wedge^2V$ by itself. Since~$U,U^{\sigma}$ and consequently~$V$
are all self-dual, and since~$V\otimes V\cong\Hom(V,V)$ 
contains exactly one copy of~$k$ in its socle by Schur's lemma, it
suffices to prove that the socle of~$\wedge^2V$ is~$k$.

Now considering~\eqref{eqn:WW-decomposition} and~\eqref{eqn:Sym2W-decomposition} with~$W=U$ we
see that~$U\otimes U$ admits a filtration with successive graded
pieces $k,U^{\sigma},k$. Similarly, $U^{\sigma}\otimes U^{\sigma}$ has
a filtration with graded pieces $k,U,k$. Tensoring these together, we
see that the Jordan--H\"older factors (with multiplicity) of
$\wedge^2V$ are $k,k,U,U^{\sigma}$.

Since~$\wedge^2V$ is~$\Gal(\F_4/\F_2)$-invariant,  if either~$U$ or~$U^{\sigma}$ occurs
in the socle of~$\wedge^2V$, then they both do.
Since ~$\wedge^2V$ is self dual, however, if~$U$ and~$U^{\sigma}$ appear
in the socle of~$\wedge^2V$, then they also appear in the cosocle, and 
 therefore occur as direct summands. 
 If this occurs,
 then since $H^1(A_5,k)=0$, the
representation~$\wedge^2V$ would be semi-simple, contradicting the
presence of exactly one copy of~$k$ in its socle. This contradiction
completes the proof.
\end{proof}

The following lemma gives our second proof that~$V$ is weakly adequate (the
first was in the proof of Lemma~\ref{toquoteone}).
\begin{lem}%
  \label{lem:A5-is-weakly-adequate}The representation~$V$ of~$A_5$ is
  weakly adequate.
\end{lem}
\begin{proof}Let~$M$ be the subspace of~$\End(V) = \Hom(V,V) \cong V\otimes V$ generated by the
  semi-simple elements of~$A_5$.
  Note that~$M$ is an~$A_5$-module: if~$[h]$ is semi-simple, then so is
  $g.[h] = [ghg^{-1}]$.  To prove that~$V$
  is weakly adequate, it suffices (by definition) to verify any of the equivalent
  conditions of Lemma~\ref{itsclear}; we shall verify condition~(\ref{spansirreducible}), namely,
  that~$M=\End(V)$.
If~$M \rightarrow \End(V)$ is not surjective, then~$\End(V)$ has a simple quotient~$Q$ such 
that the composite map~$M \rightarrow Q$ is zero. The map from~$\End(V)$ to any simple
quotient factors through the cosocle of~$\End(V)$, hence it suffices to show that~$M$ surjects
  onto the cosocle of~$\End(V)$. %
  Since~$V$ is self-dual, this cosocle
  is isomorphic to $k\oplus V$ by Lemma~\ref{lem:soc-V-otimes-V}. The
  corresponding map $\End(V)\to k$ is the trace map, and any
  non-trivial element of~$A_5$ of odd order has nonzero
  trace on~$V$, so~$M$ surjects onto~$k$.

  It remains to show that~$M$ meets the direct summand~$V$, which we will do by showing that each
  element of~$A_5$ of order~$3$ contributes to this summand, using the description
  of this summand coming
  from~\eqref{eqn:Sym2W-decomposition} with~$W=V$. Let~$e_1,e_2$ be the standard
  basis of~$U$, and let $e_1^{\sigma},e_2^{\sigma}$ be the
  corresponding basis of~$U^{\sigma}$, so that the $e_i\otimes
  e_j^{\sigma}$ for~$i=1,2$ give a basis~$v_{i,j}$ for~$V$. The
  two elements of order~$3$ in~$A_5$ correspond to the diagonal
  matrices~$(\omega,\omega^{-1})\in\SL_2(\F_4)$, where~$\omega^3=1$,
  and each~$v_{i,j}$ is an eigenvector for these matrices (with
  eigenvalues $1,1,\omega,\omega^{-1}$). The same is true for
  the~$v_{i,j}\otimes v_{i,j}$, and these give a basis
  modulo~$\wedge^2V$ for the direct summand~$V\subset V\otimes V$, so
  we are done.
\end{proof}
\subsection{A \texorpdfstring{$2$}{2}-adic ordinary modularity theorem}\label{subsec: modularity theorem}
In this section we will establish our  $2$-adic modularity theorem
(Theorem~\ref{thm: residually A5 implies modular}). We begin by proving the
following lemma which  establishes
residual modularity in our situation; a closely related result was also obtained by Tsuzuki
and Yamauchi, see~\cite[Thm.\ 4.7]{tsuzuki2022automorphy}.
\begin{lemma}\label{lem: A5 solvable potential modularity B[2]}%
Let~$F^+$ be a totally real field, and let
\[\rhobar: G_{F^+} \rightarrow \GSp_4(\F_2) \simeq S_6\]
be a continuous Galois representation with the following  properties:
\begin{enumerate}
\item The image of~$\rhobar$ is either~$S_5(b)$ or~$A_5(b)$.
\item The image of each complex conjugation has order~$2$ and lands
  in~$A_5(b)$.%
\end{enumerate}
Then there exists a solvable extension of totally real
fields~$E^+/F^+$, and an imaginary CM quadratic extension~$E/E^+$,
such that:
\begin{itemize}
\item $\rhobar(G_E)=\rhobar(G_{E^+})= A_5(b)$, and
\item there is an ordinary RACSDC representation~$\pi$ 
of~$\GL_4/E$ with
$\rbar_{\pi,2}\cong\rhobar|_{G_E}$. \end{itemize}
\end{lemma}

\begin{proof}%
  Let~ $E^+/F^+$ denote the extension of
  degree at most~$2$ corresponding to the kernel of the
  composite~$G_{F^+}\to S_5(b) \rightarrow \Z/2\Z$. Then
  $\rhobar(G_{E^+})=A_5$, and~$E^+$ is totally real by the assumption
  on complex conjugations. Making a further solvable base change, we
  can and do assume that for each place~$v|2$, $\rhobar|_{G_{E^{+}_v}}$
  is trivial.

Let
$$\varrhobar:G_{E^+} \rightarrow \SL_2(\F_4) \simeq A_5$$
denote the residual $2$-dimensional Galois representation associated to
this~$A_5$-extension. (There are two such representations which are
permuted by the outer automorphism; choose either.) Note that
$\rhobar=\Sym^3 \varrhobar$
by Lemma~\ref{irreps modulo
    2}~\eqref{four}.

By a theorem of Tate~\cite[Thm~4]{Serrelifting},
the composite~$\varrhobar: G_{E^{+}} \rightarrow A_5 \hookrightarrow
\PGL_2(\C)$
lifts to a representation
$$\varrho:G_{E^+}\to
\GL_2(\C)$$
with finite image (which will be some central extension of~$A_5$).
Since the image of complex conjugation under~$\varrhobar$ is non-trivial in~$\PGL_2(\C)$,
the image in~$\GL_2(\C)$ is non-scalar and hence~$\varrho$ is odd.
By the odd Artin conjecture for~$\GL_2$ (i.e.\ by the main results
of~\cite{MR3581178, MR3904451}), $\varrho$ is modular. More precisely,
$\varrho$ is the Galois representation associated to an ordinary
Hilbert modular eigenform~$f$ of parallel weight~$1$ (the ordinarity
being a consequence of local-global compatibility, and the assumption
that $\rhobar|_{G_{E^{+}_v}}$ is trivial for all~$v|2$.). In
particular~$f$ is contained in a Hida family. Specializing this Hida
family to parallel weight~$2$, and making a further solvable base
change if necessary, we obtain an ordinary cuspidal
automorphic representation~$\pi_{E^+}$ of~$\GL_2/E^{+}$ of
weight~$0$ and trivial central character, with
$\rbar_{\pi,2}\cong\varrhobar|_{G_{E^+}}$.

Let~$E/E^{+}$ be an imaginary quadratic CM extension, and let~$\pi$ be the base change of~$\pi_{E^+}$ to~$E$. Then the symmetric
cube $\Sym^3\pi$ (which exists by~\cite{MR1923967}) is an
ordinary RACSDC automorphic representation
of~$\GL_4/E$. Since $\rhobar=\Sym^3 \varrhobar$, we are done.
\end{proof}

We now prove the main result of this section.
\begin{thm}
  \label{thm: residually A5 implies modular}Suppose that~$A/\Q$ is an
  abelian surface such that %
  \begin{enumerate}
  \item\label{item:6}  $A_5(b)\subseteq \rhobar_{A,2}(G_\Q)\subseteq S_5(b)$.
  \item\label{item:7}  The image of complex conjugation has order~$2$ and lands
  in~$A_5(b)$.
  \item\label{item:8}   $A$ has good ordinary or semistable reduction at~$2$, and $\rho_{A,2}|_{G_{\Q_2}}$ is ordinary and
    $2$-distinguished. %
  \end{enumerate}
  Then~$A$ is modular. More precisely, there is a
    weight~$2$ cuspidal automorphic representation~$\pi$ for~$\GSp_4 /\Q$ 
which is ordinary at~2, and satisfies $\rho_{\pi,p}\cong\rho_{A,p}$ for all~$p$.
\end{thm}
\begin{proof}
  As recalled in Section~\ref{sec:galoisintro}, the representation $\rho_{A,2}$ unramified at all but finitely many
  primes, is pure, and
  $\nu\circ\rho_{A,2}=\varepsilon^{-1}$.  %
  By Lemma~\ref{lem: A5 solvable potential modularity B[2]}, there is
  an imaginary CM field ~$E$
  and an ordinary RACSDC automorphic
       representation~$\pi$ of 
       $\GL_4/E$ such that~$E/\Q$ is solvable, $\rhobar_{A,2}(G_{E^+})=A_5(b)$, and
       $\rbar_{\pi,2}\cong\rhobar_{A,2}|_{G_E}$. Making a
       further solvable extension, we can and do assume that
       furthermore $E/E^+$ is everywhere unramified, and all of
  the places~$v|2$ of~$E^+$ split in~$E$, as do  all places lying
  under a place at which~$\pi$ is ramified, and all places lying over
  a place at which~$A$ does not have good reduction. %

  By Theorem~\ref{thm:rho-is-modular-from-mult-one-and-classicity} it
  therefore suffices to check that:

  \begin{enumerate}[(a)]
      \item  \label{modularfrommultonea} the
Zariski closure
of~$\rho_{A,2}(G_\Q)$ contains~$\Sp_4$.
\item \label{modularfrommultoneb} $\rhobar_{A,2}(G_{\Q(i)})=\rhobar_{A,2}(G_{\Q})$.
  \item \label{modularfrommultonec}  $\rho_{A,2}(G_{\Q(\zeta_{2^\infty})})$ is integrally enormous.
     \item \label{item:nearlyadequatesecond} $\rhobar_{A,2}(G_E)$ is nearly adequate.
     \item \label{modularfrommultonee} $\rhobar_{A,2}(G_E)$ contains a regular semi-simple element.
     \item \label{modularfrommultonef} There exists an infinite place~$v$
       of~$E^+$ such that the polarized pair 
       $(\rhobar_{A,2}|_{G_E},1)$ is strongly
       residually odd at~$v$.
        \end{enumerate}
Since $A_5(b)\subseteq \rhobar_{A,2}(G_\Q)$, it follows
from~\cite[Thm.\ 2.1]{MR1748293} that~$\End(A_{\Qbar}) = \Z$. By~\cite[Thm~3]{MR1730973}
(see also~\cite[Thm~5.14]{MR1603865}), this implies that the
Zariski closure
of~$\rho_{A,2}(G_\Q)$ contains~$\Sp_4$, which verifies condition~\ref{modularfrommultonea}.
To see that~$\rhobar_{A,2}(G_{\Q(i)})=\rhobar_{A,2}(G_{\Q})$ 
(condition~\ref{modularfrommultoneb}), we note that the
assumptions~\eqref{item:6} and~\eqref{item:7} imply that
$\rhobar_{A,2}(G_{\Q})$ contains at most one normal subgroup of index~$2$, and
that such a subgroup corresponds to a real quadratic field, and in particular
not to~$\Q(i)$.
For condition~\ref{modularfrommultonec},
note that %
$A_5(b)\subseteq
\rhobar_{A,2}(G_{\Q})$, so we see that $\rhobar_{A,2}(G_{\Q})$ contains a
regular semi-simple element (of order $5$), and then by
Corollary~\ref{cor: still integrally enormous up the pro p tower} %
we deduce that $\rho_{A,2}(G_{\Q(\zeta_{2^\infty})})$ is integrally
enormous. Conditions~\ref{item:nearlyadequatesecond},~\ref{modularfrommultonee}, and~\ref{modularfrommultonef}
are immediate from
Proposition~\ref{prop: A5 2 torsion tolerable}. 
\end{proof}

\section{Local Geometry of curves with a Weierstrass point}
\label{sec:twothreeswitchsection}

The goal of this section is to complete the proofs of our main modularity theorems (Theorem~\ref{first} and~\ref{second})
 using a $2$-$3$ switch. By Theorems~\ref{thm:rho-is-modular-from-mult-one-and-classicity} and ~\ref{thm: residually A5 implies modular}, we have established the following (we omit  the full list of hypotheses):

\begin{enumerate}
\item \label{secondpoint} A $2$-adic ordinary modularity theorem in weight~$2$ under
  a hypothesis on the residual image: the image $\Gamma$ of
  $G_{\Q}$ in $\GSp_4(\F_2)$ contains a copy of $A_5$ with index at
  most two acting absolutely irreducibly, and moreover complex conjugation is non-trivial and lies in $A_5\subseteq\Gamma$.

\item \label{firstpoint} A $3$-adic ordinary modularity lifting theorem in weight~$2$.

\end{enumerate}

We combine these two results as follows. Given an abelian surface
$A/\Q$ with good ordinary reduction at $3$ (satisfying our
supplementary hypotheses), we construct a second abelian surface
$B/\Q$ with $A[3]=B[3]$, such that~$B$ also has good ordinary
reduction at $3$, and so that the result described in point~(\ref{secondpoint}) can be applied to establish the modularity of ~$B$. This implies that the $3$-adic representation associated to ~$A$ is residually modular and hence that ~$A$ is modular using point~(\ref{firstpoint}).

The construction of $B$ uses the rationality of a certain twisted moduli space
$P(\rhobar)$ of principally polarized abelian surfaces
(introduced
 in~\cite[\S 10.2]{BCGP}, see Definition~\ref{MandP} below).  The space $P(\rhobar)$ is closely related to the moduli space
$\MMM^w_2(\rhobar)$ of genus two
curves $X$ with a fixed Weierstrass point and fixed level $3$ structure
$\rhobar = \rhobar_{\Jac(X),3}$. Concretely,
the Torelli map $\MMM^w_2(\rhobar) \rightarrow P(\rhobar)$ is an isomorphism onto its (open) image.
This relationship suggests a natural approach to understanding points on $P(\rhobar)$ with
suitable local properties at $p=2$ and $p=3$, including having good ordinary
reduction when $p=3$ and good reduction when $p=2$.
Namely, we can consider genus $2$
curves over $\F_p$ and $\Q_p$ with a rational Weierstrass point
 with the corresponding local properties. 
We carry out this analysis in
\S\ref{localattwo} for $p=2$ and
\S\ref{localatthree} for $p=3$.

This would suffice to
prove \emph{some} version of our main theorem (with a more restrictive hypothesis
at $2$ but still applying to a positive proportion of genus $2$ curves). 
However,
we can push these arguments further by exploiting the fact that the Jacobian
of a genus $2$ curve can have good reduction even when the original curve does not.
Moreover, a principally polarized abelian surface need not even be a Jacobian.
We carry out such auxiliary constructions (for $p=2$) in
\S\ref{surfacesattwo}. It also follows from the results of that section
that weakening the hypothesis at~ $2$ any further would require some new ideas.
Note that many of the arguments in this section could potentially become much simpler (and stronger)
once our modularity lifting theorems are generalized to totally real fields (since passing
to finite extensions makes finding local points much easier). It seemed potentially
useful, however, to push our current methods as far as possible until such results are available.

Finally, in 
\S\ref{subsec: 2 3 switch}, we carry out the details of the $2$-$3$ switch
using the results in the previous three sections and then complete
the proof of our main theorems in~\S\ref{sec:proofsofAandB}.

\subsection{Genus~$2$ curves  locally at~$2$} \label{localattwo}

 In this section, we discuss some explicit computations with genus~$2$ curves (with or without
 rational Weierstrass points) over~$\F_2$  and also over local fields. 
 Note that if~$B/\Q_2$
 is an abelian surface with good reduction, then~$\rhobar_{B,3}(\Frob_2)
 \in \GSp_4(\F_3)$ has similitude character~$\varepsilon^{-1}(\Frob_2) = -1$,
 and thus the image of~$\rhobar_{B,3}(\Frob_2)$ in the group~$\PGSp_4(\F_3)$
 does not land in~$\PSp_4(\F_3)$.

\begin{rem}[Reminder concerning conventions]  \label{rem:reminder} If~$B$ is an abelian surface over a
field of characteristic prime to~$p$,
our convention (see Section~\ref{sec:galoisintro}) is that $T_pB$ and $B[p]$ correspond to~$\rho^{\vee}_{B,p}$  and~$\rhobar^{\vee}_{B,p}$.
Note, however, that since these representations are self-dual up to twist, the associated projective
representations are independent of this choice.  
\end{rem}
\begin{defn}%
  \label{defn: p-distinguished}Let~$B/\Qp$ be an abelian surface 
  for which the associated Galois representation on~$T_pB$ is ordinary.
  We say that $B$ is
  \emph{$p$-distinguished} if
 the corresponding representation
  $\rho_{B,p}:G_{\Qp}\to\GSp_{4}(\Qp)$ is $p$-distinguished in the
  sense of Definition~\ref{defn:ordinary-for-Galois-reps}.
    For example, if~$B$ has good ordinary reduction, then~$B$
    is~$p$-distinguished if and only     if the characteristic polynomial~$Q(x)$
  of~$\Frob_p$ on~$T_\ell B$, $\ell\not=p$, has pairwise distinct roots, or equivalently 
  if~$Q(x)$ is not a square. If~$B_0/\F_p$ is ordinary,
  we say that~$B_0$ is~$p$-distinguished if the characteristic polynomial~$Q(x)$
  of Frobenius has pairwise distinct roots.  
     If~$X/\Qp$ is a genus~$2$
  curve,
    we say that~$X$ is
  $p$-distinguished if~$\Jac(X)$ is $p$-distinguished.
\end{defn}

  We begin with some basic group--theoretic
 facts, which can easily be extracted from~\cite[pp.\ 26--27]{Atlas}:
 
 \begin{lemma} \label{atlas}
 There are~$10$ conjugacy classes of elements in~$\PGSp_4(\F_3) \setminus \PSp_4(\F_3)$.
 Their orders and the characteristic polynomials of any lift to~$\GSp_4(\F_3)$
 are given by the following table. 
Here the name of the conjugacy class \emph{(}with the first number
indicating the order of the element\emph{)} follows
 the same convention as the Atlas~\cite{Atlas}: %
 \medskip
 
 \begin{center}
 \begin{tabular}{*3c}
\toprule
 $\langle g \rangle$ & $P(x) \in \F_3[x]$ &  Size \\
\midrule
 $2C$  &  $x^4 + 2 x^2 + 1$ &  $36$\\
 $2D$  &  $x^4 + x^2 + 1$ & $540$  \\
$4C$ &  $x^4 \pm x^3 + 2 x^2 \mp  x + 1$ &$540$ \\
 $4D$ &  $x^4 + 1$ & $1620$ \\ 
$6G$ & $x^4 + 2 x^2 + 1$ &  $1440$  \\ 
 $6H$ & $x^4 + 2 x^2 + 1$ &  $1440$  \\ 
 $6I$ & $x^4 + x^2 + 1$ & $4320$ \\ 
$8A$ & $x^4 \pm x^3 +  x^2 \mp  x + 1$ &  $6480$  \\
$10A$ & $x^4 \pm x^3  \mp  x + 1$  &  $5184$ \\
 $12C$ & $x^4 \pm x^3 +2 x^2  \mp  x + 1$ &  $4320$ \\
\bottomrule
 \end{tabular}
 \end{center}

 There are, in particular, $9$ different possible characteristic polynomials of elements of $\GSp_4(\F_3)\setminus\Sp_4(\FF_3)$.
There is a natural permutation representation
\[\PGSp_4(\F_3) \rightarrow S_{40}\]
coming from the action on the~$40$ points~$(\F^4_3 - 0)/\{\pm1\}$, and there is also
a unique transitive action~$\PGSp_4(\F_3) \rightarrow
S_{27}$ whose stabilizer is the maximal subgroup~$2^4:S_5 \subset \PGSp_4(\F_3)$; see \cite[p. 26]{Atlas}. %
The conjugacy class of~$g \in \PGSp_4(\F_3) \setminus \PSp_4(\F_3)$
is determined by the conjugacy classes of its images in~$S_{40}$  and~$S_{27}$,
and even the image in~$S_{40}$ suffices except for the following classes:

 \medskip
 
   \begin{center}
 \begin{tabular}{*3c}
\toprule
 $g$ & $S_{40}$ & $S_{27}$ \\
 \midrule
  $4C$ & $(4,4,4,4,4,4,4,4,4,4)$ & $(1,2,2,2,4,4,4,4,4)$ \\
 $4D$ & $(4,4,4,4,4,4,4,4,4,4)$ & $(1,1,1,1,1,2,4,4,4,4,4)$ \\
 $6G$ & $(2,2,6,6,6,6,6,6)$ &
 $(1,1,1,2,2,2,3,3,3,3,6)$\\
 $6H$ &  $(2,2,6,6,6,6,6,6)$ &
 $(3,3,3,3,3,6,6)$\\
\bottomrule
  \end{tabular}
 \end{center}
 \end{lemma}

  \begin{lemma} \label{curveexamples}
 Consider the genus~$2$ curves~$C_i$ over~$\Q$:
 $$\begin{aligned}
 C_1: y^2 + (x^3 + 1) y = & \ x^2 + x \\
 C_2: y^2 + (x^3 + 1) y = & \ - x^5 + x^3 + 2 x^2 + x - 1 \\
 C_3: y^2 + (x^2 + x) y =  & \ x^5 + 2x^4 + x^3 - 16x^2 - 8x - 1 \\
 \end{aligned}$$
Then the~$C_i$
have good ordinary reduction at~$2$ and a rational Weierstrass point.
Moreover, the~$C_i$ are~$2$-distinguished. %
 The corresponding characteristic polynomials of~$\Frob_2$ %
 are as follows:
 $$
 \begin{aligned}
P_1: & \ x^4 + 2 x^3 + 3 x^2 + 4 x + 4, \\
  P_2: &  \ x^4 + x^2 + 4, \\
 P_3: &  \ x^4 - x^2 + 4, \\
 \end{aligned}$$
and the conjugacy classes of~$\rhobar_{C_i,3}(\Frob_2)$
in~$\PGSp_4(\F_3)$ have type~$10A$, ~$6I$,
and~$6H$ respectively.
 \end{lemma}

 \begin{proof}
 These three curves have conductors~$249$, $975$, and~$1947$ respectively
 (they are taken from the LMFDB~\cite{LMFDB}).
The characteristic polynomials of Frobenius at~$2$ can be obtained by
an explicit point count;
 they
are irreducible over~$\Q$, which proves they are~$2$-distinguished.
One can determine directly from the characteristic polynomial
 that the conjugacy class of~$\rhobar_{C_i,3}(\Frob_2)$ must be~$10A$
 for~$i=1$; $2D$ or~$6I$ for~$i=2$; and~$2C$, $6G$, or~$6H$ for~$i = 3$.
 Given a genus two curve with a rational Weierstrass point,
 one can write down the general degree~$40$ polynomial
 whose splitting field is~$\PGSp_4(\F_3)$
 (see~\cite[\S 3]{CCR}), and then compute a degree~$27$
 resolvent. From this one can determine the correct conjugacy class using Lemma~\ref{atlas}.
 (For~$C_2$, the image in~$S_{40}$ is already enough to determine that the element
 has order~$6$ and so must be~$6I$.)
\end{proof}

 \begin{rem} \label{magmacite} The computations in this section and the next are all done
using \texttt{magma}~\cite{MR1484478}, and the explicit code
with documentation can be found at~\cite{magma}).
\end{rem}

We now consider what happens as we loop over all ordinary genus two curves over~$\F_2$.

\begin{lemma} \label{exhaust}
    Let~$X_0/\F_2$ be an ordinary smooth genus two curve. \begin{enumerate}\item
    \label{exhaust implies weaker condition at two}
   The action of~$\Frob_2$
    on~$\mathrm{Jac}(X_0)[3]/\{\pm1\}$ 
    has conjugacy class of type~$6G$, $6H$, $6I$, $8A$, or~$10A$.
 \item \label{item: exhaust Weierstrass implies 6H6I10A}    If~$X_0$ has  a smooth lift ~$X/\Z_2$
    with a~$\Q_2$-rational Weierstrass point, then the action
    of~$\Frob_2$ on~$\mathrm{Jac}(X_0)[3]/\{\pm1\}$ has conjugacy class of
    type~$6H$, $6I$, or~$10A$.%
  \item \label{item: exhaust 6H6I10A implies Weierstrass 2
      distinguished lift} \label{toolong} If~$\rhobar:G_{\Q_2}\to\GSp_4(\F_3)$ is an unramified
    representation with similitude~$\varepsilonbar^{-1}$, and the
    image of~$\rhobar(\Frob_2)$ in~$\PGSp_4(\F_3)$ has conjugacy class of
    type~$6H$, $6I$, or~$10A$, then there is an ordinary smooth genus
    2 curve  ~$X/\Z_2$ with
    ~$\rhobar\cong\rhobar_{X,3}$, such that 
$X$ has a~$\Q_2$-rational Weierstrass point, and is $2$-distinguished. %
  \end{enumerate}
\end{lemma}

\begin{proof} We may write any smooth genus two curve~$X_0$ over~$\F_2$ in the form
\begin{equation} \label{weierstrasstwo}
y^2 + h(x) y = f(x),
\end{equation}
where~$\deg(f(x)) \le 6$ and~$\deg(h(x)) \le 3$. 
 We may  enumerate
all such equations. We do not concern ourselves with identifying
either isomorphism classes of curves or of their Jacobians, and so in particular
when we talk of ``curves'' below we really mean curves
with a given Weierstrass equation as in~\eqref{weierstrasstwo}. Let~$A = \Jac(X_0)$.
We find that: 
\begin{enumerate}
\item There are~$2^{11}$ possible pairs of~$h(x)$ and~$f(x)$.
\item There are~$768$ curves which are smooth of genus~$2$.
\item \label{distribution} There are~$384$ ordinary curves, of which:
\begin{enumerate}
\item $32$ have~$\rhobar_{A,3}(\Frob_2)$ in~$\PGSp_4(\F_3)$ of type~$6G$,
\item $16$ have~$\rhobar_{A,3}(\Frob_2)$ in~$\PGSp_4(\F_3)$ of type~$6H$,
\item $48$ have~$\rhobar_{A,3}(\Frob_2)$  in~$\PGSp_4(\F_3)$  of type~$6I$,
\item $96$ have~$\rhobar_{A,3}(\Frob_2)$  in~$\PGSp_4(\F_3)$  of type~$8A$,
\item $192$ have~$\rhobar_{A,3}(\Frob_2)$  in~$\PGSp_4(\F_3)$  of type~$10A$.
\end{enumerate}
\item \label{semidistribution} The are are~$384$ non-ordinary curves, of which:
\begin{enumerate}
\item $48$ have~$\rhobar_{A,3}(\Frob_2)$ in~$\PGSp_4(\F_3)$ of type $6G$,
\item $144$ have~$\rhobar_{A,3}(\Frob_2)$ in~$\PGSp_4(\F_3)$ of type~$12C$,
\item $48$ have~$\rhobar_{A,3}(\Frob_2)$  in~$\PGSp_4(\F_3)$  of type~$4D$,
\item $48$ have~$\rhobar_{A,3}(\Frob_2)$  in~$\PGSp_4(\F_3)$  of type~$8A$,
\item $96$ have~$\rhobar_{A,3}(\Frob_2)$  in~$\PGSp_4(\F_3)$  of type~$10A$.
\end{enumerate}
\item \label{classify}  If~$X_0$ is ordinary and is additionally the reduction of a smooth curve
 over~$\Z_2$ with a~$\Q_2$-rational Weierstrass
point,
then~$\rhobar_{A,3}(\Frob_2)$ in~$\PGSp_4(\F_3)$ has type~$6H$, $6I$, or~$10A$.
\end{enumerate}
We first explain how to distinguish between the various conjugacy classes
in parts~(\ref{distribution}) and~(\ref{semidistribution}), and then we explain part~(\ref{classify}).
\begin{enumerate}
\item By point counting, we can compute the characteristic polynomial~$Q(x)$ of Frobenius.
The characteristic polynomials of the conjugacy classes~$8A$ and~$10A$ are not
congruent modulo~$3$ to the characteristic polynomial of any other class
in~$\PGSp_4(\F_3) \setminus \PSp_4(\F_3)$, so in these cases we are done.
This is also enough to
determine the counts in part~(\ref{semidistribution}).
\item  The conjugacy classes~$\{2C,6G,6H\}$, $\{2D,6I\}$, and~$\{4C,12C\}$ are complete
sets of conjugacy classes in~$\PGSp_4(\F_3) \setminus \PSp_4(\F_3)$ with the same
characteristic polynomial.
We now show how to distinguish the elements of order~$2$ (respectively, $4$) from the elements of order~$6$ (respectively, $12$)
(none of the elements of order~$2$ or~$4$ of this type actually occur).
 If~$g \in \PGSp_4(\F_3)$ is any element, and we choose a lift in~$\GSp_4(\F_3)$,
then~$g^2$ is independent of the choice of lift.
If~$g$ is of type~$2C$ or~$2D$, the square of any lift is scalar and
given by~$I$ and~$-I$ respectively,  and so~$g^4$ will be trivial.
If we start with~$g \in \PGSp_4(\F_3)$ of order~$6$, however, then for any lift, the element~$g^4 \in \GSp_4(\F_3)$
will not be trivial, since it will have order divisible by~$3$.
This allows us to distinguish the classes of types~$2C$ and~$2D$  from
the classes of order divisible by~$3$ by 
computing~$\mathrm{Jac}(X_0)(\F_{16})[3]$. %
  Similarly, if~$g$ is of type~$4C$ or~$12C$, then, for any lift, the element~$g^8 \in \GSp_4(\F_3)$ will
  be trivial if and only if~$g$ has type~$4C$.
We have the following table (where as above~$Q(x)$ denotes the characteristic
polynomial of Frobenius):
\begin{center}
\begin{tabular}{*3c}
\toprule
$Q(x) \mod 3$ & $\langle g \rangle$ & $\dim \ker(g^8 - 1)$ \\
\midrule
$x^4 + 2 x^2 + 1$ & $2C$ & $4$   \\
$x^4 + 2 x^2 + 1$ & $6G$ & $2$   \\
$x^4 + 2 x^2 + 1$ & $6H$ & $2$   \\
\midrule
$x^4 +  x^2 + 1$ & $2D$ & $4$   \\
$x^4 +  x^2 + 1$ & $6I$ & $2$   \\
\midrule
 $x^4 \pm x^3 +2 x^2  \mp  x + 1$  & $4C$ & $4$   \\
$x^4 \pm x^3 +2 x^2  \mp  x + 1$ & $12C$ & $2$   \\
\bottomrule
\end{tabular}
\end{center}

We find in all the~$96$ ordinary cases and~$192$ non-ordinary cases when
$Q(x) \bmod 3$ is a polynomial corresponding
to one of the conjugacy classes in this table,
 there
is an
 isomorphism
 \[\Jac(X_0)(\F_{256})[3] \simeq (\Z/3\Z)^2.\]
This rules out the
 case that~$\rhobar_{A,3}(\Frob_2)$  in~$\PGSp_4(\F_3)$ has order either~$2$ or~$4$.
   When~$Q(x) \bmod 3$ is either $x^4 + x^2 + 1$ or
   $x^4 \pm x^3 +2 x^2  \mp  x + 1$, 
   this suffices to determine
the conjugacy class exactly
for the~$48$ ordinary curves lying in~$\{2D,6I\}$, and the~$144$ non-ordinary  curves lying in~$\{4C,12C\}$.
\item For the remaining~$48$ ordinary curves and~$48$ non-ordinary curves where
the conjugacy class is either of type~$6G$ or~$6H$,
we first compute the degree~$40$ polynomial corresponding to
the~$\PGSp_4(\F_3)$ representation, and then compute the degree~$27$
resolvent, and then use the table in  Lemma~\ref{atlas}.
\item To establish point~(\ref{classify}), we need to show that  all of the ordinary
  curves where~$\rhobar_{B,3}(\Frob_2)$ has conjugacy class~$6G$
  or~$8A$ do not lift to a smooth curve~$X/\Z_2$ with a rational
  Weierstrass point.  By Lemma~\ref{boxerstrass} below, the Jacobian of such a
  curve~$X$ has a rational~$2$-torsion point, so that in particular, the
  polynomial~$Q(x)\pmod{2}$ would need to have~$1$ as a root.
  However, the~$32$ curves of type~$6G$
  and the~$96$ curves of type~$8A$ have the property that
  $Q(x)\equiv x^{2}(x^2+x+1)\pmod{2}$, so no such lift can exist.
\end{enumerate}
It remains to note that if~$\rhobar(\Frob_2)$ is of type~$6H$, $6I$, or~$10A$,
then  an appropriate~$X$ exists by Lemma~\ref{curveexamples}.
\end{proof}

\begin{df}
Say that an abelian variety~$A/\Q_p$ has \emph{semistable ordinary 
reduction} if it has semistable reduction and the abelian part of the special
fibre of the N\'{e}ron model is ordinary. (In particular, good ordinary reduction is a special
case of semistable ordinary reduction.)
\end{df}
Recall from Section~\ref{standardfactstorsion} that we have fixed an
identification $S_6 = \GSp_4(\F_2)$.
\begin{lemma} \label{boxerstrass}
Let~$B/\Q_2$ be an abelian surface with semistable ordinary 
reduction. %
Suppose that the image of~$\rhobar_{B,2}$ lands inside~$S_5(b) \subset S_6 = \GSp_4(\F_2)$.
Then:
\begin{enumerate}
\item The image of~$\rhobar_{B,2}$ is a~$2$-group. \label{boxerstrassadditional}
\item There exists a rational~$2$-torsion point~$P \in B[2](\Q_2)$. \label{boxerstrassold}
\end{enumerate}
In particular, this holds if~$B = \Jac(X)$, and~$X/\Q_2$ has good
ordinary reduction and a rational Weierstrass point.
\end{lemma}

\begin{proof} 
The ordinary assumption implies that the image of
$G_{\Q_2}$ in~$\GSp_4(\F_2) = \Aut(B[2])$ lands (up to conjugation) in the Siegel
parabolic: %
\numequation
\label{siegelparabolicfinite}
\left( \begin{matrix} * & * & * & * \\ * & * & * & * \\ 0 & 0 & * & * \\ 0 & 0 & * & * \end{matrix} \right)
\cap \GSp_4(\F_2).
\end{equation}%
To see this, it suffices to show that~$G_{\Q_2}$ preserves an (isotropic) subspace~$(\Z_2)^2$ inside~$T_2(B)$.
If~$B$ has good ordinary reduction,  the subspace is the kernel of the mod-$2$ reduction.
More generally, if~$B$
has semistable reduction, we can use the description of the Tate module
given in~\cite[Exp.9, IX]{MR0354656}.
There is a~$G_{\Q_2}$-equivariant
 filtration~$T_2(B)_t \subset T_2(B)_f \subset T_2(B)$
of (saturated)~$\Z_2$-modules of ranks~$t>0$ and~$t+2a$ where~$2(t+a)=4$. 
Moreover,
by the orthogonality theorem~\cite[Thm~2.4, Exp.9, IX]{MR0354656},
$T_2(B)^{\perp}_t = T_2(B)_f$.
If~$B$ is purely toric, then~$t=2$ and~$T_2(B)$ gives the desired space.
If~$t=1$, then the abelian part of~$B$ is an abelian variety, and the kernel of reduction gives a rank 
one~$G_{\Q_2}$-stable submodule of~$T_2(B)_f/T_2(B)_t$, and the inverse image of this in~$T_2(B)_f$
is the desired submodule.

This group~(\ref{siegelparabolicfinite}) is a subgroup of~$\GSp_4(\F_2)$ order~$48$ which is isomorphic to~$S_4 \times S_2$, and is the normalizer
of the element:
\numequation
\label{central}
\left( \begin{matrix} 1 & 0 & 1 & 0 \\ 0 & 1 & 0 & 1 \\ 0 & 0 & 1 & 0 \\ 0 & 0 & 0 & 1 \end{matrix} \right).
\end{equation}

There are two
 (non-conjugate) subgroups of order~$48$ in~$S_6$; one given by the centralizer of~$({*}{*})({*}{*})({*}{*})$ in~$S_6$ and
the other by the centralizer of~$({*}{*})$; they are permuted by the outer automorphism.
Under our fixed isomorphism, the element~(\ref{central}) is conjugate
to~$({*}{*})({*}{*})({*}{*})$ by Lemma~\ref{explicit}. (In fact the other conjugacy class of subgroups %
of order~$48$ is given by the Klingen parabolic.)
But now the assumption that the image of~$\rhobar_{B,2}$
lands inside~$S_5(b)$
 implies that the image of~$G_{\Q_2}$ lands inside
 the intersection of~$S_5(b)$ with the centralizer of an element of the form~$({*}{*})({*}{*})({*}{*})$.
 If that intersection is not a~$2$-group, then it contains an element of order~$3$.
 But the conjugacy class of elements of order~$3$ inside the normalizer of~$({*}{*})({*}{*})({*}{*})$ consists
 of elements with cycle
shape~$({*}{*}{*})({*}{*}{*})$, whereas the conjugacy class of elements of order~$3$ in~$S_5(b) \subset S_6$
consists of elements with cycle shape~$({*}{*}{*})$, and thus the intersection is a~$2$-group,
proving~(\ref{boxerstrassadditional}).
 Hence the image of~$\rhobar_{B,2}$ is 
certainly contained within the~$2$-Sylow
of~$\GSp_4(\F_2)$, so the action of $G_{\Q_2}$
fixes a~$2$-torsion point, %
proving part~(\ref{boxerstrassold}).
\end{proof}

\subsection{Abelian surfaces with  semistable ordinary reduction
  at~$2$}%
\label{surfacesattwo}
In this section, we study abelian surfaces~$A/\Q_2$ with either good ordinary or semistable ordinary reduction at~$2$.

Let~$F$ be a number field or a local field of characteristic zero, and 
suppose that~$\rhobar: G_F \rightarrow \GSp_4(\F_3)$  has similitude
character~$\varepsilonbar^{-1}$.
We now recall some rational varieties associated to~$\rhobar$ constructed in~\cite[\S 10.2]{BCGP}.

\begin{df} \label{MandP} Let~$P = P(\rhobar)$  be
 the fine moduli space over~$F$ parametrizing
principally polarized abelian surfaces~$A$ with a given symplectic isomorphism~$A[3] \simeq \rhobar^{\vee}$
and a fixed odd theta characteristic.

Let~$\MMM^w_2(\rhobar)$ be the fine moduli space over~$F$
parametrizing genus two curves~$X/F$ with a fixed Weierstrass point and a fixed
symplectic isomorphism~$\Jac(X)[3] \simeq \rhobar^{\vee}$.
\end{df}More explicitly (see~\cite[Defn.\ 10.2.2]{BCGP}) the
space~$P(\rhobar)$ can be defined as follows: we let~$B$ be the moduli
space of principally polarized abelian surfaces~$A$ with a given
symplectic isomorphism~$A[3] \simeq \rhobar^{\vee}$, and
let~$B(2)\to B$ be the $S_{6}\cong \PSp_4(\F_2)$-cover corresponding to
a full level~$2$ structure. Then~$P$ is the intermediate cover
corresponding to the subgroup~$S_5(b)\subset S_{6}$. In particular we note
that a principally polarized abelian surface~$A/F$ gives rise to a
point in~$P(\rhobar)(F)$ if, in additional to having a symplectic
isomorphism~$A[3] \simeq \rhobar^{\vee}$, the image of~$\rhobar_{A,2}$
is conjugate to a subgroup of~$S_5(b)$.

The space~$P(\rhobar)$ is smooth and rational~\cite[Thm~10.2.3]{BCGP}.
The Torelli 
map~$\MMM^w_2(\rhobar) \rightarrow P(\rhobar)$ is an open immersion,
and hence~$\MMM^w_2(\rhobar)$ is also smooth and rational, and dense in~$P(\rhobar)$. %

 An unramified representation
$$\rhobar: G_{\Q_2} \rightarrow \GSp_4(\F_3)$$
with similitude character~$\varepsilonbar^{-1}$ is given
up to conjugation and unramified twist by a conjugacy class
in~$\PGSp_4(\F_3) \setminus \PSp_4(\F_3)$. Given such a class,
 the goal of this section
is  (when possible) to find a point~$A \in P(\rhobar)(\Q_2)$
which either has good ordinary or semistable ordinary reduction and is in addition~$2$-distinguished.
Naturally, one such source of representations
comes from a point~$X \in \MMM^w_2(\rhobar)(\Q_2)$ where~$X$ has good  ordinary reduction at~$2$, however,
this turns out not to exhaust the list of possibilities. There are three reasons 
for this.  The first is that~$A = \Jac(X)$ can have good reduction even when~$X$ does
not. The second is that~$A$ can have bad reduction and 
yet~$\rhobar_{A,3}$ can still be unramified (although
such~$A$ will necessarily be semistable).  The third is that some of the most accessible points
of~$P$ lie on the complement of the image of~$\MMM^w_2(\rhobar)$, namely, direct products
of elliptic curves.
We exploit a number of these  phenomena to find points for various different representations~$\rhobar$.

Since we shall only consider~$\rhobar: G_{\Q_2} \rightarrow \GSp_4(\F_3)$ which are unramified, 
we begin with
following, which is a (specialization of a) standard result:
\begin{lemma} \label{groth} Let~$A/\Q_2$ be an abelian variety.
Suppose that~$\rhobar_{A,3}$ is unramified. Then~$A$ has semistable reduction.
\end{lemma}

\begin{proof}  The assumption that~$\rhobar_{A,3}$ is unramified implies that the action of inertia on~$T_3(A)$ is
  unipotent, %
  so 
the claim follows from Grothendieck's semi-stability Theorem~\cite[Exp.9, IX]{MR0354656}.
\end{proof}

The ultimate goal of this section is to prove the following theorem:
\begin{theorem} \label{restrictionsattwo}
Let~$\rhobar: G_{\Q_2} \rightarrow \GSp_4(\F_3)$ be unramified with similitude
character~$\varepsilonbar^{-1}$.
There exists a point~$A \in P(\rhobar)(\Q_2)$ which has 
semistable ordinary reduction %
and is~$2$-distinguished 
if and only if the conjugacy class of the image of~$\rhobar(\Frob_2)$ in~$\PGSp_4(\F_3) \setminus \PSp_4(\F_3)$ is not of type~$4C$
or~$12C$. Moreover, one can additionally take~$A$ to have good reduction
if and only if the conjugacy class of the image of~$\rhobar(\Frob_2)$ in~$\PGSp_4(\F_3) \setminus \PSp_4(\F_3)$ is 
of type~$4D$, $6H$, $6I$, or~$10A$. 
This is summarized by the table below.

 \begin{center}
 \begin{tabular}{*3c}
\toprule
 $\langle g \rangle \subset \PGSp_4(\F_3) \setminus \PSp_4(\F_3)$ & good ordinary  &  semistable ordinary \\
\midrule
 $2C$  &\FFFFF & \TTTT \\
 $2D$  &  \FFFFF & \TTTT \\
$4C$ &  \FFFFF &  \FFFFF \\
 $4D$ & \TTTT  &  \TTTT  \\ 
$6G$ &  \FFFFF & \TTTT  \\ 
 $6H$ &  \TTTT & \TTTT \\ 
 $6I$ &  \TTTT & \TTTT \\ 
$8A$ & \FFFFF &  \TTTT   \\
$10A$ & \TTTT &  \TTTT   \\
 $12C$ & \FFFFF &  \FFFFF \\
\bottomrule
 \end{tabular}
 \end{center}
\end{theorem}

\begin{proof}
(Most of) the proof is carried out in the remainder of this section,
and we give the proof by the order in which the argument occurs, namely: %
\begin{enumerate}
\item When~$\langle g \rangle$ is one of the conjugacy classes~$6H$, $6I$, and~$10A$, the result follows directly from
the fact that there exist genus two curves~$X/\Q_2$ with a  rational
Weierstrass point, good ordinary reduction, and with~$2$-distinguished Jacobians,
by  Lemma~\ref{exhaust}(\ref{toolong}).
\item When~$\langle g \rangle$ has the form~$4C$ or~$12C$, the result follows by Lemma~\ref{4C12C}.
\item When~$\langle g \rangle$ has the form~$8A$, the good reduction case is covered by Lemma~\ref{goodcase},
and the   semistable  reduction case by Lemma~\ref{firstconstruction}.
\item When~$\langle g \rangle$ has the form~$2D$ or~$4D$, the semistable reduction case follows
from Lemma~\ref{firstconstruction}, which also covers the good reduction case for the conjugacy class~$4D$.
\item When~$\langle g \rangle$ has the form~$2C$ or~$6G$, the good reduction case follows
from Lemma~\ref{secondgoodcase} and Lemma~\ref{6Hgood} (together with an
examination of Table~\ref{tab:reductions}).
\item The semistable reduction case for the conjugacy class~$2C$
is Lemma~\ref{secondconstruction}.
\item The semistable reduction case for the conjugacy class~$6G$
is Lemma~\ref{thirdconstruction}.
\item  The good reduction case for the conjugacy class~$2D$ is
  Lemma~\ref{2Dgood}. \qedhere
\end{enumerate}
\end{proof}

\begin{df} \label{ordinarypoly} An ordinary Weil polynomial of weight one for~$p$ is
a degree~$4$ polynomial~$X^4 + a X^3 + b X^2 + p a X + p^2 \in \Z[X]$
all of whose roots have absolute value~$p^{1/2}$ and for which~$(b,p) =
1$. 
\end{df}

If~$A/\Q_2$ has good ordinary reduction, then certainly the characteristic 
polynomial~$Q(x)$ of Frobenius at~$2$ will (by the Weil
conjectures) be an ordinary Weil polynomial of weight one for~$p=2$.

 There are~$16$ possible ordinary Weil
polynomials of weight one for~$p=2$, listed in factored form in Table~\ref{tab:reductions},
together with the list of corresponding conjugacy classes in~$\PGSp_4(\F_3) \setminus \PSp_4(\F_3)$
(as described in Lemma~\ref{atlas})
whose conjugacy class admits a lift to~$\GSp_4(\F_3)$ with the given characteristic
polynomial over~$\F_3[x]$.

\numtable
 \begin{center}
 \begin{tabular}{*3c}
\toprule
 $Q(x)$ &    $Q(x) \bmod 3$ & $\langle g \rangle$ \\
\midrule
$x^4 - x^2  + 4$ & $x^4 - x^2 + 1$ & $2C$, $6G$, $6H$ \\
$x^4 - 3 x^3 + 5 x^2 - 6 x + 4$ & $x^4 - x^2 + 1$ \\
$x^4 + 3 x^3 + 5 x^2 + 6 x + 4$ & $x^4 - x^2 + 1$ \\
\hline
$x^4 + x^2 + 4$ & $x^4 + x^2 + 1$ & $2D$, $6I$ \\
\hline
$(x^2 - x + 2)^2$ & $x^4 + x^3 + 2 x^2 + 2 x + 1$ &  $4C$, $12C$  \\
$x^4 + x^3 - x^2 + 2 x + 4$ &  $x^4 + x^3 + 2 x^2 + 2 x + 1$ \\
\hline
$(x^2 + x + 2)^2$ & $x^4 + 2 x^3 + 2 x^2 + x + 1$  &  $4C$, $12C$ \\
$x^4 - x^3 - x^2 - 2 x + 4$ & $x^4 + 2 x^3 + 2 x^2 + x + 1$ \\
\hline
$x^4 - 3 x^2 + 4$ & $x^4 + 1$ & $4D$ \\
$(x^2 + x +2)(x^2 - x + 2)$ & $x^4 + 1$ &  \\
\hline
$x^4 + x^3 +  x^2 + 2 x  +4$ & $x^4 + x^3 + x^2 + 2 x + 1$  & $8A$ \\
\hline
$x^4 - x^3 + x^2 - 2 x + 4$ & $x^4 + 2 x^3 + x^2 + x + 1$ & $8A$ \\
\hline
$x^4 - 2 x^3 + 3 x^2 - 4 x + 4$ & $x^4 + x^3 + 2 x + 1$ & $10A$ \\
$x^4 + x^3 + 3  x^2 + 2 x + 4$ & $x^4 + x^3 + 2 x + 1$ \\
\hline
$x^4 + 2 x^3 + 3 x^2 + 4 x + 4$ & $x^4 + 2 x^3 + x + 1$ & $10A$ \\
$x^4 - x^3 + 3 x^2 - 2 x + 4$ & $x^4 + 2 x^3 + x + 1$ \\
\bottomrule
 \end{tabular}
 \caption{Mod~$3$ reduction of ordinary Weil polynomials of weight one for~$p=2$}
 \label{tab:reductions}
 \end{center}
 \end{table}

\begin{lemma} \label{goodcase}
There does not exist a principally polarized abelian surface~$A/\Q_2$ with good ordinary reduction at~$2$ such
that~$\rhobar_{A,3}(\Frob_2)$ has conjugacy class~$8A$, and~$\rhobar_{A,2}$
has image inside~$S_5(b)$.
If one further insists that~$A$ is~$2$-distinguished, then~$\rhobar_{A,3}(\Frob_2)$
can not have conjugacy class~$4C$ and~$12C$.
\end{lemma}

\begin{proof} Consider first the case of~$4C$ and~$12C$. Up to unramified twist,
the characteristic polynomial of~$\rhobar_{A,3}(\Frob_2)$ is then
$$(x^2 - x + 2)^2 \bmod 3.$$
If~$Q(x) \equiv (x^2 - x + 2)^2 \bmod 3$ is an ordinary Weil polynomial, this forces
(by Table~\ref{tab:reductions})
either the equality~$Q(x) = (x^2 - x + 2)^2$ or
$$Q(x) = x^4 + x^3 - x^2 + 2 x + 4 \equiv x^2 (x^2+x+1) \bmod 2.$$
The first case is ruled out by the~$2$-distinguished condition. For the second, it implies
that the action of~$\Frob_2$ on~$A[2](\Fbar_2)$ has order~$3$.  But
this contradicts Lemma~\ref{boxerstrass}.

Now consider~$8A$. The characteristic polynomial up to twist is
$$(x+1)(x-1)(x^2+x+2) \bmod 3,$$ 
in which case (see Table~\ref{tab:reductions}) there is a unique possibility
$$Q(x) = x^4 + x^3 +  x^2 + 2 x  +4 \equiv x^2(x^2+x+1) \bmod 2,$$
and we are again done by Lemma~\ref{boxerstrass}.
\end{proof}

We can upgrade this lemma as follows:

\begin{lemma} \label{4C12C}
Suppose that~$A/\Q_2$ is a principally polarized abelian surface
with \emph{(}potentially\emph{)} semistable ordinary reduction 
and such that:
\begin{enumerate}
\item $\rhobar_{A,3}$ is unramified. %
\item $\rhobar_{A,3}(\Frob_2)$ has projective conjugacy class~$4C$ or~$12C$.
\item $\rhobar_{A,2}$
has image inside some conjugate of~$S_5(b)$.
\end{enumerate}
Then~$A$ has good reduction at~$2$ and is not $2$-distinguished.
\end{lemma}

\begin{proof} %
 The conditions imply that~$\rhobar_{A,3}(\Frob_2)$ has, up to unramified twist,
 characteristic polynomial
 $$(x^2 + x + 2)^2 \bmod 3.$$
 If~$\alpha$ is a root of this polynomial, then~$-\alpha$ is clearly not a root. But that implies
 that no two roots have ratio~$2 \equiv -1 \bmod 3$. In particular we
 see that~$H^2(G_{\Q_2},\rhobar_{A,3})=0$, so all lifts
 of~$\rhobar_{A,3}$ are unramified.   %
In particular~$\rho_{A,3}|_{G_{\Q_2}}$ is unramified, and by N\'{e}ron--Ogg--Shafarevich, we deduce that~$A$ must 
have good ordinary reduction at~$2$.
The result  now follows from Lemma~\ref{goodcase}.
 \end{proof}
 
We now move on to the classes~$2D$, $4D$, and~$8A$,
which we can construct directly using products of elliptic curves.

\begin{lemma} \label{firstconstruction} 
Let~$X/\Q_2$ be an elliptic curve with
split multiplicative reduction and such that the Tate parameter~$q$
is a perfect~$6$th power.
Let~$Y/\Q_2$ be an elliptic curve with good ordinary reduction with characteristic
polynomial~$Q(x) = x^2+x+2$ and with~$\rhobar_{Y,2}$ trivial.
Let~$X'$ and~$Y'$ denote the
unramified quadratic twists of~$X$ and~$Y$ respectively.
Let~$A = X \times Y$, and let~$B = Y \times Y'$, and~$C = X \times X'$.
Then~$A$, $B$, and~$C$ are
principally polarized abelian surfaces with the following properties:
\begin{enumerate}
\item $A$ and~$C$ have semistable ordinary reduction,
and~$B$ has good ordinary reduction.
\item $\Q(A[3])$, $\Q(B[3])$, and~$\Q(C[3])$ are unramified at~$2$.
\item $\rhobar_{A,3}(\Frob_2)$ has projective conjugacy class~$8A$,
$\rhobar_{B,3}(\Frob_2)$ has projective conjugacy class~$4D$,
and~$\rhobar_{C,3}(\Frob_2)$ has projective conjugacy class~$2D$.
\item $\rhobar_{A,2}$, $\rhobar_{B,2}$, and~$\rhobar_{C,2}$
have image inside~$S_5(b) \subset \GSp_4(\F_2)$ up to conjugacy.
\item $A$, $B$, and~$C$ are~$2$-distinguished.
\end{enumerate}
\end{lemma}

\begin{proof} First we note that both~$X$ and~$Y$ exist; there exists a Tate
curve for any~$q \in \Q_2$ with~$v(q) > 0$, and one can take~$Y$ to
be~$y^2 + x y  + y = x^3 - x^2 - 6 x - 4$,
which is actually the base change to~$\Q_2$ of an elliptic
curve~$E$ over~$\Q$ of conductor~$17$ with~$E[2] \simeq (\Z/2\Z)^2$
as a~$G_{\Q}$-module.

The surfaces~$A$, $B$, and~$C$ have a principal
polarization coming from the principal polarization on each
elliptic curve. They clearly all have
 semistable ordinary reduction, and~$B$
in addition has good ordinary reduction. The assumption that the Tate parameter~$q$ is a cube implies
that
the action of~$G_{\Q_2}$ on~$X[3]$ is isomorphic to~$\mu_3 \oplus \Z/3\Z$ as a Galois
representation, and hence is unramified. Moreover, we deduce
that there is also an isomorphism~$X'[3] \simeq \Z/3\Z \oplus \mu_3$,
since the unique unramified quadratic character is the cyclotomic character. %
 The second claim then follows
since~$Y$ and~$Y'$ have good reduction.

 The element~$\rhobar_{Y,3}(\Frob_2) \in \GL_2(\F_3)$
has characteristic polynomial~$x^2 + x + 2 \bmod 3$, and~$\rhobar_{X,3}(\Frob_2)$
has characteristic polynomial~$(x^2 - 1)$ from the description above.
Thus~$\rhobar_{A,3}(\Frob_2)$ has characteristic polynomial
$$(x^2 - 1)(x^2 + x + 2) = x^4 + x^3 + x^2 - x + 1 \bmod 3.$$
The only elements with this characteristic polynomial modulo~$3$
have conjugacy class~$8A$ in~$\PGSp_4(\F_3)$.
For~$B$, we see that~$Q(x) = (x^2 + x + 2)(x^2 - x + 2)$,
and thus (from Table~\ref{tab:reductions}) the only possibility
is that~$\rhobar_{B,3}(\Frob_2)$
has projective conjugacy class~$4D$.
For~$C$, we see that the characteristic polynomial
of Frobenius on~$C[3]$ is~$(x^2-1)^2 \bmod 3$ and that Frobenius
clearly has order~$2$, so the projective conjugacy class is~$2D$.

We now show that the mod~$2$ reductions
are conjugate to a subgroup of~$S_5(b)$. The assumption that~$q$ is a square implies
that~$\rhobar_{X,2}$ (and its quadratic twist) are trivial. On the other hand, by construction,~$\rhobar_{Y,2}$ and thus its quadratic twist are also trivial.
So~$\rhobar_{A,2}$, $\rhobar_{B,2}$, and~$\rhobar_{C,2}$ are also trivial and the claim
follows.

It remains to show that~$A$, $B$, and~$C$ are~$2$-distinguished.
In each case, we can compute directly the  unit Frobenius eigenvalues
on the semi-simplification of the Tate module. If~$\alpha$ denotes
the unit root of~$x^2 + x + 2 = 0$, then for~$A$, $B$, and~$C$
they are given by~$\{1,\alpha\}$, $\{\alpha,-\alpha\}$, and~$\{1,-1\}$
respectively. Since~$\alpha \ne 1$ and~$\alpha \ne - \alpha$, these pairs  all consist
of distinct numbers
and we are done.
\end{proof}

\subsubsection{The cases~$2C$, $6G$, and~$6H$.}
We now turn to the cases of~$2C$, $6G$, and~$6H$, where (see Table~\ref{tab:reductions})
the characteristic polynomial of~$\rhobar_{A,3}(\Frob_2)$
is
\numequation
\label{2Cpoly}
x^4 - x^2 + 1 = (x^2 + 1)^{2} \bmod 3
\end{equation}

We have the following:

\begin{lemma} \label{secondgoodcase}
If~$A/\Q_2$ has good ordinary reduction
and
$$Q(x) = x^4 \pm 3 x^3 + 5 x^2 \pm 6x + 4
\equiv x^2(x^2+x+1) \bmod 2,$$
then the image of~$\rhobar_{A,2}$ is not
conjugate to a subgroup of~$S_5(b)$.
\end{lemma}

\begin{proof} This follows from  Lemma~\ref{boxerstrass}, exactly as in the proof of Lemma~\ref{goodcase}.
\end{proof}

Lemma~\ref{secondgoodcase} implies that, in the good reduction case
with conjugacy classes~$2C$, $6G$, and~$6H$,
the only possibility for~$Q(x)$ is~$x^4 - x^2 + 4$.

\begin{lemma} \label{6Hgood} Let~$A/\F_2$ be an
 abelian surface
with~$Q(x) = (x^4 - x^2 + 4)$.
Then the action of~$\Frob_2$ on~$A[3](\Fbar_2)$
is not semi-simple and the projective
image of~$\rhobar_{A,3}(\Frob_2)$ has conjugacy class~$6H$.
\end{lemma}

\begin{proof}Note firstly that the abelian surface~$A = \Jac(C_3)$ %
where ~$C_3$ is as in
 Lemma~\ref{curveexamples} satisfies the hypothesis and conclusions of the
 lemma. Let~$B$ be another abelian surface with~$Q(x)= x^4 - x^2 +
 4$; it suffices to show that~$A[3]\cong B[3]$. There is an inclusion~$\Z[\phi] \subset \End(A)$
where~$\phi$ is the Frobenius endomorphism which
satisfies~$\phi^4 - \phi^2 + 4 = 0$. Let~$\psi = \phi^2 + 1$,
so~$\psi^{\vee} = 3 - \psi$ and~$\psi^{\vee} \circ \psi = [6]$.
 Note that~$\Frob_2$ acts non-semi-simply
on~$A[3]$.  Since~$\Frob_2$ on~$A[3]$ has characteristic polynomial~$(x^2+1)^2 \bmod 3$, and~$x^2+1 \bmod 3$ is irreducible, it follows
that the only proper~$\Gal(\Fbar_2/\F_2)$-equivariant submodule
of~$A[3]$ is $\ker(\psi)\cap A[3]$.

Since~$Q(x)$ determines  ~$A$ up to
isogeny, there is an isogeny ~$\chi:A\to B$. %
Either~$\chi$ has order prime to~$3$ or
$\ker(\psi) \cap A[3] \subset \ker(\chi)$, in which case
there is a factorization:
 \[\xymatrix{A \ar[r]^{\chi}  \ar[d]^{\psi} & B \ar[d]^{[2]}  \\
    A \ar[r]^{\chi'} & B}  \] 
where the~$3$-part of the degree of~$\chi'$
is less than that of~$\chi$. By induction, there exists an isogeny of~$A$ to~$B$ of order prime to~$3$,
which implies that~$A[3]$ and~$B[3]$ are isomorphic, as required.
\end{proof}

\begin{lemma} \label{badisgood}
Suppose that~$A/\Q_2$ is an abelian surface such that:
\begin{enumerate}
\item \label{firstconditionbadisgood} $\rhobar_{A,3}|_{G_{\Q_2}}$ is unramified. %
\item \label{secondconditionbadisgood}  $\rhobar_{A,3}(\Frob_2)$ has projective conjugacy class~$2C$
or~$6G$. 
\item \label{thirdconditionbadisgood} The image of~$\rhobar_{A,2}$ is conjugate to a subgroup
of~$S_5(b)$.
\end{enumerate}
Then~$A/\Q_2$  has semistable ordinary reduction with  purely toric reduction.
\end{lemma}

\begin{proof}
Suppose that~$A$ had good reduction. Condition~(\ref{secondconditionbadisgood}) implies
that~$Q(x) \equiv x^4 - x^2 + 1 \bmod 3$.
By Lemmas~\ref{secondgoodcase} and~\ref{6Hgood}, 
 assumptions~(\ref{firstconditionbadisgood}),  and~(\ref{thirdconditionbadisgood}) 
(taking into account Table~\ref{tab:reductions})   imply that $\rhobar_{A,3}(\Frob_2)$
 has projective conjugacy class~$6H$, 
contradicting condition~(\ref{secondconditionbadisgood}). Thus~$A$ cannot have good reduction.
By Lemma~\ref{groth}, $A$ has semistable reduction.
As in the proof of Lemma~\ref{boxerstrass} there is a~$G_{\Q_2}$-equivariant
 filtration~$T_3(A)_t \subset T_3(A)_f \subset T_3(A)$
of (saturated)~$\Z_3$-modules of ranks~$t$ and~$t+2a$
where~$2(t+a)=4$. Since~$A$ does not  have good reduction, we have~$t>0$.
If~$t=1$, then~$T_3(A)_t/3$ is a Galois invariant line
inside~$A[3]$, but this is not compatible with the fact
that $\rhobar_{A,3}(\Frob_2)$ has characteristic polynomial~$(x^2+1)^2$,
and~$(x^2+1)$ has no roots over~$\F_3$.
So~$t=2$ and~$A$ has purely toric (and hence semistable ordinary) reduction.
\end{proof}

On the other hand, we have the
following variation
on Lemma~\ref{firstconstruction}.

\begin{lemma} \label{secondconstruction} 
There exists a principally polarized abelian surface $A/\Q_2$ satisfying the following:
\begin{enumerate}
\item $A$ has semistable ordinary reduction.
\item $\rhobar_{A,3}$ is unramified. %
\item $\rhobar_{A,3}(\Frob_2)$ has projective conjugacy class~$2C$,
\item $\rhobar_{A,2}$ has image inside~$S_5(b) \subset \GSp_4(\F_2)$ up to conjugacy.
\item $A$ is~$2$-distinguished.
\end{enumerate}
\end{lemma}

\begin{proof}
Let~$X/\Q_2$ be an elliptic curve with
split multiplicative reduction and such that the Tate parameter~$q \in \Q^{\times}_2$
is a perfect cube, and that~$q \in 5 \cdot (\Q^{\times}_2)^2$.  With these choices
$\rhobar_{X,3}$ is the trivial representation while $\rhobar_{X,2}$ is unramified and $\rhobar_{X,2}(\Frob_2)$ is conjugate to $\left(\begin{matrix} 1 & 1 \\ 0 & 1 \end{matrix} \right)$.

We now take $A=X\otimes\Z^2$ where $G_{\Q_2}$ acts on $\Z^2$ via an unramified
quotient with~$\Frob_2$
acting by
\numequation
\label{isselfdual}
\left(\begin{matrix} 0 & -1 \\ 1 & 0 \end{matrix} \right)
\end{equation}
or in other words $A$ is descended from $(X\times X)_{\Q_{16}}$ where Frobenius is twisted by the automorphism $(a,b)\mapsto (b,-a)$.

Certainly~$A$ is semistable ordinary. Since~$X$
is principally polarized and since the map~$G_{\Q_2} \rightarrow \GL_2(\Z)$
associated to~\eqref{isselfdual} is  self-dual,
it follows that~$A$ is isomorphic to its dual and hence 
is also principally polarized.
Moreover we have, up to conjugation,
$$\rhobar_{A,3}(\Frob_2)
= \left(\begin{matrix} 1 & 0 \\ 0 & 1 \end{matrix} \right)
\otimes  \left(\begin{matrix} 0 & -1 \\ 1 & 0 \end{matrix} \right)
= \left( \begin{matrix} 0 & 1 & 0 & 0 \\ -1 & 0 & 0 & 0 
\\ 0 & 0 & 0& -1 \\ 0 & 0 & 1 & 0 \end{matrix}\right),$$
with order~$4$ and projective image~$2C$.
On the other hand,
$$\rhobar_{A,2}(\Frob_2)
= \left(\begin{matrix} 1 & 1 \\ 0 & 1 \end{matrix} \right)
\otimes  \left(\begin{matrix} 0 & -1 \\ 1 & 0 \end{matrix} \right)
=\left( \begin{matrix} 0 & 1 & 0 & 1 \\ 1 & 0 & 1 & 0 
\\ 0 & 0 & 0& 1 \\ 0 & 0 & 1 & 0 \end{matrix}\right).$$
The centralizer of this element
 has order~$16$, which means
 (with respect to the isomorphism
in Lemma~\ref{explicit}, although this conjugacy
class is preserved by the outer automorphism) that it is conjugate
to~$({*}{*})({*}{*})$ in~$S_6$ and so is conjugate to an element of~$S_5(b)$.
\end{proof}

This leads us to~$6G$ as the last remaining semistable
case of Theorem~\ref{restrictionsattwo}. %
There is a construction in this case along the lines
of Lemma~\ref{secondconstruction}, although
the details are more cumbersome. Instead, we use a different
idea motivated by Lemma~\ref{badisgood}. 
  Consider the genus two curve
$$Y: y^2 + y =x^5 - x^4 + x^3$$
with good reduction at~$2$ (it has conductor~$797$).
We have~$Q(x) \equiv x^4 - x^2 + 1 \bmod 3$,
and so (the projective image of)~$\rhobar_{\Jac(Y),3}(\Frob_2)$ 
has conjugacy class~$2C$,
$6G$, or~$6H$. By computing the corresponding
degree~$40$ and degree~$27$ polynomials, we find
that it has conjugacy class~$6G$. This does not
contradict Lemma~\ref{badisgood} because~$Q(x) = x^4 + 2 x^2 + 1$
and~$Y$ is not ordinary.
We may write~$Y/\Q$ as
\numequation
\label{Xactually}
Y: y^2 = x^5 + a x^3 + b x^2 + c x + d,
\end{equation}
with
$$(a,b,c,d) =  \left(\frac{48}{5}, \frac{704}{25}, \frac{3072}{125},
\frac{821504}{3125} \right).$$
We can think of~$\MMM^{w}_2(\rhobar)$ explicitly as the moduli space of genus two
curves~$X$ given by~$y^2 = x^5 + A x^3 + B x^2 + C x + D$
with a (symplectic) isomorphism~$\rhobar_{\Jac(X),3} \simeq \rhobar_{\Jac(Y),3}$.
In~\cite[Thm~2]{CCR}, an explicit parametrization~$\mathbf{P}^3_{\Q}(s,t,u,v)
\rightarrow \MMM^{w}_2(\rhobar)$ is given; that is, $A$, $B$, $C$, and~$D$
are explicit polynomials in~$(s,t,u,v)$ whose specialization to~$(1,0,0,0)$
gives the parameters~$(a,b,c,d)$ of equation~(\ref{Xactually}). 
By Lemma~\ref{badisgood}, any specialization of this family
which does \emph{not} have good reduction is necessarily
semistable ordinary with purely toric reduction.
Moreover,  it will also necessarily be~$2$-distinguished;
the pair of eigenvalues of Frobenius on the  unramified quotient will
be Galois invariant and yet be roots of~$(x^2+1) \bmod 3$.
Thus in practice we can choose random points on this family
to find one which does not have good reduction,
and then we are done.
The specialization of this family
to the point~$(0,0,4,1)$ is the curve~$y^2=x^5+ A x^3 + B x^2 + C x + D$, with
(after scaling  down by~$(4^{36},8^{36},16^{36},32^{36})$ from the formulas in~\cite{CCR})
$$\begin{aligned}
A = & \ 672315215064342/5, \\
B= & \ -197745818620367722373332/25, \\
C = & \ -3038748471428312132304651799323/125, \\
D =  & \ 405130036222076498453650257209001453372/3125. \end{aligned}$$
Thus we have  produced a curve of the required form.
As a sanity check, the conductor has the form~$2^2 \cdot N$ where~$(2,N)=1$, and the Euler factor
at~$2$ is~$x^2+1$, and one can indeed compute the~$3$-torsion division polynomial of degree~$40$
and its resolvent of degree~$27$ and find that they give an extension unramified at~$2$ 
with~$\rhobar(\Frob_2)$
of conjugacy class~$6G$ (as they should).
On the other hand, after reducing this modulo a large enough power of~$2$ ($2^{12}$ in this case),
we get a more manageable example
(except now a different \emph{global} representation):
\begin{lemma} \label{thirdconstruction} Let~$X$ be the curve
$$y^2 + (-x^2 - x - 1)y = x^5 - 4x^4 + 47x^3 - 
    43x^2 - 2x + 8.$$
    Then~$A =\Jac(X)/\Q_2$ has purely toric reduction, and
    is~$2$-distinguished. Furthermore
    $A[3]$ is unramified, $\rhobar_{A,2}$ has image conjugate
    to a subgroup of~$S_5(b)$, and~$\rhobar_{A,3}(\Frob_2)$ has conjugacy class~$6G$.
    \end{lemma}
    
    \begin{proof}
We may compute directly using division polynomials that~$A[3]$
is unramified at~$2$ and~$\rhobar_{A,3}(\Frob_2)$ has conjugacy class~$6G$.
Since~$X$ has a~$\Q_2$ (even a~$\Q$) Weierstrass point, the image of~$\rhobar_{A,2}$
is conjugate to a subgroup of~$S_5(b)$. The conductor at~$2$ of~$A$ is~$2^2$,
so by Lemma~\ref{badisgood}, $A$ has purely toric reduction,
and the eigenvalues of Frobenius on the unramified quotient are~$\pm i$, so~$A$
is $2$-distinguished.
    \end{proof}

\subsubsection{The case~$2D$} We finish the proof
of Theorem~\ref{restrictionsattwo} by ruling out~$2D$ in the case of good ordinary reduction.

\begin{lemma} \label{2Dgood}
There does not exist 
an abelian surface $A/\Q_2$
with good ordinary reduction and
$\rhobar_{A,3}(\Frob_2)$ projectively
conjugate to $2D$.
\end{lemma}

\begin{proof}
From Table~\ref{tab:reductions}, we see that such an~$A$ must
satisfy~$Q(x) = x^4 + x^2 + 4 \equiv (x-1)^2 (x+1)^2 \bmod 3$. From 
Lemma~\ref{curveexamples} (in particular the curve~$C_2$),
 we see that
there exists a smooth ordinary~$X/\F_2$  with the
same~$Q(x)$, and thus~$\Jac(X)$ is isogenous to~$A/\F_2$. 
Let~$\chi: \Jac(X) \rightarrow A$ be an isogeny, which we may assume
is not divisible by~$[3]$. 
Since~$\rhobar_{\Jac(X),3}(\Frob_2)$ is conjugate to~$6I$, the
\emph{minimal} polynomial of $\Frob_2 $ on~$A[3]$ is~$(x-1)^2(x+1)^2$.
Thus the only Galois invariant subspaces of~$A[3]$ contain either the
intersection of~$A[3]$ with the kernel of~$\psi = \phi - 1$ or~$\psi = \phi+1$ respectively.
 Thus (as in Lemma~\ref{6Hgood}) we may reduce to the case when~$\chi$
 has degree prime to~$3$, which implies that~$\rhobar_{A,3}(\Frob_2)$ is conjugate
 to~$\rhobar_{\Jac(X),3}(\Frob_2)$ of class~$6I$.
\end{proof}

\subsection{Genus~$2$ curves  locally at~$3$} \label{localatthree}
In this section, we carry out some computations similar to~\S\ref{localattwo}
except now over~$\F_3$. (The \texttt{magma} files for these computations
can also be found in~\cite{magma}, as noted in Remark~\ref{magmacite}.)

Recall that an irreducible monic polynomial~$Q(x) \in \Z[x]$  with roots of absolute value~$q^{1/2}$ 
(for~$q$ some power of~$p$) 
 corresponds  (by Honda--Tate theory) to an isogeny class of simple
abelian varieties~$A$ of dimension~$\frac{\deg(Q)}{2} \cdot [E:F]^{1/2}$ over~$\F_q$, where~$F = \Q(\alpha) = \Q[x]/Q(x)$
 and~$E$ is a certain division algebra whose centre
is~$F$ and whose invariants (also determined by~$Q(x)$) are trivial away from primes dividing~$p$ and $\infty$. However, if
one also assumes that~$Q(x)$ is ordinary in the sense that it has degree $2g$ and for~$g$ of the embeddings~$F \hookrightarrow \Qbar_p$,
the valuation of~$\alpha$ is zero, this
forces $F$ to be totally complex and the invariants of~$E$ to be trivial at~$v|p$, which implies
that~$E=F$ and~$\dim(A) = g$.

Specializing to the case~$g=2$, 
recall (Definition~\ref{ordinarypoly}) that by an ordinary Weil polynomial of weight one for~$p$, we mean
a degree~$4$ polynomial~$X^4 + a X^3 + b X^2 + p a X + p^2 \in \Z[X]$
all of whose roots have absolute value~$p^{1/2}$ and for which~$(b,p) =
1$. 

\begin{lemma} \label{compute}
Table~\ref{table1} contains the following data
concerning all pairs consisting of a  smooth genus two curve together with an explicit Weierstrass equation~$X: y^2 = f(x)$
with~$f(x) \in \F_3[x]$.
The columns indicate:
\begin{itemize}
\item All~$40$ ordinary Weil polynomials~$Q(x)$ of weight one for~$p = 3$.
By Honda--Tate theory, these correspond to 
 isogeny classes of ordinary abelian 
surfaces~$A/\F_3$.
\item The reduction of~$Q(x) \bmod 3$.
\item Whether the isogeny class of~$A/\F_3$ contains 
the Jacobian of an ordinary curve~$X/\F_3$ with 
 a rational
Weierstrass point.
\item Whether the isogeny class of~$A/\F_3$ contains 
the Jacobian of an ordinary curve~$X/\F_3$.%
  \item How many such~$X$ have a Jacobian with the corresponding~$Q(x)$.
\end{itemize}
Of the~$3^7$ possible~$f(x)$ of degree~$\le 6$, we find that:
\begin{enumerate}
\item There are~$1296$ curves which are smooth of genus~$2$.
\item There are~$864$ ordinary curves.
\item \label{troll} Exactly~$10$ of these ordinary curves are
  not~$3$-distinguished; equivalently, the polynomial~$Q(x)$ is a
  square.
Moreover, these are precisely the curves for which:
$$\Jac(X)(\Fbar_3)[3] \simeq \chi \otimes (\Z/3\Z)^2$$
as a~$G_{\F_3}$-representation, where~$\chi^2 = 1$. None of these curves have a rational Weierstrass point.
\end{enumerate}
If one enumerates  curves together with a generalized Weierstrass equation
$$y^2 + h(x) y = f(x),$$
where~$\deg(f(x)) \le 6$ and~$\deg(h(x)) \le 3$, all the relative ratios
remain unchanged. 
\end{lemma}

\begin{proof}%
This is a straightforward computation, although we explain
  point~(\ref{troll}). %
If there is an isomorphism~$\Jac(X)(\Fbar_3)[3] \simeq \chi \otimes (\Z/3\Z)^2$ as a
~$G_{\F_3}$-representation, then $Q(x)\equiv x^{2}(x\pm 1)^2\pmod{3}$.
Similarly, if $Q(x)$ is a square then it is a square modulo~$3$, and thus
$Q(x)\equiv x^{2}(x\pm 1)^2\pmod{3}$.

We may make a quadratic twist to reduce to the case $\chi=1$ and $Q(x)\equiv x^2(x-1)^2$.  We are reduced to checking that if $\Jac(X)$ is in the isogeny class corresponding
to~$Q(x)$ with $Q(x)\equiv x^{2}(x-1)^2\pmod{3}$, then ~$Q(x)$ is a square if
and only if $\Jac(X)(\F_3)[3] =  (\Z/3\Z)^2$.  One checks this directly for each of the $4+24+1+8+48+8+48=141$ curves~$X$ corresponding to
such a~$Q(x)$.  One further checks that the 5 such $X$ where $Q(x)$ is a square do not have rational Weierstrass points. %
\end{proof}
\begin{rem}
  The fact that some~$Q(x)$ in Table~\ref{table1} do not arise from
any~$X$  means that there exist ordinary abelian surfaces over~$\F_3$ which are not
isogenous to Jacobians of genus two curves. The simple (although not absolutely simple)
examples
in our table (with~$Q(x) = x^4 - 5 x^2+9$ and~$Q(x) = x^4 - 4 x^2 + 9$)
actually generalize to similar examples over~$\F_q$ for any odd~$q$, see~\cite{MR2041770}.
\end{rem}

\numtable
\begin{tabular}{*5c}
\toprule
$Q(x)$ & $Q(x) \bmod 3$ & 
 $\mathrm{Jac}(X_{\mathrm{WP}})$ & $\mathrm{Jac}(X)$
& $\#X$ \\
\midrule
$9 - 5x^2 +x^4 $ & $x^2 +x^4$ & \FFFF & \FFFF & $0$ \\
$9 - 2x^2 +x^4$ &  $x^2 +x^4$  &\TTTT  &  \TTTT &  $30$ \\ 
$9 +x^2 +x^4$ &  $ x^2 +x^4$  &  \FFFF &  \TTTT & $24$  \\
   $9 + 4x^2 +x^4$ &  $x^2 +x^4$ &  \TTTT   &  \TTTT & $24$ \\ 
 $9 - 9x + 7x^2 - 3x^3 +x^4$ &  $x^2 +x^4$    &\TTTT   &  \TTTT & $24$ \\ 
 $9 + 9x + 7x^2 + 3x^3 +x^4$ &  $x^2 +x^4$  & \TTTT  &  \TTTT & $24$  \\ 
 \midrule
 $9 - 4x^2 +x^4$ & $2x^2 +x^4$ &   \FFFF & \FFFF  & $0$ \\
 $(3 -x +x^2)(3 +x +x^2)$ &  $ 2x^2 +x^4$   & \FFFF  &  \TTTT & $24$ \\ 
 $(3 - 2x +x^2)(3 + 2x +x^2)$ &  $ 2x^2 +x^4$   &\TTTT   &  \TTTT & $36$  \\ 
  $(3 - 2x +x^2)(3 -x +x^2)$ &  $ 2x^2 +x^4$ &  \FFFF   &  \FFFF & $0$ \\ 
 $(3 + 2x +x^2)(3 +x +x^2)$ &  $ 2x^2 +x^4$ &  \FFFF   & \FFFF & $0$  \\ 
  $9 -x^2 +x^4$ &  $ 2x^2 +x^4$ &  \FFFF &  \TTTT  &    $24$ \\ 
 $9 - 9x + 5x^2 - 3x^3 +x^4$ &  $ 2x^2 +x^4$   &\TTTT   &  \TTTT & $24$ \\ 
 $9 + 9x + 5x^2 + 3x^3 +x^4$ &  $ 2x^2 +x^4$   & \TTTT  &  \TTTT & $24$ \\ 
  \midrule
 $(3 -x +x^2)^2$ &  $ x^2 + x^3 + x^4$ & \FFFF  &  \TTTT & $4$  \\ 
 $(3 - x + x^2)(3 + 2x + x^2)$ &  $ x^2 + x^3 + x^4$    & \TTTT  &  \TTTT & $24$  \\ 
 $(3 + 2x + x^2)^2$ &  $ x^2 + x^3 + x^4$    & \FFFF  &  \TTTT & $1$  \\ 
 $9 - 6x + x^2 - 2x^3 + x^4$ &  $ x^2 + x^3 + x^4$  & \FFFF  &  \TTTT & $8$ \\ 
 $9 - 6x + 4x^2 - 2x^3 + x^4$ &  $ x^2 + x^3 + x^4$   & \TTTT  &  \TTTT & $48$ \\ 
 $9 + 3x - 2x^2 + x^3 + x^4$ &  $ x^2 + x^3 + x^4$  & \TTTT  &  \TTTT & $8$  \\ 
 $9 + 3x + x^2 + x^3 + x^4$ &  $ x^2 + x^3 + x^4$   &  \TTTT  &  \TTTT & $48$  \\ 
  \midrule
 $(3 - 2x + x^2)^2$ &  $ x^2 + 2x^3 + x^4$   & \FFFF  &  \TTTT & $1$ \\ 
 $(3 - 2x + x^2)(3 + x + x^2)$ &  $ x^2 + 2x^3 + x^4$    & \TTTT  &  \TTTT & $24$  \\ 
 $(3 + x + x^2)^2$ &  $ x^2 + 2x^3 + x^4$ &   \FFFF  &  \TTTT & $4$  \\ 
 $9 - 3x - 2x^2 - x^3 + x^4$ &  $ x^2 + 2x^3 + x^4$  & \TTTT   &  \TTTT & $8$ \\ 
 $9 - 3x + x^2 - x^3 + x^4$ &  $ x^2 + 2x^3 + x^4$   & \TTTT   &  \TTTT & $48$  \\ 
 $9 + 6x + x^2 + 2x^3 + x^4$ &  $ x^2 + 2x^3 + x^4$  &  \FFFF  &  \TTTT & $8$ \\ 
 $9 + 6x + 4x^2 + 2x^3 + x^4$ &  $ x^2 + 2x^3 + x^4$   & \TTTT   &  \TTTT & $48$  \\ 
  \midrule
 $9 - 12x + 8x^2 - 4x^3 + x^4$ &  $ 2x^2 + 2x^3 + x^4$  &   \FFFF  &  \TTTT & $6$  \\ 
 $9 - 3x - x^2 - x^3 + x^4$ &  $ 2x^2 + 2x^3 + x^4$   & \TTTT  &  \TTTT & $24$ \\ 
 $9 - 3x + 2x^2 - x^3 + x^4$ &  $   2x^2 + 2x^3 + x^4$   & \TTTT   &  \TTTT & $48$ \\
  $9 - 3x + 5x^2 - x^3 + x^4$ &  $  2x^2 + 2x^3 + x^4$   & \TTTT  &  \TTTT & $24$ \\ 
 $9 + 6x + 2x^2 + 2x^3 + x^4$ &  $ 2x^2 + 2x^3 + x^4$   & \TTTT  &  \TTTT & $36$ \\ 
 $9 + 6x + 5x^2 + 2x^3 + x^4$ &  $ 2x^2 + 2x^3 + x^4$   & \FFFF  &  \TTTT & $24$ \\
\midrule 
 $9 - 6x + 2x^2 - 2x^3 + x^4$ &  $ 2x^2 + x^3 + x^4$  & \TTTT  &  \TTTT & $36$ \\ 
 $9 - 6x + 5x^2 - 2x^3 + x^4$ &  $ 2x^2 + x^3 + x^4$   & \FFFF   &  \TTTT & $24$ \\ 
 $9 + 3x - x^2 + x^3 + x^4$ &  $ 2x^2 + x^3 + x^4$   & \TTTT  &  \TTTT & $24$ \\ 
 $9 + 3x + 2x^2 + x^3 + x^4$ &  $ 2x^2 + x^3 + x^4$   & \TTTT  &  \TTTT & $48$ \\ 
 $9 + 3x + 5x^2 + x^3 + x^4$ &  $ 2x^2 + x^3 + x^4$   & \TTTT  &  \TTTT & $24$  \\ 
 $9 + 12x + 8x^2 + 4x^3 + x^4$ &  $ 2x^2 + x^3 + x^4$   & \FFFF  &  \TTTT & $6$ \\ 
\bottomrule
\end{tabular}
\caption{Data from Lemma~\ref{compute}} \label{table1}
\end{table}

Given a finite flat ordinary mod~$3$ representation~$\rhobar^{\vee}: G_{\Q_3}
\rightarrow \GSp_4(\F_3)$, we would like to realize it as the~$3$-torsion
in the Jacobian of a genus~$2$ curve with good ordinary reduction and a rational Weierstrass point. We shall do this
(under some restrictions)
in Lemma~\ref{belowthree}, using the following lemma.

\begin{lem}  \label{lem: lifting genus 2 curves and finite flat group
    schemes}Let~$p>2$ be prime, and let~$\cO$ be the ring of integers
  in a finite extension of~$\Qp$ with residue
  field~$\F$.  Let~$G_1/\cO$ be a principally quasi-polarized 
  finite flat group scheme of rank~$4$.
Suppose that~$A_0/\F$ is a principally polarized  abelian surface
  with~$A_0[p] \simeq G_{1,\F}$ compatibly with the quasi-polarization.
Then there exists a lift  of~$A_0$ 
 to  a principally
polarized abelian surface ~$A/\cO$
with ~$A[p] \simeq G_1$.
\end{lem}
\begin{proof}
  By~\cite[(2.17)]{MR1827029}, there
  is a lift of~$A_{0}[p^{\infty}]$ to a principally quasi-polarized $p$-divisible group
  $G/\cO$ with~$G[p]\simeq G_1$. By Serre--Tate
  theory~\cite[Thm~1.2.1]{localmoduli}, there exists a principally
  polarized formal abelian surface~$A/\cO$ with~$A_{\F} = A_0$
  and~$A[p^{\infty}] \simeq G$. Since deformations of polarized
  abelian surfaces are effective, we are done.
\end{proof}%

\begin{rem}
  \label{rem: could lift to Jacobians of curves if we needed.}While we
  do not use this fact, we note that if
  $A_{0}$ in Lemma~\ref{lem: lifting genus 2 curves and finite flat group
    schemes} is of the form~$\Jac(C_0)$ for
  a smooth genus~$2$ curve~$C_0/\F$, then necessarily~$A=\Jac(C)$ for a lift~$C$ of~$C_{0}$
  to~$\cO$. To see this, note that since deformations of curves are effective, it suffices
  to show that taking the functor taking a formal lift of~$C_{0}$ to
  its Jacobian is an isomorphism to the deformation problem of
  lifting~$J_0=\mathrm{Jac}(C_0)$ to a principally polarized abelian variety. Since both
  deformation problems are formally smooth of dimension~$3$, it is enough to
  show that the morphism on tangent spaces is injective. This is
  classical; see \cite[\S 2.1]{MR4241766} for an exposition.
\end{rem}

\begin{cor} \label{findlift} Let
$\rhobar^{\vee}: G_{\Q_3} \rightarrow \GSp_4(\F_3)$
be a finite flat representation with similitude
character~$\varepsilonbar$. Assume that~$\rhobar^{\vee}$ is ordinary,
so it 
  is an extension of an unramified
 2-dimensional representation~$\Vbar$ by its Cartier dual.

Suppose that there exists an ordinary principally polarized abelian surface~$A_0/\F_3$ 
with the following properties:
\begin{enumerate}
\item There is an isomorphism of~$G_{\F_3}$-representations~$\Vbar \simeq
 A_0[3]^{\textup{\'{e}t}}$, and
\item The image of~$\rhobar_{A_0,2}$ is conjugate to a subgroup
of~$S_5(b)$.
\end{enumerate}
Then there exists a genus~$2$ curve~$X/\Q_3$ with a~$\Q_3$-rational
Weierstrass point such that~$J = \Jac(X)$ has good ordinary reduction,
and~$\rhobar_{J,3} \simeq \rhobar$.
Moreover, if~$A_0$ is~$3$-distinguished, then so is~$J$.%
\end{cor}

\begin{proof}
  The~$p$-divisible group~$A_0[3^{\infty}]$ is the direct product of an
  \'{e}tale part~$V$ and its Cartier dual. By abuse of notation, we
  may also consider~$V$ as an unramified representation of~$G_{\Q_3}$
  (equivalently, the generic fibre of an \'{e}tale~$3$-divisible group).
We are assuming 
  that~$\Vbar$ is isomorphic to the unramified quotient of~$\rhobar$.

  We can therefore apply Lemma~\ref{lem: lifting genus 2 curves and
    finite flat group schemes} to~$A_0$ where~$G_0$ taken to be
  the (unique) finite flat group scheme with generic
  fibre~$\rhobar^{\vee}$. Let~$A$ be the resulting lift of~$A_0$. 
  Since~$\rhobar_{A_0,2}$ has image inside a conjugate of~$S_5(b)$,
  so does~$\rhobar_{A,2}$ (since~$A$ has good reduction, these two representations are the same). Hence~$A$ gives a~$\Q_3$-rational
  point of~$P(\rhobar)$ (see Definition~\ref{MandP}). 
 By  a version of Krasner's Lemma due to Kisin ~\cite[Thm.\ 5.1]{kisin-krasner},
all the properties listed hold in any open ball around~$A \in
P(\rhobar)$, 
and hence there exists a~$\Q_3$-point~$J = \Jac(X)$
in the corresponding (dense) open subscheme~$\MMM^w_2(\rhobar)$,  and we are done.
\end{proof}

We use this to deduce the following:

\begin{lemma} \label{belowthree}
Let~$\rhobar: G_{\Q_3} \rightarrow \GSp_4(\F_3)$ be an ordinary
representation with similitude factor~$\varepsilonbar^{-1}$,
and suppose that $\rhobar^{\vee}$ is finite flat. 
Then there exists a genus two curve~$X/\Q_3$ with a rational
Weierstrass point such that%
~$\rhobar_{\Jac(X),3} \cong \rhobar$,
and~$\Jac(X)$ has good ordinary reduction and is~$3$-distinguished.
\end{lemma}

\begin{proof}Write~$\overline{q}(x)\in\F_3[x]$ for the characteristic
  polynomial of the Frobenius on the  unramified
  2-dimensional quotient of~$\rhobar^{\vee}$.  Note that~$\overline{q}(x)$
  determines the unramified quotient~$\Vbar$ of~$\rhobar$
  unless it has repeated roots, in which case there are two possible~$\Vbar$;
  one semi-simple and one non-semi-simple.
By Corollary~\ref{findlift}, it suffices to  find for each such~$\Vbar$ an~$A_0/\F_3$
satisfying the hypotheses of Corollary~\ref{findlift}.
    We first consider Jacobians~$A_0 = \Jac(X_0)$ of smooth
  ordinary genus~$2$ curves~$X_0/\F_3$ with an~$\F_3$-rational
  Weierstrass point
   such that the characteristic polynomial~$Q(x)$ of
  Frobenius at~$3$ on~$T_{3}X_{0}$ lifts~$x^{2}\overline{q}(x)$.

There are six possibilities~$x^2 \pm 1$ and~$x^2 \pm x \pm 1$ for
$\overline{q}(x)$, and the existence of such an~$X_0$
follows immediately from Lemma~\ref{compute}, in particular from Table~\ref{table1}.
We can also give explicit examples of such curves~$X_0: y^2 = f(x)$ as follows,
noting that 
(after taking into account unramified quadratic twists) 
we only need to consider four of the six cases.
\begin{center}
\begin{tabular}{*3c}
\toprule
$Q(x)$ & $\overline{q}(x)$ & $f(x)$  \\
\midrule
$x^4 + 3 x^3 + 7 x^2 + 9 x + 9$ & $x^2 + 1$ &    $x^5+2x+1 $ \\
$x^4 + 3 x^3 + 5 x^2 + 9 x + 9$  & $x^2 - 1$ &    $x^5+x^3+x+1$ \\
 $x^4 - x^3 + 2 x^2 - 3 x + 9$ & $x^2 - x - 1$ &  $x^5+x^2+x$ \\
 $x^4 +  x^3 +  x^2 + 3 x + 9$  & $x^2 + x + 1$ &   $x^5 + x^4 + x^2 + 1$  \\
\bottomrule
 \end{tabular}
 \end{center}
 In the ambiguous case where~$\overline{q}(x)$ has repeated roots, 
 it follows from Lemma~\ref{compute}(\ref{troll}) that in all examples
 which arise (including the final example above) the representation~$\Vbar$
 is not semi-simple.
 Hence it remains to consider the case when~$\Frob_3$ acts on~$\Vbar$ by a scalar.
  In this case, we shall construct~$A_0$ directly.
After an unramified quadratic twist (if necessary), we may assume that~$\Vbar$ is trivial.
Let~$E_{-1}/\F_3$ and~$E_{2}/\F_3$ denote the elliptic
curves with~$a_3 = -1$ and~$a_3 = 2$ respectively. Note that they are both ordinary and  they each
have a rational point over~$\F_3$.
Let~$A_0 = E_{-1} \times E_{2}$. Then~$A_0/\F_3$ is principally polarized and~$3$-distinguished,
since~$Q(x) = (x^2 + x + 3)(x^2 - 2 x + 3)$. 
Moreover, $\rhobar_{E_{-1},2}(\Frob_3)$ has order~$2$ and~$\rhobar_{E_{2},2}(\Frob_3)$ has order~$3$.
It follows that~$\rhobar_{A,2}(\Frob_3)$ has order~$6$ and characteristic polynomial~$(x-1)^2(x^2+x+1) \bmod 2$.
This uniquely identifies the conjugacy class as the element~$({*}{*}{*})({*}{*}) \in S_5(b) \subset \GSp_4(\F_2)$, since
the characteristic polynomial of the other conjugacy class of order six elements~$({*}{*}{*}{*}{*}{*}) \in S_6 \simeq \GSp_4(\F_2)$ is equal to~$(x^2+x+1)^2$.

One can verify these claims directly using Lemma~\ref{explicit}, where, for example,
$$(12)(345) \mapsto
 \left( \begin{matrix} 1 & 0 & 0 & 1 \\ 0 & 0 & 1 & 0 \\ 0 & 1 & 1 & 0 \\ 0 & 0 & 0 & 1
\end{matrix}\right), \quad (125346) \mapsto
\left( \begin{matrix} 0 & 0 & 0 & 1\\ 0 & 0 & 1 & 0 \\ 0 & 1 & 0 & 1 \\ 1 & 0 & 1 &0
\end{matrix}\right),$$
 Alternatively,
 the claim about conjugacy classes is equivalent to the claim that the eigenvalues of the semi-simple element~$({*}{*}{*}) \in A_5(b)$ are~$1$, $1$, $\omega$, $\omega^{-1}$ for a primitive third root of unity~$\omega$ and not~$\omega$ and~$\omega^{-1}$
 with multiplicity two;  equivalently that the Brauer character of~$V$ on~$({*}{*}{*})$
evaluates to~$1+1+\omega+\omega^{-1} = 1$ rather than~$2\omega +2\omega^{-1} = -2$,
and this follows from Lemma~\ref{irreps modulo 2}. (The~$4$-dimensional
representation of~$A_5(a)$ coming from~$\GSp_4(\F_2)$, in contrast,  is isomorphic over~$\F_4$
to~$U \oplus U^{\sigma}$.) See also~\cite[Lemma~5.1.7]{BPVY}.\end{proof}

\begin{rem} Although~$A_0/\F_3 =  E_{-1} \times E_{2}$ is  principally polarized, it is not a Jacobian
of an ordinary curve~$X/\F_3$ (by Lemma~\ref{compute}~\eqref{troll}), and so the~$X$ whose existence
is proven in  Lemma~\ref{belowthree} must have bad reduction at~$3$, even though the Jacobian of~$X$
has good reduction at~$3$.
From Table~\ref{table1}, we see there are exactly five %
isogeny classes of principally polarized ordinary abelian surfaces~$A_0/\F_3$ which are~$3$-distinguished
and with~$Q(x) \equiv x^2(x^2+x+1) = x^2+x^3+ x^4 \bmod 3$. 
It turns out that in four out of these five examples, it is \emph{not} possible
to find an~$A_0/\F_3$ in the corresponding isogeny class with~$A_0[3](\F_3) = (\Z/3\Z)^2$.
This can be proved by an argument similar to Lemma~\ref{6Hgood}; 
for each of the five isogeny classes there exists a Jacobian~$B_0/\F_3$ with~$B_0[3](\F_3) = \Z/3\Z$.
Suppose there exists an isogeny~$\chi: B_0 \rightarrow A_0$ with~$A_0[3](\F_3) = (\Z/3\Z)^2$. The kernel of~$\chi$
must contain~$B_0[3](\F_3)$.  Now suppose that  the characteristic polynomial~$Q(x)$
of Frobenius satisfies~$Q(1) \equiv \pm 3 \bmod 9$, which occurs
 in precisely four of these cases.
It follows that (up to isogenies of
degree prime to~$3$) the map~$\chi$ will factor through~$1 - \phi$ where~$\phi$ is the Frobenius morphism,
and since this reduces the power of three dividing the degree,  we
reduce to the case when~$\chi$ has degree prime to three and we obtain a contradiction.
In the remaining case (which we exploited above), we have $Q(x) = (3-x+x^2)(3+2x+x^2)$ and so~$Q(1)  \equiv 0 \bmod 9$,
and now such an isogeny is possible.
\end{rem}

\subsection{A 2-3 switch}\label{subsec: 2 3 switch}

We begin with the following approximation lemma.

\begin{lemma} \label{rationalmb} Let~$Z$ be a rational variety
  over~$\Q$. Let ~$S$ be a finite set of places of~$\Q$, and for each
  $v\in S$, let~$\Omega_v$ be a non-empty open subset of~$Z(\Q_v)$
  \emph{(}for the $v$-adic topology\emph{)}.  Then there exists a rational
  point~$P \in Z(\Q)$ with~$P_v \in \Omega_v$ for all~$v \in S$, and
  such that~$P$ avoids any fixed thin subset of~$Z(\Q)$.
\end{lemma}

\begin{proof}Apart from the statement that we may avoid any fixed thin
  subset of~$Z(\Q)$, this is a special case of~\cite[Lem.\
  3.5.5]{MR2363329}, and our proof is an obvious variation on the
  arguments of \cite[\S 3.4, \S 3.5]{MR2363329}.  We may assume that~$S$ is nonempty. After shrinking~$Z$ if necessary, we may assume
 that~$Z \hookrightarrow \mathbf{P}^n_{\Q}$
is an open immersion. Here we use that~$Z$ is smooth; 
this guarantees that, for any open~$U \subset Z$,
$\Omega_v \cap U(\Q_v) \subset U(\Q_v)$ is non-empty.
Since~$Y = \mathbf{P}^n_{\Q} \setminus Z$ is
closed, $Y(\Q_v) \subset \mathbf{P}^n_{\Q}(\Q_v)$ is closed, and
so~$\Omega_v \subset Z(\Q_v) \subset \mathbf{P}^n_{\Q}(\Q_v)$ is
open. Since
$$\mathbf{P}^n_{\Q}(\Q) \cap \Omega_v \subset \mathbf{P}^n_{\Q}(\Q) \cap Z(\Q_v) = Z(\Q),$$
 we can and do assume that $Z=\mathbf{P}^n_{\Q}$. %

The number of points in~$\mathbf{P}^n_{\Q}(\Q)$
which are of height at most~$H$ and are contained in~$\Omega_v$ for
all~$v\in S$  grows at the rate of a positive constant
 times~$H^n$ (the precise constant depending on the open
 sets~$\{\Omega_v\}$), whereas the number of points in any fixed thin set is bounded
by~$O(H^{n-1/2} \log H)$ by~\cite[Thm.\ 3.4.4]{MR2363329}, and the
result follows.  %
\end{proof}

We now construct a suitable abelian surface through which to do our~$2$-$3$-switch.

\begin{lemma}  \label{switching}
Suppose that 
$$\rhobar: G_{\Q} \rightarrow \GSp_4(\F_3)$$ has
similitude~$\varepsilonbar^{-1}$, that $\rhobar^{\vee}|_{G_{\Q_{3}}}$ is
ordinary and finite flat, and that $\rhobar|_{G_{\Q_2}}$ is unramified.
\begin{enumerate} %
\item The following conditions are equivalent: \label{cattwoswitch}
\begin{enumerate}
\item \label{frob} The image of~$\rhobar(\Frob_2)$ in~$\PGSp_4(\F_3)
  \setminus \PSp_4(\F_3)$ is not conjugate to~$4C$ or~$12C$
\emph{(}see
Lemma~\ref{atlas}\emph{)}.
\item \label{curve} $\rhobar |_{G_{\Q_2}} \simeq \rhobar_{A,3}$, where~$A$ is the Jacobian of a genus~$2$ curve~$Y/\Q_2$ with 
a rational Weierstrass point, and where~$A$ has either good ordinary
 or semistable ordinary reduction at~$2$ and is~$2$-distinguished.
\end{enumerate}
\item\label{item: actually building a curve over Q} Assume that the equivalent
  conditions in~\eqref{cattwoswitch} hold.
Then there exists a genus two curve~$X/\Q$ with a rational Weierstrass point,
with~$B = \mathrm{Jac}(X)$ having the following properties:
\begin{enumerate}
\item $\rhobar_{B,3}\cong\rhobar$. %
\item $B$ has good ordinary  or semistable ordinary reduction at~$2$,
and is~$2$-distinguished.
\item $B$ has good ordinary reduction at~$3$. %
\item The representation
$$\rhobar_{B,2}: G_{\Q} \rightarrow \GSp_4(\F_2)$$
has image $S_5(b)$, and the image of complex conjugation has
conjugacy class~$({*}{*})({*}{*})$.
\end{enumerate}
Moreover, $\End(B_{\Qbar}) = \Z$.
\end{enumerate}
\end{lemma}
\begin{proof}
We recall from Definition~\ref{MandP}; that~$\MMM^w_2(\rhobar)$ and~$P
= P(\rhobar)$ are the fine moduli spaces over~$\Q$ parametrizing respectively
genus 2 curves~$X$ with a rational Weierstrass point together with a symplectic
isomorphism~$\rhobar_{\Jac(X),3}\cong\rhobar$, and principally
polarized abelian surfaces with a fixed odd theta characteristic and a
symplectic  isomorphism~$\rhobar_{A,3} \cong \rhobar$.%

We claim that condition~(\ref{frob}) is equivalent to condition~(\ref{curve}) by 
Theorem~\ref{restrictionsattwo}. More precisely, that theorem shows that
~(\ref{frob}) implies there exists a point~$A \in P(\rhobar)(\Q_2)$
with either good ordinary or semistable ordinary reduction (and which is~$2$-distinguished), whereas
condition~(\ref{curve}) shows that there is a point~$A = \Jac(Y) \in P(\rhobar)(\Q_2)$
which lies in the image of~$\MMM^w_2(\rhobar)$. The variety~$P(\rhobar)$ is smooth
and the map~$\MMM^w_2(\rhobar) \rightarrow P(\rhobar)$ is an open immersion.
Moreover, all the properties listed hold in any open ball around any such point~$A$
by~\cite[Thm.\ 5.1]{kisin-krasner}. Hence given~$A \in P(\rhobar)(\Q_2)$
there exists a point~$B \in P(\rhobar)(\Q_2)$ with the same properties but lying
in the image of~$\MMM^w_2(\rhobar)$.

Having established this equivalence, we now turn to the proof of
part~\eqref{item: actually building a curve over Q}, so we in
particular assume that condition~(\ref{curve}) holds. We now use Lemma~\ref{rationalmb} (applied to~$Z = \MMM^w_2(\rhobar)$)
to produce a suitable point~$X/\Q$.
Our set~$S$ will consist of the primes~$2$, $3$, $\infty$.
The corresponding thin set inside~$Z(\Q) \subset P(\Q)$ is the union of the rational points
in the images of~$P_G(\Q)$,  where~$P_G \rightarrow P$ is the cover
corresponding to imposing that the image of~$\rhobar_{B,2}$
lands inside a strict subgroup~$G \subset S_5(b)$. There are finitely
many such~$G$ and the degree of~$P_G$ over~$P$ is~$[S_5(b):G] > 1$, so
this is indeed a thin set.

\begin{enumerate}
\item Suppose that~$p=2$. Condition~(\ref{curve})
implies that there exists a point~$X$ in~$Z(\Q_2)$ with
the required properties ($\Jac(X)/\Q_2$ with good ordinary reduction or
semistable ordinary reduction and~$2$-distinguished
in characteristic zero).
By~\cite[Thm.\ 5.1]{kisin-krasner},
there exists an open ball~$\Omega_2 \subset P(\Q_2)$
around~$X$ consisting of points which also have good
ordinary reduction and are~$2$-distinguished.
\item 
Suppose that~$p=3$. 
Then there exists a suitable point~$X \in Z(\Q_3)$  by
Lemma~\ref{belowthree}. 
As above, we take~$\Omega_3$ to be a suitable open ball around~$X$.
\item For~$p  = \infty$, we choose~$\Omega_{\infty} \subset Z(\R)$
to be a sufficiently small open ball around any point with the correct local properties,
namely~$y^2 = f(x)$ for any separable  $f(x)\in\R[x]$ of degree~$5$
with exactly one real root.
\end{enumerate}

The existence of~$X$ and~$B$ then follows from Lemma~\ref{rationalmb}.
Since the image of~$\rhobar_{B,2}$ is~$S_5(b)$, it follows from~\cite{MR1748293}
that~$\End(B_{\Qbar}) = \Z$.
\end{proof}
\begin{rem}
  \label{rem:don't-care-3-distinguished}The application of
  Lemma~\ref{belowthree} in the proof of Lemma~\ref{switching} can obviously
  additionally be used to show that~$B$ can be chosen to be~$3$-distinguished,
  but we shall not use this fact, so we have not explicitly recorded it.
\end{rem}
\begin{rem} Suppose that~$\rhobar: G_{\Q} \rightarrow \GSp_4(\F_3)$ has
  multiplier~$\varepsilon^{-1}$, and $\rhobar|_{G_{\Q_3}}$ is 
ordinary and (dual to) finite flat.
Then, exactly as in the proof of
 Lemma~\ref{switching} (now ignoring
the conditions at~$2$) obtains infinitely many
genus two curves~$X/\Q$ with a rational Weierstrass point and such that~$A = \Jac(X)$ has good ordinary
reduction at~$3$, such that $\rhobar_{X,3}\cong\rhobar$.
This was implicitly assumed in the proof of~\cite[Theorem~10.2.1]{BCGP}.
\end{rem}

\subsection{Proof of Theorems~\ref{first} and~\ref{second}}
\label{sec:proofsofAandB}

In this section, we prove Theorem~\ref{first}, which we restate as Theorem~\ref{firstlater} below, except that
the hypothesis on the image of~$\rhobar_{A,3}$ has been relaxed. (Note
that if~$\rhobar_{A,3}$ is surjective, then~$\End(A_{\Qbar}) = \Z$ is
automatic, so Theorem~\ref{firstlater} really does imply
Theorem~\ref{first}.)
We begin, however, with the following %
modularity lifting theorem.

\begin{thm}%
  \label{thm: p>2 modularity from one abelian surface to
    another}Suppose that $p>2$, and that~$A/\Q$ and~$B/\Q$ are abelian surfaces such
  that:
  \begin{enumerate}
    \item \label{conditionequalimage} $\rhobar_{A,p}\cong \rhobar_{B,p}$.
  \item \label{conditiondistinguished} $A$ and $B$ both have good ordinary reduction at~$p$, and
    $\rho_{A,p}|_{G_{\Qp}}$ %
    is $p$-distinguished.
    \item \label{conditionmodular} $B$ is modular; more precisely, there is a
    weight~$2$ cuspidal automorphic representation~$\pi$  for~$\GSp_4 /\Q$ of
    level prime to~$p$,
    which is ordinary at~$p$ and satisfies $\rho_{\pi,p}\cong\rho_{B,p}$.
   \item \label{conditionzariski} The Zariski closure of $\rho_{A,p}(G_\Q)$
    contains~$\Sp_4$. %
  \item \label{conditionreasonable} $\rhobar_{A,p}$ is $\GSp_4$-reasonable, in the sense
          of~\cite[Defn.\
          3.19]{Whitmore}. 
        \item \label{conditiontidy} $\rhobar_{A,p}$ is tidy, in the sense of~\cite[Defn.\
           7.5.11]{BCGP}. %
         \item \label{conditionregularsemi-simple} $\rhobar_{A,p}(G_{\Q(\zeta_p)})$ contains a regular semi-simple
           element. %
           \item\label{condition-regular-semi-simple-central-character}   $\rhobar_{A,p}(G_{\Q})\setminus \Sp_4 (\F_p)$ contains a regular semi-simple element.
         \end{enumerate}
         Then~$A$ is modular.
         More precisely, there exists a cuspidal automorphic representation~$\pi$
for~$\GL_4/\Q$ \emph{(}the transfer of a cuspidal
automorphic representation of~$\GSp_4/\Q$ of weight~$2$\emph{)} %
such that~$L(s,H^1(A),s) = L(s,\pi)$.
\end{thm}
\begin{proof}
  We deduce the theorem from
  Theorem~\ref{thm:rho-is-modular-from-mult-one-and-classicity} (taking~$\rho$
  there to be~$\rho_{A,p}$). By
  our assumption~\eqref{conditionzariski}, it suffices to check that hypotheses (1)--(5)
  and \ref{item: this needs a label}--\ref{item: same component} of
  Proposition~\ref{prop: the modularity result for p equals 2 or 3} hold.

  Most of these conditions hold
    either explicitly by our assumptions, or by the purity of Galois representations associated to abelian surfaces.
    The only remaining conditions are:
     \begin{enumerate}[(a)]
  \item \label{item:enormous in proof of B to A}
    $\rho(G_{\Q(\zeta_{p^\infty})})$ is integrally enormous.
           \item\label{item:same component in proof of B to A}We can choose
             $p$-stabilizations of $\rho_{A,p}|_{G_{\Qp}}$ and~$\rho_{B,p}|_{G_{\Qp}}$ such
             that the representations ~$\rho_{B,p}|_{G_{\Qp}}$ lies on a unique irreducible component of $\Spec R_p^\triangle$ and $\rho_{A,p}|_{G_{\Q_p}}$ lies on the same component.%
         \end{enumerate}
Part~\ref{item:enormous in proof of B to A}
follows from Corollary~\ref{cor: still integrally enormous up the pro p tower}.
For Part~\ref{item:same component in proof of B to A}, we firstly choose a
$p$-stabilization of~$\rho_{B,p}|_{G_{\Qp}}$, and thus
of~$\rhobar_{B,p}|_{G_{\Qp}}=\rhobar_{A,p}|_{G_{\Qp}}$. Then at least one of the
 $p$-stabilizations of~$\rho_{A,p}|_{G_{\Qp}}$ is compatible with this fixed
 choice, and we conclude by Lemma~\ref{lem: good ordinary reduction points are on a single
    component} and the assumption that~$A$ and~$B$ both have good
  ordinary reduction at~$p$.
   \end{proof}

We are now ready to prove our main theorem.

\begin{thm} \label{firstlater}
Let~$A/\Q$ be an abelian surface with a polarization of degree prime to~$3$. Suppose that
the following conditions hold:
\begin{enumerate}
\item \label{surjectivelater} The image of the mod~$3$ representation:
$$\rhobar_{A,3}: \Gal(\Qbar/\Q) \rightarrow \GSp_4(\F_3)$$
is one of the~$15$ subgroups listed in Lemma~\ref{listofsubgroups}, and ~$\End(A_{\Qbar}) = \Z$.
\item  \label{conditionattwolater} $\rhobar_{A,3}|_{{G_{\Q_{2}}}}$ is unramified, and 
the image 
of~$\rhobar_{A,3}(\Frob_2)$ inside $\PGSp_4(\F_3) \setminus
\PSp_4(\F_3)$  does not have conjugacy class~$4C$ or~$12C$
\emph{(}see
Lemma~\ref{atlas}\emph{)}. Equivalently, the characteristic polynomial 
of $\rhobar_{A,3}(\Frob_2)$ is not $(x^2\pm x+2)^2$.
\item \label{conditionatthreelater}
 $A$ has good ordinary reduction at~$3$ and is $3$-distinguished.
\end{enumerate}
Then~$A$ is modular. More precisely, there exists a cuspidal automorphic representation~$\pi$
for~$\GL_4/\Q$ \emph{(}the transfer of a cuspidal
automorphic representation of~$\GSp_4/\Q$ of weight~$2$\emph{)} %
such that~$L(s,H^1(A)) = L(s,\pi)$.
\end{thm}
\begin{proof}By Lemma~\ref{switching} \eqref{item: actually building a
    curve over Q} (which applies to~$\rhobar_{A,3}$, since
  condition~(\ref{frob}) of  Lemma~\ref{switching}  holds by our assumption~\eqref{conditionattwolater}),
 there exists a genus two curve~$X/\Q$ with a rational Weierstrass point,
with~$B = \mathrm{Jac}(X)$ having the following properties:
\begin{itemize}
\item $\rhobar_{A,3}\cong\rhobar_{B,3}$. %
\item $B$ has semistable ordinary or good ordinary reduction at~$2$,
and is~$2$-distinguished.
\item $B$ has  good ordinary  reduction at~$3$. %
\item The representation
$$\rhobar_{B,2}: G_{\Q} \rightarrow \GSp_4(\F_2)$$
has image ~$S_5(b)$, and the image of complex conjugation has
conjugacy class~$({*}{*})({*}{*})$.%
\item $\End(B_{\Qbar}) = \Z$.
\end{itemize}

We shall first apply Theorem~\ref{thm: residually A5 implies modular} (a
$2$-adic modularity theorem; we take
~$A$ there to be our~$B$) to deduce that~$B$ is modular. To recall, the hypotheses of Theorem~\ref{thm: residually A5 implies modular} 
are as follows:
  \begin{enumerate}[(i)]
  \item \label{imageagain} $A_5(b)\subseteq \rhobar_{B,2}(G_\Q)\subseteq S_5(b)$.
  \item  \label{infinityagain} The image of each complex conjugation has order~$2$ and lands
  in~$A_5(b)$.
  \item  \label{distinguishedagain} $\rho_{B,2}|_{G_{\Q_2}}$ is ordinary and $2$-distinguished.
  \end{enumerate}
All of these conditions are guaranteed by the properties of~$B$ listed
above, noting that %
  ~$({*}{*})({*}{*})$
  is a non-trivial conjugacy class contained in~$A_5(b)$. %
  Thus~$B$ is modular.  More precisely, there is a
    weight~$2$ cuspidal automorphic representation~$\pi$  for~$\GSp_4 /\Q$
    which in particular satisfies $\rho_{\pi,3}\cong\rho_{B,3}$;
    furthermore~$\pi$ is necessarily of level prime to~$3$
    and is ordinary at~$3$ by local-global compatibility.%
We now wish to use Theorem~\ref{thm: p>2 modularity from one abelian surface to
    another} at~$p=3$ to deduce that~$A$ is modular, 
    so we need to check the conditions
    of that theorem. %
We established that~$B$ is modular above, and
thus condition~\eqref{conditionmodular} holds. The isomorphism $\rhobar_{A,3}\cong\rhobar_{B,3}$ 
holds by the construction of~$B$, hence we have
condition~\eqref{conditionequalimage}. Both $A$ and~$B$ have good ordinary
reduction at~$3$ and $A$ is furthermore~$3$-distinguished
 (by assumption for~$A$
  and by construction for~$B$), and thus we have
  condition~\eqref{conditiondistinguished}. We are assuming that~$\End(A_{\Qbar}) = \Z$,
so condition~\eqref{conditionzariski} holds
by~\cite[Thm~3]{MR1730973}. Finally conditions~\eqref{conditionreasonable},
\eqref{conditiontidy}, \eqref{conditionregularsemi-simple} and \eqref{condition-regular-semi-simple-central-character} 
hold by Lemma~\ref{listofsubgroups} and our assumptions on~$A$.
\end{proof}

\begin{thm} \label{thirdlater}
Let~$X$ be a smooth genus two curve over~$\Q$.
Suppose that:
\begin{enumerate}
\item The image of~$\rhobar_{\Jac(X),3}: G_{\Q} \rightarrow \GSp_4(\F_3)$ is 
 one of the~$15$ subgroups listed in Lemma~\ref{listofsubgroups},
and~$\End(\Jac(X)_{\Qbar}) = \Z$.
\item $X$ has good ordinary reduction at~$2$.
\item $X$ has good ordinary reduction at~$3$.
\item $\Jac(X)$ is~$3$-distinguished.
\end{enumerate}
Then~$X$ is modular.
\end{thm}
\begin{proof}
Let~$A = \Jac(X)$ (so that~$A$ is in particular principally
polarized).  It suffices to verify the conditions \eqref{surjectivelater}--\eqref{conditionatthreelater} of Theorem~\ref{firstlater}.
Condition~\eqref{surjectivelater} is identical to our first condition. Since
we are assuming that~ $X$  has good ordinary
reduction at~$2$ and~$3$, so does~$A$.
Condition~\eqref{conditionattwolater} follows from Lemma~\ref{exhaust}\eqref{exhaust implies weaker condition at two}.
Finally condition~\eqref{conditionatthreelater} is immediate from our assumptions.
 Hence~$A$ (and thus~$X$) is modular.\end{proof}

 We now deduce Theorem~\ref{second}, which we restate here, again with
 a weakening of the assumption that~$\rhobar_{X,3}$ is surjective.
\begin{thm} \label{secondlater}
Let~$X: y^2 = f(x)$ with~$\deg(f(x)) = 5$ be a smooth genus two curve over~$\Q$.
Suppose that:
\begin{enumerate}
\item The image of~$\rhobar_{X,3}: G_{\Q} \rightarrow \GSp_4(\F_3)$ is 
 one of the~$15$ subgroups listed in Lemma~\ref{listofsubgroups},
and~$\End(\Jac(X)_{\Qbar}) = \Z$.
\item $X$ has good ordinary reduction at~$2$.
\item $X$ has good ordinary reduction at~$3$.
\end{enumerate}
Then~$X$ is modular.
\end{thm}
\begin{proof} %
By Theorem~\ref{thirdlater}, it suffices to show that~$\Jac(X)$ is $3$-distinguished.
Since~$X$ has a rational Weierstrass point (by our assumption that
$f(x)$ has degree~$5$), this follows immediately from Lemma~\ref{compute}\eqref{troll}. %
\end{proof}

Note that Theorem~\ref{first} and Theorem~\ref{firstlater} do not require that~$A$ has good reduction at~$2$,
only that~$\rhobar_{A,3}$ is unramified at~$2$.
Here we answer a question of Drew Sutherland, who asks if the conditions
of our main theorem are easy to verify computationally if~$\Q(A[3])$ is unramified
at~$2$ but~$A$ has \emph{bad} reduction at~$2$. It turns out
that the answer is surprisingly simple.

\begin{theorem} \label{drew}
Let~$A/\Q$ be an abelian surface with a polarization of degree prime to~$3$. Suppose 
the following  holds:
\begin{enumerate}
\item  The image of the mod~$3$ representation:
$$\rhobar_{A,3}: \Gal(\Qbar/\Q) \rightarrow \GSp_4(\F_3)$$
is one of the~$15$ subgroups listed in Lemma~\ref{listofsubgroups}, and ~$\End(A_{\Qbar}) = \Z$.
\item  $\rhobar_{A,3}|_{G_{{\Q_2}}}$ is unramified.
\item  $A$ has good ordinary reduction at~$3$ and
 the characteristic polynomial
of  Frobenius at~$3$ does not have repeated roots.
\item $A$ has bad reduction at~$2$.
\end{enumerate}
Then~$A$ is modular. \end{theorem}

\begin{proof}
We shall apply Theorem~\ref{firstlater}.
It suffices to show that, under the assumption that~$A$ has \emph{bad}
reduction at~$2$, that the action of~$\Frob_2$ on~$A[3]$
does not have characteristic polynomial~$(x^2 \pm x + 2)^2$.
By Lemma~\ref{groth}, we deduce that~$A$ has semistable reduction.
Hence, as in the proof of Lemma~\ref{boxerstrass}, we
deduce the existence 
of a~$G_{\Q_2}$-equivariant
 filtration~$T_2(B)_t \subset T_2(B)_f \subset T_2(B)$
of (saturated)~$\Z_3$-modules of ranks~$t>0$ and~$t+2a$ where~$2(t+a)=4$. 
If~$t=1$, then~$A[3]$ has a~$G_{\Q_2}$-stable line. But~$x^2 \pm x + 2$
has no eigenvalues in~$\F_3$, which concludes the proof in this case.
Assume that~$t=2$, so~$A$ has purely multiplicative reduction.
It follows that~$A$ has split multiplicative reduction over some minimal unramified extension~$K/\Q_2$.
There is a corresponding action of~$\Gal(K/\Q_2)$ on~$\mathbf{G}_m \times \mathbf{G}_m$ which gives
the descent data to~$\Q_2$; this determines a finite order element of~$\GL_2(\Z)$,
and such elements can only have orders~$1$, $2$, $3$, $4$, or~$6$.
(The characteristic polynomial of this element will be, up to normalization,
 the~$L$-factor of~$A$ at~$p = 2$.)
On the other hand, 
 the action
of~$G_{\Q_2}$ on~$T_3(B) \otimes \Q_3$ factors through~$\Gal(K/\Q_2)$.
By considering the action on the unramified quotient~$A[3]$, we deduce that $\Gal(K/\Q_2)$ has
order divisible by~$8$, since~$x^4 + 1 = (x^2 + x + 2)(x^2 - x + 2) \bmod 3$.
This is a contradiction.
\end{proof}

\begin{remark} An alternative argument is to note that if the characteristic polynomial of~$\Frob_2$
on~$\rhobar_{A,3}$ is~$(x^2 \pm x + 2)^2$, then the ratio of any two eigenvalues is
never equal to~$2 = -1 \bmod 3$, and so~$\rhobar_{A,3}|_{G_{\Q_{2}}}$ has no ramified %
lifts.
\end{remark}

\section{Complements}
\label{sec:complements}

This final section includes a number of results which are complementary to the main theorems  of our paper (and in particular are not used elsewhere).

In~\S\ref{sub:examples}, we explain how our main theorems apply (relative to a certain  natural way of enumerating genus two curves) to slightly over 10\% of all such curves, and we compare this to the data
in the LMDFB~\cite{LMFDB}.
In~\S\ref{nongeneric}, we prove the automorphy  of any abelian surface $A/\Q$ which falls into~$32$ of the~$34$ possible
Galois types.
In~\S\ref{residualmod2}, we prove some residual modularity
theorems (Serre's conjecture) for mod-$2$ representations~$\rhobar: G_{\Q} \rightarrow \GSp_4(\F_2)$ 
with image~$A_6$ or~$S_6$. Finally, in~\S\ref{regularimpliesirregular}, we point out that a sufficiently strong version of Serre's conjecture for~$\GSp_4$ in \emph{regular} weight would be enough to prove the modularity of all abelian surfaces~$A/\Q$.

\subsection{Examples}
\label{sub:examples}

Suppose one samples genus two curves
$$X:y^2 + h(x) y = f(x)$$
with~$h(x), f(x) \in \Z[x]$ of degrees~$\le 3$ and~$\le 6$ in any way in which the distributions
modulo~$2$ and~$3$  are equidistributed, 
and considers curves~$X$ with the following properties:
\begin{enumerate}
\item $X$ has good reduction at~$2$,
\item $X$ has good ordinary reduction at~$3$,
\item $\rhobar_{\Jac(X),3}(\Frob_2)$ does not have characteristic polynomial
$x^4 \pm x^3 + 2 x^2 \pm x + 1$, equivalently, is not projectively conjugate
to~$4C$ or~$12C$,
\item The characteristic polynomial~$Q_3(x)$ of~$\Frob_3$
 has distinct eigenvalues.
 \end{enumerate}
 Then from Lemmas~\ref{exhaust} and~\ref{compute} these~$X$ form a subset of density
$$\frac{768 - 144}{2^{11}} \cdot \frac{864 - 10}{3^7} =\frac{13}{16} \cdot \frac{3}{8} \cdot
\frac{854}{2187} = 
\frac{5551}{46656} = 0.1189\ldots $$
(Note that~$13/16$ is the density of allowable elements for~$\rhobar_{X,3}(\Frob_2)$,
and~$3/8$ is the density of curves with good reduction at~$2$.) Since~$\End(A_{\Qbar}) = \Z$ holds for a set of density one,
we see in particular that Theorem~\ref{firstlater} applies to a positive (although not so large,
slightly over~10\%) proportion
of all genus~$2$ curves. (The main theorem of~\cite{MR1333035} also applies to a positive
but strictly less than one proportion of all genus~$1$ curves by any natural counting.)

Another point of comparison is with the curves in the database~\cite{LMFDB}:
There are~$66158$ genus two curves~$X$ in~\cite{LMFDB,database}:
\begin{enumerate}
\item Of those, $63107$ have~$\End A_{\Qbar} = \Z$,
where~$A = \Jac(X)$. 
\item Of those, $22158$ have good reduction at~$2$ and~$3$. (In the
  range of the data, a genus~$2$ curve~$X$
has good reduction at~$p$ (for any~$p$) if and only if~$A = \Jac(X)$ has good reduction at~$p$.)
\item Of those, $21552$ have surjective mod~$3$ representations.
\item Of those, $14856$ have ordinary reduction at~$p=3$.
\item Of those, our theorem applies to~$11384$ curves, where
the distribution of various conjugacy classes and~$3$-distinguishedness
conditions is indicated in  Table~\ref{table: counting LMFDB examples}.
\end{enumerate}
\numtable[!htbp]\label{table: counting LMFDB examples}
\centering
\caption{$\rhobar_{A,3}(\Frob_2)$ distributions
of $A = \Jac(X)$ where~$X$ has good 
reduction at~$2$, good ordinary reduction at~$3$, $\rhobar_{A,3}$ is
surjective, 
and~$X$ is taken from~\cite{database}, together with a count
of those to whom Theorem~\ref{firstlater} applies.}
\begin{tabular}{*5c}
\toprule
$\rhobar_{A,3}(\Frob_2)$ &  \multicolumn{2}{c}{ordinary at~$2$} & \multicolumn{2}{c}{non-ordinary at~$2$}\\
\midrule
{}   & $3$-dist   & not $3$-dist    & $3$-dist    & not $3$-dist \\
$2C/6G/6H$   &  1048 & 8  & 840 & 7\\
$4D$ & 0 & 0 & 890 & 2 \\
$2D/6I$    &  825  &  1   & 0  & 0\\
$8A$    &  1233  &  9   & 854  & 6\\
$10A$    &  3407  &  48   & 2287  & 4 \\
$4C/12C$ & 0 & 0 & 3369 & 18 \\
\bottomrule
All     &  6513 & 66 &  8240  & 37  \\ 
Theorem~\ref{firstlater} applies   &  6513 & 0 &  4871 & 0  \\ 
\end{tabular}
\end{table}

If one 
allows~$\rhobar_{A,3}$ to be
any of the~$15$ subgroups listed
 in Lemma~\ref{listofsubgroups},
there are three additional curves, precisely one of which
we can deduce is modular by Theorem~\ref{firstlater}.
This is
the curve \texttt{7889.b.55223.1} 
of conductor~$7^3 \cdot 23$. 
The representation~$\rhobar_{A,3}$ in this case (with image of  order~$2304$)
 is induced from a  representation~$\rhobar_{E,3}: G_F \rightarrow \GL_2(\F_3)$,
where~$E$ is a modular elliptic curve over~$F = \Q(\sqrt{-7})$
with~\cite{LMFDB} label \texttt{2.0.7.1-322.1-a1} and
conductor of norm $322 = 2 \cdot 7 \cdot 23$.

\medskip

There are $41324$ curves in the~\cite{LMFDB} database with $\rhobar_{A,3}$ surjective,
 $\End_{\Qbar}(A) = \Z$,  
 and such that~$A$ has good reduction at~$3$.
Of those,  $19772$ of these curves have \emph{bad} reduction at~$2$.
 We find that for precisely~$360$ of these curves, $A$
 is ordinary at~$3$ and $\rhobar_{A,3}$ is unramified at~$2$.
Of those, $359$ are~$3$-distinguished
and thus modular by Theorem~\ref{drew}. 
In particular, Theorem~\ref{first} applies
to precisely~$11384 + 359 = 11743$ of the~$66158$ curves in the~\cite{LMFDB}.
(One can verify modularity for more of the curves in~\cite{LMFDB} by 
including quadratic twists.)
 The smallest conductor of such
an~$A$ with bad reduction at~$2$ is 1982; 
the corresponding $A$ is the Jacobian of the curve
$$y^2 + (x + 1)y = -x^5 + x^4 - x^3 + x^2.$$

\subsection{Automorphy for abelian surfaces with small Sato--Tate group}
\label{nongeneric}

In this section, we prove the automorphy
(in the sense of Definition~\ref{definitionautomorphic})  for~$32$
of the~$34$ Galois types (in the sense of~\cite{MR2982436}) of
abelian surfaces~$A$ over~$\Q$. 
We closely follow~\cite[\S9.2]{BCGP} and use freely the notation of that section,
as well as the results summarized there from~\cite{MR2982436,MR3660222} (see
also~\cite{MR4082247}). We say that a Galois representation is
``finite up to twist'' if it is a twist by a character of a
representation with finite image.

Recall that the Galois type of~$A/\Q$ is~$\mathbf{A}$
precisely when~$\End(A_{\Qbar}) = \Z$,
and~$A/\Q$ has type~$\mathbf{B}[C_2]$ if there
exists a quadratic field~$K/\Q$ so that~$\End(A) = \Z$ 
but $\End(A_K) \otimes \Q$ is either~$\Q \oplus \Q$ or a real
quadratic field. (In~\cite[\S9.2]{BCGP}, we call an abelian
surface~$A/\Q$ ``challenging'' precisely when it is one of these two types.)

The main theorem of this section is as follows.
\begin{theorem} \label{remainingcases} Let~$A/\Q$ be an abelian surface.
Suppose that the Galois type of~$A$ is neither~$\mathbf{A}$ nor~$\mathbf{B}[C_2]$.
Then~$A$ is modular.
\end{theorem}

\begin{rem}[Abelian surfaces of Galois type~$\mathbf{B}[C_2{]}$]
A natural source of abelian surfaces of type~$\mathbf{B}[C_2]$ are those of the form~$\Res_{K/\Q}(E)$
for a non-CM elliptic curve~$E$ which is not isogenous to its~$\Gal(K/\Q)$-conjugate.
If~$K/\Q$ is  real quadratic, then $E$ is automorphic
 for~$\GL_2/K$ by~\cite{FLS} and then~$A$ is automorphic for~$\GL_4/\Q$.
The modularity of elliptic curves~$E$ over imaginary quadratic
fields~$K$ is known in many cases (but not yet all) by~\cite{caraiani-newton}.
On the other hand,  for~$A$ of type~$\mathbf{B}[C_2]$, the endomorphism algebra
 $\End(A_K) \otimes \Q$ could also be a real quadratic field~$E$ rather than~$\Q \times \Q$, in which case~$A/K$
will be a simple abelian surface of~$\GL_2$-type. 
This happens, for example, when~$A$ is the Jacobian of the genus~$2$
 curve
 $$y^2 + (x^3 + 1)y = x^6 + 2x^3 - x$$
with~$E=K=\Q(\sqrt{5})$ 
\cite[genus~$2$ curve \texttt{12500.a.12500.1}]{LMFDB}.
The modularity of such abelian surfaces  remains open in general even
for real quadratic fields~$K$.
\end{rem}

\begin{proof}[Proof of Theorem~\ref{remainingcases}]
Following the discussion in~\cite[\S9.2]{BCGP}  and~\cite[Prop~9.2.1]{BCGP},
all abelian surfaces~$A/\Q$ can be divided up into a number of possible Galois types,
which, writing~$\{\rho_{A,p}\}$ for the compatible system of Galois
representations $\{H^1(A_{\overline{\Q}},\Qbar_p)\}$,
fall into the following categories independently of~$p$:
\begin{enumerate}
\item strongly irreducible (type~$\mathbf{A}$),
\item reducible (type~$\mathbf{B}[C_1]$, $\mathbf{C}$, $\mathbf{E}[C_n]$, some~$\mathbf{D}$, some~$\mathbf{F}$),
\item potentially abelian but not reducible (of type the remaining~$\mathbf{D}$ and~$\mathbf{F}$ cases),
\item induced from a quadratic extension~$K/\Q$ but not potentially abelian, in which case either:
\begin{enumerate}
\item the two~$2$-dimensional representations over~$K$ are equivalent
  up to twist
(type~$\mathbf{E}[D_n]$), or
\item the two~$2$-dimensional representations over~$K$ are not
 equivalent up to twist (type~$\mathbf{B}[C_2]$).
\end{enumerate}
\end{enumerate}

We will prove automorphy in all cases except those of type~$\mathbf{A}$
and those of type~$\mathbf{B}[C_2]$.

In the reducible cases, it follows from~\cite{MR2982436}
that
the compatible system associated to~$A$ can be written as a direct
sum of two irreducible, odd, regular,
weakly compatible systems of
Galois representations over~$\Q$.
 These are modular by~\cite{MR2551763}.
 In case~$\mathbf{E}[D_n]$, we see (as in the proof
 of~\cite[Prop~9.2.1]{BCGP}) that there exists
 a quadratic extension~$K/\Q$
 and an odd irreducible regular weakly compatible system~${\cS} = \{s_p\}$  of~$G_{\Q}$ such that
 $\rho_p \simeq s_p \otimes   \Ind^{G_{\Q}}_{G_K} \psi^{-1} = \Ind^{G_{\Q}}_{G_K} s_p |_{G_K} \otimes \psi^{-1}$
 for some fixed finite order character~$\psi$. Once more~$\cS$ is automorphic for~$\GL_2/\Q$
 by~\cite{MR2551763}, and then~$A$ is is automorphic for~$\GL_4/\Q$,
 as required.

 It remains to consider the cases where~$\rho_p$ is (absolutely) irreducible but potentially
 abelian. Since the representations~$\rho_p$ have similitude character~$\varepsilon^{-1}$ 
 and~$\rho_p$ is not finite up to twist (since it has distinct Hodge--Tate weights), 
this last case follows from Lemma~\ref{lemma:potabmod} below.
\end{proof}

In the remainder of this section we prove Lemma~\ref{lemma:potabmod}, which was
used in the proof of Theorem~\ref{remainingcases}. %
We begin with some preliminary lemmas,
the first of which concerns
representations which have potentially abelian image.

\begin{lemma} \label{grouparguments}
Let~$F$ be a number field, and let~$\rho: G_{F} \rightarrow \GL_n(\Qbar_p)$
be a continuous irreducible representation which is de Rham at all
places dividing~$p$
and potentially abelian over a finite Galois extension~$L/F$.
Then there exist integers~$a$, $b$, with~$ab = n$,  and~$b$ pairwise distinct characters~$\chi_i:G_L\to\Qpbartimes$
such that
$$\rho |_{G_L} \simeq \bigoplus_{i=1}^{b}  (\chi_i)^{\oplus a}.$$
The action of~$\Gal(L/F)$ on the characters~$\chi_i$
via~$\chi^{\sigma}_i(g) = \chi_i(\sigma g \sigma^{-1})$  induces
a map
$$\Gal(L/F) \rightarrow S_b$$
with transitive image. Let~$\Gal(L/K_i)$ be the stabilizer of~$\chi_i$.
Then there exists an irreducible representation~$V_i$ of~$G_{K_i}$
such that~$V_i|_{G_L}\cong (\chi_i)^{\oplus a}$,
and
$$\rho \simeq \Ind^{G_F}_{G_{K_i}} V_i.$$
If~$\rho$ is not finite up to twist,
then:
\begin{enumerate}
\item the characters~$\chi_i$ are associated to algebraic Hecke characters of
  non-parallel weight. 
\item $b>1$. %
  \item If~$a=1$, then each~$K_i$ contains an imaginary CM field. 
\end{enumerate}
\end{lemma}

\begin{proof}
Since the~$\chi_i$
eigenspace is mapped to the~$\chi^{\sigma}_i$
eigenspace under the action of~$\rho(\sigma)$
for any lift of~$\sigma \in \Gal(L/F)$ to~$G_F$, the group~$\Gal(L/F)$ acts transitively
on the characters (since otherwise the direct sum
of the eigenspaces for~$\chi^{\sigma}_i$
for any given~$i$ would be a non-trivial~$G_F$-invariant
subspace of~$\rho$, and we are assuming that~$\rho$
is irreducible).
Similarly, the multiplicity of each~$\chi_i$
is independent of~$i$.
Let~$V$ be the vector space underlying the representation~$\rho$, and 
let~$V_i$ denote the subspace on which~$G_L$ acts by~$\chi_i$.
Since~$G_F$ preserves the decomposition~$V =\bigoplus V_i$, it
follows that~$V_i$ extends to a representation of~$G_{K_i}$
where~$\Gal(L/K_i)$ stabilizes~$\chi_i$.  By the orbit--stabilizer
theorem, $[K_i:F] = b$. By Frobenius
reciprocity, there is a non-trivial map~$V \rightarrow \Ind^{G_{F}}_{G_{K_i}} V_i$,
which (because~$V$ is irreducible) is an isomorphism, and hence~$V_i$
must also be irreducible.

Assume  for the remainder of the proof that~$\rho$ is not finite up to twist.
Each character~$\chi_i$ is de Rham and thus is either
a finite order character times an integer power of the cyclotomic
character or  has non-parallel weight. In the first case, after twisting~$\rho$
by a power of the cyclotomic character, we may assume that~$\chi_i$
has finite image. But then~$V_i$ and hence~$\rho$ also have finite
image, contrary to our assumption.
If~$b= 1$, then the projective image of~$\rho$ restricted to~$G_L$
is trivial,  and hence the image of
the projective representation~$P\rho$
is finite.  From the vanishing of~$H^2(F,\Qbar^{\times}_p)$~\cite[Thm.~4]{Serrelifting},  it follows $P\rho$ lifts to a genuine representation~$\widehat{\rho}: G_F \rightarrow \GL_n(\Qbar_p)$
which has finite image. Since~$\rho$ is irreducible, $\widehat{\rho}
\simeq \rho \otimes \chi$ for some~$\chi$,
and thus~$\rho$ is finite up to twist, once more contrary to our assumption.

If~$a=1$, then the action of~$G_{K_i}$ on~$V_i$
is via a character~$\psi_i$ which restricts to~$\chi_i$
over~$G_L$. We have already shown that~$\chi_i$ has 
non-parallel weight, 
 and thus~$\psi_i$ also corresponds to an algebraic Hecke character
of non-parallel weight, which implies that~$K_i$ contains an imaginary
CM field.
\end{proof}

\begin{lemma} \label{dimtwo}
Let~$F$ be a number field, and let~$\rho: G_F \rightarrow \GL_2(\Qbar_p)$
be a continuous irreducible representation which is de Rham at~$p$
and potentially abelian over a finite extension. Then either:
\begin{enumerate}
\item $\rho$ is finite up to twist, or
\item $\rho$ is automorphic for a cuspidal automorphic representation~$\pi$ for~$\GL_2/F$
which is the automorphic induction of an algebraic
Hecke character.
\end{enumerate}
\end{lemma}

\begin{proof}
Suppose that~$\rho |_{G_L} $ becomes reducible over a Galois extension~$L/F$,
and write~$\rho |_{G_L} \simeq \chi_1 \oplus \chi_2$. If~$\chi_1 = \chi_2$,
 then~$\rho$ is finite up to twist by Lemma~\ref{grouparguments}.
Thus we may assume that~$\chi_1 \ne \chi_2$, and then by Lemma~\ref{grouparguments},
we see that~$K = K_i/F$ is cyclic of degree~$2$ and~$\rho =
\Ind^{G_F}_{G_K} \psi_p$ where~$\psi_p$ corresponds to an algebraic Hecke character~$\psi$ of~$G_K$.
\end{proof}

\begin{lemma} \label{lemma:potabmod} Let~$\rho: G_{\Q} \rightarrow \GSp_4(\Qbar_p)$
be a continuous irreducible representation which is de Rham at~$p$. Suppose that:
\begin{enumerate}
\item There exists a  Galois extension~$L/\Q$ so that the image of~$\rho |_{G_L}$ is abelian.
\item If~$\nu$ is the similitude character, then~$\nu(c) = -1$, where~$c$ is complex conjugation.
\item $\rho$ is not finite up to twist.
\end{enumerate}
Then~$\rho$ is modular. %
\end{lemma}

\begin{proof} By Lemma~\ref{grouparguments},
there is a decomposition~$\rho |_{G_L} \simeq \bigoplus_{i=1}^{b}  (\chi_i)^{\oplus a}$
with~$ab = 4$. Since~$\rho$ is not finite up to twist, it follows from
Lemma~\ref{grouparguments} that~$b>1$. Suppose that~$b = 2$.
Then~$\rho \simeq \Ind^{G_{\Q}}_{G_{K_i}} V_i$ where~$[K_i:\Q] = 2$ and~$\dim(V_i) = 2$.
Note that since~$\Gal(L/K_i)$ is the stabilizer of a point with respect to the map~$\Gal(L/\Q) \rightarrow S_2$,
the field~$F = K_i$ does not depend on~$i$. 
We know that~$V_1$ and~$V_2$  are irreducible and potentially abelian, so by Lemma~\ref{dimtwo}
either both~$V_i$ are automorphic for~$\GL_2/F$, in which case~$\rho$
is modular, and we are done; %
or both~$V_i$ are finite up to twist, which we assume from now on.

If we write~$V = \Ind^{G_{\Q}}_{G_F} V_i$, then
$$\wedge^2 V = \Ind^{G_{\Q}}_{G_F} \det(V_i) \oplus \mathrm{Asai}_{F/\Q}(V_i),$$
where~$\mathrm{Asai}_{F/\Q}(V_i) |_{G_F} \simeq V_1 \otimes V_2$ (this only
 characterizes the representation over~$G_{\Q}$ up to quadratic twist but it is all that we will use
 in this argument).
Since~$V$ admits a symplectic form which is Galois invariant
up to a similitude character, we know that~$\wedge^2 V$ must contain a character.
We deduce either that~$\det(V_i)$
extends to~$\Q$ or~$\mathrm{Asai}(F/\Q)(V_i)$ is reducible. Suppose firstly
that~$\det(V_i)$ is the restriction to~$G_G$ of a character~$\chi$ of~$G_{\Q}$. The image of~$\chi$ lands in~$\OL^{\times}_E$ for some finite extension~$E/\Q_p$,
and so~$\chi$
factors through a quotient of~$G_{\Q}$
of the form~$\Z^r_p \oplus T$ for some finite group~$T$.
Define a new character~$\psi:G_{\Q}\to\cO_E^{\times}$ by sending topological
generators~$\sigma$ of each of these~$\Z_p$ factors to any square root
of~$\chi(\sigma)$, so that~$\chi \psi^{-2}$ has finite order. Then $\det(V_i \otimes
\psi|_{G_F}^{-1})$ has finite order, and $V_i \otimes
\psi|_{G_F}^{-1}$ is finite up to twist, so 
$V_i \otimes  \psi|_{G_F}^{-1}$ has finite image. Since \[V =  \Ind^{G_{\Q}}_{G_F} V_i 
     = \psi \otimes \Ind^{G_{\Q}}_{G_F} (V_i \otimes  \psi|_{G_F}^{-1}), 
    \]and $V_i \otimes  \psi|_{G_F}^{-1}$ has finite image, we see that~ $\rho$ is
    finite up to twist, contradicting our assumptions.

Hence we may assume that the Asai representation contains a character as a constituent, and in particular its 
restriction~$V_1 \otimes V_2$ to~$G_F$ does as well, and thus (since the~$V_i$
are irreducible over~$F$) we have~$V_1 \simeq V_2 \otimes \phi$
for some~$\phi$. This implies that the projective representations associated to~$V_1$ and~$V_2$
are isomorphic. Since~$V_1$ and~ $ V_2$ are $G_{\Q}$-conjugate, it follows that the projective
representation associated to~$V_i$ extends to~$\Q$, and thus (by Tate's theorem)
$V_i$ itself lifts (up to twist) to a representation~$V$ of~$G_{\Q}$. It then follows that
$$\rho \simeq V \otimes \Ind^{G_{\Q}}_{G_F} \chi$$
where now~$F/\Q$ is an imaginary quadratic field,~$\chi$ is an algebraic Hecke character,
and~$V$ is a representation of~$G_{\Q}$ of finite image. 
If~$V$ is induced, then~$V$  and~$\rho$ are automorphic, so
we may assume that~$V$ is not
induced. The~$G_{\Q}$-module $V$ admits a unique symplectic form invariant up to  a similitude
character which is given by~$\det(V)$, but~$V$ does not admit any corresponding orthogonal form, 
since~$V$ is not induced. On the other hand, we also see that~$W = \Ind^{G_{\Q}}_{G_F} \chi$ admits 
a symplectic form which is invariant up to similitude character~$\det(W)$,
and an orthogonal form which is invariant up to similitude character~$\det(W) \eta_{F/\Q}$,
where~$\eta_{F/\Q}$ is the quadratic character associated to the imaginary quadratic
field~$F/\Q$. We deduce that the unique symplectic form
on~$\rho = V \otimes W$ has similitude character~$\det(V) \det(W) \eta_{F/\Q}$, which is
odd if and only if~$\det(V)$ is odd, since~$\det(W) \eta_{F/\Q}$ is even.
Thus the oddness assumption implies
that~$V$ is an \emph{odd} Artin representation, and thus~$V$ is modular by
known cases of the Artin conjecture~\cite{MR3581178, MR3904451}, and the automorphy
of~$\rho$ follows.

Finally, suppose that~$b=4$, so~$a=1$. By Lemma~\ref{grouparguments},
 there exists a degree~$4$ field~$K/\Q$ such
that~$\rho \simeq \Ind^{G_{\Q}}_{G_K} \chi$, where~$\chi$ corresponds
to an algebraic Hecke character of non-parallel weight, so that
~$K/\Q$ contains an imaginary CM field. In particular, either~$K$ is itself an imaginary CM field, and thus contains a real
quadratic subfield~$E = K^{+}$, or~$K$ %
contains an imaginary
quadratic subfield~$E$. In either case, we see that~$\varrho =
\Ind^{G_E}_{G_K} \chi$ corresponds to a cuspidal automorphic
representation of~$\GL_2/E$, and by another application of automorphic
induction we deduce that ~$\rho = \Ind^{G_{\Q}}_{G_E} \varrho$ is automorphic.
\end{proof}

\begin{rem} Various rationality considerations (see~\cite{MR2982436}) imply that, if~$\rho$ is a potentially
abelian Galois representation associated to an abelian surface~$A/\Q$,
then~$\rho$ is actually potentially abelian over a \emph{solvable} extension of~$\Q$, which
can be used to simplify the argument in this case, On the other hand.
Lemma~\ref{lemma:potabmod} is conjecturally
 still true without either
the oddness assumption, or the finiteness up to twist condition, although presumably 
extremely difficult.
In the first case,
 it would include the automorphy of representations of the 
form~$V \otimes  \Ind^{G_{\Q}}_{G_F} \chi$ where~$\chi$ is an algebraic Hecke character
of a CM field~$F$ and~$V$ is an even Galois representation with projective image~$A_5$,
which would imply the automorphy of such a~$V$.

Similarly, assume oddness holds but drop the finiteness up to twist condition.
The group~$S_6$ is a subgroup of~$\PGSp_4(\C)$; this~$S_6$ can also be seen inside~$\PGSp_4(\F_3)$.
Let~$L/\Q$ be any~$S_6$ extension such that complex conjugation is odd.
Then there is  a projective
representation~$G_{\Q} \rightarrow \Gal(L/\Q) \simeq S_6 \hookrightarrow \PGSp_4(\C)$.
Any lift to an Artin representation~$\rho: G_{\Q}  \rightarrow \GSp_4(\C)$ will be odd, but the modularity of~$\rho$ 
is unknown for any representation with projective image~$S_6$ regardless of the image of complex conjugation.
Fortunately, neither case is relevant for the automorphy of abelian surfaces over~$\Q$.
\end{rem}

\subsection{Residual modularity theorems (modulo $2$)} \label{residualmod2} 
The goal in this section is to prove some residual modularity theorems for mod-$2$ representations
with image~$A_6$ or~$S_6$.
We will need the following variation of Lemma~\ref{lem: lifting genus 2 curves and finite flat group
    schemes} for the prime~$p=2$.

\begin{lem}  \label{lem:liftingvariant}%
 Let ~$\cO$ be the ring of integers
  in a finite extension~$K$ 
    of~$\Q_2$ with residue
  field~$\F$.  Let~$G_1/\cO$ be a $2$-torsion
  finite flat group scheme of order~$16=2^{4}$, together with an isomorphism~$\lambda: G_1 \rightarrow G^{\vee}_1$
  such that~$\lambda^{\vee} = - \lambda$. 
Suppose that~$A_0/\F$ is a principally polarized ordinary abelian surface
  with~$A_0[2]/\F \simeq G_1/\F$.
Then there exists a lift of
~$A_0$  to  a principally
polarized  abelian surface ~$A/\cO$
with ~$A[p] \simeq G_1$.
\end{lem}

\begin{remark}[Remarks on the proof and the statement
 of Lemma~\ref{lem:liftingvariant}]
Let~$p=2$. By Serre--Tate theory, we are reduced to finding an appropriate
lifting~$A_0[p^{\infty}]$ to a Barsotti--Tate group,
together with a lifting of~$\lambda$ to make~$A_0[p^{\infty}]$ a quasi-polarized BT. By a result of Grothendieck~\cite{MR801922}, there is no issue in lifting~$A_0[p^{\infty}]$ as
a Barsotti--Tate group, so the subtlety is imposing the polarization. 
We proved an analogous statement in Lemma~\ref{lem: lifting genus 2 curves and finite flat group schemes}
(without any ordinary hypothesis) using results from~\cite{MR1827029}.
Wedhorn's argument in~\cite[(2.17)]{MR1827029} involves certain constructions in which one obtains
a polarization by an averaging procedure involving dividing by~$2$ --- this naturally causes issues when~$p = 2$.
One difficulty  is that, when~$p=2$, one needs to  decide
what it means for a pairing on a finite flat
group scheme~$G$ to be alternating rather
than skew symmetric.
If~$G = \alpha_2/\F_2$, then~$G$ is Cartier self-dual via
the map~$\alpha_2 \times \alpha_2 \rightarrow \mathbf{G}_m$
given on points by~$(x,y) \mapsto 1+xy$. This is alternating
on points (since~$x^2 = 0$ for~$x \in \alpha_2(A)$) but
one does not want to regard it as an alternating pairing.
Instead, following~\cite[\S3.2]{Lurie} (where the idea is attributed in part to de Jong),
one could define a pairing~$G \times G \rightarrow \mathbf{G}_m$ of finite flat
group schemes 
to be \emph{strongly alternating} if, fpqc locally on the base,
 there is
is a central extension 
\numequation
\label{heisenberg}
1 \rightarrow  \mathbf{G}_m
\rightarrow \widetilde{G} \rightarrow G \rightarrow 1
\end{equation}
such that the pairing arises from the commutator pairing on this extension.

One strategy would be to determine the precise conditions for any~$G = G_1$
of exponent~$2$ 
admitting an  isomorphism~$\lambda: G_1 \isoto G^{\vee}_1$ with~$\lambda^{\vee} =  \lambda$ (since~$p = 2$  there
is no choice of sign) to give
rise to 
 a corresponding Heisenberg
group extension~$\widetilde{G}$ of the form~(\ref{heisenberg}), 
and then to prove a version of this lemma
without any ordinary hypothesis on~$A_0/\F$, but with a suitably modified definition
of what it means for~$G/\OL$ to be quasi-polarized.

Alternatively, instead of addressing any of the more subtle issues which 
might arise in the general case, we 
exploit the assumption that~$A_0/\F$ is ordinary. In this case,  the assumption that a commutative
finite flat group  scheme~$G$
is an extension of an \'{e}tale group scheme by a
multiplicative (dual to \'{e}tale) group scheme simplifies the situation considerably: There is a  connected-\'{e}tale sequence
$$1 \rightarrow G^0 \rightarrow G \rightarrow  G^{\et} \rightarrow 1,$$
where~$G^0$ is multiplicative and  so its Cartier dual is \'{e}tale.
Now any isomorphism~$\lambda: G \simeq G^{\vee}$
clearly has the property that the induced map
 \numequation
\label{ordinaryOK}
G^0 \rightarrow G \stackrel{\lambda}{\simeq} G^{\vee} \rightarrow (G^0)^{\vee} = ((G^0)^{\vee})^{\et}
\end{equation}
is trivial, and so any such~$\lambda$ %
will automatically be alternating 
on the generic fibre. %
In particular, 
 we only work with the assumption that
there exists an isomorphism~$\lambda: G_1 \rightarrow G^{\vee}_1$ with~$\lambda^{\vee} =  - \lambda$ %
(the
 sign makes no difference for finite flat group schemes annihilated by~$2$),
even though 
one expects
this may will be the ``wrong'' definition in the non-ordinary case. 
  \end{remark}
  
  \begin{proof}[Proof of Lemma~\ref{lem:liftingvariant}]
The assumption that~$A_0/\F$ is ordinary implies that the corresponding~$2$-divisible group splits
into toroidal  (ind-multiplicative)
and (ind)-\'{e}tale parts which are Cartier
dual to each other. 
Over~$\OL/\pi^n_K$ or over~$\OL$, these factors have unique lifts,
and the lifts 
of~$A_0[2^{\infty}]$ are equivalent to the category of extensions of these factors
(\cite[Prop~2.1]{MessingThesis}).
But over~$\OL$, $2$-divisible groups are determined by their generic fibres, and so the lifts are classified  in terms of
Galois cohomology. More precisely,
 if~$V$ denotes the free (rank~$2$) $\Z_2$ module
  corresponding to the Pontryagin dual  of the (ind)-\'{e}tale
 part of~$A_0[2^{\infty}]$, and~$W = V^{\vee}$ the~$\Z_2$-dual of~$V$, then the extensions of interest are computed by the group
$$H^1_f(K,W \otimes W(1)).$$
The result~\cite{MR801922} then implies the surjectivity of the reduction map:
$$H^1_f(K,W \otimes W(1)) \rightarrow H^1_f(K,\Wbar \otimes \Wbar(1)),$$
where~$\Wbar =W/2$. 
Now we wish to impose  the condition that there exists a suitable polarization~$\lambda$.
There is an exact sequence of flat~$\Z_2$-modules %
$$0 \rightarrow S(W) \rightarrow W \otimes W \rightarrow \wedge^2(W) \rightarrow 0,$$
where~$S(W)$ is the submodule generated by~$x \otimes x$ for all~$x \in
W$. The vector space~$S(W) \otimes K$ is isomorphic to~$\Sym^2(W) \otimes K$, but
this is not used below.
By purity, the (generalized) eigenvalues of Frobenius on~$W \otimes W$
cannot have absolute value~$1$
and so in particular are~$\ne 1$.
It follows that for ~$M$ equal to any of~$S(W)$, $W \otimes W$,
or~$\wedge^2 W$, we have  $H^2(K,M(1)) = 0$,
and $H^1_f(K,M(1)) = H^1(K,M(1))$.
We say that a class~$\eta$ in~$H^1_f(K,W \otimes W(1))
\subset H^1(K,W \otimes W(1))$ 
is alternating if it lies in the image of~$H^1(K,S(W)(1))$, and we denote the alternating classes by~$H^1_f(K,W \otimes W(1))^{\Alt}$.
We similarly write~$H^1_f(K,\Wbar \otimes \Wbar(1))^{\Alt}$ for the 
same condition modulo~$2$. (The reason this is the correct choice
is ultimately explained by equation~\eqref{meaningoflife} below.)
We  have a commutative diagram as follows.
\[\xymatrix{H^1_f(K,S(W)(1)) \ar[r] \ar[d] & H^1_f(K,W \otimes W(1))  \ar[r]  \ar@{->>}[d] &H^1_f(K,\wedge^2(W)(1))  \ar[d] \ar[r] & 0  \\
H^1_f(K,S(\Wbar)(1)) \ar[r] & H^1_f(K,\Wbar \otimes \Wbar(1))  \ar[r]   &H^1_f(K,\wedge^2(\Wbar)(1)) & 
}\]
Because the~$H^2$ groups vanish, the kernel of each vertical map consists of classes divisible by~$2$. 
Now take an alternating class~$\etabar \in H^1_f(K,\Wbar \otimes \Wbar(1))^{\Alt}$. It lifts to~$\eta \in H^1_f(K,W \otimes W(1))$,
which then maps to a class~$\nu \in H^1_f(K,\wedge^2(W)(1))$ whose reduction~$\nubar$ is trivial. But that implies that~$\nu$ is divisible by~$2$,
and thus, writing~$\nu = 2 \gamma$, 
and lifting~$\gamma$ to~$\widetilde{\gamma} \in H^1_f(K, W \otimes W(1))$,
we see that~$\eta - 2 \widetilde{\gamma}  \in
H^1_f(K, W \otimes W(1))^{\Alt}$, 
and hence there is a surjection
$$H^1_f(K, W \otimes W(1))^{\Alt} \rightarrow H^1_f(K,\Wbar \otimes \Wbar(1))^{\Alt}.$$
It now suffices to show that $H^1_f(K, W \otimes W(1))^{\Alt}$
 classifies possible 
principally quasi-polarized
Barsotti--Tate groups~$G/\OL$ lifting~$A_0[2^{\infty}]$, equivalently,
a BT~$G/\OL$ together with an isomorphism~$\lambda: G \rightarrow G^{\vee}$
with~$\lambda^{\vee} = -\lambda$,
whereas
 $H^1_f(K,\Wbar
\otimes \Wbar(1))^{\Alt}$ classifies
finite flat group schemes over~$\OL$ lifting~$A_0[2]$ together with an isomorphism~$\lambda: G \rightarrow G^{\vee}$
with~$\lambda^{\vee} = -\lambda$.
In both settings, the corresponding lifts are determined by extensions of (fixed) \'{e}tale by multiplicative group schemes,
and these extensions
 are determined by their generic fibres.
In either case, $\lambda$ induces a skew-symmetric pairing~$J$ on the generic fibre which 
 (as explained in the
 discussion surrounding ~\eqref{ordinaryOK} using the ordinary hypothesis) is 
alternating. Moreover, the generic fibre of the connected part (respectively, \'{e}tale part) is isotropic
with respect to this pairing.
That implies that
the action on the generic fibre factors through the generalized symplectic
group, and in particular that the extension class of the \'{e}tale by
multiplicative part is
alternating in the sense described above.
(This is equivalent to the computation that
\begin{equation}
\label{meaningoflife}
\left( \begin{matrix} I & A \\ 0 & I \end{matrix} \right) 
\left( \begin{matrix} 0 & I \\ -I & 0 \end{matrix} \right) \left( \begin{matrix} I & A \\ 0 & I \end{matrix} \right)^T = \left( \begin{matrix} 0 & I \\ -I & 0 \end{matrix} \right) 
\end{equation}
if and only if~$A$ is symmetric.)
Conversely, once the image of the Galois representation lies in the generalized symplectic group
compatible with the connected part being isotropic, the
Barsotti--Tate group (or finite flat group scheme) admits a suitable~$\lambda$. %
We deduce that~$G_1$ corresponds to a class in~$H^1_f(K,\Wbar \otimes \Wbar(1))^{\Alt}$, and that there exists a lift of~$A_{0}[2^{\infty}]$ to a principally quasi-polarized $\BT$
  $G/\cO$ with~$G[2]\simeq G_1$, and we conclude as in
  the proof of
   Lemma~\ref{lem: lifting genus 2 curves and finite flat group
    schemes}.
\end{proof}

We also offer the following alternative proof using stacks, for those who are gripped to the pages of this manuscript
and don't wish it to end:

\begin{proof}[Alternate proof of Lemma~\ref{lem:liftingvariant}]
 Let us consider the $p$-divisible group $\mu_{p^\infty}^{2} \oplus \qq_p/\ZZ_p^{2}$ equipped with its standard polarization over $\Fbar_p$. 
By Serre--Tate theory, the moduli space of polarized extension  of $0 \rightarrow \mu_{p^\infty}^{2} \rightarrow \mathcal{G} \rightarrow \qq_p/\ZZ_p^{2} \rightarrow 0$ on local Artinian $W(\Fbar_p)$-algebras is represented by 
$$\widehat{\mathbf{G}}_m^3 = \Spf W(\Fbar_p)\llb X_1,X_2,X_3\rrb.$$
Identify $\F = \F_q$ and let  $\phi$ be the $q$-th power  Frobenius which topologically generates  $\mathrm{Gal}(\Fbar_p/\F)$. The $p$-divisible group $A_0[p^\infty]$ is isomorphic over   $\Fbar_p$ to $\mu_{p^\infty}^{2} \oplus \qq_p/\ZZ_p^{2}$. This implies that the moduli stack of deformations (as polarized extensions) of  $A_0[p^\infty]$ to local $W(\F)$-Artinian algebras   is represented by 
$$\widehat{\mathbf{G}}_m^3 / \phi^{\hat{\Z}}$$ where $\phi$ acts naturally on $W(\Fbar_p)$ and the action on the Serre--Tate parameters is induced  by the action of Frobenius on  $T_p A_0^{{\et}} \otimes T_p (A_0^m)^D$. Let $ W(\Fbar_p)/p^n= W_n(\Fbar_p)$. 
We can think of  $\widehat{\mathbf{G}}_m^3 / \phi^{\hat{\Z}}$ as the ind-stack $\colim \Spec W_n(\Fbar_p)\llb X_1,X_2,X_3\rrb / \phi^{\hat{\Z}}$ and each  $\Spec W_n(\Fbar_p)\llb X_1,X_2,X_3\rrb / \phi^{\hat{\Z}}$ is the inverse limit on $r$ of:
 $$\Spec W_n(\Fbar_{q^r})\llb X_1,X_2,X_3\rrb / \phi^{\Z/q^{r}\Z}$$ where $r$ is large enough so that the action on the Serre--Tate parameters of $\phi^{q^{r}}$ is trivial modulo $p^n$. The map $$\Spec W_n(\Fbar_{q^r})\llb X_1,X_2,X_3\rrb  \rightarrow \Spec W_n(\Fbar_{q^{r}})\llb X_1,X_2,X_3\rrb / \phi^{\Z/q^{r}\Z}$$ is formally \'etale as the group $\Z/q^{r}\Z$ is \'etale. 
The moduli of polarized extensions 
$$0 \rightarrow \mu_{p}^{2} \rightarrow G \rightarrow (\Z/p\Z)^{2} \rightarrow 0$$
  on local Artinian $W(\Fbar_p)$-algebras is the quotient stack $\widehat{\mathbf{G}}_m^3/ \widehat{\mathbf{G}}_m^3$ where
  each copy of $\widehat{\mathbf{G}}_m$ acts on itself through the map $(1+x,1+y) \mapsto (1+x)^p(1+y)$. 
  (The moduli of all extensions is given by~$\widehat{\mathbf{G}}_m^4/ \widehat{\mathbf{G}}_m^4$, and
  the inclusion of the polarized extensions into all extensions corresponds 
  in the previous argument to the inclusion of~$S(W)$ into~$W \otimes W$.)
The map $\widehat{\mathbf{G}}_m^3 \rightarrow  \widehat{\mathbf{G}}_m^3/ \widehat{\mathbf{G}}_m^3$ is formally smooth. 
We now consider  the map
\begin{equation} \label{thismap}
\widehat{\mathbf{G}}_m^3/\phi^{\hat{\Z}}  \rightarrow  \widehat{\mathbf{G}}_m^3/ \widehat{\mathbf{G}}_m^3 /\phi^{\hat{\Z}}. 
\end{equation}
The map~\eqref{thismap}  sends a deformation (as a polarized extension) of $A_0$ to a deformation (as a polarized extension) of $A_0[p]$. 
The map~\eqref{thismap} is moreover the inductive limit of the maps 
\begin{equation} \label{thesemaps}
\widehat{\mathbf{G}}_m^3\vert_{W_n(\Fbar_p)}/\phi^{\hat{\Z}}  \rightarrow  \widehat{\mathbf{G}}_m^3\vert_{W_n({\F}_{q^r})}/ \widehat{\mathbf{G}}_m^3\vert_{W_n({\F}_{q^r})}/\phi^{\hat{\Z}}.
\end{equation}
 In turn, these maps~\eqref{thesemaps} are the inverse limit of the maps:
$$\widehat{\mathbf{G}}_m^3\vert_{W_n({\F}_{q^r})}/\phi^{{\Z}/q^{r}\Z}  \rightarrow  \widehat{\mathbf{G}}_m^3\vert_{W_n({\F}_{q^r})}/ \widehat{\mathbf{G}}_m^3\vert_{W_n({\F}_{q^r})}/\phi^{{{\Z}/q^{r}\Z}} $$
This last  map is  formally smooth because the groups $\widehat{\mathbf{G}}_m^3\vert_{W_n({\F}_{q^r})}$ and ${\Z}/q^{r}\Z$ 
are both formally smooth.  

By assumption, we begin with~$G_1/\cO$, which  is a point $\Spf \ocal \rightarrow \widehat{\mathbf{G}}_m^3/ \widehat{\mathbf{G}}_m^3/\phi^{\hat{\Z}} $, and our goal is to lift  this to an $\Spf~\ocal$-point of $\widehat{\mathbf{G}}_m^3/\phi^{\hat{\Z}}$. Let $\varpi$ be a uniformizer of~$\cO$. Assume that we have found a lift of  $G_1\vert_{ \Spec(\cO/\varpi^k)}$ to $G\rightarrow \Spec(\cO/\varpi^k)$;  we shall upgrade it to  a lift on $\Spec(\cO/\varpi^{k+1})$ of $G_1\vert_{ \Spec(\cO/\varpi^{k+1})}$, and
then we are done by induction. 

There exists, $n, r$ such that we have a commutative diagram (given by the solid arrows):
\[
\begin{tikzcd}
\Spec (\mathcal{O}/\varpi^k) \ar[r] \ar[d] & \widehat{\mathbf{G}}_m^3\vert_{W_n({\mathbb{F}}_{q^r})}/\phi^{{\mathbb{Z}/q^{r}\mathbb{Z}}} \ar[d] \\
\Spec (\mathcal{O}/\varpi^{k+1}) \ar[r] \ar[ur,dotted]& \widehat{\mathbf{G}}_m^3\vert_{W_n({\mathbb{F}}_{q^r})}/ \widehat{\mathbf{G}}_m^3\vert_{W_n(\mathbb{F}_{q^r})}/\phi^{{\mathbb{Z}/q^{r}\mathbb{Z}}}
\end{tikzcd}
\]
 By formal smoothness, we can produce the lift given by the dotted arrow, completing the proof.
\end{proof}

\subsubsection{Moduli of~$X$ with fixed~$\Jac(X)[2]$}

Let~$F$ be a global or local field of characteristic zero, and let
$$\rhobar: G_{F} \rightarrow \GSp_4(\F_2)$$
be a continuous representation. Under the identification (Lemma~\ref{explicit}) of~$S_6$ with~$\GSp_4(\F_2)$,
there exists a corresponding degree~$6$ separable polynomial~$f(x)$ such that~$\rhobar$ is 
isomorphic to the~$2$-torsion representation on the Jacobian of~$y^2 = f(x)$.
Let~$K =  F[x]/f(x) \simeq \prod F_i$, which is a degree~$6$ \'{e}tale~$F$-algebra.
Given~$\theta \in K$, the multiplication by~$\theta$ map~$K \rightarrow K$ naturally has a characteristic 
polynomial of degree~$6$ with roots we denote by~$\sigma \theta$. If we fix a basis 
for~$K$ over~$F$, for example given by the powers of~$x$, this map is compatible
with extensions of~$F$.
If we identify~$K$ with~$F^6$, then
by the primitive element theorem, the~$\sigma \theta$ will be distinct for~$\theta$ outside
a finite number of hyperplanes (which are defined over the splitting field of~$f(x)$ and
compatible with field homomorphisms~$F\to F'$).
If  the~$\sigma \theta$ are distinct, then
$$X: y^2 = \prod (x - \sigma \theta),$$
will be a smooth genus two curve over~$F$
with~$\Jac(X)[2] \simeq \rhobar$. We have
therefore constructed
a smooth rational variety~$Z(\rhobar) \subset \mathbf{P}^6$ over~$F$ given by the complement
of finitely many hyperplanes, whose~$F$-rational points give
smooth genus two curves~$X$ with $2$-torsion given by~$\rhobar$. 
Moreover, the construction of~$Z(\rhobar)$
(having fixed~$K$) is compatible with both extensions of~$F$ and completions at primes of~$F$.
There is a map from~$Z(\rhobar)$ to the corresponding moduli stack~$\MMM_2(\rhobar)$,
but to avoid any issues concerning fields of moduli versus fields of definition it is fine for
our purposes to work directly with~$Z(\rhobar)$.

\begin{lemma} \label{localmod2}
Let~$\rhobar:G_{\Q} \rightarrow \GSp_4(\F_2)$ be a continuous
representation unramified at~$3$. Assume that there exists
a finite flat model~$\WW/\Z_2$ for~$\rhobar$ over~$\Z_2$ which is
isomorphic to its Cartier dual and which is ordinary, that is,
the extension of an \'{e}tale group scheme by a multiplicative group scheme. 
Suppose that
$\rhobar(\Frob_3)$ is non-trivial.
There exists an abelian surface~$A/\Q$ such that:
\begin{enumerate}
\item $A$ has good ordinary reduction at~$3$ and is~$3$-distinguished.
\item $A$ has good ordinary reduction at~$2$. There is an
isomorphism of finite flat group schemes~$A[2]/\Z_2 \simeq \WW$, and
the characteristic polynomial~$Q(x)$ of Frobenius at~$2$ satisfies
\numequation \label{not4C12C}
Q(x) \not\equiv 
 x^4 \pm x^3 + 2 x^2 \mp  x + 1 \bmod 3.
 \end{equation}
 \item %
   $\rhobar_{A,3}$ is surjective.
 \item $\rhobar_{A,2}\cong\rhobar$.
\end{enumerate}
\end{lemma}

\begin{proof}
Let~$Z = Z(\rhobar)$.
The conditions we are imposing at~$2$ and~$3$ are open conditions
in~$Z(\Q_2)$ and~$Z(\Q_3)$ respectively. The condition that~$\rhobar_{A,3}$
is surjective holds outside a thin set.
Since~$Z$ is smooth and rational,
by Lemma~\ref{rationalmb}  once we find suitable points on~$Z(\Q_2)$ and~$Z(\Q_3)$,
we obtain an~$A$ which has the required properties.

Let us consider~$Z(\Q_3)$.
With~$F = \Q$, let~$K$ denote the 
degree~$6$ \'{e}tale~$F$-algebra corresponding to~$\rhobar$ as described above. Since we are assuming that~$\rhobar$
is unramified at~$3$, we certainly have that~$K$ 
is unramified at~$3$, so the possible %
completions~$K \otimes \Q_3$ are determined by a partition of~$6$. There
are exactly~$11$ such partitions, the partition~$6=1+1+1+1+1+1$ corresponding
to the case when~$\rhobar(\Frob_3)$ is trivial, which we are excluding.
For the remaining~$10$ partitions, we
now produce an explicit~$X: y^2 = f(x)$ with good ordinary reduction at~$3$ which is~$3$-distinguished.
We actually write down~$X/\Q$ with these properties. Note that for~$y^2 = f(x)$ where one
of the Weierstrass points is at~$\infty$, the corresponding partition corresponds
to the factorization of~$f(x)$ over~$\Q_3$ plus another copy of~$\Q_3$.
In other words, the partition corresponds to the degrees of the (unramified) fields of definition of
the~$6$ Weierstrass points of~$X$ over~$\Q_3$.
\begin{center}
 \begin{tabular}{*3c}
\toprule
 $f(x)$ & partition & $N$ \\
\midrule
$4 x^5 + 32 x^4 + 64 x^3 + x^2 + 4 x$ &
$[ 1, 1, 1, 1, 2 ]$ & 1051 \\
$4 x^5 - 7 x^2 + 4 x$ &
$[ 1, 1, 1, 3 ]$  & 709 \\
$4 x^5 - 11 x^4 + 6 x^3 + 3 x^2 - 2 x + 1$ &
$[ 1, 1, 2, 2 ]$ & 1415 \\
$x^6 + 4 x^5 - 6 x^4 - 32 x^3 + x^2 + 64 x + 28$ &
$[ 1, 1, 4 ]$ & 389 \\
$x^6 + 2 x^5 + 5 x^4 + 4 x^3 - 4 x - 8$ &
$[ 1, 2, 3 ]$  & 847 \\
 $x^6 + 2 x^5 + 3 x^4 - x^2 + 2 x + 1$ & 
$[ 1, 5 ]$ & 349   \\
$x^6 + 4 x^4 + 6 x^3 - 8 x^2 + 1$ &
$[ 2, 2, 2 ]$ & 7165 \\
$x^6 + 2 x^4 + 2 x^3 + 5 x^2 + 2 x + 1$ &
$[ 2, 4 ]$  & 353 \\
$x^6 + 2 x^5 + 5 x^4 - 10 x^3 + 10 x^2 - 4 x + 1$ &
$[ 3, 3 ]$ & 4889 \\
$x^6 + 4 x^5 - 6 x^4 + 2 x^3 + x^2 - 2 x + 1$ &
$[ 6 ]$  & 1343 \\
\bottomrule
 \end{tabular}
 \end{center}
 
 Let us now turn to the prime~$2$. 
By Lemma~\ref{lem:liftingvariant},   the  required  abelian
   surface~$A/\Z_2$ will exist provided that there is an~$A_0/\F_2$ with~$A_0[2] \simeq \WW/\F_2$ (also
    satisfying equation~\eqref{not4C12C}). 
 The finite  flat group scheme~$\WW/\Z_2$ 
 is an extension of 
 an \'{e}tale group scheme~$\Vbar$
 by its Cartier dual~$\Vbar^{\vee}(1)$, and in particular
~$\WW/\F_2$ is determined by~$\Vbar$. There are three possibilities for~$\Vbar$:
  \begin{enumerate}
  \item $\Vbar$ is trivial as a~$G_{\F_2}$-module,
  \item $G_{\F_2}$ acts on~$\Vbar$ via a (non-semi-simple) element of order~$2$,
  \item $G_{\F_2}$ acts on~$\Vbar$ via a (semi-simple) element of order~$3$.
  \end{enumerate}
  It now suffices to find an~$A_0/\F_2$ of each form.

  One subtlety is that, given~$\rhobar$, the finite
 flat group scheme~$\WW/\Z_2$ is not determined by~$\rhobar$.
 Consider the following two examples:
 \begin{enumerate}
  \item The action of~$G_{\F_2}$ on~$\Vbar$ is trivial, 
 and the extension class
 of~$\Vbar = (\Z/2\Z)^2$ 
 by its Cartier dual~$\Vbar^{\vee}(1) = (\mu_2)^2$ is a direct sum of two extensions  corresponding to
 the \emph{unramified} class in~$H^1_f(\Q,\mu_2) = \Z^{\times}_2/\Z^{\times 2}_2$.
 \item The action of~$G_{\F_2}$ on~$\Vbar$ has order~$2$, and the group scheme
 is~$\Vbar \oplus \Vbar^{{\vee}}(1)$.
 \end{enumerate}
 In both cases, the representation~$\rhobar$ is unramified of order~$2$, and the
 non-trivial elements in the images are given by
 $$\left(\begin{matrix} 1 & 0 & 1 & 0 \\  0 & 1 & 0& 1 \\ 0 & 0 & 1 & 0 \\ 0 & 0 & 0 & 1 \end{matrix}\right), \quad \left(\begin{matrix} 1 & 1 & 0 & 0 \\  0 & 1 & 0& 0 \\ 0 & 0 & 1 & 1 \\ 0 & 0 & 0 & 1 \end{matrix}\right),$$
 respectively. Under the isomorphism of Lemma~\ref{explicit}, these are equal
 to~$(12)(34)(56) \in S_6$ and~$(12)(35)(46) \in S_6$, and so are conjugate.
 (They are not, however, conjugate inside the Siegel parabolic~$C((12)(34)(56)) = S_4 \times S_2$
 described in Lemma~\ref{boxerstrass}.)
 
  We now consider the following examples of genus~$2$ curves
  given by the equations~$y^2 = f(x)$ where~$f(x)$ is
  listed in the table below. One can check that the corresponding minimal models have good ordinary reduction at~$2$
  and compute the corresponding polynomial~$Q(x)$. Since these curves
  are defined over~$\Q$, one can also compute the global conductor, which is indicated in the table by~$N$.
  \begin{center}
 \begin{tabular}{*4c}
\toprule
 $f(x)$ & partition & $Q(x) $ &  $N$ \\
\midrule
$-4 x^5 + x^4 + 6 x^3 - 3 x^2 - 4 x$ &
$[ 1, 1, 1, 1, 1,1 ]$ & $ x^4 - x^2 + 4$ & 3451 \\
$x^6 - 12 x^4 + 2 x^3 + 16 x^2 + 8 x + 1$ &$ [2,2,2]$ & 
$ x^4 + 2 x^3 + 3 x^2 + 4 x + 4$  & 2225 \\
$x^6 - 4 x^5 + 2 x^4 + 2 x^3 + x^2 - 2 x + 1$ &  
$[ 3, 3 ] $&   $ x^4 + 3 x^3 + 5 x^2 + 6 x + 4$ &
713 \\
\bottomrule
 \end{tabular}
 \end{center}
 Here the partition indicates the factorization of~$f(x)$ over~$\Q_2$; the
corresponding Galois extension is cyclic and unramified of degree~$1$,
$2$, and~$3$ (in that order).
It follows immediately in the first and last cases that~$A[2]/\Z_2$ 
is the \emph{split} extension of~$\Vbar$ by~$\Vbar^{\vee}(1)$,
where~$G_{\F_2}$ acts on~$\Vbar$ through a cyclic group of order~$1$ 
and~$3$ respectively.
In the second case, we still need to check (in light of the example above)
that~$G_{\F_2}$ acts on~$\Vbar$ through a cyclic group of order~$2$.
 In this case, a  smooth model~$C/\Z_2 $  is given by the equation
 $$y^2 + (x^3+1)y = -3 x^4 +  4 x^2 + 2x,$$
 from which we find that~$\Jac(C)(\F_2) \simeq  \Z/14 \Z$.
 If~$\Vbar$ was trivial, then~$\Jac(C)(\F_2)$ would contain~$(\Z/2 \Z)^2$ as a subgroup, which it does not.
 Hence we deduce that
 the action of~$G_{\F_2}$ on~$\Vbar$ is through a cyclic group
 of order~$2$ (which suffices for our purposes, but  is not sufficient to determine~$\WW$ as
 an 
 extension). Thus we obtain a suitable~$A_0/\F_2$ in all possible cases.
\end{proof}

 \begin{rem}
Note that the method of proof of Lemma~\ref{localmod2} fails when~$\rhobar(\Frob_3)$ is trivial. %
This would imply that~$X$ has~$6$ Weierstrass points over~$\F_3$,
but if~$X$ has good reduction at~$3$ these points are distinct, %
and they are exactly the ramification
points over the map to~$\mathbf{P}^1$. But~$\mathbf{P}^1(\F_3)$ has only~$4 < 6$ points,
so this is impossible.
It seems unlikely one can avoid this using~$X$ for which~$A = \Jac(X)$ has good reduction
but~$X$ does not; at least the idea of using a product of elliptic curves does not work,
since if $\# E[2](\F_3) = 4$,
 then by the Hasse bounds~$\#E(\F_3) = 4$ and~$a_3 = 0$, and~$E$ is supersingular. 
\end{rem}
 
 By combining this with our main modularity theorem for abelian
 surfaces, we  deduce the following:
 
 \begin{theorem} \label{theoremWisneeded}
 Let~$\rhobar: G_{\Q} \rightarrow \GSp_4(\F_2)$ be a continuous %
 representation which is unramified at~$3$ and such that~$\rhobar(\Frob_3)$ is non-trivial. 
 Suppose as in Lemma~\ref{localmod2} that
 there exists
 a finite flat model~$\WW/\Z_2$ for~$\rhobar$ over~$\Z_2$ which is
isomorphic to its Cartier dual and which is ordinary.
 Then~$\rhobar$ is ordinarily modular
 of weight~$2$ and level prime to~$6$.
  \end{theorem}

 \begin{proof} Consider the abelian surface~$A/\Q$  with~$\rhobar_{A,2} = \rhobar$
 whose existence
 follows from Lemma~\ref{localmod2}. It suffices to show that~$A$ is modular.
 Condition~\eqref{not4C12C} %
 guarantees (by  Lemma~\ref{atlas}) that~$\rhobar_{A,3}(\Frob_2)$
is not of conjugacy class~$4C$ or~$12C$ in~$\PGSp_4(\F_3) \setminus \PSp_4(\F_3)$.
 Moreover, $A$ has  good ordinary  reduction %
 and~$3$,
is~$3$-distinguished, and~$\rhobar_{A,3}$ is surjective.
Thus~$A$ satisfies the conditions of Theorem~\ref{firstlater} and hence~$A$ is
modular.
\end{proof}

\begin{remark} If one  proved a version of Lemma~\ref{lem:liftingvariant}
  in the non-ordinary case, one could %
  improve the statement of
Theorem~\ref{theoremWisneeded}.  But note that in either case the required assumption on~$\rhobar$
is stronger than merely the assumption that~$\rhobar$ is finite flat, that is, arises as the generic fibre
of some~$\WW/\Z_2$ without any duality assumption. For example, if~$\rhobar |_{\Q_2}$ is unramified
with image of order~$5$, then~$\rhobar$ is both ordinary and finite flat,
and yet there does not exist any abelian surface~$A/\Z_2$ with~$A[2] \simeq \rhobar$. The issue
is that the only finite flat~$\WW/\Z_2$ with generic fibre~$\rhobar$
are either \'{e}tale or multiplicative, and so certainly not Cartier self-dual.
This is analogous to the fact that an unramified representation~$\rhobar: G_{\Q_2} \rightarrow \GL_2(\F_2)$
with image of order~$3$ is ordinary and finite flat in the usual sense but does not come from the~$2$-torsion
of an elliptic curve with good reduction, although after one extends the coefficients of~$\rhobar$ to~$\F_4$
it does come from the~$2$-torsion of an abelian surface with endomorphisms by~$\Z[(1+\sqrt{5})/2]$ with good reduction at~$2$
(for example, the modular abelian surface $J_0(23)/\Z_2$.) 
Note that this subtlety only arises (over~$\Q$) for~$p=2$, since for~$p>2$ any finite flat~$\WW/\Z_p$ is determined
by its generic fibre and so the Cartier self-duality of~$\WW/\Z_p$ follows from the corresponding property of~$\rhobar$. 
\end{remark}

\subsection{Consequences of Serre's Conjecture in regular weight}
\label{regularimpliesirregular}

The odd Artin conjecture for odd~$2$-dimensional
\emph{complex} representations of~$G_{\Q}$   is a %
consequence (\cite[Cor~10.2]{MR2551763})
of Serre's Conjecture for odd~$2$-dimensional mod-$p$ representations of~$G_{\Q}$
 (this implication was proved by Khare in~\cite{MR1434905}, using the weight
 lowering results~\cite{MR1074305,MR1185584}).
 It seems worthwhile remarking here that as a consequence of our main theorems,
 an  analogous deduction is valid for abelian surfaces over~$\Q$.

\begin{lemma}[Serre's Conjecture in regular weight implies modularity]\label{higherserre} Suppose that for every residual representation:
$$\rhobar: G_{\Q}  \rightarrow \GSp_4(\F_p)$$
satisfying the following conditions:
\begin{enumerate}
\item $\rhobar$ has multiplier~$\varepsilonbar^{-1}$,
\item %
  $\rhobar$ is absolutely irreducible,
\item the semi-simplification of $\rhobar |_{G_{\Qp}}$ is a direct sum of characters,
\end{enumerate}
there exists an ordinary cuspidal automorphic
           representation~$\pi$ of~$\GSp_4/\Q$ of regular weight, level prime to~$p$,
           and
           central character $|\cdot|^2$, such that \[\rhobar_{\pi,p}\cong\rhobar.\]
Then all abelian surfaces~$A/\Q$ are modular.
\end{lemma}

\begin{remark}
There are several possible natural variations on the hypotheses of this Lemma;
for example, one could only demand the statement for~$p$ sufficiently large. %
We have simply chosen one such version for illustrative purposes.
\end{remark}

\begin{proof}[Proof of Lemma~\ref{higherserre}]
By Theorem~\ref{remainingcases}, we may assume that~$A$ is ``challenging''
in the terminology of~\cite[\S9]{BCGP}, i.e.\ that either~$\End(A_{\Qbar})=\Z$,
or that~$A$ has Galois type~$\mathbf{B}[C_2]$, so there exists
a quadratic extension~$K/\Q$ so that~$E = \End(A_K) \otimes \Q
= \End(A_{\Qbar}) \otimes \Q$ is either~$\Q \oplus \Q$ or a real quadratic field.
By~\cite[Lemma~9.2.5]{BCGP}, there is a density one set of primes~$p>2$ 
 such that~$A$
is ordinary at~$p$ and residually~$p$-distinguished
in the
           sense of~\cite[Def~7.3.1]{BCGP},  and moreover that~$\rhobar =
           \rhobar_{A,p}$ satisfies the hypotheses listed in the statement of
           this lemma, as well as being vast in the sense of~\cite[Defn.\
           7.5.6]{BCGP} (and in particular reasonable in the sense of~\cite[Defn.\
          3.19]{Whitmore}) and  tidy  in the sense
           of~\cite[Defn.\
           7.5.11]{BCGP}. Furthermore, we may assume that if ~$A$ has Galois
           type~$\mathbf{B}[C_2]$ then ~$p$ splits in~$E$.
           
        We now deduce the modularity of~$A$ as a consequence of
        Theorem~\ref{thm:rho-is-modular-from-mult-one-and-classicity} for~$\rho=\rho_{A,p}$.
        It suffices to check the conditions of that theorem;  and in particular
         it suffices to check that hypotheses (1)--(5)
  and \ref{item: this needs a label}--\ref{item: same component} of
  Proposition~\ref{prop: the modularity result for p equals 2 or 3} hold (the other condition of Theorem~\ref{thm:rho-is-modular-from-mult-one-and-classicity} being immediate
        from the definition of ``challenging''). 
        Since~$A$  is an abelian surface with good ordinary reduction, and
        is additionally residually $p$-distinguished, the only conditions that
        need to be checked are that:
        \begin{enumerate}[(a)]
        \item \label{checkitenormous} $\rho_{A,p}(G_{\Q(\zeta_{p^{\infty}})})$ is integrally enormous,
        \item  \label{checkitregular} $\rhobar_{A,p}(G_{\Q})\setminus \Sp_4 (\Fp)$ contains a regular
          semi-simple element,
        \item \label{checkitcomponent}There are choices of  $p$-stabilizations
          such that $\rho_{\pi,p}|_{G_{\Q_p}}$ lies on a unique irreducible component of $\Spec R_p^\triangle$ and $\rho_{A,p}|_{G_{\Q_p}}$ lies on the same component.
         \end{enumerate}
        Suppose firstly that~ $A$ has Galois type $\mathbf{B}[C_2]$. %
         Then for sufficiently large primes~$p$ (splitting in~$E$)
         the mod~$p$ representations
         $$
         \begin{aligned}
         \rhobar_{A,p}: & \ G_{K(\zeta_{p^{\infty}})} \rightarrow
         \SL_2(\OL_E/p) = 
         \SL_2(\F_p) \times \SL_2(\F_p) \\
         \rhobar_{A,p}: & \ G_K \rightarrow \{(A,B) \in \GL_2(\F_p) \times \GL_2(\F_p), \ \det(A) = \det(B)\}
         \end{aligned}
         $$
          are surjective. (This
          follows from~\cite[Lem.\ 9.1.10(3)]{BCGP}, and for~$E=\Q \oplus \Q$
          goes back to~\cite{MR387283}.)
                     Thus we may additionally
         assume that~$p$ is chosen so that~$K \not\subseteq
	\Q(\zeta_{p^{\infty}})$ and
   	~$\rho_{A,p}(G_{\Q(\zeta_{p^{\infty}})})$ is
        precisely~$\SL_2(\F_p) \wr \Z/2\Z$.
	As explained in the proof of~\cite[Lemma~7.5.18]{BCGP}, for~$p \ge 3$, the
	 set~$\rhobar_{A,p}(G_{\Q(\zeta_{p^{\infty}})}) \setminus \rhobar(G_{G(\zeta_{p^{\infty}})})
	=  \SL_2(\F_p) \wr \Z/2\Z \setminus \SL_2(\F_p)^2$  always
	contains a regular semi-simple element %
        with eigenvalues~$(\zeta_8,\zeta_8^{-1},-\zeta_8 ,-\zeta_8 ^{-1})$.
	 Thus  condition~\ref{checkitenormous}  follows from Corollary~\ref{cor: still integrally enormous up the pro p tower}.
	 For~$p > 5$, $\rhobar_{A,p}(G_{\Q})$ contains~$(A,B)=(\diag(1,6),\diag(2,3))$,
	 which
	 is regular semi-simple and which does not lie in~$\Sp_4(\F_p)$,
	 which verifies
         condition~\ref{checkitregular} in this case as well.
         
	 Finally if~$\End(A_{\Qbar}) = \Z$, then  for sufficiently large~$p$ we
         have ~$\rhobar_{A,p}(G_{\Q})=\GSp_4(\F_p)$
	  and~$\rhobar_{A,p}(G_{\Q(\zeta_{p^{\infty}}})=\Sp_4(\F_p)$.
	  Hence both $\rhobar_{A,p}(G_{\Q}) \setminus \Sp_4(\F_p)$
	  and $\rhobar_{A,p}(G_{\Q(\zeta_{p^{\infty}}})$ contain regular semi-simple
	  elements for large enough~$p$ (for example, the same elements that were used above).
	  This verifies condition~\ref{checkitregular} in this case, and~$\rho_{A,p}(G_{\Q(\zeta_{p^{\infty}})})$ is
        integrally enormous by  Corollary~\ref{cor: still integrally enormous up
          the pro p tower}, verifying condition~\ref{checkitenormous}.

     For condition~\ref{checkitcomponent}, we firstly choose a
$p$-stabilization of~$\pi_p$, and thus
of~$\rhobar_{A,p}|_{G_{\Qp}}$. Then at least one of the
 $p$-stabilizations of~$\rho_{A,p}|_{G_{\Qp}}$ is compatible with this fixed
 choice. Having made this choice, since
     $\rhobar_{A,p}|_{G_{\Qp}}$ is residually $p$-distinguished it follows from
     \cite[Prop. 7.3.4]{BCGP} that $\Spec
     R_p^\triangle[1/p]$ is irreducible, and we are done.
	 \end{proof}

\bibliographystyle{amsalpha}
\bibliography{modularityabeliansurfaces}

\end{document}